# Stability Analysis and Stabilization of Continuous-Time Linear Systems with Distributed Delays

*Qian Feng*,
School of Control and Computer Engineering,
North China Electric Power University Beijing, China.

*Sing Kiong Nugnag*†,
The Department of Electrical and Computer Engineering
The University of Auckland, New Zealand.

*Alexandre Seuret*
Departamento de Ingeniería de
Sistemas y Automática, Universidad de Sevilla

# Preface

This monograph is devoted to the exploration of new methods for the stability and stabilization of linear systems with nontrivial distributed delays through the construction of Lyapunov-Krasovskiĭ functionals. We begin by proposing methods for the design of a dissipative state feedback controller for a linear distributed delay system in which the delay is known and the distributed delay kernels belong to a specific class of functions. The stabilization problem is addressed by constructing a functional associated with the delay kernels of the system through the use of an integral inequality of general form. We subsequently extend these results to accommodate uncertain linear distributed delay systems, where the embedded linear fractional uncertainties are handled by a newly proposed lemma. Given the crucial role that integral inequalities have played in applying the Lyapunov-Krasovskiĭ functional approach, we introduce two general classes of novel integral inequalities, establishing their relations and properties among their lower bounds. These proposed inequalities have significant general structures that can generalize a wide variety of existing integral inequalities found in delay-related literature. Next, we propose a new method for the dissipativity and stability analysis of linear coupled differential-difference systems, where the distributed delays can contain an unlimited number of $\mathcal{L}^2$ functions. The delay kernels are approximated by a class of functions that are the solution to linear homogeneous differential equations with constant coefficients. In addition, approximate errors are included in the resulting dissipativity (stability) condition via a matrix framework based on the application of a novel proposed integral inequality via the construction of an LKF. The previous results are then followed by a study of delay range analysis where the problem of dissipativity and stability analysis of coupled differential-difference systems with distributed delays is considered with an unknown but bounded delay. By constructing a functional whose matrix parameters are dependent polynomially on the delay value, we can derive a dissipativity and stability condition as sum-of-squares constraints. Finally, we propose new methods for the dissipative synthesis of linear systems with nontrivial distributed delays where the delay is time-varying, measurable and bounded. The problem is solved by utilizing the Lyapunov-Krasovskiĭ functional approach based on the application of a novel integral inequality.

ΙΣ ΧΣ ΘΥ ΥΣ ΣΗΡ

# List of Notation and Acronyms

The following symbols are utilized throughout this book that remain consistent and unambiguous. However, it is important to note that the notation of scalars, functions, and matrices may be overloaded in distinct chapters without introducing ambiguities.

## Symbols:

universal quantifier $\forall$

existential quantifier $\exists$

unique existential quantifier $\exists!$

$\mathbb{N} := \{1, 2, 3 \cdots\}$

$\mathbb{N}_0 := \mathbb{N} \cup \{0\}$

$\mathbb{R} := \{\text{All Real Numbers}\}$

$\mathbb{R}_{\geq a} := \{x \in \mathbb{R} : x \geq a\}$

$\mathbb{C} := \{\text{All Complex Numbers}\}$

$\mathbb{R}^{n\times m} := \{\text{All Real } n \times m \text{ Matrices}\}$

$\mathbb{C}^{n\times m} := \{\text{All Complex } n \times m \text{ Matrices}\}$

$\mathbb{S}^n := \{X \in \mathbb{R}^{n\times n} : X = X^\top\}$

$\mathbb{R}^{n\times m}_{[r]} := \{X \in \mathbb{R}^{n\times m} : \text{rank}(X) = r\}$

$\mathsf{Sy}(X) := X + X^\top$ where $X \in \mathbb{R}^{n\times n}$

$\mathcal{Y}^{\mathcal{X}} := \{f(\cdot) : f(\cdot) \text{ is a function from } \mathcal{X} \text{ onto } \mathcal{Y}\}$

$\boldsymbol{\mathcal{L}}(\mathcal{X}) = \{\mathcal{Y} : \mathcal{Y} \subseteq \mathcal{X} \;\; \& \;\; \mathcal{Y} \text{ is Lebesgue measurable}\}$

$\boldsymbol{\mathcal{B}}(\mathbb{R}^m)$: Borel Algebra on $\mathbb{R}^m$

$\mathcal{M}(\mathbb{R}^n; \mathbb{R}^m) := \left\{f(\cdot) \in (\mathbb{R}^m)^{\mathbb{R}^n} : \forall \mathcal{Y} \in \boldsymbol{\mathcal{B}}(\mathbb{R}^m),\; f^{-1}(\mathcal{Y}) \in \boldsymbol{\mathcal{L}}(\mathbb{R}^n)\right\}$

$\mathbb{R}^n \ni \mathbf{x} \to \|\mathbf{x}\|_p := \left(\sum_{i=1}^n |x_i|^p\right)^{\frac{1}{p}}$

$\|f(\cdot)\|_p := \left(\int_{\mathcal{X}} |f(x)|^p \mathrm{d}x\right)^{\frac{1}{p}}$ $\;\;p$–seminorm of $f(\cdot)$

$\|\boldsymbol{f}(\cdot)\|_p := \left(\int_{\mathcal{X}} \|\boldsymbol{f}(x)\|_2^p \mathrm{d}x\right)^{\frac{1}{p}}$ $p$–seminorm of $\boldsymbol{f}(\cdot)$

$\|\boldsymbol{f}(\cdot)\|_\infty := \sup_{x\in\mathcal{X}} \|\boldsymbol{f}(x)\|_2$

$\mathcal{L}^p(\mathcal{X};\mathbb{R}) := \{f(\cdot) \in \mathcal{M}(\mathcal{X};\mathbb{R}) : \|f(\cdot)\|_p < +\infty\}$

$\mathcal{L}^p(\mathcal{X};\mathbb{R}^n) := \{\boldsymbol{f}(\cdot) \in \mathcal{M}(\mathcal{X};\mathbb{R}^n) : \|\boldsymbol{f}(\cdot)\|_p < +\infty\}$

$\widehat{\mathcal{L}}_p(\mathcal{X};\mathbb{R}^n)$ The local integrable version of $\mathcal{L}^p(\mathcal{X};\mathbb{R}^n)$

$\mathcal{C}(\mathcal{X};\mathbb{R}^n) := \left\{\boldsymbol{f}(\cdot) \in (\mathbb{R}^n)^{\mathcal{X}} : \boldsymbol{f}(\cdot) \text{ is continuous on } \mathcal{X}\right\}$

$\mathcal{C}^k([a,b];\mathbb{R}^n) := \left\{\boldsymbol{f}(\cdot) \in \mathcal{C}([a,b];\mathbb{R}^n) : \frac{\mathsf{d}^k \boldsymbol{f}(x)}{\mathsf{d}x^k} \in \mathcal{C}([a,b];\mathbb{R}^n)\right\}$, the derivatives at $a$ and $b$ are one sided.

$\gamma(\cdot)$ Gamma function

$\mathbf{Col}_{i=1}^n X_i := \left[\mathbf{Row}_{i=1}^n X_i^\top\right]^\top = \left[X_1^\top \cdots X_i^\top \cdots X_n^\top\right]^\top$

$[*]YX = X^\top YX$ or $X^\top Y[*] = X^\top YX$

$\mathsf{O}_{n,m}$ $n \times m$ zero matrix

$\mathbf{vec}(A) = \mathbf{vec}\left(\mathbf{Row}_{i=1}^m(\mathbf{a}_i)\right) := \mathbf{Col}_{i=1}^m \mathbf{a}_i,$ vectorization of $A = \mathbf{Row}_{i=1}^m \mathbf{a}_i \in \mathbb{R}^{n\times m}$

$\mathsf{O}_n$ $n \times n$ zero matrix

$\mathbf{0}_n$ zeros column vector with dimension $n$.

$x \vee y = \max(x, y)$

$x \wedge y = \min(x, y)$

$X \oplus Y := \begin{bmatrix} X & \mathsf{O} \\ \mathsf{O} & Y \end{bmatrix}$

$\bigoplus_{i=1}^n X_i = \begin{bmatrix} X_1 & \mathsf{O} & \cdots & \mathsf{O} \\ \mathsf{O} & X_2 & \cdots & \mathsf{O} \\ \vdots & \vdots & \ddots & \vdots \\ \mathsf{O} & \mathsf{O} & \mathsf{O} & X_n \end{bmatrix}$

$[*]YX := X^\top YX$ or $X^\top Y[*] := X^\top YX$

$\otimes$ the Kronecker product

$A \succeq 0 \;:\Longleftrightarrow\; \forall \mathbf{x} \in \mathbb{R}^n \setminus \{\mathbf{0}_n\},\;\; \mathbf{x}^\top A\mathbf{x} \geq 0, \quad A = A^\top$

$A \preceq 0 \;:\Longleftrightarrow\; \forall \mathbf{x} \in \mathbb{R}^n \setminus \{\mathbf{0}_n\},\;\; \mathbf{x}^\top A\mathbf{x} \leq 0, \quad A = A^\top$

$X \succeq Y \;:\Longleftrightarrow\; X - Y \succeq 0$

$X \preceq Y \;:\Longleftrightarrow\; X - Y \preceq 0$

$A \succ 0 \;:\Longleftrightarrow\; \forall \mathbf{x} \in \mathbb{R}^n \setminus \{\mathbf{0}_n\},\;\; \mathbf{x}^\top A\mathbf{x} > 0, \quad A = A^\top$

$A \prec 0 \;:\Longleftrightarrow\; \forall \mathbf{x} \in \mathbb{R}^n \setminus \{\mathbf{0}_n\},\;\; \mathbf{x}^\top A\mathbf{x} < 0, \quad A = A^\top$

$X \succ Y \;:\Longleftrightarrow\; X - Y \succ 0$

$X \prec Y \;:\Longleftrightarrow\; X - Y \prec 0$

$\mathbb{S}^n_{\preceq 0} := \{X \in \mathbb{S}^n : X \preceq 0\}$

$\mathbb{S}^n_{\prec 0} := \{X \in \mathbb{S}^n : X \prec 0\}$

$\mathbb{S}^n_{\succeq 0} := \{X \in \mathbb{S}^n : X \succeq 0\}$

$\mathbb{S}^n_{\succ 0} := \{X \in \mathbb{S}^n : X \succ 0\}$

The symbol $*$ is used as abbreviations for $[*]YX = X^\top YX$ or $X^\top Y[*] = X^\top YX$ or $\left[\begin{smallmatrix} A & B \\ * & C \end{smallmatrix}\right] = \left[\begin{smallmatrix} A & B \\ B^\top & C \end{smallmatrix}\right]$. We use $\sqrt{X}$ to represent[1] the unique square root of $X \succ 0$. The order of matrix operations in this chapter is *matrix (scalars) multiplications* $> \otimes > \oplus > +$. Finally, the notation [ ], is used to represent empty matrices [1, See I.7] following the same definition and rules in MATLAB© . We assume $I_0 = [\,]_{0,0}$, $\mathsf{O}_{0,m} = [\,]_{0,m}$ and $\mathbf{Col}_{i=1}^n x_i = [\,]_{0,m}$, $\mathbf{Row}_{i=1}^n x_i = [\,]_{m,0}$ if $n < 1$, where $[\,]_{0,m}$, $[\,]_{m,0}$ are an empty matrices with $m \in \mathbb{N}$ depending on specific contexts.

## Acronyms:

DDs: Distributed Delays

LKFs: Lyapunov-Krasovskiĭ Functionals

CLKFs: Complete Lyapunov-Krasovskiĭ Functionals

FDEs: Functional Differential Equations

TDSs: Time-Delay Systems

LTDSs: Linear Time-Delay Systems

LDDSs: Linear Distributed Delay System

SA: Spectral Abscissa

LMIs: Linear Matrix Inequalities

SDP: Semidefinite Programming

BMIs: Bilinear Matrix Inequality

SoSs: Sums-of-Squares

SRF: Supply Rate Function

IIs: Integral Inequalities

LBs: Lower Bounds

[1]Note that $\sqrt{X^{-1}} = \left(\sqrt{X}\right)^{-1}$ for any $X \succ 0$ based on the application of eigendecomposition of $X \succ 0$

# Chapter 1

# Introduction

## 1.1 Background

Time-delay systems (TDSs) [2–6] are instrumental in characterizing real-time processes influenced by transport, propagation, or aftereffects [5, 7, 8]. TDSs can be mathematically represented via functional differential equations (FDEs) [9–13], coupled differential-difference equations [14–22], or posed as infinite-dimensional systems (IDSs) [23–33]. Practical and engineering examples[1] of systems with delays can be found in the modeling of

- mechanical vibration [35–38], industrial processes [39–41], drilling [42–45]
- robotic systems [46–51], UAV [52–56], vehicle systems [57–61],
- biological processes [62–67], cell dynamics [68–70], SIR epidemic model [71–75],
- financial market [76–80], supply chain management [81–85],
- multi-agents systems [86–97], network control systems [98–101],
- smart grid [102–105], power systems [106–111], renewable energy systems [112–116].

Delays in real-time engineering systems can arise from actuators, sensors, wires, or communication channels. It has been shown that the presence of delay in a system can result in positive [117–121] or negative [122, 123] system behaviors. As a result, the analysis and control of TDSs remain a vibrant research area in system engineering.

Generally speaking, two types of delays, pointwise and distributed, are commonly utilized to model aftereffects in dynamical systems. Pointwise delays can be modeled by classical transport equations with appropriate boundary conditions [5], which can intuitively explain the connections between delays and transport phenomena, and why systems with delays may have infinite dimension [23, 33, 124]. Delays can also be introduced by media of propagation with more complex structures. Distributed delay (DD) [125] is represented by an integral operator integrating over a delay interval, which takes into account a segment of past dynamics information. This indicates that DDs contain more delay information than a simple pointwise-delay but also makes them more difficult to analyze due to the various forms an integral operator can take. Systems with DDs can be found in the modeling of biological processes [64, 71, 72, 126–130], population dynamics [131–138], traffic flows [139–141], complex and neural networks [142–152], modeling actuator saturation [153, 154], characterizing event-triggered mechanisms [100, 114, 115], matching dynamics [155, 156], machine tool vibration [36], shimmy dynamics of wheels [157], chemical kinetics [39, 40] and milling processes [158]. Distributed delays can also exist at the input of a system [113, 159, 160], or within controllers [161–165], which may produce positive effects [165] on the behavior of the system's dynamics.

[1] For more examples, see [34] and [10]

Stability is crucial in qualitatively describing the asymptotic behavior of system dynamics. The characteristics and delay values in a system [166] can fundamentally affect its stability. Even for linear time delay systems (LTDSs) [8], stability analysis and stabilization require the use of sophisticated mathematical tools, in contrast to linear time-invariant systems with finite dimensions where achieving these objectives is more straightforward by constructing a quadratic Lyapunov function [167]. Moreover, the presence of DDs in an LTDS may further complicate stability analysis and stabilization, making it of great research interest to propose effective solutions to these problems.

This book is devoted to the presentation of the latest approaches for stability analysis and stabilization of linear systems with DDs that possess finite delay values and nontrivial kernels, considering dissipativity constraints. Specifically, the DDs considered in this book are addressed using various approaches, and the resulting synthesis conditions are derived from the construction of Lyapunov-Krasovskiĭ functionals (LKFs) [9, Section 5.2]. The remainder of this chapter introduces mathematical preliminaries concerning systems with delays, emphasizing the Krasovskiĭ stability criterion [168]. A literature review section follows, discussing existing methods for the stability analysis and stabilization of LTDSs based on the construction of LKFs functionals. Finally, the research motivations and an outline of subsequent chapters are provided. Some content in this book is based on or related to the results reported in [169–172, 539? ].

## 1.2 A Review in the Stability Analysis of Systems with Delays

In this section, we investigate general models for systems with delays, based on the framework of coupled differential-functional equations (CDFEs). CDFEs are capable of characterizing various types of systems affected by delay effects, including those denoted by retarded and neutral type FDEs [9, 10, 173]. We also mathematically present the Lyapunov-Krasovskiĭ stability criteria for CDFEs, as a generalization of the Lyapunov direct method for ordinary differential equations [174].

It should be noted that this chapter does not aim to provide a comprehensive and systematic study of the mathematical theories of coupled differential-functional equations or FDEs. For further information on the extensive theories underlying this topic, readers are referred to [9, 22, 175]. Additionally, the monograph [176] provides an excellent treatise on the theoretical fundamentals of dynamical systems and mathematical control theories.

### 1.2.1 Systems with delays denoted by Coupled Differential-Functional Equations

A general coupled differential-functional equation [22] is denoted by

$$\begin{aligned} &\dot{\boldsymbol{x}}(t) := \lim_{\eta \downarrow 0} \frac{\boldsymbol{x}(t+\eta) - \boldsymbol{x}(t)}{\eta} = \boldsymbol{f}\left(t, \boldsymbol{x}(t), \mathbf{y}_t(\cdot)\right), \quad \boldsymbol{y}(t) = \boldsymbol{g}(t, \boldsymbol{x}(t), \mathbf{y}_t(\cdot)), \quad t \geq t_0 \in \mathbb{R} \\ &\boldsymbol{x}(t_0) = \boldsymbol{\omega} \in \mathbb{R}^n, \;\; \forall \theta \in [-r, 0), \;\; \boldsymbol{y}(t_0 + \theta) = \mathbf{y}_{t_0}(\theta) = \boldsymbol{\psi}(\theta), \;\; \boldsymbol{\psi}(\cdot) \in \widehat{\mathcal{C}}\left([-r, 0); \mathbb{R}^\nu\right) \end{aligned} \tag{1.1a}$$

where functionals $\boldsymbol{f} : \mathbb{R} \times \mathbb{R}^n \times \widehat{\mathcal{C}}([-r, 0); \mathbb{R}^n) \to \mathbb{R}^n$ and $\boldsymbol{g} : \mathbb{R} \times \mathbb{R}^n \times \widehat{\mathcal{C}}([-r, 0); \mathbb{R}^\nu) \to \mathbb{R}^\nu$ are continuous and satisfy

$$\forall t \in \mathbb{R}, \; \mathbf{0}_n = \boldsymbol{f}\left(t, \mathbf{0}_n, \mathbf{0}_\nu\right), \quad \mathbf{0}_\nu = \boldsymbol{g}(t, \mathbf{0}_n, \mathbf{0}_\nu), \tag{1.1b}$$

and $\mathbf{y}_t(\cdot) \in \widehat{\mathcal{C}}\,([-r,0);\mathbb{R}^{\nu})$ is defined as

$$\forall t \geq t_0, \;\; \forall \theta \in [-r,0), \;\; \mathbf{y}_t(\theta) = \boldsymbol{y}(t+\theta). \tag{1.1c}$$

Notation $\widehat{\mathcal{C}}([-r,0);\mathbb{R}^{\nu})$ denotes the Banach space of bounded right piecewise continuous functions endowed with a uniform norm $\|\boldsymbol{f}(\cdot)\|_{\infty} := \sup_{\tau\in[-r,0)} \|\boldsymbol{f}(\tau)\|_2$. As detailed in [22, Section 2], the uniqueness and existence of the solution to (1.1a) can be established via the procedures in [177], provided that certain properties are satisfied by the functional $\boldsymbol{f} : \mathbb{R} \times \mathbb{R}^n \times \widehat{\mathcal{C}}([-r,0);\mathbb{R}^n) \to \mathbb{R}^n$ and $\boldsymbol{g} : \mathbb{R} \times \mathbb{R}^n \times \widehat{\mathcal{C}}([-r,0);\mathbb{R}^{\nu}) \to \mathbb{R}^{\nu}$.

Now let $\boldsymbol{g}(t,\boldsymbol{x}(t),\mathbf{y}_t(\cdot)) := \boldsymbol{x}(t)$ with $n = \nu$ which satisfies (1.1b), then system (1.1a) becomes

$$\begin{aligned} &\dot{\boldsymbol{x}}(t) = \boldsymbol{f}(t,\boldsymbol{x}(t),\mathbf{x}_t(\cdot)), \\ &\boldsymbol{x}(t_0) = \boldsymbol{\omega} \in \mathbb{R}^n, \;\; \forall \theta \in [-r,0), \;\; \boldsymbol{x}(t_0+\theta) = \mathbf{x}_{t_0}(\theta) = \boldsymbol{\psi}(\theta), \end{aligned} \tag{1.2a}$$

where $\mathbf{x}_t(\cdot) \in \widehat{\mathcal{C}}\,([-r,0);\mathbb{R}^n)$ is defined as

$$\forall t \geq t_0, \;\; \forall \theta \in [-r,0), \;\; \mathbf{x}_t(\theta) = \boldsymbol{x}(t+\theta). \tag{1.2b}$$

The differential equation in (1.2a) can be further simplified as

$$\forall t \geq t_0 \in \mathbb{R}, \;\; \dot{\boldsymbol{\xi}}(t) = \widehat{\boldsymbol{f}}\,(t,\boldsymbol{\xi}_t(\cdot))\,, \;\; \forall \theta \in [-r,0], \;\; \boldsymbol{\xi}(t_0+\theta) = \boldsymbol{\phi}(\theta), \;\; \boldsymbol{\phi}(\cdot) \in \mathcal{C}([-r,0];\mathbb{R}^n) \tag{1.3a}$$

where $\widehat{\boldsymbol{f}} : \mathbb{R} \times \mathbb{R}^n \times \mathcal{C}([-r,0];\mathbb{R}^n) \to \mathbb{R}^n$ is continuous and $\boldsymbol{\xi}_t(\cdot)$ is defined as

$$\forall t \geq t_0, \;\; \forall \theta \in [-r,0], \;\; \boldsymbol{\xi}_t(\theta) = \boldsymbol{\xi}(t+\theta). \tag{1.3b}$$

It should be noted that in equation (1.3a) we have chosen a continuous function $\boldsymbol{\phi}(\cdot) \in \mathcal{C}([-r,0];\mathbb{R}^n)$ as the initial condition of (1.3a), which is a simpler choice compared to the piecewise continuous function $\boldsymbol{\psi}(\cdot)$ in (1.2a). Equation (1.3a) represents the standard expression for general FDEs, as described in references [9, 10]. Additionally, let us consider the following neutral FDE [9, 10]

$$\begin{aligned} &\dot{\boldsymbol{x}}(t) = \acute{\boldsymbol{f}}(t,\mathbf{y}_t(\cdot)), \quad \boldsymbol{x}(t) = \boldsymbol{y}(t) - \boldsymbol{h}(t,\mathbf{y}_t(\cdot)), \quad t \geq t_0 \\ &\forall \theta \in [-r,0], \;\; \mathbf{y}_{t_0}(\theta) = \boldsymbol{\phi}(\theta) \end{aligned} \tag{1.4a}$$

with $t_0 \in \mathbb{R}$ and $\boldsymbol{\phi}(\cdot) \in \mathcal{C}([-r,0];\mathbb{R}^n)$ and $\mathbf{y}_t(\cdot)$ defined as

$$\forall t \geq t_0, \;\; \forall \theta \in [-r,0], \;\; \mathbf{y}_t(\theta) = \boldsymbol{y}(t+\theta), \tag{1.4b}$$

where $\acute{\boldsymbol{f}} : \mathbb{R} \times \mathcal{C}([-r,0];\mathbb{R}^n) \to \mathbb{R}^{\nu}$ and $\boldsymbol{h} : \mathbb{R} \times \mathcal{C}([-r,0];\mathbb{R}^n) \to \mathbb{R}^n$ are continuous and satisfy

$$\forall t \in \mathbb{R}, \;\; \mathbf{0}_n = \acute{\boldsymbol{f}}(t,\mathbf{0}_n), \quad \mathbf{0}_n = \boldsymbol{h}(t,\mathbf{0}_n). \tag{1.4c}$$

Now (1.4a) can be further reformulated as

$$\begin{aligned} &\dot{\boldsymbol{x}}(t) = \grave{\boldsymbol{f}}\,(t,\boldsymbol{y}(t),\acute{\mathbf{y}}_t(\cdot)) = \grave{\boldsymbol{f}}\left(t,\left[\boldsymbol{x}(t)+\boldsymbol{h}(t,\mathbf{y}_t(\cdot))\right],\acute{\mathbf{y}}_t(\cdot)\right) = \widetilde{\boldsymbol{f}}\,(t,\boldsymbol{x}(t),\acute{\mathbf{y}}_t(\cdot))\,, \\ &\boldsymbol{y}(t) = \boldsymbol{x}(t) + \boldsymbol{h}(t,\mathbf{y}_t(\cdot)) = \widetilde{\boldsymbol{h}}\,(t,\boldsymbol{x}(t),\acute{\mathbf{y}}_t(\cdot))\,, \\ &\mathbb{R}^n \ni \boldsymbol{x}(t_0) = \boldsymbol{\phi}(0) - \boldsymbol{h}(t_0,\boldsymbol{\phi}(\cdot)), \;\; \forall \theta \in [-r,0), \;\; \acute{\mathbf{y}}_{t_0}(\theta) = \boldsymbol{\phi}(\theta) \end{aligned} \tag{1.5a}$$

where $\acute{\mathbf{y}}_t(\cdot)$ satisfies

$$\forall t \geq t_0, \;\; \forall \theta \in [-r,0), \;\; \acute{\mathbf{y}}_t(\theta) = \boldsymbol{y}(t+\theta), \tag{1.5b}$$

and $\widetilde{\boldsymbol{f}} : \mathbb{R} \times \mathbb{R}^n \times \mathcal{C}([-r,0);\mathbb{R}^n) \to \mathbb{R}^n$ and $\widetilde{\boldsymbol{h}} : \mathbb{R} \times \mathbb{R}^n \times \mathcal{C}([-r,0);\mathbb{R}^\nu) \to \mathbb{R}^n$ are continuous and satisfy $\forall t \in \mathbb{R},\ \mathbf{0} = \widetilde{\boldsymbol{h}}(t,\mathbf{0},\mathbf{0}),\ \mathbf{0} = \widetilde{\boldsymbol{f}}(t,\mathbf{0},\mathbf{0})$. Upon closer examination, it becomes evident that (1.5a) can be effectively analyzed through (1.1a). This demonstrates the versatility of (1.1a) in its ability to model neutral FDEs.

In the field of control engineering, determining the stability of dynamical systems is a fundamental property of crucial importance. In the subsequent subsection, we will present the Lyapunov direct approach for CDFEs, which serves to generalize the Lyapunov's second method for FDEs.

### 1.2.2 Stability analysis of CDFEs via the LKF criterion

The subsequent Lyapunov-Krasovskiĭ stability theorem for CDFEs is derived from [22, Theorem 3], and has been paraphrased using our established notation. For a more comprehensive understanding of the precise mathematical definitions pertaining to the various types of stability, we recommend consulting [22, Definitions 1 and 2].

**Theorem 1.1.** *Let $\boldsymbol{f} : \mathbb{R} \times \mathbb{R}^n \times \widehat{\mathcal{C}}([-r,0);\mathbb{R}^n) \to \mathbb{R}^\nu$ and $\boldsymbol{g} : \mathbb{R} \times \mathbb{R}^n \times \widehat{\mathcal{C}}([-r,0);\mathbb{R}^\nu) \to \mathbb{R}^\nu$ in (1.1a) to satisfy the prerequisites of [22, Theorem 3], and assume that $\boldsymbol{y}(t) = \boldsymbol{g}(t,\boldsymbol{x}(t),\mathbf{y}_t(\cdot))$ in (1.1a) is uniformly input to state stable[2]. Moreover, let $\alpha_1(\cdot),\alpha_2(\cdot),\alpha_3(\cdot) \in \mathcal{C}(\mathbb{R}_{\geq 0};\mathbb{R}_{\geq 0})$ be non-decreasing function and $\forall \theta > 0,\ \alpha_1(\theta) > 0,\ \alpha_2(\theta) > 0$ with $\alpha_1(0) = \alpha_2(0) = 0$. Then the trivial solution $\boldsymbol{x}(t) \equiv \mathbf{0}_n, \boldsymbol{y}(t) \equiv \mathbf{0}_\nu$ of (1.1a) is uniformly stable, if there exist differentiable functionals $\mathsf{v} : \mathbb{R} \times \mathbb{R}^n \times \widehat{\mathcal{C}}([-r,0);\mathbb{R}^\nu) \to \mathbb{R}_{\geq 0}$ such that $\forall t \in \mathbb{R}, \mathsf{v}(t,\mathbf{0}_n,\mathbf{0}_\nu) = 0$ and*

$$\alpha_1\left(\|\boldsymbol{\omega}\|_2\right) \leq \mathsf{v}(t_0,\boldsymbol{\omega},\boldsymbol{\phi}(\cdot)) \leq \alpha_2\left(\|\boldsymbol{\omega}\|_2 \vee \|\boldsymbol{\phi}(\cdot)\|_\infty\right) \tag{1.6a}$$

$$\dot{\mathsf{v}}(t_0,\boldsymbol{\omega},\boldsymbol{\phi}(\cdot)) = \frac{\mathsf{d}^+}{\mathsf{d}t}\mathsf{v}(t,\boldsymbol{x}(t),\mathbf{y}_t(\cdot))\Big|_{t=t_0,\boldsymbol{x}(t_0)=\boldsymbol{\omega},\mathbf{y}_{t_0}(\cdot)=\boldsymbol{\phi}(\cdot)} \leq -\alpha_3\left(\|\boldsymbol{\omega}\|_2\right) \tag{1.6b}$$

*for any initial condition $\boldsymbol{\omega} \in \mathbb{R}^n$ and $\boldsymbol{\phi}(\cdot) \in \widehat{\mathcal{C}}([-r,0);\mathbb{R}^n)$ in (1.1a), where $t_0 \in \mathbb{R}$, $\|\boldsymbol{\omega}\|_2 \vee \|\boldsymbol{\phi}(\cdot)\|_\infty := \max\left(\|\boldsymbol{\omega}\|_2, \|\boldsymbol{\phi}(\cdot)\|_\infty\right)$ and $\frac{\mathsf{d}^+}{\mathsf{d}x}f(x) = \limsup_{\eta\downarrow 0}\frac{f(x+\eta)-f(x)}{\eta}$. Furthermore, if $\alpha_3(\theta) > 0$ for all $\theta > 0$, then the trivial solution $\boldsymbol{x}(t) \equiv \mathbf{0}_n$, $\boldsymbol{y}(t) \equiv \mathbf{0}_\nu$ of (1.1a) is globally uniformly asymptotically stable. In addition, the trivial solution $\boldsymbol{x}(t) \equiv \mathbf{0}_n$, $\boldsymbol{y}(t) \equiv \mathbf{0}_\nu$ of (1.1a) is globally uniformly asymptotically stable if $\lim_{\theta\to\infty}\alpha_1(\theta) = \infty$.*

It is important to note that if $\widehat{\boldsymbol{x}}(t)$ and $\widehat{\boldsymbol{y}}(t)$ are any solution to (1.1a), then their stability can be determined by the trivial solution $\boldsymbol{z}(t) \equiv \mathbf{0}_n, \boldsymbol{\zeta}(t) \equiv \mathbf{0}_\nu$ of

$$\begin{aligned}
\dot{\boldsymbol{z}}(t) &= \boldsymbol{f}\Big(t,\boldsymbol{z}(t)+\widehat{\boldsymbol{x}}(t),\mathbf{z}_t(\cdot)+\widehat{\mathbf{y}}_t(\cdot)\Big) - \boldsymbol{f}\left(t,\widehat{\boldsymbol{x}}(t),\widehat{\mathbf{y}}_t(\cdot)\right), \\
\boldsymbol{\zeta}(t) &= \boldsymbol{g}\Big(t,\boldsymbol{z}(t)+\widehat{\boldsymbol{x}}(t),\mathbf{z}_t(\cdot)+\widehat{\mathbf{y}}_t(\cdot)\Big) - \boldsymbol{g}(t,\widehat{\boldsymbol{x}}(t),\widehat{\mathbf{y}}_t(\cdot)),
\end{aligned} \tag{1.7}$$

where $\mathbf{z}_t(\cdot)$ is defined as $\forall t \geq t_0,\ \ \forall \theta \in [-r,0),\ \ \mathbf{z}_t(\theta) = \boldsymbol{z}(t+\theta)$.

In the case of a system described by equations (1.2a)–(1.3b), Theorem 1.1 can be modified to specifically address (1.3a), where Theorem 1.1 becomes equivalent to [9, Theorem 2.1] or [178, Theorem 1.3]. Furthermore, this adaptation illustrates that the stability of certain types of neutral FDEs [44] can be effectively analyzed through (1.1) with Theorem 1.1. This is demonstrated by the procedures outlined in (1.4a)–(1.5a).

[2] For the mathematical definition of the uniformly input to state stability of $\boldsymbol{y}(t) = \boldsymbol{g}(t,\boldsymbol{x}(t),\mathbf{y}_t(\cdot))$, see [22, Definition 2]

**Remark 1.1.** Theorem 1.1 can be viewed as an extension of the direct Lyapunov method for systems with finite dimensions. It offers a powerful means of verifying the stability of (1.1a) without the need for explicit knowledge of the analytic expressions of its solution.

Within the field of control engineering, it is often desirable to numerically construct Lyapunov functions or functionals, where conditions such as those presented in (1.6a) and (1.6b) are implied by the feasible solution to specific optimization programs. In this book, we focus on the synthesis and stability analysis of linear systems with DDs. In contrast to employing a quadratic Lyapunov function to analyze the stability of an LTI delay-free system, sufficient and necessary stability conditions cannot typically be derived through the construction of LKFs for linear systems with delays. In the ensuing section, we will examine existing methods for the stability analysis and stabilization of LTDSs, with a specific focus on linear systems with DDs. We will also discuss some recent advancements in methodologies and expound upon the technical difficulties inherent to this topic.

## 1.3 Review of Stability Analysis Methods for Linear Systems with DDs

In this section, we present a comprehensive review of two primary branches of existing methodologies for the stability analysis and stabilization of LTDSs: frequency and time-domain approaches. Particular emphasis is placed on linear systems with DDs. The scenarios proposed in this treatise are predicated upon the construction of LKFs, which fall within the purview of time-domain approaches. Additionally, we provide a concise overview of the development of optimization methods via semidefinite programming (SDP), through which our proposed stability and stabilization conditions can be numerically computed.

In order to provide a comprehensive description of existing works pertaining to the stability analysis and stabilization of linear systems with DDs, we utilize the following linear distributed delay system (LDDS)

$$\dot{\boldsymbol{x}}(t) = A_1\boldsymbol{x}(t) + A_2\boldsymbol{x}(t-r) + \int_{-r}^{0} A_3(\tau)\boldsymbol{x}(t+\tau)\mathrm{d}\tau, \quad r > 0 \tag{1.8}$$

as a reference in time-domain, where $A_3(\cdot) \in \mathcal{L}^2([-r,0];\mathbb{R}^{n\times n})$ with $r > 0$ and $\boldsymbol{x}(t) \in \mathbb{R}^n$. The spectrum of (1.8):

$$\left\{s \in \mathbb{C} : p(s) = 0\right\}, \quad p(s) = \det\left(sI_n - A_1 - A_2\mathrm{e}^{-rs} - \int_{-r}^{0} A_3(\tau)\mathrm{e}^{\tau s}\mathrm{d}\tau\right) \tag{1.9}$$

is used as a reference for the discussion of existing frequency-domain-based methods.

### 1.3.1 Frequency-domain approaches

To analyze the spectrum in (1.9), it is necessary to possess at least some knowledge concerning the zeros of $p(\cdot)$, which contains complex transcendental functions. Unlike the spectrum of finite-dimensional LTI system, which have a finite number of characteristic roots, the number of roots of an LTDS is typically countably infinite [178, 179]. Thus, more advanced mathematical theories are required for the analysis. It has been proved that (see [178, Theorem 1.5]) (1.8) is asymptotically (exponentially) stable if and only if the spectral abscissa is negative. The spectral abscissa can be computed by finding the zeros of $p(s)$ within a critical region in the complex plane, as only a finite number of zeros [8] can exist in the region $\{s \in \mathbb{C} : \Re(s) > 0\}$. As a result, the stability of (1.8) can be determined by computing [8] the spectral abscissa of $\{s \in \mathbb{C} : p(s) = 0\}$.

Numerous methodologies have been developed for computing spectral abscissae, including those for finding the zeros of analytic functions [180–186], and algorithms for solving nonlinear eigenvalue problems [187–201]. Notably, several approaches with available MATLAB© code, can handle a general[3] DD $\int_{-r}^{0} A_3(\tau)e^{\tau s}d\tau$ with $s \in \mathbb{C}$. For example, the numerical method proposed in [202–205] (MATLAB© code in http://cdlab.uniud.it/software#eigAM-eigTMN) can compute the spectral abscissa of retarded (renewal) linear systems with an unlimited number of pointwise and distributed delays, provided that the expressions of DDs and delay values are known. Additionally, algorithms for computing the zeros of quasi-polynomials [206–212] have been proposed in [213–215], with accompanying MATLAB© package. These algorithms can be used to compute $\max\{\Re(s) : p(s) = 0\}$, as $p(s)$ is a quasi-polynomial in $s \in \mathbb{C}$. However, they may not be applicable if $\int_{-r}^{0} A_3(\tau)e^{\tau s}d\tau$ in (1.9) does not have an analytic closed-form expression.

It should be noted that this chapter is not intended to be a complete survey of frequency-domain approaches [216] for the stability analysis and stabilization of LTDSs. For further information on this topic, readers may refer to monographs [17, 178, 179, 217–221] and the references therein. Additional approaches, such as algebraic factorizations on rings [222–224], the Nyquist criteria [225] and root locii [226], frequency-sweeping framework [227–236], cluster treatment of characteristic Roots (CTCR) [237–248], linear multistep (LMS) methods for computing characteristic roots [249, 250], spectral manifold [251, 252], geometric methods [253–262], $\mathcal{H}^2$ norm computation [263–265], etc, can be found in [254, 266–282].

Over the past decades, significant progress has been made in the research on the stabilization of LTDSs in the frequency-domain. A non-exhaustive list of methods developed for this purpose is presented below:

- Stabilization of linear systems with pointwise input/output delays [283–302].
- Delay reduction approach via transcendental matrix equations [303–306]
- Stabilization of linear system with a state distributed delay: [307–309]
- Stabilization of linear systems with state pointwise-delays [310–314]
- (Quasi) Pole placement for LTDSs [310, 315–317]
- PID controller for linear system with pointwise-delays [318–326]
- Stabilization of linear systems with pointwise-delays and performance constraints [327–330]
- Stabilization of SISO LTDSs [331, 332]
- Stabilization of SISO IDSs systems (IDSs [333–338]
- Stabilization of MIMO IDSs [26, 33, 339–347]

To the best of our knowledge, the latest trends in the frequency-domain-based methods for stabilizing linear systems with delays are represented by the results in [8, 327, 329] and [348–350][4]. The rapid development of these methods is predominantly nourished by the recent advancement of numerical algorithms for non-smooth optimization such as gradient sampling [351–355], BFGS-SQP [356–358] and bundle trust region [350, 359–366] algorithms.

[3]"General" here means that $A(\tau)$ is not a constant over $[-r, 0]$

[4]The synthesis solutions in [348, 349] are developed for infinite-dimensional linear systems, hence it may be applied to design a controller for a linear system with delays whenever it is applicable.

### 1.3.2 Time-domain approaches

Research on systems with delays in the time-domain has a long history [6, 9–11, 17, 367–372], with various perspectives being explored. By extending classical optimal control theory, numerous results on the linear-quadratic control of LTDSs with finite [373–387] and infinite [388–392] time horizon have been established since the 1970s by researchers such as M.C. Delfour [393–398], H.T. Banks [391, 399], J.S. Gibson [389, 400] and others [399, 401]. However, these approaches require solving matrix ODE-PDE/integral boundary-valued problems [18, 402–405] coupled with algebraic Riccati equations, or operator Riccati equations with infinite dimension [24, 25, 399, 401, 406–416], both of which are extremely challenging to compute numerically.

The LKF approach [178, 417–419] is one of the major time-domain methods for analyzing the stability and stabilization of systems with delays. By exploiting the properties of functionals with predefined structures, stability (synthesis) conditions can be constructed using the unknown parameters of the functionals. These conditions can then be posed as optimization constraints and efficiently solved using numerical algorithms. The main challenges of the LKF approach lie in its inherent conservatism, which is significantly influenced by the predefined structures and mathematical handling of the functionals to be constructed. As a matter of fact, the LKF approach for delay systems can be considered as an extension of the second Lyapunov method for finite-dimensional systems, where the latter was successfully applied to characterize the stability of $\dot{\boldsymbol{x}}(t) = A\boldsymbol{x}(t)$ via a quadratic Lyapunov function $\mathsf{v}(\boldsymbol{\xi}) = \boldsymbol{\xi}^\top P\boldsymbol{\xi}$. General LKFs such as

$$\mathsf{v}(\boldsymbol{\phi}(\cdot)) := \boldsymbol{\phi}^\top(0)P_1\boldsymbol{\phi}(0) + 2\boldsymbol{\phi}^\top(0)\int_{-r}^{0} P_2(\tau)\boldsymbol{\phi}(\tau)\mathsf{d}\tau + \int_{-r}^{0}\int_{-r}^{0} \boldsymbol{\phi}^\top(\theta)P_3(\tau,\theta)\boldsymbol{\phi}(\tau)\mathsf{d}\tau\mathsf{d}\theta + \int_{-r}^{0} \boldsymbol{\phi}^\top(\tau)Q(\tau)\boldsymbol{\phi}(\tau)\mathsf{d}\tau \tag{1.10}$$

were previously investigated in [420, 421] to carry out the stability analysis of a simple delay system $\dot{\boldsymbol{x}}(t) = A_1\boldsymbol{x}(t) + A_2\boldsymbol{x}(t-r)$, where $P \in \mathbb{S}^n$, $P_2(\tau) \in \mathbb{R}^{n\times n}$, $P_3(\tau,\theta) \in \mathbb{R}^{n\times n}$, $Q(\tau) \in \mathbb{S}^n$ with $P_3(\tau,\theta) = P_3^\top(\theta,\tau)$, $\forall\tau;\theta \in [-r,0]$. It has been shown in [422] that the form of (1.10) with $Q(\tau) = Q_1 + (\tau+r)Q_2$, which can render (1.10) to admit a quadratic lower bound [422], is adequate to provide sufficient and necessary stability conditions for the asymptotic stability. Thus (1.10) with $Q(\tau) = Q_1 + (\tau+r)Q_2$ can be considered as an example of complete Lyapunov-Krasovskiĭ functional[5]. It has been shown in [423] that a complete LKF exists that can determine the stability of a linear time-delay system (TDS) with general delay structure. Furthermore, the concept of a complete LKF has been extended in [22] as

$$\mathsf{v}(\boldsymbol{\phi}(\cdot)) := \boldsymbol{\xi}^\top P_1\boldsymbol{\xi} + 2\boldsymbol{\xi}^\top\int_{-r}^{0} P_2(\tau)\boldsymbol{\phi}(\tau)\mathsf{d}\tau + \int_{-r}^{0}\int_{-r}^{0} \boldsymbol{\phi}^\top(\theta)P_3(\tau,\theta)\boldsymbol{\phi}(\tau)\mathsf{d}\tau\mathsf{d}\theta + \int_{-r}^{0} \boldsymbol{\phi}^\top(\tau)Q(\tau)\boldsymbol{\phi}(\tau)\mathsf{d}\tau \tag{1.11}$$

with $P \in \mathbb{S}^n$, $P_2(\tau) \in \mathbb{R}^{n\times\nu}$, $P_3(\tau,\theta) = P_3^\top(\theta,\tau) \in \mathbb{R}^{\nu\times\nu}$ and $Q(\tau) \in \mathbb{S}^\nu$, which can provide sufficient and necessary conditions for the stability of a coupled differential-difference equation

$$\dot{\boldsymbol{x}}(t) = A_1\boldsymbol{x}(t) + A_2\boldsymbol{y}(t-r), \quad \boldsymbol{y}(t) = A_3\boldsymbol{x}(t) + A_4\boldsymbol{y}(t-r) \tag{1.12}$$

with $r > 0$ and $\boldsymbol{x}(t) \in \mathbb{R}^n$ and $\boldsymbol{y}(t) \in \mathbb{R}^\nu$.

[5] For a thorough study on the theory of complete LKFs for LTDSs, refer to [419]

**Remark 1.2.** Note that for a system with multiple delays such as

$$\dot{\boldsymbol{x}}(t) = A\boldsymbol{x}(t) + \sum_{i=1}^{p} A_i \boldsymbol{x}(t - r_i), \quad r_{i+1} > r_i, \ \ \forall i = 1 \dots p, \tag{1.13}$$

the structure of the corresponding the complete LKF can be easily determined by adding more delay related integral terms in (1.11). See the functional in Chapter 7 of [178] for details.

In contrast to finding a quadratic function $\boldsymbol{x}^\top(t)Q\boldsymbol{x}(t)$ for $\dot{\boldsymbol{x}}(t) = A\boldsymbol{x}$ where $Q$ has finite dimensions, the dimensions of the decision variables $P_2(\cdot), P_3(\cdot,\cdot)$ in (1.11) are infinite-dimensional and significantly more challenging to construct numerically. As a matter of fact, the handling of the DD kernel $A_3(\cdot)$ in (1.8) inherits similar difficulties, since $A_3(\cdot)$ is also infinite dimensional. As a result, early results based on the LKF approach generally assumed simple structures for the matrix parameters in (1.11) and with a constant DD kernel $A_3(\tau) = D \in \mathbb{R}^{n\times n}$. For instance, by setting $P_2(\tau) = P_3(\tau, \theta) = \mathsf{O}_n$ and $Q(\tau) = Q_1$, delay-independent conditions can be derived for $\dot{\boldsymbol{x}}(t) = A_1\boldsymbol{x}(t) + A_2\boldsymbol{x}(t - r)$ with (1.11), ensuring that the system is asymptotically stable regardless of the value of delay $r > 0$. However, such stability conditions can be overly restrictive, as many systems may only be stable within specific ranges of delay values. Therefore, it is imperative to derive delay-dependent stability conditions taking into account the value of delays.

Since the mid-1990s, there has been a tremendous amount of literature dedicated to stability analysis [424–437] and stabilization [438–449], supported by efficient numerical algorithms for solving linear matrix inequalities [167]. For comprehensive collections of existing literature on this topic, see monographs [178, 405, 450–452]. While assuming $P_2(\cdot), P_3(\cdot,\cdot)$ to be constants can lead to stability (synthesis) conditions with finite dimensions, the resulting conservatism can be significant since the original structures of $P_2(\cdot), P_3(\cdot,\cdot)$ in (1.11) are functions. Thus, it is preferable to consider more general structures for $P_2(\cdot), P_3(\cdot,\cdot)$ as functions. On the other hand, simply adding more terms to the functionals may not necessarily result in less conservative stability conditions. One critical factor in deriving stability conditions using the LKF approach is the application of integral inequalities (IIs), whose lower bounds conservatism can significantly impact the conservatism of the resulting stability (synthesis) conditions. Finally, it is worth remarking that only sufficient conditions can be obtained generally via the LKF approach due to the intrinsic mathematical structures encountered in the procedure of deriving stability (synthesis) conditions that are generally denoted by matrix inequalities. Therefore, constructing the stability (synthesis) conditions with less conservatism and fewer variables is a paramount goal when considering the LKF approach.

Since the dimensions of $P_2(\cdot)$ and $P_3(\cdot,\cdot)$ in (1.11) are infinite, we may need to apply certain discretization or approximation schemes to the variables so that they can be constructed via finite-dimensional SDP. In [453], Gu initiated the use of piecewise linear functions over subregions of a delay interval to discretize the variables $P_2(\cdot), P_3(\cdot,\cdot)$ and $Q(\cdot)$ in (1.11). His idea has produced fruitful results over the past decade [454–458] and has been successfully extended to tackle linear coupled differential-difference systems [22, 459, 460], and linear systems with DD terms having piecewise constant integrands [461, 462]. However, using piecewise linear functions as the basis for $P_2(\cdot), P_3(\cdot,\cdot)$ and $Q(\cdot)$ can be too conservative, and it is not clear how to handle DDs with more general kernels. In [463], a novel approach was presented using the full-block S-procedure [464] to handle linear systems with rational DD kernels. This method could be applied to a system with general DDs by approximating the DD kernels via rational functions [465], although the relationship between the stability of the original and approximated systems is unclear. On the other hand, the stabilization condition in [463] requires that $(A, B)$ is stabilizable[6] where $A$ is the delay-free state space matrix and $B$ is the input gain matrix, indicating that the induced conservatism cannot be ignored.

Research on finding new IIs for constructing Lyapunov-Krasovskiĭ functionals (LKFs) has become increasingly popular. In [436], a general form of the Wirtinger inequality was derived, generalizing the integral form of the Jensen inequality [466]. A significant breakthrough was made in [467] with the derivation of the Legendre-Bessel integral inequality, which generalizes both the Jensen [466] and Wirtinger [436] inequalities and is well-suited for constructing stability conditions using functionals like (1.11) with polynomial kernels for $P_2(\cdot)$ and $P_3(\cdot, \cdot)$. A distinct feature of the results in [467] is that the feasibility of the stability conditions is hierarchical with respect to the degree of Legendre polynomials [468, 22.2.10] in the LKF. Moreover, the numerical examples tested in [467] and the subsequent literature [469] had clearly demonstrated the advantage of the reported methods over existing results in terms of reduced conservatism and lower numerical complexity. Following this idea, the Legendre-Bessel integral inequality was further extended in [470] to analyze the stability of linear systems with a general DD whose kernels are continuous, where the kernels are approximated by Legendre polynomials. Although the DDs in [470] are handled via approximation, the stability of the original DD system was analyzed, and the results are not obtained from an approximated system that neglects certain information. This is reflected by the fact that the approximation error of the DD function is included in the hierarchical stability condition in [470]. However, a potential issue with the method developed in [470] is that a very large degree of Legendre polynomials may be required to approximate a DD kernel function if its form is not well-suited for polynomial approximations.

In addition to the LKF approach, there are various time-domain methods that can be applied to analyze the stability and stabilization of systems with delays. For example, necessary or both sufficient and necessary conditions for the stability of LTDSs have been proposed in [471–473] and [474–478], respectively, denoted using delay Lyapunov matrix. However, as remarked in [479], how to solve the delay Lyapunov matrix numerically remains a critical challenge for the applicability of this approach. A very interesting method for stabilizing LTDSs with DDs was proposed in [480], combining time and frequency-domain methods based on the concept of smoothed spectral abscissa [481] and delay Lyapunov matrix [474, 476, 477]. However, this method requires numerical computation of the delay Lyapunov matrix and its derivatives, and it is not clear how this can be carried out for an LTDS with general DDs where all delays are non-commensurate. A classical approach [482], based on the utilization of spectral decomposition [9, 483], established that if an LTDS is stabilizable/detectable [484–486], a controller/observer can always be constructed by stabilizing its "unstable" component using methods for simple LTI systems, given that the spectrum of an LTDS can contain only a finite number of unstable roots [8] as emphasized in Section 1.3.1. Nonetheless, this approach is feasible only if all unstable roots can be computed numerically, which in itself represents a quite challenging aspect in the research of LTDSs. Another representative approach in addressing the presence of time delays is the construction of predictor controllers [487, 488], also known as the reduction approach [489, 490], for systems with input/output delays [491–496], and for systems with both state and input delays [305, 306, 495, 497–502]. Notably, the classical prediction scheme has been further integrated with constructive synthesis approaches, aiming to design controllers for delay systems with specific structures (backstepping, forwarding) [503–505]. Finally, the semi-discretization scheme[7] in [219] can be applied to analyze the stability of linear-periodic TDSs with DDs.

[6]For the case of stability analysis, this means that the stability condition can be feasible only if the delay-free matrix $A$ is Hurwitz

[7]This approach might be interpreted as a mixed combination of both time and frequency-domain approaches

### 1.3.3 Computational tools for the analysis and synthesis of delay systems

Many control-related problems [167, 506] can be recast as SDP problems [507], which can then be efficiently solved using numerical algorithms [508–513]. For example, by using a quadratic Lyapunov function for a finite-dimensional linear time-invariant (LTI) system, the problem of stability analysis or stabilization can be transformed into an optimization problem denoted by linear matrix inequalities (LMIs). A significant advantage of using the SDP approach is that stability (stabilization) conditions can be extended (modified) to incorporate additional information, such as dissipativity [514] and uncertainties [515], without introducing intractable mathematical complications. There are numerous references dedicated to the application of LMIs in the field of systems and control [464, 516, 517].

Recently, three papers [518–520] have demonstrated similar methodologies for solving polynomial optimization problems based on the principles of algebraic geometry [521–523]. The so-called sum-of-squares (SoS) programming can impose positive constraints on polynomial decision variables, which can ultimately be solved using finite-dimensional LMIs. This provides an effective approach to solve infinite-dimensional robust LMIs that are polynomially dependent on uncertainties over compact sets. A series of results in [519, 524–529] have addressed a wide variety of potential applications from different mathematical perspectives. In terms of delay systems, the construction of general LKF using the SoS approach has been considered in [530–534] where the functionals, with structures similar to (1.11), contain integral terms with polynomial kernels. However, the numerical complexity when using SoS to construct Krasovskiĭ functionals still needs to be addressed compared to standard LMI approaches.

## 1.4 Research Motivations and Outline of Contents

This section presents the research motivations and outline of this book. The first subsection discusses the theoretical and practical motivations for investigating the stability and stabilization of LDDSs. Moreover, we provide a summary in the second subsection to outline the works of this book by each chapter.

### 1.4.1 Research motivations

The first research motivation for this book is the limited availability of solutions for analyzing the stability and stabilization of linear systems with general DDs. Mathematically, a linear system with both discrete and general DDs can represent a significantly general class of LTDSs.[8] The development of methods for analyzing the stability and stabilization of linear systems with DDs is particularly important for advancing the theory of general linear time delay systems. The results developed in [467, 469, 470] suggest that the LKF approach holds promise for analyzing the stability of delay systems. However, there is still considerable room for developing new approaches and addressing many outstanding questions concerning DDs. For example, can nonconservative stability (synthesis) conditions be derived for linear systems with nontrivial DDs that can detect delay margins? Can DDs be directly addressed in the context of stability analysis (stabilization) without using approximations? Can functions other than Legendre polynomials [470] be used to handle general DDs in stability analysis? Due to the mathematical nature of constructing stability conditions using the LKF approach, these questions imply that new IIs must be developed that can handle non-polynomial integral kernels. Finally, stabilization for systems with delays require further attention, as it may not be trivial to construct convex synthesis conditions using the LKF approach. All these questions are worthy of consideration and are addressed

[8] For instance, see Chapter 7 in [9] for the general model of autonomous linear functional differential systems.

by the results presented in this book.

The second motivation for studying linear systems with delays arises from their ubiquitous theoretical and practical applications. For instance, if the transfer function of a linear distributed parameter system (i.e., a linear IDS) is identical to that of a system with delays, the former stability can be analyzed via the construction of LKFs using the corresponding delay model. Indeed, it has been shown in [535–538] that coupled PDE-ODE systems, which can be applied to drilling mechanisms modeling [44, 45], can be equivalently expressed in terms of systems with delays. Therefore, if a linear distributed parameter system contains a DD, this system might be analyzed by an equivalent linear system with a DD. Furthermore, the LKF approach offers several advantages over common frequency-domain solutions. For example, it does not necessarily introduce intractable mathematical difficulties when extended to consider dissipativity [452], uncertainties [577], or uncertain bounded time-varying delays [539]. The rich existing results on the stability analysis and control of networked control [101, 540–544] and sampled-data systems [545–550] advanced by the LKF approach further demonstrate the importance of developing effective solutions to linear systems with time-varying delays. Additionally, it has been shown that the digital communication channel [551, 552] of a networked control system, which experiences stochastic packet delay and loss, can be modeled using DDs. Consequently, new methods based on the LKF approach for systems with DDs may lead to significant advancements in the modeling and control of networked control systems, which has become one of the major subjects in the field of control engineering. Finally, given the similarity [553] between analyzing delay systems and distributed parameter systems, results derived by the Krasovskiĭ functional approach for delay systems can be extended to handle general PDE systems [554–559], which have very rich practical applications.

### 1.4.2 Outline of contents

The contents of the rest of the chapters in this book are summarized as follows.

In Chapter 2, we examine the problem of stabilizing a linear system with DDs subject to dissipativity constraints. The DDs exist in the states, inputs, and outputs of the system, and the DD kernels can include a certain class of elementary functions such as polynomials, trigonometric and exponential functions. Sufficient conditions for the existence of a stabilizing state feedback controller under the dissipativity constraints are derived in terms of linear and bilinear matrix inequalities (BMIs) via constructing a LKF related to the distributed kernels. The construction of the functional is achieved through the application of a novel general integral inequality. To tackle the non-convexity induced by the BMI in the synthesis conditions, we employ the Projection Lemma to convexify the BMI so that a sufficient convex condition is obtained. Moreover, an iterative algorithm is constructed to further enhance the feasibility of our methods. Finally, the testing results of numerical examples are presented to demonstrate the strength and effectiveness of the proposed methodology.

In Chapter 3, the solutions to stability analysis and stabilization presented in Chapter 2 are further extended to deal with the problem of stabilizing an uncertain linear system with DDs subject to dissipativity constraints. The uncertainties are of linear fractional forms and subject to the constraints of full block scaling structures. To handle these complex uncertainties, we propose a novel lemma that establishes a relation between finite dimensional LMIs and a robust LMI with linear fractional uncertainties, where the well-posedness of the uncertainties can be determined by a matrix inequality condition with finite dimensions. Based on the results derived in Chapter 2 without considering the presence of system uncertainties, two theorems containing synthesis conditions can be derived based on the application of the aforementioned lemma where the conditions of the second theorem are convex. Similar to the paradigm utilized in Chapter 2, an iterative algorithm is also proposed to solve the BMI in the first theorem to further reduce conservatism. A distinct feature of the result in this chapter is that the proposed method is further modified to design a resilient dynamical state feedback for a linear system with input delays, where both the plant and resulting controller are robust against uncertainties. More importantly, the corresponding iterative algorithm for the design of resilient dynamical state feedback can be initiated simply by the gain of a constructible predictor controller, without appealing to a separate theorem. A battery of numerical examples was tested to demonstrate the effectiveness of our proposed methods.

Knowing the crucial importance of using optimal IIs with the LKF approach, two general classes of integral inequality are specifically developed in Chapter 4 for the stability analysis and stabilization of linear systems with delay. Our inequalities exhibit very general structures in terms of the generality of weight functions and integral kernels of the lower quadratic bounds. Almost all existing inequalities in the peer-reviewed literature, including those with free matrix variables, are special cases of our proposed inequalities. Moreover, relations are established concerning the lower bounds between our proposed inequalities. For specific applications concerning systems with delays, our inequalities are applied to construct a stability condition via the LKF approach for a linear continuous differential difference system (CDDS) with a DD. It is shown that the resulting stability condition is invariant with respect to a parameter in the LKF and equivalent stability conditions can be derived by means of different types of inequalities. Finally, it is worth noting that the proposed inequalities can be applied in more general contexts such as the stability analysis of PDE-related systems or sampled-data systems.

Chapter 5 presents a new method for the dissipativity and stability analysis of a CDDS with general distributed-delays at both state and output. More precisely, the DD under consideration can contain any $\mathscr{L}^2$ which are approximated by a class of elementary functions that includes Legendre polynomials [468, 22.2.10]. By taking this broader class of functions compared with the existing approach in [470] where only Legendre polynomials are utilized for approximations, one can construct LKFs with more general structures as compared with the existing approach in [470], where the functional is parameterized via Legendre polynomials. Furthermore, a novel generalized integral inequality is also proposed to incorporate approximation error in our stability (dissipativity) conditions. Based on the proposed approximation scenario with the proposed integral inequality, sufficient conditions determining the dissipativity and stability of a CDDS are derived in terms of LMIs. In addition, several hierarchies in terms of the feasibility of the proposed conditions are derived under certain constraints. Finally, the testing results of several numerical examples are presented in this chapter to show the effectiveness of our proposed methodologies.

The problem of delay range stability analysis for a CDDS with DDs is investigated in Chapter 6 where the system is also subject to a dissipative constraint. Polynomial DD kernels are considered in this chapter so that tractable stability conditions can be obtained. A LKF with non-constant matrix parameters, which are related to the delay value polynomially, is applied to construct sufficient conditions that guarantee range stability of a linear CDDS subject to a dissipative constraint. The proposed sufficient conditions for the stability and dissipativity of the system are denoted by sum-of-squares constraints which are constructed based on the application of a matrix relaxation technique for robust LMIs without introducing any potential conservatism. Furthermore, the proposed methods can be extended to solve delay margin estimation problems for a linear CDDS subject to prescribed dissipative constraint. Finally, numerical examples are tested to demonstrate the effectiveness of the proposed methodologies.

In Chapter 7, we develop new methods to stabilize a linear DD system where the delay function is time-varying and bounded. The DDs exist in the states, inputs, and outputs of the system, and the DD kernels in this chapter can be functions belonging to a class of elementary functions. Furthermore, a novel integral inequality is proposed to construct synthesis (stability) conditions via LKFs where the conditions are expressed in terms of matrix inequalities. The proposed synthesis (stability) conditions, which can determine the stability and dissipativity of the system with a supply rate function, are related to the values of the bounds of the time-varying delay where the information of the derivatives of the time-varying delays is absent. Given the results in Chapter 2 concerning the handling of BMIs, the resulting synthesis (stability) conditions in this chapter can be either solved directly by standard solvers for SDP problems if they are convex or reshaped into linear matrix inequalities (LMIs), or solved by an iterative algorithm. Finally, the testing results of various numerical examples are provided to demonstrate the effectiveness of the proposed methodologies.

At the end of this book in Chapter 8, some future works and suggestions are discussed and evaluated.

# Chapter 2

# Dissipative Stabilization of Linear Delay Systems with Distributed Delays

## 2.1 Introduction

As previously noted in Chapter 1, DDs can be very difficult to deal with because of their infinite dimensions. If one assumes a linear system with DDs of constant integral kernels and analyzes this system via a functional with constant decision variables, it is possible to obtain constraints denoted by LMIs with finite dimensions. There has been a significant body of literature in this direction aimed at performing either stability analysis or controller synthesis for linear DDS [445, 448, 449]. For a collection of the previous works on this topic, consult the monographs [405, 452].

In terms of addressing DDs with non-constant kernels, the results in [560] have demonstrated that certain linear DDS can be transformed into a system with only pointwise-delays. However, this method introduces noticeable conservatism due to the additional dynamics required by the introduction of new states. To mitigate this conservatism, it is advantageous to employ CLKFs as the primary tool for handling systems with DDs where the functionals are parameterized by unknowns of infinite dimensions. This approach was first implemented in [462, 561] using a discretization scheme with piecewise linear functions, where the LDDS contains piecewise constant kernels. Furthermore, by applying the full-block S-procedure [464], a novel synthesis method was set forth in [463] for addressing systems with rational DD kernels, which is capable of dealing with general distributed terms via approximations [465]. However, the derived stabilization conditions require $(A, B)$ to be stabilizable in [463] ($A$ is the delay-free state space matrix and $B$ is the input gain matrix), indicating that the induced conservatism cannot be disregarded. Finally, a systematic methodology for constructing controllers for linear systems with discrete and DDs, featuring forwarding or backstepping structures, has been explored in [503].

In this chapter, we present methods for stabilizing linear DDS with DDs in states, inputs, and outputs. The DDs considered here may have non-constant structures, with delay kernel functions that can be polynomials, trigonometric, or exponential functions. Our synthesis schemes incorporate a quadratic supply rate function [452, 514], providing a comprehensive characterization of controller performance. Additionally, we derive a new integral inequality for formulating synthesis conditions, which can be viewed as a generalization of the recently proposed Bessel-Legendre inequality [467, 469]. By constructing a general LKF using this inequality, we derive sufficient conditions for the existence of a stabilizing state feedback controller that accounts for dissipativity constraints. These conditions are expressed as matrix inequalities containing a non-convex bilinear matrix inequality.

We then utilize [562–564, Projection Lemma] on the BMI to construct convex synthesis conditions, while minimizing conservatism by simplifying a portion of the slack variables. To further reduce conservatism, we develop an iterative algorithm based on the inner approximation approach in [565]. In contrast to existing methods, our proposed synthesis solutions do not require $(A, B)$ to be stabilizable as in [463], nor demand forwarding or backstepping structures as in [503]. In addition, the integral kernels of our LKF are not necessarily polynomials compared with the existing results in [467, 469]. Regarding the performance of stability analysis, a numerical example is tested which can demonstrate that our approach can outperform the method in [470] in terms of numerical complexity. This mainly achieved thanks to the direct handling of DDs without resorting to polynomial approximations [470], as some functions are difficult to approximate with insufficient degrees.

This chapter is structured as follows: The synthesis problem is first formulated in Section 2.2. Secondly, we present vital mathematical results in Section 2.3 for the derivation of the synthesis solutions in this chapter. Next, the main results on controller synthesis are given in Section 2.4. To illustrate the capacity and effectiveness of our methodologies, we investigated several numerical examples in Section 2.5.

## 2.2 Problem formulation

The following property of the Kronecker product is made significant use throughout the whole work.

**Lemma 2.1.** $\forall X \in \mathbb{R}^{n\times m}, \ \ \forall Y \in \mathbb{R}^{m\times p}, \ \ \forall Z \in \mathbb{R}^{q\times r}$,

$$(X \otimes I_q)(Y \otimes Z) = (XY) \otimes (I_q Z) = (XY) \otimes Z = (XY) \otimes (Z I_r) = (X \otimes Z)(Y \otimes I_r). \tag{2.1a}$$

*Moreover,* $\forall X \in \mathbb{R}^{n\times m}$, *we have*

$$\begin{bmatrix} A & B \\ C & D \end{bmatrix} \otimes X = \begin{bmatrix} A\otimes X & B\otimes X \\ C\otimes X & D\otimes X \end{bmatrix} \tag{2.1b}$$

*for any* $A, B, C, D$ *with appropriate dimensions for the partition of the block matrix.*

Consider the following model of LDDS

$$\begin{aligned}
\dot{\boldsymbol{x}}(t) &= A_1\boldsymbol{x}(t) + A_2\boldsymbol{x}(t-r) + \int_{-r}^{0} \widetilde{A}_3(\tau)\boldsymbol{x}(t+\tau)\mathsf{d}\tau + B_1\boldsymbol{u}(t) + B_2\boldsymbol{u}(t-r) \\
&\qquad + \int_{-r}^{0} \widetilde{B}_3(\tau)\boldsymbol{u}(t+\tau)\mathsf{d}\tau + D_1\boldsymbol{w}(t), \ \ \forall t \geq t_0 \\
\boldsymbol{z}(t) &= C_1\boldsymbol{x}(t) + C_2\boldsymbol{x}(t-r) + \int_{-r}^{0} \widetilde{C}_3(\tau)\boldsymbol{x}(t+\tau)\mathsf{d}\tau + B_4\boldsymbol{u}(t) + B_5\boldsymbol{u}(t-r) \\
&\qquad + \int_{-r}^{0} \widetilde{B}_6(\tau)\boldsymbol{u}(t+\tau)\mathsf{d}\tau + D_2\boldsymbol{w}(t) \\
\dot{\boldsymbol{x}}(t) &:= \lim_{\eta\downarrow 0} \frac{\boldsymbol{x}(t+\eta)-\boldsymbol{x}(t)}{\eta}, \qquad \forall\theta\in[-r,0], \ \ \boldsymbol{x}(t_0+\theta) = \boldsymbol{\phi}(\theta),
\end{aligned} \tag{2.2}$$

where $t_0 \in \mathbb{R}$ and $\boldsymbol{\phi}(\cdot) \in \mathcal{C}\left([-r,0];\mathbb{R}^n\right)$, and $\boldsymbol{x}(t) \in \mathbb{R}^n$ satisfies the delay equation in (2.2), $\boldsymbol{u}(t) \in \mathbb{R}^p$ denotes input signals, $\boldsymbol{w}(\cdot) \in \mathcal{L}^2([t_0,\infty);\mathbb{R}^q)$ represents the disturbance and $\boldsymbol{z}(t) \in \mathbb{R}^m$ is the regulated output. Note that $\boldsymbol{\phi}(\cdot) \in \mathcal{C}\left([-r,0];\mathbb{R}^n\right)$ is the initial condition for (2.2) at $t = t_0$. The size of the state spaces matrices in (2.2) is determined by the values of $n; m; p; q \in \mathbb{N}$. Here, we assume delay $r > 0$ is known. Finally, $\widetilde{A}_3(\tau), \widetilde{B}_3(\tau), \widetilde{C}_3(\tau)$ and $\widetilde{B}_6(\tau)$ satisfy the following assumption:

**Assumption 2.1.** There exist $\mathbf{Col}_{i=1}^{d} f_i(\tau) = \boldsymbol{f}(\cdot) \in \mathcal{C}^1([-r,0];\mathbb{R}^d)$ with $d \in \mathbb{N}$, and constant matrices $A_3 \in \mathbb{R}^{n\times dn}$, $B_3 \in \mathbb{R}^{n\times dp}$, $C_3 \in \mathbb{R}^{m\times dn}$, $B_6 \in \mathbb{R}^{m\times dp}$ such that $\forall \tau \in [-r,0]$, $\mathbb{R}^{n\times n} \ni \widetilde{A}_3(\tau) = A_3 F(\tau)$ and $\mathbb{R}^{n\times p} \ni \widetilde{B}_3(\tau) = B_3\left(\boldsymbol{f}(\tau)\otimes I_p\right)$ and $\mathbb{R}^{m\times n} \ni \widetilde{C}_3(\tau) = C_3 F(\tau)$ and $\mathbb{R}^{m\times p} \ni \widetilde{B}_6(\tau) = B_6\left(\boldsymbol{f}(\tau)\otimes I_p\right)$ where $\mathbb{R}^{dn\times n} \ni F(\tau) := \boldsymbol{f}(\tau)\otimes I_n$ and $\boldsymbol{f}(\cdot)$ satisfies

$$\exists M \in \mathbb{R}^{d\times d},\ \frac{\mathsf{d}\boldsymbol{f}(\tau)}{\mathsf{d}\tau} = M\boldsymbol{f}(\tau), \tag{2.3a}$$

$$\int_{-r}^{0} \boldsymbol{f}(\tau)\boldsymbol{f}^\top(\tau) \succ 0. \tag{2.3b}$$

**Remark 2.1.** The condition in (2.3a) indicates that the functions in $\boldsymbol{f}(\cdot)$ are the solutions to linear homogeneous differential equations with constant coefficients such as polynomials, trigonometric and exponential functions. In other words, $\boldsymbol{f}(\cdot)$ contains functions that are entries of $\mathsf{e}^{X\tau}$ with $X \in \mathbb{R}^{d\times d}$. In addition, there is no restriction on the dimension of $\boldsymbol{f}(\tau)$ as long as it can encompass all kernels in the DDs of (2.2). The inequality in (2.3b) indicates that the functions in $\boldsymbol{f}(\cdot)$ are linearly independent in a Lebesgue sense according to [566, Theorem 7.2.10]. With respect to the generality of $\boldsymbol{f}(\tau)$, a wide range of applications can be modeled by (2.2) based on the structure in Assumption 2.1. For instance, examples include the compartmental dynamic systems with DDs examined in [567], and the LDDS with gamma distributions in [463] and [568] with a finite delay range.

Assume that all states are available to feedback in (2.2), and our goal is to construct a state feedback controller $\boldsymbol{u}(t) = K\boldsymbol{x}(t)$ with $K \in \mathbb{R}^{p\times n}$ to stabilize the open-loop system. By utilizing $\boldsymbol{u}(t) = K\boldsymbol{x}(t)$ with (2.2) and taking into account Assumption 2.1 with (2.1a), we find that the input DD matrices become

$$\begin{aligned} B_3F(\tau)K &= B_3\left(\boldsymbol{f}(\tau)\otimes I_p\right)K = B_3\left(I_d\otimes K\right)\left(\boldsymbol{f}(\tau)\otimes I_n\right), \\ B_6F(\tau)K &= B_6\left(\boldsymbol{f}(\tau)\otimes I_p\right)K = B_6\left(I_d\otimes K\right)\left(\boldsymbol{f}(\tau)\otimes I_n\right). \end{aligned} \tag{2.4}$$

As a result, the resulting closed-loop system (CLS) can be written as

$$\begin{aligned} &\dot{\boldsymbol{x}}(t) = \big(\mathbf{A} + \mathbf{B}_1\left[(I_{2+d}\otimes K)\oplus \mathsf{O}_q\right]\big)\boldsymbol{\chi}(t),\ t \geq t_0 \\ &\boldsymbol{z}(t) = \big(\mathbf{C} + \mathbf{B}_2\left[(I_{2+d}\otimes K)\oplus \mathsf{O}_q\right]\big)\boldsymbol{\chi}(t) \\ &\forall\theta \in [-r,0],\quad \boldsymbol{x}(t_0+\theta) = \boldsymbol{\phi}(\theta) \end{aligned} \tag{2.5}$$

with $t_0 \in \mathbb{R}$ and $\boldsymbol{\phi}(\cdot) \in \mathcal{C}\left([-r,0];\mathbb{R}^n\right)$ in (2.2), where

$$\mathbf{A} = \begin{bmatrix} A_1 & A_2 & A_3 & D_1 \end{bmatrix},\quad \mathbf{B}_1 = \begin{bmatrix} B_1 & B_2 & B_3 & \mathsf{O}_{n,q} \end{bmatrix} \tag{2.6a}$$

$$\mathbf{C} = \begin{bmatrix} C_1 & C_2 & C_3 & D_2 \end{bmatrix},\quad \mathbf{B}_2 = \begin{bmatrix} B_4 & B_5 & B_6 & \mathsf{O}_{m,q} \end{bmatrix} \tag{2.6b}$$

$$\boldsymbol{\chi}(t) := \begin{bmatrix} \boldsymbol{x}^\top(t) & \boldsymbol{x}^\top(t-r) & \int_{-r}^{0}\boldsymbol{x}^\top(t+\tau)F^\top(\tau)\mathsf{d}\tau & \boldsymbol{w}^\top(t) \end{bmatrix}^\top,\quad F(\tau) = \boldsymbol{f}(\tau)\otimes I_n. \tag{2.6c}$$

## 2.3 Important Mathematical Tools

In this section, we present several essential theoretical tools for deriving our synthesis conditions in this chapter. These include the Krasovskiĭ stability criteria for determining the stability of delay systems and the definition of dissipativity. Additionally, we propose a new integral inequality that will be utilized in constructing an LKF in the following section.

**Lemma 2.2.** *Given $r > 0$, the CLS in (2.5) with $\boldsymbol{w}(t) \equiv \mathbf{0}_q$ is uniformly globally asymptotically stable at its origin if there exist $\epsilon_1; \epsilon_2; \epsilon_3 > 0$ and a differentiable functional $\mathsf{v} : \mathcal{C}([-r, 0]; \mathbb{R}^n) \to \mathbb{R}$ with $\mathsf{v}(\mathbf{0}_n) = 0$ such that*

$$\epsilon_1 \|\boldsymbol{\phi}(0)\|_2^2 \leq \mathsf{v}(\boldsymbol{\phi}(\cdot)) \leq \epsilon_2 \|\boldsymbol{\phi}(\cdot)\|_\infty^2 , \tag{2.7a}$$

$$\left.\frac{\mathrm{d}^+}{\mathrm{d}t}\mathsf{v}(\mathbf{x}_t(\cdot))\right|_{t=t_0, \mathbf{x}_{t_0}(\cdot)=\boldsymbol{\phi}(\cdot)} \leq -\epsilon_3 \|\boldsymbol{\phi}(0)\|_2^2 \tag{2.7b}$$

*hold for any $\boldsymbol{\phi}(\cdot) \in \mathcal{C}([-r, 0]; \mathbb{R}^n)$ in (2.5), where $t_0 \in \mathbb{R}$ and $\|\boldsymbol{\phi}(\cdot)\|_\infty := \sup_{-r \leq \tau \leq 0} \|\boldsymbol{\phi}(\tau)\|_2$ and $\frac{\mathrm{d}^+}{\mathrm{d}x} f(x) = \operatorname{limsup}_{\eta \downarrow 0} \frac{f(x+\eta) - f(x)}{\eta}$. Moreover, the term $\mathbf{x}_t(\cdot)$ in (2.7b) is given by $\forall t \geq t_0$, $\forall \theta \in [-r, 0]$, $\mathbf{x}_t(\theta) = \boldsymbol{x}(t + \theta)$ where $\boldsymbol{x} : [t_0 - r, \infty) \to \mathbb{R}^n$ satisfies (2.5) with $\boldsymbol{w}(t) \equiv \mathbf{0}_q$.*

*Proof.* Let $u(\cdot), v(\cdot), w(\cdot)$ in [178, Theorem 3] be quadratic functions with positive parameters $\epsilon_1, \epsilon_2$ and $\epsilon_3$, then Lemma 2.2 can be obtained accordingly since (2.5) with $\boldsymbol{w}(t) \equiv \mathbf{0}_q$ is a particular case of the system $\dot{\boldsymbol{x}}(t) = \boldsymbol{f}(t, \mathbf{x}_t(\cdot)))$, $\mathbf{x}_t(\cdot) \in \mathcal{C}([-r, 0]; \mathbb{R}^n)$ considered in [178, Theorem 3]. ■

The following definition of dissipativity is presented based on the general definition of dissipativity in [569].

**Definition 2.1.** Given $r > 0$, the CLS in (2.5) with a supply rate function $\mathsf{s}(\boldsymbol{z}(t), \boldsymbol{w}(t))$ is dissipative if there exists a differentiable functional $\mathsf{v} : \mathcal{C}([-r, 0]; \mathbb{R}^n) \to \mathbb{R}_{\geq 0}$ such that

$$\forall t \geq t_0, \quad \dot{\mathsf{v}}(\mathbf{x}_t(\cdot)) - \mathsf{s}(\boldsymbol{z}(t), \boldsymbol{w}(t)) \leq 0 \tag{2.8}$$

with $t_0 \in \mathbb{R}$ in (2.5), where $\dot{\mathsf{v}}(\mathbf{x}_t(\cdot))$ is well defined for all $t \geq t_0$ and $\mathbf{x}_t(\cdot)$ satisfies $\forall t \geq t_0$, $\forall \theta \in [-r, 0]$, $\mathbf{x}_t(\theta) = \boldsymbol{x}(t + \theta)$ with $\boldsymbol{x}(t)$ and $\boldsymbol{z}(t)$ in (2.5) with $\boldsymbol{w}(\cdot) \in \mathcal{L}^2([t_0, \infty); \mathbb{R}^q)$.

If (2.8) holds, then we have

$$\forall t \geq t_0, \;\; \mathsf{v}(\mathbf{x}_t(\cdot)) - \mathsf{v}(\mathbf{x}_{t_0}(\cdot)) \leq \int_{t_0}^{t} \mathsf{s}(\boldsymbol{z}(\theta), \boldsymbol{w}(\theta)) \mathrm{d}\theta \tag{2.9}$$

which now aligns the original definition of dissipativity in [569], given that $\dot{\mathsf{v}}(\mathbf{x}_t(\cdot))$ is well defined for all $t \geq t_0$.

To denote dissipativity in this chapter, we apply a quadratic supply rate function

$$\mathsf{s}(\boldsymbol{z}(t), \boldsymbol{w}(t)) = \begin{bmatrix} \boldsymbol{z}(t) \\ \boldsymbol{w}(t) \end{bmatrix}^\top \mathbf{J} \begin{bmatrix} \boldsymbol{z}(t) \\ \boldsymbol{w}(t) \end{bmatrix}, \quad \mathbf{J} = \begin{bmatrix} \widetilde{J}^\top J_1^{-1} \widetilde{J} & J_2 \\ * & J_3 \end{bmatrix} \in \mathbb{S}^{(m+q)}, \;\; \widetilde{J}^\top J_1^{-1} \widetilde{J} \preceq 0, \;\; J_1^{-1} \prec 0 \tag{2.10}$$

where the form of $\mathbf{J}$ is constructed considering the general quadratic constraints applied in [514] together with the idea of factorizing the matrix $U_j$ in [514]. Note that (2.10) is able to characterize numerous performance criteria such as

- $\mathcal{L}^2$ gain performance: $J_1 = -\gamma I_m$, $\widetilde{J} = I_m$, $J_2 = \mathsf{O}_{m,q}$, $J_3 = \gamma I_q$ where $\gamma > 0$.
- Strict Passivity: $J_1 \prec 0$, $\widetilde{J} = \mathsf{O}_m$, $J_2 = I_m$, $J_3 = \mathsf{O}_m$ with $m = q$.
- Various sector constraints in [570, Table 1].

The following general integral inequality is derived to be utilized for the construction of KFs in the next section.

**Lemma 2.3.** *Given $\mathcal{K} \subseteq \mathbb{R} \cup \{\pm\infty\}$ where the Lebesgue measure [571] of $\mathcal{K}$ is nonzero, and $U \in \mathbb{S}^n_{\succeq 0}$ and $\boldsymbol{g}(\tau) = \mathbf{Col}_{i=1}^d g_i(\tau) \in \mathscr{L}_2(\mathcal{K}; \mathbb{R}^d)$ with $d \in \mathbb{N}$ which satisfies*

$$\int_{\mathcal{K}} \boldsymbol{g}(\tau)\boldsymbol{g}^\top(\tau)\mathrm{d}\tau \succ 0, \tag{2.11}$$

*then we have*

$$\forall \boldsymbol{x}(\cdot) \in \mathscr{L}_2(\mathcal{K}; \mathbb{R}^n), \quad \int_{\mathcal{K}} \boldsymbol{x}^\top(\tau)U\boldsymbol{x}(\tau)\mathrm{d}\tau \geq \int_{\mathcal{K}} \boldsymbol{x}^\top(\tau)G^\top(\tau)\mathrm{d}\tau\,(\mathsf{G} \otimes U)\int_{\mathcal{K}} G(\tau)\boldsymbol{x}(\tau)\mathrm{d}\tau, \tag{2.12}$$

*where $\mathsf{G}^{-1} = \int_{\mathcal{K}} \boldsymbol{g}(\tau)\boldsymbol{g}^\top(\tau)\mathrm{d}\tau \succ 0$ and $G(\tau) := \boldsymbol{g}(\tau) \otimes I_n$ such that $\boldsymbol{g}(\tau) := \mathbf{Col}_{i=1}^d g_i(\tau)$.*

*Proof.* The proof is inspired by the results in [467, 470]. To begin with, we can conclude that the matrix $\mathsf{G}$ is invertible given (2.11). Let $\boldsymbol{y}(\cdot) \in \mathscr{L}^2(\mathcal{K}; \mathbb{R}^n)$ and

$$\boldsymbol{y}(\tau) := \boldsymbol{x}(\tau) - G^\top(\tau)(\mathsf{G} \otimes I_n)\int_{\mathcal{K}} G(\theta)\boldsymbol{x}(\theta)\mathrm{d}\theta. \tag{2.13}$$

Substituting (2.13) for $\boldsymbol{y}(\tau)$ in $\int_{\mathcal{K}} \boldsymbol{y}^\top(\tau)U\boldsymbol{y}(\tau)\mathrm{d}\tau$ gives

$$\begin{aligned}\int_{\mathcal{K}} \boldsymbol{y}^\top(\tau)U\boldsymbol{y}(\tau)\mathrm{d}\tau = \int_{\mathcal{K}} \boldsymbol{x}^\top(\tau)U\boldsymbol{x}(\tau)\mathrm{d}\tau + \mathbf{z}^\top \int_{\mathcal{K}} (\mathsf{G} \otimes I_n)^\top G(\tau)UG^\top(\tau)(\mathsf{G} \otimes I_n)\mathrm{d}\tau\mathbf{z}\\ - 2\int_{\mathcal{K}} \boldsymbol{x}^\top(\tau)UG^\top(\tau)\mathrm{d}\tau(\mathsf{G} \otimes I_n)\mathbf{z}\end{aligned} \tag{2.14}$$

where $\mathbf{z} := \int_{\mathcal{K}} G(\theta)\boldsymbol{x}(\theta)\mathrm{d}\theta$.

Now apply (2.1a) to the product term $UG^\top(\tau)$ in (2.14) and consider the fact that $G(\tau) = \boldsymbol{g}(\tau) \otimes I_n$, we have

$$U\left(\boldsymbol{g}^\top(\tau) \otimes I_n\right) = U\left(\boldsymbol{g}^\top(\tau) \otimes I_n\right)\left(\boldsymbol{g}^\top(\tau) \otimes I_n\right)(I_d \otimes U) = G^\top(\tau)(I_d \otimes U). \tag{2.15a}$$

Applying (2.15a) to the term $\int_{\mathcal{K}} \boldsymbol{x}^\top(\tau)UG^\top(\tau)\mathrm{d}\tau(\mathsf{G} \otimes I_n)\mathbf{z}$ in (2.14) yields

$$\int_{\mathcal{K}} \boldsymbol{x}^\top(\tau)UG^\top(\tau)\mathrm{d}\tau(\mathsf{G} \otimes I_n)\mathbf{z} = \mathbf{z}^\top(I_d \otimes U)(\mathsf{G} \otimes I_n)\mathbf{z} = \mathbf{z}^\top(\mathsf{G} \otimes U)\mathbf{z}. \tag{2.15b}$$

Furthermore, by (2.15a) and the fact that $\mathsf{G} = \mathsf{G}^\top$, the term $\int_{\mathcal{K}} (\mathsf{G} \otimes I_n)^\top G(\tau)UG^\top(\tau)(\mathsf{G} \otimes I_n)\mathrm{d}\tau$ in (2.14) can be reformulated as

$$\begin{aligned}\int_{\mathcal{K}} (\mathsf{G} \otimes I_n)G(\tau)G^\top(\tau)(I_d \otimes U)(\mathsf{G} \otimes I_n)\mathrm{d}\tau &= \int_{\mathcal{K}} (\mathsf{G} \otimes I_n)G(\tau)G^\top(\tau)(\mathsf{G} \otimes U)\mathrm{d}\tau\\ &= (\mathsf{G} \otimes I_n)\int_{\mathcal{K}} (\boldsymbol{g}(\tau) \otimes I_n)\left(\boldsymbol{g}^\top(\tau) \otimes I_n\right)\mathrm{d}\tau(\mathsf{G} \otimes U)\\ &= (\mathsf{G} \otimes I_n)\left(\int_{\mathcal{K}} \boldsymbol{g}(\tau)\boldsymbol{g}^\top(\tau)\mathrm{d}\tau \otimes I_n\right)(\mathsf{G} \otimes U) = \mathsf{G} \otimes U.\end{aligned} \tag{2.15c}$$

By (2.15c) and (2.15b), (2.14) is rewritten as the form

$$\int_{\mathcal{K}} \boldsymbol{y}^\top(\tau)U\boldsymbol{y}(\tau)\mathrm{d}\tau = \int_{\mathcal{K}} \boldsymbol{x}^\top(\tau)U\boldsymbol{x}(\tau)\mathrm{d}\tau - \int_{\mathcal{K}} \boldsymbol{x}^\top(\tau)G^\top(\tau)\mathrm{d}\tau\,(\mathsf{G} \otimes U)\int_{\mathcal{K}} G(\tau)\boldsymbol{x}(\tau)\mathrm{d}\tau. \tag{2.16}$$

Given $U \succeq 0$, (2.12) can now be derived from (2.16). This finishes the proof. ∎

**Remark 2.2.** By the conclusions of [566, Theorem 7.2.10], the condition in (2.11) indicates that the functions in $\boldsymbol{g}(\cdot)$ are linearly independent in a Lebesgue sense.

**Remark 2.3.** Given $\mathcal{K} = [-r, 0], r > 0$ and $\boldsymbol{g}(\tau) = \boldsymbol{\ell}_d(\tau) = \mathbf{Col}_{i=0}^{d}\, \ell_i(\tau)$ with the Legendre polynomials[1] [467, 469, 470, 572]

$$\ell_d(\tau) := \sum_{k=0}^{d} \binom{d}{k}\binom{d+k}{k}\left(\frac{\tau}{r}\right)^k, \quad \forall d \in \mathbb{N}\cup\{0\}, \quad \forall \tau \in [-r, 0], \tag{2.17}$$

then (2.12) holds with $\mathsf{G}^{-1} = \bigoplus_{i=0}^{d} \frac{r}{2i+1}$ in this case. This demonstrates the fact that (2.12) can be regarded as a generalization of the Bessel-Legendre Inequality proposed in [467, 469]. Furthermore, (2.12) can also be considered as a generalization of the results in [573]. Finally, it is worth remarking that a discrete version of (2.12)

$$\begin{gathered}\sum_{k\in\mathcal{J}} \boldsymbol{x}^\top(k) U \boldsymbol{x}(k) \geq [*]\left(\mathsf{G}\otimes U\right)\left(\sum_{k\in\mathcal{J}}\left(\boldsymbol{g}(k)\otimes I_n\right)\boldsymbol{x}(k)\right), \\ \mathsf{G}^{-1} = \sum_{k\in\mathcal{J}} \boldsymbol{g}(k)\boldsymbol{g}^\top(k) \succ 0, \quad \mathcal{J}\subseteq\mathbb{Z}\end{gathered} \tag{2.18}$$

can be easily derived given the connections between integrations and summations.

The following Projection Lemma is applied in the derivation of convex synthesis conditions in throughout this book.

**Lemma 2.4** (Projection Lemma). *[562–564] Given $n; p; q \in \mathbb{N}$, $\Pi \in \mathbb{S}^n, P \in \mathbb{R}^{q\times n}, Q \in \mathbb{R}^{p\times n}$, there exists $\Theta \in \mathbb{R}^{p\times q}$ such that the following two propositions are equivalent:*

$$\Pi + P^\top \Theta^\top Q + Q^\top \Theta P \prec 0, \tag{2.19a}$$

$$P_\perp^\top \Pi P_\perp \prec 0 \;\; \text{and} \;\; Q_\perp^\top \Pi Q_\perp \prec 0, \tag{2.19b}$$

*where the columns of $P_\perp$ and $Q_\perp$ contain bases of null space of matrix $P$ and $Q$, respectively, which means that $PP_\perp = \mathsf{O}$ and $QQ_\perp = \mathsf{O}$.*

*Proof.* Refer to the proof of [564, Lemma 3.1] and [452, Lemma C.12.1]. ■

## 2.4 Main Results on the Dissipative Stabilization

In this section, the main results on controller synthesis in view of dissipativity are set forth in Theorem 2.1 and Theorem 2.2 and Algorithm 1. The proposed methods are based on the construction of

$$\mathsf{v}(\boldsymbol{\chi}_t(\cdot)) := \begin{bmatrix} \boldsymbol{x}(t) \\ \int_{-r}^{0} F(\tau)\boldsymbol{x}(t+\tau)\mathrm{d}\tau \end{bmatrix}^\top \begin{bmatrix} P & Q \\ * & R \end{bmatrix} \begin{bmatrix} \boldsymbol{x}(t) \\ \int_{-r}^{0} F(\tau)\boldsymbol{x}(t+\tau)\mathrm{d}\tau \end{bmatrix} + \int_{-r}^{0} \boldsymbol{x}^\top(t+\tau)\left(S + (\tau+r)U\right)\boldsymbol{x}(t+\tau)\mathrm{d}\tau \tag{2.20}$$

via the application of (2.12), where $\boldsymbol{x}(t)$ satisfies (2.5) and $F(\tau) = \boldsymbol{f}(\tau)\otimes I_n$ with $\boldsymbol{f}(\cdot)$ in Assumption 2.1, and the values of $P \in \mathbb{S}^n, Q \in \mathbb{R}^{n\times dn}, R \in \mathbb{S}^{dn}$, $S; U \in \mathbb{S}^n$ are to be constructed numerically. Note that (2.20) can be considered as a particular case of the Complete Lyapunov-Krasovskiĭ Functional (CKLF) in [178] whose infinite-dimensional variables are parameterized by $Q$ and $R$ here via the integral term $\int_{-r}^{0} F(\tau)\boldsymbol{x}(t+\tau)\mathsf{d}\tau$.

[1] The expression of Legendre polynomials here is obtained by setting $\alpha = \beta = 0$ for the expersion of Jacobi polynomials over $[-r, 0]$. See (4.6) for more details.

In the following theorem, sufficient conditions for the existence of a state feedback controller which stabilizes (2.5) and takes into account the supply rate (2.10) are derived in terms of matrix inequalities.

**Theorem 2.1.** *Given $\boldsymbol{f}(\cdot)$ and $M$ in (2.1), then the CLS in (2.5) with supply rate function (2.10) is dissipative and the trivial solution $\boldsymbol{x}(t) \equiv \mathbf{0}_n$ of (2.5) with $\boldsymbol{w}(t) \equiv \mathbf{0}_q$ is globally asymptotically stable if there exist $K \in \mathbb{R}^{n\times p}$ and $P \in \mathbb{S}^n$, $Q \in \mathbb{R}^{n\times dn}$, $R \in \mathbb{S}^{dn}$ and $S; U \in \mathbb{S}^n$ such that*

$$\begin{bmatrix} P & Q \\ * & R+\mathsf{F}\otimes S \end{bmatrix} \succ 0, \quad S \succ 0, \quad U \succ 0, \tag{2.21a}$$

$$\boldsymbol{\Phi} + \mathsf{Sy}\left(\mathbf{P}^\top \boldsymbol{\Pi}\right) \prec 0, \tag{2.21b}$$

*where* $\mathbf{P} := \begin{bmatrix} P & \mathsf{O}_n & Q & \mathsf{O}_{n,q} & \mathsf{O}_{n,m} \end{bmatrix}$, $\boldsymbol{\Pi} := \begin{bmatrix} \mathbf{A} + \mathbf{B}_1\left[(I_{2+d}\otimes K)\oplus \mathsf{O}_q\right] & \mathsf{O}_{n,m} \end{bmatrix}$ *and*

$$\boldsymbol{\Phi} := \mathsf{Sy}\left( \begin{bmatrix} Q \\ \mathsf{O}_{n,dn} \\ R \\ \mathsf{O}_{(q+m),dn} \end{bmatrix} \begin{bmatrix} \mathbf{F} & \mathsf{O}_{dn,m} \end{bmatrix} \right) + \mathsf{Sy}\left( \begin{bmatrix} \mathsf{O}_{(2n+dn),m} \\ -J_2^\top \\ \widetilde{J} \end{bmatrix} \begin{bmatrix} \boldsymbol{\Sigma} & \mathsf{O}_m \end{bmatrix} \right)$$
$$+ [S+rU]\oplus[-S]\oplus[-\mathsf{F}\otimes U]\oplus(-J_3)\oplus J_1 \tag{2.22a}$$

*with* $\mathsf{F}^{-1} := \int_{-r}^{0} \boldsymbol{f}(\tau)\boldsymbol{f}^\top(\tau)\mathsf{d}\tau$ *and*

$$\mathbf{F} = \begin{bmatrix} F(0) & F(-r) & -M\otimes I_n & \mathsf{O}_{dn,q} \end{bmatrix}, \quad \boldsymbol{\Sigma} = \mathbf{C} + \mathbf{B}_2\left[(I_{2+d}\otimes K)\oplus \mathsf{O}_q\right] \tag{2.22b}$$

*with $F(\tau)$ given in (2.6c). Moreover, the matrices $\mathbf{A}, \mathbf{B}_1, \mathbf{B}_2$ and $\mathbf{C}$, which contains the parameters of open-loop system (2.2), are given in (2.6a) and (2.6b).*

*Proof.* To begin with, note that $\mathsf{s}(\boldsymbol{z}(t), \boldsymbol{w}(t))$ in (2.10) can be written as

$$\mathsf{s}(\boldsymbol{z}(t), \boldsymbol{w}(t)) = \boldsymbol{z}^\top(t)\widetilde{J}^\top J_1^{-1}\widetilde{J}\boldsymbol{z}(t) + \mathsf{Sy}\left[\boldsymbol{z}^\top(t)J_2\boldsymbol{w}(t)\right] + \boldsymbol{w}^\top(t)J_3\boldsymbol{w}(t), \tag{2.23a}$$

where only the term $\boldsymbol{z}^\top(t)\widetilde{J}^\top J_1^{-1}\widetilde{J}\boldsymbol{z}(t)$ introduces nonlinearity with respect to the undetermined variables in (2.5). Substituting (2.5) for $\boldsymbol{z}(t)$ in $\boldsymbol{z}^\top(t)\widetilde{J}^\top J_1^{-1}\widetilde{J}\boldsymbol{z}(t)$ produces

$$\boldsymbol{z}^\top(t)\widetilde{J}^\top J_1^{-1}\widetilde{J}\boldsymbol{z}(t) = \boldsymbol{\chi}^\top(t)\boldsymbol{\Sigma}^\top\widetilde{J}^\top J_1^{-1}\widetilde{J}\boldsymbol{\Sigma}\boldsymbol{\chi}(t), \tag{2.23b}$$

where $\boldsymbol{\Sigma}$ is given in (2.22b) and $\boldsymbol{\chi}(t)$ is defined in (2.6c).

Now differentiate the functional $\mathsf{v}(\mathbf{x}_t(\cdot))$ in (2.20) along the trajectory of the CLS in (2.5) and consider (2.3a), (2.23a), (2.23b) and the relation

$$\frac{\mathsf{d}}{\mathsf{d}t}\int_{-r}^{0} F(\tau)\boldsymbol{x}(t+\tau)\mathsf{d}\tau = F(0)\boldsymbol{x}(t) - F(-r)\boldsymbol{x}(t-r) - (M\otimes I_n)\int_{-r}^{0} F(\tau)\boldsymbol{x}(t+\tau)\mathsf{d}\tau = \mathbf{F}\boldsymbol{\chi}(t) \tag{2.24}$$

with $\mathbf{F}$ in (2.22b). Then it yields

$$\forall t \geq t_0, \quad \dot{\mathsf{v}}(\mathbf{x}_t(\cdot)) - \mathsf{s}(\boldsymbol{z}(t), \boldsymbol{w}(t))$$
$$= \boldsymbol{\chi}^\top(t)\,\mathsf{Sy}\left( \begin{bmatrix} I_n & \mathsf{O}_{n,dn} \\ \mathsf{O}_n & \mathsf{O}_{n,dn} \\ \mathsf{O}_{dn,n} & I_{dn} \\ \mathsf{O}_{q,n} & \mathsf{O}_{q,dn} \end{bmatrix} \begin{bmatrix} P & Q \\ * & R \end{bmatrix} \begin{bmatrix} \mathbf{A} + \mathbf{B}_1\left[(I_{2+d}\otimes K)\oplus \mathsf{O}_q\right] \\ \mathbf{F} \end{bmatrix} \right)\boldsymbol{\chi}(t)$$

$$
+\boldsymbol{x}^\top(t)\left(S+rU\right)\boldsymbol{x}(t)-\boldsymbol{x}^\top(t-r)S\boldsymbol{x}(t-r)-\int_{-r}^{0}\boldsymbol{x}^\top(t+\tau)U\boldsymbol{x}(t+\tau)\mathrm{d}\tau-\boldsymbol{w}^\top(t)J_3\boldsymbol{w}(t)
$$

$$
-\boldsymbol{\chi}^\top(t)\,\mathsf{Sy}\left(\boldsymbol{\Sigma}^\top\widetilde{J}^\top J_1^{-1}\widetilde{J}\boldsymbol{\Sigma}+\begin{bmatrix}\mathsf{O}_{(2n+dn),m}\\ J_2^\top\end{bmatrix}\boldsymbol{\Sigma}\right)\boldsymbol{\chi}(t). \tag{2.25}
$$

Assume $U\succeq 0$, then applying (2.12) to the integral $\int_{-r}^{0}\boldsymbol{x}^\top(t+\tau)U\boldsymbol{x}(t+\tau)\mathrm{d}\tau$ in (2.25) with $\boldsymbol{g}(\tau)=\boldsymbol{f}(\tau)$ produces

$$
\forall t\geq t_0,\quad \int_{-r}^{0}\boldsymbol{x}^\top(t+\tau)U\boldsymbol{x}(t+\tau)\mathrm{d}\tau\geq\int_{-r}^{0}\boldsymbol{x}^\top(t+\tau)F^\top(\tau)\mathrm{d}\tau\,(\mathsf{F}\otimes U)\int_{-r}^{0}F(\tau)\boldsymbol{x}(t+\tau)\mathrm{d}\tau, \tag{2.26}
$$

where $\mathsf{F}^{-1}=\int_{-r}^{0}\boldsymbol{f}(\tau)\boldsymbol{f}^\top(\tau)\mathrm{d}\tau$. Now apply (2.26) to (2.25), then we have

$$
\forall t\geq t_0,\quad \dot{\mathsf{v}}(\mathbf{x}_t(\cdot))-\mathsf{s}(\boldsymbol{z}(t),\boldsymbol{w}(t))\leq\boldsymbol{\chi}^\top(t)\left(\boldsymbol{\Psi}-\boldsymbol{\Sigma}^\top\widetilde{J}^\top J_1^{-1}\widetilde{J}\boldsymbol{\Sigma}\right)\boldsymbol{\chi}(t) \tag{2.27a}
$$

with $\boldsymbol{\chi}(t)$ defined in (2.6c), where

$$
\begin{aligned}
\boldsymbol{\Psi}:=\mathsf{Sy}&\left(\begin{bmatrix}P & Q\\ \mathsf{O}_n & \mathsf{O}_{n,nd}\\ Q^\top & R\\ \mathsf{O}_{q,n} & \mathsf{O}_{q,nd}\end{bmatrix}\begin{bmatrix}\mathbf{A}+\mathbf{B}_1\left[(I_{2+d}\otimes K)\oplus\mathsf{O}_q\right]\\ \mathbf{F}\end{bmatrix}\right)\\
&-\left([-S-rU]\oplus S\oplus(\mathsf{F}\otimes U)\oplus J_3\right)-\mathsf{Sy}\left(\begin{bmatrix}\mathsf{O}_{(2n+dn),m}\\ J_2^\top\end{bmatrix}\left(\mathbf{C}+\mathbf{B}_2\left[(I_{2+d}\otimes K)\oplus\mathsf{O}_q\right]\right)\right)
\end{aligned} \tag{2.27b}
$$

containing all the matrix terms induced by $\dot{\mathsf{v}}(\mathbf{x}_t(\cdot))-\mathsf{s}(\boldsymbol{z}(t),\boldsymbol{w}(t))$ in (2.25) excluding (2.23b). Now considering the property of positive definite matrices, it becomes evident that if

$$
\boldsymbol{\Psi}-\boldsymbol{\Sigma}^\top\widetilde{J}^\top J_1^{-1}\widetilde{J}\boldsymbol{\Sigma}\prec 0,\quad U\succeq 0 \tag{2.28}
$$

is satisfied then the dissipative inequality in (2.8) : $\forall t\geq t_0,\ \dot{\mathsf{v}}(\mathbf{x}_t(\cdot))-\mathsf{s}(\boldsymbol{z}(t),\boldsymbol{w}(t))\leq 0$ holds. Moreover, utilizing the Schur complement to the first inequality in (2.28) concludes that

$$
\begin{bmatrix}\boldsymbol{\Psi} & \boldsymbol{\Sigma}^\top\widetilde{J}^\top\\ * & J_1\end{bmatrix}=\mathbf{P}^\top\boldsymbol{\Pi}+\boldsymbol{\Phi}\prec 0,\quad U\succeq 0 \tag{2.29}
$$

holds if and only if (2.28) holds given $J_1^{-1}\prec 0$ in (2.10), where $\boldsymbol{\Phi}$ is defined in (2.22a). For this reason, it follows that there exists a functional (2.20) satisfying (2.8) if (2.29) has feasible solutions. By considering the properties of positive definite matrices, it is obvious that given (2.29) holds then $\exists\epsilon_3>0$, $\forall t\geq t_0$ $\dot{\mathsf{v}}(\mathbf{x}_t(\cdot))\leq-\epsilon_3\left\|\boldsymbol{x}(t)\right\|_2^2=-\epsilon_3\left\|\mathbf{x}_t(0)\right\|_2^2$ where $\boldsymbol{x}(t)$ here satisfies (2.5) with $\boldsymbol{w}(t)\equiv\mathbf{0}_q$. Let $t=t_0$ with the fact[2] that $\forall\theta\in[-r,0]$, $\mathbf{x}_{t_0}(\theta)=\boldsymbol{x}(t_0+\theta)=\boldsymbol{\phi}(\theta)$ in (2.5), then we have $\exists\epsilon_3>0$, $\dot{\mathsf{v}}(\mathbf{x}_{t_0}(\cdot))=\dot{\mathsf{v}}(\boldsymbol{\phi}(\cdot))\leq-\epsilon_3\left\|\mathbf{x}_{t_0}(0)\right\|_2^2=-\epsilon_3\left\|\boldsymbol{\phi}(0)\right\|_2^2$ for any $\boldsymbol{\phi}(\cdot)\in\mathcal{C}\left([-r,0];\mathbb{R}^n\right)$ in (2.5), which gives (2.7b). Thus we can conclude that if (2.29) holds then (2.20) satisfies (2.8) and (2.7b).

Now we start to prove that (2.20) satisfies (2.7a) if (2.21a) holds. Knowing that $\forall\theta\in[-r,0]$, $\mathbf{x}_{t_0}(\theta)=\boldsymbol{x}(t_0+\theta)=\boldsymbol{\phi}(\theta)$ in (2.5), let $t=t_0$ in (2.20) and $S\succeq 0$. By applying (2.12) with $\boldsymbol{g}(\tau)=\boldsymbol{f}(\tau)$ and $U=S$ to the integral $\int_{-r}^{0}\boldsymbol{x}^\top(t+\tau)S\boldsymbol{x}(t+\tau)\mathrm{d}\tau$ in (2.20) with $t=t_0$, it follows that

$$
\int_{-r}^{0}\boldsymbol{\phi}^\top(\tau)S\boldsymbol{\phi}(\tau)\mathrm{d}\tau\geq\left[\int_{-r}^{0}F(\tau)\boldsymbol{\phi}(\tau)\mathrm{d}\tau\right]^\top(\mathsf{F}\otimes S)\int_{-r}^{0}F(\tau)\boldsymbol{\phi}(\tau)\mathrm{d}\tau \tag{2.30}
$$

[2] Again $\boldsymbol{x}(t)$ here satisfies (2.5) with $\boldsymbol{w}(t)\equiv\mathbf{0}_q$

holds for any $\phi(\cdot) \in \mathcal{C}([-r,0];\mathbb{R}^n)$ in (2.5). Apply (2.30) to (2.20) with $t = t_0$ and $S \succ 0$, we have

$$\mathsf{v}(\mathbf{x}_{t_0}(\cdot)) = \mathsf{v}(\phi(\cdot)) \geq \begin{bmatrix} \phi(0) \\ \int_{-r}^{0} F(\tau)\phi(\tau)\mathsf{d}\tau \end{bmatrix}^\top \begin{bmatrix} P & Q \\ * & R + \mathsf{F} \otimes S \end{bmatrix} \begin{bmatrix} \phi(0) \\ \int_{-r}^{0} F(\tau)\phi(\tau)\mathsf{d}\tau \end{bmatrix} + \int_{-r}^{0} (\tau + r)\phi^\top(\tau) U \phi(\tau)\mathsf{d}\tau \quad (2.31)$$

for any $\phi(\cdot) \in \mathcal{C}([-r,0];\mathbb{R}^n)$ in (2.5). According to the property of positive definite matrices considering the structure of (2.31), it is obvious to conclude that if (2.21a) are satisfied, then $\exists \epsilon_1 > 0$ : $\mathsf{v}(\phi(\cdot)) \geq \epsilon_1 \|\phi(0)\|_2$ for any $\phi(\cdot) \in \mathcal{C}([-r,0];\mathbb{R}^n)$ in (2.5). Furthermore, for (2.20) with $t = t_0$ it follows that $\exists \lambda_1, \lambda_2 > 0$ such that

$$\begin{aligned} &\mathsf{v}(\phi(\cdot)) \leq \begin{bmatrix} \phi(0) \\ \int_{-r}^{0} F(\tau)\phi(\tau)\mathsf{d}\tau \end{bmatrix}^\top \lambda_1 \begin{bmatrix} \phi(0) \\ \int_{-r}^{0} F(\tau)\phi(\tau)\mathsf{d}\tau \end{bmatrix} + \lambda_1 \int_{-r}^{0} \sup_{-r \leq \tau \leq 0} \|\phi(\tau)\|_2^2 \,\mathsf{d}\tau \\ &\leq \phi^\top(0)\lambda_1\phi(0) + \int_{-r}^{0} \phi^\top(\tau) F^\top(\tau)\mathsf{d}\tau (\lambda_2 \mathsf{F} \otimes I_n) \int_{-r}^{0} F(\tau)\phi(\tau)\mathsf{d}\tau + r\lambda_1 \|\phi(\cdot)\|_\infty^2 \\ &\leq \phi^\top(0)\lambda_1\phi(0) + \int_{-r}^{0} \phi^\top(\tau)\lambda_2\phi(\tau)\mathsf{d}\tau + r\lambda_1 \|\phi(\cdot)\|_\infty^2 \leq (\lambda_1 + r\lambda_2 + r\lambda_1) \|\phi(\cdot)\|_\infty^2 \end{aligned} \quad (2.32)$$

for any $\phi(\cdot) \in \mathcal{C}([-r,0];\mathbb{R}^n)$ in (2.5), which is derived via (2.12) with $\boldsymbol{g}(\tau) = \boldsymbol{f}(\tau)$ and the property $\forall X \in \mathbb{S}^n, \exists \lambda > 0, \lambda I_n - X \succ 0$. This shows that (2.7a) can be satisfied for any functional in the form of (2.20). Therefore, we have shown that if (2.21a) holds then the functional in (2.20) satisfies (2.7a).

As we have proved that (2.20) satisfies (2.8) and (2.7b) for some $\epsilon_3 > 0$ if (2.29) holds, we can conclude that feasible solutions to (2.21a) and (2.21b) imply the existence of (2.20) and $\epsilon_1; \epsilon_2; \epsilon_3 > 0$ satisfying the dissipativity condition in (2.8) and the stability criteria in (2.7a) and (2.7b). ■

As (2.21b) is bilinear (non-convex) with nonzeros $B_1, B_2, B_3$, it cannot be computed using convex SDP solvers. To address this issue, we derive the following theorem based on the application of the Projection Lemma, which has been utilized in the context of SDP to decouple the products between decision variables [574–577]. The theorem provides a convex synthesis condition for dissipative stabilization with some given parameters[3], which can be solved using standard SDP solvers.

**Theorem 2.2.** *Let $\boldsymbol{f}(\cdot)$ and $M$ in (2.1) and $\{\alpha_i\}_{i=1}^{2+d} \subset \mathbb{R}$ be given, then the CLS in (2.5) with supply rate function (2.10) is dissipative and the trivial solution $\boldsymbol{x}(t) \equiv \mathbf{0}_n$ of (2.5) with $\boldsymbol{w}(t) \equiv \mathbf{0}_q$ is exponentially stable if there exist $\acute{P} \in \mathbb{S}^n$, $X \in \mathbb{R}^{n\times n}$, $V \in \mathbb{R}^{p\times n}$, $\acute{Q} \in \mathbb{R}^{n\times dn}$, $\acute{R} \in \mathbb{S}^{dn}$, $\acute{S}; \acute{U} \in \mathbb{S}^n$ such that*

$$\begin{bmatrix} \acute{P} & \acute{Q} \\ * & \acute{R} + \mathsf{F} \otimes \acute{S} \end{bmatrix} \succ 0, \quad \acute{S} \succ 0, \quad \acute{U} \succ 0, \quad (2.33a)$$

$$\acute{\boldsymbol{\Theta}} := \mathsf{Sy}\left( \begin{bmatrix} I_n \\ \mathbf{Col}_{i=1}^{2+d}\, \alpha_i I_n \\ \mathsf{O}_{(q+m),n} \end{bmatrix} \begin{bmatrix} -X & \acute{\boldsymbol{\Pi}} \end{bmatrix} \right) + \begin{bmatrix} \mathsf{O}_n & \acute{\mathbf{P}} \\ * & \acute{\boldsymbol{\Phi}} \end{bmatrix} \prec 0, \quad (2.33b)$$

[3]It is illustrated later in Remark 2.4 that it possible to only adjust the value of a single parameter with other parameters assigned to be zeros

*hold, where* $\mathsf{F}^{-1} := \int_{-r}^{0} \boldsymbol{f}(\tau)\boldsymbol{f}^{\top}(\tau)\mathsf{d}\tau$ *and*

$$
\begin{aligned}
&\acute{\mathbf{P}} := \begin{bmatrix} \mathsf{O}_n & \acute{P} & \mathsf{O}_n & \acute{Q} & \mathsf{O}_{n,q} & \mathsf{O}_{n,m} \end{bmatrix}, \quad \acute{\mathbf{\Pi}} := \begin{bmatrix} \mathbf{A}\left(I_{2+d} \otimes X\right) + \mathbf{B}_1\left(I_{2+d} \otimes V\right) & \mathsf{O}_{n,m} \end{bmatrix}, \\
&\acute{\mathbf{\Phi}} := \mathsf{Sy}\left( \begin{bmatrix} \acute{Q} \\ \mathsf{O}_{n,dn} \\ \acute{R} \\ \mathsf{O}_{q,dn} \end{bmatrix} \begin{bmatrix} \mathbf{F} & \mathsf{O}_{dn,m} \end{bmatrix} \right) + \left(\acute{S} + r\acute{U}\right) \oplus \left(-\acute{S}\right) \oplus \left(-\mathsf{F} \otimes \acute{U}\right) \oplus (-J_3) \oplus J_1 \\
&\quad + \mathsf{Sy}\left( \begin{bmatrix} \mathsf{O}_{(2n+dn),m} \\ -J_2^{\top} \\ \widetilde{J} \end{bmatrix} \begin{bmatrix} \acute{\mathbf{\Sigma}} & \mathsf{O}_m \end{bmatrix} \right).
\end{aligned} \tag{2.34}
$$

*with* $\acute{\mathbf{\Sigma}} = \begin{bmatrix} C_1X + B_4V & C_2X + B_5V & C_3(I_d \otimes X) + B_2\left(I_d \otimes V\right) & D_2 \end{bmatrix}$ *and* $\mathbf{F}$ *in* (2.22b). *Furthermore the controller parameter is calculated via* $K = VX^{-1}$ *where the invertibility of* $X$ *is guaranteed by* (2.33b).

*Proof.* It is obvious to see that the root of non-convexity in (2.21b) is the matrix products in $\mathsf{Sy}\left(\mathbf{P}^{\top}\mathbf{\Pi}\right)$ concerning $P$ and $Q$. Now we can reformulate (2.21b) into the form

$$
\begin{bmatrix} \mathbf{\Pi} \\ I_{2n+dn+q+m} \end{bmatrix}^{\top} \begin{bmatrix} \mathsf{O}_n & \mathbf{P} \\ * & \mathbf{\Phi} \end{bmatrix} \begin{bmatrix} \mathbf{\Pi} \\ I_{2n+dn+q+m} \end{bmatrix} \prec 0. \tag{2.35}
$$

It is obvious that the structure of (2.35) is consistent with the structure of the inequalities in (2.19b) as part of the statements in Lemma 2.4. However, since there are two matrix inequalities in (2.19b), we must construct a new matrix inequality. Consider the matrix inequality

$$
\Upsilon^{\top} \begin{bmatrix} \mathsf{O}_n & \mathbf{P} \\ * & \mathbf{\Phi} \end{bmatrix} \Upsilon \prec 0 \tag{2.36a}
$$

with $\Upsilon^{\top} := \begin{bmatrix} \mathsf{O}_{(q+m),(3n+dn)} & I_{q+m} \end{bmatrix}$, which is further simplified into

$$
\Upsilon^{\top} \begin{bmatrix} \mathsf{O}_n & \mathbf{P} \\ * & \mathbf{\Phi} \end{bmatrix} \Upsilon = \begin{bmatrix} -J_3 - \mathsf{Sy}(D_2^{\top}J_2) & D_2^{\top}\widetilde{J} \\ * & J_1 \end{bmatrix} \prec 0. \tag{2.36b}
$$

As a matter of fact, (2.36b) is identical to the matrix obtained by extracting the $2 \times 2$ block matrix at the bottom right of (2.35). As a consequence, one can conclude that (2.36b) is automatically satisfied if condition (2.35) holds. So (2.36b) introduces no additional conservatism to the original condition.

To utilize Lemma 2.4, we meed to construct two matrices $\mathbf{U}$, $\mathbf{Y}$ in line with the inequalities in (2.19b) satisfying

$$
\mathbf{U}_{\perp} = \Upsilon, \quad \mathbf{Y}_{\perp} = \begin{bmatrix} \mathbf{\Pi} \\ I_{2n+dn+q+m} \end{bmatrix} \tag{2.37}
$$

where $\Upsilon$ and $\begin{bmatrix} \mathbf{\Pi}^{\top} & I_{2n+dn+q+m} \end{bmatrix}^{\top}$ contain bases of the null spaces of $\mathbf{U}$ and $\mathbf{Y}$, respectively. For $\Upsilon^{\top} := \begin{bmatrix} \mathsf{O}_{(q+m),(3n+dn)} & I_{q+m} \end{bmatrix}$, we have $\operatorname{rank}(\Upsilon) = q + m$ which gives $\operatorname{rank}(\mathbf{U}) = 3n + dn$ by the rank nullity theorem. Likewise, we can conclude that $\operatorname{rank}(\mathbf{Y}) = n$. Without loss of generality., let

$$
\mathbf{Y} := \begin{bmatrix} -I_n & \mathbf{\Pi} \end{bmatrix}, \quad \mathbf{U} := \begin{bmatrix} I_{3n+dn} & \mathsf{O}_{(3n+dn),(q+m)} \end{bmatrix} \tag{2.38}
$$

by which we have

$$
\mathbf{Y}\mathbf{Y}_{\perp} = \mathbf{Y} \begin{bmatrix} \mathbf{\Pi} \\ I_{2n+dn+q+m} \end{bmatrix} = \mathsf{O}_{n,(2n+dn+q+m)}, \quad \mathbf{U}\mathbf{U}_{\perp} = \mathsf{O}_{(3n+dn),(q+m)}.
$$

Now the choice of $\mathbf{U}$ and $\mathbf{Y}$ in (2.38) satisfies $\text{rank}(\mathbf{U}) = 3n+dn$, $\text{rank}(\mathbf{Y}) = n$, and $\Upsilon$, $\begin{bmatrix}\mathbf{\Pi}^\top & I_{2n+dn+q+m}\end{bmatrix}^\top$ contain bases of the null spaces of $\mathbf{U}$ and $\mathbf{Y}$, respectively. Now one can apply Lemma 2.4 to (2.35) and (2.36b) yields that (2.35) and (2.36b) hold if and only if

$$\exists \mathbf{W} \in \mathbb{R}^{(3n+dn)\times n} : \mathsf{Sy}\begin{bmatrix}\mathbf{U}^\top \mathbf{W}\mathbf{Y}\end{bmatrix} + \begin{bmatrix} \mathsf{O}_n & \mathbf{P} \\ * & \mathbf{\Phi} \end{bmatrix} \prec 0. \tag{2.39}$$

with $\mathbf{U}$ and $\mathbf{Y}$ defined in (2.38). To ultimately produce convex synthesis conditions, let $\{\alpha_i\}_{i=1}^{d+2} \subset \mathbb{R}$ be given and

$$\mathbf{W} := \begin{bmatrix} W^\top & \mathsf{Col}_{i=1}^{d+2}\, \alpha_i W \end{bmatrix}^\top \tag{2.40}$$

with $W \in \mathbb{R}^{n\times n}$. Substituting the expressions in (2.40) and (2.38) for $\mathbf{W}, \mathbf{Y}, \mathbf{U}$ in (2.39) produces

$$\mathbf{\Theta} = \mathsf{Sy}\left(\mathbf{U}^\top \begin{bmatrix} W \\ \mathsf{Col}_{i=1}^{d+2}\, \alpha_i W \end{bmatrix} \mathbf{Y}\right) + \begin{bmatrix} \mathsf{O}_n & \mathbf{P} \\ * & \mathbf{\Phi} \end{bmatrix} \prec 0. \tag{2.41}$$

where the condition in (2.41) can be convexified via congruent transformations [578, page 12].

It should be noted that the inequality in (2.41) is no longer equivalent to (2.35), but rather it is a sufficient condition due to the structural constraints in (2.40). It is important to recognize that the invertibility of $W$ is implied by (2.41) as the expression $-W - W^\top$ is the only element in the first diagonal block of (2.41). Now by applying congruent transformations to (2.41) and (2.21a) with the fact that $W^{-1}$ is well defined, we have

$$\begin{gathered} \left((I_{3+d}\otimes X)\oplus I_{q+m}\right)^\top \mathbf{\Theta}\left((I_{3+d}\otimes X)\oplus I_{q+m}\right) \prec 0, \\ \left(I_{d+1}\otimes X^\top\right)\begin{bmatrix} P & Q \\ * & R + \mathsf{F}\otimes S \end{bmatrix}(I_{d+1}\otimes X) \succ 0,\quad X^\top S X \succ 0,\quad X^\top U X \succ 0 \end{gathered} \tag{2.42}$$

holds if and only if (2.41) and (2.21a) hold, where $X^\top := W^{-1}$. Moreover, letting

$$\begin{bmatrix} \acute{P} & \acute{Q} \\ * & \acute{R} \end{bmatrix} := \left(I_{1+d}\otimes X^\top\right)\begin{bmatrix} P & Q \\ * & R \end{bmatrix}(I_{1+d}\otimes X),\quad \begin{bmatrix}\acute{S} & \acute{U}\end{bmatrix} := X^\top\begin{bmatrix} SX & UX \end{bmatrix} \tag{2.43a}$$

and considering (2.42) with (2.1a) yields (2.33a) and

$$[*]\mathbf{\Theta}\left((I_{3+d}\otimes X)\oplus I_{q+m}\right) = \acute{\mathbf{\Theta}} = \mathsf{Sy}\left(\begin{bmatrix} I_n \\ \mathsf{Col}_{i=1}^{2+d}\, \alpha_i I_n \\ \mathsf{O}_{(q+m),n} \end{bmatrix}\begin{bmatrix} -X & \acute{\mathbf{\Pi}} \end{bmatrix}\right) + \begin{bmatrix} \mathsf{O}_n & \acute{\mathbf{P}} \\ * & \acute{\mathbf{\Phi}} \end{bmatrix} \prec 0 \tag{2.43b}$$

where $\acute{\mathbf{P}} = X\mathbf{P}\left[(I_{3+d}\otimes X)\oplus I_{q+m}\right] = \begin{bmatrix} \mathsf{O}_n & \acute{P} & \mathsf{O}_n & \acute{Q} & \mathsf{O}_{n,(q+m)} \end{bmatrix}$ and

$$\begin{aligned} \acute{\mathbf{\Pi}} = \mathbf{\Pi}\left[(I_{3+d}\otimes X)\oplus I_{q+m}\right] &= \begin{bmatrix} \mathbf{A}\left(I_{3+d}\otimes X\right) + \mathbf{B}_1\left(I_{3+d}\otimes KX\right) & \mathsf{O}_{n,m} \end{bmatrix} \\ &= \begin{bmatrix} \mathbf{A}\left(I_{3+d}\otimes X\right) + \mathbf{B}_1\left(I_{3+d}\otimes V\right) & \mathsf{O}_{n,m} \end{bmatrix} \end{aligned} \tag{2.43c}$$

with $V = KX$, and

$$\acute{\mathbf{\Phi}} = \left[\left(I_{2+d}\otimes X^\top\right)\oplus I_{q+m}\right]\mathbf{\Phi}\left[(I_{2+d}\otimes X)\oplus I_{q+m}\right] =$$

$$
\mathsf{Sy}\left(\begin{bmatrix} \acute{Q} \\ \mathsf{O}_{n,dn} \\ \acute{R} \\ \mathsf{O}_{(q+m),dn} \end{bmatrix} \begin{bmatrix} \mathbf{F} & \mathsf{O}_{dn,m} \end{bmatrix}\right) + \left(\left[\acute{S} + r\acute{U}\right] \oplus \left[-\acute{S}\right] \oplus \left[-\mathsf{F} \otimes \acute{U}\right] \oplus J_3 \oplus (-J_1)\right)
$$

$$
+ \mathsf{Sy}\left(\begin{bmatrix} \mathsf{O}_{(2n+dn),m} \\ -J_2^\top \\ \widetilde{J} \end{bmatrix} \begin{bmatrix} \acute{\mathbf{\Sigma}} & \mathsf{O}_{n,m} \end{bmatrix}\right) \tag{2.43d}
$$

with $\acute{\mathbf{\Sigma}} = \begin{bmatrix} C_1X + B_4V & C_2X + B_5V & C_3(I_d \otimes X) + B_2(I_d \otimes V) & D_2 \end{bmatrix}$, are the same given in (2.34). Note that (2.43b) is equivalent to (2.33b), and the form of $\acute{\mathbf{\Phi}}$ in (2.43d) is derived by considering the relation

$$
\begin{aligned}
\mathbf{F}\left[(I_{2+d} \otimes X) \oplus I_{q+m}\right] = \begin{bmatrix} \widehat{\mathbf{F}} \otimes I_n & \mathsf{O}_{dn,q} \end{bmatrix} \left[(I_{2+d} \otimes X) \oplus I_{q+m}\right] &= \begin{bmatrix} I_d\widehat{\mathbf{F}} \otimes XI_n & \mathsf{O}_{dn,q} \end{bmatrix} \\
&= (I_d \otimes X)\begin{bmatrix} \widehat{\mathbf{F}} \otimes I_n & \mathsf{O}_{dn,q} \end{bmatrix}.
\end{aligned}
$$

Moreover, since $-X - X^\top = 2X$ is the only term in the first diagonal block of $\acute{\mathbf{\Theta}}$ in (2.33b), it follows that $X$ is invertible if (2.33b) holds. This is consistent with the fact that $W$ being full rank is inferred from the inequality in (2.41).

As a consequence, the equivalence between (2.21a) and (2.33a) has been established. Furthermore, since (2.41) is equivalent to (2.33b) from which (2.21b) is inferred , it follows that there exist matrices such that (2.21a) and (2.21b) are satisfied if there exist feasible solutions to (2.33a) and (2.33b). Because feasible solutions to (2.21a) and (2.21b) always imply the existence of an LKF (2.20) satisfying (2.8) and (2.7a), (2.7b), this indicates that (2.33a) and (2.33b) also necessitates the existence of (2.20) satisfying the corresponding stability and dissipativity criteria. This finishes the proof. ■

**Remark 2.4.** Given the structure of (2.33b), some values of $\{\alpha_i\}_{i=1}^{2+d} \subset \mathbb{R}$ may have a greater impact on the feasibility of (2.33b) than others. For instance, $\alpha_1$ may be particularly critical as it directly influences the feasibility of the diagonal related to $A_1$ in (2.33b). A practical rule for assigning values to $\{\alpha_i\}_{i=1}^{2+d} \subset \mathbb{R}$ can be $\alpha_1 \in \mathbb{R}$ and $\alpha_i = 0, i = 2 \ldots 2 + d$, allowing one to adjust only the value $\alpha_1$ when using Theorem 2.2.

**Remark 2.5.** Even without considering dissipativity constraints, it is still possible to introduce slack variables in the same manner as in (2.33b) to solve a synthesis problem. However, in such a situation, the Projection Lemma may not be applicable as it may not be possible to construct two matrix inequalities as in (2.36b). Instead, a particular version of the Projection Lemma, called the Finsler Lemma [578], which only requires one inequality analogous to the structure of the inequalities in (2.19b), can be applied. By utilizing the notation of empty matrices, the corresponding synthesis condition derived through the application of the Finsler Lemma can be obtained by setting $m = q = 0$ in (2.33b).

**Remark 2.6.** It is important to stress that one can directly apply the Finsler lemma at the step of (2.35) to obtain a similar condition with more extra variables than (2.39). This observation suggests demonstrates that there is an advantage to applying the Projection Lemma over the Finsler Lemma, especially when (2.36b) can be constructed without introducing additional conservatism.

### 2.4.1 An inner convex approximation algorithm for Theorem 2.1

By specifying the values of $\{\alpha_i\}_{i=1}^{2+d} \subset \mathbb{R}$, a dissipative stabilizing controller gain $K$ can be obtained by solving the constraints in Theorem 2.2 using standard SDP solvers. However, the structure in (2.40) may introduce conservatism compared to the synthesis condition in Theorem 2.1. For this reason, it is preferable to have methods that can directly compute BMI (2.33b). For tackling generic non-convex SDP problems, various approaches [579] exist, such as augmented Lagrangian [580–588], sequential SDP [589–592] and primal-dual interior point methods [593–601], each with unique advantages and weaknesses.

In this subsection, an iterative algorithm is proposed based on the inner convex approximation scheme in [565] to solve Theorem 2.1, where each iteration involves a convex SDP problem. The algorithm can be considered as an example of sequential SDP, with initial values provided by feasible solutions to Theorem 2.2.

First of all, the non-convex inequality (2.21b) can be rewritten as

$$\mathfrak{D}(\mathbf{\Lambda}, K) = \mathsf{Sy}\left[\mathbf{P}^\top \mathbf{\Pi}\right] + \mathbf{\Phi} = \mathsf{Sy}\left[\mathbf{P}^\top \mathbf{B}\left[(I_{2+d} \otimes K) \oplus \mathsf{O}_{p+m}\right]\right] + \widehat{\mathbf{\Phi}} \prec 0 \tag{2.44}$$

given the structure of $\mathbf{\Pi}$ in (2.21b), where $\mathbf{B} = \begin{bmatrix} \mathbf{B}_1 & \mathsf{O}_{n,m} \end{bmatrix}$ and $\widehat{\mathbf{\Phi}} = \mathsf{Sy}\left(\mathbf{P}^\top \begin{bmatrix} \mathbf{A} & \mathsf{O}_{n,m} \end{bmatrix}\right) + \mathbf{\Phi}$ and $\mathbf{\Lambda} = \begin{bmatrix} P & Q \end{bmatrix}$ with $P$ and $Q$ in Theorem 2.1. Now consider the expression

$$\begin{aligned} \Delta\left(\mathbf{\Omega}, \widetilde{\mathbf{\Omega}}, \mathbf{\Gamma}, \widetilde{\mathbf{\Gamma}}\right) := \begin{bmatrix} \mathbf{\Omega}^\top - \widetilde{\mathbf{\Omega}}^\top & \mathbf{\Gamma}^\top - \widetilde{\mathbf{\Gamma}}^\top \end{bmatrix} \left[Z \oplus (I_n - Z)\right]^{-1} \begin{bmatrix} \mathbf{\Omega} - \widetilde{\mathbf{\Omega}} \\ \mathbf{\Gamma} - \widetilde{\mathbf{\Gamma}} \end{bmatrix} \\ + \mathsf{Sy}\left(\widetilde{\mathbf{\Omega}}^\top \mathbf{\Gamma} + \mathbf{\Omega}^\top \widetilde{\mathbf{\Gamma}} - \widetilde{\mathbf{\Omega}}^\top \widetilde{\mathbf{\Gamma}}\right) + \mathbf{T} \end{aligned} \tag{2.45}$$

where $\mathbf{\Omega}; \widetilde{\mathbf{\Omega}} \in \mathbb{R}^{n \times \mu}$, $\mathbf{\Gamma}; \widetilde{\mathbf{\Gamma}} \in \mathbb{R}^{n \times \mu}$ and $\mathbf{T} \in \mathbb{S}^{\mu}$, $Z \in \{X \in \mathbb{S}^n : X \oplus (I_n - X) \succ 0\}$. By Example 3 in [565], it becomes evident that (2.45) satisfies

$$\mathbf{T} + \mathsf{Sy}\left(\mathbf{\Omega}^\top \mathbf{\Gamma}\right) \preceq \Delta\left(\mathbf{\Omega}, \widetilde{\mathbf{\Omega}}, \mathbf{\Gamma}, \widetilde{\mathbf{\Gamma}}\right), \quad \mathbf{T} + \mathsf{Sy}\left(\mathbf{\Omega}^\top \mathbf{\Gamma}\right) = \Delta(\mathbf{\Omega}, \mathbf{\Omega}, \mathbf{\Gamma}, \mathbf{\Gamma}) \tag{2.46a}$$

for all $\mathbf{\Omega}; \widetilde{\mathbf{\Omega}} \in \mathbb{R}^{n \times \mu}$ and for all $\mathbf{\Gamma}; \widetilde{\mathbf{\Gamma}} \in \mathbb{R}^{n \times \mu}$ with $\widetilde{\mathbf{T}} \in \mathbb{S}^{\mu}$ and $Z \oplus (I_n - Z) \succ 0$, which indicates that $\Delta\left(\bullet, \widetilde{\mathbf{\Omega}}, \bullet, \widetilde{\mathbf{\Gamma}}\right)$ in (2.45) is a PSD-overestimate of $\acute{\Delta}(\mathbf{\Omega}, \mathbf{\Gamma}) = \mathbf{T} + \mathsf{Sy}\left[\mathbf{\Omega}^\top \mathbf{\Gamma}\right]$ with respect to the parameterization

$$\begin{bmatrix} \mathbf{vec}(\widetilde{\mathbf{\Omega}}) \\ \mathbf{vec}(\widetilde{\mathbf{\Gamma}}) \end{bmatrix} = \begin{bmatrix} \mathbf{vec}(\mathbf{\Omega}) \\ \mathbf{vec}(\mathbf{\Gamma}) \end{bmatrix}. \tag{2.46b}$$

Now let $\mu = 2n + dn + q + m$ and $Z \oplus (I_n - Z) \succ 0$ and

$$\begin{aligned} &\mathbf{T} = \widehat{\mathbf{\Phi}}, \quad \mathbf{\Omega} = \mathbf{P} = \begin{bmatrix} P & \mathsf{O}_n & Q & \mathsf{O}_{n,q} & \mathsf{O}_{n,m} \end{bmatrix}, \quad \mathbf{\Lambda} = \begin{bmatrix} P & Q \end{bmatrix} \\ &\widetilde{\mathbf{\Omega}} = \widetilde{\mathbf{P}} = \begin{bmatrix} \widetilde{P} & \mathsf{O}_n & \widetilde{Q} & \mathsf{O}_{n,q} & \mathsf{O}_{n,m} \end{bmatrix}, \ \widetilde{\mathbf{\Lambda}} = \begin{bmatrix} \widetilde{P} & \widetilde{Q} \end{bmatrix}, \ \widetilde{P} \in \mathbb{S}^n, \ \widetilde{Q} \in \mathbb{R}^{n \times dn} \\ &\mathbf{\Gamma} = \mathbf{B}\mathbf{K} = \mathbf{B}\left[(I_{2+d} \otimes K) \oplus \mathsf{O}_{p+m}\right], \quad \widetilde{\mathbf{\Gamma}} = \mathbf{B}\widetilde{\mathbf{K}} = \mathbf{B}\left[\left(I_{2+d} \otimes \widetilde{K}\right) \oplus \mathsf{O}_{p+m}\right], \quad \widetilde{K} \in \mathbb{R}^{p \times n} \end{aligned} \tag{2.47}$$

with $\widehat{\mathbf{\Phi}}$, $\mathbf{B}$, $\mathbf{\Lambda}$ and $K$ in (2.44). By (2.46a) with (2.47), we have

$$\begin{aligned} \mathfrak{D}(\mathbf{\Lambda}, K) = \widehat{\mathbf{\Phi}} + \mathsf{Sy}\left[\mathbf{P}^\top \mathbf{B}\mathbf{K}\right] \preceq \mathfrak{J}\left(\mathbf{\Lambda}, \widetilde{\mathbf{\Lambda}}, K, \widetilde{K}\right) := \widehat{\mathbf{\Phi}} + \mathsf{Sy}\left(\widetilde{\mathbf{P}}^\top \mathbf{B}\mathbf{K} + \mathbf{P}^\top \mathbf{B}\widetilde{\mathbf{K}} - \widetilde{\mathbf{P}}^\top \mathbf{B}\widetilde{\mathbf{K}}\right) \\ + \begin{bmatrix} \mathbf{P}^\top - \widetilde{\mathbf{P}}^\top & \mathbf{K}^\top \mathbf{B}^\top - \widetilde{\mathbf{K}}^\top \mathbf{B}^\top \end{bmatrix} \left[Z \oplus (I_n - Z)\right]^{-1} \begin{bmatrix} \mathbf{P} - \widetilde{\mathbf{P}} \\ \mathbf{B}\mathbf{K} - \mathbf{B}\widetilde{\mathbf{K}} \end{bmatrix} \end{aligned} \tag{2.48a}$$

where $\mathfrak{J}\left(\bullet, \widetilde{\mathbf{\Lambda}}, \bullet, \widetilde{K}\right)$ is a PSD-convex overestimate of $\mathfrak{D}(\mathbf{\Lambda}, K)$ in (2.44) with respect to the parameterization

$$\begin{bmatrix} \mathbf{vec}(\widetilde{\mathbf{\Lambda}}) \\ \mathbf{vec}(\widetilde{K}) \end{bmatrix} = \begin{bmatrix} \mathbf{vec}(\mathbf{\Lambda}) \\ \mathbf{vec}(K) \end{bmatrix}. \tag{2.48b}$$

Now it is evident that one can infer $\mathfrak{D}(\mathbf{\Lambda}, K)$ in (2.44) from $\mathfrak{J}\left(\mathbf{\Lambda}, \widetilde{\mathbf{\Lambda}}, K, \widetilde{K}\right) \prec 0$. Moreover, applying the Schur complement to the inequality $\mathfrak{J}\left(\mathbf{\Lambda}, \widetilde{\mathbf{\Lambda}}, K, \widetilde{K}\right) \prec 0$ concludes that $\mathfrak{J}\left(\mathbf{\Lambda}, \widetilde{\mathbf{\Lambda}}, K, \widetilde{K}\right) \prec 0$ with $Z \in \{X \in \mathbb{S}^n : X \oplus (I_n - X) \succ 0\}$ if and only if

$$\begin{bmatrix} \widehat{\mathbf{\Phi}} + \mathsf{Sy}\left(\widetilde{\mathbf{P}}^\top \mathbf{B}\mathbf{K} + \mathbf{P}^\top \mathbf{B}\widetilde{\mathbf{K}} - \widetilde{\mathbf{P}}^\top \mathbf{B}\widetilde{\mathbf{K}}\right) & \mathbf{P}^\top - \widetilde{\mathbf{P}}^\top & \mathbf{K}^\top \mathbf{B}^\top - \widetilde{\mathbf{K}}^\top \mathbf{B}^\top \\ * & -Z & \mathsf{O}_n \\ * & * & Z - I_n \end{bmatrix} \prec 0 \tag{2.49}$$

which now can be handled by standard interior algorithms of SDP provided that the values of $\widetilde{\mathbf{P}}$ and $\widetilde{\mathbf{K}}$ are given. To apply the methods in [565], one has to determine an initial value for $\widetilde{\mathbf{P}}$ and $\widetilde{\mathbf{K}}$ which must be included by the corresponding elements in the relative interior of the feasible set of (2.21a)–(2.21b) in Theorem 2.1. In other words, one may use $\widetilde{P} \leftarrow P$, $\widetilde{Q} \leftarrow Q$ and $\widetilde{K} \leftarrow K$ as the initial data for (2.49) where $P$, $Q$ and $K$ is a feasible solution to Theorem 2.1.

By compiling all the aforementioned procedures according to the expositions in [565], Algorithm 1 can be formulated as follows where $\mathbf{x}$ comprises all the decision variables of $R \in \mathbb{S}^{dn}$, $S; U \in \mathbb{S}^n$ in Theorem 2.1 and $Z \in \mathbb{S}^n$ in (2.49), while $\mathbf{\Lambda}$, $\widetilde{\mathbf{\Lambda}}$, $K$, $\widetilde{K}$ in Algorithm 1 are in line with (2.47). Furthermore, $\rho_1, \rho_2$ and $\varepsilon$ are given constants for achieving regularizations and designating error tolerance, respectively.

---

**Algorithm 1:** An inner convex approximation solution to Theorem 2.1

---

**begin**

Find values for $\mathbf{\Lambda} = \begin{bmatrix} P & Q \end{bmatrix}$ and $K$ included by the corresponding elements in the relative interior of the feasible set of Theorem 2.1. This may be achieved by the feasible solution to Theorem 2.2.

**update** $\widetilde{\mathbf{\Lambda}} \longleftarrow \mathbf{\Lambda}$, $\quad \widetilde{K} \longleftarrow K$,

**solve** $\min\limits_{\mathbf{x}, \mathbf{\Lambda}, K} \mathsf{tr}\left[\rho_1 [*](\mathbf{\Lambda} - \widetilde{\mathbf{\Lambda}}) + \rho_2 [*](K - \widetilde{K})\right]$ subject to (2.21a) and (2.49) to obtain the values of $\mathbf{\Lambda}$ and $K$

**while** $\dfrac{\left\| \begin{bmatrix} \mathbf{vec}(\mathbf{\Lambda}) \\ \mathbf{vec}\,(K) \end{bmatrix} - \begin{bmatrix} \mathbf{vec}(\widetilde{\mathbf{\Lambda}}) \\ \mathbf{vec}(\widetilde{K}) \end{bmatrix} \right\|_\infty}{\left\| \begin{bmatrix} \mathbf{vec}(\widetilde{\mathbf{\Lambda}}) \\ \mathbf{vec}(\widetilde{K}) \end{bmatrix} \right\|_\infty + 1} \geq \varepsilon$ **do**

  **update** $\widetilde{\mathbf{\Lambda}} \longleftarrow \mathbf{\Lambda}$, $\quad \widetilde{K} \longleftarrow K$;

  **solve** $\min\limits_{\mathbf{x}, \mathbf{\Lambda}, K} \mathsf{tr}\left[\rho_1 [*](\mathbf{\Lambda} - \widetilde{\mathbf{\Lambda}}) + \rho_2 [*](K - \widetilde{K})\right]$ subject to (2.21a) and (2.49) to obtain $\mathbf{\Lambda}$ and $K$;

**end**

**end**

---

**Remark 2.7.** When the synthesis condition in Theorem 2.1 contains a convex objective function, for instance $\mathcal{L}^2$ gain $\gamma$ minimization, a termination condition might be added to Algorithm 1 concerning the values of objective functions between two successive iterations [565]. Nonetheless, this condition has not been applied in our numerical examples in this chapter. Moreover, note that a termination condition in terms of the number of the iterations in the while loop can be added in Algorithm 1. Finally, note that the regularization term $\mathsf{tr}\left[\rho_1 [*](\mathbf{\Lambda} - \widetilde{\mathbf{\Lambda}}) + \rho_2 [*](K - \widetilde{K})\right]$ in Algorithm 1 is a special case[4] of the general regularization term

$$\frac{1}{2}\begin{bmatrix}\mathbf{x}-\widetilde{\mathbf{x}}\\ \mathbf{vec}\left(\mathbf{\Lambda}-\widetilde{\mathbf{\Lambda}}\right)\\ \mathbf{vec}\left(K-\widetilde{K}\right)\end{bmatrix}^{\top} Q_k \begin{bmatrix}\mathbf{x}-\widetilde{\mathbf{x}}\\ \mathbf{vec}\left(\mathbf{\Lambda}-\widetilde{\mathbf{\Lambda}}\right)\\ \mathbf{vec}\left(K-\widetilde{K}\right)\end{bmatrix}, \quad Q_k \succeq 0 \tag{2.50}$$

corresponding to the one proposed in CSDP [565].

**Remark 2.8.** The most challenging step in using Algorithm 1 is its initialization if one only considers the synthesis condition in Theorem 2.1. Notwithstanding, given what has been proposed in Theorem 2.2, one way to acquire the initial values of $\widetilde{P}$, $\widetilde{Q}$ and $\widetilde{K}$ is to find a feasible solution to (2.33a) and (2.33b) with given values of $\{\alpha_i\}_{i=1}^{2+d}$. [5]

## 2.5 Numerical Examples

We experimented two numerical examples to demonstrate the effectiveness of the proposed methods. All semidefinite programs are coded in MATLAB© with the aid of the optimization parser Yalmip [603]. We used SeDuMi [604] and SDPT3 [512, 605, 606] as the numerical solvers of SDP.

### 2.5.1 Stability analysis of an LDDS

The example in this subsection was previously utilized as the first numerical example in [169] where SDP problems are solved via SeDuMi [604].

Consider an LDDS

$$\dot{x}(t) = 0.395x(t) - 5\int_{-r}^{0}\cos(12\tau)x(t+\tau)\mathsf{d}\tau \tag{2.51}$$

with $t_0 \in \mathbb{R}$, which corresponds to $A_1 = 0.395$, $A_2 = 0$ and $\widetilde{A}_3(\tau) = -5\cos(12\tau)$ with $n = 1$ and $p = m = q = 0$ in (2.2) and the remaining parts of the state space matrices in (2.2) corresponding to (2.51) are empty matrices. Since $0.395 > 0$ and the distributed term contains a trigonometric function, the methodologies in [607] and [608] cannot analyze the stability of (2.51).

In order to apply the methodology in this chapter, $\boldsymbol{f}(\cdot)$ in Assumption 2.1 is chosen as

$$\boldsymbol{f}(\tau) = \begin{bmatrix}1\\ \sin(12\tau)\\ \cos(12\tau)\end{bmatrix} \quad \text{with} \quad M = \begin{bmatrix}0 & 0 & 0\\ 0 & 0 & 12\\ 0 & -12 & 0\end{bmatrix} \tag{2.52}$$

with $A_3 = \begin{bmatrix}0 & 0 & -5\end{bmatrix}$ which satisfies the conditions in Assumption 2.1 with $d = 3, n = 1$. Furthermore, applying the spectrum methods in [202, 203] with $M = 20$ as the discretization index yields Figure 2.1 as a stability diagram. This diagram plots the values of $\operatorname{sign}\left[\Re(\lambda)\right]$, where $\lambda$ represents the rightmost characteristic roots of the system (2.51). Specifically, by testing sufficiently large $r$ using the code in [203], it is found that $[0.104, 0.1578]$, $[0.6276, 0.6814]$, $[1.1512, 1.205]$, $[1.6748, 1.7286]$ and $[2.1984, 2.2522]$ are the stable delay intervals of (2.51).

[4] This can be comprehended by the relation $\mathrm{tr}(A^{\top}B) = \mathbf{vec}(A)^{\top}\mathbf{vec}(B)$, see the vectorization section in `http://www.ee.ic.ac.uk/hp/staff/dmb/matrix/property.html` which is part of [602].

[5] As elucidated in Remark 2.4, one can apply $\alpha_1 = \alpha_2 = 0$ and $\alpha_i = 0, i = 4 \ldots 2+d$ which allows users to only adjust the value of $\alpha_3$ to apply Theorem 2.2.

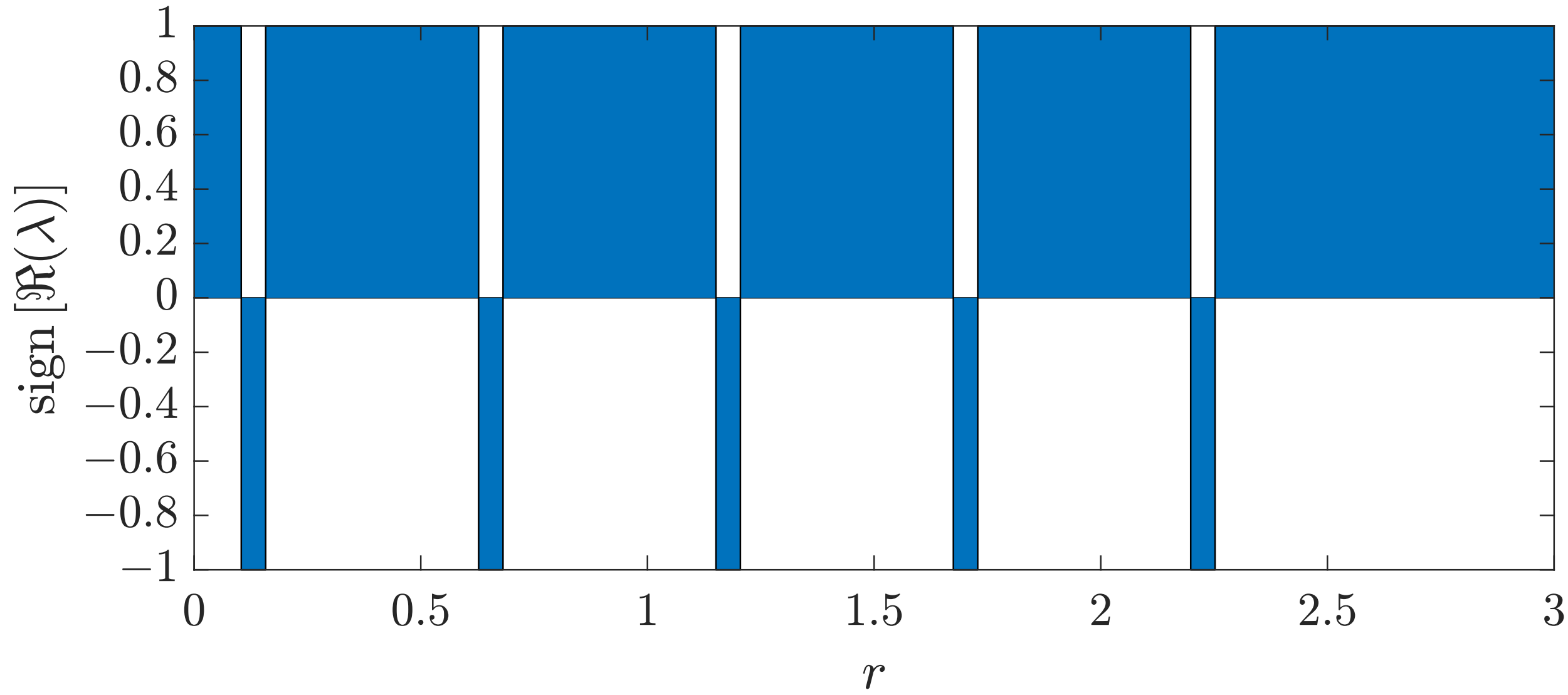

**Figure 2.1:** Diagram showing stability regions of (2.51)

In Table 2.1, we present a comparison between our proposed methodology and the approximation approach in [470] with respect to their ability to detect the boundaries of the stable delay intervals of (2.51). All the semidefinite programs corresponding to the results in Table 2.1 are solved via SeDuMi [604]. Our analysis reveals that $d = 3$ (which is equivalent to the $N$ in [470]) with 12 variables is sufficient to detect all boundaries of the first stability intervals. As $r$ increases, a higher degree of Legendre polynomials is necessary to produce feasible results. It follows that $d = 22$ with 302 variables is required to detect the upper stability boundary 2.2522. In contrast, applying Theorem 2.1 with (2.52) and $A_3 = \begin{bmatrix} 0 & 0 & -5 \end{bmatrix}$ to (2.51), we are able to detect all boundaries of the stable intervals using only 12 decision variables.

| Methodology | first interval | second interval | third interval | forth interval | fifth interval | NDVs |
|---|---|---|---|---|---|---|
| [608] | – | – | – | – | – | – |
| [463] | – | – | – | – | – | – |
| [470], $d = 2$ | $[0.104, 0.1578]$ | – | – | – | – | 12 |
| [470], $d = 8$ | $[0.104, 0.1578]$ | $[0.6276, 0.6814]$ | – | – | – | 57 |
| [470], $d = 13$ | $[0.104, 0.1578]$ | $[0.6276, 0.6814]$ | $[1.1512, 1.205]$ | – | – | 122 |
| [470], $d = 17$ | $[0.104, 0.1578]$ | $[0.6276, 0.6814]$ | $[1.1512, 1.205]$ | $[1.6748, 1.7286]$ | – | 192 |
| [470], $d = 22$ | $[0.104, 0.1578]$ | $[0.6276, 0.6814]$ | $[1.1512, 1.205]$ | $[1.6748, 1.7286]$ | $[2.1984, 2.2522]$ | 302 |
| Theorem 2.1, $d = 2$ | $[0.104, 0.1578]$ | $[0.6276, 0.6814]$ | $[1.1512, 1.205]$ | $[1.6748, 1.7286]$ | $[2.1984, 2.2522]$ | 12 |

**Table 2.1:** Feasible Stability Testing Intervals (NDVs stands for the number of decision variables).

**Remark 2.9.** The boundaries of the stable intervals of (2.51) can be accurately detected using both the approaches in [470] and Theorem 2.1. This demonstrates that both methods are consistent with the reliable calculations in [202]. However, a significant contribution of our methodology is its ability to require fewer variables for distributed kernels exhibiting intensive oscillation patterns, as is the case with (2.51).

**Remark 2.10.** For practical systems, the proposed methods in this chapter can be applied to the models of hematopoietic cell maturation in [127] or SIR Epidemic in [71].

### 2.5.2 Dissipative stabilization via a static state feedback controller

The numerical example in this subsection is based on subsection 4.2.1 in the author's chapter [169]. Note that all semidefinite programs in this subsection were solved via SDPT3 [512], whereas the resulting controllers in [169] are solved via SeDuMi [604] instead.

Consider (2.2) with delay value $r = 1$ and state space parameters

$$
\begin{aligned}
&A_1 = \begin{bmatrix} 0 & 0 \\ 0 & 0.1 \end{bmatrix},\ A_2 = \begin{bmatrix} -1 & -1 \\ 0 & 0.9 \end{bmatrix},\ B_1 = \begin{bmatrix} 0 \\ 1 \end{bmatrix},\ B_2 = B_3(\tau) = \mathbf{0}_2,\ D_1 = \begin{bmatrix} 0.1 & -0.11 \\ 0.21 & 0.1 \end{bmatrix},\\
&\widetilde{A}_3(\tau) = \begin{bmatrix} -0.4 - 0.1\mathrm{e}^{\tau}\sin(20\tau) + 0.3\mathrm{e}^{\tau}\cos(20\tau) & 1 + 0.2\mathrm{e}^{\tau}\sin(20\tau) + 0.2\mathrm{e}^{\tau}\cos(20\tau) \\ -1 + 0.01\mathrm{e}^{\tau}\sin(20\tau) - 0.2\mathrm{e}^{\tau}\cos(20\tau) & 0.4 + 0.3\mathrm{e}^{\tau}\sin(20\tau) + 0.4\mathrm{e}^{\tau}\cos(20\tau) \end{bmatrix},\\
&C_1 = \begin{bmatrix} -0.1 & 0.2 \\ 0 & 0.1 \end{bmatrix},\ C_2 = \begin{bmatrix} -0.1 & 0 \\ 0 & 0.2 \end{bmatrix},\ B_4 = \begin{bmatrix} 0.2 \\ 0.3 \end{bmatrix},\ B_5 = B_6(\tau) = \mathbf{0}_2,\\
&\widetilde{C}_3(\tau) = \begin{bmatrix} 0.2\mathrm{e}^{\tau}\sin(20\tau) & 0.1 + 0.1\mathrm{e}^{\tau}\cos(20\tau) \\ 0.1\mathrm{e}^{\tau}\sin(20\tau) - 0.1\mathrm{e}^{\tau}\cos(20\tau) & -0.2 + 0.3\mathrm{e}^{\tau}\sin(20\tau) \end{bmatrix},\ D_2 = \begin{bmatrix} 0.1 & 0.2 \\ 0.12 & 0.1 \end{bmatrix}
\end{aligned}
\tag{2.53}
$$

which correspond to $n = q = m = 2$ and $p = 1$. Since $(A_1, B_1)$ is not stabilizable, the method in [463] cannot be applied here, regardless of whether $\widetilde{A}_3(\tau)$ and $\widetilde{C}_3(\tau)$ can be approximated by rational functions [465]. Moreover, based on the structures of $\widetilde{A}_3(\tau)$ and $\widetilde{C}_3(\tau)$ in (2.53), the corresponding delay system does not possess either forwarding or backstepping structures. Thus, the constructive approaches outlined in [503] may not be applicable. Now by using the spectrum method in [202, 203] with a testing grid vector of different delay values, one can make the estimation[6] that the delay system with (2.53) is unstable for $0 \leq r \leq 10$. Additionally, $\mathscr{L}^2$ attenuation is considered here as the performance objective for the system with parameters (2.53), which corresponds to $J_3 = -J_1 = \gamma I_2$, $\widetilde{J} = I_2$ and $J_2 = \mathsf{O}_2$ for the supply rate function (2.10).

Observing the elements inside of $\widetilde{A}_3(\tau)$, $\widetilde{C}_3(\tau)$, we choose

$$
\boldsymbol{f}(\tau) = \begin{bmatrix} 1 & 10\mathrm{e}^{\tau}\sin(20\tau) & 10\mathrm{e}^{\tau}\cos(20\tau) \end{bmatrix}^{\top}
\tag{2.54}
$$

which gives $M = \begin{bmatrix} 0 & 0 & 0 \\ 0 & 1 & 20 \\ 0 & -20 & 1 \end{bmatrix}$ and

$$
A_3 = 0.1\begin{bmatrix} -4 & 10 & -0.1 & 0.2 & 0.3 & 0.2 \\ -10 & 4 & 0.01 & 0.3 & -0.2 & 0.4 \end{bmatrix},\quad C_3 = 0.1\begin{bmatrix} 0 & 1 & 0.2 & 0 & 0 & 0.1 \\ 0 & -2 & 0.1 & 0.3 & -0.1 & 0 \end{bmatrix}
\tag{2.55a}
$$

in accordance with the relations in Assumption 2.1 with $d = 3$, $n = m = q = 2$.

**Remark 2.11.** For $\widetilde{A}_3(\tau)$, $\widetilde{C}_3(\tau)$ in (2.53), constant matrices $A_3$ and $C_3$ can be constructed via direct observations of the elements in $\widetilde{A}_3(\tau), \widetilde{C}_3(\tau)$, given that the structure of $\boldsymbol{f}(\tau)$ in (2.54) is well ordered. Moreover, applying substitutions $\theta_1 = 10\mathrm{e}^{\tau}\sin(20\tau)$, $\theta_2 = 10\mathrm{e}^{\tau}\cos(20\tau)$ to $\widetilde{A}_3(\tau)$, $\widetilde{C}_3(\tau)$ yields $\widetilde{A}_3(\tau) = \widehat{A}_3(\theta_1, \theta_2) = A_3 M(\theta_1, \theta_2)$, $\widetilde{C}_3(\tau) = \widehat{C}_3(\theta_1, \theta_2) = C_3 M(\theta_1, \theta_2)$, where $A_3, C_3$ can be constructed based on the corresponding coefficients with respect to $\theta_1, \theta_2$. Finally, one can use the formulas $A_3 = \begin{bmatrix} \widehat{A}_3(0,0) & \frac{\partial \widehat{A}_3(\theta_1,\theta_2)}{\partial \theta_1} & \frac{\partial \widehat{A}_3(\theta_1,\theta_2)}{\partial \theta_2} \end{bmatrix}$, $C_3 = \begin{bmatrix} \widehat{C}_3(0,0) & \frac{\partial \widehat{C}_3(\theta_1,\theta_2)}{\partial \theta_1} & \frac{\partial \widehat{C}_3(\theta_1,\theta_2)}{\partial \theta_2} \end{bmatrix}$ to construct $A_3, C_3$. (This approach is suggested by Kwin in `https://www.mathworks.com/matlabcentral/answers/309797-extracting-the-coefficient-of-a-polynomials-matrix#answer_241277`)

[6] It is an estimation since only a finite number of pointwise delay values can be tested by the method in [203].

By applying Theorem 2.2 to (2.5) with the parameters in (2.53)–(2.55a) and $\alpha_1 = 1$, $\{\alpha_i\}_{i=2}^{8} = 0$, it is shown that the system (2.2) with the system parameters in (2.53) is stabilized by

$$\boldsymbol{u}(t) = \begin{bmatrix} 1.78 & -6.3569 \end{bmatrix} \boldsymbol{x}(t) \tag{2.56}$$

with performance $\min \gamma = 0.3468$. It should be noted that SDPT3 [512] is used as the solver for SDP following this step. Due to the simplification applied in step (2.40), the value of $\min \gamma$ calculated by Theorem 2.2 may be more conservative compared to Theorem 2.1 with a given controller gain $K$. In consequence, we apply Theorem 2.1 with $K = \begin{bmatrix} 1.78 & -6.3569 \end{bmatrix}$ to (2.5) with the parameters in (2.53)–(2.55a), which shows that (2.56) is able to achieve $\min \gamma = 0.27078$. To verify (2.56) is a stabilizing controller for the system with (2.53), we once again apply the spectrum method described in [202] to the resulting CLS. This yields $-0.16051 < 0$ as the real part of the rightmost characteristic root pair, thereby proving that the CLS is stable.

To calculate controllers with further performance improvement using Algorithm 1, one can first solve Theorem 2.1 with the controller gain specified in (2.56) to compute a feasible solution for $P$ and $Q$. The controller gain in (2.56), with the aforementioned $P$ and $Q$, can then be used as the initial values for $\widetilde{K}$, $\widetilde{P}$, $\widetilde{Q}$ in Algorithm 1. The results produced by Algorithm 1 with $\rho_1 = \rho_2 = 10^{-8}$ and $\varepsilon = 10^{-12}$ are summarized in Table 2.2 where NoIs represents the number of iterations executed by the while loop. Additionally, the spectral abscissas of the CLSs in Table 2.2 are calculated by the method in [202, 203]. From Table 2.2, it is observed that more iterations can lead to better $\min \gamma$ values at the expense of increased computational burden.

| Controller gain $K$ | $\begin{bmatrix} 2.4692 \\ -9.3734 \end{bmatrix}^\top$ | $\begin{bmatrix} 2.9890 \\ -11.6754 \end{bmatrix}^\top$ | $\begin{bmatrix} 3.3877 \\ -13.4389 \end{bmatrix}^\top$ | $\begin{bmatrix} 3.7316 \\ -14.9638 \end{bmatrix}^\top$ |
|---|---|---|---|---|
| $\min \gamma$ | 0.27041 | 0.27031 | 0.27027 | 0.27025 |
| Number of Iterations | 10 | 20 | 30 | 40 |
| Spectral abscissa of the CLSs | $-0.15943$ | $-0.1587$ | $-0.15818$ | $-0.15781$ |

**Table 2.2:** $\min \gamma$ produced by different iterations

**Remark 2.12.** Some functions in the $\boldsymbol{f}(\cdot)$ in (2.55a) have been scaled compared to the form in (2.52). This is because a $\mathsf{F}$ with a large condition number may affect the numerical solvability of the corresponding SDP constraints.

Now consider the CLS stabilized by the controller with the gain $K = \begin{bmatrix} 3.7316 & -14.9638 \end{bmatrix}$ in Table 2.2, assume $\boldsymbol{z}(t) = 0, t < 0$ with initial condition $\forall \tau \in [-1, 0]$, $\boldsymbol{\phi}(\tau) = \begin{bmatrix} 5 & 3 \end{bmatrix}^\top$ and disturbance $\boldsymbol{w}(t) = \begin{bmatrix} 10 \sin 20\pi t(\mathit{1}(t) - \mathit{1}(t-5)) & 5 \sin 20\pi t(\mathit{1}(t) - \mathit{1}(t-5)) \end{bmatrix}$, where $\mathit{1}(t)$ is the Heaviside step function. Then the corresponding state trajectories and outputs of the CLS are presented in Figure 2.2–2.3, where all the figures are plotted over $[-1, 50]$ which takes into account the initial condition $\boldsymbol{\phi}(\tau) = \begin{bmatrix} 5 & 3 \end{bmatrix}^\top$, $\tau \in [-1, 0]$, and the fact that $\boldsymbol{z}(t) = 0, t < 0$.

The results in Figure 2.2–2.3 clearly show the effectiveness of $K = \begin{bmatrix} 3.7316 & -14.9638 \end{bmatrix}$ in Table 2.2 with respect to stabilization and disturbance attenuation ($\mathcal{L}^2$ gain control)

**Remark 2.13.** The plots in Figure 2.2–2.3 are obtained from Simulink©, where the DDs were implemented numerically by discretizing the integrals using the trapezoidal rule (200 sampling points).

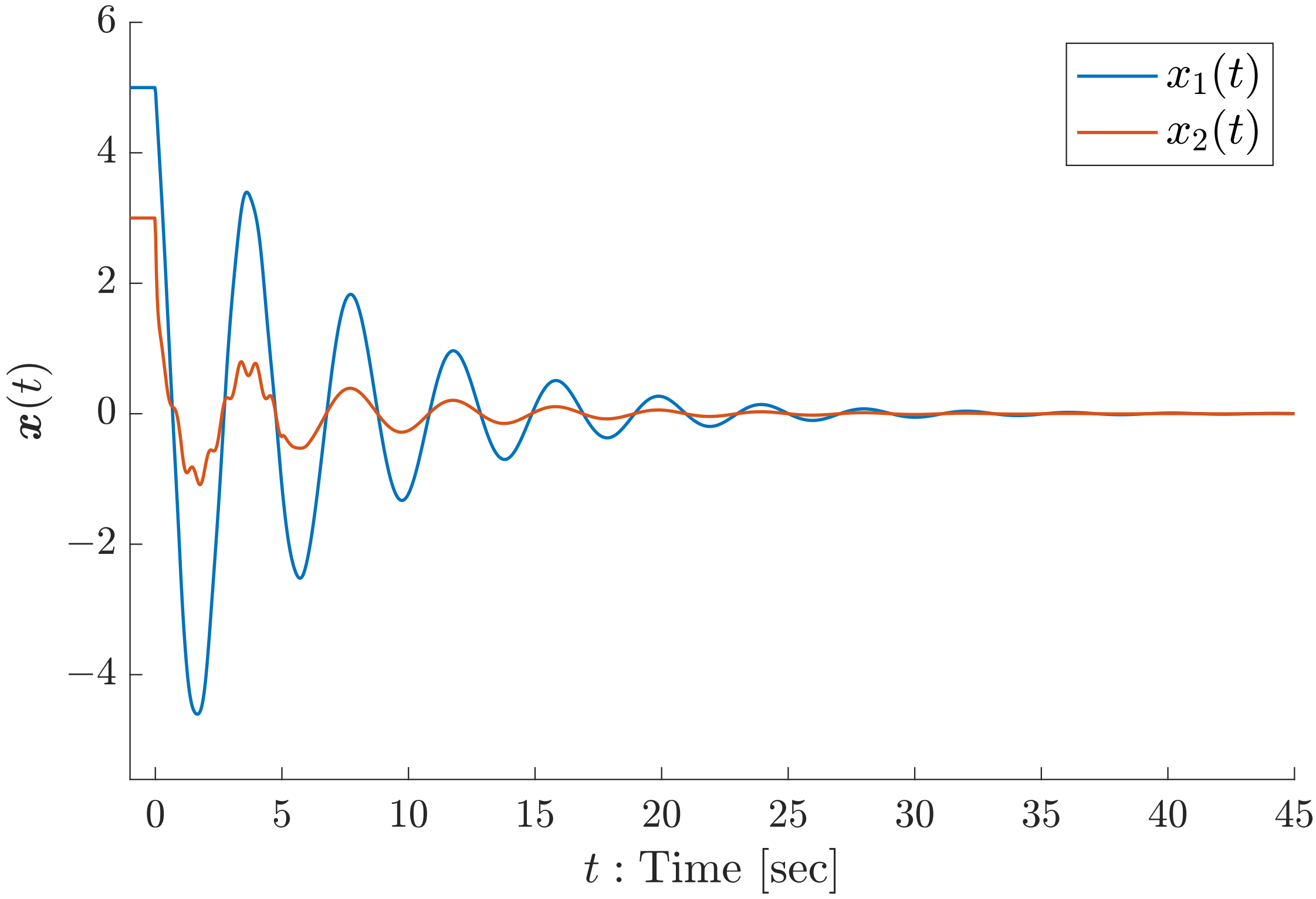

**Figure 2.2:** The close-loop system's trajectory with $K = \begin{bmatrix} 3.7316 & -14.9638 \end{bmatrix}$ in Table 7.2

Since no numerical solvers for delay systems are available in Simulink©, an ODE solver[7] was used for computations. However, due to the unpredictability of the potential issues that may arise from using an ODE solver for delay systems, the results in Figure 2.2–2.3 should only be considered as an estimation of the actual behavior of the system's trajectories and output.

[7]For our simulation, we used ode8 with fundamental sampling time 0.001.

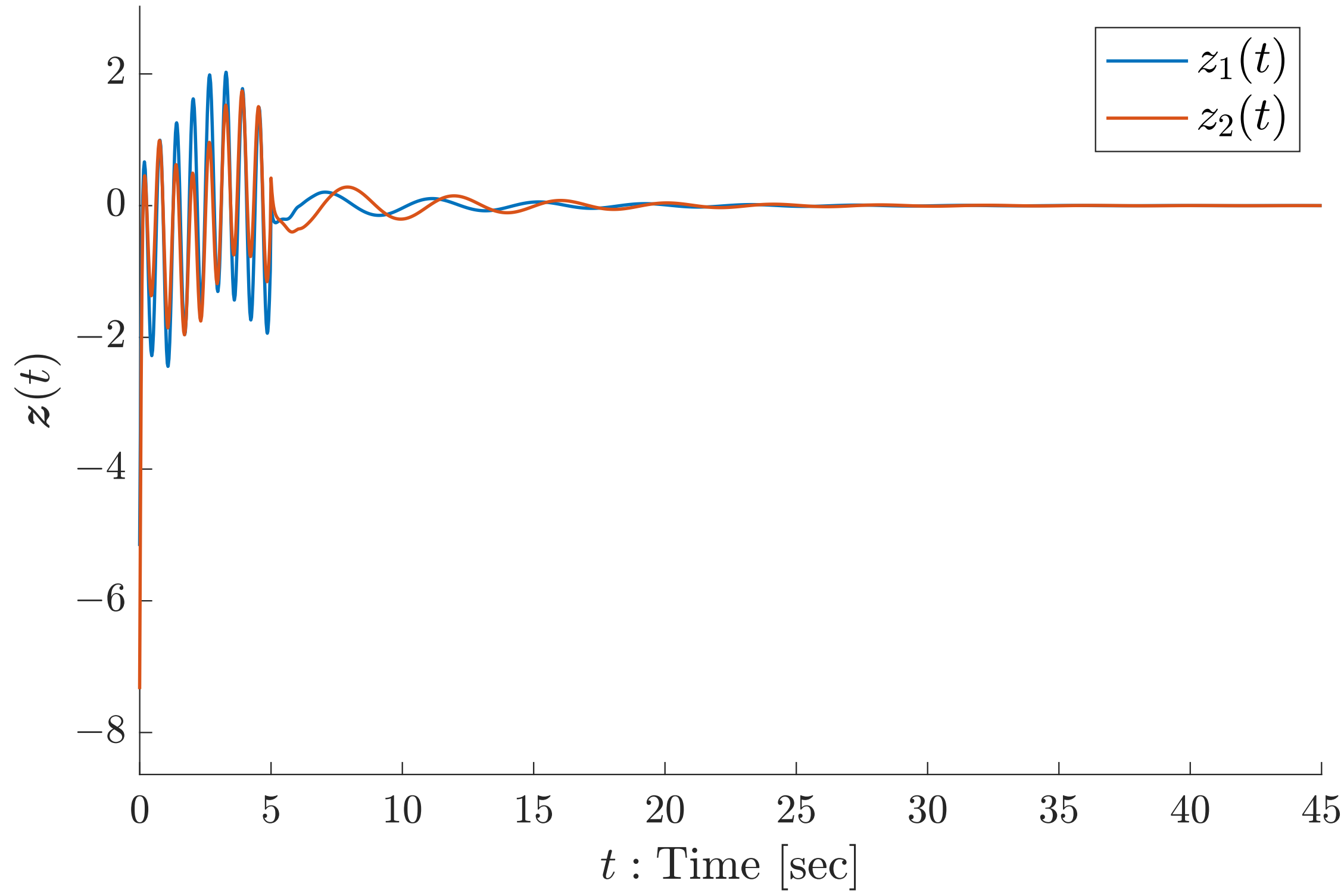

**Figure 2.3:** The output of the CLS with $K = \begin{bmatrix} 3.7316 & -14.9638 \end{bmatrix}$ in Table 7.2

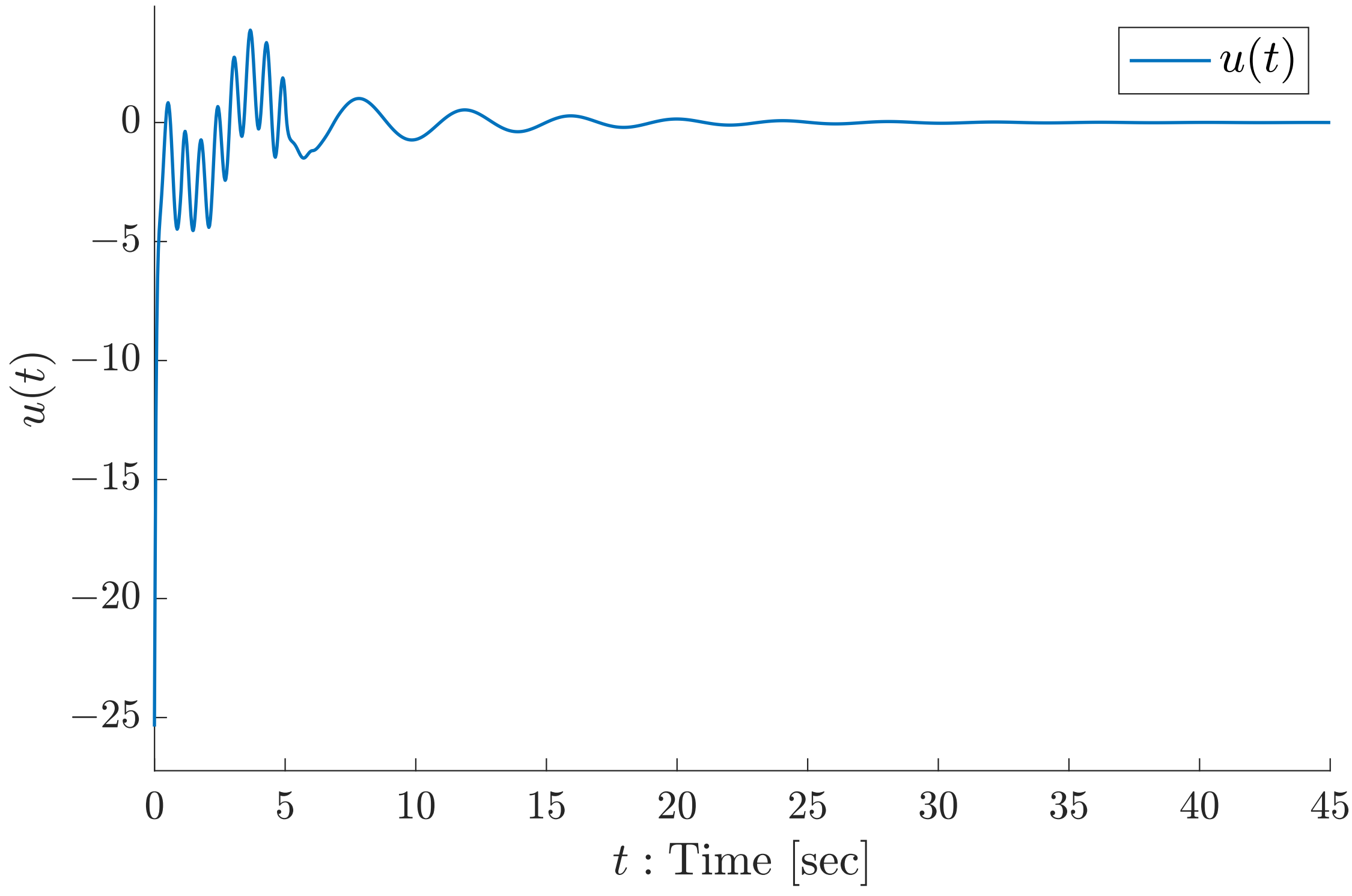

**Figure 2.4:** Control action of the CLS with $K = \begin{bmatrix} 3.7316 & -14.9638 \end{bmatrix}$ in Table 7.2

# Chapter 3

# Dissipative Stabilization of Uncertain Linear Distributed Delay Systems

## 3.1 Introduction

To characterize modeling errors and the impact of operating environment, it is preferable to incorporate uncertainties into the systems mathematical model. Addressing uncertainties in dynamical systems has been extensively researched over the past several decades [609–613], and these methodologies may also be applied to problems pertaining to linear parameter varying systems when Lyapunov approaches are considered [452, 516, 614, 615].

The mathematical description of a systems robustness is directly affected by the complexity of its uncertainties. Therefore, incorporating uncertainties with more general structures into a system can lead to more general results in terms of the systems robustness. One common structure for uncertainties is norm-bounded uncertainty [610, 616], represented as $G\Delta H$ with $\Delta^\top\Delta \preceq I$, [1] where $G$ and $H$ are given and $\Delta$ can be a function of time $t$ or a function of other variables. To handle norm-bounded uncertainties in the context of using LMIs, the Petersen Lemma [610] (see the summary in [452, Lemma C.10.1] also) was introduced to provide tractable conditions for an LMI term with norm-bounded uncertainties. This idea was further extended in [617] to handle uncertainties of linear fractional form. Additionally, the handling of norm-bounded uncertainties subject to a full block scaling constraint is expounded in [452, Lemma C.10.4] with its proof supported by the linearization procedure described in [618]. Nevertheless, it is certainly more beneficial to consider linear fractional uncertainties subject to full block scaling constraints.

In Chapter 3, we investigated the problem of stabilizing an uncertain linear DDS with DDs in its states, inputs, and outputs, where the structure of the DD kernels follows the same class considered in Chapter 2. We propose methodologies for calculating the gains of a static state feedback controller for an LDDS with uncertainties in linear fractional form. Specifically, the uncertainties of state-space matrices are represented as $G(I-\Delta F)^{-1}\Delta H$, where[2] $\Delta$ subject to the full block scaling constraint

$$\Delta \in \left\{ \widehat{\Delta} \;\middle|\; \begin{bmatrix} I \\ \widehat{\Delta} \end{bmatrix}^\top \begin{bmatrix} \Xi^{-1} & \Lambda \\ * & \Gamma \end{bmatrix} \begin{bmatrix} I \\ \widehat{\Delta} \end{bmatrix} \succeq 0 \right\}, \quad \Xi^{-1} \succ 0, \quad \Gamma \preceq 0.$$

Additionally, the coefficients of the DD kernels are also assumed to be affected by uncertainties, mak-ing the model of the uncertain system considered in this chapter sufficiently general. Furthermore, the proposed scenarios for static controller design can be modified to enable the design of a resilient dynamic state feedback controller for an uncertain LTI system with input delays. Interestingly, compared to finding a robust static state feedback controller, computing a resilient dynamic state feedback controller may be easier due to the availability of results pertaining to the construction of predictor controllers for input delay systems.

[1] Note that the constraints on $\Delta$ can be chosen as $\Delta^\top\Delta \preceq R$ where $R \succ 0$.

[2] Note that each state space matrix of the system in this chapter has its own uncertainty parameter (constraint), the term $G(I-\Delta F)^{-1}\Delta H$ and its constraint only provide a common characteristic of the uncertainties considered in this chapter

The rest of the chapter is organized as follows. Section 3.2 presents the plant of an open-loop uncertain system, where the uncertainties are denoted in linear fractional form subject to full block scaling constraints. Since the synthesis solutions for an LDDS have already been presented in Chapter 2, these solutions can be adapted to handle the uncertain extension as long as the uncertainties can be addressed by a mathematical solution. To handle uncertainties in the form of $G(I-\Delta F)^{-1}\Delta H$ in the context of SDP constraints, Lemma 3.1 is introduced and later applied to derive optimization constraints for robust synthesis. Utilizing this lemma, the main results on the synthesis of robust static state feedback controllers are set forth in Theorem 3.1–3.2 in Section 3.3 expressed in terms of matrix inequalities of finite dimensions. Additionally, an iterative algorithm for solving the BMI in Theorem 3.1 is derived based on the approach used in Chapter 2 for the derivation of Algorithm 1. Next, the whole Section 3.4 is dedicated to the study of designing a resilient dynamic state feedback controller for an uncertain LTI system with input delays, where both the models of both the plant and controller incorporate uncertainties terms. Following the same strategy proposed in the previous section, the corresponding synthesis condition for the existence of a resilient dynamic state feedback controller is presented in Theorem 3.3, which can be solved by the iterative algorithm enunciated in Algorithm 3. Interestingly, this iterative algorithm for designing a robust dynamic state feedback controller can be initiated directly using the gains of the constructed predictor controllers. Finally, the testing results of several numerical examples are presented in Section 3.5 to demonstrate the effectiveness and capabilities of our proposed methodologies.

## 3.2 Problem Formulation

Consider the following uncertain LDDS

$$\begin{aligned}
\dot{\boldsymbol{x}}(t) &= \grave{A}_1\boldsymbol{x}(t)+\grave{A}_2\boldsymbol{x}(t-r)+\int_{-r}^{0}\grave{A}_3F(\tau)\boldsymbol{x}(t+\tau)\mathsf{d}\tau+\grave{B}_1\boldsymbol{u}(t)+\grave{B}_2\boldsymbol{u}(t-r)\\
&\qquad +\int_{-r}^{0}\grave{B}_3F(\tau)\boldsymbol{u}(t+\tau)\mathsf{d}\tau+\grave{D}_1\boldsymbol{w}(t),\quad \forall t\geq t_0\\
\boldsymbol{z}(t) &= \grave{C}_1\boldsymbol{x}(t)+\grave{C}_2\boldsymbol{x}(t-r)+\int_{-r}^{0}\grave{C}_3(\tau)\boldsymbol{x}(t+\tau)\mathsf{d}\tau+\grave{B}_4\boldsymbol{u}(t)+\grave{B}_5\boldsymbol{u}(t-r)\\
&\qquad +\int_{-r}^{0}\grave{B}_6F(\tau)\boldsymbol{u}(t+\tau)\mathsf{d}\tau+\grave{D}_2\boldsymbol{w}(t)\\
&\forall\theta\in[-r,0],\quad \boldsymbol{x}(t_0+\theta)=\boldsymbol{\phi}(\theta)
\end{aligned}\tag{3.1}$$

to be stabilized, where $t_0\in\mathbb{R}$, $\boldsymbol{x}(t)\in\mathbb{R}^n$ satisfies (3.1), $\boldsymbol{u}(t)\in\mathbb{R}^p$ denotes input signals, $\boldsymbol{w}(\cdot)\in\mathcal{L}^2([t_0,\infty);\mathbb{R}^q)$ represents the disturbance, $\boldsymbol{z}(t)\in\mathbb{R}^m$ is the regulated output, and $\boldsymbol{\phi}(\cdot)\in\mathcal{C}([-r,0];\mathbb{R}^n)$ denotes initial condition. Moreover, the DD kernel matrix is $\mathbb{R}^{n\times dn}\ni F(\tau)=\boldsymbol{f}(\tau)\otimes I_n$ with

$\boldsymbol{f}(\cdot) \in C^1\left([-r,0];\mathbb{R}^d\right)$ which satisfies

$$
\begin{gathered}
\int_{-r}^{0} \boldsymbol{f}(\tau)\boldsymbol{f}^{\top}(\tau)\mathsf{d}\tau \succ 0, \\
\exists M \in \mathbb{R}^{d\times d} : \frac{\mathsf{d}\boldsymbol{f}(\tau)}{\mathsf{d}\tau} = M\boldsymbol{f}(\tau).
\end{gathered}
\tag{3.2}
$$

The state space parameters in (3.1) are defined as

$$
\begin{aligned}
&\begin{bmatrix} \grave{A}_1 & \grave{B}_1 & \grave{A}_2 & \grave{B}_2 & \grave{A}_3 & \grave{B}_3 & \grave{D}_1 \end{bmatrix} \\
&\qquad = \begin{bmatrix} A_1 & B_1 & A_2 & B_2 & A_3 & B_3 & D_1 \end{bmatrix} + \mathop{\mathsf{Row}}_{i=1}^{7} \left( G_i (I - \Delta_i F_i)^{-1} \Delta H_i \right)
\end{aligned}
\tag{3.3a}
$$

and

$$
\begin{aligned}
&\begin{bmatrix} \grave{C}_1 & \grave{B}_4 & \grave{C}_2 & \grave{B}_5 & \grave{C}_3 & \grave{B}_6 & \grave{D}_2 \end{bmatrix} \\
&\qquad = \begin{bmatrix} C_1 & B_4 & C_2 & B_5 & C_3 & B_6 & D_2 \end{bmatrix} + \mathop{\mathsf{Row}}_{i=8}^{14} \left[ G_i (I - \Delta_i F_i)^{-1} \Delta H_i \right]
\end{aligned}
\tag{3.3b}
$$

where the dimensions of the parameters $A_i$; $C_i$, $i = 1,2,3$ and $B_j$, $j = 1,\ldots,6$ and $D_1$, $D_2$ are identical to the parameters in (3.1) and Assumption 2.1 in the previous chapter. The given matrices $G_i$, $F_i$, $H_i$ with $i = 1,\ldots,14$ determine the configuration of the uncertainty $\Delta_i$ which subject to the full block constraints [452, 515]

$$
\Delta_i \in \left\{ \widehat{\Delta}_i \;\middle|\; \begin{bmatrix} I \\ \widehat{\Delta}_i \end{bmatrix}^{\top} \begin{bmatrix} \Xi_i^{-1} & \Lambda_i \\ * & \Gamma_i \end{bmatrix} \begin{bmatrix} I \\ \widehat{\Delta}_i \end{bmatrix} \succeq 0 \right\}, \quad \forall i = 1,\cdots,14, \quad \Xi_i^{-1} \succ 0, \quad \Gamma_i \preceq 0
\tag{3.3c}
$$

where $\Xi_i$, $\Lambda_i$ and $\Gamma_i$ are given. Finally, all matrices in (3.3a)–(3.3c) are supposed to have compatible dimensions.

**Remark 3.1.** We have omitted the decomposition procedure for the DDs as in Assumption 2.1 of Chapter 1. It is a well-established fact that any DD containing functions as entries of a matrix exponential $\mathsf{e}^{X\tau}$, $X \in \mathbb{R}^{d\times d}$ can be expressed using the forms of DDs specified in (3.3a).

The constraints in (3.3c) can be rewritten as $\widehat{\Delta}_i^{\top}\Gamma_i\widehat{\Delta}_i + \mathsf{Sy}(\Lambda_i\widehat{\Delta}_i) + \Xi_i^{-1} \succeq 0, \forall i = 1,\cdots,14$, which is equivalent to a single inequality $\bigoplus_{i=1}^{14}\left(\widehat{\Delta}_i^{\top}\Gamma_i\widehat{\Delta}_i + \mathsf{Sy}(\Lambda_i\widehat{\Delta}_i) + \Xi_i^{-1}\right) \succeq 0$. Moreover, using the property of diagonal matrices $(X+Y)\oplus(X+Y) = (X\oplus X)+(Y\oplus Y)$ and $XY\oplus XY = (X\oplus X)(Y\oplus Y)$ shows that (3.3c) is equivalent to

$$
\begin{gathered}
\bigoplus_{i=1}^{14}\Delta_i \in \mathcal{U} := \left\{ \bigoplus_{i=1}^{14}\widehat{\Delta}_i \;\middle|\; \begin{bmatrix} I \\ \bigoplus_{i=1}^{14}\widehat{\Delta}_i \end{bmatrix}^{\top} \begin{bmatrix} \bigoplus_{i=1}^{14}\Xi_i^{-1} & \bigoplus_{i=1}^{14}\Lambda_i \\ * & \bigoplus_{i=1}^{14}\Gamma_i \end{bmatrix} \begin{bmatrix} I \\ \bigoplus_{i=1}^{14}\widehat{\Delta}_i \end{bmatrix} \succeq 0 \right\} \\
\bigoplus_{i=1}^{14}\Xi_i^{-1} \succ 0, \quad \bigoplus_{i=1}^{14}\Gamma_i \preceq 0.
\end{gathered}
\tag{3.4}
$$

Having demonstrated the equivalence relation between (3.4) and (3.3c), representation (3.4) is utilized in the next section to derive our robust synthesis conditions due to its compact structure.

**Remark 3.2.** The uncertainties in (3.3a) with constraint (3.4) provide a comprehensive characterization of uncertainties in an LDDS. Note that the robust terms $\grave{A}_3 = A_3 + G_5(I-\Delta_5F_5)^{-1}\Delta H_5$, $\grave{B}_3 = B_3 + G_6(I-\Delta_6F_6)^{-1}\Delta_6H_6$, $\grave{C}_3 = C_3 + G_{12}(I-\Delta_{12}F_{12})^{-1}\Delta_{12}H_{12}$ and $\grave{B}_6 = B_6 + G_{13}(I -$ $\Delta_{13}F_{13})^{-1}\Delta H_{13}$ lead to the distributed terms $\widetilde{A}_3(\tau) = A_3F(\tau) + G_5(I - \Delta_5F_5)^{-1}\Delta_5H_5F(\tau)$ and $\widetilde{B}_3(\tau) = B_3F(\tau) + G_6(I - \Delta_6F_6)^{-1}\Delta H_6F(\tau)$ and $\widetilde{C}_3(\tau) = C_3F(\tau) + G_{12}(I - \Delta_{12}F_{12})^{-1}\Delta H_{12}F(\tau)$ and $\widetilde{B}_6(\tau) = B_6F(\tau) + G_{13}(I - \Delta_{13}F_{13})^{-1}\Delta H_{13}F(\tau)$, respectively. This further demonstrates that the uncertainties associated with the distributed terms in (3.3a), as all coefficients of the functions in $\widetilde{A}_3(\tau), \widetilde{B}_3(\tau), \widetilde{C}_3(\tau)$ and $\widetilde{B}_6(\tau)$ are subject to variations of $G_5(I-\Delta_5F_5)^{-1}\Delta H_5$, $G_6(I-\Delta_6F_6)^{-1}\Delta H_6$, $G_{12}(I - \Delta_{12}F_{12})^{-1}\Delta H_{12}$ and $G_{13}(I - \Delta_{13}F_{13})^{-1}\Delta H_{12}$, respectively.

Now substitute $\boldsymbol{u}(t) = K\boldsymbol{x}(t)$ into (3.1), one can derive the expression of the closed-loop uncertain system

$$\begin{bmatrix}\dot{\boldsymbol{x}}(t)\\ \boldsymbol{z}(t)\end{bmatrix} = \left(\begin{bmatrix}\mathbf{A} + \mathbf{B}_1\left[(I_{2+d}\otimes K)\oplus O_q\right]\\ \mathbf{C} + \mathbf{B}_2\left[(I_{2+d}\otimes K)\oplus O_q\right]\end{bmatrix} + \begin{bmatrix}\mathsf{Row}_{i=1}^{7}\, G_i & O\\ O & \mathsf{Row}_{i=8}^{14}\, G_i\end{bmatrix}\left(I - \bigoplus_{i=1}^{14}\Delta_iF_i\right)^{-1}\boldsymbol{\Delta}\begin{bmatrix}\mathbf{H}_1\\ \mathbf{H}_2\end{bmatrix}\right)\boldsymbol{\chi}(t) \tag{3.5}$$
$$\forall\theta \in [-r, 0],\quad \boldsymbol{x}(t_0 + \theta) = \boldsymbol{\phi}(\theta)$$

with $\boldsymbol{\chi}(t)$ in (2.6c), where $\boldsymbol{\Delta} := \bigoplus_{i=1}^{14}\Delta_i$ and

$$\begin{aligned}\mathbf{H}_1 &= \left[\begin{bmatrix}H_1\\ H_2K\end{bmatrix}\oplus\begin{bmatrix}H_3\\ H_4K\end{bmatrix}\oplus\begin{bmatrix}H_5\\ H_6\left(I_d\otimes K\right)\end{bmatrix}\oplus H_7\right],\\ \mathbf{H}_2 &= \left[\begin{bmatrix}H_8\\ H_9K\end{bmatrix}\oplus\begin{bmatrix}H_{10}\\ H_{11}K\end{bmatrix}\oplus\begin{bmatrix}H_{12}\\ H_{13}\left(I_d\otimes K\right)\end{bmatrix}\oplus H_{14}\right]\end{aligned} \tag{3.6}$$

In order to handle the uncertainties structure in (3.3a) and (3.4), the following lemma is derived.

### 3.2.1 A Lemma dealing with uncertainties

**Lemma 3.1.** *For any $n; m; p; q \in \mathbb{N}$, $\Theta_1 \in \mathbb{S}^p_{\succ 0}$, $\Theta_3 \in \mathbb{S}^m_{\preceq 0}$, $\Theta_2 \in \mathbb{R}^{p\times m}$, $\Phi \in \mathbb{S}^n$, $G \in \mathbb{R}^{n\times m}$, $H \in \mathbb{R}^{p\times n}$, $F \in \mathbb{R}^{p\times m}$ if*

$$\exists\alpha > 0 : \begin{bmatrix}I_m & -I_m - \alpha F^\top\Theta_2 & \alpha F^\top\\ * & I_m - \alpha\Theta_3 & O_{m,p}\\ * & * & \alpha\Theta_1\end{bmatrix} \succ 0, \tag{3.7}$$

*then*

$$\Phi + \mathsf{Sy}\left[G(I_m - \Delta F)^{-1}\Delta H\right] \prec 0,\ \forall\Delta \in \mathcal{F} \subseteq \mathcal{D} := \left\{\widehat{\Delta} \in \mathbb{R}^{m\times p}\ \middle|\ [*]\begin{bmatrix}\Theta_1^{-1} & \Theta_2\\ * & \Theta_3\end{bmatrix}\begin{bmatrix}I_p\\ \widehat{\Delta}\end{bmatrix} \succeq 0\right\} \tag{3.8}$$

*holds if*

$$\exists\kappa > 0 : \begin{bmatrix}\Phi & G + \kappa H^\top\Theta_2 & \kappa H^\top\\ * & \kappa F^\top\Theta_2 + \kappa\Theta_2^\top F + \kappa\Theta_3 & \kappa F^\top\\ * & * & -\kappa\Theta_1\end{bmatrix} \prec 0. \tag{3.9}$$

*Moreover, for the situation when $\Theta_1^{-1} = O_p$, (3.9) and (3.7) become*

$$\exists\kappa > 0 : \begin{bmatrix}\Phi & G + \kappa H^\top\Theta_2\\ * & \kappa F^\top\Theta_2 + \kappa\Theta_2^\top F + \kappa\Theta_3\end{bmatrix} \prec 0, \tag{3.10}$$

$$\exists\alpha > 0 : \begin{bmatrix}I_m & -I_m - \alpha F^\top\Theta_2\\ * & I_m - \alpha\Theta_3\end{bmatrix} \succ 0, \tag{3.11}$$

*respectively.*

*Proof.* See Appendix A. ■

**Remark 3.3.** The hypothesis $\Theta_1 \succ 0$ in this lemma is motivated by the expectation to derive convex conditions (3.7) and (3.9) by the application of the Schur complement at the steps of (A.15) and (A.7).

Lemma 3.1 can cover a wide range of uncertainty configurations. Specifically, let $\mathbb{S}^p \ni \Theta_1^{-1} = R \succ 0$, $\Theta_2 = \mathsf{O}_{p,m}$, $\Theta_3 = -I_m$, $F = \mathsf{O}_{p,m}$ and $\mathcal{D} = \left\{\hat{\Delta} \in \mathbb{R}^{m\times p} : \hat{\Delta}^\top \hat{\Delta} \preceq R\right\}$ and

$$\exists \kappa > 0 : \begin{bmatrix} \Phi & G & \kappa H^\top \\ * & -\kappa I_m & \mathsf{O}_m \\ * & * & -\kappa R^{-1} \end{bmatrix} \prec 0 \tag{3.12}$$

which corresponds to (3.9). Now it is obvious that (3.12) is equivalent to

$$\exists \kappa > 0 : \Phi - \begin{bmatrix} G & \kappa H^\top \end{bmatrix} \begin{bmatrix} -\kappa I_m & \mathsf{O}_m \\ * & -\kappa R^{-1} \end{bmatrix}^{-1} \begin{bmatrix} G^\top \\ \kappa H \end{bmatrix} = \Phi + \kappa^{-1} G G^\top + \kappa H^\top R H, \tag{3.13}$$

where (3.13) with $m = p$ is equivalent to the result of the Petersen Lemma (see [452, Lemma C.10.1]). Note that in this case the well-posedness condition (3.7) does not need to be considered given $F = \mathsf{O}_{p,m}$.

Furthermore, consider the case of $F \neq \mathsf{O}_{p,m}$, $\Theta_1^{-1} = R \in \mathbb{S}^p_{\succ 0}$, $\Theta_2 = \mathsf{O}_{p,m}$ and $\Theta_3 = -I_m$ with $\mathcal{D} = \left\{\hat{\Delta} \in \mathbb{R}^{m\times p} : \hat{\Delta}^\top \hat{\Delta} \preceq R\right\}$. In this scenario, Lemma 3.1 generalizes the rational versions of Petersens Lemma in [617] and [452], with the corresponding conditions of (3.7) and (3.9) being

$$\exists \alpha > 0 : \begin{bmatrix} I_m & -I_m & \alpha F^\top \\ * & I_m + \alpha I_m & \mathsf{O}_{m,p} \\ * & * & \alpha R^{-1} \end{bmatrix} \succ 0, \tag{3.14a}$$

$$\exists \kappa > 0 : \begin{bmatrix} \Phi & G & \kappa H^\top \\ * & -\kappa I_m & \kappa F^\top \\ * & * & -\kappa R^{-1} \end{bmatrix} \prec 0. \tag{3.14b}$$

For (3.14a), it is equivalent to

$$I_m - \begin{bmatrix} -I_m & \alpha F^\top \end{bmatrix} \begin{bmatrix} (\alpha+1)^{-1} I_m & \mathsf{O}_{m,p} \\ * & \alpha^{-1} R \end{bmatrix} \begin{bmatrix} -I_m \\ \alpha F \end{bmatrix} = I_m - (\alpha+1)^{-1} I_m - \alpha F^\top R F \succ 0. \tag{3.15}$$

which in turn is equivalent to the proposition: there exists $\alpha > 0$ such that

$$\begin{aligned} \alpha(\alpha+1) F^\top R F + I_m \prec (\alpha+1) I_m \iff \alpha(\alpha+1) F^\top R F - \alpha I_m \prec 0 \\ \iff (\alpha+1)\, F^\top R F - I_m \prec 0. \end{aligned} \tag{3.16}$$

Furthermore, for any $\alpha > 0$, we have $F^\top R F - I_m \preceq (\alpha+1)\, F^\top R F - I_m \prec 0$, since $F^\top R F \preceq (\alpha+1)\, F^\top R F$ for any $\alpha > 0$ based on the fact that $\alpha > 0$ and $F^\top R F \succeq 0$ with $R \succ 0$. Hence we can infer $F^\top R F - I_m \prec 0$[3] from (3.16). This demonstrates that $F^\top R F - I_m \prec 0$ is implied by (3.14a) which further indicates that (3.14a) implies the satisfaction of the well-posedness conditions described in [452] and [617] corresponding to the cases of $m = p$ and $R = I_p$, respectively, are satisfied. Additionally, we can conclude that (3.14b) holds if and only if

[3] It is possible to infer $(\alpha+1) F^\top R F - I_m \prec 0$ from $F^\top R F - I_m \prec 0$ based on the application of eigendecomposition and the property of real numbers.

$$
\exists \kappa > 0: \begin{bmatrix} I_n & \mathsf{O}_{n,m} & \mathsf{O}_{n,p} \\ \mathsf{O}_{p,n} & \mathsf{O}_{p,m} & I_p \\ \mathsf{O}_{m,n} & I_m & \mathsf{O}_{m,p} \end{bmatrix} \begin{bmatrix} \Phi & G & \kappa H^\top \\ G^\top & -\kappa I_m & \kappa F^\top \\ \kappa H & \kappa F & -\kappa R^{-1} \end{bmatrix} \begin{bmatrix} I_n & \mathsf{O}_{n,p} & \mathsf{O}_{n,m} \\ \mathsf{O}_{m,n} & \mathsf{O}_{m,p} & I_m \\ \mathsf{O}_{p,n} & I_p & \mathsf{O}_{p,m} \end{bmatrix}
$$

$$
= \begin{bmatrix} \Phi & G & \kappa H^\top \\ \kappa H & \kappa F & -\kappa R^{-1} \\ G^\top & -\kappa I_m & \kappa F^\top \end{bmatrix} \begin{bmatrix} I_n & \mathsf{O}_{n,p} & \mathsf{O}_{n,m} \\ \mathsf{O}_{m,n} & \mathsf{O}_{m,p} & I_m \\ \mathsf{O}_{p,n} & I_p & \mathsf{O}_{p,m} \end{bmatrix} = \begin{bmatrix} \Phi & \kappa H^\top & G \\ \kappa H & -\kappa R^{-1} & \kappa F \\ G^\top & \kappa F^\top & -\kappa I_m \end{bmatrix} \prec 0 \quad (3.17)
$$

holds given the properties of congruent transformations. By applying the Schur complement to (3.17), it follows that (3.17) is equivalent to

$$
\exists \kappa > 0,\ \ \Phi + \begin{bmatrix} \kappa H^\top & G \end{bmatrix} \kappa^{-1} \begin{bmatrix} R^{-1} & -F \\ -F^\top & I_m \end{bmatrix}^{-1} \begin{bmatrix} \kappa H \\ G^\top \end{bmatrix}
$$

$$
= \Phi + \begin{bmatrix} \sqrt{\kappa} H^\top & \sqrt{\kappa^{-1}} G \end{bmatrix} \begin{bmatrix} R^{-1} & -F \\ -F^\top & I_m \end{bmatrix}^{-1} \begin{bmatrix} \sqrt{\kappa} H \\ \sqrt{\kappa^{-1}} G^\top \end{bmatrix} \prec 0. \quad (3.18)
$$

Now let $m = p$, then (3.18) is equivalent to the result in (C.30) of Lemma C.10.2 in [452], which further shows the equivalence between (3.14b) and (C.31) of Lemma C.10.2 with $m = p$. Moreover, with $\sqrt{\kappa} = \varepsilon^{-1} > 0$ and $R = I_q$, then (3.18) is equivalent to the inequality in (11) of [617].

## 3.3 Main Results on Controller Synthesis

By combining the synthesis results in Theorem 2.1 with (3.3a) and (3.4) and applying Lemma 3.1, we can derive the following theorem, which provides sufficient conditions for the existence of a state feedback controller that ensures both robust dissipativity and stability.

**Theorem 3.1.** *Given $\boldsymbol{f}(\cdot)$ and $M$ in (3.2), the uncertain CLS in (3.5) with supply rate function (2.10) is robustly dissipative subject to the uncertainty constraints in (3.4), and the trivial solution $\boldsymbol{x}(t) \equiv \mathbf{0}_n$ of (3.5) with $\boldsymbol{w}(t) \equiv \mathbf{0}_q$ is robustly globally asymptotically stable subject to the uncertainty constraints in (3.4), if there exist $P \in \mathbb{S}^n$, $Q \in \mathbb{R}^{n \times dn}$, $R \in \mathbb{S}^{dn}$ and $S; U \in \mathbb{S}^n$ and $\varkappa_1, \varkappa_2 > 0$ such that the inequalities in (2.21a) and*

$$
\begin{bmatrix} I & -I - \varkappa_1 \mathbf{F}^\top \mathbf{J}_2 & \varkappa_1 \mathbf{F}^\top \\ * & I - \varkappa_1 \mathbf{J}_3 & \mathsf{O} \\ * & * & \varkappa_1 \mathbf{J}_1 \end{bmatrix} \succ 0, \quad (3.19a)
$$

$$
\begin{bmatrix} \mathbf{P}^\top \mathbf{\Pi} + \mathbf{\Phi} & \mathbf{G} + \varkappa_2 \mathbf{H}^\top \mathbf{J}_2 & \varkappa_2 \mathbf{H}^\top \\ * & \varkappa_2 \mathbf{F}^\top \mathbf{J}_2 + \varkappa_2 \mathbf{J}_2^\top \mathbf{F} + \varkappa_2 \mathbf{J}_3 & \varkappa_2 \mathbf{F}^\top \\ * & * & -\varkappa_2 \mathbf{J}_1 \end{bmatrix} \prec 0, \quad (3.19b)
$$

*are satisfied, where $\mathbf{P}^\top \mathbf{\Pi} + \mathbf{\Phi}$ is given in Theorem 2.1 with the nominal state space parameters in (3.1), and*

$$
\mathbf{G} := \begin{bmatrix} P & \mathsf{O}_{n,m} \\ \mathsf{O}_n & \mathsf{O}_{n,m} \\ Q^\top & \mathsf{O}_{dn,m} \\ \mathsf{O}_{q,n} & -J_2^\top \\ \mathsf{O}_{m,n} & \tilde{J} \end{bmatrix} \begin{bmatrix} \operatorname*{\mathsf{Row}}_{i=1}^{7} G_i & \mathsf{O} \\ \mathsf{O} & \operatorname*{\mathsf{Row}}_{i=8}^{14} G_i \end{bmatrix}, \quad \begin{bmatrix} \mathbf{F} \\ \mathbf{J}_1 \\ \mathbf{J}_2 \\ \mathbf{J}_3 \end{bmatrix} := \begin{bmatrix} \bigoplus_{i=1}^{14} F_i \\ \bigoplus_{i=1}^{14} \Xi_i \\ \bigoplus_{i=1}^{14} \Lambda_i \\ \bigoplus_{i=1}^{14} \Gamma_i \end{bmatrix}, \quad (3.20a)
$$

$$
\mathbf{H} := \begin{bmatrix} \begin{bmatrix} H_1 \\ H_2 K \end{bmatrix} \oplus \begin{bmatrix} H_3 \\ H_4 K \end{bmatrix} \oplus \begin{bmatrix} H_5 \\ H_6 (I_d \otimes K) \end{bmatrix} \oplus H_7 \\ \begin{bmatrix} H_8 \\ H_9 K \end{bmatrix} \oplus \begin{bmatrix} H_{10} \\ H_{11} K \end{bmatrix} \oplus \begin{bmatrix} H_{12} \\ H_{13} (I_d \otimes K) \end{bmatrix} \oplus H_{14} \end{bmatrix} \begin{bmatrix} I_{2n+dn+q} & \mathsf{O}_{m,(2n+dn+q)} \end{bmatrix} \quad (3.20b)
$$

*with $G_i$, $F_i$, $H_i$ and $\Xi_i$, $\Lambda_i$, $\Gamma_i$ given in (3.5) and (3.3c), respectively.*

*Proof.* By substituting uncertain CLS (3.5) into (2.21b), then we have

$$\forall \boldsymbol{\Delta} \in \mathcal{U}, \;\; \mathsf{Sy}\left(\mathbf{P}^\top \boldsymbol{\Pi}\right) + \boldsymbol{\Phi} + \mathsf{Sy}\left[\mathbf{G}\left(\mathbf{I} - \boldsymbol{\Delta}\mathbf{F}\right)^{-1}\boldsymbol{\Delta}\mathbf{H}\right] \prec 0, \tag{3.21}$$

where $\mathcal{U}$ is given in (3.4), and $\mathbb{R}^{\nu_1 \times \nu_2} \ni \boldsymbol{\Delta} := \bigoplus_{i=1}^{14} \Delta_i$ whose dimensions $\nu_1$, $\nu_2$ are determined by the dimensions of $\Delta_i$, $i = 1, \cdots, 14$. Note that (3.21) can be derived based on the structure of (2.21b) and (3.5). A crucial property utilized for the following derivations is that inequality

$$\forall \boldsymbol{D} \in \mathcal{W}, \;\; \mathsf{Sy}\left(\mathbf{P}^\top \boldsymbol{\Pi}\right) + \boldsymbol{\Phi} + \mathsf{Sy}\left[\mathbf{G}\left(\mathbf{I} - \boldsymbol{D}\mathbf{F}\right)^{-1}\boldsymbol{D}\mathbf{H}\right] \prec 0 \tag{3.22a}$$

implies (3.21) with

$$\mathcal{W} := \left\{ \widetilde{\boldsymbol{\Delta}} \in \mathbb{R}^{\nu_1 \times \nu_2} \;\middle|\; \begin{aligned} &\begin{bmatrix} I_{\nu_2} \\ \widetilde{\boldsymbol{\Delta}} \end{bmatrix}^\top \begin{bmatrix} \bigoplus_{i=1}^{14} \Xi_i^{-1} & \bigoplus_{i=1}^{14} \Lambda_i \\ * & \bigoplus_{i=1}^{14} \Gamma_i \end{bmatrix} \begin{bmatrix} I_{\nu_2} \\ \widetilde{\boldsymbol{\Delta}} \end{bmatrix} \succeq 0 \\ &\bigoplus_{i=1}^{14} \Xi_i^{-1} \succ 0, \; \bigoplus_{i=1}^{14} \Gamma_i \preceq 0 \end{aligned} \right\} \tag{3.22b}$$

since $\mathcal{U} \subseteq \mathcal{W}$. Now applying Lemma 3.1 to (3.22a) with (3.4) shows that we can infer (3.22a) from (3.19a) and (3.19b), which further indicates that (3.21) is implied by (3.19a) and (3.19b). This shows that the existence of feasible solutions to (2.21a), (3.19a) and (3.19b) ensures the existence of (2.20) satisfying the corresponding robust version of the stability criteria in (2.7a), (2.7b) and the robust version of the dissipativity in (2.8). These results further imply the robust stability of the trivial solution $\boldsymbol{x}(t) \equiv \mathbf{0}_n$ of (3.5) with $\boldsymbol{w}(t) \equiv \mathbf{0}_q$ and its robust dissipativity with (2.10). ■

**Remark 3.4.** The inequality in (3.19b) is bilinear (non-convex) if $K$ and $B_1, B_2, B_3$ are of nonzero values. On the other hand, the constraints in Theorem 3.1 are convex for a robust dissipativity (stability) analysis problem with $K = \mathsf{O}_{n,q}$ or $B_1 = B_2 = \mathsf{O}_n$ and $B_3 = \mathsf{O}_{n,dn}$.

Similar to Theorem 2.2 in the previous chapter, we specifically derive the following theorem providing a convex optimization-based solution to a genuine robust dissipative control problem.

**Theorem 3.2.** *Given $\boldsymbol{f}(\cdot)$ and $M$ in (3.2) and $\{\alpha_i\}_{i=1}^{2+d} \subset \mathbb{R}$ with $d \in \mathbb{N}$, then the uncertain CLS in (3.5) with supply rate function (2.10) is robustly dissipative subject to the uncertainty constraints in (3.4), and the trivial solution $\boldsymbol{x}(t) \equiv \mathbf{0}_n$ of (3.5) with $\boldsymbol{w}(t) \equiv \mathbf{0}_q$ is robustly globally asymptotically stable subject to (3.4), if there exist $\acute{P} \in \mathbb{S}^n$, $X \in \mathbb{R}^{n \times n}_{[n]}$, $V \in \mathbb{R}^{p \times n}$, $\acute{Q} \in \mathbb{R}^{n \times dn}$, $\acute{R} \in \mathbb{S}^{dn}$, $\acute{S}; \acute{U} \in \mathbb{S}^n$ and $\varkappa_1, \varkappa_2 > 0$ such that (2.33a) and*

$$\begin{bmatrix} I & -I - \varkappa_1 \mathbf{F}^\top \mathbf{J}_2 & \varkappa_1 \mathbf{F}^\top \\ * & I - \varkappa_1 \mathbf{J}_3 & \mathsf{O} \\ * & * & \varkappa_1 \mathbf{J}_1 \end{bmatrix} \succ 0, \quad \begin{bmatrix} \acute{\boldsymbol{\Theta}} & \acute{\mathbf{G}} + \varkappa_2 \mathbf{H}^\top \mathbf{J}_2 & \varkappa_2 \acute{\mathbf{H}}^\top \\ * & \varkappa_2 \mathbf{F}^\top \mathbf{J}_2 + \varkappa_2 \mathbf{J}_2^\top \mathbf{F} + \varkappa_2 \mathbf{J}_3 & \varkappa_2 \mathbf{F}^\top \\ * & * & -\varkappa_2 \mathbf{J}_1 \end{bmatrix} \prec 0, \tag{3.23}$$

*are satisfied, where $\acute{\boldsymbol{\Theta}}$ is defined in (2.33b) and*

$$\acute{\mathbf{G}} := \begin{bmatrix} I_n & \mathsf{O}_{n,m} \\ \mathsf{Col}_{i=1}^{d+2} \alpha_i I_n & \mathsf{O}_{(d+2)n,m} \\ \mathsf{O}_{q,n} & -J_2^\top \\ \mathsf{O}_{m,n} & \widetilde{J} \end{bmatrix} \begin{bmatrix} \operatorname*{\mathsf{Row}}\limits_{i=1}^{7} G_i & \mathsf{O} \\ \mathsf{O} & \operatorname*{\mathsf{Row}}\limits_{i=8}^{14} G_i \end{bmatrix}, \tag{3.24a}$$

$$\acute{\mathbf{H}} := \begin{bmatrix} \begin{bmatrix} H_1X \\ H_2V \end{bmatrix} \oplus \begin{bmatrix} H_3X \\ H_4V \end{bmatrix} \oplus \begin{bmatrix} H_5(I_d \otimes X) \\ H_6\,(I_d \otimes V) \end{bmatrix} \oplus H_7 \\ \begin{bmatrix} H_8X \\ H_9V \end{bmatrix} \oplus \begin{bmatrix} H_{10}X \\ H_{11}V \end{bmatrix} \oplus \begin{bmatrix} H_{12}(I_d \otimes X) \\ H_{13}(I_d \otimes V) \end{bmatrix} \oplus H_{14} \end{bmatrix} \begin{bmatrix} \mathsf{O}_{(2n+dn+q),n} & I_{2n+dn+q} & \mathsf{O}_{(2n+dn+q),m} \end{bmatrix} \tag{3.24b}$$

*with $G_i$, $F_i$, $H_i$ and $\Xi_i$, $\Lambda_i$, $\Gamma_i$ given in (3.5) and (3.3c), respectively.*

*Proof.* By substituting (3.5) into (2.33b) and considering the proof procedure of Theorem 2.2, the robust version of (2.33b) can be derived as

$$\forall \mathbf{\Delta} \in \mathcal{U}, \quad \acute{\mathbf{\Theta}} := \mathsf{Sy}\left( \begin{bmatrix} I_n \\ \mathsf{Col}_{i=1}^{2+d}\, \alpha_i I_n \\ \mathsf{O}_{(q+m),n} \end{bmatrix} \begin{bmatrix} -X & \widehat{\mathbf{\Pi}} \end{bmatrix} \right) + \begin{bmatrix} \mathsf{O}_n & \acute{\mathbf{P}} \\ * & \widehat{\mathbf{\Phi}} \end{bmatrix} \prec 0, \tag{3.25}$$

where $\mathbb{R}^{\nu_1 \times \nu_2} \ni \mathbf{\Delta} := \bigoplus_{i=1}^{14} \Delta_i$ and

$$\begin{aligned} \widehat{\mathbf{\Pi}} = \acute{\mathbf{\Pi}} + \left( \underset{i=1}{\overset{7}{\mathsf{Row}}}\, G_i \right) \left( I - \bigoplus_{i=1}^{7} \Delta_i F_i \right)^{-1} \left( \bigoplus_{i=1}^{7} \Delta_i \right) \\ \times \left( \begin{bmatrix} H_1X \\ H_2V \end{bmatrix} \oplus \begin{bmatrix} H_3X \\ H_4V \end{bmatrix} \oplus \begin{bmatrix} H_5(I_d \otimes X) \\ H_6\,(I_d \otimes V) \end{bmatrix} \oplus H_7 \right) \begin{bmatrix} I_{2n+dn+q} & \mathsf{O}_{(2n+dn+q),m} \end{bmatrix} \end{aligned} \tag{3.26}$$

with $\acute{\mathbf{\Pi}}$, $\acute{\mathbf{P}}$ in line with the definitions in (2.34), and

$$\begin{aligned} \widehat{\mathbf{\Phi}} := \mathsf{Sy}\left( \begin{bmatrix} Q \\ \mathsf{O}_{n,dn} \\ R \\ \mathsf{O}_{q,dn} \end{bmatrix} \begin{bmatrix} \mathbf{F} & \mathsf{O}_{dn,m} \end{bmatrix} \right) + (S + rU) \oplus (-S) \oplus (-\mathsf{F} \otimes U) \oplus (-J_3) \oplus J_1 \\ + \mathsf{Sy}\left( \begin{bmatrix} \mathsf{O}_{(2n+dn),m} \\ -J_2^\top \\ \widetilde{J} \end{bmatrix} \begin{bmatrix} \grave{C}_1X + \grave{B}_4V & \grave{C}_2X + \grave{B}_5V & \grave{C}_3(I_d \otimes X) + \grave{B}_2\,(I_d \otimes V) & \grave{D}_2 & \mathsf{O}_m \end{bmatrix} \right). \end{aligned} \tag{3.27}$$

Now considering the expressions of (3.25)–(3.27), (3.25) can be reformulated as

$$\forall \mathbf{\Delta} \in \mathcal{U}, \quad \acute{\mathbf{\Theta}} + \mathsf{Sy}\left[ \acute{\mathbf{G}} \left( \mathbf{I}_\nu - \mathbf{\Delta}\mathbf{F} \right)^{-1} \mathbf{\Delta} \acute{\mathbf{H}} \right] \prec 0 \tag{3.28}$$

with $\mathbb{R}^{\nu_1 \times \nu_2} \ni \mathbf{\Delta} := \bigoplus_{i=1}^{14} \Delta_i$ and $\acute{\mathbf{G}}$ in (3.24a) and $\acute{\mathbf{H}}$ in (3.24b) and $\acute{\mathbf{\Theta}}$ in (2.33b). Similar to the proof of Theorem 3.1, it is obvious that

$$\forall \acute{\mathbf{\Delta}} \in \mathcal{W}, \quad \acute{\mathbf{\Theta}} + \mathsf{Sy}\left[ \acute{\mathbf{G}} \left( \mathbf{I}_{\nu_1} - \acute{\mathbf{\Delta}}\mathbf{F} \right)^{-1} \acute{\mathbf{\Delta}} \acute{\mathbf{H}} \right] \prec 0 \tag{3.29}$$

infers (3.28) since $\mathcal{U} \subseteq \mathcal{W}$, where $\mathcal{W}$ is given in (3.22b). Also note that the well-posedness of the uncertainty in (3.29) guarantees the well-posedness of the uncertainties in (3.28) since $\mathcal{U} \subseteq \mathcal{W}$. Applying Lemma 3.1 to (3.29) with $\mathcal{W}$ in (3.22b) shows that the existence of feasible solutions to (3.23) ensures the satisfaction of (3.29), which further implies (3.28). The conditions in (3.23) also ensure the well-posedness of the linear fractional uncertainties in (3.29) and (3.28) and ultimately (3.5). From the proof of Theorem 2.2, it shows that the existence of feasible solutions to (2.33a) and (3.23) guarantees the existence of (2.20) satisfying the corresponding robust version of (2.7a), (2.7b) and (2.8). This further ensures the robust stability of the origin of (3.5) with $\boldsymbol{w}(t) \equiv \mathbf{0}_q$, and its robust dissipativity with supply function (2.10). Finally, one can infer the well-posedness of all uncertainties in (3.5) from the first inequality in (3.23). ∎

**Remark 3.5.** It is worth mentioning that all the uncertainties in (3.3a) can be pulled out that results in an interconnection form as demonstrated in [464, 515]. However, in order to generate a single convex condition without requiring the application of the dualization lemma [452, 464], the uncertainties have been chosen in the form of (3.3a).

### 3.3.1 An inner convex approximation solution to Theorem 3.1

Similar to the approach presented in Section 2.4.1, an iterative algorithm is presented in this subsection based on the algorithm proposed in [565] to solve Theorem 3.1. The resulting iterative algorithm can be initiated using the feasible solutions to the synthesis condition in Theorem 2.2.

We first set forth the following lemma to derive the condition for inner convex approximation.

**Lemma 3.2.** *Given* $A \in \mathbb{S}^n$, $B \in \mathbb{R}^{n\times m}$, $C \in \mathbb{S}^m$, $D \in \mathbb{S}^p_{\prec 0}$ *and* $E \in \mathbb{R}^{n\times p}$, *we have*

$$\begin{bmatrix} A - ED^{-1}E^\top & B \\ B^\top & C \end{bmatrix} \prec 0 \iff \begin{bmatrix} A & E & B \\ E^\top & D & \mathsf{O}_{p,n} \\ B^\top & \mathsf{O}_{n,p} & C \end{bmatrix} \tag{3.30}$$

*Proof.* Note that

$$\begin{bmatrix} A - ED^{-1}E^\top & B \\ B^\top & C \end{bmatrix} = \begin{bmatrix} A & B \\ B^\top & C \end{bmatrix} - \begin{bmatrix} E \\ \mathsf{O}_{m,p} \end{bmatrix} D^{-1} \begin{bmatrix} E^\top & \mathsf{O}_{p,m} \end{bmatrix}. \tag{3.31}$$

Applying the Schur complement to (3.31) concludes that (3.31) holds if and only if

$$\begin{bmatrix} A & B & E \\ B^\top & C & \mathsf{O}_{m,p} \\ E^\top & \mathsf{O}_{p,m} & D \end{bmatrix} \prec 0 \tag{3.32}$$

Now apply congruent transformations to (3.3.1), it follows that (3.3.1) holds if and only

$$\begin{aligned}
&\begin{bmatrix} I_n & \mathsf{O}_{n,m} & \mathsf{O}_{n,p} \\ \mathsf{O}_{m,n} & \mathsf{O}_{m,p} & I_m \\ \mathsf{O}_{p,n} & I_p & \mathsf{O}_{p,m} \end{bmatrix} \begin{bmatrix} A & B & E \\ B^\top & C & \mathsf{O}_{m,p} \\ E^\top & \mathsf{O}_{p,m} & D \end{bmatrix} \begin{bmatrix} I_n & \mathsf{O}_{n,m} & \mathsf{O}_{n,p} \\ \mathsf{O}_{m,n} & \mathsf{O}_{m,p} & I_m \\ \mathsf{O}_{p,n} & I_p & \mathsf{O}_{p,m} \end{bmatrix} \\
&= \begin{bmatrix} A & B & E \\ E^\top & \mathsf{O}_{p,m} & D \\ B^\top & C & \mathsf{O}_{m,p} \end{bmatrix} \begin{bmatrix} I_n & \mathsf{O}_{n,m} & \mathsf{O}_{n,p} \\ \mathsf{O}_{m,n} & \mathsf{O}_{m,p} & I_m \\ \mathsf{O}_{p,n} & I_p & \mathsf{O}_{p,m} \end{bmatrix} = \begin{bmatrix} A & E & B \\ E^\top & D & \mathsf{O}_{p,n} \\ B^\top & \mathsf{O}_{n,p} & C \end{bmatrix}
\end{aligned} \tag{3.33}$$

This finishes the proof. ■

Similar to the structure of (2.44), the inequality in (3.19b) can be rewritten as

$$\begin{bmatrix} \mathsf{Sy}\left[\mathbf{P}^\top \mathbf{B}\left[(I_{2+d} \otimes K) \oplus \mathsf{O}_{p+m}\right]\right] + \widehat{\boldsymbol{\Phi}} & \mathbf{G} + \varkappa_2 \mathbf{H}^\top \mathbf{J}_2 & \varkappa_2 \mathbf{H}^\top \\ * & \varkappa_2 \mathbf{F}^\top \mathbf{J}_2 + \varkappa_2 \mathbf{J}_2^\top \mathbf{F} + \varkappa_2 \mathbf{J}_3 & \varkappa_2 \mathbf{F}^\top \\ * & * & -\varkappa_2 \mathbf{J}_1 \end{bmatrix} \prec 0, \tag{3.34}$$

where $\mathbf{B}$ and $\widehat{\boldsymbol{\Phi}}$ are defined in (2.44), and the other parameters are detailed in Theorem 3.1. Now by applying (2.46a) with (2.47) and the conclusion in (3.30) to (3.34), one can infer (3.34) from

$$\begin{bmatrix} \widehat{\boldsymbol{\Phi}} + \mathsf{Sy}\left(\widetilde{\mathbf{P}}^\top \mathbf{B}\mathbf{K} + \mathbf{P}^\top \mathbf{B}\widetilde{\mathbf{K}} - \widetilde{\mathbf{P}}^\top \mathbf{B}\widetilde{\mathbf{K}}\right) & \mathbf{P}^\top - \widetilde{\mathbf{P}}^\top & \mathbf{K}^\top \mathbf{B}^\top - \widetilde{\mathbf{K}}^\top \mathbf{B}^\top & \mathbf{G} + \varkappa_2 \mathbf{H}^\top \mathbf{J}_2 & \varkappa_2 \mathbf{H}^\top \\ * & -Z & \mathsf{O}_n & \mathsf{O} & \mathsf{O} \\ * & * & Z - I_n & \mathsf{O} & \mathsf{O} \\ * & * & * & \varkappa_2 \mathbf{F}^\top \mathbf{J}_2 + \varkappa_2 \mathbf{J}_2^\top \mathbf{F} + \varkappa_2 \mathbf{J}_3 & \varkappa_2 \mathbf{F}^\top \\ * & * & * & * & -\varkappa_2 \mathbf{J}_1 \end{bmatrix} \prec 0, \tag{3.35}$$

where all the matrices in (3.35) align with the definitions in Theorem 3.1 and (2.47).

By using the procedures expounded in [565], Algorithm 2 is constructed analogous to Algorithm 1, where $\mathbf{x}$ contains all the decision variables of $R \in \mathbb{S}^{nd}$ and $S; U; Z \in \mathbb{S}^n$. Moreover, let $\mathbf{\Lambda} := \begin{bmatrix} P & Q \end{bmatrix}$ and $\widetilde{\mathbf{\Lambda}} := \begin{bmatrix} \widetilde{P} & \widetilde{Q} \end{bmatrix}$. Furthermore, $\rho_1$, $\rho_2$ and $\varepsilon$ are given constants utilized for regularizations and determining error tolerance, respectively.

---

**Algorithm 2:** An inner convex approximation solution to Theorem 3.1

---

**begin**

Find initial values for $\widetilde{\mathbf{\Lambda}} = \begin{bmatrix} \widetilde{P} & \widetilde{Q} \end{bmatrix}$ and $\widetilde{K}$ belonging to the relative interior of the feasible set of Theorem 3.1.

**solve** $\min\limits_{\mathbf{x},\mathbf{\Lambda},K} \operatorname{tr}\left[\rho_1[*](\mathbf{\Lambda} - \widetilde{\mathbf{\Lambda}}) + \rho_2[*](K - \widetilde{K})\right]$ subject to (2.21a), (3.19a) and (3.35) to obtain $\mathbf{\Lambda}$ and $K$

**while** $\dfrac{\left\| \begin{bmatrix} \mathbf{vec}(\mathbf{\Lambda}) \\ \mathbf{vec}\,(K) \end{bmatrix} - \begin{bmatrix} \mathbf{vec}(\widetilde{\mathbf{\Lambda}}) \\ \mathbf{vec}(\widetilde{K}) \end{bmatrix} \right\|_\infty}{\left\| \begin{bmatrix} \mathbf{vec}(\widetilde{\mathbf{\Lambda}}) \\ \mathbf{vec}(\widetilde{K}) \end{bmatrix} \right\|_\infty + 1} \geq \varepsilon$ **do**

**update** $\widetilde{\mathbf{\Lambda}} \longleftarrow \mathbf{\Lambda}, \quad \widetilde{K} \longleftarrow K;$

**solve** $\min\limits_{\mathbf{x},\mathbf{\Lambda},K} \operatorname{tr}\left[\rho_1[*](\mathbf{\Lambda} - \widetilde{\mathbf{\Lambda}}) + \rho_2[*](K - \widetilde{K})\right]$ subject to (2.21a), (3.19a) and (3.35) to obtain $\mathbf{\Lambda}$ and $K$;

**end**

**end**

---

**Remark 3.6.** When Theorem 3.1 contains a convex objective function, for instance $\mathcal{L}^2$ gain $\gamma$ minimization, a termination condition can be added to Algorithm 2 to regulate the minimization of objective function between two successive iterations [565]. Nonetheless, this condition was not implemented in our numerical examples in this chapter, as we explicitly controlled the number of iterations in our numerical examples.

**Remark 3.7.** The Initial values for $\widetilde{P}$, $\widetilde{Q}$ and $\widetilde{K}$ in Algorithm 2 can be provided by some feasible solutions to (2.33a) and (3.23) with given $\{\alpha_i\}_{i=1}^{2+d}$. As explained in Remark 2.4 in the previous chapter, one can apply $\alpha_1 = \alpha_2 = 0$ and $\alpha_i = 0, i = 4, \cdots, 2+d$, which allows users to adjust only the value of $\alpha_3$ when Theorem 3.2 is employed.

## 3.4 Application to the Dissipative Resilient Stabilizations of A Linear System with A Pointwise Input Delay

In this section, we demonstrate that the idea presented in Section 3.3 can be adapted to handle the dissipative stabilization of a linear input delay system by employing a dynamical state feedback controller [619] under both general uncertainties and a dissipative constraint. Due to the mathematical structure of the CLS with a dynamical state feedback controller, we can use the synthesis conditions proposed in this section to compute the gains of a resilient controller that is robust against uncertainties. Moreover, unlike Algorithm 2 that may require a feasible solution to Theorem 3.2 for its initialization, Algorithm 3 may be initialized based on the gains of predictor controllers that are always constructible.

### 3.4.1 Problem formulation

Consider the following system

$$\dot{\boldsymbol{x}}(t) = \grave{A}\boldsymbol{x}(t) + \grave{B}\boldsymbol{u}(t-r) + \grave{D}_1\boldsymbol{w}(t), \quad t \geq t_0$$
$$\boldsymbol{z}(t) = \grave{C}_1 \begin{bmatrix} \boldsymbol{x}(t) \\ \boldsymbol{u}(t) \end{bmatrix} + \grave{C}_2 \begin{bmatrix} \boldsymbol{x}(t-r) \\ \boldsymbol{u}(t-r) \end{bmatrix} + \int_{-r}^{0} \grave{C}_3 \left(\sqrt{\mathsf{F}}\boldsymbol{f}(\tau) \otimes I_\nu\right) \begin{bmatrix} \boldsymbol{x}(t+\tau) \\ \boldsymbol{u}(t+\tau) \end{bmatrix} \mathsf{d}\tau + \grave{D}_2\boldsymbol{w}(t)$$

where $r > 0$ and $t_0 \in \mathbb{R}$ and $\boldsymbol{x}(t) \in \mathbb{R}^n$, $\boldsymbol{u}(t) \in \mathbb{R}^p$ and $\boldsymbol{w}(t) \in \mathbb{R}^q$ and $\boldsymbol{f}(\cdot) \in \mathcal{C}^1\left([-r, 0]; \mathbb{R}^d\right)$ with $\mathsf{F}^{-1} = \int_{-r}^{0} \boldsymbol{f}(\tau)\boldsymbol{f}^\top(\tau)\mathsf{d}\tau \succ 0$ and $\sqrt{\mathsf{F}}$ stands for the unique[4] square root of $\mathsf{F}$. The state space matrices in (3.36) contain uncertainties. Specifically, $A \in \mathbb{R}^{n\times n}$ and $B \in \mathbb{R}^{n\times p}$ are the nominal part of the state space matrices $\grave{A}$ and $\grave{B}$ without uncertainties[5], and we assume that there exists $K \in \mathbb{R}^{p\times n}$ such that $A + BK$ is Hurwitz. The matrix parameters in the output are defined as

$$\begin{bmatrix} \grave{C}_1 & \grave{C}_2 & \grave{C}_3 & \grave{D}_2 \end{bmatrix} = \begin{bmatrix} C_1 & C_2 & C_3 & D_2 \end{bmatrix} + \underset{i=5}{\overset{8}{\mathsf{Row}}} \begin{bmatrix} G_i(I - \Delta_i F_i)^{-1}\Delta H_i \end{bmatrix} \tag{3.36}$$

where $C_1; C_2 \in \mathbb{R}^{m\times\nu}$, $C_3 \in \mathbb{R}^{m\times d\nu}$ with $\nu = n+p$. Furthermore, $\boldsymbol{f}(\cdot) \in \mathcal{C}^1\left(\mathbb{R}; \mathbb{R}^d\right)$ in (3.36) satisfies

$$\exists K \in \mathbb{R}^{p\times n}, \exists X \in \mathbb{R}^{p\times p},\ \exists \widehat{K} \in \mathbb{R}^{p\times dp},\ \begin{bmatrix} \mathsf{O}_{p,n} & (KA - XK)\,\mathsf{e}^{-A\tau}B \end{bmatrix} = \widehat{K}\left(\sqrt{\mathsf{F}}\boldsymbol{f}(\tau) \otimes I_\nu\right) \tag{3.37}$$
$$\exists M \in \mathbb{R}^{d\times d} : \frac{\mathsf{d}\boldsymbol{f}(\tau)}{\mathsf{d}\tau} = M\boldsymbol{f}(\tau) \tag{3.38}$$

where the motivation of having the structure in (3.37) is explained later in light of the structure of predictor controllers.

Now, in order to stabilize (3.36), we consider a dynamical state feedback controller

$$\dot{\boldsymbol{u}}(t) = K_1 \begin{bmatrix} \boldsymbol{x}(t) \\ \boldsymbol{u}(t) \end{bmatrix} + K_2 \begin{bmatrix} \boldsymbol{x}(t-r) \\ \boldsymbol{u}(t-r) \end{bmatrix} + \int_{-r}^{0} K_3 \left(\sqrt{\mathsf{F}}\boldsymbol{f}(\tau) \otimes I_\nu\right) \begin{bmatrix} \boldsymbol{x}(t+\tau) \\ \boldsymbol{u}(t+\tau) \end{bmatrix} \mathsf{d}\tau \tag{3.39}$$

and assume that all the states are measurable, where $K_1; K_2 \in \mathbb{R}^{p\times\nu}$ and $K_3 \in \mathbb{R}^{p\times d\nu}$ are controller gains with $\nu = n+p$. Owing to the uncertainties caused by implementation environments and disturbances, it is more realistic to consider

$$\dot{\boldsymbol{u}}(t) = \grave{K}_1 \begin{bmatrix} \boldsymbol{x}(t) \\ \boldsymbol{u}(t) \end{bmatrix} + \grave{K}_2 \begin{bmatrix} \boldsymbol{x}(t-r) \\ \boldsymbol{u}(t-r) \end{bmatrix} + \int_{-r}^{0} \grave{K}_3 \left(\sqrt{\mathsf{F}}\boldsymbol{f}(\tau) \otimes I_\nu\right) \begin{bmatrix} \boldsymbol{x}(t+\tau) \\ \boldsymbol{u}(t+\tau) \end{bmatrix} \mathsf{d}\tau + \grave{D}_3\boldsymbol{w}(t), \tag{3.40}$$

as the mathematical model of (3.39) for the theoretical analysis taking into account uncertainties, where $\grave{K}_1$,$\grave{K}_2$,$\grave{K}_3$ and $\grave{D}_3$ contain uncertainties.

**Remark 3.8.** It is important to note that (3.40) serves only as a mathematical model for (3.39), which is utilized in the theoretical derivation of the robust synthesis condition in this section. When the values of $K_1$,$K_2$ and $K_3$ are obtained using the proposed scenarios, the resulting controller is implemented using the model in (3.39). Thus no uncertainties or $\boldsymbol{w}(t)$ are included in the actual implementation of the resulting controller in (3.39). On the other hand, since uncertainties and disturbances are taken into account in (3.40), the resulting controller in (3.39) is resilient and can withstand external perturbations.

[4]For the uniqueness of the square root of a positive semidefinite matrix, see [566, Theorem 7.2.6]

[5]The exact structures of the uncertainties in $\grave{A}$ and $\grave{B}$ are specified later

Now combining (3.36) and (3.40) with (3.36) produces the following CLS

$$
\begin{aligned}
&\begin{bmatrix}\dot{\boldsymbol{x}}(t)\\ \dot{\boldsymbol{u}}(t)\end{bmatrix} = \begin{bmatrix}\begin{bmatrix}\grave{A} & \mathsf{O}_{n,p}\\ \multicolumn{2}{c}{\grave{K}_1}\end{bmatrix} & \begin{bmatrix}\mathsf{O}_n & \grave{B}\\ \multicolumn{2}{c}{\grave{K}_2}\end{bmatrix} & \begin{bmatrix}\mathsf{O}_{n,d\nu}\\ \grave{K}_3\end{bmatrix} & \begin{bmatrix}\grave{D}_1\\ \grave{D}_3\end{bmatrix}\end{bmatrix}\boldsymbol{\vartheta}(t), \quad t \geq t_0\\
&\boldsymbol{z}(t) = \begin{bmatrix}\grave{C}_1 & \grave{C}_2 & \grave{C}_3 & \grave{D}_2\end{bmatrix}\boldsymbol{\vartheta}(t)\\
&\forall\theta\in[-r,0], \quad \begin{bmatrix}\boldsymbol{x}(t_0+\theta)\\ \boldsymbol{u}(t_0+\theta)\end{bmatrix} = \widehat{\boldsymbol{\phi}}(\theta)
\end{aligned}
\tag{3.41}
$$

where $r > 0$ and $t_0 \in \mathbb{R}$ and $\widehat{\boldsymbol{\phi}}(\cdot) \in \mathcal{C}\left([-r,0];\mathbb{R}^{n+p}\right)$ and

$$
\boldsymbol{\vartheta}(t) = \mathbf{Col}\begin{bmatrix}\begin{bmatrix}\boldsymbol{x}(t)\\ \boldsymbol{u}(t)\end{bmatrix}, & \begin{bmatrix}\boldsymbol{x}(t-r)\\ \boldsymbol{u}(t-r)\end{bmatrix}, & \displaystyle\int_{-r}^{0}\left(\sqrt{\mathsf{F}}\boldsymbol{f}(\tau)\otimes I_\nu\right)\begin{bmatrix}\boldsymbol{x}(t+\tau)\\ \boldsymbol{u}(t+\tau)\end{bmatrix}, & \boldsymbol{w}(t)\end{bmatrix}
\tag{3.42}
$$

and the matrices in (3.41) are defined as

$$
\begin{aligned}
&\begin{bmatrix}\begin{bmatrix}\grave{A} & \mathsf{O}_{n,p}\\ \multicolumn{2}{c}{\grave{K}_1}\end{bmatrix} & \begin{bmatrix}\mathsf{O}_n & \grave{B}\\ \multicolumn{2}{c}{\grave{K}_2}\end{bmatrix} & \begin{bmatrix}\mathsf{O}_{n,d\nu}\\ \grave{K}_3\end{bmatrix} & \begin{bmatrix}\grave{D}_1\\ \grave{D}_3\end{bmatrix}\end{bmatrix} = \begin{bmatrix}\begin{bmatrix}A & \mathsf{O}_{n,p}\\ \multicolumn{2}{c}{K_1}\end{bmatrix} & \begin{bmatrix}\mathsf{O}_n & B\\ \multicolumn{2}{c}{K_2}\end{bmatrix} & \begin{bmatrix}\mathsf{O}_{n,d\nu}\\ K_3\end{bmatrix} & \begin{bmatrix}D_1\\ D_3\end{bmatrix}\end{bmatrix}\\
&\qquad\qquad + \mathop{\mathbf{Row}}_{i=1}^{4}\left(G_i(I-\Delta_i F_i)^{-1}\Delta H_i\right)\\
&\begin{bmatrix}\grave{C}_1 & \grave{C}_2 & \grave{C}_3 & \grave{D}_2\end{bmatrix} = \begin{bmatrix}C_1 & C_2 & C_3 & D_2\end{bmatrix} + \mathop{\mathbf{Row}}_{i=5}^{8}\begin{bmatrix}G_i(I-\Delta_i F_i)^{-1}\Delta H_i\end{bmatrix}.
\end{aligned}
\tag{3.43}
$$

The uncertainties $\Delta_i$, $i = 1,\cdots,8$ in (3.43) and (3.36) are subject to the constraints

$$
\Delta_i \in \left\{\widehat{\Delta}_i \;\middle|\; \begin{bmatrix}I\\ \widehat{\Delta}_i\end{bmatrix}^\top\begin{bmatrix}\Xi_i^{-1} & \Lambda_i\\ * & \Gamma_i\end{bmatrix}\begin{bmatrix}I\\ \widehat{\Delta}_i\end{bmatrix} \succeq 0\right\}, \quad \forall i = 1,\cdots,8, \quad \Xi_i^{-1} \succ 0, \quad \Gamma_i \preceq 0.
\tag{3.44}
$$

Note that the matrices $G_i$, $H_i$, $F_i$, $\Delta_i$ in (3.43) and $\Xi_i$, $\Lambda_i$, $\Gamma_i$ in (3.44) are locally defined in this section and are not the same as those in (3.3a), (3.3b) and (3.3c). Meanwhile, the value of $D_3 \in \mathbb{R}^{p,q}$ in (3.43) is given. Finally, similar to (3.4), the constraints in (3.44) can be reformulated as

$$
\begin{gathered}
\bigoplus_{i=1}^{8}\Delta_i \in \mathcal{T} := \left\{\bigoplus_{i=1}^{8}\widehat{\Delta}_i \;\middle|\; [*]\begin{bmatrix}\bigoplus_{i=1}^{8}\Xi_i^{-1} & \bigoplus_{i=1}^{8}\Lambda_i\\ * & \bigoplus_{i=1}^{8}\Gamma_i\end{bmatrix}\begin{bmatrix}I\\ \bigoplus_{i=1}^{8}\widehat{\Delta}_i\end{bmatrix} \succeq 0\right\},\\
\bigoplus_{i=1}^{8}\Xi_i^{-1} \succ 0, \quad \bigoplus_{i=1}^{8}\Gamma_i \preceq 0.
\end{gathered}
\tag{3.45}
$$

**Remark 3.9.** The constraints in (3.37) indicate that the elements in $\boldsymbol{f}(\cdot)$ cover the functions in $\mathrm{e}^{-A\tau}B$ satisfying $\frac{\mathrm{d}\mathrm{e}^{-A\tau}B}{\mathrm{d}\tau} = -A\mathrm{e}^{-A\tau}B$ where the functions in $-A\mathrm{e}^{-A\tau}B$ are inherently consistent with the property described in (3.38).[6] Moreover, this also means that (3.37) does not impose any additional restriction on the structure $A$, $B$ as long as a $K$ can be found to make $A + BK$ Hurwitz. Finally, we stress that the use of $\sqrt{\mathsf{F}}$ does not affect the existence of $C_3$ and $\Gamma$, given that $\sqrt{\mathsf{F}}$ is a symmetrical and full rank matrix.

**Remark 3.10.** Unlike the scenario of addressing a standard state feedback synthesis problem in (3.5), the controller parameters in (3.43) are not directly multiplied by the system parameters in (3.36). For this reason, it is possible to produce resilient dissipative stabilization results with (3.43) by leveraging the conclusions in (3.3).

[6]This can be understood from the property of $\mathrm{e}^{-A\tau}$ in light of the Putzer Algorithm for Matrix Exponentials

**Remark 3.11.** In [620], DDs with uncertainties are investigated based on the structure of matrix exponential functions $(C + \Delta C)\mathsf{e}^{(\Delta A + A)\tau}(B + \Delta B)$. Obviously, the structures of the uncertainties in [620] can be more general than the structures in (3.43) and (3.36), albeit the constraints on the uncertainties in [620, eq.16.5] are simpler than the characterization in (3.45). However, it can be argued that the structure of $(C+\Delta C)\mathsf{e}^{(\Delta A+A)\tau}(B+\Delta B)$ is unnecessary given that the implementation of a predictor type controller is not subject to the uncertainties of the open-loop plant parameters $\Delta A, \Delta B, \Delta C$. Therefore, the linear fractional uncertainties in (3.43) and (3.36) are sufficiently general given the discussions in Remark 3.2. Moreover, one may argue that the uncertainties in (3.43) and (3.36) are more general than those in [620] since $\Delta_i$ here are independent of each other.

**Remark 3.12.** The CLSs in (3.41) can be addressed by the methodologies proposed in the previous section. In addition, since the system in (3.41) is of the retarded type, it adheres to the properties stated in [163, Theorem 2].This implies that if the trivial solution $\boldsymbol{x}(t) \equiv \mathbf{0}_n$ of (3.41) with $\boldsymbol{w}(t) \equiv 0$ is robustly globally uniformly asymptotically stable with the controller in (3.40), then the integral in (3.39) can always be numerically implemented in real-time without potential stability problems, as long as the accuracy of the approximation of the integral in (3.39) reaches certain thresholds.

### 3.4.2 An example of dynamical state feedback controllers

The motivation behind the assumption in (3.37) can be elucidated in relation to the expression in (3.40). Consider system $\dot{\boldsymbol{x}}(t) = A\boldsymbol{x}(t) + B\boldsymbol{u}(t-r)$ with $r > 0$ and $\boldsymbol{x}(t) \in \mathbb{R}^n$ and $\boldsymbol{u}(t) \in \mathbb{R}^p$. This system can always be exponentially stabilized by the predictor controller

$$\boldsymbol{u}(t) = K\boldsymbol{x}(t+r) = K\left(\mathsf{e}^{Ar}\boldsymbol{x}(t) + \int_{-r}^{0} \mathsf{e}^{-A\tau}B\boldsymbol{u}(t+\tau)\mathsf{d}\tau\right) \tag{3.46}$$

for any $r > 0$, as long as $A + BK$ is Hurwitz for some $K \in \mathbb{R}^{p\times n}$. However, the form of (3.46) may not guarantee safe numerical implementation. A solution to this issue was proposed in [163], where the author demonstrated that a dynamical state feedback controller with special forms can achieve the same objective as a predictor controller. By substituting appropriate parameters into eq.(10) of [163], one can conclude from the analysis in the frequency domain that the dynamical state feedback controller

$$\dot{\boldsymbol{u}}(t) = (KB + X)\,\boldsymbol{u}(t) + (KA - XK)\left(\mathsf{e}^{Ar}\boldsymbol{x}(t) + \int_{-r}^{0} \mathsf{e}^{-A\tau}B\boldsymbol{u}(t+\tau)\mathsf{d}\tau\right) \tag{3.47}$$

can asymptotically stabilize $\dot{\boldsymbol{x}}(t) = A\boldsymbol{x}(t) + B\boldsymbol{u}(t-r)$ for any $r > 0$, where $A + BK$ and $X \in \mathbb{R}^{p\times p}$ are Hurwitz, and the form of (3.47) is a particular case of the controller structure in (10) of [163]. Specifically, the spectrum of $\dot{\boldsymbol{x}}(t) = A\boldsymbol{x}(t) + B\boldsymbol{u}(t-r)$ with (3.47) is

$$\left\{s \in \mathbb{C} : \det\left(sI_n - A - BK\right)\det\left(sI_p - X\right) = 0\right\}, \tag{3.48}$$

which is exponentially stabilizable for some $K \in \mathbb{R}^{p\times n}$ given $X \in \mathbb{R}^{p\times p}$ is Hurwitz.

By examining the matrices in (3.47), it becomes clear that (3.37) signifies that function $\boldsymbol{f}(\cdot)$ in (3.36) is able to cover all the functions in $\mathsf{e}^{-A\tau}B$ and $\widetilde{C}_3(\tau)$, where $A$ and $B$ are given in (3.36). Thus, the condition in (3.37) indicates that it is always possible to construct a predictor controller in the form of (3.47) to stabilize the nominal system of (3.36) without taking into account uncertainties and disturbance $\boldsymbol{w}(\cdot)$. Additionally, the structure of (3.39) demonstrates that it encompasses (3.47) as a special case. This implies that we can always find controller gains in (3.39) that can exponentially stabilize the nominal system of (3.36) with $\boldsymbol{w}(t) \equiv 0$. Specifically, (3.47) can be expressed in the form of (3.39) as

$$K_1 = \begin{bmatrix} (KA - XK)\,\mathsf{e}^{Ar} & KB + X \end{bmatrix}, \quad K_2 = \mathsf{O}_{p\times\nu}, \quad K_3 = \widehat{K} \tag{3.49}$$

where $\widehat{K}$ is given in (3.37) and both $A + BK$ and $X$ are Hurwitz.

Based on Theorem 3.1, the following theorem can be derived for the CLS in (3.41).

**Theorem 3.3.** *Let the parameters $G_i$, $F_i$, $H_i$ and $\Xi_i$, $\Lambda_i$, $\Gamma_i$ in (3.43) and (3.44) be given. Given $\boldsymbol{f}(\cdot)$ with $\mathsf{F}^{-1} = \int_{-r}^{0} \boldsymbol{f}(\tau)\boldsymbol{f}^\top(\tau)\mathsf{d}\tau$ and $M$ in (3.37)–(3.38), if there exist $P \in \mathbb{S}^\nu$, $Q \in \mathbb{R}^{\nu\times d\nu}$, $R \in \mathbb{S}^{d\nu}$, $S; U \in \mathbb{S}^\nu$ and $\begin{bmatrix} K_1 & K_2 & K_3 \end{bmatrix} \in \mathbb{R}^{p\times(2\nu+d\nu)}$ and $\varkappa_1, \varkappa_2 > 0$ such that the following matrix inequalities*

$$\begin{bmatrix} P & Q \\ * & R + I_d \otimes S \end{bmatrix} \succ 0, \quad S \succeq 0, \quad U \succeq 0 \tag{3.50}$$

$$\begin{bmatrix} I & -I - \varkappa_1 \mathbf{F}^\top \mathbf{J}_2 & \varkappa_1 \mathbf{F}^\top \\ * & I - \varkappa_1 \mathbf{J}_3 & \mathsf{O} \\ * & * & \varkappa_1 \mathbf{J}_1 \end{bmatrix} \succ 0 \tag{3.51}$$

$$\begin{bmatrix} \mathbf{\Phi} + \mathsf{Sy}\left[\mathbf{P}^\top(\mathbf{A} + \mathbf{B}\mathbf{K})\right] & \mathbf{G} + \varkappa_2 \mathbf{H}^\top \mathbf{J}_2 & \varkappa_2 \mathbf{H}^\top \\ * & \varkappa_2 \mathbf{F}^\top \mathbf{J}_2 + \varkappa_2 \mathbf{J}_2^\top \mathbf{F} + \varkappa_2 \mathbf{J}_3 & \varkappa_2 \mathbf{F}^\top \\ * & * & -\varkappa_2 \mathbf{J}_1 \end{bmatrix} \prec 0 \tag{3.52}$$

*are satisfied with*

$$\mathbf{G} := \begin{bmatrix} P & \mathsf{O}_{\nu,m} \\ \mathsf{O}_\nu & \mathsf{O}_{\nu,m} \\ Q^\top & \mathsf{O}_{d\nu,m} \\ \mathsf{O}_{q,\nu} & -J_2^\top \\ \mathsf{O}_{m,\nu} & \widetilde{J} \end{bmatrix} \left( \begin{bmatrix} \mathop{\mathbf{Row}}\limits_{i=1}^{4} G_i & \mathsf{O} \\ \mathsf{O} & \mathop{\mathbf{Row}}\limits_{i=4}^{8} G_i \end{bmatrix} \right), \quad \begin{bmatrix} \mathbf{F} \\ \mathbf{J}_1 \\ \mathbf{J}_2 \\ \mathbf{J}_3 \end{bmatrix} := \begin{bmatrix} \bigoplus_{i=1}^{8} F_i \\ \bigoplus_{i=1}^{8} \Xi_i \\ \bigoplus_{i=1}^{8} \Lambda_i \\ \bigoplus_{i=1}^{8} \Gamma_i \end{bmatrix} \tag{3.53}$$

$$\mathbf{H} := \begin{bmatrix} \bigoplus_{i=1}^{4} H_i \\ \bigoplus_{i=4}^{8} H_i \end{bmatrix} \begin{bmatrix} I_{2\nu+d\nu+q} & \mathsf{O}_{(2\nu+d\nu+q),m} \end{bmatrix} \tag{3.54}$$

$$\mathbf{P} := \begin{bmatrix} P & \mathsf{O}_\nu & Q & \mathsf{O}_{\nu,q} & \mathsf{O}_{\nu,m} \end{bmatrix}, \quad \mathbf{A} := \begin{bmatrix} A & \mathsf{O}_{n,\nu} & B & \mathsf{O}_{n,d\nu} & D_1 & \mathsf{O}_{n,m} \\ \mathsf{O}_{p,n} & \mathsf{O}_{p,\nu} & \mathsf{O}_p & \mathsf{O}_{p,d\nu} & D_2 & \mathsf{O}_{p,m} \end{bmatrix} \tag{3.55}$$

$$\mathbf{B} = \mathbf{Col}\left[\mathsf{O}_{n,p}, \; I_p\right], \quad \mathbf{K} := \begin{bmatrix} K_1 & K_2 & K_3 & \mathsf{O}_{p,(q+m)} \end{bmatrix} \tag{3.56}$$

$$\begin{aligned} \mathbf{\Phi} := \mathsf{Sy}&\left( \begin{bmatrix} Q \\ \mathsf{O}_{\nu,d\nu} \\ R \\ \mathsf{O}_{(q+m),d\nu} \end{bmatrix} \begin{bmatrix} \sqrt{\mathsf{F}}\boldsymbol{f}(0) \otimes I_\nu & -\sqrt{\mathsf{F}}\boldsymbol{f}(-r) \otimes I_\nu & -\left(\sqrt{\mathsf{F}} M \sqrt{\mathsf{F}^{-1}}\right) \otimes I_\nu & \mathsf{O}_{d\nu,(q+m)} \end{bmatrix} \right) \\ &+ (S + rU) \oplus (-S) \oplus (-I_d \otimes U) \oplus J_3 \oplus J_1 \\ &+ \mathsf{Sy}\left( \begin{bmatrix} \mathsf{O}_{(2\nu+d\nu),m} \\ -J_2^\top \\ \widetilde{J} \end{bmatrix} \begin{bmatrix} C_1 & C_2 & C_3 & D_3 & \mathsf{O}_m \end{bmatrix} \right), \end{aligned} \tag{3.57}$$

*then the CLS in (3.41) with supply rate function (2.10) is dissipative, and the origin of (3.41) with $\boldsymbol{w}(t) \equiv \mathbf{0}_q$ is robustly globally asymptotically stable subject to the uncertainty constraints in (3.44).*

*Proof.* The proof is straightforward given the results presented in Theorem 3.1 and Theorem 2.1. ∎

Based on the derivation of Theorem 3.1 was derived, an iterative algorithm can be constructed to solve Theorem 3.3. Given the results derived in Section 3.3.1, we see that (3.52) holds if

$$\begin{bmatrix} \mathbf{\Phi} + \mathsf{Sy}\left(\widetilde{\mathbf{P}}^\top\mathbf{B}\mathbf{K} + \mathbf{P}^\top\mathbf{B}\widetilde{\mathbf{K}} - \widetilde{\mathbf{P}}^\top\mathbf{B}\widetilde{\mathbf{K}}\right) & \mathbf{P}^\top - \widetilde{\mathbf{P}}^\top & \mathbf{K}^\top\mathbf{B}^\top - \widetilde{\mathbf{K}}^\top\mathbf{B}^\top & \mathbf{G} + \varkappa_2\mathbf{H}^\top\mathbf{J}_2 & \varkappa_2\mathbf{H}^\top \\ * & -Z & \mathsf{O}_\nu & \mathsf{O} & \mathsf{O} \\ * & * & Z - I_\nu & \mathsf{O} & \mathsf{O} \\ * & * & * & \varkappa_2\mathbf{F}^\top\mathbf{J}_2 + \varkappa_2\mathbf{J}_2^\top\mathbf{F} + \varkappa_2\mathbf{J}_3 & \varkappa_2\mathbf{F}^\top \\ * & * & * & * & -\varkappa_2\mathbf{J}_1 \end{bmatrix} \prec 0, \quad (3.58)$$

where $Z \in \mathbb{S}^\nu$ and

$$\begin{aligned} \widetilde{\mathbf{P}} &:= \begin{bmatrix} \widetilde{P} & \mathsf{O}_\nu & \widetilde{Q} & \mathsf{O}_{\nu,q} & \mathsf{O}_{\nu,m} \end{bmatrix} \quad \text{with} \quad \widetilde{P} \in \mathbb{S}^\nu \quad \text{and} \quad \widetilde{Q} \in \mathbb{R}^{\nu \times d\nu} \\ \widetilde{\mathbf{K}} &:= \begin{bmatrix} \widetilde{K}_1 & \widetilde{K}_2 & \widetilde{K}_3 & \mathsf{O}_{p,(q+m)} \end{bmatrix} \quad \text{with} \quad \widetilde{K} = \begin{bmatrix} \widetilde{K}_1 & \widetilde{K}_2 & \widetilde{K}_3 \end{bmatrix} \in \mathbb{R}^{p,(2\nu+d\nu)} \end{aligned} \quad (3.59)$$

and all other matrices in (3.58) are in line with the definitions in Theorem 3.3. With the derivations in Section 3.3.1 in mind, Algorithm 3 can be constructed analogously to Algorithm 2, where $\mathbf{x}$ contains all the decision variables of $R \in \mathbb{S}^{d\nu}$ and $S; U; Z \in \mathbb{S}^\nu$, while $\mathbf{\Lambda} := \begin{bmatrix} P & Q \end{bmatrix} \in \mathbb{R}^{\nu \times \nu(d+1)}$ and $\widetilde{\mathbf{\Lambda}} := \begin{bmatrix} \widetilde{P} & \widetilde{Q} \end{bmatrix} \in \mathbb{R}^{\nu \times \nu(d+1)}$. Furthermore, $\rho_1$, $\rho_2$ and $\varepsilon$ in Algorithm 3 are given constants for regularizations and determining error tolerance, respectively. A distinguishing feature of Algorithm 3 in contrast to Algorithm 2 is that the gains of a predictor controller in (3.49) can be utilized to obtain initial values for the iterative algorithm in Algorithm 3. This is because (3.47) can asymptotically stabilize the nominal system of (3.36) with $\boldsymbol{w}(t) \equiv \mathbf{0}_q$, as shown in Section 3.4.2. Moreover, one can always find $K$ and $X$ for (3.49) since $(A, B)$ is stabilizable. Because of this, an initial value of $\widetilde{\mathbf{K}}$ and $\widetilde{\mathbf{P}}$ for Algorithm 3 can be constructed by solving Theorem 3.1 using the gains in (3.49) without the need for a separate theorem such as Theorem 3.2. As for the value of $K$ in (3.49), we suggest that it can be obtained by solving the convex synthesis condition of state feedback stabilization for the system excluding all the delays in (3.36).

**Algorithm 3:** An inner convex approximation solution to Theorem 3.3

**begin**

Given $K \in \mathbb{R}^{p \times n}$ such that $A + BK$ is Hurwitz

**solve** Theorem 3.3 with (3.49) to produce $P$ and $Q$.

**solve** Theorem 3.3 with the previous $P$ and $Q$ to produce $K_1$, $K_2$ and $K_3$.

**let** $\widetilde{P} \leftarrow P, \ \widetilde{Q} \leftarrow Q, \ \widetilde{K}_1 \leftarrow K_1, \ \widetilde{K}_2 \leftarrow K_2, \ \widetilde{K}_3 \leftarrow K_3,$

**solve** $\min_{\mathbf{x},\mathbf{\Lambda},\mathbf{K}} \mathrm{tr}\left[\rho_1[*](\mathbf{\Lambda} - \widetilde{\mathbf{\Lambda}}) + \rho_2[*](\mathbf{K} - \widetilde{\mathbf{K}})\right]$ subject to (3.50), (3.51) and (3.58) to obtain $\mathbf{\Lambda}$ and $\mathbf{K}$

**while** $\dfrac{\left\| \begin{bmatrix} \mathbf{vec}(\mathbf{\Lambda}) \\ \mathbf{vec}\,(\mathbf{K}) \end{bmatrix} - \begin{bmatrix} \mathbf{vec}(\widetilde{\mathbf{\Lambda}}) \\ \mathbf{vec}(\widetilde{\mathbf{K}}) \end{bmatrix} \right\|_\infty}{\left\| \begin{bmatrix} \mathbf{vec}(\widetilde{\mathbf{\Lambda}}) \\ \mathbf{vec}(\widetilde{\mathbf{K}}) \end{bmatrix} \right\|_\infty + 1} \geq \varepsilon$ **do**

**update** $\widetilde{\mathbf{\Lambda}} \longleftarrow \mathbf{\Lambda}, \ \widetilde{\mathbf{K}} \longleftarrow \mathbf{K};$

**solve** $\min_{\mathbf{x},\mathbf{\Lambda},\mathbf{K}} \mathrm{tr}\left[\rho_1[*](\mathbf{\Lambda} - \widetilde{\mathbf{\Lambda}}) + \rho_2[*](\mathbf{K} - \widetilde{\mathbf{K}})\right]$ subject to (3.50), (3.51) and (3.58) to obtain $\mathbf{\Lambda}$ and $\mathbf{K}$;

**end**

**end**

## 3.5 Numerical Examples

We tested two examples numerically in this section to illustrate the strength of the proposed methods in Chapter 3. The numerical examples are coded in MATLAB© with the help of Yalmip [603] as the optimization parser. Furthermore, delay margins of the examples are examined by the spectral method in [202, 203] with code from http://cdlab.uniud.it/software#eigAM-eigTMN.

### 3.5.1 Robust stabilization of an uncertain LDDs with disspative constraints

Semidefinite programs in this subsection are solved by SDPT3 [512], with the exception of the programs corresponding to the controller in (3.64), which are solved using SeDuMi and reported in subsection 4.2.2 in [169].

Consider a system in the form of (3.1) with delay value $r = 1$ and the state space matrices

$$
\begin{aligned}
&A_1 = \begin{bmatrix} 0 & 0 \\ 0 & 0.1 \end{bmatrix}, \; A_2 = \begin{bmatrix} -1 & -1 \\ 0 & 0.9 \end{bmatrix}, \; B_1 = \begin{bmatrix} 0 \\ 1 \end{bmatrix}, \; B_2 = \begin{bmatrix} 0.1 & -0.11 \\ 0.21 & 0.1 \end{bmatrix}, \; D_1 = \begin{bmatrix} 0.2 \\ 0.3 \end{bmatrix} \\
&\widetilde{A}_3(\tau) = \begin{bmatrix} -0.4 - 0.1\mathrm{e}^{\tau}\sin(20\tau) + 0.3\mathrm{e}^{\tau}\cos(20\tau) & 1 + 0.2\mathrm{e}^{\tau}\sin(20\tau) + 0.2\mathrm{e}^{\tau}\cos(20\tau) \\ -1 + 0.01\mathrm{e}^{\tau}\sin(20\tau) - 0.2\mathrm{e}^{\tau}\cos(20\tau) & 0.4 + 0.3\mathrm{e}^{\tau}\sin(20\tau) + 0.4\mathrm{e}^{\tau}\cos(20\tau) \end{bmatrix} \\
&D_2 = \begin{bmatrix} 0.1 & 0.2 \\ 0.12 & 0.1 \end{bmatrix}, \; C_1 = \begin{bmatrix} -0.1 & 0.2 \\ 0 & 0.1 \end{bmatrix}, \; C_2 = \begin{bmatrix} -0.1 & 0 \\ 0 & 0.2 \end{bmatrix} \\
&\widetilde{C}_3(\tau) = \begin{bmatrix} 0.2\mathrm{e}^{\tau}\sin(20\tau) & 0.1 + 0.1\mathrm{e}^{\tau}\cos(20\tau) \\ 0.1\mathrm{e}^{\tau}\sin(20\tau) - 0.1\mathrm{e}^{\tau}\cos(20\tau) & -0.2 + 0.3\mathrm{e}^{\tau}\sin(20\tau) \end{bmatrix}
\end{aligned}
\tag{3.60}
$$

with uncertainties parameters $\Delta_i \in \mathbb{R}^{n\times n}$ subject to (3.4) with $\Lambda_i = F_i = \mathsf{O}_n, \forall i = 1, \cdots, 10$ and

$$
\begin{aligned}
&H_1 = \begin{bmatrix} -0.1 & -0.7 \\ -0.3 & 0.3 \end{bmatrix}, \; H_2 = \begin{bmatrix} 0.1 \\ -0.3 \end{bmatrix}, \; H_3 = \begin{bmatrix} 0.1 & 0.2 \\ 0.1 & 0.1 \end{bmatrix}, \; H_5 = \left[\begin{smallmatrix} 0.2 & 0.1 & -0.1 & 0.3 & 0.12 & -0.2 \\ 0.14 & 0.25 & 0.19 & -0.11 & -0.1 & -0.23 \end{smallmatrix}\right] \\
&H_7 = \begin{bmatrix} 0.2 & 0.2 \\ 0.21 & 0.21 \end{bmatrix}, \; H_8 = \begin{bmatrix} -0.12 & -0.14 \\ 0.01 & 0.2 \end{bmatrix}, \; H_4 = H_9 = H_{11} = \begin{bmatrix} 0.4 \\ -0.2 \end{bmatrix}, \; H_{10} = \begin{bmatrix} 0.12 & -0.14 \\ 0.01 & 0.2 \end{bmatrix} \\
&H_{12} = \left[\begin{smallmatrix} 0.2 & 0.1 & -0.1 & -0.14 & 0.1 & -0.1 \\ -0.1 & 0.3 & 0.2 & -0.1 & -0.1 & -0.1 \end{smallmatrix}\right], \; H_{14} = \begin{bmatrix} 0.22 & 0.23 \\ 0.22 & 0.23 \end{bmatrix}, G_i = \begin{bmatrix} 0.04 & 0.04 \\ 0.11 & 0.11 \end{bmatrix}, \forall i = 1, \ldots, 5 \\
&G_i = \begin{bmatrix} 0.17 & 0.17 \\ 0.14 & 0.14 \end{bmatrix}, \forall i = 6, \ldots 10, \; \Xi_1 = \begin{bmatrix} 2.3 & 1 \\ * & 2.4 \end{bmatrix}, \Xi_2 = \begin{bmatrix} 1.5 & -0.5 \\ * & 2.9 \end{bmatrix} \\
&\Xi_3 = \begin{bmatrix} 1.7 & 0.48 \\ * & 1.6 \end{bmatrix}, \; \Xi_4 = \begin{bmatrix} 2.5 & 0.51 \\ * & 2 \end{bmatrix}, \; \Xi_5 = \begin{bmatrix} 1.7 & 0.44 \\ * & 1.7 \end{bmatrix}, \; \Xi_6 = \begin{bmatrix} 1.6 & 0.15 \\ * & 1.4 \end{bmatrix}, \; \Xi_7 = \begin{bmatrix} 3.37 & -1.1 \\ * & 1.8 \end{bmatrix} \\
&\Xi_8 = \begin{bmatrix} 1.54 & 0.13 \\ * & 1.34 \end{bmatrix}, \; \Xi_9 = \begin{bmatrix} 2.7 & -0.65 \\ * & 1.87 \end{bmatrix}, \; \Xi_{10} = \begin{bmatrix} 1.7 & 0.44 \\ * & 1.7 \end{bmatrix}, \; \Gamma_i = -I_2, \forall i = 1 \cdots 10.
\end{aligned}
\tag{3.61}
$$

Note that (3.60) is identical to (2.53) when uncertainties are not taken into account. Thus the methods in [463] and [503] still cannot handle (3.60) as argued in the previous chapter, even without considering the presence of uncertainties in (3.60). Moreover, it is shown in Chapter 2 by means of the spectrum method that (3.60) is unstable for $0 \leq r \leq 10$ without considering the presence of uncertainties. For the controller objective, again we choose to minimize the value of $\mathcal{L}^2$ attenuation factor $\gamma$ which corresponds to $J_3 = -J_1 = \gamma I_2$, $\widetilde{J} = I_2$, $J_2 = O_2$ in (2.10). Now consider the parameters

$$\boldsymbol{f}(\tau) = \begin{bmatrix} 1 \\ 10e^{\tau}\sin(20\tau) \\ 10e^{\tau}\cos(20\tau) \end{bmatrix} \otimes I_2, \quad M = \begin{bmatrix} 0 & 0 & 0 \\ 0 & 1 & 20 \\ 0 & -20 & 1 \end{bmatrix} \tag{3.62}$$

and

$$A_3 = 0.1 \begin{bmatrix} -4 & 10 & -0.1 & 0.2 & 0.3 & 0.2 \\ -10 & 4 & 0.01 & 0.3 & -0.2 & 0.4 \end{bmatrix}, \quad C_3 = 0.1 \begin{bmatrix} 0 & 1 & 0.2 & 0 & 0 & 0.1 \\ 0 & -2 & 0.1 & 0.3 & -0.1 & 0 \end{bmatrix} \tag{3.63}$$

which are the same as in (2.55a) with $d = 3, \quad n = m = q = 2$.

Now we apply Theorem 3.1 to the system with the parameters in (3.60)–(3.63). It follows that the corresponding uncertain system is robustly stabilized by the controller

$$\boldsymbol{u}(t) = \begin{bmatrix} 3.2847 & -16.7739 \end{bmatrix} \boldsymbol{x}(t) \tag{3.64}$$

for any uncertainties in the corresponding $\mathcal{D}$ with $\min \gamma = 0.62$. To reduce the potential conservatism of the value of $\min \gamma$ calculated by Theorem 3.2, we apply Theorem 3.1 to the previous resulting CLS with controller gain $K = \begin{bmatrix} 3.2847 & -16.7739 \end{bmatrix}$. It shows that $K = \begin{bmatrix} 3.2847 & -16.7739 \end{bmatrix}$ can achieve $\min \gamma = 0.45941$. Next, we apply the spectrum method again to the resulting CLS without considering the uncertainties. This produces $-0.1773 < 0$ as the real part of the rightmost characteristic root pair, indicating that the nominal resulting CLS is stable.

We can now apply Algorithm 2 to compute controller gains with improved performance. The results produced by Algorithm 2 with $\rho_1 = \rho_2 = 10^{-7}$ and $\varepsilon = 10^{-12}$ are cataloged in Table 3.1, where Number of Iterations represents the number of iterations executed by the while loop in Algorithm 2. Furthermore, denotes the spectral abscissas of the sets containing the characteristic roots of the nominal CLSs (without uncertainties) listed in Table 3.1. These values were computed using the method described in [202, 203].

| Controller gains $K$ | $\begin{bmatrix} 3.6897 \\ -19.0393 \end{bmatrix}^\top$ | $\begin{bmatrix} 4.0571 \\ -21.1008 \end{bmatrix}^\top$ | $\begin{bmatrix} 4.4230 \\ -23.1575 \end{bmatrix}^\top$ | $\begin{bmatrix} 4.7878 \\ -25.2108 \end{bmatrix}^\top$ |
|---|---|---|---|---|
| $\min \gamma$ | 0.459357 | 0.459321 | 0.459292 | 0.459269 |
| Number of Iterations | 10 | 20 | 30 | 40 |
| SPA | $-0.177$ | $-0.1768$ | $-0.1766$ | $-0.1764$ |

**Table 3.1:** $\min \gamma$ produced by different iterations

The results in Table 3.1 demonstrate that more iterations lead to better $\min \gamma$ value at the cost of increased numerical complexity. Note that by employing Theorem 3.1 with the numerical $K$ in Table 3.1, the same values of $\min \gamma$ can be obtained by solving the corresponding convex SDP constraints. This indicates that the numerical results produced by the iterative Algorithm 2 in terms of $\min \gamma$ are reliable.

### 3.5.2 Resilient dynamical state feedback design for an uncertain linear system with an input delay

We apply Mosek [513] as the numerical solver for SDP in this subsection.

Consider the uncertain open-loop system (3.36) with $r = 3$ and the state space matrices

$$\begin{aligned} A &= \begin{bmatrix} -1 & 1 \\ 0 & 0.1 \end{bmatrix}, \ B = \begin{bmatrix} 0 \\ 1 \end{bmatrix}, \ D_1 = \begin{bmatrix} 0.1 \\ -0.1 \end{bmatrix}, \ C_1 = \begin{bmatrix} -0.3 & 0.4 & 0.1 \\ -0.3 & 0.1 & -0.1 \end{bmatrix}, \ C_2 = \begin{bmatrix} 0 & 0.2 & 0 \\ -0.2 & 0.1 & 0 \end{bmatrix}, \\ \widetilde{C}_3(\tau) &= \begin{bmatrix} 0.2 + 0.1\mathrm{e}^{\tau} & 0.1 & 0.12\mathrm{e}^{3\tau} \\ -0.2 & 0.3 + 0.14\mathrm{e}^{2\tau} & 0.11\mathrm{e}^{3\tau} \end{bmatrix}, \ D_2 = 0.12, \ \ D_3 = \begin{bmatrix} 0.14 \\ 0.1 \end{bmatrix} \end{aligned} \tag{3.65}$$

with $\Delta_i \in \mathbb{R}^{\nu \times \nu}$ subject to (3.45) with $\Lambda_i = F_i = \mathsf{O}_\nu, \forall i = 1, \ldots, 8$ and

$$H_1 = \begin{bmatrix} -0.1 & -0.7 & 0 \\ -0.3 & 0.3 & 0 \\ 0 & 0 & 0.1 \end{bmatrix}, \ H_2 = \begin{bmatrix} -0.1 & 0 & 0 \\ -0.3 & 0.3 & 0 \\ 0 & 0 & 0.1 \end{bmatrix}, \ H_3 = \begin{bmatrix} 0 & 0 & 0 & 0 & 0 & 0 & \mathbf{0}_9^\top \\ 0.14 & 0.25 & 0.19 & -0.11 & -0.1 & -0.23 & \mathbf{0}_9^\top \\ 0.1 & 0.14 & 0.07 & 0.12 & 0.17 & 0.15 & \mathbf{0}_9^\top \end{bmatrix} \tag{3.66}$$

$$H_4 = \begin{bmatrix} 0.01 \\ 0.02 \\ 0.01 \end{bmatrix}, \ H_5 = \begin{bmatrix} 0.12 & -0.14 & 0.13 \\ 0.01 & 0.2 & 0.05 \\ 0.12 & 0.18 & 0.15 \end{bmatrix}, \ H_6 = \begin{bmatrix} 0.12 & -0.14 & 0.11 \\ 0.01 & 0.2 & 0.15 \\ 0.11 & 0.12 & 0.011 \end{bmatrix} \tag{3.67}$$

$$H_7 = \begin{bmatrix} 0.2 & 0.1 & -0.1 & -0.14 & 0.1 & -0.1 & \mathbf{0}_9^\top \\ -0.1 & 0.3 & 0.2 & 0.1 & -0.1 & 0.1 & \mathbf{0}_9^\top \\ 0.2 & 0.1 & -0.1 & 0.3 & 0.12 & -0.2 & \mathbf{0}_9^\top \end{bmatrix}, \ H_8 = \begin{bmatrix} 0.2 \\ 0.1 \\ 0.1 \end{bmatrix} \tag{3.68}$$

$$G_1 = \begin{bmatrix} 0.04 & 0.04 & 0 \\ 0.11 & 0.11 & 0.1 \\ 0.02 & 0.01 & 0.03 \end{bmatrix}, \ G_2 = \begin{bmatrix} 0 & 0 & 0.04 \\ 0.11 & 0.11 & 0.12 \\ 0.05 & 0.07 & 0.03 \end{bmatrix}, \ G_3 = \begin{bmatrix} 0 & 0 & 0 \\ 0.11 & 0.11 & 0.1 \\ 0.02 & 0.01 & 0.03 \end{bmatrix} \tag{3.69}$$

$$G_4 = \begin{bmatrix} 0.04 & 0.04 & 0.05 \\ 0.11 & 0.11 & 0.1 \\ 0.02 & 0.01 & 0.03 \end{bmatrix}, \ G_i = \begin{bmatrix} 0.11 & 0.14 & 0.12 \\ 0.1 & 0.11 & 0.12 \end{bmatrix}, \forall i = 6, \ldots, 10, \ \Xi_1 = \begin{bmatrix} 2.3 & 1 & 0 \\ 1 & 2.4 & 0 \\ 0 & 0 & 1 \end{bmatrix} \tag{3.70}$$

$$\Xi_2 = \begin{bmatrix} 1.5 & -0.5 & 0.1 \\ -0.5 & 2.9 & 0.1 \\ 0.1 & 0.1 & 0.11 \end{bmatrix}, \ \Xi_3 = \begin{bmatrix} 1.7 & 0.48 & 0 \\ 0.48 & 1.6 & 0 \\ 0 & 0 & 1.2 \end{bmatrix}, \ \Xi_4 = \begin{bmatrix} 2.5 & 0.51 & 0.2 \\ 0.51 & 2 & 0.2 \\ 0.2 & 0.2 & 1.5 \end{bmatrix} \tag{3.71}$$

$$\Xi_5 = \begin{bmatrix} 1.7 & 0.44 & 0 \\ 0.44 & 1.7 & 0 \\ 0 & 0 & 1.6 \end{bmatrix}, \ \Xi_6 = \begin{bmatrix} 1.6 & 0.15 & 0.12 \\ 0.15 & 1.4 & 0.2 \\ 0.12 & 0.2 & 1.8 \end{bmatrix}, \ \Xi_7 = \begin{bmatrix} 3.37 & -1.1 & 0 \\ -1.1 & 1.8 & 0 \\ 0 & 0 & 2.2 \end{bmatrix} \tag{3.72}$$

$$\Xi_8 = \begin{bmatrix} 1.54 & 0.13 & 0 \\ 0.13 & 1.34 & 0 \\ 0 & 0 & 1.7 \end{bmatrix}, \ \Gamma_i = \begin{bmatrix} -1.75 & -0.58 & 0 \\ -0.58 & -1.75 & 0 \\ 0 & 0 & -1.75 \end{bmatrix}, \ i = 1, \ldots, 8. \tag{3.73}$$

Now consider the functions in

$$\mathrm{e}^{-A\tau} B = \mathbf{Col}\left[ (10/11)\, \mathrm{e}^{-0.1\tau} - (10/11)\, \mathrm{e}^{\tau}, \ \mathrm{e}^{-0.1\tau} \right] \tag{3.74}$$

and the functions inside of $\widetilde{C}_3(\tau)$, we apply

$$\boldsymbol{f}(\tau) = \mathbf{Col}\left(1, \mathrm{e}^{\tau}, \mathrm{e}^{2\tau}, \mathrm{e}^{3\tau}, \mathrm{e}^{-0.1\tau}\right), \ \frac{\mathrm{d}\boldsymbol{f}(\tau)}{\mathrm{d}\tau} = M \boldsymbol{f}(\tau) \tag{3.75}$$

with $M = 0 \oplus 1 \oplus 2 \oplus 3 \oplus (-0.1)$ as the basis function to denote all the distributed terms in (3.74) and (3.65). Now let $K = \begin{bmatrix} -0.5249 & -0.4173 \end{bmatrix}$ and $X = -0.1$ in (3.49) with (3.75), where $A + BK$ and $X$ are Hurwitz. Then we have

$$\begin{aligned} C_3 &= \begin{bmatrix} 0.2 & 0.1 & 0 & 0.1 & 0 & 0 & 0 & 0 & 0 & 0 & 0 & 0.12 & 0 & 0 & 0 \\ -0.2 & 0.3 & 0 & 0 & 0 & 0 & 0 & 0.14 & 0 & 0 & 0 & 0.11 & 0 & 0 & 0 \end{bmatrix} \left( \sqrt{\mathsf{F}^{-1}} \otimes I_3 \right), \\ \widehat{K} &= \begin{bmatrix} \mathbf{0}_5^\top & -0.4295 & \mathbf{0}_8^\top & -0.1789 \end{bmatrix} \left( \sqrt{\mathsf{F}^{-1}} \otimes I_3 \right), \end{aligned} \tag{3.76}$$

where $\mathsf{F}^{-1} = \int_{-r}^{0} \boldsymbol{f}(\tau)\boldsymbol{f}^{\top}(\tau)\mathrm{d}\tau$ with $\boldsymbol{f}(\tau)$ in (3.75), and

$$K_1 = \begin{bmatrix} 0.0235 & -0.2629 & -0.5173 \end{bmatrix}, \; K_2 = \mathbf{0}_3^{\top}, \tag{3.77}$$

where $K_1$, $K_2$ and $\widehat{K}$ can be used to initiate the iterative algorithm in Algorithm 3.

Given the parameters in (3.65)–(3.73), applying Theorem 3.3 to (3.41) with the controller gains in (3.76) and (3.77) yields a feasible solution and shows the controller can achieve $\min \gamma = 0.99295$. Following the procedures in Algorithm 3, the resulting $P$ and $Q$ produced by the previous optimization problem concerning $K_1$, $K_2$ and $\widehat{K}$ in (3.76)–(3.77) are substituted into Theorem 3.3, which can now be solved as a convex optimization program to compute a new $\mathbf{K}$. By doing so, we obtain

$$\begin{aligned} K_1 &= \begin{bmatrix} -0.0250 & -0.1808 & -0.5636 \end{bmatrix}, \; K_2 = \begin{bmatrix} -0.0098 & 0.0067 & -0.0394 \end{bmatrix}, \\ K_3 &= \Big[ 0.0294 \;\; -0.0610 \;\; -0.4612 \;\; 0.0455 \;\; -0.1037 \;\; -0.2287 \;\; -0.0011 \;\; -0.0932 \cdots \\ &\qquad -0.1429 \;\; -0.0181 \;\; -0.0767 \;\; -0.0931 \;\; -0.0071 \;\; -0.0328 \;\; -0.5428 \Big], \end{aligned} \tag{3.78}$$

which can achieve the performance $\gamma = 0.92315$. To proceed, use (3.78) with the values of the associated $P$ and $Q$ for $\grave{\mathbf{Y}}$ and $\grave{\mathbf{K}}$ to follow the steps in Algorithm 3. The results produced by Algorithm 3 with $\rho_1 = 0.01$, $\rho_2 = 0.01$ and $\varepsilon = 10^{-12}$ are given in Table 3.2, where Number of Iterations represents the number of iterations executed by the while loop in Algorithm 3, and SPA denotes the spectral abscissas of the sets containing the characteristic roots of the resulting nominal CLSs (without uncertainties). Note that the values of SPA are calculated via the methods in [202, 203].

| $\min \gamma$ | 0.88264 | 0.8562911 | 0.831113 | 0.807242 |
|---|---|---|---|---|
| Number of Iterations | 100 | 200 | 300 | 400 |
| SPA | $-0.2749$ | $-0.2635$ | $-0.2608$ | $-0.2874$ |

**Table 3.2:** $\min \gamma$ produced by Algorithm 3 with different numbers of iterations

Moreover, the controller gains corresponding to the results in Table 3.2 are presented as follows.

$$\begin{aligned} \mathrm{NoI} = 100 : K_1 &= \begin{bmatrix} -0.0288 & -0.206 & -0.6271 \end{bmatrix}, \; K_2 = \begin{bmatrix} -0.0132 & 0.0124 & -0.0411 \end{bmatrix} \\ K_3 &= \Big[ 0.016 \;\; -0.0383 \;\; -0.4609 \;\; 0.0395 \;\; -0.11 \;\; -0.2298 \;\; -0.0175 \;\; -0.1108 \cdots \\ &\qquad -0.1579 \;\; -0.0225 \;\; -0.0872 \;\; -0.1107 \;\; -0.0167 \;\; -0.0236 \;\; -0.5467 \Big] \end{aligned} \tag{3.79}$$

$$\begin{aligned} \mathrm{NoI} = 200 : K_1 &= \begin{bmatrix} -0.0274 & -0.2203 & -0.6585 \end{bmatrix}, \; K_2 = \begin{bmatrix} -0.0109 & 0.012 & -0.0351 \end{bmatrix} \\ K_3 &= \Big[ 0.0016 \;\; -0.0174 \;\; -0.4539 \;\; 0.0288 \;\; -0.1194 \;\; -0.2313 \;\; -0.0235 \;\; -0.1306 \cdots \\ &\qquad -0.1657 \;\; -0.0283 \;\; -0.1051 \;\; -0.1186 \;\; -0.0306 \;\; -0.0083 \;\; -0.5541 \Big] \end{aligned} \tag{3.80}$$

$$\begin{aligned} \mathrm{NoI} = 300 : K_1 &= \begin{bmatrix} -0.0294 & -0.2350 & -0.6847 \end{bmatrix}, \; K_2 = \begin{bmatrix} -0.0072 & 0.0069 & -0.0289 \end{bmatrix} \\ K_3 &= \Big[ -0.0015 \;\; -0.0006 \;\; -0.4439 \;\; 0.0271 \;\; -0.1219 \;\; -0.2308 \;\; -0.0286 \;\; -0.1433 \cdots \\ &\qquad -0.1741 \;\; -0.0325 \;\; -0.1206 \;\; -0.1217 \;\; -0.0381 \;\; -0.0019 \;\; -0.562 \Big] \end{aligned} \tag{3.81}$$

$$\begin{aligned} \mathrm{NoI} = 400 : K_1 &= \begin{bmatrix} -0.0479 & -0.2484 & -0.6754 \end{bmatrix}, \; K_2 = \begin{bmatrix} 0.0089 & 0.0057 & -0.0155 \end{bmatrix} \\ K_3 &= \Big[ -0.0055 \;\; 0.0258 \;\; -0.4477 \;\; 0.0205 \;\; -0.1192 \;\; -0.2471 \;\; -0.0431 \;\; -0.1603 \cdots \\ &\qquad -0.1909 \;\; -0.0404 \;\; -0.1226 \;\; -0.1304 \;\; -0.0532 \;\; 0.0069 \;\; -0.569 \Big] \end{aligned} \tag{3.82}$$

Note that all the controller gains in the steps of Algorithm 3 can stabilize (3.36) with (3.65).

Clearly, the results in Table 3.2 demonstrate that more iterations lead to smaller values of $\min \gamma$ at the expense of greater numbers of iterations. In addition, it also shows an example of using Algorithm 3 to compute controller gains with better performance. It is important to note that all of the resulting CLSs mentioned above are of the retarded type and thus satisfy the properties emphasized in Remark 3.12.

# Chapter 4

# Two General Classes of Integral Inequalities Including Weight Functions

## 4.1 Introduction

The derivations of the synthesis conditions in Chapter 2 clearly demonstrate the importance of utilizing capable IIs for the construction of synthesis conditions at key steps such as (2.30) and (2.26), when using the LKF approach. Without the IIs in (2.30) and (2.26), the resulting synthesis condition (2.1) would become exceptionally conservative, and could even be infeasible in all situations. The IIs in (2.30), (2.26) provide the necessary lower bounds (LBs) that match the structures of the "large" quadratic forms in $\mathsf{v}(\cdot)$ of (2.20) and $\dot{\mathsf{v}}(\cdot)$ of (2.25), which are also strategically located in the diagonal blocks in (2.21a) and (2.21b). As researchers have come to realize that IIs are instrumental in improving the effectiveness of the LKFs method, tremendous efforts have been made to reduce the induced conservatisms by developing novel IIs with tighter LBs. Broadly speaking, two classes of IIs have been proposed with different characteristics. The first class can be referred to as the Bessel type IIs [169, 469, 572, 621], as they share similarities with the structure of the Legendre-Bessel inequality initially proposed in [467]. These IIs contain no additional matrices other than the original matrix in the quadratic term. In contrast, free matrix type inequalities with extra matrices have been proposed in [622–625], motivated by their applications to the stability analysis of time-varying delay systems. Based on our literature review, it becomes evident that the applicable structures of LKFs that can be employed are intrinsically determined by the IIs available to use, which in turn dictates the level of conservatism inherent in the resulting stability/synthesis conditions. Thus, it is imperative to develop IIs with optimal LBs, as these bounds serve a crucial role in the effective application of the LKF approach.

In this chapter, we investigate two general classes of IIs and analyze the properties of their LBs in detail. The proposed IIs can be applied to a wide range of different situations, including but not limited to the analysis of delay (time-varying) related systems [424, 426–429, 431, 433, 434, 436–438, 440–442, 444, 445, 447–449], PDE-related systems [45, 536–538, 556, 626], and sampled-data systems [550, 627–630]. The first class of IIs is proposed in Section 4.2, which can generalize numerous existing Bessel-type inequalities in [169, 467, 469, 572, 621, 631–633]. The second class, formulated based on the application of [634, Lemma 4.1], is presented in Section 4.3, where we prove an important relationship concerning the equivalence of the LBs between the first and second classes. Furthermore, another inequality, which is a special case of the second class, is also derived in Section 4.4 generalizing the existing inequalities in [623–625]. We then establish a significant conclusion concerning the LBs between all the proposed IIs, through which equivalence between many existing inequalities can be confirmed. To demonstrate some applications of the proposed IIs, we apply them in Section 4.5 to derive a stability condition for a linear CDDS [22] with a DD by constructing a parameterized complete LKF. We show that equivalent stability conditions, whose feasibility is invariant with respect to some parameters in the LKF, can be obtained by making use of the proposed IIs. The core contributions of this chapter are to establish a general framework for the IIs utilized in the analysis of TDSs and other IDSs. By clearly demonstrating the optimality of their LBs and their mutual relationships, the proposed IIs offer great potential in addressing the analysis of a variety of different systems that are related to the structures of TDSs or IDSs.

## 4.2 The Inequality of The First Class

The first class of IIs is derived in this section, where we also mathematically demonstrate its generality in comparison to existing results in the literature.

We first define the weighted Lebesgue function space

$$\mathscr{L}^2_{\varpi}(\mathcal{K};\mathbb{R}^d) := \left\{\phi(\cdot)\in \mathcal{M}_{\mathcal{L}(\mathcal{X})/\mathcal{B}(\mathbb{R}^d)}(\mathcal{X};\mathbb{R}^d) : \|\phi(\cdot)\|_{\varpi} < \infty\right\} \tag{4.1}$$

with $d\in\mathbb{N}$ and the seminorm $\|\phi(\cdot)\|_{\varpi} := \int_{\mathcal{K}} \varpi(\tau)\phi^\top(\tau)\phi(\tau)\mathrm{d}\tau$ where $\varpi(\cdot)\in \mathcal{M}_{\mathcal{L}(\mathcal{K})/\mathcal{B}(\mathbb{R}_{\geq 0})}(\mathcal{K};\mathbb{R}_{\geq 0})$ and the function $\varpi(\cdot)$ has only countably infinite or finite number of zero values. Furthermore, $\mathcal{K}\subseteq \mathbb{R}\cup\{\pm\infty\}$ and the Lebesgue measure [571] of $\mathcal{K}$ is nonzero.

**Theorem 4.1.** *Given $\varpi(\cdot)$ in (4.1) and $U\in\mathbb{S}^n_{\succ 0}$ and $\boldsymbol{f}(\cdot)\in\mathscr{L}^2_{\varpi}(\mathcal{K};\mathbb{R}^d)$ which satisfies*

$$\int_{\mathcal{K}} \varpi(\tau)\boldsymbol{f}(\tau)\boldsymbol{f}^\top(\tau)\mathrm{d}\tau \succ 0, \tag{4.2}$$

*then for all $\boldsymbol{x}(\cdot)\in\mathscr{L}^2_{\varpi}(\mathcal{K};\mathbb{R}^n)$ we have*

$$\int_{\mathcal{K}} \varpi(\tau)\boldsymbol{x}^\top(\tau)U\boldsymbol{x}(\tau)\mathrm{d}\tau \geq \int_{\mathcal{K}}\varpi(\tau)\boldsymbol{x}^\top(\tau)F^\top(\tau)\mathrm{d}\tau\,(\mathsf{F}\otimes U)\int_{\mathcal{K}}\varpi(\tau)F(\tau)\boldsymbol{x}(\tau)\mathrm{d}\tau \tag{4.3}$$

*where $F(\tau) := \boldsymbol{f}(\tau)\otimes I_n$ and $\mathsf{F}^{-1} = \int_{\mathcal{K}}\varpi(\tau)\boldsymbol{f}(\tau)\boldsymbol{f}^\top(\tau)\mathrm{d}\tau$.*

*Proof.* See Appendix B for details. ∎

Inequality (4.3) holds for any $\boldsymbol{f}(\cdot)\in\mathscr{L}^2_{\varpi}(\mathcal{K};\mathbb{R}^d)$ satisfying (4.2) with a given $\varpi(\cdot)$. Note that the constraint (4.2) in Theorem 4.1 indicates the functions in $\boldsymbol{f}(\cdot)$ are linearly independent (see the [566, Theorem 7.2.10]) in a Lebesgue sense. Thus we can select general functions for $\boldsymbol{f}(\cdot)$ when using (4.3). This includes the situation of $\boldsymbol{f}(\tau)$ containing orthogonal functions, elementary functions or the other type of function as long as they are linearly independent in a Lebesgue sense. As a result, the structure of (4.3) is by far the most general Bessel type inequality in terms of the applicable integral kernels of the lower quadratic bound. Finally, the inequality in (4.3) still holds if $U\succeq 0$. Note that $U\succ 0$ in (4.3) is taken as the prerequisite of Theorem 4.1 to ensure that relations can be established between (4.3) and the IIs proposed in later sections.

The generality of (4.3) is demonstrated with mathematical details as follows. To do so, let us first provide the standard expression of Jacobi polynomials (see [468, 22.3.2])

$$j_d^{\alpha,\beta}(\tau)_{-1}^{1} := \frac{\gamma(d+1+\alpha)}{d!\,\gamma(d+1+\alpha+\beta)}\sum_{k=0}^{d}\binom{d}{k}\frac{\gamma(d+k+1+\alpha+\beta)}{\gamma(k+1+\alpha)}\left(\frac{\tau-1}{2}\right)^k,\quad \tau\in[-1,1] \tag{4.4}$$

over $[-1,1]$, where $\gamma(\cdot)$ stands for the standard gamma function with $d\in\mathbb{N}_0$ and $\alpha>-1,\beta>-1$. The polynomials represented in (4.6) follow the orthogonal property (see [468, 22.2.1])

$$\int_{-1}^{1}(1-\tau)^{\alpha}(\tau+1)^{\beta}\boldsymbol{j}_d^{\alpha,\beta}(\tau)_{-1}^{1}\left[\boldsymbol{j}_d^{\alpha,\beta}(\tau)_{-1}^{1}\right]^{\top}\mathsf{d}\tau=\bigoplus_{k=0}^{d}\frac{2^{\alpha+\beta+1}\gamma(k+\alpha+1)\gamma(k+\beta+1)}{k!(2k+\alpha+\beta+1)\gamma(k+\alpha+\beta+1)} \tag{4.5}$$

with $\boldsymbol{j}_d^{\alpha,\beta}(\tau)_{-1}^{1}=\mathbf{Col}_{i=0}^{d}\, j_i^{\alpha,\beta}(\tau)_{-1}^{1}$, where the polynomials are orthogonal with respect to $(1-\tau)^{\alpha}(\tau+1)^{\beta}$ over $[-1,1]$. However, it is preferable to derive a general expression for Jacobi polynomials defined over $[a,b]$ with $b>a$. Specifically, consider the affine transformation $\frac{2\tau-a-b}{b-a}\to\tau$ where the affine function $\frac{2\tau-a-b}{b-a}$ satisfies $-1\leq\frac{2\tau-a-b}{b-a}\leq 1$ for $\tau\in[a,b]$ with $b>a$. The shift-scaled Jacobi polynomials $j_d^{\alpha,\beta}\left(\frac{2\tau-a-b}{b-a}\right)_{-1}^{1}$ is denoted as

$$j_d^{\alpha,\beta}(\tau)_a^b:=j_d^{\alpha,\beta}\left(\frac{2\tau-a-b}{b-a}\right)_{-1}^{1}=\frac{\gamma(d+1+\alpha)}{d!\gamma(d+1+\alpha+\beta)}\sum_{k=0}^{d}\binom{d}{k}\frac{\gamma(d+k+1+\alpha+\beta)}{\gamma(k+1+\alpha)}\left(\frac{\tau-b}{b-a}\right)^{k} \tag{4.6}$$

with $\tau\in[a,b]$. Now by using the affine transformation $\frac{2\tau-a-b}{b-a}\to\tau$ to (4.5), it yields

$$\begin{aligned}\int_{a}^{b}\left(\frac{-2\tau+2b}{b-a}\right)^{\alpha}\left(\frac{2\tau-2a}{b-a}\right)^{\beta}\boldsymbol{j}_d^{\alpha,\beta}\left(\frac{2\tau-a-b}{b-a}\right)_{-1}^{1}\left[\boldsymbol{j}_d^{\alpha,\beta}\left(\frac{2\tau-a-b}{b-a}\right)_{-1}^{1}\right]^{\top}\mathsf{d}\left(\frac{2\tau-a-b}{b-a}\right)\\ =\frac{2^{\alpha+\beta+1}}{(b-a)^{\alpha+\beta+1}}\int_{a}^{b}(b-\tau)^{\alpha}(\tau-a)^{\beta}\boldsymbol{j}_d^{\alpha,\beta}(\tau)_a^b\left[\boldsymbol{j}_d^{\alpha,\beta}(\tau)_a^b\right]^{\top}\mathsf{d}\tau\\ =\bigoplus_{k=0}^{d}\frac{2^{\alpha+\beta+1}\gamma(k+\alpha+1)\gamma(k+\beta+1)}{k!(2k+\alpha+\beta+1)\gamma(k+\alpha+\beta+1)}\end{aligned} \tag{4.7}$$

where $\boldsymbol{j}_d^{\alpha,\beta}(\tau)_a^b:=\mathbf{Col}_{i=0}^{d}\, j_k^{\alpha,\beta}(\tau)_a^b$ for $\tau\in[a,b]$. Moreover, eq.(4.7) can be rewritten as

$$\int_{a}^{b}(b-\tau)^{\alpha}(\tau-a)^{\beta}\boldsymbol{j}_d^{\alpha,\beta}(\tau)_a^b\left[\boldsymbol{j}_d^{\alpha,\beta}(\tau)_a^b\right]^{\top}\mathsf{d}\tau=\bigoplus_{k=0}^{d}\frac{(b-a)^{\alpha+\beta+1}\gamma(k+\alpha+1)\gamma(k+\beta+1)}{k!(2k+\alpha+\beta+1)\gamma(k+\alpha+\beta+1)}$$

which now is the expression for the orthogonality of (4.4) with respect to $\varpi(\tau)=(b-\tau)^{\alpha}(\tau-a)^{\beta}$. In fact, for $\alpha=\beta=0$, (4.6) becomes Legendre polynomials (see [468, 22.2.10])

$$\ell_d(\tau):=\sum_{k=0}^{d}\binom{d}{k}\binom{d+k}{k}\left(\frac{\tau-b}{b-a}\right)^{k} \tag{4.8}$$

with $d\in\mathbb{N}_0$ and $\tau\in[a,b]$, which satisfies $\int_a^b\boldsymbol{\ell}_d(\tau)\boldsymbol{\ell}_d^{\top}(\tau)\mathsf{d}\tau=\bigoplus_{k=0}^{d}\frac{b-a}{(2k+1)}$.

**Remark 4.1.** Let $\alpha=0$ and $\beta\in\mathbb{N}_0$, then (4.6) becomes identical to the orthogonal hyper-geometric polynomials defined in [621, eq.(13) and (14)]. By using the Cauchy formula for repeated integrations[1], it becomes evident that (4.6) with $\alpha=0,\beta\in\mathbb{N}_0$ and $\alpha\in\mathbb{N}_0,\beta=0$ are also equivalent to the polynomials defined in [632, eq.(3) and (4)], respectively. This conclusion is also corroborated by the structure in [635, eq.(10)].

Having presented the definitions of the Jacobi and Legendre polynomials, we provide a list of existing IIs in Table 4.1 as the special cases of (4.3) with appropriate $\mathcal{K}$, $\varpi(\tau)$, $\boldsymbol{f}(\tau)$ and $\boldsymbol{x}(\cdot)$. Note that some of the results in Table 4.1 require the application of the Cauchy formula for repeated integration.

[1]See [621, eq.(5),(6) and eq.(25),(26)] and the [625, Lemma 1] for concrete examples

| (4.3) | $\mathcal{K}$ | $\varpi(\tau)$ | $\boldsymbol{f}(\tau)$ | $\boldsymbol{x}(\cdot)$ |
|---|---|---|---|---|
| [469, eq.(5)] | $[-r, 0]$ | $1$ | $\boldsymbol{j}_d^{0,0}(\tau)_{-r}^{0}$ | $\boldsymbol{x}(\cdot)$ |
| [469, eq.(6)] | $[-r, 0]$ | $1$ | $\boldsymbol{j}_d^{0,0}(\tau)_{-r}^{0}$ | $\dot{\boldsymbol{x}}(\cdot)$ |
| [573, eq.(5)] | $[a, b]$ | $1$ | $\begin{bmatrix} 1 \\ \boldsymbol{p}(\tau) \end{bmatrix}$ | $\boldsymbol{x}(\cdot)$ |
| [621, eq.(27)] | $[a, b]$ | $\frac{(\tau-a)^p}{(a-b)^p}$ | $\boldsymbol{j}_d^{0,p}(\tau)_a^b$ | $\boldsymbol{x}(\cdot)$ |
| [621, eq.(34)] | $[a, b]$ | $\frac{(\tau-a)^p}{(a-b)^p}$ | $\boldsymbol{j}_d^{0,p}(\tau)_a^b$ | $\dot{\boldsymbol{x}}(\cdot)$ |
| [632, eq.(1)] | $[a, b]$ | $(\tau-a)^{m-1}$ | $\boldsymbol{j}_d^{0,m-1}(\tau)_a^b$ | $\boldsymbol{x}(\cdot)$ |
| [632, eq.(2)] | $[a, b]$ | $(b-\tau)^{m-1}$ | $\boldsymbol{j}_d^{m-1,0}(\tau)_a^b$ | $\boldsymbol{x}(\cdot)$ |
| [631, eq.(2)] | $[a, b]$ | $(\tau-a)^k$ | $\boldsymbol{p}_k(\tau)$ | $\boldsymbol{x}(\cdot)$ |
| (2.12) | $\mathcal{K}$ | $1$ | $\boldsymbol{g}(\tau)$ | $\boldsymbol{x}(\cdot)$ |
| [633, eq.(9)] | $[0, +\infty]$ | $K(\tau)$ | $\begin{bmatrix} 1 \\ g(\tau) \end{bmatrix}$ | $\boldsymbol{x}(\cdot)$ |

**Table 4.1:** List of IIs encompassed by (4.3)

For $\boldsymbol{p}_k(\cdot)$ in Table 4.1, it is defined as

$$\boldsymbol{p}_k(\cdot) \in \left\{ \boldsymbol{f}(\cdot) = \operatorname*{\mathsf{Col}}_{i=1}^{d} f_i(\cdot) \in \mathscr{L}^2_{(\tau-a)^k}\left([a,b];\mathbb{R}^d\right) : \int_{\mathcal{K}} (\tau-a)^k \boldsymbol{f}(\tau)\boldsymbol{f}^\top(\tau)\mathsf{d}\tau = \bigoplus_{i=1}^{d} \int_a^b (\tau-a)^k f_i^2(\tau)\mathsf{d}\tau \right\}$$

with $k \in \mathbb{N}_0$, which means all the functions in $\boldsymbol{p}_k(\cdot)$ are orthogonal functions with respect to the corresponding weight functions. Furthermore, $\boldsymbol{p}(\cdot)$ in Table 4.1 is defined as

$$\boldsymbol{p}(\cdot) \in \left\{ \boldsymbol{f}(\cdot) = \operatorname*{\mathsf{Col}}_{i=1}^{d-1} f_i(\cdot) \in \mathscr{L}^2\left([a,b];\mathbb{R}^{d-1}\right) : \int_a^b \boldsymbol{f}(\tau)\boldsymbol{f}^\top(\tau)\mathsf{d}\tau = \bigoplus_{i=1}^{d-1} \int_a^b f_i^2(\tau)\mathsf{d}\tau \quad \& \quad \int_a^b \boldsymbol{f}(\tau)\mathsf{d}\tau = \mathbf{0}_d \right\}$$

where it contains the auxiliary functions generated by the process in [573, Lemma 1]. The terms $\varpi(\tau) = K(\tau)$ and $\boldsymbol{f}(\tau) = \mathsf{Col}\,[1, g(\tau)]$ are in line with the definitions in [633, Theorem 1]. Finally, the inclusions by (4.3) in Table 4.1 concerning the inequalities in [621] are demonstrated as follows.

Let $\alpha = 0$, $\beta = p \in \mathbb{N}_0$ with $\varpi(\tau) = (\tau-a)^p$ and $\boldsymbol{f}(\tau) = \boldsymbol{j}_d^{0,p}(\tau)_a^b$, then (4.3) becomes

$$\begin{aligned} \int_a^b (\tau-a)^p \boldsymbol{x}^\top(\tau) U \boldsymbol{x}(\tau)\mathsf{d}\tau &\geq \int_a^b (\tau-a)^p \boldsymbol{x}^\top(\tau)\left(\boldsymbol{j}_d^{0,p,\top}(\tau)_a^b \otimes I_n\right)\mathsf{d}\tau\,(\mathsf{D}_d \otimes U)\,[*] \\ &= \int_a^b \boldsymbol{x}^\top(\tau)\left(\boldsymbol{\ell}_{d+p}^\top(\tau) \otimes I_n\right)\mathsf{d}\tau\, \Xi^\top(\mathsf{D}_d \otimes U)\,\Xi \int_a^b \left(\boldsymbol{\ell}_{d+p}(\tau) \otimes I_n\right)\boldsymbol{x}(\tau)\mathsf{d}\tau \end{aligned} \tag{4.9}$$

for all $\boldsymbol{x}(\cdot) \in \mathscr{L}^2_{\varpi}\left(\mathcal{K};\mathbb{R}^n\right)$, where

$$\mathsf{D}_d^{-1} = \bigoplus_{k=0}^{d} \frac{(b-a)^{p+1}}{2k+1+p}, \quad \Xi := \mathbf{Y}(P \otimes I_n), \quad \mathbf{Y} := \left(\mathsf{J}_{0,p}\mathsf{L}^{-1}\right) \otimes I_n, \quad \boldsymbol{j}_d^{0,p}(\tau)_a^b = \mathsf{J}_{0,p} \operatorname*{\mathsf{Col}}_{i=0}^{d} \tau^i, \quad \boldsymbol{\ell}_d(\tau) = \mathsf{L} \operatorname*{\mathsf{Col}}_{i=0}^{d} \tau^i$$

with $\mathsf{J}_{0,p} \in \mathbb{R}^{(d+1)\times(d+1)}_{[d+1]}$ and $\mathsf{L} \in \mathbb{R}^{(d+1)\times(d+1)}_{[d+1]}$. Moreover, the matrix $P \in \mathbb{R}^{(d+1)\times(d+p+1)}$ satisfies $(\tau-a)^p\boldsymbol{\ell}_d(\tau) = P\boldsymbol{\ell}_{d+p}(\tau)$. By using the multiplier $(b-a)^{-p}$, it is evident that (27) in [621] is equivalent to (4.9). Now considering (4.9) with substitution $\dot{\boldsymbol{x}}(\tau) \to \boldsymbol{x}(\tau)$, we have

$$
\begin{aligned}
\int_a^b (\tau-a)^p \dot{\boldsymbol{x}}^\top(\tau) U \dot{\boldsymbol{x}}(\tau) \mathrm{d}\tau &\ge [*] \left[\Xi^\top \left(\mathsf{D}_d \otimes U\right) \Xi\right] \int_a^b \left(\boldsymbol{\ell}_{d+p}(\tau) \otimes I_n\right) \dot{\boldsymbol{x}}(\tau) \mathrm{d}\tau \\
&= \boldsymbol{\eta}_1^\top \Omega_1^\top \Xi^\top (\mathsf{D}_d \otimes U)\, \Xi \Omega_1 \boldsymbol{\eta}_1 = \boldsymbol{\eta}_2^\top \Omega_2^\top \Xi^\top (\mathsf{D}_d \otimes U)\, \Xi \Omega_2 \boldsymbol{\eta}_2
\end{aligned} \tag{4.10}
$$

where

$$
\begin{aligned}
\boldsymbol{\eta}_1 &:= \begin{bmatrix} \boldsymbol{x}(b) \\ \boldsymbol{x}(a) \\ \int_a^b \left(\boldsymbol{\ell}_{d+p}(\tau) \otimes I_n\right) \boldsymbol{x}(\tau) \mathrm{d}\tau \end{bmatrix}, \quad \boldsymbol{\eta}_2 := \begin{bmatrix} \boldsymbol{x}(b) \\ \boldsymbol{x}(a) \\ \int_a^b \left(\boldsymbol{\ell}_{d+p-1}(\tau) \otimes I_n\right) \boldsymbol{x}(\tau) \mathrm{d}\tau \end{bmatrix} \\
\Omega_1 &= \begin{bmatrix} \boldsymbol{\ell}_{d+p}(0), & \boldsymbol{\ell}_{d+p}(-r) & \Lambda_1 \end{bmatrix} \otimes I_n \quad \Omega_2 = \begin{bmatrix} \boldsymbol{\ell}_{d+p}(0) & \boldsymbol{\ell}_{d+p}(-r) & \Lambda_2 \end{bmatrix} \otimes I_n
\end{aligned} \tag{4.11}
$$

with $\Lambda_1 \in \mathbb{R}^{(d+p)\times(d+p)}$ and $\Lambda_2 \in \mathbb{R}^{(d+p)\times(d+p-1)}$ satisfying $\dot{\boldsymbol{\ell}}_{d+p}(\tau) = \Lambda_1 \boldsymbol{\ell}_{d+p}(\tau) = \Lambda_2 \boldsymbol{\ell}_{d+p-1}(\tau)$. Again by adjusting the factor $(b-a)^{-p}$ with (4.10), one can conclude that the result in [621, eq.(34)] can be obtained by (4.3).

It is clear that the structure of $\boldsymbol{f}(\cdot)$ in Theorem 4.1 may significantly affect the optimality of the LB in (4.3). In the following corollary, however, we show that a substitution $G\boldsymbol{f}(\tau) \to \boldsymbol{f}(\tau)$ for (4.3) with an invertible $G$ does not change the value of the LB in (4.3). This shows that when $d$ and $f_i(\tau)$ in $\boldsymbol{f}(\tau) = \mathbf{Col}_{i=1}^d f_i(\tau)$ are fixed, using linear combinations of $f_i(\tau)$ for (4.3) results in equivalent inequalities.

**Corollary 4.1.** *Given the same $\varpi(\cdot)$, $U \in \mathbb{S}^n_{\succ 0}$ and $\boldsymbol{f}(\cdot) \in \mathscr{L}^2_{\varpi}(\mathcal{K};\mathbb{R}^d)$ in Theorem 4.1, we have*

$$
\begin{aligned}
\int_{\mathcal{K}} \varpi(\tau) \boldsymbol{x}^\top(\tau) U \boldsymbol{x}(\tau) \mathrm{d}\tau &\ge \int_{\mathcal{K}} \varpi(\tau) \boldsymbol{x}^\top(\tau) F^\top(\tau) \mathrm{d}\tau \left(\mathsf{F} \otimes U\right) \int_{\mathcal{K}} \varpi(\tau) F(\tau) \boldsymbol{x}(\tau) \mathrm{d}\tau \\
&= \int_{\mathcal{K}} \varpi(\tau) \boldsymbol{x}^\top(\tau) \Phi^\top(\tau) \mathrm{d}\tau \left(\mathsf{\Phi} \otimes U\right) \int_{\mathcal{K}} \varpi(\tau) \Phi(\tau) \boldsymbol{x}(\tau) \mathrm{d}\tau
\end{aligned} \tag{4.12}
$$

*for all $\boldsymbol{x}(\cdot) \in \mathscr{L}^2_{\varpi}(\mathcal{K};\mathbb{R}^n)$ and $G \in \mathbb{R}^{n\times n}_{[n]}$, where $\Phi(\tau) = \boldsymbol{\varphi}(\tau) \otimes I_n$ with $\boldsymbol{\varphi}(\tau) = G\boldsymbol{f}(\tau)$, and $\mathsf{\Phi}^{-1} := \int_{\mathcal{K}} \varpi(\tau) \boldsymbol{\varphi}(\tau) \boldsymbol{\varphi}^\top(\tau) \mathrm{d}\tau$ and $F(\tau)$, $\mathsf{F}$ are the same defined in Theorem 4.1.*

*Proof.* Note that $\mathsf{\Phi}$ is calculated by the expression

$$
\begin{aligned}
\mathsf{\Phi}^{-1} = \int_{\mathcal{K}} \varpi(\tau) \boldsymbol{\varphi}(\tau) \boldsymbol{\varphi}^\top(\tau) \mathrm{d}\tau &= \int_{\mathcal{K}} \varpi(\tau) G \boldsymbol{f}(\tau) \boldsymbol{f}^\top(\tau) G^\top \mathrm{d}\tau \\
&= G \int_{\mathcal{K}} \varpi(\tau) \boldsymbol{f}(\tau) \boldsymbol{f}^\top(\tau) \mathrm{d}\tau G^\top = G \mathsf{F}^{-1} G^\top
\end{aligned} \tag{4.13}
$$

where $\mathsf{\Phi}^{-1}$ is well defined given the fact that $G \in \mathbb{R}^{d\times d}_{[d]}$. By (4.13) with the property of the Kronecker product in (2.1a), we have

$$
\begin{aligned}
&\int_{\mathcal{K}} \varpi(\tau) \boldsymbol{x}^\top(\tau) \Phi^\top(\tau) \mathrm{d}\tau \left(\mathsf{\Phi} \otimes U\right) \int_{\mathcal{K}} \varpi(\tau) \Phi(\tau) \boldsymbol{x}(\tau) \mathrm{d}\tau = \int_{\mathcal{K}} \boldsymbol{x}^\top(\tau) \left(\boldsymbol{\varphi}^\top(\tau) \otimes I_n\right) \mathrm{d}\tau \left[(G^{-1})^\top \mathsf{F} G^{-1} \otimes U\right] \times \\
&\int_{\mathcal{K}} \left(\boldsymbol{\varphi}(\tau) \otimes I_n\right) \boldsymbol{x}(\tau) \mathrm{d}\tau = \int_{\mathcal{K}} \boldsymbol{x}^\top(t+\tau) \left(\boldsymbol{f}^\top(\tau) G^\top \otimes I_n\right) \mathrm{d}\tau \left([*] \left(\mathsf{F} \otimes U\right) \left(G^{-1} \otimes I_n\right)\right) \int_{\mathcal{K}} \left(G \boldsymbol{f}(\tau) \otimes I_n\right) \boldsymbol{x}(t+\tau) \mathrm{d}\tau \\
&\qquad = \int_{\mathcal{K}} \boldsymbol{x}^\top(t+\tau) F^\top(\tau) \mathrm{d}\tau \left(\mathsf{F} \otimes U\right) \int_{\mathcal{K}} F(\tau) \boldsymbol{x}(t+\tau) \mathrm{d}\tau
\end{aligned}
$$

which gives (4.12) based on the inequality in (4.3). ■

**Remark 4.2.** Let $G$ in (4.12) be given by $G^{-2} = \int_{\mathcal{K}} \boldsymbol{f}(\tau) \boldsymbol{f}^\top(\tau) \mathrm{d}\tau$, which implies $G\boldsymbol{f}(\tau)$ only contains functions which are mutually orthogonal. Corollary 4.1 has an important implication: using function $\boldsymbol{f}(\cdot)$, which may not only contain orthogonal functions, for (4.3) does not degenerate the LB in (4.3) compared to using the orthogonal option $G\boldsymbol{f}(\tau)$.

## 4.3 The Second Class of IIs with Slack Variables

Inspired by [634, Theorem 4.1], where the maximum LB value is obtained when the slack variable is chosen to make the second inequality identical to the first in (4.3). The structure of the proposed inequality indicates that it can be useful for the analysis of sampled-data systems as pointed out in [634] or other potential problems that are subject to future research. We also show that our proposed inequality, in fact, generalizes most existing results in the literature.

We first present the following lemma which can be considered as an extension of [634, Lemma 4.1], where only real matrices are considered.

**Lemma 4.1.** *Given matrices $C \in \mathbb{S}^m_{\succ 0}$, $B \in \mathbb{R}^{m\times n}$, then*

$$\forall M \in \mathbb{R}^{m\times n}, \quad B^\top C^{-1} B \succeq M^\top B + B^\top M - M^\top C M \tag{4.14}$$

*where $B^\top C^{-1}B = M^\top B + B^\top M - M^\top C M$ can be obtained with $M = C^{-1}B$.*

*Proof.* It is obvious to see that

$$M^\top B + B^\top M - M^\top C M = B^\top C^{-1} B - \begin{bmatrix} -B^\top & M^\top C \end{bmatrix} \begin{bmatrix} C^{-1} & C^{-1} \\ C^{-1} & C^{-1} \end{bmatrix} \begin{bmatrix} -B \\ CM \end{bmatrix} \preceq B^\top C^{-1} B \tag{4.15}$$

since

$$\begin{aligned}\begin{bmatrix} -B^\top & M^\top C \end{bmatrix} \begin{bmatrix} C^{-1} & C^{-1} \\ C^{-1} & C^{-1} \end{bmatrix} \begin{bmatrix} -B \\ CM \end{bmatrix} &= \begin{bmatrix} -B^\top C^{-1} + M^\top & -B^\top C^{-1} + M^\top \end{bmatrix} \begin{bmatrix} -B \\ CM \end{bmatrix} \\ &= B^\top C^{-1} B - M^\top B - B^\top M + M^\top C M \end{aligned} \tag{4.16}$$

and $\begin{bmatrix} C^{-1} & C^{-1} \\ C^{-1} & C^{-1} \end{bmatrix} \succeq 0$. Moreover, $B^\top C^{-1} B = M^\top B + B^\top M - M^\top C M$ if $M = C^{-1}B$. ■

**Remark 4.3.** Compared to the proofs in [634], our proof of Lemma 4.1 is much simpler and can clearly demonstrate that (4.14) can be interpreted as a generalization of the idea of "completion of squares".

**Theorem 4.2.** *Given the same $\varpi(\cdot)$, $U$ and $\boldsymbol{f}(\cdot)$ defined in Theorem 4.1, then*

$$\forall \boldsymbol{x}(\cdot) \in \mathscr{L}^2_{\varpi}(\mathcal{K};\mathbb{R}^n), \quad \int_{\mathcal{K}} \varpi(\tau)\boldsymbol{x}^\top(\tau) U \boldsymbol{x}(\tau)\mathsf{d}\tau \geq \boldsymbol{\zeta}^\top \left[\mathsf{Sy}(H^\top \Omega) - H^\top \left(\mathsf{F}^{-1} \otimes U^{-1}\right) H\right] \boldsymbol{\zeta} \tag{4.17}$$

*with $H \in \mathbb{R}^{dn\times\nu}$, where $\Omega\boldsymbol{\zeta} = \int_{\mathcal{K}} \varpi(\tau)F(\tau)\boldsymbol{x}(\tau)\mathsf{d}\tau$ with $\boldsymbol{\zeta} \in \mathbb{R}^\nu$ and $\Omega \in \mathbb{R}^{dn\times\nu}$. Finally, with $H = (\mathsf{F} \otimes U)\,\Omega$, then (4.17) and (4.3) become identical and this is the case that the LB in (4.17) become the maximum.*

*Proof.* Since $\Omega\boldsymbol{\zeta} = \int_{\mathcal{K}} \varpi(\tau)F(\tau)\boldsymbol{x}(\tau)\mathsf{d}\tau$ with $\Omega \in \mathbb{R}^{dn\times\nu}$ and $\boldsymbol{\zeta} \in \mathbb{R}^\nu$, (4.3) can be rewritten as

$$\int_{\mathcal{K}} \varpi(\tau)\boldsymbol{x}^\top(\tau) U \boldsymbol{x}(\tau)\mathsf{d}\tau \geq \boldsymbol{\zeta}^\top \Omega^\top \left(\mathsf{F} \otimes U\right) \Omega \boldsymbol{\zeta}. \tag{4.18}$$

Now since $\mathsf{F} \succ 0$ and $U \succ 0$ with $\mathsf{F}^{-1} \otimes U^{-1} = \left(\mathsf{F} \otimes U\right)^{-1}$, applying Lemma 4.1 to the LB of (4.18) yields the results of Theorem 4.2. ■

**Remark 4.4.** Note that the definition $\Omega\boldsymbol{\zeta} = \int_{\mathcal{K}} \varpi(\tau)F(\tau)\boldsymbol{x}(\tau)\mathsf{d}\tau$ does not impose extra constraints on $\boldsymbol{f}(\cdot)$ or $\boldsymbol{x}(\cdot)$, since one can always find certain values of $\boldsymbol{\zeta}$ and $\Omega$ to make the equality hold. When analyzing the stability of systems with delays, one can choose a fixed $\Omega$ with appropriate $\boldsymbol{\zeta}$ to render $\Omega\boldsymbol{\zeta} = \int_{\mathcal{K}} \varpi(\tau)F(\tau)\boldsymbol{x}(\tau)\mathsf{d}\tau$ as an identity that is valid for all $\boldsymbol{x}(\cdot) \in \mathscr{L}^2_{\varpi}(\mathcal{K};\mathbb{R}^n)$.

**Remark 4.5.** Theorem 4.2 generalizes [634, Theorem 4.1] with $\boldsymbol{f}(\tau) = 1$ and $\varpi(\tau) = 1$. Furthermore, let $\mathcal{K} = [a, b]$ and $\varpi(\cdot) = 1$ and $\boldsymbol{f}(\tau)$ to contain Legendre polynomials over $[a, b]$, then [635, Lemma 1] can be obtained from (4.17) with appropriate $\boldsymbol{\zeta}$ and $\Omega$ using the substitution $\dot{\boldsymbol{x}}(\cdot) \to \boldsymbol{x}(\cdot)$. Now consider the fact that the left-hand side of [635, eq.(9)] can be rewritten as a single-fold integral with a weight function by way of the Cauchy formula for repeated integrations. Let $\mathcal{K} = [a, b]$, $\varpi(\tau) = (\tau - a)^m$ and $\boldsymbol{f}(\tau)$ to contain Jacobi polynomials associated with the weighted function $(\tau - a)^m$ over $[a, b]$, hence the integral inequality in [635] can be obtained by (4.17) with appropriate $\boldsymbol{\zeta}$ and $\Omega$ using the substitution $\dot{\boldsymbol{x}}(\cdot) \to \boldsymbol{x}(\cdot)$. Thus, all inequalities in [635] are the particular examples of (4.17). Finally, since (4.17) is equivalent to (4.3), it also indicates that equivalence relations can be established between the inequalities in [632, 635].

Similar to Corollary 4.1, we show in the following corollary that using a substitution $G\boldsymbol{f}(\tau) \to \boldsymbol{f}(\tau)$ for (4.17) with an invertible $G$ does not alter the maximum LB value in (4.17).

**Corollary 4.2.** *Given the same $\varpi(\cdot)$, $U$ and $\boldsymbol{f}(\cdot)$ as in Theorem 4.1, then for any $G \in \mathbb{R}^{n\times n}_{[n]}$ we have*

$$\forall \boldsymbol{x}(\cdot) \in \mathscr{L}^2_{\varpi}(\mathcal{K}; \mathbb{R}^n), \quad \int_{\mathcal{K}} \varpi(\tau)\boldsymbol{x}^\top(\tau) U \boldsymbol{x}(\tau) \mathrm{d}\tau \geq \widehat{\boldsymbol{\zeta}}^\top \left[ \mathsf{Sy}\left( \widehat{H}^\top \widehat{\Omega} \right) - H^\top \left( \Phi^{-1} \otimes U^{-1} \right) H \right] \widehat{\boldsymbol{\zeta}} \tag{4.19}$$

*with $H \in \mathbb{R}^{dn\times\nu}$, where $\Phi^{-1} := \int_{\mathcal{K}} \varpi(\tau)\boldsymbol{\varphi}(\tau)\boldsymbol{\varphi}^\top(\tau)\mathrm{d}\tau$ with $\boldsymbol{\varphi}(\tau) = G\boldsymbol{f}(\tau)$, and $\widehat{\Omega}\widehat{\boldsymbol{\zeta}} = \int_{\mathcal{K}} \varpi(\tau)\Phi(\tau)\boldsymbol{x}(\tau)\mathrm{d}\tau$ with $\widehat{\boldsymbol{\zeta}} \in \mathbb{R}^\nu$ and $\widehat{\Omega} \in \mathbb{R}^{dn\times\nu}$. Finally, (4.19) and (4.3) become identical with the assignment $\widehat{H} = (\mathsf{F} \otimes U)\widehat{\Omega}$, and this is the case that the maximum LB value in (4.19) is obtained which is invariant to the value of $G \in \mathbb{R}^{n\times n}_{[n]}$ and identical to the maximum LB value in (4.17).*

*Proof.* Let $G \in \mathbb{R}^{n\times n}_{[n]}$. Since $\widehat{\Omega}\widehat{\boldsymbol{\zeta}} = \int_{\mathcal{K}} \varpi(\tau)\Phi(\tau)\boldsymbol{x}(\tau)\mathrm{d}\tau$ with $\widehat{\boldsymbol{\zeta}} \in \mathbb{R}^\nu$ and $\widehat{\Omega} \in \mathbb{R}^{dn\times\nu}$. Then (4.12) can be rewritten as

$$\int_{\mathcal{K}} \varpi(\tau)\boldsymbol{x}^\top(\tau) U \boldsymbol{x}(\tau)\mathrm{d}\tau \geq \widehat{\boldsymbol{\zeta}}^\top \widehat{\Omega}^\top (\Phi \otimes U) \widehat{\Omega}\widehat{\boldsymbol{\zeta}}. \tag{4.20}$$

From (4.13), we obtain $\Phi \succ 0$ for all $G \in \mathbb{R}^{n\times n}_{[n]}$. Now since $\Phi \succ 0$ and $U \succ 0$, applying Lemma 4.1 to the LB of (4.18) yields the results in (4.2). ■

## 4.4 Another Integral Inequality of The Free Matrix Type

This section is dedicated to presenting another general integral inequality that includes additional matrix variables, concerning the analysis of the optimality of its LB with other properties. The proposed inequality belongs to the class of free matrix type, which has been previously researched in [623–625]. They have been utilized in addressing the stability analysis of systems with time-varying delays using the LKF approach. As mentioned in [624, Remark 7], the utilization of a free matrix type inequality can avoid making use of reciprocally convex combinations when analyzing the stability of a system with a time-varying delay. See [624, Remark 7] for further references therein. Finally, we can prove that the maximum LB value of the proposed inequality in this section is the same as (4.3) and (4.17) under the same $\varpi(\cdot)$, $U$ and $\boldsymbol{f}(\cdot)$, and it is invariant for any $G \in \mathbb{R}^{n\times n}_{[n]}$ if $G\boldsymbol{f}(\cdot)$ is considered.

For the derivations of the integral inequality, we make use of the following lemma, which can be straightforwardly obtained from the definition of matrix multiplication.

**Lemma 4.2.** *Given a matrix $X := \mathbf{Row}_{i=1}^d X_i \in \mathbb{R}^{n\times d\rho n}$ with $n; d; \rho \in \mathbb{N}$ and a function $\boldsymbol{f}(\tau) = \mathbf{Col}_{i=1}^d f_i(\tau) \in \mathbb{R}^d$, we have*

$$X(\boldsymbol{f}(\tau) \otimes I_{\rho n}) = \sum_{i=1}^{d} f_i(\tau) X_i = \left( \boldsymbol{f}^\top(\tau) \otimes I_n \right) \widehat{X} \tag{4.21}$$

*where* $\widehat{X} := \mathbf{Col}_{i=1}^{d} X_i \in \mathbb{R}^{dn \times \rho n}$.

**Theorem 4.3.** *Let* $\varpi(\cdot)$ *be given as in* (4.1) *and* $U \in \mathbb{S}_{\succ 0}^{n}$ *and* $\boldsymbol{f}(\cdot) \in \mathscr{L}_{\varpi}^{2}(\mathcal{K};\mathbb{R}^{d})$ *which satisfies the inequality in* (4.2)*. For any* $Y \in \mathbb{S}^{\rho dn}$ *and* $X = \mathbf{Row}_{i=1}^{d} X_i \in \mathbb{R}^{n \times \rho dn}$ *satisfying*

$$\begin{bmatrix} U & -X \\ * & Y \end{bmatrix} \succeq 0, \tag{4.22}$$

*we have*

$$\forall \boldsymbol{x}(\cdot) \in \mathscr{L}_{\varpi}^{2}(\mathcal{K};\mathbb{R}^{n}), \quad \int_{\mathcal{K}} \varpi(\tau)\boldsymbol{x}^{\top}(\tau)U\boldsymbol{x}(\tau)\mathrm{d}\tau \geq \boldsymbol{z}^{\top}\left[\mathsf{Sy}\left(\Upsilon^{\top}\widehat{X}\right) - \mathsf{W}\right]\boldsymbol{z}, \tag{4.23}$$

*where* $\rho \in \mathbb{N}$ *and* $\mathsf{W} := \int_{\mathcal{K}} \varpi(\tau)(\boldsymbol{f}^{\top}(\tau) \otimes I_{\rho n})Y(\boldsymbol{f}(\tau) \otimes I_{\rho n})\mathrm{d}\tau \in \mathbb{S}^{\rho n}$ *and* $\widehat{X} = \mathbf{Col}_{i=1}^{d} X_i \in \mathbb{R}^{dn \times \rho n}$*, and* $\Upsilon\boldsymbol{z} = \int_{\mathcal{K}} \varpi(\tau)F(\tau)\boldsymbol{x}(\tau)\mathrm{d}\tau$ *with* $\Upsilon \in \mathbb{R}^{dn \times \rho n}$ *and* $\boldsymbol{z} \in \mathbb{R}^{\rho n}$.

*Proof.* Given (4.22) and $\Upsilon\boldsymbol{z} = \int_{\mathcal{K}} \varpi(\tau)F(\tau)\boldsymbol{x}(\tau)\mathrm{d}\tau$ with $\Upsilon \in \mathbb{R}^{dn \times \rho n}$ and $\boldsymbol{z} \in \mathbb{R}^{\rho n}$, $\rho \in \mathbb{N}$, we have

$$\begin{aligned}\int_{\mathcal{K}} \varpi(\tau)\,[*]^{\top}\begin{bmatrix} U & -X \\ * & Y \end{bmatrix}\begin{bmatrix} \boldsymbol{x}(\tau) \\ \boldsymbol{f}(\tau) \otimes \boldsymbol{z} \end{bmatrix}\mathrm{d}\tau = \int_{\mathcal{K}} \varpi(\tau)\boldsymbol{x}^{\top}(\tau)U\boldsymbol{x}(\tau)\mathrm{d}\tau - \mathsf{Sy}\left[\int_{\mathcal{K}} \varpi(\tau)\boldsymbol{x}^{\top}(\tau)X\left(\boldsymbol{f}(\tau) \otimes \boldsymbol{z}\right)\mathrm{d}\tau\right] \\ + \int_{\mathcal{K}} \varpi(\tau)\left(\boldsymbol{f}(\tau) \otimes \boldsymbol{z}\right)^{\top} Y \left(\boldsymbol{f}(\tau) \otimes \boldsymbol{z}\right)\mathrm{d}\tau \geq 0.\end{aligned} \tag{4.24}$$

Now applying the property of the Kronecker product in (2.1a) with (4.21) to the terms in (4.24) yields

$$\begin{aligned}\int_{\mathcal{K}} \varpi(\tau)\boldsymbol{x}^{\top}(\tau)X\left(\boldsymbol{f}(\tau) \otimes \boldsymbol{z}\right)\mathrm{d}\tau = \int_{\mathcal{K}} \varpi(\tau)\boldsymbol{x}^{\top}(\tau)X(\boldsymbol{f}(\tau) \otimes I_{\rho n})\mathrm{d}\tau\boldsymbol{z} = \left(\sum_{i=1}^{d}\int_{\mathcal{K}} \varpi(\tau)\boldsymbol{x}^{\top}(\tau)f_i(\tau)\mathrm{d}\tau X_i\right)\boldsymbol{z} \\ = \int_{\mathcal{K}} \varpi(\tau)\boldsymbol{x}^{\top}(\tau)(\boldsymbol{f}^{\top}(\tau) \otimes I_n)\mathrm{d}\tau\widehat{X}\boldsymbol{z} = \boldsymbol{z}^{\top}\Upsilon^{\top}\widehat{X}\boldsymbol{z}\end{aligned} \tag{4.25}$$

and

$$\int_{\mathcal{K}} \varpi(\tau)\left(\boldsymbol{f}(\tau) \otimes \boldsymbol{z}\right)^{\top} Y \left(\boldsymbol{f}(\tau) \otimes \boldsymbol{z}\right)\mathrm{d}\tau = \boldsymbol{z}^{\top}\int_{\mathcal{K}} \varpi(\tau)\left(\boldsymbol{f}^{\top}(\tau) \otimes I_{\rho n}\right) Y \left(\boldsymbol{f}(\tau) \otimes I_{\rho n}\right)\mathrm{d}\tau\boldsymbol{z} = \boldsymbol{z}^{\top}\mathsf{W}\boldsymbol{z} \tag{4.26}$$

where $X = \mathbf{Row}_{i=1}^{d} X_i \in \mathbb{R}^{n \times \rho dn}$ and $\widehat{X} = \mathbf{Col}_{i=1}^{d} X_i$. With (4.25) and (4.26) at hand, (4.24) gives (4.23). This finishes the proof. ■

As previously noted, one potential motivation for proposing an inequality with extra matrix parameters is its utility in analyzing the stability of systems with time-varying delays. However, the free matrix structure of (4.23) may also present potential advantages, which can be explored in future works.

**Remark 4.6.** Since $\boldsymbol{f}(\cdot)$ in Theorem 4.3 is subject to the same constraint $\int_{\mathcal{K}} \varpi(\tau)\boldsymbol{f}(\tau)\boldsymbol{f}^{\top}(\tau)\mathrm{d}\tau \succ 0$ as the $\boldsymbol{f}(\cdot)$ in Theorem 4.1, the structure of (4.23) is more general than existing free matrix type inequalities in the literature. Moreover, let $\mathcal{K} = [a, b]$ and $\varpi(\cdot) = 1$ and $\boldsymbol{f}(\tau)$ comprising the Legendre polynomials over $[a, b]$, then one can obtain [624, Lemma 3] by Theorem 4.3 with appropriate $\Upsilon$ and $\boldsymbol{z}$ and the substitution $\dot{\boldsymbol{x}}(\cdot) \to \boldsymbol{x}(\cdot)$. Note that this also means that Theorem 4.3 covers the special cases of [624, Lemma 3] such as [623] mentioned therein. Furthermore, let $\mathcal{K} = [a, b]$ and $\varpi(\cdot) = 1$, then by (4.23) with appropriate $\boldsymbol{f}(\cdot)$ and without considering $\Upsilon\boldsymbol{z} = \int_{\mathcal{K}} \varpi(\tau)F(\tau)\boldsymbol{x}(\tau)\mathrm{d}\tau$, one can derive an inequality, which is equivalent to the [625, Lemma 3], given the proofs of Theorem 4.3.

The following theorem shows the relation among the maximum values of the LBs of (4.3) and (4.17) and (4.23).

**Theorem 4.4.** *By choosing the same $\varpi(\cdot)$, $U$ and $\boldsymbol{f}(\cdot)$ for Theorem 4.1, Theorem 4.2 and Theorem 4.3, it is always possible to construct appropriate $X$ and $Y$ for (4.22) such that (4.23) becomes (4.3). In addition, the maximum LB value in (4.23) is identical to that in (4.17), which in turn is the same as the value of the LB in (4.3).*

*Proof.* See Appendix C. ■

**Remark 4.7.** Let $\mathcal{K} = [a, b]$, $\varpi(\cdot) = 1$ and let $\boldsymbol{f}(\tau)$ includes Legendre polynomials over $[a, b]$, then [624, Theorem 1] can be established from Theorem 4.4 with appropriate $\Upsilon$ and $\boldsymbol{z}$ using substitution $\dot{\boldsymbol{x}}(\cdot) \to \boldsymbol{x}(\cdot)$. Since we have proved that (4.17) is equivalent to (4.3), inequality (4.23) is also equivalent to (4.17). Consequently, it is possible to show that there are comprehensive equivalence relations[2] between the inequalities in [624, 632, 635], as discussed in Remark 4.5.

Theorem 4.4 serves as a pivotal point in bridging the relations between (4.3), (4.17) and (4.23). Since all three inequalities are equivalent with respect to the maximum value of their LBs, if one finds a special example of one of these three inequalities, then two equivalent inequalities can be immediately constructed.

The following Corollary 4.3 can be established for (4.23) in a manner similar to Corollary 4.1.

**Corollary 4.3.** *Given the same $\varpi(\cdot)$, $U$ and $\boldsymbol{f}(\cdot)$ in Theorem 4.3, then for any $Y \in \mathbb{S}^{\rho dn}$ and $X = \mathbf{Row}_{i=1}^{d} X_i \in \mathbb{R}^{n\times \rho dn}$ satisfying*

$$\begin{bmatrix} U & -X \\ * & Y \end{bmatrix} \succeq 0, \tag{4.27}$$

*we have*

$$\forall \boldsymbol{x}(\cdot) \in \mathscr{L}^2_{\varpi}(\mathcal{K}; \mathbb{R}^n), \quad \int_{\mathcal{K}} \varpi(\tau)\boldsymbol{x}^\top(\tau)U\boldsymbol{x}(\tau)\mathrm{d}\tau \geq \widehat{\boldsymbol{z}}^\top \left[\mathsf{Sy}\left(\Pi^\top \widehat{X}\right) - \mathsf{V}\right]\widehat{\boldsymbol{z}} \tag{4.28}$$

*for all $G \in \mathbb{R}^{n\times n}_{[n]}$, where $\widehat{X} = \mathbf{Col}_{i=1}^{d} X_i \in \mathbb{R}^{dn\times \rho n}$ and*

$$\begin{aligned} &\mathbb{S}^{\rho n} \ni \mathsf{V} := \int_{\mathcal{K}} \varpi(\tau)(\boldsymbol{\varphi}^\top(\tau) \otimes I_{\rho n})Y(\boldsymbol{\varphi}(\tau) \otimes I_{\rho n})\mathrm{d}\tau, \quad \boldsymbol{\varphi}(\tau) = G\boldsymbol{f}(\tau) \\ &\Pi\widehat{\boldsymbol{z}} = \int_{\mathcal{K}} \varpi(\tau)\Phi(\tau)\boldsymbol{x}(\tau)\mathrm{d}\tau, \quad \Pi \in \mathbb{R}^{dn\times \rho n}, \quad \widehat{\boldsymbol{z}} \in \mathbb{R}^{\rho n}, \quad \Phi(\tau) = \boldsymbol{\varphi}(\tau) \otimes I_n. \end{aligned} \tag{4.29}$$

***Finally, given the same** $\varpi(\cdot)$, $U$ **and** $\boldsymbol{f}(\cdot)$, the maximum LB value in (4.28) is identical to that in (4.23). The maximum value is equal to the value of the LB in (4.3), and it is invariant to $G \in \mathbb{R}^{n\times n}_{[n]}$.*

*Proof.* Let $\varpi(\cdot)$, $U$ and $\boldsymbol{f}(\cdot)$ from Theorem 4.3 be given throughout the entire proof. The inequality in (4.28) can be obtained based on the substitution $G\boldsymbol{f}(\cdot) = \boldsymbol{\varphi}(\cdot) \to \boldsymbol{f}(\cdot)$ in (4.23). Now consider the inequality

$$\int_{\mathcal{K}} \varpi(\tau)\boldsymbol{x}^\top(\tau)U\boldsymbol{x}(\tau)\mathrm{d}\tau \geq \int_{\mathcal{K}} \varpi(\tau)\boldsymbol{x}^\top(\tau)\Phi^\top(\tau)\mathrm{d}\tau\,(\Phi \otimes U)\int_{\mathcal{K}} \varpi(\tau)\Phi(\tau)\boldsymbol{x}(\tau)\mathrm{d}\tau \tag{4.30}$$

in (4.12). By applying Theorem 4.4 with the fact that $\Pi\widehat{\boldsymbol{z}} = \int_{\mathcal{K}} \varpi(\tau)\Phi(\tau)\boldsymbol{x}(\tau)\mathrm{d}\tau$, it follows that the maximum LB value in (4.28) is identical to that in (4.30). Since both (4.3) and (4.30) are part of (4.12), we can conclude that maximum LB value in (4.28) is equal to the one in (4.3) which is invariant to the values of $G \in \mathbb{R}^{n\times n}_{[n]}$. Furthermore, as Theorem 4.4 has shown that the maximum LB value in (4.3) is identical to the one in (4.23), it proves the results in Corollary 4.3. ■

[2] The equivalence relations here are understood by considering the structure of inequalities irrespective of using $\boldsymbol{x}(\cdot)$ or $\dot{\boldsymbol{x}}(\cdot)$.

**Remark 4.8.** Upon a careful examination of the results we have proposed, we can establish a coherent chain of relations among the inequalities in [636]– [633] (see Table 4.1) and their "slack variables" counterpart obtained from the inequalities we have presented with appropriate $\Omega, \boldsymbol{\zeta}$ in (4.17), and $\widehat{\Omega}, \widehat{\boldsymbol{\zeta}}$ in (4.19), and $\Upsilon, \boldsymbol{z}$ in (4.23), and $\Pi, \widehat{\boldsymbol{z}}$ in (4.28).

## 4.5 Applications of IIs to the Stability Analysis of A System with Delays

To demonstrate the usefulness of the results we have established, we derive a stability condition in this section for a linear CDDS with a DD via constructing a parameterized version of the complete LKF [22] based on the application of (4.3). We show that the feasibility of the resulting stability condition remains unchanged if two different, yet linearly independent, groups of integral kernel functions are utilized for the construction of LKFs. Furthermore, equivalent stability conditions that preserve this invariance property can also be derived through the application of (4.17) and (4.23).

Consider a linear coupled differential-difference system in the form of

$$\begin{aligned} \dot{\boldsymbol{x}}(t) &= A_1\boldsymbol{x}(t) + A_2\boldsymbol{y}(t-r) + \int_{-r}^{0} \widetilde{A}_3(\tau)\boldsymbol{y}(t+\tau)\mathsf{d}\tau \\ \boldsymbol{y}(t) &= A_4\boldsymbol{x}(t) + A_5\boldsymbol{y}(t-r) \\ \boldsymbol{x}(t_0) &= \boldsymbol{\xi}, \quad \forall\theta \in [-r, 0], \quad \boldsymbol{y}(t_0+\theta) = \boldsymbol{\phi}(\theta) \end{aligned} \tag{4.31}$$

where $t_0 \in \mathbb{R}$ and $\boldsymbol{\xi} \in \mathbb{R}^n$ and $\boldsymbol{\phi}(\cdot) \in \widehat{C}([-r, 0); \mathbb{R}^\nu)$, and the notation $\widehat{C}([-r, 0); \mathbb{R}^n)$ stands for the space of bounded right piecewise continuous functions endowed with norm $\|\boldsymbol{\phi}(\cdot)\|_\infty = \sup_{\tau\in\mathcal{X}} \|\boldsymbol{\phi}(\tau)\|_2$. Furthermore, $\boldsymbol{x}(t) \in \mathbb{R}^n$ with $\boldsymbol{y}(t) \in \mathbb{R}^\nu$ is the solution to (4.31) and the size of the state space parameters in (4.31) are determined by $n, \nu \in \mathbb{N}$. We also assume that $\rho(A_5) < 1$ which ensures the input to state stability of $\boldsymbol{y}(t) = A_4\boldsymbol{x}(t) + A_5\boldsymbol{y}(t-r)$ [22], where $\rho(A_5)$ is the spectral radius of $A_5$. Since $\rho(A_5) < 1$ is independent of $r$, hence this condition ensures the input to state stability of $\boldsymbol{y}(t) = A_4\boldsymbol{x}(t) + A_5\boldsymbol{y}(t-r)$ for all $r > 0$. Finally, $\widetilde{A}_3(\tau)$ satisfies the following assumption.

**Assumption 4.1.** There exist function $\mathbf{Col}_{i=1}^d f_i(\tau) = \boldsymbol{f}(\cdot) \in C^1(\mathbb{R}; \mathbb{R}^d)$ with $d \in \mathbb{N}$, and $A_3 \in \mathbb{R}^{n\times d\nu}$ such that for all $\tau \in [-r, 0]$ we have $\widetilde{A}_3(\tau) = A_3 F(\tau) \in \mathbb{R}^{n\times\nu}$ where $F(\tau) := \boldsymbol{f}(\tau) \otimes I_\nu \in \mathbb{R}^{d\nu\times\nu}$. In addition, we assume that $\boldsymbol{f}(\cdot)$ satisfies $\int_{-r}^{0} \boldsymbol{f}(\tau)\boldsymbol{f}^\top(\tau)\mathsf{d}\tau \succ 0$ and

$$\int_{-r}^{0} \boldsymbol{f}(\tau)\boldsymbol{f}^\top(\tau)\mathsf{d}\tau \succ 0 \tag{4.32}$$

$$\exists M \in \mathbb{R}^{d\times d}, \quad \frac{\mathsf{d}\boldsymbol{f}(\tau)}{\mathsf{d}\tau} = M\boldsymbol{f}(\tau), \tag{4.33}$$

$$\exists N_1 \in \mathbb{R}_{[\delta]}^{\delta\times d}, \quad \exists N_2 \in \mathbb{R}_{[\delta]}^{\delta\times d}, \quad (\tau + r)N_1\boldsymbol{f}(\tau) = N_2\boldsymbol{f}(\tau) \tag{4.34}$$

where $0 \le \delta \le d \in \mathbb{N}$.

Examples of $\boldsymbol{f}(\cdot)$ in Assumption 4.1 are the solutions to homogeneous differential equations. For instance, $\boldsymbol{f}(\tau) = \mathbf{Col}\,(1, \mathsf{e}^\tau, \tau, \tau\mathsf{e}^\tau)$ and $N_1\boldsymbol{f}(\tau) = \mathbf{Col}\,(1, \mathsf{e}^\tau)$ with parameters

$$M = \begin{bmatrix} 0 & 0 & 0 & 0 \\ 0 & 1 & 0 & 0 \\ 1 & 0 & 0 & 0 \\ 0 & 1 & 0 & 1 \end{bmatrix}, \quad N_1 = \begin{bmatrix} I_2 & \mathsf{O}_2 \end{bmatrix}, \quad N_2 = \begin{bmatrix} r & 1 & 0 & 0 \\ 0 & r & 0 & 1 \end{bmatrix}. \tag{4.35}$$

Note that for any $\boldsymbol{f}(\cdot)$ satisfying (4.33), one can always enlarge the dimension of $\boldsymbol{f}(\cdot)$ by adding new functions to ensure that $\boldsymbol{f}(\cdot)$ satisfies (4.34). Moreover, (4.34) is satisfied for any $\boldsymbol{f}(\cdot) \in \mathcal{C}^1(\mathbb{R};\mathbb{R}^d)$ if $N_1$ and $N_2$ are empty matrices, implying that the constraint in (4.34) can be omitted on demand. Note that the rank constraint on $N_1$ ensures that $N_1 \int_{-r}^{0}(\tau+r)\boldsymbol{f}(\tau)\boldsymbol{f}^\top(\tau)\mathrm{d}\tau N_1^\top \succ 0$.

To prove the results in this section, we introduce the subsequent lemma that is a Lyapunov-Krasovskiĭ stability criteria for (4.31).

**Lemma 4.3.** *Given $r > 0$, the system (4.31) is globally uniformly asymptotically (exponentially)*[3] *stable at its origin, if there exist $\epsilon_1;\epsilon_2;\epsilon_3 > 0$ and a differentiable functional $\mathsf{v} : \mathbb{R}^n \times \widehat{\mathcal{C}}([-r,0);\mathbb{R}^\nu) \to \mathbb{R}_{\geq 0}$ such that $\mathsf{v}(\mathbf{0}_n,\mathbf{0}_\nu) = 0$ and*

$$\epsilon_1 \left\|\boldsymbol{\xi}\right\|_2^2 \leq \mathsf{v}(\boldsymbol{\xi},\boldsymbol{\phi}(\cdot)) \leq \epsilon_2 \left(\left\|\boldsymbol{\xi}\right\|_2 \vee \left\|\boldsymbol{\phi}(\cdot)\right\|_\infty\right)^2 \tag{4.36}$$

$$\dot{\mathsf{v}}(\boldsymbol{\xi},\boldsymbol{\phi}(\cdot)) := \left.\frac{\mathsf{d}^+}{\mathsf{d}t}\mathsf{v}(\boldsymbol{x}(t),\boldsymbol{\mathfrak{y}}_t(\cdot))\right|_{t=t_0,\boldsymbol{x}(t_0)=\boldsymbol{\xi},\boldsymbol{\mathfrak{y}}_{t_0}(\cdot)=\boldsymbol{\phi}(\cdot)} \leq -\epsilon_3 \left\|\boldsymbol{\xi}\right\|_2^2 \tag{4.37}$$

*for any $\boldsymbol{\xi} \in \mathbb{R}^n$ and $\boldsymbol{\phi}(\cdot) \in \widehat{\mathcal{C}}([-r,0);\mathbb{R}^\nu)$ in (4.31), where $t_0 \in \mathbb{R}$ and $\frac{\mathsf{d}^+}{\mathsf{d}x}f(x) = \operatorname{limsup}_{\eta\downarrow 0}\frac{f(x+\eta)-f(x)}{\eta}$. Furthermore, $\boldsymbol{\mathfrak{y}}_t(\cdot)$ in (4.37) is defined by the equality $\forall t \geq t_0$, $\forall\theta \in [-r,0)$, $\boldsymbol{\mathfrak{y}}_t(\theta) = \boldsymbol{y}(t+\theta)$ where $\boldsymbol{x}(t)$ and $\boldsymbol{y}(t)$ satisfying (4.31).*

*Proof.* Let $u(\cdot),\mathsf{v}(\cdot),w(\cdot)$ in [22, Theorem 3] be quadratic functions with multiplier factors $\epsilon_1;\epsilon_2;\epsilon_3 > 0$. Since (4.31) is a particular case of the general system considered in [22, Theorem 3], then Lemma 4.3 is obtained accordingly. ■

Consider the following parameterized LKF

$$\mathsf{v}(\boldsymbol{\xi},\boldsymbol{\phi}(\cdot)) := \begin{bmatrix}\boldsymbol{\xi}\\ \int_{-r}^{0}\widehat{G}(\tau)\boldsymbol{\phi}(\tau)\mathsf{d}\tau\end{bmatrix}^\top \widehat{P}\begin{bmatrix}\boldsymbol{\xi}\\ \int_{-r}^{0}\widehat{G}(\tau)\boldsymbol{\phi}(\tau)\mathsf{d}\tau\end{bmatrix} + \int_{-r}^{0}\boldsymbol{\phi}^\top(\tau)\left[S+(\tau+r)U\right]\boldsymbol{\phi}(\tau)\mathsf{d}\tau \tag{4.38}$$

with

$$\widehat{G}(\tau) = \boldsymbol{g}(\tau)\otimes I_\nu,\quad \boldsymbol{g}(\tau) = G\boldsymbol{f}(\tau),\qquad G \in \mathbb{R}_{[d]}^{d\times d}$$

and $\boldsymbol{f}(\cdot)$ in Assumption 4.1, where $\boldsymbol{\xi} \in \mathbb{R}^n$, $\boldsymbol{\phi}(\cdot) \in \widehat{\mathcal{C}}([-r,0);\mathbb{R}^\nu)$ in (4.38) are the initial conditions in (4.31), and $\widehat{P} \in \mathbb{S}^{n+d\nu}$ and $S;U \in \mathbb{S}^\nu$ are unknown parameters need to be determined. Note that $\widehat{G}(\tau)$ can be rewritten as $\widehat{G}(\tau) = \boldsymbol{g}(\tau)\otimes I_\nu = G\boldsymbol{f}(\tau)\otimes I_\nu = (G\otimes I_\nu)F(\tau)$ with $F(\tau) := \boldsymbol{f}(\tau)\otimes I_\nu$ based on the property of the Kronecker product in (2.1a). Note that also (4.38) can be interpreted as a parameterized version of the complete LKF proposed in [22].

In the subsequent theorem, we demonstrate that the feasibility of the resulting stability condition remains unchanged for any $G \in \mathbb{R}_{[d]}^{d\times d}$ in (4.38), irrespective of whether (4.3) or (4.23) is employed for the derivation.

**Theorem 4.5.** *Given $G \in \mathbb{R}_{[d]}^{d\times d}$ and $\boldsymbol{f}(\cdot)$ with $M$, $N_1$ and $N_2$ in Assumption 4.1, then then the origin of (4.31) under Assumption 4.1 is globally uniformly asymptotically stable if there exists $\widehat{P} \in \mathbb{S}^{n+d\nu}$ and $S;U \in \mathbb{S}^\nu$ such that*

$$\widehat{P} + \left[\mathsf{O}_n \oplus \left(G^{-1\top}\mathsf{F}G^{-1}\otimes S + \left(G^{-1\top}N_2^\top\widetilde{\mathsf{F}}N_2G^{-1}\right)\otimes U\right)\right] \succ 0, \tag{4.39}$$

$$S \succ 0,\ \ U \succ 0,\ \ \boldsymbol{\Phi} \prec 0, \tag{4.40}$$

[3] See [22] for the exposition on the equivalence bewteen uniform asymptotic and exponential stability for a linear coupled differential functional system

*hold, where* $\mathsf{F}^{-1} = \int_{-r}^{0} \boldsymbol{f}(\tau)\boldsymbol{f}^{\top}(\tau)\mathsf{d}\tau$ *and* $\widetilde{\mathsf{F}}^{-1} = \int_{-r}^{0}(\tau+r)N_1\boldsymbol{f}(\tau)\boldsymbol{f}^{\top}(\tau)N_1^{\top}\mathsf{d}\tau$ *and*

$$\boldsymbol{\Phi} := \mathsf{Sy}\left(H\widehat{P}\begin{bmatrix}\mathbf{A}\\ \mathbf{G}\end{bmatrix}\right) + \Gamma^{\top}(S+rU)\Gamma - \left(\mathsf{O}_n \oplus S \oplus \left(G^{-1\top}\mathsf{F}G^{-1}\otimes U\right)\right) \tag{4.41}$$

*with*

$$H = \begin{bmatrix} I_n & \mathsf{O}_{n,d\nu} \\ \mathsf{O}_{\nu,n} & \mathsf{O}_{\nu,d\nu} \\ \mathsf{O}_{d\nu,n} & I_{d\nu} \end{bmatrix}, \quad \Gamma := \begin{bmatrix} A_4 & A_5 & \mathsf{O}_{\nu,d\nu} \end{bmatrix}, \tag{4.42}$$

$$\mathbf{A} = \begin{bmatrix} A_1 & A_2 & A_3(G^{-1}\otimes I_\nu) \end{bmatrix}, \tag{4.43}$$

$$\mathbf{G} = \begin{bmatrix} \widehat{G}(0)A_4 & \widehat{G}(0)A_5 - \widehat{G}(-r) & -\widehat{M} \end{bmatrix}, \tag{4.44}$$

*where* $\widehat{G}(0) = (G\otimes I_\nu)F(0)$ *and* $\widehat{G}(-r) = (G\otimes I_\nu)F(-r)$ *and* $\widehat{M} = (G\otimes I_\nu)(M\otimes I_\nu)(G^{-1}\otimes I_\nu)$. *Moreover, the feasibility of* (4.39) *and* (4.40) *is invariant for any* $G\in\mathbb{R}_{[d]}^{d\times d}$. *Note that* (4.39) *and* (4.40) *are derived by the application of* (4.3). *In addition, if* (4.17) *or* (4.23) *are applied instead of* (4.3) *for the derivation of stability conditions, then the corresponding stability conditions are equivalent to* (4.39) *and* (4.40)*, respectively, and the feasibility of the resulting conditions remain invariant for any* $G\in\mathbb{R}_{[d]}^{d\times d}$.

*Proof.* Let $G\in\mathbb{R}_{[d]}^{d\times d}$ and $\boldsymbol{f}(\cdot)$ with $M$, $N_1$ and $N_2$ in Assumption 4.1 be given. Given the fact that the eigenvalues of $S+(\tau+r)U$, $\tau\in[-r,0]$ are bounded and $\widehat{G}(\tau) = (G\otimes I_n)F(\tau)$, it follows that (4.38) satisfies

$$\begin{aligned}
&\exists\lambda>0,\ \exists\eta>0,\ \mathsf{v}(\boldsymbol{\xi},\boldsymbol{\phi}(\cdot)) \le \begin{bmatrix}\boldsymbol{\xi}\\ \int_{-r}^{0}F(\tau)\boldsymbol{\phi}(\tau)\mathsf{d}\tau\end{bmatrix}^{\top}\lambda\begin{bmatrix}\boldsymbol{\xi}\\ \int_{-r}^{0}F(\tau)\boldsymbol{\phi}(\tau)\mathsf{d}\tau\end{bmatrix} + \int_{-r}^{0}\boldsymbol{\phi}^{\top}(\tau)\lambda\boldsymbol{\phi}(\tau)\mathsf{d}\tau \\
&\le \lambda\left\|\boldsymbol{\xi}\right\|_2^2 + \int_{-r}^{0}\boldsymbol{\phi}^{\top}(\tau)F^{\top}(\tau)\mathsf{d}\tau\lambda\int_{-r}^{0}F(\tau)\boldsymbol{\phi}(\tau)\mathsf{d}\tau + r\lambda\left\|\boldsymbol{\phi}(\cdot)\right\|_\infty^2 \\
&\le \lambda\left\|\boldsymbol{\xi}\right\|_2^2 + \lambda r\left\|\boldsymbol{\phi}(\cdot)\right\|_\infty^2 + \int_{-r}^{0}\boldsymbol{\phi}^{\top}(\tau)F^{\top}(\tau)\mathsf{d}\tau\left(\eta\mathsf{F}\otimes I_n\right)\int_{-r}^{0}F(\tau)\boldsymbol{\phi}(\tau)\mathsf{d}\tau \\
&\le \lambda\left\|\boldsymbol{\xi}\right\|_2^2 + \lambda r\left\|\boldsymbol{\phi}(\cdot)\right\|_\infty^2 + \int_{-r}^{0}\boldsymbol{\phi}^{\top}(\tau)\eta I_n\boldsymbol{\phi}(\tau)\mathsf{d}\tau \le \lambda\left\|\boldsymbol{\xi}\right\|_2^2 + (\lambda r+\eta r)\left\|\boldsymbol{\phi}(\cdot)\right\|_\infty^2 \\
&\le (\lambda r+\eta r)\left\|\boldsymbol{\xi}\right\|_2^2 + (\lambda r+\eta r)\left\|\boldsymbol{\phi}(\cdot)\right\|_\infty^2 \le 2\left(\lambda r+\eta r\right)\left[\max\left(\left\|\boldsymbol{\xi}\right\|_2, \left\|\boldsymbol{\phi}(\cdot)\right\|_\infty\right)\right]^2
\end{aligned} \tag{4.45}$$

for any $\boldsymbol{\xi}\in\mathbb{R}^n$ and $\boldsymbol{\phi}(\cdot)\in\widehat{\mathscr{C}}\left([-r_2,0);\mathbb{R}^\nu\right)$ in (4.31), where (4.45) is derived via the property of quadratic forms: $\forall X\in\mathbb{S}^n, \exists\lambda>0:\forall\mathbf{x}\in\mathbb{R}^n\setminus\{\mathbf{0}\}, \mathbf{x}^{\top}\left(\lambda I_n - X\right)\mathbf{x}>0$ together with the application of (4.3) with $\boldsymbol{f}(\cdot)$ in (4.1). Consequently, we have shown from (4.45) that (4.38) satisfies the upper bound property in (4.36).

Now, when we apply (4.3) with $\varpi(\tau)=1$ to the integral term $\int_{-r}^{0}\boldsymbol{\phi}^{\top}(\tau)S\boldsymbol{\phi}(\tau)\mathsf{d}\tau$ in (4.38) given $S\succ 0$ and $\boldsymbol{f}(\cdot)$ in Assumption 4.1 and the fact that $\boldsymbol{\phi}(\cdot)\in\widehat{\mathscr{C}}([-r,0);\mathbb{R}^\nu)\subset\mathscr{L}^2([-r,0);\mathbb{R}^\nu)$. It yields

$$\int_{-r}^{0}\boldsymbol{\phi}^{\top}(\tau)S\boldsymbol{\phi}(\tau)\mathsf{d}\tau \ge \int_{-r}^{0}\boldsymbol{\phi}^{\top}(\tau)\widehat{G}^{\top}(\tau)\mathsf{d}\tau\left[(G^{-1})^{\top}\mathsf{F}G^{-1}\otimes S\right]\int_{-r}^{0}\widehat{G}(\tau)\boldsymbol{\phi}(\tau)\mathsf{d}\tau \tag{4.46}$$

for any initial condition $\boldsymbol{\xi}\in\mathbb{R}^n$ and $\boldsymbol{\phi}(\cdot)\in\widehat{\mathscr{C}}\left([-r_2,0);\mathbb{R}^\nu\right)$ in (4.31). On the other hand, apply (4.3) with $\varpi(\tau)=\tau+r$ to the term $\int_{-r}^{0}(\tau+r)\boldsymbol{\phi}^{\top}(\tau)U\boldsymbol{\phi}(\tau)\mathsf{d}\tau$ in (4.38) with $U\succ 0$ and $\boldsymbol{f}(\cdot)$ in Assumption 4.1. It yields

$$\int_{-r}^{0}(\tau+r)\boldsymbol{\phi}^{\top}(\tau)S\boldsymbol{\phi}(\tau)\mathsf{d}\tau \ge [*]\left(\widetilde{\mathsf{F}}\otimes U\right)\left[\int_{-r}^{0}(\tau+r)\left(N_1\boldsymbol{f}(\tau)\otimes I_n\right)\boldsymbol{\phi}(\tau)\mathsf{d}\tau\right]$$

$$= \int_{-r}^{0} \boldsymbol{\phi}^\top(\tau)\left(\boldsymbol{f}^\top(\tau)N_2^\top \otimes I_n\right) \mathrm{d}\tau \left(\widetilde{\mathsf{F}} \otimes U\right) \int_{-r}^{0} \left(N_2 \boldsymbol{f}(\tau)\otimes I_n\right)\boldsymbol{\phi}(\tau)\mathrm{d}\tau$$
$$= \int_{-r}^{0} \boldsymbol{\phi}^\top(\tau)F^\top(\tau)\mathrm{d}\tau \left[\left(N_2^\top \widetilde{\mathsf{F}} N_2\right)\otimes U\right]\int_{-r}^{0} F(\tau)\boldsymbol{\phi}(\tau)\mathrm{d}\tau$$
$$= \int_{-r}^{0} \boldsymbol{\phi}^\top(\tau)\widehat{G}^\top(\tau)\mathrm{d}\tau \left[\left(G^{-1\top}N_2^\top\widetilde{\mathsf{F}}N_2G^{-1}\right)\otimes U\right]\int_{-r}^{0} \widehat{G}(\tau)\boldsymbol{\phi}(\tau)\mathrm{d}\tau \quad (4.47)$$

for any $\boldsymbol{\xi}$ and $\boldsymbol{\phi}(\cdot)$ in (4.31), where $\widetilde{\mathsf{F}}^{-1} = \int_{-r}^{0}(\tau + r)N_1\boldsymbol{f}(\tau)\boldsymbol{f}^\top(\tau)N_1^\top \mathrm{d}\tau$ with $N_1$, $N_2$ in (4.34).

Now by utilizing (4.46) and (4.47) on (4.38), we see that if (4.39) is feasible, then it guarantees the existence of (4.38) satisfying (4.36) given what is shown in (4.45). In addition, considering the property of congruent transformations knowing that $G \in \mathbb{R}_{[d]}^{d\times d}$, we can deduce that (4.39) holds if and only if

$$\left[I_n \oplus \left(G^\top \otimes I_\nu\right)\right] \widehat{P}\left[I_n \oplus (G\otimes I_\nu)\right]$$
$$+ \left[I_n \oplus \left(G^\top \otimes I_\nu\right)\right]\left[\mathsf{O}_n \oplus \left(G^{-1\top}\mathsf{F}G^{-1}\otimes S + [*]\widetilde{\mathsf{F}}\left(N_2G^{-1}\right)\otimes U\right)\right]\left[I_n \oplus (G\otimes I_\nu)\right]$$
$$= P + \left[\mathsf{O}_n \oplus \left(\mathsf{F}\otimes S + N_2^\top\widetilde{\mathsf{F}}N_2 \otimes U\right)\right] \succ 0, \quad (4.48)$$

where $P = \left[I_n \oplus \left(G^\top \otimes I_\nu\right)\right] \widehat{P}\left[I_n \oplus (G\otimes I_\nu)\right]$. Treating the matrix $P$ as a new variable, we observe that the feasibility of the last matrix inequality in (4.48) is unchanged for any $G \in \mathbb{R}_{[n]}^{n\times n}$. As a consequence, we have shown that (4.39) has the same feasibility for any invertible $G$.

Now we employ (4.38) to construct conditions which can imply (4.37). Differentiate $\mathsf{v}(\boldsymbol{x}(t), \mathfrak{y}_t(\cdot))$ along the trajectory of (4.31) at $t = t_0$ and consider the relation

$$\frac{\mathrm{d}}{\mathrm{d}t}\int_{-r}^{0}\widehat{G}(\tau)\boldsymbol{\phi}(\tau)\mathrm{d}\tau = \frac{\mathrm{d}}{\mathrm{d}t}\int_{-r}^{0}(G\otimes I_\nu)F(\tau)\boldsymbol{\phi}(\tau)\mathrm{d}\tau = (G\otimes I_\nu)F(0)\boldsymbol{\phi}(0) - (G\otimes I_\nu)F(-r)\boldsymbol{\phi}(-r)$$
$$-\widehat{M}\int_{-r}^{0}\left(G\boldsymbol{f}(\tau)\otimes I_\nu\right)\boldsymbol{\phi}(\tau)\mathrm{d}\tau = \widehat{G}(0)A_4\boldsymbol{\xi} + \left[\widehat{G}(0)A_5 - \widehat{G}(-r)\right]\boldsymbol{\phi}(-r) - \widehat{M}\int_{-r}^{0}\widehat{G}(\tau)\boldsymbol{\phi}(\tau)\mathrm{d}\tau \quad (4.49)$$

where $\widehat{M} = (G\otimes I_\nu)(M\otimes I_\nu)(G^{-1}\otimes I_\nu)$ and (4.49) can be obtained by the relation in (4.31). Then we have

$$\left.\frac{\mathrm{d}^+}{\mathrm{d}t}\mathsf{v}(\boldsymbol{x}(t), \mathfrak{y}_t(\cdot))\right|_{t=t_0, \boldsymbol{x}(t_0)=\boldsymbol{\xi}, \mathfrak{y}_{t_0}(\cdot)=\boldsymbol{\phi}(\cdot)} = \boldsymbol{\chi}^\top \mathsf{Sy}\left(H\widehat{P}\begin{bmatrix}\mathbf{A}\\ \mathbf{G}\end{bmatrix}\right)\boldsymbol{\chi}$$
$$+ \boldsymbol{\chi}^\top\left[\Gamma^\top(S+rU)\Gamma - \left(\mathsf{O}_n \oplus S \oplus \mathsf{O}_{d\nu}\right)\right]\boldsymbol{\chi} - \int_{-r}^{0}\boldsymbol{\phi}^\top(\tau)U\boldsymbol{\phi}(\tau)\mathrm{d}\tau, \quad (4.50)$$

where $H$, $\mathbf{A}$, $\mathbf{G}$ and $\Gamma$ have been defined in the statement of Theorem 4.5 and

$$\boldsymbol{\chi} := \mathbf{Col}\left(\boldsymbol{\xi},\ \boldsymbol{\phi}(-r),\ \int_{-r}^{0}\widehat{G}(\tau)\boldsymbol{\phi}(\tau)\mathrm{d}\tau\right). \quad (4.51)$$

Given $U \succ 0$ in (4.39) and apply (4.3) with $\varpi(\tau) = 1$ to the integral $\int_{-r}^{0}\boldsymbol{\phi}^\top(\tau)U\boldsymbol{\phi}(\tau)\mathrm{d}\tau$ in (4.50) analogous to the procedure in (4.46). It produces

$$\int_{-r}^{0}\boldsymbol{\phi}^\top(\tau)U\boldsymbol{\phi}(\tau)\mathrm{d}\tau \geq \int_{-r}^{0}\widehat{G}(\tau)\boldsymbol{\phi}(\tau)\mathrm{d}\tau\left[(G^{-1})^\top\mathsf{F}G^{-1}\otimes U\right]\int_{-r}^{0}\widehat{G}(\tau)\boldsymbol{\phi}(\tau)\mathrm{d}\tau \quad (4.52)$$

for any $\boldsymbol{\xi} \in \mathbb{R}^n$ and $\boldsymbol{\phi}(\cdot) \in \widehat{C}\left([-r_2, 0); \mathbb{R}^\nu\right)$ in (4.31). By applying (4.52) to (4.50), we have

$$\left.\frac{\mathrm{d}^+}{\mathrm{d}t}\mathsf{v}(\boldsymbol{x}(t), \mathfrak{y}_t(\cdot))\right|_{t=t_0, \boldsymbol{x}(t_0)=\boldsymbol{\xi}, \mathfrak{y}_{t_0}(\cdot)=\boldsymbol{\phi}(\cdot)} \leq \boldsymbol{\chi}^\top\boldsymbol{\Phi}\boldsymbol{\chi} \quad (4.53)$$

given that $U \succ 0$ in (4.39), where $\mathbf{\Phi}$ is defined in (4.41). By (4.53) and (4.51), it becomes evident that the feasible solutions to (4.40) guarantees the existence of (4.38) satisfying (4.37).

Considering the property of congruent transformations knowing that $G \in \mathbb{R}^{d\times d}_{[d]}$, we have

$$\mathbf{\Phi} \prec 0 \iff \left[I_{n+\nu} \oplus (G^\top \otimes I_\nu)\right] \mathbf{\Phi} \left[I_{n+\nu} \oplus (G \otimes I_\nu)\right] = \mathbf{\Theta} \prec 0 \tag{4.54}$$

where

$$\mathbf{\Theta} := \mathsf{Sy}\left(HP\Psi\right) + \Gamma^\top (S + rU)\Gamma - \left[\mathsf{O}_n \oplus S \oplus (\mathsf{F} \otimes U)\right] \tag{4.55}$$

with $P = \left[I_n \oplus (G^\top \otimes I_\nu)\right] \widehat{P} \left[I_n \oplus (G \otimes I_\nu)\right]$ and $H$ defined in (4.42) and

$$\Psi = \begin{bmatrix} A_1 & A_2 & A_3 \\ F(0)A_4 & F(0)A_5 - F(-r) & -M \otimes I_\nu \end{bmatrix}, \tag{4.56}$$

which is derived via (2.1a) and (4.49). By treating the matrix $P$ as a new variable, it is clear to see that the feasibility of (4.54) is invariant from $G \in \mathbb{R}^{n\times n}_{[n]}$, which indicates the feasibility of (4.40) remains unchanged for any invertible $G$.

Finally, based on the conclusion in Theorem 4.4, we know that if (4.17) or (4.23) are employed during the steps at (4.46) or (4.52) instead of (4.3), then the resulting conditions with the extra constraints induced by (4.22) has the same feasibility as (4.39) and (4.40) for any $G \in \mathbb{R}^{d\times d}_{[d]}$. Since the feasibility of (4.39) and (4.40) is invariant with respect to $G \in \mathbb{R}^{n\times n}_{[n]}$, hence the feasibility of the conditions derived via (4.17) and (4.23) is also invariant with respect to $G \in \mathbb{R}^{n\times n}_{[n]}$ $G$. This finishes the proof of this theorem. ■

By Theorem 4.5, we have established that applying a linear transformation to $\boldsymbol{f}(\cdot)$, namely $\boldsymbol{g}(\tau) = G\boldsymbol{f}(\tau)$ in (4.38), cannot change the feasibility of the stability conditions derived via (4.3) or (4.17) or (4.23). In fact, if $\boldsymbol{f}(\cdot)$ contains only orthogonal functions in (4.38) with $G = I_d$, such option cannot render the corresponding stability condition more feasible compared to the case of $G \neq I_d$ that gives a non-orthogonal structure for $\boldsymbol{g}(\cdot)$. However, although selecting orthogonal functions for $\boldsymbol{g}(\cdot)$ cannot generate less conservative stability conditions, as we have demonstrated, it may still be beneficial to do so. Specifically, the matrix $\left(\int_{-r}^{0} \boldsymbol{g}(\tau)\boldsymbol{g}^\top(\tau)\mathsf{d}\tau\right)^{-1} = (G^{-1})^\top \mathsf{F} G^{-1}$ in (4.41) is always diagonal if $\{g_i(\cdot)\}_{i=1}^{d}$ only contains mutually orthogonal functions, which might be a positive factor for numerical computations.

# Chapter 5

# Stability and Dissipativity Analysis of Linear Coupled Differential-Difference Systems with Distributed Delays

## 5.1 Introduction

Coupled differential-functional equations (CDFEs), which generalize the expressions of TDSs [452], can express a broad class of models concerning delay or propagation effects [637]. CDFEs can model systems such as standard or neutral TDSs or certain singular delay systems [638]. For more information on the topic of CDFEs, see [22, 175] and the references therein.

Over the past few decades, a series of significant results on the stability of CDESs [639, 640] has been proposed based on the approach of constructing LKFs. In particular, the idea of the complete LKF of LTDSs [452] has been extended in [22] to formulate a complete functional for a linear coupled differential-difference system (CDDS)[1] [641], which can be constructed numerically [642] via SDP. To the best of our knowledge, however, methods on linear CDDSs with nontrivial (non-constant) DDs remain underdeveloped. Generally speaking, analyzing DDs may necessitate much more effort due to the complexities induced by different types of DD kernels. For the latest time-domain-based results on the subject of DDs, see [159, 169, 419, 462, 463, 470].

In [470], an approximation scheme is introduced to address DDs with continuous kernels based on the application of Legendre polynomials [468, 22.2.10]. Despite the fact that only one or two DD kernels are considered, the stability conditions derived in [470] are notably effective and exhibit a pattern of hierarchical feasibility enhancement with respect to the degree of the Legendre polynomials utilized for the approximation. In this chapter, we propose a new approach generalizing the results in [470]. In contrast to the approximation scheme in [470] where approximations are solely carried out by the application of Legendre orthogonal polynomials, our proposed approximation solution employs a class of elementary functions such as polynomials and trigonometric functions. Our approach offers a unified framework for handling situations where multiple distributed matrix kernels are individually approximated over two distinct delay intervals with general matrix structures. Additionally, unified measures for approximation errors are formulated via a matrix framework and incorporated into our proposed stability and dissipativity condition.

In this chapter, we propose solutions to the dissipativity and stability analysis of a linear CDDS

[1] A CDDS can be considered as a special case of the systems denoted by CDFEs

with DDs in both the state and output equation. Specifically, the DD kernels considered can be any $\mathscr{L}^2$ function and all kernel functions are approximated by a class of elementary functions. Many existing models with delays, such as the ones in [22, 169, 463, 470, 642] are the special cases of the system model under consideration. Meanwhile, we also analyze the behavior of approximation errors using matrix representations, generalizing existing results in [470]. Our analysis shows that when orthogonal functions are used as approximators, the approximation errors converge to zero in an $\mathscr{L}^2$ sense as $d \to \infty$. Furthermore, a quadratic supply rate function is also considered for the dissipative analysis. Additionally, we consider a quadratic supply rate function for dissipative analysis. To incorporate approximation errors into the optimization constraints for dissipativity and stability analysis, we derive a general integral inequality that introduces error-related terms into its LB. Using this inequality to construct a Lyapunov-Krasovskiĭ functional (LKF), we derive sufficient conditions for dissipativity and asymptotic (exponential) stability in terms of LMIs. These conditions exhibit hierarchical feasibility enlargement when only orthogonal functions are used to approximate the DD kernels, offering an significant extension of the result in [469]. Finally, several numerical examples were experimented to demonstrate the effectiveness and capacity of the proposed methodologies.

## 5.2 Problem Formulation

The following linear CDDS

$$
\begin{aligned}
\dot{\boldsymbol{x}}(t) &= A_1\boldsymbol{x}(t) + A_2\boldsymbol{y}(t-r_1) + A_3\boldsymbol{y}(t-r_2) \\
&\qquad + \int_{-r_1}^{0} \widetilde{A}_4(\tau)\boldsymbol{y}(t+\tau)\mathrm{d}\tau + \int_{-r_2}^{-r_1} \widetilde{A}_5(\tau)\boldsymbol{y}(t+\tau)\mathrm{d}\tau + D_1\boldsymbol{w}(t) \\
\boldsymbol{y}(t) &= A_6\boldsymbol{x}(t) + A_7\boldsymbol{y}(t-r_1) + A_8\boldsymbol{y}(t-r_2), \qquad t \geq t_0 \\
\boldsymbol{z}(t) &= C_1\boldsymbol{x}(t) + C_2\boldsymbol{y}(t-r_1) + C_3\boldsymbol{y}(t-r_2) + \int_{-r_1}^{0} \widetilde{C}_4(\tau)\boldsymbol{y}(t+\tau)\mathrm{d}\tau \\
&\qquad + \int_{-r_2}^{-r_1} \widetilde{C}_5(\tau)\boldsymbol{y}(t+\tau)\mathrm{d}\tau + C_6\dot{\boldsymbol{y}}(t-r_1) + C_7\dot{\boldsymbol{y}}(t-r_2) + D_2\boldsymbol{w}(t) \\
\boldsymbol{x}(t_0) &= \boldsymbol{\xi} \in \mathbb{R}^n,\ \ \forall\theta \in [-r_2, 0),\ \boldsymbol{y}(t_0+\theta) = \boldsymbol{\psi}(\theta), \quad \boldsymbol{\psi}(\cdot) \in \mathcal{A}\left([-r_2,0);\mathbb{R}^\nu\right)
\end{aligned}
\tag{5.1}
$$

with DDs is considered in this chapter, where $r_2 > r_1 > 0$ and $t_0 \in \mathbb{R}$. The notation $\mathcal{A}\left([-r_2,0);\mathbb{R}^\nu\right)$ in (5.1) stands for

$$
\mathcal{A}\left([-r_2,0);\mathbb{R}^\nu\right) := \left\{\boldsymbol{\psi}(\cdot) \in \mathcal{C}\left([-r_2,0);\mathbb{R}^\nu\right) :\ \dot{\boldsymbol{\psi}}(\cdot) \in \mathscr{L}^2\left([-r_2,0);\mathbb{R}^\nu\right)\ \&\ \ \|\boldsymbol{\psi}(\cdot)\|_\infty + \|\dot{\boldsymbol{\psi}}(\cdot)\|_2 < +\infty\right\}
$$

the space of absolutely continuous functions, where $\|\boldsymbol{\psi}(\cdot)\|_\infty := \sup_{\tau\in\mathcal{X}} \|\boldsymbol{\psi}(\tau)\|_2$ and $\dot{\boldsymbol{\psi}}(\cdot)$ stands for the weak derivatives of $\boldsymbol{\psi}(\cdot)$. Furthermore, $\boldsymbol{x}(t) \in \mathbb{R}^n$, $\boldsymbol{y}(t) \in \mathbb{R}^\nu$ satisfy (5.1), and $\boldsymbol{w}(\cdot) \in \mathscr{L}^2([t_0,\infty);\mathbb{R}^q)$, $\boldsymbol{z}(t) \in \mathbb{R}^m$ are the disturbance and output of (5.1), respectively. The size of the state space matrices in (5.1) are determined by the given dimensions $n;\nu \in \mathbb{N}$ and $m;q \in \mathbb{N}_0$. All the functions in the entries of the matrix-valued DD kernels $\widetilde{A}_4(\cdot), \widetilde{C}_4(\cdot)$ and $\widetilde{A}_5(\cdot), \widetilde{C}_5(\cdot)$ are the elements of $\mathscr{L}^2\left(\mathcal{I}_1;\mathbb{R}\right)$ and $\mathscr{L}^2\left(\mathcal{I}_2;\mathbb{R}\right)$, respectively, where $\mathcal{I}_1 = [-r_1, 0]$ and $\mathcal{I}_2 = [-r_2, -r_1]$. Finally, $A_7$ and $A_8$ satisfy

$$
\sup\left\{s \in \mathbb{C} : \det\left(I_\nu - A_7\mathrm{e}^{-r_1 s} - A_8\mathrm{e}^{-r_2 s}\right) = 0\right\} < 0,
\tag{5.2}
$$

which ensures input to state stability for the associated difference equation [459] of (5.1).

In order to address the DDs in (5.1), we first define $\acute{\boldsymbol{f}}(\cdot) \in C^1(\mathcal{I}_1; \mathbb{R}^d)$ and $\grave{\boldsymbol{f}}(\cdot) \in C^1(\mathcal{I}_2; \mathbb{R}^\delta)$ which satisfy

$$\exists! M_1 \in \mathbb{R}^{d\times d},\ \ \exists! M_2 \in \mathbb{R}^{\delta\times\delta} : \frac{\mathsf{d}\acute{\boldsymbol{f}}(\tau)}{\mathsf{d}\tau} = M_1 \acute{\boldsymbol{f}}(\tau) \ \text{ and } \ \frac{\mathsf{d}\grave{\boldsymbol{f}}(\tau)}{\mathsf{d}\tau} = M_2 \grave{\boldsymbol{f}}(\tau), \tag{5.3}$$

$$\begin{aligned} &\exists \acute{\boldsymbol{\phi}}(\cdot) \in C^1(\mathcal{I}_1; \mathbb{R}^{\kappa_1}),\ \exists \grave{\boldsymbol{\phi}}(\cdot) \in C^1(\mathcal{I}_2; \mathbb{R}^{\kappa_2}),\ \exists! M_3 \in \mathbb{R}^{\kappa_1\times d},\ \ \exists! M_4 \in \mathbb{R}^{\kappa_2\times\delta} : \\ &\qquad\qquad\qquad\qquad \frac{\mathsf{d}\acute{\boldsymbol{\phi}}(\tau)}{\mathsf{d}\tau} = M_3 \acute{\boldsymbol{f}}(\tau) \ \text{ and } \ \frac{\mathsf{d}\grave{\boldsymbol{\phi}}(\tau)}{\mathsf{d}\tau} = M_4 \grave{\boldsymbol{f}}(\tau), \end{aligned} \tag{5.4}$$

$$\mathbb{S}^d \ni \acute{\mathsf{F}}_d^{-1} = \int_{-r_2}^{-r_1} \grave{\boldsymbol{f}}(\tau)\grave{\boldsymbol{f}}^\top(\tau)\mathsf{d}\tau \succ 0, \quad \mathbb{S}^\delta \ni \grave{\mathsf{F}}_\delta^{-1} = \int_{-r_2}^{-r_1} \grave{\boldsymbol{f}}(\tau)\grave{\boldsymbol{f}}^\top(\tau)\mathsf{d}\tau \succ 0, \tag{5.5}$$

$$\mathbb{S}^{\kappa_1} \ni \acute{\Phi}_{\kappa_1}^{-1} = \int_{-r_1}^{0} \acute{\boldsymbol{\phi}}(\tau)\acute{\boldsymbol{\phi}}^\top(\tau)\mathsf{d}\tau \succ 0, \quad \mathbb{S}^{\kappa_2} \ni \grave{\Phi}_{\kappa_2}^{-1} = \int_{-r_2}^{-r_1} \grave{\boldsymbol{\phi}}(\tau)\grave{\boldsymbol{\phi}}^\top(\tau)\mathsf{d}\tau \succ 0, \tag{5.6}$$

where $d; \delta \in \mathbb{N}$. The conditions in (5.6) indicate that the functions in $\acute{\boldsymbol{f}}(\cdot)$, $\grave{\boldsymbol{f}}(\cdot)$, $\acute{\boldsymbol{\phi}}(\cdot)$ and $\grave{\boldsymbol{\phi}}(\cdot)$ are linearly independent in a Lebesgue sense, respectively. See [566, Theorem 7.2.10] for further explanation on the mathematical foundation of (5.6).

**Remark 5.1.** The constraint in (5.3) specifies that the functions in $\acute{\boldsymbol{f}}(\cdot)$, $\grave{\boldsymbol{f}}(\cdot)$ are the solutions to homogeneous differential equations with constant coefficients. (polynomials, exponential, trigonometric functions, etc) Note that (5.4) does not impose additional constraints on $\acute{\boldsymbol{f}}(\cdot)$, $\grave{\boldsymbol{f}}(\cdot)$. This is because for any given $\acute{\boldsymbol{f}}(\cdot)$, $\grave{\boldsymbol{f}}(\cdot)$ satisfying (5.3), the one can always select $\acute{\boldsymbol{\phi}}(\tau) = \acute{\boldsymbol{f}}(\tau)$ and $\grave{\boldsymbol{\phi}}(\tau) = \grave{\boldsymbol{f}}(\tau)$ with $M_3 = M_1$ and $M_4 = M_2$ to satisfy (5.4).

Now given $\acute{\boldsymbol{f}}(\cdot) \in C^1(\mathcal{I}_1; \mathbb{R}^d)$ and $\grave{\boldsymbol{f}}(\cdot) \in C^1(\mathcal{I}_2; \mathbb{R}^\delta)$ satisfying (5.3), we see that for any $\widetilde{A}_4(\cdot); \widetilde{A}_5(\cdot)$ and $\widetilde{C}_4(\cdot); \widetilde{C}_5(\cdot)$ in (5.1), there exist constant matrices $A_4 \in \mathbb{R}^{n\times(d+\mu_1)\nu}$, $A_5 \in \mathbb{R}^{n\times(\delta+\mu_2)\nu}$, $C_4 \in \mathbb{R}^{m\times(d+\mu_1)\nu}$, $C_5 \in \mathbb{R}^{m\times(\delta+\mu_2)\nu}$ and the functions $\boldsymbol{\varphi}_1(\cdot) \in \mathcal{L}^2(\mathcal{I}_1; \mathbb{R}^{\mu_1})$, $\boldsymbol{\varphi}_2(\cdot) \in \mathcal{L}^2(\mathcal{I}_2; \mathbb{R}^{\mu_2})$ such that

$$\begin{aligned} \widetilde{A}_4(\tau) &= A_4 \left( \begin{bmatrix} \boldsymbol{\varphi}_1(\tau) \\ \acute{\boldsymbol{f}}(\tau) \end{bmatrix} \otimes I_\nu \right), \qquad & \widetilde{A}_5(\tau) &= A_5 \left( \begin{bmatrix} \boldsymbol{\varphi}_2(\tau) \\ \grave{\boldsymbol{f}}(\tau) \end{bmatrix} \otimes I_\nu \right), \\ \widetilde{C}_4(\tau) &= C_4 \left( \begin{bmatrix} \boldsymbol{\varphi}_1(\tau) \\ \acute{\boldsymbol{f}}(\tau) \end{bmatrix} \otimes I_\nu \right), \qquad & \widetilde{C}_5(\tau) &= C_5 \left( \begin{bmatrix} \boldsymbol{\varphi}_2(\tau) \\ \grave{\boldsymbol{f}}(\tau) \end{bmatrix} \otimes I_\nu \right), \end{aligned} \tag{5.7}$$

$$\int_{-r_1}^{0} \begin{bmatrix} \boldsymbol{\varphi}_1(\tau) \\ \acute{\boldsymbol{f}}(\tau) \end{bmatrix} \begin{bmatrix} \boldsymbol{\varphi}_1^\top(\tau) & \acute{\boldsymbol{f}}^\top(\tau) \end{bmatrix} \mathsf{d}\tau \succ 0, \qquad \int_{-r_2}^{-r_1} \begin{bmatrix} \boldsymbol{\varphi}_2(\tau) \\ \grave{\boldsymbol{f}}(\tau) \end{bmatrix} \begin{bmatrix} \boldsymbol{\varphi}_2^\top(\tau) & \grave{\boldsymbol{f}}^\top(\tau) \end{bmatrix} \mathsf{d}\tau \succ 0, \tag{5.8}$$

where $\mu_1, \mu_2 \in \mathbb{N}_0$ and (5.8) indicates that the functions in $\mathbf{Col}\left[\boldsymbol{\varphi}_1(\tau), \acute{\boldsymbol{f}}(\tau)\right]$ and $\mathbf{Col}\left[\boldsymbol{\varphi}_2(\tau), \grave{\boldsymbol{f}}(\tau)\right]$ are linearly independent in a Lebesgue sense, respectively. Thus, (5.7) can be applied to equivalently denote the DD kernels in (5.1). Finally, note that (5.5) is satisfied if (5.8) holds.

**Remark 5.2.** The functions inside $\acute{\boldsymbol{f}}(\cdot)$ and $\grave{\boldsymbol{f}}(\cdot)$ in (5.7) are selected in view of the kernel functions in $\widetilde{A}_4(\cdot)$, $\widetilde{A}_5(\cdot)$, $\widetilde{C}_4(\cdot)$ and $\widetilde{C}_5(\cdot)$. Note that we can always let $\acute{\boldsymbol{f}}(\cdot)$ and $\grave{\boldsymbol{f}}(\cdot)$ to only contain orthogonal functions since one can always adjust the elements in $\boldsymbol{\varphi}_1(\cdot) \in \mathcal{L}^2(\mathcal{I}_1; \mathbb{R}^{\mu_1})$ and $\boldsymbol{\varphi}_2(\cdot) \in \mathcal{L}^2(\mathcal{I}_2; \mathbb{R}^{\mu_2})$ to satisfy (5.7). Note that $\boldsymbol{\varphi}_1(\cdot)$ and $\boldsymbol{\varphi}_2(\cdot)$ can become a $0 \times 1$ empty vector if $\mu_1 = \mu_2 = 0$. Finally, the matrix inequalities in (5.8) can be verified numerically[2] with given $\acute{\boldsymbol{f}}(\cdot)$, $\grave{\boldsymbol{f}}(\cdot)$ and $\boldsymbol{\varphi}_1(\cdot), \boldsymbol{\varphi}_2(\cdot)$.

**Remark 5.3.** It is worth remarking that (5.1) generalizes all the models in [169, 462, 470] without considering uncertainties. Similar to the strategies outlined in the previous chapters, the decomposition in (5.7) is employed to handle the DD in (5.1), enabling us to derive well-posed dissipativity and stability conditions. On the other hand, we will employ an approximation scheme in this chapter to deal with functions $\boldsymbol{\varphi}_1(\cdot), \boldsymbol{\varphi}_2(\cdot)$, which is fundamentally different from the solutions in the previous chapters.

[2] One option is to use `vpaintegral` in MATLAB© which performs high-precision numerical integration.

**Remark 5.4.** A neutral delay system

$$\frac{\mathsf{d}}{\mathsf{d}t}\left(\boldsymbol{y}(t) - A_4\boldsymbol{y}(t-r)\right) = A_1\boldsymbol{y}(t) + A_2\boldsymbol{y}(t-r) + \int_{-r}^{0} A_3(\tau)\boldsymbol{x}(t+\tau)\mathsf{d}\tau$$

can be equivalently expressed by CDDS

$$\begin{aligned}
\dot{\boldsymbol{x}}(t) &= A_1\boldsymbol{x}(t) + (A_2 + A_1A_4)\boldsymbol{y}(t-r) + \int_{-r}^{0} A_3(\tau)\boldsymbol{x}(t+\tau)\mathsf{d}\tau \\
\boldsymbol{y}(t) &= \boldsymbol{x}(t) + A_4\boldsymbol{y}(t-r).
\end{aligned}$$

On the other hand, if there is rank redundancy in the delay matrices, namely,

$$\frac{\mathsf{d}}{\mathsf{d}t}\left(\boldsymbol{y}(t) - A_4N\boldsymbol{y}(t-r)\right) = A_1\boldsymbol{y}(t) + A_2N\boldsymbol{y}(t-r) + \int_{-r}^{0} A_3(\tau)N\boldsymbol{y}(t+\tau)\mathsf{d}\tau, \tag{5.9}$$

then we can first reformulate (5.9) as

$$\frac{\mathsf{d}}{\mathsf{d}t}\left(\boldsymbol{y}(t) - A_4\boldsymbol{z}(t-r)\right) = A_1\boldsymbol{y}(t) + A_2\boldsymbol{z}(t-r) + \int_{-r}^{0} A_3(\tau)\boldsymbol{z}(t+\tau)\mathsf{d}\tau, \quad \boldsymbol{z}(t) = N\boldsymbol{y}(t). \tag{5.10}$$

Then, by setting $\boldsymbol{x}(t) = \boldsymbol{y}(t) - A_4\boldsymbol{z}(t-r)$ considering (5.10), we can obtain the equivalent CDDS representation

$$\begin{aligned}
\dot{\boldsymbol{x}}(t) &= A_1\boldsymbol{x}(t) + \left(A_1A_4 + A_2\right)\boldsymbol{z}(t-r) + \int_{-r}^{0} A_3(\tau)\boldsymbol{z}(t+\tau)\mathsf{d}\tau \\
\boldsymbol{z}(t) &= N\boldsymbol{x}(t) + NA_4\boldsymbol{z}(t-r).
\end{aligned}$$

This representation is advantageous in terms of reducing dimensionality if $\dim\left[\boldsymbol{z}(t)\right] \ll \dim\left[\boldsymbol{y}(t)\right]$. For further details on exploiting the rank redundancies among state-space variables in retarded cases, see [22].

In this chapter, the functions $\acute{\boldsymbol{f}}(\cdot)$ and $\grave{\boldsymbol{f}}(\cdot)$ in (5.7) are applied to approximate the functions $\boldsymbol{\varphi}_1(\cdot) \in \mathscr{L}^2(\mathcal{I}_1;\mathbb{R}^{\mu_1})$ and $\boldsymbol{\varphi}_2(\cdot) \in \mathscr{L}^2(\mathcal{I}_2;\mathbb{R}^{\mu_2})$ in (5.7), respectively, where $\boldsymbol{\varphi}_1(\cdot)$ and $\boldsymbol{\varphi}_2(\cdot)$ might not satisfy (5.3). Specifically, the approximations are denoted by the relations

$$\boldsymbol{\varphi}_1(\tau) = \acute{\Gamma}_d\acute{\boldsymbol{f}}(\tau) + \acute{\boldsymbol{\varepsilon}}_d(\tau), \quad \boldsymbol{\varphi}_2(\tau) = \grave{\Gamma}_\delta\grave{\boldsymbol{f}}(\tau) + \grave{\boldsymbol{\varepsilon}}_\delta(\tau) \tag{5.11}$$

where $\acute{\Gamma}_d$ and $\grave{\Gamma}_\delta$ are given coefficients. Furthermore, $\acute{\boldsymbol{\varepsilon}}_d(\tau) = \boldsymbol{\varphi}_1(\tau) - \acute{\Gamma}_d\acute{\boldsymbol{f}}(\tau)$ and $\grave{\boldsymbol{\varepsilon}}_\delta(\tau) = \boldsymbol{\varphi}_2(\tau) - \grave{\Gamma}_\delta\grave{\boldsymbol{f}}(\tau)$ contain the errors of approximations. In addition, we define matrices

$$\mathbb{S}^{\mu_1} \ni \acute{\mathsf{E}}_d := \int_{-r_1}^{0} \acute{\boldsymbol{\varepsilon}}_d(\tau)\acute{\boldsymbol{\varepsilon}}_d^\top(\tau)\mathsf{d}\tau, \quad \mathbb{S}^{\mu_2} \ni \grave{\mathsf{E}}_\delta := \int_{-r_2}^{-r_1} \grave{\boldsymbol{\varepsilon}}_\delta(\tau)\grave{\boldsymbol{\varepsilon}}_\delta^\top(\tau)\mathsf{d}\tau \tag{5.12}$$

to measure the error residues of (5.11). Inspired by the idea of orthogonal approximation in Hilbert space [643], one option for the values of $\acute{\Gamma}_d$ and $\grave{\Gamma}_\delta$ in (5.11) is to use the least squares method

$$\begin{aligned}
\mathbb{R}^{\mu_1\times d} \ni \acute{\Gamma}_d &:= \int_{-r_1}^{0} \boldsymbol{\varphi}_1(\tau)\acute{\boldsymbol{f}}^\top(\tau)\mathsf{d}\tau\acute{\mathsf{F}}_d, \; \acute{\mathsf{F}}_d^{-1} = \int_{-r_2}^{-r_1} \acute{\boldsymbol{f}}(\tau)\acute{\boldsymbol{f}}^\top(\tau)\mathsf{d}\tau \\
\mathbb{R}^{\mu_2\times \delta} \ni \grave{\Gamma}_\delta &:= \int_{-r_2}^{-r_1} \boldsymbol{\varphi}_2(\tau)\grave{\boldsymbol{f}}^\top(\tau)\mathsf{d}\tau\grave{\mathsf{F}}_\delta, \; \grave{\mathsf{F}}_\delta^{-1} = \int_{-r_2}^{-r_1} \grave{\boldsymbol{f}}(\tau)\grave{\boldsymbol{f}}^\top(\tau)\mathsf{d}\tau.
\end{aligned} \tag{5.13}$$

**Remark 5.5.** The expressions in (5.13) can be interpreted as a vector form of the standard approximations (Least Squares) in Hilbert space, as detailed in section 10.2 in [643]) Specifically, let $\{\acute{f}_i(\cdot)\}_{i=1}^{\infty} \subset \mathcal{L}^2(\mathcal{I}_1; \mathbb{R})$ and $\{\grave{f}_i(\cdot)\}_{i=1}^{\infty} \subset \mathcal{L}^2(\mathcal{I}_2; \mathbb{R})$ be the orthogonal bases (not necessary normal here) of the corresponding Hilbert space and satisfying (5.3) (Trigonometric series, for instance). Then one can conclude that by using (5.13) with $\{\acute{f}_i(\cdot)\}_{i=1}^{d}$ and $\{\grave{f}_i(\cdot)\}_{i=1}^{\delta}$, each entry of the error vectors $\acute{\boldsymbol{\varepsilon}}_d(\cdot)$ and $\grave{\boldsymbol{\varepsilon}}_\delta(\cdot)$ in (5.12) converges to zero in a $\mathcal{L}^2$ sense, respectively, as $d \to \infty$ and $\delta \to \infty$. If $\acute{\boldsymbol{f}}(\cdot)$ and $\grave{\boldsymbol{f}}(\cdot)$ in (5.13) contains only Legendre polynomials, then (5.11)–(5.13) generalizes the polynomials approximation scheme proposed in [470] via a matrix framework. Finally, it is very crucial to stress that (5.11) does not restrict one only to apply (5.13) for the values of $\acute{\Gamma}_d$ and $\grave{\Gamma}_\delta$. Other appropriate options for $\acute{\Gamma}_d$ and $\grave{\Gamma}_\delta$ can also be considered based on specific contexts.

The system (5.1) can be rewritten as

$$
\begin{aligned}
&\dot{\boldsymbol{x}}(t) = \mathbf{A}\boldsymbol{\vartheta}(t), \quad \boldsymbol{y}(t) = \begin{bmatrix} \mathsf{O}_{\nu,(2\nu+q)} & \Xi & \mathsf{O}_{\nu,\nu\mu} \end{bmatrix} \boldsymbol{\vartheta}(t), \quad \boldsymbol{z}(t) = \Sigma\boldsymbol{\vartheta}(t) \\
&\boldsymbol{x}(t_0) = \boldsymbol{\xi} \in \mathbb{R}^n, \quad \forall \theta \in [-r_2, 0], \; \boldsymbol{y}(t_0 + \theta) = \boldsymbol{\psi}(\theta)
\end{aligned}
\tag{5.14}
$$

with

$$
\begin{aligned}
\mathbf{A} = \Bigg[ \mathsf{O}_{n,2\nu} \quad D_1 \quad A_1 \quad A_2 \quad A_3 \quad A_4 \left( \begin{bmatrix} \acute{\Gamma}_d \\ I_d \end{bmatrix} \otimes I_\nu \right) \quad A_5 \left( \begin{bmatrix} \grave{\Gamma}_\delta \\ I_\delta \end{bmatrix} \otimes I_\nu \right) \cdots \\
\cdots A_4 \left( \begin{bmatrix} \acute{\mathsf{E}}_d \\ \mathsf{O}_{d,\mu_1} \end{bmatrix} \otimes I_\nu \right) \quad A_5 \left( \begin{bmatrix} \grave{\mathsf{E}}_\delta \\ \mathsf{O}_{\delta,\mu_2} \end{bmatrix} \otimes I_\nu \right) \Bigg]
\end{aligned}
\tag{5.15}
$$

$$
\Xi = \begin{bmatrix} A_6 & A_7 & A_8 & \mathsf{O}_{\nu,\varrho\nu} \end{bmatrix}
\tag{5.16}
$$

$$
\begin{aligned}
\Sigma = \Bigg[ C_6 \quad C_7 \quad D_2 \quad C_1 \quad C_2 \quad C_3 \quad C_4 \left( \begin{bmatrix} \acute{\Gamma}_d \\ I_d \end{bmatrix} \otimes I_\nu \right) \quad C_5 \left( \begin{bmatrix} \grave{\Gamma}_\delta \\ I_\delta \end{bmatrix} \otimes I_\nu \right) \cdots \\
\cdots C_4 \left( \begin{bmatrix} \acute{\mathsf{E}}_d \\ \mathsf{O}_{d,\mu_1} \end{bmatrix} \otimes I_\nu \right) \quad C_5 \left( \begin{bmatrix} \grave{\mathsf{E}}_\delta \\ \mathsf{O}_{\delta,\mu_2} \end{bmatrix} \otimes I_\nu \right) \Bigg]
\end{aligned}
\tag{5.17}
$$

$$
\boldsymbol{\vartheta}(t) := \mathbf{Col}\left( \begin{bmatrix} \dot{\boldsymbol{y}}(t-r_1) \\ \dot{\boldsymbol{y}}(t-r_2) \end{bmatrix}, \begin{bmatrix} \boldsymbol{w}(t) \\ \boldsymbol{x}(t) \end{bmatrix}, \begin{bmatrix} \boldsymbol{y}(t-r_1) \\ \boldsymbol{y}(t-r_2) \end{bmatrix}, \begin{bmatrix} \int_{-r_1}^{0} \acute{F}_d(\tau)\boldsymbol{y}(t+\tau)\mathsf{d}\tau \\ \int_{-r_2}^{-r_1} \grave{F}_\delta(\tau)\boldsymbol{y}(t+\tau)\mathsf{d}\tau \end{bmatrix}, \begin{bmatrix} \int_{-r_1}^{0} \acute{E}_d(\tau)\boldsymbol{y}(t+\tau)\mathsf{d}\tau \\ \int_{-r_2}^{-r_1} \grave{E}_\delta(\tau)\boldsymbol{y}(t+\tau)\mathsf{d}\tau \end{bmatrix} \right),
\tag{5.18}
$$

where $\mathbb{R}^{d\nu\times\nu} \ni \acute{F}_d(\tau) := \acute{\boldsymbol{f}}(\tau) \otimes I_\nu$ and $\mathbb{R}^{\delta\nu\times\nu} \ni \grave{F}_\delta(\tau) := \grave{\boldsymbol{f}}(\tau) \otimes I_\nu$ and $\acute{E}_d(\tau) := \acute{\mathsf{E}}_d^{-1}\acute{\boldsymbol{\varepsilon}}_d(\tau) \otimes I_\nu$ and $\grave{E}_\delta(\tau) := \grave{\mathsf{E}}_\delta^{-1}\grave{\boldsymbol{\varepsilon}}_\delta(\tau) \otimes I_\nu$ with $\acute{\mathsf{E}}_d$ and $\grave{\mathsf{E}}_\delta$ in (5.12). Note that $\acute{\mathsf{E}}_d$ and $\grave{\mathsf{E}}_\delta$ in (5.12) are invertible according to what is enunciated in Remark 5.8 based on (5.23) and (D.1). Note that also the DD in (5.14) are derived based on the identities

$$
\begin{aligned}
&\left( \begin{bmatrix} \boldsymbol{\varphi}_1(\tau) \\ \acute{\boldsymbol{f}}(\tau) \end{bmatrix} \otimes I_\nu \right) \boldsymbol{y}(t+\tau) = \left( \begin{bmatrix} \acute{\Gamma}_d \acute{\boldsymbol{f}}(\tau) + \acute{\boldsymbol{\varepsilon}}_d(\tau) \\ \acute{\boldsymbol{f}}(\tau) \end{bmatrix} \otimes I_\nu \right) \boldsymbol{y}(t+\tau) = \left( \begin{bmatrix} \acute{\Gamma}_d \\ I_d \end{bmatrix} \acute{\boldsymbol{f}}(\tau) \otimes I_\nu \right) \boldsymbol{y}(t+\tau) \\
&\quad + \left( \begin{bmatrix} I_{\mu_1} \\ \mathsf{O}_{d,\mu_1} \end{bmatrix} \acute{\boldsymbol{\varepsilon}}_d(\tau) \otimes I_\nu \right) \boldsymbol{y}(t+\tau) = \left( \begin{bmatrix} \acute{\Gamma}_d \\ I_d \end{bmatrix} \otimes I_\nu \right) \acute{F}_d(\tau)\boldsymbol{y}(t+\tau) + \left( \begin{bmatrix} \acute{\mathsf{E}}_d \\ \mathsf{O}_{d,\mu_1} \end{bmatrix} \otimes I_\nu \right) \acute{E}_d(\tau)\boldsymbol{y}(t+\tau) \\
&\left( \begin{bmatrix} \boldsymbol{\varphi}_2(\tau) \\ \grave{\boldsymbol{f}}(\tau) \end{bmatrix} \otimes I_\nu \right) \boldsymbol{y}(t+\tau) = \left( \begin{bmatrix} \grave{\Gamma}_\delta \grave{\boldsymbol{f}}(\tau) + \grave{\boldsymbol{\varepsilon}}_\delta(\tau) \\ \grave{\boldsymbol{f}}(\tau) \end{bmatrix} \otimes I_\nu \right) \boldsymbol{y}(t+\tau) = \left( \begin{bmatrix} \grave{\Gamma}_\delta \\ I_d \end{bmatrix} \grave{\boldsymbol{f}}(\tau) \otimes I_\nu \right) \boldsymbol{y}(t+\tau) \\
&\quad + \left( \begin{bmatrix} I_{\mu_2} \\ \mathsf{O}_{\delta,\mu_2} \end{bmatrix} \grave{\boldsymbol{\varepsilon}}_\delta(\tau) \otimes I_\nu \right) \boldsymbol{y}(t+\tau) = \left( \begin{bmatrix} \grave{\Gamma}_\delta \\ I_\delta \end{bmatrix} \otimes I_\nu \right) \grave{F}_\delta(\tau)\boldsymbol{y}(t+\tau) + \left( \begin{bmatrix} \grave{\mathsf{E}}_\delta \\ \mathsf{O}_{\delta,\mu_2} \end{bmatrix} \otimes I_\nu \right) \grave{E}_\delta(\tau)\boldsymbol{y}(t+\tau)
\end{aligned}
$$

which themselves are obtained via the property of the Kronecker product in (2.1a).

## 5.3 Mathematical Preliminaries

We present some essential lemmas and definitions that are necessary for subsequent derivations. These include a novel integral inequality that will be applied later to derive our dissipative stability condition.

The following lemma provides sufficient conditions that can verify the stability of (5.1). It can be interpreted as a special case of [22, Theorem 3] with certain modifications.

**Lemma 5.1.** *Given $r_2 \geq r_1 > 0$, the origin of system (5.1) with $\boldsymbol{w}(t) \equiv \mathbf{0}_q$ is globally uniformly asymptotically stable if there exist $\epsilon_1; \epsilon_2; \epsilon_3 > 0$ and a differentiable functional $\mathsf{v} : \mathbb{R}^n \times \mathcal{A}([-r_2, 0); \mathbb{R}^\nu) \to \mathbb{R}_{\geq 0}$ such that $\mathsf{v}(\mathbf{0}_n, \mathbf{0}_\nu) = 0$ and*

$$\epsilon_1 \left\|\boldsymbol{\xi}\right\|_2^2 \leq \mathsf{v}(\boldsymbol{\xi}, \boldsymbol{\psi}(\cdot)) \leq \epsilon_2 \left[\left\|\boldsymbol{\xi}\right\|_2 \vee \left(\left\|\boldsymbol{\psi}(\cdot)\right\|_\infty + \|\dot{\boldsymbol{\psi}}(\cdot)\|_2\right)\right]^2, \tag{5.19}$$

$$\dot{\mathsf{v}}(r, \boldsymbol{\xi}, \boldsymbol{\psi}(\cdot)) := \left.\frac{\mathsf{d}^+}{\mathsf{d}t}\mathsf{v}(\boldsymbol{x}(t), \mathbf{y}_t(\cdot))\right|_{t=t_0, \boldsymbol{x}(t_0)=\boldsymbol{\xi}, \mathbf{y}_{t_0}(\cdot)=\boldsymbol{\psi}(\cdot)} \leq -\epsilon_3 \left\|\boldsymbol{\xi}\right\|_2^2 \tag{5.20}$$

*for any $\boldsymbol{\xi} \in \mathbb{R}^n$ and $\boldsymbol{\psi}(\cdot) \in \mathcal{A}([-r_2, 0); \mathbb{R}^\nu)$ in (5.1), where $t_0 \in \mathbb{R}$ and $\frac{\mathsf{d}^+}{\mathsf{d}x} f(x) = \limsup_{\eta \downarrow 0} \frac{f(x+\eta) - f(x)}{\eta}$. Furthermore, $\mathbf{y}_t(\cdot)$ in (5.20) is defined by the equality $\forall t \geq t_0$, $\forall \theta \in [-r, 0)$, $\mathbf{y}_t(\theta) = \boldsymbol{y}(t + \theta)$ where $\boldsymbol{x}(t)$ and $\boldsymbol{y}(t)$ here satisfying (5.1) with $\boldsymbol{w}(t) \equiv \mathbf{0}_q$.*

**Definition 5.1** (Dissipativity)**.** Given $r_2 \geq r_1 > 0$, the coupled differential functional system (5.1) with a supply rate function $\mathsf{s}(\boldsymbol{z}(t), \boldsymbol{w}(t))$ is dissipative if there exists a differentiable functional $\mathsf{v} : \mathbb{R}^n \times \mathcal{A}([-r_2, 0); \mathbb{R}^\nu) \to \mathbb{R}$ such that

$$\forall t \geq t_0 : \ \dot{\mathsf{v}}(\boldsymbol{x}(t), \mathbf{y}_t(\cdot)) - \mathsf{s}(\boldsymbol{z}(t), \boldsymbol{w}(t)) \leq 0 \tag{5.21}$$

with $t_0 \in \mathbb{R}$ in (5.1), where $\mathbf{y}_t(\cdot)$ is defined by the equality $\forall t \geq t_0$, $\forall \theta \in [-r_2, 0)$, $\mathbf{y}_t(\theta) = \boldsymbol{y}(t + \theta)$, and $\boldsymbol{x}(t)$, $\boldsymbol{y}(t)$ and $\boldsymbol{z}(t)$ satisfy the equalities in (5.1) with $\boldsymbol{w}(\cdot) \in \mathscr{L}^2([t_0, \infty); \mathbb{R}^q)$.

To incorporate performance objectives (dissipativity) into the analysis of (5.14), we still apply a quadratic form

$$\mathsf{s}(\boldsymbol{z}(t), \boldsymbol{w}(t)) = \begin{bmatrix} \boldsymbol{z}(t) \\ \boldsymbol{w}(t) \end{bmatrix}^\top \mathbf{J} \begin{bmatrix} \boldsymbol{z}(t) \\ \boldsymbol{w}(t) \end{bmatrix}, \quad \mathbf{J} = \begin{bmatrix} J_1 & J_2 \\ * & J_3 \end{bmatrix} \in \mathbb{S}^{m+q}, \quad \widetilde{J}^\top J_1^{-1} \widetilde{J} \preceq 0, \quad J_1^{-1} \prec 0, \quad \widetilde{J} \in \mathbb{R}^{m \times m} \tag{5.22}$$

just as in (2.10) for the supply rate function.

The following generalized new integral inequality is proposed for the derivation of our major results on the dissipativity and stability analysis in this section.

**Lemma 5.2.** *Given $\mathcal{K}$ and $\varpi(\cdot)$ in (4.1) and $U \in \mathbb{S}^n_{\succeq 0} := \{X \in \mathbb{S}^n : X \succeq 0\}$ with $n \in \mathbb{N}$. Let $\mathsf{f}(\cdot) := \mathbf{Col}_{i=1}^d \mathsf{f}_i(\cdot) \in \mathscr{L}^2_\varpi(\mathcal{K}; \mathbb{R}^d)$ and $\mathsf{g}(\cdot) := \mathbf{Col}_{i=1}^\delta \mathsf{g}_i(\cdot) \in \mathscr{L}^2_\varpi(\mathcal{K}; \mathbb{R}^\delta)$ with $d \in \mathbb{N}$ and $\delta \in \mathbb{N}_0$, in which $\mathscr{L}^2_\varpi(\mathcal{K}; \mathbb{R}^d)$ is in line with the definition in (4.1) and the functions $\mathsf{f}(\cdot)$ and $\mathsf{g}(\cdot)$ satisfy*

$$\int_{\mathcal{K}} \varpi(\tau) \begin{bmatrix} \mathsf{g}(\tau) \\ \mathsf{f}(\tau) \end{bmatrix} \begin{bmatrix} \mathsf{g}^\top(\tau) & \mathsf{f}^\top(\tau) \end{bmatrix} \mathsf{d}\tau \succ 0. \tag{5.23}$$

*Then we have,*

$$\begin{aligned}\int_{\mathcal{K}} \varpi(\tau) \boldsymbol{x}^\top(\tau) U \boldsymbol{x}(\tau) \mathsf{d}\tau \geq & \int_{\mathcal{K}} \varpi(\tau) \boldsymbol{x}^\top(\tau) \mathsf{F}^\top(\tau) \mathsf{d}\tau \left(\mathcal{F}_d \otimes U\right) \int_{\mathcal{K}} \varpi(\tau) \mathsf{F}(\tau) \boldsymbol{x}(\tau) \mathsf{d}\tau \\ & + \int_{\mathcal{K}} \varpi(\tau) \boldsymbol{x}^\top(\tau) \mathsf{E}^\top(\tau) \mathsf{d}\tau \left(\mathcal{E}_d^{-1} \otimes U\right) \int_{\mathcal{K}} \varpi(\tau) \mathsf{E}(\tau) \boldsymbol{x}(\tau) \mathsf{d}\tau \end{aligned} \tag{5.24}$$

*for all* $\boldsymbol{x}(\cdot) \in \mathcal{L}^2_{\varpi}(\mathcal{K}; \mathbb{R}^n)$ *, where*

$$
\begin{aligned}
\mathsf{F}(\tau) &= \mathbf{f}(\tau) \otimes I_n \in \mathbb{R}^{dn \times n}, \quad \mathcal{F}_d^{-1} = \int_{\mathcal{K}} \varpi(\tau)\mathbf{f}(\tau)\mathbf{f}^\top(\tau)\mathsf{d}\tau \in \mathbb{S}^d_{\succ 0}, \\
\mathsf{E}(\tau) &= \boldsymbol{e}(\tau) \otimes I_n \in \mathbb{R}^{\delta n \times n}, \quad \mathcal{E}_d = \int_{\mathcal{K}} \varpi(\tau)\boldsymbol{e}(\tau)\boldsymbol{e}^\top(\tau)\mathsf{d}\tau \in \mathbb{S}^\delta, \\
\boldsymbol{e}(\tau) &= \mathbf{g}(\tau) - \mathsf{A}\mathbf{f}(\tau) \in \mathbb{R}^\delta, \quad \mathsf{A} = \int_{\mathcal{K}} \varpi(\tau)\mathbf{g}(\tau)\mathbf{f}^\top(\tau)\mathsf{d}\tau\mathcal{F}_d \in \mathbb{R}^{\delta \times d}.
\end{aligned} \tag{5.25}
$$

*Proof.* See Appendix D for details. ■

**Remark 5.6.** By [566, Theorem 7.2.10] and seeing that

$$
\left(\mathcal{L}^2_{\varpi}\left(\mathcal{K}; \mathbb{R}^d\right) / \operatorname{Ker}\left(\|\bullet\|_{\varpi}\right), \quad \int_{\mathcal{K}} \varpi(\tau)\bullet_1^\top(\tau)\bullet_2(\tau)\mathsf{d}\tau\right)
$$

constitutes an inner product space[3], it can be deduced from (5.23) that the functions in $\mathbf{f}(\cdot)$ and $\mathbf{g}(\cdot)$ are linearly independent in a Lebesgue sense.

By setting $\delta = 0$ in Lemma 5.2 and utilizing the concept of empty matrices, the following inequality can be derived. Furthermore, we show that the subsequent corollary is equivalent to Theorem 4.1 in Chapter 4.

**Corollary 5.1.** *Given* $\mathcal{K}$ *and* $\varpi(\cdot)$ *in* (4.1) *and* $U \in \mathbb{S}^n_{\succeq 0} := \{X \in \mathbb{S}^n : X \succeq 0\}$ *with* $n \in \mathbb{N}$*. Let* $\mathbf{f}(\cdot) := \mathbf{Col}_{i=1}^d \mathsf{f}_i(\cdot) \in \mathcal{L}^2_{\varpi}\left(\mathcal{K}; \mathbb{R}^d\right)$ *with* $d \in \mathbb{N}$ *where* $\mathbf{f}(\cdot)$ *satisfies*

$$
\int_{\mathcal{K}} \varpi(\tau)\mathbf{f}(\tau)\mathbf{f}^\top(\tau)\mathsf{d}\tau \succ 0. \tag{5.26}
$$

*Then the inequality*

$$
\int_{\mathcal{K}} \varpi(\tau)\boldsymbol{x}^\top(\tau)U\boldsymbol{x}(\tau)\mathsf{d}\tau \geq \int_{\mathcal{K}} \varpi(\tau)\boldsymbol{x}^\top(\tau)\mathsf{F}^\top(\tau)\mathsf{d}\tau\left(\mathcal{F}_d \otimes U\right)\int_{\mathcal{K}} \varpi(\tau)\mathsf{F}(\tau)\boldsymbol{x}(\tau)\mathsf{d}\tau \tag{5.27}
$$

*holds for all* $\boldsymbol{x}(\cdot) \in \mathcal{L}^2_{\varpi}(\mathcal{K}; \mathbb{R}^n)$ *, where* $\mathsf{F}(\tau) = \mathbf{f}(\tau) \otimes I_n \in \mathbb{R}^{dn \times n}$ *and* $\mathcal{F}_d^{-1} = \int_{\mathcal{K}} \varpi(\tau)\mathbf{f}(\tau)\mathbf{f}^\top(\tau)\mathsf{d}\tau \in \mathbb{S}^d_{\succ 0}$*.*

**Remark 5.7.** The inequality in (5.24) reduces to [470, Lemma 1] if $\mathbf{f}(\cdot)$ is composed solely of Legendre polynomials [468, 22.2.10]. Moreover, all the particular cases of (4.3) in Chapter 4 are the special cases of (5.24), as Corollary 5.1 is equivalent to Theorem 4.1.

**Remark 5.8.** In (5.24), $\mathbf{f}(\cdot)$ can be interpreted as to approximate $\mathbf{g}(\cdot)$. By setting $\mathbf{f}(\tau) = \acute{\boldsymbol{f}}(\tau)$ and $\mathbf{g}(\tau) = \boldsymbol{\varphi}_1(\tau)$ with $\varpi(\tau) = 1$ in Lemma 5.2, it follows that $\mathcal{E}_d = \acute{\mathsf{E}}_d$, where the matrix $\acute{\mathsf{E}}_d$ is defined in (5.12). Similar procedures can be applied with $\mathbf{f}(\tau) = \grave{\boldsymbol{f}}(\tau)$ and $\mathbf{g}(\tau) = \boldsymbol{\varphi}_2(\tau)$ and $\varpi(\tau) = 1$. Furthermore, if $\mathbf{f}(\cdot)$ consists only of functions that are orthogonal with respect to $\varpi(\cdot)$, then the behavior of $\mathcal{E}_d$ can be quantitatively characterized with respect to $d$, as elaborated in the following corollary. Note that (D.1) holds for any $\mathsf{A} \in \mathbb{R}^{\delta \times d}$ as long as (5.23) is satisfied, even if $\mathsf{A}$ is not defined as $\mathsf{A} = \int_{\mathcal{K}} \varpi(\tau)\mathbf{g}(\tau)\mathbf{f}^\top(\tau)\mathsf{d}\tau\,\mathcal{F}_d \in \mathbb{R}^{\delta \times d}$. This is a noteworthy conclusion as it shows that the error matrices $\acute{\mathsf{E}}_d$ and $\grave{\mathsf{E}}_\delta$ in (5.12) are invertible since (5.8) holds.

An interesting corollary of Lemma 5.2 is presented as follows which can be interpreted as a generalization of [470, Lemma 1].

[3] $\operatorname{Ker}\left(\|\bullet\|_{\varpi}\right) := \left\{\boldsymbol{f}(\cdot) \in \mathcal{L}^2_{\varpi}\left(\mathcal{K}; \mathbb{R}^d\right) : \|\boldsymbol{f}(\cdot)\|_{\varpi} = \mathbf{0}_d\right\}$

**Corollary 5.2.** *Given all the parameters defined in Lemma 5.2 with* $\{\mathsf{f}_i(\tau)\}_{i=1}^{\infty}$ *and* $\mathbf{f}(\cdot) = \mathbf{Col}_{i=1}^{d}\,\mathsf{f}_i(\cdot)$ *satisfying*

$$\forall d \in \mathbb{N},\quad \mathcal{F}_d^{-1} = \int_{\mathcal{K}} \mathbf{f}(\tau)\mathbf{f}^{\top}(\tau)\mathrm{d}\tau = \bigoplus_{i=1}^{d}\left(\int_{\mathcal{K}} \varpi(\tau)\mathsf{f}_i^2(\tau)\mathrm{d}\tau\right), \tag{5.28}$$

*then we have that*

$$\forall d \in \mathbb{N},\quad 0 \prec \mathcal{E}_{d+1} = \mathcal{E}_d - \left(\int_{\mathcal{K}} \varpi(\tau)\mathsf{f}_{d+1}^2(\tau)\mathrm{d}\tau\right)\boldsymbol{\alpha}_{d+1}\boldsymbol{\alpha}_{d+1}^{\top} \preceq \mathcal{E}_d \tag{5.29}$$

*where* $\mathcal{E}_d$ *is given in Lemma 5.2 and* $\boldsymbol{\alpha}_{d+1} := \left(\int_{\mathcal{K}} \varpi(\tau)\mathbf{g}(\tau)\mathsf{f}_{d+1}(\tau)\mathrm{d}\tau\right)\left(\int_{\mathcal{K}} \varpi(\tau)\mathsf{f}_{d+1}^2(\tau)\mathrm{d}\tau\right)^{-1} \in \mathbb{R}^{\delta}$ *and* $\mathsf{f}_{d+1}(\cdot) \in \mathcal{L}^2_{\varpi}(\mathcal{K};\mathbb{R})$.

*Proof.* Note that only the dimension of $\mathbf{f}(\cdot)$ is related to $d$, whereas $\delta$ as the dimension of $\mathbf{g}(\cdot)$ is independent of $d$. It is obvious to see that given $\mathbf{f}(\cdot)$ satisfying (5.28), we have $\mathcal{F}_{d+1} = \mathcal{F}_d \oplus \left(\int_{\mathcal{K}} \varpi(\tau)\mathsf{f}_{d+1}^2(\tau)\mathrm{d}\tau\right)^{-1}$ (see the Definition 1 in [572]). By using this property, it follows that

$$\begin{aligned}
\mathbf{e}_{d+1}(\tau) &= \mathbf{g}(\tau) - \left(\int_{\mathcal{K}} \varpi(\tau)\mathbf{g}(\tau)\begin{bmatrix}\mathbf{f}^{\top}(\tau) & \mathsf{f}_{d+1}(\tau)\end{bmatrix}\mathrm{d}\tau\right)\begin{bmatrix}\mathcal{F}_d \oplus \left(\int_{\mathcal{K}} \varpi(\tau)\mathsf{f}_{d+1}^2(\tau)\mathrm{d}\tau\right)^{-1}\end{bmatrix}\begin{bmatrix}\mathbf{f}(\tau) \\ \mathsf{f}_{d+1}(\tau)\end{bmatrix} \\
&= \mathbf{g}(\tau) - \begin{bmatrix}\mathsf{A}_d & \boldsymbol{\alpha}_{d+1}\end{bmatrix}\begin{bmatrix}\mathbf{f}(\tau) \\ \mathsf{f}_{d+1}(\tau)\end{bmatrix} = \mathbf{e}_d(\tau) - \mathsf{f}_{d+1}(\tau)\boldsymbol{\alpha}_{d+1}, \quad \forall d \in \mathbb{N}
\end{aligned} \tag{5.30}$$

where $\boldsymbol{\alpha}_{d+1}$ is defined in (5.29) and $\mathbf{e}_d(\tau) = \mathbf{g}(\tau) - \mathsf{A}_d\mathbf{f}(\tau)$. Note that the index $d$ is added to the symbols $\mathsf{A}$ and $\mathbf{e}(\tau)$ in Lemma 5.2 without causing ambiguity. By (5.30) and (D.1), we have

$$\begin{aligned}
0 \prec \mathcal{E}_{d+1} = \int_{\mathcal{K}} \varpi(\tau)\mathbf{e}_{d+1}(\tau)\mathbf{e}_{d+1}^{\top}(\tau)\mathrm{d}\tau &= \mathcal{E}_d - \mathsf{Sy}\left(\boldsymbol{\alpha}_{d+1}\int_{\mathcal{K}} \varpi(\tau)\mathsf{f}_{d+1}(\tau)\mathbf{e}_d^{\top}(\tau)\mathrm{d}\tau\right) \\
&\quad + \left(\int_{\mathcal{K}} \varpi(\tau)\mathsf{f}_{d+1}^2(\tau)\mathrm{d}\tau\right)\boldsymbol{\alpha}_{d+1}\boldsymbol{\alpha}_{d+1}^{\top}.
\end{aligned} \tag{5.31}$$

By (D.2) and the fact that $\int_{\mathcal{K}} \varpi(\tau)\mathsf{f}_{d+1}(\tau)\mathbf{f}(\tau)\mathrm{d}\tau = \mathbf{0}_d$ due to (5.28), we have

$$\begin{aligned}
\mathsf{O}_{\delta,d+1} &= \int_{\mathcal{K}} \varpi(\tau)\mathbf{e}_{d+1}(\tau)\begin{bmatrix}\mathbf{f}^{\top}(\tau) & \mathsf{f}_{d+1}(\tau)\end{bmatrix}\mathrm{d}\tau \\
&= \int_{\mathcal{K}} \varpi(\tau)\left(\mathbf{e}_d(\tau) - \boldsymbol{\alpha}_{d+1}\mathsf{f}_{d+1}(\tau)\right)\begin{bmatrix}\mathbf{f}^{\top}(\tau) & \mathsf{f}_{d+1}(\tau)\end{bmatrix}\mathrm{d}\tau \\
&= \int_{\mathcal{K}} \varpi(\tau)\begin{bmatrix}\mathbf{e}_d(\tau)\mathbf{f}^{\top}(\tau) & \mathsf{f}_{d+1}(\tau)\mathbf{e}_d(\tau)\end{bmatrix}\mathrm{d}\tau - \boldsymbol{\alpha}_{d+1}\int_{\mathcal{K}} \varpi(\tau)\begin{bmatrix}\mathsf{f}_{d+1}(\tau)\mathbf{f}^{\top}(\tau) & \mathsf{f}_{d+1}^2(\tau)\end{bmatrix}\mathrm{d}\tau \\
&= \begin{bmatrix}\mathsf{O}_{\delta,d} & \int_{\mathcal{K}} \varpi(\tau)\mathsf{f}_{d+1}(\tau)\mathbf{e}_d(\tau)\mathrm{d}\tau\end{bmatrix} - \begin{bmatrix}\mathsf{O}_{\delta,d} & \int_{\mathcal{K}} \varpi(\tau)\mathsf{f}_{d+1}^2(\tau)\mathrm{d}\tau\,\boldsymbol{\alpha}_{d+1}\end{bmatrix} = \mathsf{O}_{\delta,d+1}.
\end{aligned} \tag{5.32}$$

Now (5.32) leads to the equality $\int_{\mathcal{K}} \varpi(\tau)\mathsf{f}_{d+1}(\tau)\mathbf{e}_d(\tau)\mathrm{d}\tau = \int_{\mathcal{K}} \varpi(\tau)\mathsf{f}_{d+1}^2(\tau)\mathrm{d}\tau\,\boldsymbol{\alpha}_{d+1}$. Utilizing this equality on (5.31) yields (5.29) given that $\int_{\mathcal{K}} \varpi(\tau)\mathsf{f}_{d+1}^2(\tau)\mathrm{d}\tau > 0$ and $\boldsymbol{\alpha}_{d+1}\boldsymbol{\alpha}_{d+1}^{\top} \succeq 0$. ■

**Remark 5.9.** The result in [470, Lemma 1] is generalized by Corollary 5.2 as $\mathbf{f}(\cdot)$ can have Legendre polynomials only with $\varpi(\tau) = 1$. Moreover, let $\acute{\boldsymbol{f}}(\cdot)$, $\grave{\boldsymbol{f}}(\cdot)$ in (5.7) consist only of orthogonal functions over the intervals $\mathcal{I}_1$ and $\mathcal{I}_2$, respectively. Then, $\acute{\mathsf{E}}_d$ and $\grave{\mathsf{E}}_\delta$ in (5.12) follow the property in (5.29) with $\varpi(\tau) = 1$.

## 5.4 Main Results on Dissipativity and Stability Analysis

The main results on the dissipativity and stability analysis of (5.1) are set forth in Theorem 5.1, where the condition for the dissipativity and stability analysis of (5.1) is expressed in terms of LMIs. Moreover, we also demonstrate in Corollary 5.3 and Corollary 5.4 that the resulting condition in Theorem 5.1 can exhibit a hierarchical pattern under certain prerequisites.

**Theorem 5.1.** *Suppose that all functions and the parameters in (5.3)–(5.12) are given with* $\mu_1; \mu_2 \in \mathbb{N}_0$ *and* $d; \delta \in \mathbb{N}$. *Assume also that there exist* $\acute{\boldsymbol{g}}(\cdot) \in C^1(\mathbb{R}; \mathbb{R}^{p_1})$, $\grave{\boldsymbol{g}}(\cdot) \in C^1(\mathbb{R}; \mathbb{R}^{p_2})$ *and* $N_1 \in \mathbb{R}^{p_1 \times d}$, $N_2 \in \mathbb{R}^{p_2 \times \delta}$ *such that*

$$
\begin{aligned}
&(\tau + r_1)\frac{\mathsf{d}\acute{\boldsymbol{g}}(\tau)}{\mathsf{d}\tau} = N_1 \acute{\boldsymbol{f}}(\tau) && (\tau + r_2)\frac{\mathsf{d}\grave{\boldsymbol{g}}(\tau)}{\mathsf{d}\tau} = N_2 \grave{\boldsymbol{f}}(\tau) \\
&\mathbb{S}^{p_1} \ni \acute{\mathsf{G}}_{p_1}^{-1} = \int_{-r_1}^{0} (\tau + r_1)\acute{\boldsymbol{g}}(\tau)\acute{\boldsymbol{g}}^\top(\tau)\mathsf{d}\tau \succ 0 && \mathbb{S}^{p_2} \ni \grave{\mathsf{G}}_{p_2}^{-1} = \int_{-r_2}^{-r_1} (\tau + r_2)\grave{\boldsymbol{g}}(\tau)\grave{\boldsymbol{g}}^\top(\tau)\mathsf{d}\tau \succ 0.
\end{aligned}
\tag{5.33}
$$

*Given* $r_2 > r_1 > 0$, *then delay system (5.14) with supply rate function (5.22) is dissipative and the origin of (5.14) is globally asymptotically stable with* $\boldsymbol{w}(t) \equiv \mathbf{0}_q$, *if there exist symmetric matrices* $P \in \mathbb{S}^l$ *and* $Q_1; Q_2; R_1; R_2; S_1; S_2; U_1; U_2 \in \mathbb{S}^\nu$ *such that the inequalities*

$$
\begin{aligned}
\mathbf{P} := P + \left(\mathsf{O}_{n+2\nu} \oplus \left[\acute{\mathsf{F}}_d \otimes Q_1\right] \oplus \left[\grave{\mathsf{F}}_\delta \otimes Q_2\right]\right) + \Pi^\top \left(G_1^\top \acute{\Phi}_{\kappa_1} G_1 \otimes S_1 + G_2^\top \grave{\Phi}_{\kappa_2} G_2 \otimes S_2\right) \Pi \\
+ \Pi^\top \left(H_1^\top \acute{\mathsf{G}}_{p_1} H_1 \otimes U_1 + H_2^\top \grave{\mathsf{G}}_{p_2} H_2 \otimes U_2\right) \Pi \succ 0,
\end{aligned}
\tag{5.34}
$$

$$
Q_1 \succeq 0, \quad Q_2 \succeq 0, \quad R_1 \succeq 0, \quad R_2 \succeq 0, \quad S_1 \succeq 0, \quad S_2 \succeq 0, \quad U_1 \succeq 0, \quad U_2 \succeq 0,
\tag{5.35}
$$

$$
\widetilde{\mathbf{\Omega}} = \begin{bmatrix} J_1 & \mathsf{O}_{m,\nu} & \widetilde{J}\Sigma \\ * & -S_1 - r_1 U_1 & (S_1 + r_1 U_1)(A_6 \mathbf{A} + Y) \\ * & * & \mathbf{\Omega} \end{bmatrix} \prec 0
\tag{5.36}
$$

*hold, where the positive definite matrices* $\acute{\mathsf{F}}_d$, $\grave{\mathsf{F}}_\delta$, $\acute{\Phi}_{\kappa_1}$ *and* $\acute{\Phi}_{\kappa_2}$ *are given in (5.5) and (5.6), and the parameters* $\mathbf{A}$ *and* $\Sigma$ *have been defined in (5.15)–(5.17). Moreover,*

$$
\Pi := \begin{bmatrix} \Xi \\ \mathsf{O}_{(2\nu+\varrho\nu),n} \quad I_{2\nu+\varrho\nu} \end{bmatrix}, \quad Y := \begin{bmatrix} A_7 & A_8 & \mathsf{O}_{\nu,(q+l+\mu\nu)} \end{bmatrix}
\tag{5.37}
$$

*with* $\varrho = d + \delta$ *and* $l = n + 2\nu + \varrho\nu$ *and* $\mu = \mu_1 + \mu_2$, *and*

$$
\begin{aligned}
\mathbf{\Omega} := &\ \mathsf{Sy}\left(\Theta_2^\top P \Theta_1 - \begin{bmatrix} \mathsf{O}_{(2\nu+q+l+\mu\nu),2\nu} & \Sigma^\top J_2 & \mathsf{O}_{(2\nu+q+l+\mu\nu),(l+\mu\nu)} \end{bmatrix}\right) \\
&- \left(\mathsf{O}_{q+2\nu} \oplus \left[\Pi^\top \left(G_1^\top \acute{\Phi}_{\kappa_1} G_1 \otimes U_1 + G_2^\top \grave{\Phi}_{\kappa_2} G_2 \otimes U_2\right) \Pi\right] \oplus \mathsf{O}_{\mu\nu}\right) \\
&- \Bigg( [S_1 - S_2 - r_3 U_2] \oplus S_2 \oplus J_3 \oplus \mathsf{O}_n \oplus [Q_1 - Q_2 - r_3 R_2] \oplus Q_2 \oplus \left[\acute{\mathsf{F}}_d \otimes R_1\right] \oplus \left[\grave{\mathsf{F}}_\delta \otimes R_2\right] \cdots \\
&\qquad \oplus \left[\acute{\mathsf{E}}_d \otimes R_1\right] \oplus \left[\grave{\mathsf{E}}_\delta \otimes R_2\right] \Bigg) + [*](Q_1 + r_1 R_1) \begin{bmatrix} \mathsf{O}_{\nu,(2\nu+q)} & \Xi & \mathsf{O}_{\nu,\nu\mu} \end{bmatrix}
\end{aligned}
\tag{5.38}
$$

*where*

$$
G_1 = \begin{bmatrix} \acute{\boldsymbol{\phi}}(0) & -\acute{\boldsymbol{\phi}}(-r_1) & \mathbf{0}_{\kappa_1} & -M_3 & \mathsf{O}_{\kappa_1,d} \end{bmatrix}, \; G_2 = \begin{bmatrix} \mathbf{0}_{\kappa_2} & \grave{\boldsymbol{\phi}}(-r_1) & -\grave{\boldsymbol{\phi}}(-r_2) & \mathsf{O}_{\kappa_2,\delta} & -M_4 \end{bmatrix}
\tag{5.39}
$$

$$
G_3 = \begin{bmatrix} \acute{\boldsymbol{f}}(0) & -\acute{\boldsymbol{f}}(-r_1) & \mathbf{0}_d & -M_1 & \mathsf{O}_d \end{bmatrix}, \; G_4 = \begin{bmatrix} \mathbf{0}_\delta & \grave{\boldsymbol{f}}(-r_1) & -\grave{\boldsymbol{f}}(-r_2) & \mathsf{O}_\delta & -M_2 \end{bmatrix},
\tag{5.40}
$$

$$
H_1 = \begin{bmatrix} r_1 \acute{\boldsymbol{g}}(0) & \mathbf{0}_{p_1} & \mathbf{0}_{p_1} & -N_1 & \mathsf{O}_{p_1} \end{bmatrix}, \; H_2 = \begin{bmatrix} \mathbf{0}_{p_2} & (r_2 - r_1)\grave{\boldsymbol{g}}(-r_1) & \mathbf{0}_{p_2} & \mathsf{O}_{p_2} & -N_2 \end{bmatrix},
\tag{5.41}
$$

$$
\Theta_1 := \begin{bmatrix} & \mathbf{A} & \\ & I_{2\nu} \quad \mathsf{O}_{2\nu,(q+l+\mu\nu)} & \\ \mathsf{O}_{d\nu,(2\nu+q)} & (G_3 \otimes I_\nu)\Pi & \mathsf{O}_{d\nu,\nu\mu} \\ \mathsf{O}_{\delta\nu,(2\nu+q)} & (G_4 \otimes I_\nu)\Pi & \mathsf{O}_{\delta\nu,\nu\mu} \end{bmatrix}, \quad \Theta_2 := \begin{bmatrix} \mathsf{O}_{l,(2\nu+q)} & I_l & \mathsf{O}_{l,\mu\nu} \end{bmatrix}.
\tag{5.42}
$$

*Proof.* Given $r_2 > r_1 > 0$, our goal is to construct LKF

$$
\begin{aligned}
&\mathsf{v}(\boldsymbol{\xi}, \boldsymbol{\psi}(\cdot)) = \boldsymbol{\eta}^\top(t_0)P\boldsymbol{\eta}(t_0) + \int_{-r_1}^{0} \boldsymbol{y}^\top(t+\tau)\left[Q_1 + (\tau + r_1)R_1\right]\boldsymbol{y}(t+\tau)\mathsf{d}\tau \\
&+ \int_{-r_2}^{-r_1} \boldsymbol{y}^\top(t+\tau)\left[Q_2 + (\tau + r_2)R_2\right]\boldsymbol{y}(t+\tau)\mathsf{d}\tau + \int_{-r_1}^{0} \dot{\boldsymbol{y}}^\top(t+\tau)\left[S_1 + (\tau + r_1)U_1\right]\dot{\boldsymbol{y}}(t+\tau)\mathsf{d}\tau \\
&+ \int_{-r_2}^{-r_1} \dot{\boldsymbol{y}}^\top(t+\tau)\left[S_2 + (\tau + r_2)U_2\right]\dot{\boldsymbol{y}}(t+\tau)\mathsf{d}\tau
\end{aligned} \tag{5.43}
$$

to prove the statements in Theorem 5.1, where $\boldsymbol{x}(t)$ and $\mathbf{y}_t(\cdot)$ follow the same definition in (5.21). Moreover,

$$
\boldsymbol{\eta}(t) := \mathbf{Col}\left[\boldsymbol{x}(t),\ \boldsymbol{y}(t-r_1),\ \boldsymbol{y}(t-r_2),\ \int_{-r_1}^{0} \acute{F}_d(\tau)\boldsymbol{y}(t+\tau)\mathsf{d}\tau,\ \int_{-r_2}^{-r_1} \grave{F}_\delta(\tau)\boldsymbol{y}(t+\tau)\mathsf{d}\tau\right] \tag{5.44}
$$

with $\acute{F}_d(\tau)$ and $\grave{F}_\delta(\tau)$ defined in (5.18), and the matrix parameters in (5.43) are defined as $P \in \mathbb{S}^l$ and $Q_1; Q_2; R_1; R_2; S_1; S_2; U_1; U_2 \in \mathbb{S}^\nu$ with $l := n + 2\nu + \varrho\nu$ and $\varrho := d + \delta$. Given that the eigenvalues of the matrix terms $Q_1 + (\tau + r_1)R_1$, $Q_2 + (\tau + r_2)R_2$ $S_1 + (\tau + r_1)U_1$ and $S_2 + (\tau + r_2)U_2$ in (5.43) are bounded, it follows that all associated quadratic integrals are well-defined since $\mathbf{y}_t(\cdot) \in \mathcal{A}\left([-r_2, 0); \mathbb{R}^\nu\right)$. In addition, since $\mathbf{y}_t(\cdot)$, $\acute{\boldsymbol{f}}(\tau)$ and $\grave{\boldsymbol{f}}(\tau)$ are bounded, the integrals in (5.44) are well defined as well.

Our first step is to demonstrate that the existence of feasible solutions to inequalities (5.35) and (5.36) implies that (5.43) satisfies both (5.20) and (5.21). We then show that one can infer (5.19) from the existence of the feasible solutions of (5.34) and (5.35). The existence of the upper bound of $\mathsf{v}(\boldsymbol{x}(t), \mathbf{y}_t(\cdot))$ can be independently proved without considering the inequalities (5.34)–(5.36).

Let $t_0 \in \mathbb{R}$, differentiate $\mathsf{v}(\boldsymbol{x}(t), \mathbf{y}_t(\cdot))$ along the trajectory of (5.14) and consider (5.22), it produces

$$
\begin{aligned}
&\dot{\mathsf{v}}(\boldsymbol{x}(t), \mathbf{y}_t(\cdot)) - \mathsf{s}(\boldsymbol{z}(t), \boldsymbol{w}(t)) = \boldsymbol{\vartheta}^\top(t)\,\mathsf{Sy}\left(\Theta_2^\top P\Theta_1\right)\boldsymbol{\vartheta}(t) + \boldsymbol{y}^\top(t)\left(Q_1 + r_1R_1\right)\boldsymbol{y}(t) \\
&+ \boldsymbol{y}^\top(t-r_1)\left(Q_2 + r_3R_2 - Q_1\right)\boldsymbol{y}(t-r_1) - \boldsymbol{y}^\top(t-r_2)Q_2\boldsymbol{y}(t-r_2) + \dot{\boldsymbol{y}}^\top(t)\left(S_1 + r_1U_1\right)\dot{\boldsymbol{y}}(t) \\
&+ \dot{\boldsymbol{y}}^\top(t-r_1)\left(S_2 + r_3U_2 - S_1\right)\dot{\boldsymbol{y}}(t-r_1) - \dot{\boldsymbol{y}}^\top(t-r_2)S_2\dot{\boldsymbol{y}}(t-r_2) - \int_{-r_1}^{0} \boldsymbol{y}^\top(t+\tau)R_1\boldsymbol{y}(t+\tau)\mathsf{d}\tau \\
&- \int_{-r_2}^{-r_1} \boldsymbol{y}^\top(t+\tau)R_2\boldsymbol{y}(t+\tau)\mathsf{d}\tau - \int_{-r_1}^{0} \dot{\boldsymbol{y}}^\top(t+\tau)U_1\dot{\boldsymbol{y}}(t+\tau)\mathsf{d}\tau - \int_{-r_2}^{-r_1} \dot{\boldsymbol{y}}^\top(t+\tau)U_2\dot{\boldsymbol{y}}(t+\tau)\mathsf{d}\tau \\
&- \boldsymbol{w}^\top(t)J_3\boldsymbol{w}(t) - \boldsymbol{\vartheta}^\top(t)\left[\Sigma^\top\widetilde{J}^\top J_1^{-1}\widetilde{J}\Sigma + \mathsf{Sy}\left(\begin{bmatrix}\mathsf{O}_{(2\nu+q+l+\mu\nu),2\nu} & \Sigma^\top J_2 & \mathsf{O}_{(2\nu+q+l+\mu\nu),(l+\mu\nu)}\end{bmatrix}\right)\right]\boldsymbol{\vartheta}(t),
\end{aligned} \tag{5.45}
$$

where $\boldsymbol{\vartheta}(t)$ and $\Theta_1; \Theta_2$ have been defined in (5.18) and (5.42), respectively, and the matrices $G_3$ and $G_4$ in (5.40) are obtained via the relations

$$
\begin{aligned}
\int_{-r_1}^{0} \acute{F}_d(\tau)\dot{\boldsymbol{y}}(t+\tau)\mathsf{d}\tau &= \acute{F}_d(0)\boldsymbol{y}(t) - \acute{F}_d(-r_1)\boldsymbol{y}(t-r_1) - (M_1 \otimes I_\nu)\int_{-r_1}^{0} \acute{F}_d(\tau)\boldsymbol{y}(t+\tau)\mathsf{d}\tau = \\
&= \begin{bmatrix}\mathsf{O}_{d\nu,(q+2\nu)} & (G_3 \otimes I_\nu)\,\Pi & \mathsf{O}_{d\nu,\nu\mu}\end{bmatrix}\boldsymbol{\vartheta}(t),
\end{aligned} \tag{5.46}
$$

$$
\begin{aligned}
\int_{-r_2}^{-r_1} \grave{F}_\delta(\tau)\dot{\boldsymbol{y}}(t+\tau)\mathsf{d}\tau &= \grave{F}_\delta(-r_1)\boldsymbol{y}(t-r_1) - \grave{F}_\delta(-r_2)\boldsymbol{y}(t-r_2) - (M_2 \otimes I_\nu)\int_{-r_2}^{-r_1} \grave{F}_\delta(\tau)\boldsymbol{y}(t+\tau)\mathsf{d}\tau \\
&= \begin{bmatrix}\mathsf{O}_{\delta\nu,(q+2\nu)} & (G_4 \otimes I_\nu)\,\Pi & \mathsf{O}_{\delta\nu,\nu\mu}\end{bmatrix}\boldsymbol{\vartheta}(t),
\end{aligned} \tag{5.47}
$$

which themselves can be derived by utilizing (2.1b), (2.1a) with (5.3).

To construct an upper bound for (5.45), let $R_1 \succeq 0$, $R_2 \succeq 0$ so that the inequalities

$$
\begin{aligned}
\int_{-r_1}^{0} \boldsymbol{y}^\top(t+\tau)R_1\boldsymbol{y}(t+\tau)\mathsf{d}\tau &\geq \int_{-r_1}^{0} \boldsymbol{y}^\top(t+\tau)\acute{F}_d^\top(\tau)\mathsf{d}\tau \left(\acute{\mathsf{F}}_d \otimes R_1\right) \int_{-r_1}^{0} \acute{F}_d(\tau)\boldsymbol{y}(t+\tau)\mathsf{d}\tau \\
&+ \int_{-r_1}^{0} \boldsymbol{y}^\top(t+\tau)\acute{E}_d^\top(\tau)\mathsf{d}\tau \left(\acute{\mathsf{E}}_d \otimes R_1\right) \int_{-r_1}^{0} \acute{E}_d(\tau)\boldsymbol{y}(t+\tau)\mathsf{d}\tau, \\
\int_{-r_2}^{-r_1} \boldsymbol{y}^\top(t+\tau)R_2\boldsymbol{y}(t+\tau)\mathsf{d}\tau &\geq \int_{-r_2}^{-r_1} \boldsymbol{y}^\top(t+\tau)\grave{F}_\delta^\top(\tau)\mathsf{d}\tau \left(\grave{\mathsf{F}}_\delta \otimes R_2\right) \int_{-r_2}^{-r_1} \grave{F}_\delta(\tau)\boldsymbol{y}(t+\tau)\mathsf{d}\tau \\
&+ \int_{-r_2}^{-r_1} \boldsymbol{y}^\top(t+\tau)\grave{E}_\delta^\top(\tau)\mathsf{d}\tau \left(\grave{\mathsf{E}}_\delta \otimes R_2\right) \int_{-r_2}^{-r_1} \grave{E}_\delta(\tau)\boldsymbol{y}(t+\tau)\mathsf{d}\tau
\end{aligned}
\tag{5.48}
$$

can be derived from (5.24) with $\mathsf{f}(\tau) = \acute{\boldsymbol{f}}(\tau)$; $\mathsf{g}(\tau) = \boldsymbol{\varphi}_1(\tau)$ and $\mathsf{f}(\tau) = \grave{\boldsymbol{f}}(\tau)$; $\mathsf{g}(\tau) = \boldsymbol{\varphi}_2(\tau)$, respectively, which matches $\acute{F}_d(\tau)$; $\grave{F}_\delta(\tau)$ in (5.18) and the expressions in (5.11). Furthermore, let $U_1 \succeq 0$ and $U_2 \succeq 0$ and apply (5.27) to the integral terms $\int_{-r_1}^{0} \dot{\boldsymbol{y}}^\top(t+\tau)U_1\dot{\boldsymbol{y}}(t+\tau)\mathsf{d}\tau$ and $\int_{-r_2}^{-r_1} \dot{\boldsymbol{y}}^\top(t+\tau)U_2\dot{\boldsymbol{y}}(t+\tau)\mathsf{d}\tau$ with $\mathsf{f}(\tau) = \acute{\boldsymbol{\phi}}(\tau)$ and $\mathsf{f}(\tau) = \grave{\boldsymbol{\phi}}(\tau)$ in (5.6), respectively, and consider the expression $\boldsymbol{y}(t) = \begin{bmatrix} \mathsf{O}_{\nu,(2\nu+q)} & \Xi & \mathsf{O}_{\nu,\nu\mu} \end{bmatrix} \boldsymbol{\vartheta}(t)$ in (5.16) with (2.1a) and (2.1b). It produces

$$
\begin{aligned}
\int_{-r_1}^{0} \dot{\boldsymbol{y}}^\top(t+\tau)U_1\dot{\boldsymbol{y}}(t+\tau)\mathsf{d}\tau &\geq [*]\left(\acute{\mathsf{\Phi}}_{\kappa_1} \otimes U_1\right)\left[\int_{-r_1}^{0} \left(\acute{\boldsymbol{\phi}}(\tau) \otimes I_\nu\right)\dot{\boldsymbol{y}}(t+\tau)\mathsf{d}\tau\right] \\
&= \boldsymbol{\vartheta}^\top(\tau)\left[\mathsf{O}_{2\nu+q} \oplus \Pi^\top\left(G_1^\top \acute{\mathsf{\Phi}}_{\kappa_1} G_1 \otimes U_1\right)\Pi \oplus \mathsf{O}_{\nu\mu}\right]\boldsymbol{\vartheta}(\tau),
\end{aligned}
\tag{5.49}
$$

$$
\begin{aligned}
\int_{-r_2}^{-r_1} \dot{\boldsymbol{y}}^\top(t+\tau)U_2\dot{\boldsymbol{y}}(t+\tau)\mathsf{d}\tau &\geq [*]\left(\grave{\mathsf{\Phi}}_{\kappa_2} \otimes U_2\right)\left[\int_{-r_2}^{-r_1} \left(\grave{\boldsymbol{\phi}}(\tau) \otimes I_\nu\right)\dot{\boldsymbol{y}}(t+\tau)\mathsf{d}\tau\right] \\
&= \boldsymbol{\vartheta}^\top(\tau)\left[\mathsf{O}_{2\nu+q} \oplus \Pi^\top\left(G_2^\top \grave{\mathsf{\Phi}}_{\kappa_2} G_2 \otimes U_2\right)\Pi \oplus \mathsf{O}_{\nu\mu}\right]\boldsymbol{\vartheta}(\tau),
\end{aligned}
\tag{5.50}
$$

where $G_1$ and $G_2$ are given in (5.39) which are established by the relations

$$
\begin{aligned}
&\int_{-r_1}^{0} \left(\acute{\boldsymbol{\phi}}(\tau) \otimes I_\nu\right)\dot{\boldsymbol{y}}(t+\tau)\mathsf{d}\tau = \left(\acute{\boldsymbol{\phi}}(0) \otimes I_\nu\right)\boldsymbol{y}(t) - \left(\acute{\boldsymbol{\phi}}(-r_1) \otimes I_\nu\right)\boldsymbol{y}(t-r_1) \\
&\quad - (M_3 \otimes I_\nu)\int_{-r_1}^{0} \left(\acute{\boldsymbol{\phi}}(\tau) \otimes I_\nu\right)\boldsymbol{y}(t+\tau)\mathsf{d}\tau = (G_1 \otimes I_\nu)\,\Pi\boldsymbol{\eta}(t) \\
&\quad = \begin{bmatrix} \mathsf{O}_{\kappa_1\nu,(q+2\nu)} & (G_1 \otimes I_\nu)\,\Pi & \mathsf{O}_{\kappa_1\nu,\nu\mu} \end{bmatrix}\boldsymbol{\vartheta}(t),
\end{aligned}
\tag{5.51}
$$

$$
\begin{aligned}
&\int_{-r_2}^{-r_1} \left(\grave{\boldsymbol{\phi}}(\tau) \otimes I_\nu\right)\dot{\boldsymbol{y}}(t+\tau)\mathsf{d}\tau = \left(\grave{\boldsymbol{\phi}}(-r_1) \otimes I_\nu\right)\boldsymbol{y}(t-r_1) - \left(\grave{\boldsymbol{\phi}}(-r_2) \otimes I_\nu\right)\boldsymbol{y}(t-r_2) \\
&\quad - (M_4 \otimes I_\nu)\int_{-r_2}^{-r_1} \left(\grave{\boldsymbol{\phi}}(\tau) \otimes I_\nu\right)\boldsymbol{y}(t+\tau)\mathsf{d}\tau = (G_2 \otimes I_\nu)\,\Pi\boldsymbol{\eta}(t) \\
&\quad = \begin{bmatrix} \mathsf{O}_{\kappa_2\nu,(q+2\nu)} & (G_2 \otimes I_\nu)\,\Pi & \mathsf{O}_{\kappa_2\nu,\nu\mu} \end{bmatrix}\boldsymbol{\vartheta}(t).
\end{aligned}
\tag{5.52}
$$

Now applying (5.48)–(5.50) with (5.35) to (5.45) yields

$$
\begin{aligned}
&\forall t \geq t_0,\ \ \dot{\mathsf{v}}(\boldsymbol{x}(t), \mathbf{y}_t(\cdot)) - \mathsf{s}(\boldsymbol{z}(t), \boldsymbol{w}(t)) \leq \\
&\qquad \boldsymbol{\vartheta}^\top(t)\left[\boldsymbol{\Omega} + (A_6\mathbf{A} + Y)^\top (S_1 + r_1U_1)(A_6\mathbf{A} + Y) - \Sigma^\top \widetilde{J}^\top J_1^{-1}\widetilde{J}\Sigma\right]\boldsymbol{\vartheta}(t),
\end{aligned}
\tag{5.53}
$$

where $\boldsymbol{\Omega}$ is defined in (5.38). It is evident that if $\boldsymbol{\Omega} + (A_6\mathbf{A}+Y)^\top (S_1 + r_1U_1)(A_6\mathbf{A}+Y) - \Sigma^\top \widetilde{J}^\top J_1^{-1}\widetilde{J}\Sigma \prec 0$ holds with (5.35), we have

$$
\exists \epsilon_3 > 0,\ \forall t \geq t_0,\ \ \dot{\mathsf{v}}(\boldsymbol{x}(t), \mathbf{y}_t(\cdot)) - \mathsf{s}(\boldsymbol{z}(t), \boldsymbol{w}(t)) \leq -\epsilon_3 \left\|\boldsymbol{x}(t)\right\|_2 .
\tag{5.54}
$$

Moreover, let $\boldsymbol{w}(t) \equiv \mathbf{0}_q$ and consider the structure of the quadratic term in (5.53) in conjunction with the properties of negative definite matrices. It can be concluded that $\boldsymbol{\Omega}+(A_6\mathbf{A}+Y)^\top (S_1+r_1U_1)(A_6\mathbf{A}+Y)-\Sigma^\top\widetilde{J}^\top J_1^{-1}\widetilde{J}\Sigma \prec 0$ and (5.35) are satisfied, it implies that

$$\exists\epsilon_3>0,\quad \frac{\mathsf{d}^+}{\mathsf{d}t}\mathsf{v}(\boldsymbol{x}(t),\mathbf{y}_t(\cdot))\bigg|_{t=t_0,\boldsymbol{x}(t_0)=\boldsymbol{\xi},\mathbf{y}_{t_0}(\cdot)=\boldsymbol{\phi}(\cdot)} \le -\epsilon_3\left\|\boldsymbol{\xi}\right\|_2, \tag{5.55}$$

where $\boldsymbol{x}(t)$ and $\mathbf{y}_t(\cdot)$ here follow the same definition in Lemma 5.1. As a result, we see that (5.35) with $\boldsymbol{\Omega}+(A_6\mathbf{A}+Y)^\top (S_1+r_1U_1)(A_6\mathbf{A}+Y)-\Sigma^\top\widetilde{J}^\top J_1^{-1}\widetilde{J}\Sigma \prec 0$ implies (5.20) and (5.21). Finally, by applying the Schur complement to $\boldsymbol{\Omega}+(A_6\mathbf{A}+Y)^\top (S_1+r_1U_1)(A_6\mathbf{A}+Y)-\Sigma^\top\widetilde{J}^\top J_1^{-1}\widetilde{J}\Sigma \prec 0$ with (5.35) and $J_1^{-1} \prec 0$, we can obtain (5.36). Hence, we have proved that the feasible solutions to (5.35) and (5.36) implies that (5.43) satisfies (5.20) and (5.21).

We now proceed to demonstrate that if (5.34) and (5.35) hold then (5.43) satisfies (5.19). Let $\left\|\boldsymbol{\psi}(\cdot)\right\|_\infty := \sup_{-r_2\le r\le 0}\left\|\boldsymbol{\psi}(\tau)\right\|_2$ and $\|\boldsymbol{\psi}(\cdot)\|_2^2 := \int_{-r_2}^{0}\boldsymbol{\psi}^\top(\tau)\boldsymbol{\psi}(\tau)\mathsf{d}\tau$. Given the structure of (5.43) with $t=t_0$, it follows that there exist $\lambda;\eta>0:$

$$\begin{aligned}
&\mathsf{v}(\boldsymbol{x}(t_0),\mathbf{y}_{t_0}(\cdot)) = \mathsf{v}(\boldsymbol{\xi},\boldsymbol{\psi}(\cdot)) \le \boldsymbol{\eta}^\top(t_0)\lambda\boldsymbol{\eta}(t_0)+\lambda\int_{-r_2}^{0}\begin{bmatrix}\boldsymbol{\psi}^\top(\tau) & \dot{\boldsymbol{\psi}}^\top(\tau)\end{bmatrix}\begin{bmatrix}\boldsymbol{\psi}(\tau)\\ \dot{\boldsymbol{\psi}}(\tau)\end{bmatrix}\mathsf{d}\tau\\
&\le \lambda\left\|\boldsymbol{\xi}\right\|_2^2+(2\lambda+\lambda r_2)\left\|\boldsymbol{\psi}(\cdot)\right\|_\infty^2+\lambda r_2\|\dot{\boldsymbol{\psi}}(\cdot)\|_2^2+\lambda\int_{-r_1}^{0}\boldsymbol{\psi}^\top(\tau)\acute{F}_d^\top(\tau)\mathsf{d}\tau\int_{-r_1}^{0}\acute{F}_d(\tau)\boldsymbol{\psi}(\tau)\mathsf{d}\tau\\
&+\lambda\int_{-r_2}^{-r_1}\boldsymbol{\psi}^\top(\tau)\grave{F}_\delta^\top(\tau)\mathsf{d}\tau\int_{-r_2}^{-r_1}\grave{F}_\delta(\tau)\boldsymbol{\psi}(\tau)\mathsf{d}\tau \le \lambda\left\|\boldsymbol{\xi}\right\|_2^2+(2\lambda+\lambda r_2)\left\|\boldsymbol{\psi}(\cdot)\right\|_\infty^2\\
&+\lambda r_2\|\dot{\boldsymbol{\psi}}(\cdot)\|_2^2+\int_{-r_1}^{0}\boldsymbol{\psi}^\top(\tau)\acute{F}_d^\top(\tau)\mathsf{d}\tau\left(\eta\acute{\mathsf{F}}_d\otimes I_\nu\right)\int_{-r_1}^{0}\acute{F}_d(\tau)\boldsymbol{\psi}(\tau)\mathsf{d}\tau\\
&+\int_{-r_2}^{-r_1}\boldsymbol{\psi}^\top(\tau)\grave{F}_\delta^\top(\tau)\mathsf{d}\tau\left(\eta\grave{\mathsf{F}}_\delta\otimes I_\nu\right)\int_{-r_2}^{-r_1}\grave{F}_\delta(\tau)\boldsymbol{\psi}(\tau)\mathsf{d}\tau \le \lambda\left\|\boldsymbol{\xi}\right\|_2^2+(2\lambda+\lambda r_2)\left\|\boldsymbol{\psi}(\cdot)\right\|_\infty^2\\
&+\lambda r_2\|\dot{\boldsymbol{\psi}}(\cdot)\|_2^2+\eta\int_{-r_1}^{0}\boldsymbol{\psi}^\top(\tau)\boldsymbol{\psi}(\tau)\mathsf{d}\tau+\eta\int_{-r_2}^{-r_1}\boldsymbol{\psi}^\top(\tau)\boldsymbol{\psi}(\tau)\mathsf{d}\tau = \lambda\left\|\boldsymbol{\xi}\right\|_2^2\\
&+(2\lambda+\lambda r_2+\eta r_2)\left\|\boldsymbol{\psi}(\cdot)\right\|_\infty^2+\lambda r_2\|\dot{\boldsymbol{\psi}}(\cdot)\|_2^2 \le (2\lambda+\lambda r_2+\eta r_2)\left(\left\|\boldsymbol{\xi}\right\|_2^2+\left\|\boldsymbol{\psi}(\cdot)\right\|_\infty^2+\|\dot{\boldsymbol{\psi}}(\cdot)\|_2^2\right)\\
&\le (2\lambda+\lambda r_2+\eta r_2)\left[\left\|\boldsymbol{\xi}\right\|_2^2+\left(\left\|\boldsymbol{\psi}(\cdot)\right\|_\infty+\|\dot{\boldsymbol{\psi}}(\cdot)\|_2\right)^2\right]\\
&\le (4\lambda+2\lambda r_2+2\eta r_2)\left[\left\|\boldsymbol{\xi}\right\|_2\vee\left(\left\|\boldsymbol{\psi}(\cdot)\right\|_\infty+\|\dot{\boldsymbol{\psi}}(\cdot)\|_2\right)\right]^2
\end{aligned}\tag{5.56}$$

for any initial condition $\boldsymbol{\xi}\in\mathbb{R}^n$ and $\boldsymbol{\psi}(\cdot)\in\mathcal{A}([-r_2,0);\mathbb{R}^\nu)$ in (5.1), which is derived via (5.27) and the property of quadratic forms: $\forall X\in\mathbb{S}^n,\exists\lambda>0:\forall\mathbf{x}\in\mathbb{R}^n\setminus\{\mathbf{0}\},\mathbf{x}^\top(\lambda I_n-X)\mathbf{x}>0$. Then (5.56) implies that (5.43) satisfies the rightmost inequality in (5.19).

Now assume the inequalities in (5.35) are satisfied. Apply (5.27) with appropriate $\mathbf{f}(\cdot)$ with $\varpi(\tau)=1$ to the integrals in (5.43) at $t=t_0$ and consider the initial conditions in (5.1), we have

$$\begin{aligned}
&\int_{-r_1}^{0}\boldsymbol{\psi}^\top(\tau)Q_1\boldsymbol{\psi}(\tau)\mathsf{d}\tau \ge \int_{-r_1}^{0}\boldsymbol{\psi}^\top(\tau)\acute{F}_d^\top(\tau)\mathsf{d}\tau\left(\acute{\mathsf{F}}_d\otimes Q_1\right)\int_{-r_1}^{0}\acute{F}_d(\tau)\boldsymbol{\psi}(\tau)\mathsf{d}\tau,\\
&\int_{-r_2}^{-r_1}\boldsymbol{\psi}^\top(\tau)Q_2\boldsymbol{\psi}(\tau)\mathsf{d}\tau \ge \int_{-r_2}^{-r_1}\boldsymbol{\psi}^\top(\tau)\grave{F}_\delta^\top(\tau)\mathsf{d}\tau\left(\grave{\mathsf{F}}_\delta\otimes Q_2\right)\int_{-r_2}^{-r_1}\grave{F}_\delta(\tau)\boldsymbol{\psi}(\tau)\mathsf{d}\tau
\end{aligned}\tag{5.57}$$

and

$$\int_{-r_1}^{0}\dot{\boldsymbol{\psi}}^\top(\tau)S_1\dot{\boldsymbol{\psi}}(\tau)\mathsf{d}\tau \ge \int_{-r_1}^{0}\dot{\boldsymbol{\psi}}^\top(\tau)\left(\acute{\boldsymbol{\phi}}^\top(\tau)\otimes I_\nu\right)\mathsf{d}\tau\left(\acute{\mathsf{\Phi}}_{\kappa_1}\otimes S_1\right)\int_{-r_1}^{0}\left(\acute{\boldsymbol{\phi}}(\tau)\otimes I_\nu\right)\dot{\boldsymbol{\psi}}(\tau)\mathsf{d}\tau$$

$$
\begin{aligned}
&= \boldsymbol{\eta}^\top(t_0)\Pi^\top \left(G_1^\top \acute{\mathsf{\Phi}}_{\kappa_1} G_1 \otimes S_1\right) \Pi\boldsymbol{\eta}(t_0), && (5.58)\\
\int_{-r_2}^{-r_1} \dot{\boldsymbol{\psi}}^\top(\tau) S_2 \dot{\boldsymbol{\psi}}(\tau)\mathsf{d}\tau \geq \int_{-r_2}^{-r_1} \dot{\boldsymbol{\psi}}^\top(\tau) \left(\grave{\boldsymbol{\phi}}^\top(\tau) \otimes I_\nu\right) \mathsf{d}\tau \left(\grave{\mathsf{\Phi}}_{\kappa_2} \otimes S_2\right) &\int_{-r_2}^{-r_1} \left(\grave{\boldsymbol{\phi}}(\tau) \otimes I_\nu\right) \dot{\boldsymbol{\psi}}(\tau)\mathsf{d}\tau \\
&= \boldsymbol{\eta}^\top(t_0)\Pi^\top \left(G_2^\top \grave{\mathsf{\Phi}}_{\kappa_2} G_2 \otimes S_2\right) \Pi\boldsymbol{\eta}(t_0) && (5.59)
\end{aligned}
$$

which are derived via the relations in (5.51) and (5.52). Furthermore, apply (5.24) again with appropriate weight functions to the integrals $\int_{-r_1}^{0}(r_1+\tau)\dot{\boldsymbol{y}}^\top(t+\tau)U_1\dot{\boldsymbol{y}}(t+\tau)\mathsf{d}\tau$ and $\int_{-r_2}^{-r_1}(r_2+\tau)\dot{\boldsymbol{y}}^\top(t+\tau)U_2\dot{\boldsymbol{y}}(t+\tau)\mathsf{d}\tau$ for $t = t_0$ in (5.43) with $\mathsf{f}(\tau) = \acute{\boldsymbol{g}}(\tau)$, $\mathsf{f}(\tau) = \grave{\boldsymbol{g}}(\tau)$, respectively. Then it yields

$$
\begin{aligned}
\int_{-r_1}^{0} (r_1+\tau)\dot{\boldsymbol{\psi}}^\top(\tau) U_1 \dot{\boldsymbol{\psi}}(\tau)\mathsf{d}\tau &\geq [*]\left(\acute{\mathsf{G}}_{p_1} \otimes U_1\right) \int_{-r_1}^{0} (\tau + r_1)\left(\acute{\boldsymbol{g}}(\tau) \otimes I_\nu\right) \dot{\boldsymbol{\psi}}(\tau)\mathsf{d}\tau \\
&= \boldsymbol{\eta}^\top(t_0)\Pi^\top \left(H_1^\top \acute{\mathsf{G}}_{p_1} H_1 \otimes U_1\right) \Pi\boldsymbol{\eta}(t_0) \\
\int_{-r_2}^{-r_1} (r_2+\tau)\dot{\boldsymbol{\psi}}^\top(\tau) U_2 \dot{\boldsymbol{\psi}}(\tau)\mathsf{d}\tau &\geq [*]\left(\grave{\mathsf{G}}_{p_2} \otimes U_2\right) \int_{-r_2}^{-r_1} (\tau + r_2)\left(\grave{\boldsymbol{g}}(\tau) \otimes I_\nu\right) \dot{\boldsymbol{\psi}}(\tau)\mathsf{d}\tau \\
&= \boldsymbol{\eta}^\top(t_0)\Pi^\top \left(H_2^\top \grave{\mathsf{G}}_{p_2} H_2 \otimes U_2\right) \Pi\boldsymbol{\eta}(t_0)
\end{aligned} \tag{5.60}
$$

for any initial condition $\boldsymbol{\xi} \in \mathbb{R}^n$ and $\boldsymbol{\psi}(\cdot) \in \mathcal{A}\left([-r_2, 0); \mathbb{R}^\nu\right)$ in (5.1), where $H_1$ and $H_2$ are given in (5.41) and obtained by the relations

$$
\begin{aligned}
&\int_{-r_1}^{0} (\tau + r_1)\left(\acute{\boldsymbol{g}}(\tau) \otimes I_\nu\right) \dot{\boldsymbol{\psi}}(\tau)\mathsf{d}\tau = r_1\left(\acute{\boldsymbol{g}}(0) \otimes I_\nu\right)\boldsymbol{\psi}(0) \\
&\quad - \left(N_1 \otimes I_\nu\right)\int_{-r_1}^{0} \left(\acute{\boldsymbol{f}}(\tau) \otimes I_\nu\right) \boldsymbol{\psi}(\tau)\mathsf{d}\tau = \left(H_1 \otimes I_\nu\right)\boldsymbol{\eta}(t_0)
\end{aligned} \tag{5.61}
$$

$$
\begin{aligned}
&\int_{-r_2}^{-r_1} (\tau + r_2)\left(\grave{\boldsymbol{g}}(\tau) \otimes I_\nu\right) \dot{\boldsymbol{\psi}}(\tau)\mathsf{d}\tau = (r_2 - r_1)\left(\grave{\boldsymbol{g}}(-r_1) \otimes I_\nu\right)\boldsymbol{\psi}(-r_1) \\
&\qquad - \left(N_2 \otimes I_\nu\right)\int_{-r_2}^{-r_1} \left(\grave{\boldsymbol{f}}(\tau) \otimes I_\nu\right) \boldsymbol{\psi}(\tau)\mathsf{d}\tau = \left(H_2 \otimes I_\nu\right)\boldsymbol{\eta}(t_0)
\end{aligned} \tag{5.62}
$$

via (5.33) and the properties of the Kronecker product in (2.1a) and (2.1b).

With (5.35), applying (5.57)–(5.60) to (5.43) with $t = t_0$ and considering the initial conditions in (5.1) can conclude that (5.19) is satisfied if (5.34) and (5.35) hold. This shows that feasible solutions to (5.34)–(5.36) imply the existence of the functional in (5.43) satisfying (5.19)–(5.21). ∎

**Remark 5.10.** By allowing $m, q$ to be zero, Theorem 5.1 can address the problem of conducting stability analysis without performance requirements. Moreover, if $\acute{\boldsymbol{f}}(\cdot)$ and $\grave{\boldsymbol{f}}(\cdot)$ contain only Legendre polynomials, then Theorem 5.1 with (5.13) generalizes the two delay channel version of the stability results in [470]. (Note that the method in [470] addresses systems with only a single delay channel)

**Remark 5.11.** If one wants to increase the values of $d$ and $\delta$ in (5.43) to incorporate more functions in the DD kernels in (5.44), then additional zeros must be introduced to the coefficient matrices $A_4, A_5$ and $C_4, C_5$ in (5.7) in order to make (5.43) consistent with (5.7). In conclusion, there are no upper bounds on the values of $d$ and $\delta$. Finally, (5.43) generalizes the LKF in [470], which only considers Legendre polynomials for the integral terms in (5.44).

**Remark 5.12.** If condition (5.33) is not imposed on $\acute{\boldsymbol{f}}(\cdot)$ and $\grave{\boldsymbol{f}}(\cdot)$, then dissipative conditions can still be derived, but the inequalities in (5.60) can no longer be considered. In that case, the constraints (5.35) and (5.36) remain the same, and (5.34) is changed into

$$
P + \left(\mathsf{O}_{n+2\nu} \oplus \left[\acute{\mathsf{F}}_d \otimes Q_1\right] \oplus \left[\grave{\mathsf{F}}_\delta \otimes Q_2\right]\right) + \Pi^\top \left(G_1^\top \acute{\mathsf{\Phi}}_{\kappa_1} G_1 \otimes S_1 + G_2^\top \grave{\mathsf{\Phi}}_{\kappa_2} G_2 \otimes S_2\right)\Pi \succ 0. \tag{5.63}
$$

**Remark 5.13.** Note that the position of the error matrices $\acute{\mathsf{E}}_d$ and $\grave{\mathsf{E}}_\delta$ in $\widetilde{\mathbf{\Omega}} \prec 0$ in (5.36) may cause numerical problems if the eigenvalues of $\acute{\mathsf{E}}_d$ and $\grave{\mathsf{E}}_\delta$ are too small. To circumvent this potential issue, we can apply congruent transformations to $\widetilde{\mathbf{\Omega}} \prec 0$ which concludes that $\widetilde{\mathbf{\Omega}} \prec 0$ holds if and only if

$$[*]\widetilde{\mathbf{\Omega}}\left[I_{m+q+n+5\nu+\varrho\nu} \oplus \left(\eta_1 \acute{\mathsf{E}}_d^{-\frac{1}{2}} \otimes I_\nu\right) \oplus \left(\eta_2 \grave{\mathsf{E}}_\delta^{-\frac{1}{2}} \otimes I_\nu\right)\right] \prec 0 \tag{5.64}$$

holds where $\eta_1; \eta_2 \in \mathbb{R}$ are given values. Note that the diagonal elements of the transformed matrix in (5.64) are no longer associated with the error terms appear at off-diagonal elements, hence one can use the inequality (5.64) instead of (5.36).

**Remark 5.14.** The hypothesis of $r_2 > r_1 > 0$ in Theorem 5.1 indicates that there are no redundant matrix parameters in (5.43) since two genuine delay channels are considered therein and (5.44) and (5.18) contain no zeros vectors. In the cases of $r_1 = 0$ or $r_2 = r_1$, one only needs to consider one delay channel so that the corresponding (5.44), (5.18) and (5.43) can be simplified. Note that we do not present the corresponding dissipativity and stability condition for $r_1 = 0$ or $r_2 = r_1$ in this chapter since it can be effortlessly derived based on the proof of Theorem 5.1 with a simplified (5.43).

In the following corollary, we show that a hierarchy of the stability condition in Theorem 5.1 can be established with respect to $\grave{\boldsymbol{\phi}}(\cdot)$ and its dimension when certain conditions are met.

**Corollary 5.3.** *Let all the functions and the parameters in (5.3)–(5.12) be given where* $\grave{\boldsymbol{\phi}}(\tau) := \mathbf{Col}_{i=1}^{\kappa_2} \grave{\phi}_i(\tau)$ *with* $\{\grave{\phi}_i(\cdot)\}_{i=1}^{\kappa_2} \subset \{\grave{\phi}_i(\cdot)\}_{i=1}^{\infty} \subset \mathcal{C}^1(\mathcal{I}_2; \mathbb{R})$ *satisfying*

$$\exists \varkappa \in \mathbb{N},\ \ \forall \kappa_2 \in \{j \in \mathbb{N} : j \leq \varkappa\},\ \exists! M_4 \in \mathbb{R}^{\kappa_2 \times \delta},\ \ \frac{\mathsf{d}}{\mathsf{d}\tau} \underset{i=1}{\overset{\kappa_2}{\mathbf{Col}}}\, \grave{\phi}_i(\tau) = M_4 \grave{\boldsymbol{f}}(\tau) \tag{5.65}$$

$$\forall \kappa_2 \in \mathbb{N},\ \ \grave{\Phi}_{\kappa_2} = \int_{-r_2}^{-r_1} \underset{i=1}{\overset{\kappa_2}{\mathbf{Col}}}\, \grave{\phi}_i(\tau) \underset{i=1}{\overset{\kappa_2}{\mathbf{Row}}}\, \grave{\phi}_i(\tau) \mathsf{d}\tau = \bigoplus_{j=1}^{\kappa_2} \grave{\varphi}_j, \quad \grave{\varphi}_j^{-1} = \int_{-r_2}^{-r_1} \grave{\phi}_j^2(\tau) \mathsf{d}\tau. \tag{5.66}$$

*Now given* $\acute{\boldsymbol{g}}(\cdot)$, $\grave{\boldsymbol{g}}(\cdot)$ *and* $N_1$, $N_2$ *in* Theorem 5.1, *we have*

$$\forall \kappa_2 \in \{j \in \mathbb{N} : j \leq \varkappa\},\ \ \mathcal{G}_{\kappa_2} \subseteq \mathcal{G}_{\kappa_2+1} \tag{5.67}$$

*where* $\varkappa \in \mathbb{N}$ *is given and*

$$\mathcal{G}_{\kappa_2} := \left\{(r_1, r_2)\ \middle|\ r_1 > 0, r_2 > r_1\ \ \&\ \ (5.34)\text{–}(5.36)\ \textit{hold}\ \ \&\ \ P \in \mathbb{S}^l, Q_1; Q_2; R_1; R_2; S_1; S_2; U_1; U_2 \in \mathbb{S}^\nu\right\}$$

*with* $l := n + 2\nu + (d + \delta)\nu$.

*Proof.* Given $r_2 > r_1 > 0$ and all the parameters in (5.3)–(5.7) and (5.11), (5.12), let $\mathbf{Col}(r_1, r_2) \in \mathcal{G}_{\kappa_2}$ with $\mathcal{G}_{\kappa_2} \neq \varnothing$ which implies that there exist feasible solutions to (5.34)–(5.36). Consider the situation when the dimensions and elements of $\acute{\boldsymbol{f}}(\tau)$, $\grave{\boldsymbol{f}}(\tau)$, $\acute{\boldsymbol{\phi}}(\tau)$, $\acute{\boldsymbol{g}}(\tau)$ and $\grave{\boldsymbol{g}}(\tau)$ are all fixed, and let $P \in \mathbb{S}^l$ and $Q_1; Q_2; R_1; R_2; S_1; S_2; U_1; U_2 \in \mathbb{S}^\nu$ be a given feasible solution to $\mathbf{P}_{\kappa_2} \succ 0$, (5.35) and $\widetilde{\mathbf{\Omega}}_{\kappa_2} \prec 0$ at $\kappa_2$. Note that the matrix $G_2$ and $\Phi_{\kappa_2}$ in (5.36) are indexed by the value of $\kappa_2$. Given (5.35), we show that the corresponding feasible solutions to (5.34) and (5.36) at $\kappa_2 + 1$ exist if feasible solutions to (5.34) and (5.36) at $\kappa_2$ exist, which proves (5.67).

The conditions in (5.66) indicate that $\grave{\phi}_i(\cdot)$ are orthogonal functions with respect to the weight function $\varpi(\tau) = 1$ over $\mathcal{I}_2$, Assume $\kappa_2 + 1 \leq \varkappa$ and by the structure of $G_2$ in (5.39) with (5.65) and (5.66), we have

$$G_{2,\kappa_2+1}^\top \grave{\Phi}_{\kappa_2+1} G_{2,\kappa_2+1} = [*]\begin{bmatrix} \grave{\Phi}_{\kappa_2} & \mathbf{0}_{\kappa_2+1} \\ * & \grave{\varphi}_{\kappa_2+1} \end{bmatrix}\begin{bmatrix} G_{2,\kappa_2} \\ \mathbf{g}_{\kappa_2+1}^\top \end{bmatrix} = G_{2,\kappa_2}^\top \grave{\Phi}_{\kappa_2} G_{2,\kappa_2} + \grave{\varphi}_{\kappa_2+1} \mathbf{g}_{\kappa_2+1} \mathbf{g}_{\kappa_2+1}^\top, \tag{5.68}$$

where $\mathbf{g}_{\kappa_2+1} \in \mathbb{R}^{3+d+\delta}$ can be easily determined by the structure of $G_2$ with (5.52) and (5.65), and $G_{2,\kappa_2+1}$ denotes the corresponding $G_2$ at $\kappa_2 + 1$. It is important to note that there is no increase in the dimension indexes $d$, $\delta$, $p_1$ and $p_2$. By (5.68) and considering the structure of the matrix inequalities in (5.34) and (5.36), we have

$$\begin{aligned}\mathbf{P}_{\kappa_2+1} &= \mathbf{P}_{\kappa_2} + \Pi^\top \left(\grave{\varphi}_{\kappa_2+1}\mathbf{g}_{\kappa_2+1}\mathbf{g}_{\kappa_2+1}^\top \otimes S_2\right)\Pi \\ \widetilde{\mathbf{\Omega}}_{\kappa_2+1} &= \widetilde{\mathbf{\Omega}}_{\kappa_2} + \left(\mathsf{O}_{q+2\nu} \oplus \Pi^\top \left(\grave{\varphi}_{\kappa_2+1}\mathbf{g}_{\kappa_2+1}\mathbf{g}_{\kappa_2+1}^\top \otimes U_2\right)\Pi \oplus \mathsf{O}_{\mu\nu}\right).\end{aligned} \tag{5.69}$$

Since $\grave{\varphi}_{\kappa_2+1} > 0$, $\mathbf{g}_{\kappa_2+1}\mathbf{g}_{\kappa_2+1}^\top \succeq 0$ and $S_2 \succeq 0$, $U_2 \succeq 0$ in (5.35), it is clearly to see that the feasible solutions to $\mathbf{P}_{\kappa_2} \succ 0$, $\widetilde{\mathbf{\Omega}}_{\kappa_2} \succ 0$ guarantee the existence of a feasible solution to $\mathbf{P}_{\kappa_2+1} \succ 0$, $\widetilde{\mathbf{\Omega}}_{\kappa_2+1} \succ 0$ given the prerequisites of Corollary 5.3. This finishes the proof. ■

**Remark 5.15.** A hierarchical pattern of the LMIs in Theorem 5.1 can also be established for the situation when $\acute{\boldsymbol{\phi}}(\cdot)$ contains orthogonal functions which satisfy appropriate constraints resembling (5.65) and (5.66). Note that the corresponding hierarchy result can be derived without the use of congruent transformations, as the dimensions of $\mathbf{P}$ in (5.34) and $\widetilde{\mathbf{\Omega}}$ in (5.36) are not related to the dimensions of $\acute{\boldsymbol{\phi}}(\tau) \in \mathbb{R}^{\kappa_1}$.

On the other hand, a hierarchy of the stability condition in Theorem 5.1 can also be established with respect to $\grave{\boldsymbol{g}}(\cdot)$ and its dimensions.

**Corollary 5.4.** *Given the functions with the parameters in (5.3)–(5.12), let $\acute{\boldsymbol{g}}(\cdot)$, $\grave{\boldsymbol{g}}(\cdot)$ and $N_1$, $N_2$ in Theorem 5.1 be given where $\grave{\boldsymbol{g}}(\tau) = \mathsf{Col}_{i=1}^{p_2}\, \grave{g}_i(\tau)$ with $\{\grave{g}_i(\cdot)\}_{i=1}^{p_2} \subset \{\grave{g}_i(\cdot)\}_{i=1}^{\infty} \subset \mathcal{C}^1(\mathcal{I}_2;\mathbb{R})$ satisfying*

$$\exists \alpha \in \mathbb{N},\ \ \forall p_2 \in \{j \in \mathbb{N} : j \leq \alpha\},\ \exists! N_2 \in \mathbb{R}^{p_2\times\delta},\ \ (r_2+\tau)\frac{\mathsf{d}}{\mathsf{d}\tau}\underset{i=1}{\overset{p_2}{\mathsf{Col}}}\, \grave{g}_i(\tau) = N_2\grave{\boldsymbol{g}}(\tau) \tag{5.70}$$

$$\forall p_2 \in \mathbb{N},\ \ \grave{\mathsf{G}}_{p_2} = \int_{-r_2}^{-r_1} \underset{i=1}{\overset{p_2}{\mathsf{Col}}}\, \grave{g}_i(\tau)\, \underset{i=1}{\overset{p_2}{\mathsf{Row}}}\, \grave{g}_i(\tau)\mathsf{d}\tau = \bigoplus_{j=1}^{p_2} \grave{\mathsf{g}}_j,\quad \grave{\mathsf{g}}_j^{-1} = \int_{-r_2}^{-r_1} (\tau+r_2)\grave{g}_j^2(\tau)\mathsf{d}\tau. \tag{5.71}$$

*Then we have*

$$\forall p_2 \in \{j \in \mathbb{N} : j \leq \alpha\},\ \ \mathcal{H}_{p_2} \subseteq \mathcal{H}_{p_2+1} \tag{5.72}$$

*where $\alpha \in \mathbb{N}$ is given and*

$$\mathcal{H}_{p_2} := \left\{(r_1, r_2)\ \middle|\ r_1 > 0, r_2 > r_1\ \ \&\ \ (5.34)\text{–}(5.36)\ \textit{hold}\ \ \&\ \ P \in \mathbb{S}^l, Q_1; Q_2; R_1; R_2; S_1; S_2; U_1; U_2 \in \mathbb{S}^\nu\right\}$$

*with $l := n + 2\nu + (d+\delta)\nu$.*

*Proof.* The proof is similar to the proof of Corollary 5.3, with the exception that for Corollary 5.4, one only needs to consider the increase of the value of $p_2$ instead of $\kappa_2$ in Corollary 5.3. Given $r_2 > r_1 > 0$ with all the parameters in (5.3)–(5.7) and (5.11) and (5.12), let $\mathsf{Col}(r_1, r_2) \in \mathcal{H}_{p_2}$ with $\mathcal{H}_{p_2} \neq \varnothing$, indicating that there exist feasible solutions to (5.34)–(5.36). Let the dimensions and elements of $\acute{\boldsymbol{f}}(\tau)$, $\grave{\boldsymbol{f}}(\tau)$, $\acute{\boldsymbol{\phi}}(\tau)$, $\grave{\boldsymbol{\phi}}(\tau)$ and $\acute{\boldsymbol{g}}(\tau)$ to be all fixed, and let $P \in \mathbb{S}^l$ and $Q_1; Q_2; R_1; R_2; S_1; S_2; U_1; U_2$ be a given feasible solution to $\mathbf{P}_{p_2} \succ 0$, (5.35) and $\widetilde{\mathbf{\Omega}} \prec 0$ at $p_2$. Note that the matrix $H_2$ and $\mathsf{G}_{p_2}$ in (5.36) are indexed by the value of $\kappa_2$ whereas $\widetilde{\mathbf{\Omega}} \prec 0$ is not related to $\acute{\boldsymbol{g}}(\tau)$ and $\grave{\boldsymbol{g}}(\tau)$ or their dimensions $p_1$, $p_2$. Given (5.35), we show that in the following that feasible solutions to (5.34) and (5.36) at $p_2 + 1$ exist if feasible solutions to (5.34) and (5.36) at $p_2$ exist, which leads to (5.67).

The constraints in (5.71) show that $\grave{g}_i(\cdot)$ contains functions which are orthogonal with respect to the weight function $\varpi(\tau) = (\tau + r_2)$ over $\mathcal{I}_2$. Suppose $p_2 + 1 \leq \alpha$. Now by the structure of $H_2$ in

(5.41) and (5.70) and (5.71), we have

$$H_{2,p_2+1}^\top \grave{\mathsf{G}}_{p_2+1} H_{2,p_2+1} = [*]\begin{bmatrix}\grave{\mathsf{G}}_{p_2} & \mathbf{0}_{p_2}\\ * & \grave{\mathsf{g}}_{p_2+1}\end{bmatrix}\begin{bmatrix}H_{2,p_2}\\ \mathbf{h}_{p_2+1}^\top\end{bmatrix} = H_{2,p_2}^\top \grave{\mathsf{G}}_{p_2} H_{2,p_2} + \grave{\mathsf{g}}_{p_2+1}\mathbf{h}_{p_2+1}\mathbf{h}_{p_2+1}^\top \tag{5.73}$$

where $\mathbf{h}_{p_2+1} \in \mathbb{R}^{3+d+\delta}$ can be easily determined by the structure of $H_2$ with (5.62) and (5.70), and $H_{2,p_2+1}$ denotes the corresponding $H_2$ at $p_2+1$. Note that here the values of the dimension indexes $d$, $\delta$, $\kappa_1$, $\kappa_2$ and $p_1$ remain unchanged.

By (5.73) and considering the structure of $\mathbf{P} \succ 0$ in (5.34), it follows that

$$\mathbf{P}_{p_2+1} = \mathbf{P}_{p_2} + \Pi^\top\left(\grave{\mathsf{g}}_{p_2+1}\mathbf{h}_{p_2+1}\mathbf{h}_{p_2+1}^\top \otimes U_2\right)\Pi. \tag{5.74}$$

Since $\grave{\mathsf{g}}_{p_2+1} > 0$, $\mathbf{h}_{p_2+1}\mathbf{h}_{p_2+1}^\top \succeq 0$ with $U_2 \succeq 0$ in (5.35), one can conclude that the feasible solutions to $\mathbf{P}_{p_2} \succ 0$ imply the existence of the feasible solution to $\mathbf{P}_{p_2+1} \succ 0$ given the prerequisites in Corollary 5.4. On the other hand, since the inequality in (5.36) is not related to $\grave{\boldsymbol{g}}(\tau)$, hence $\widetilde{\boldsymbol{\Omega}} \prec 0$ remains unchanged at $p_2+1$. This finishes the proof. ■

**Remark 5.16.** By following the strategy used in the proof of Corollary 5.4, a hierarchy of conditions in Theorem 5.1 can also be established when $\acute{\boldsymbol{g}}(\cdot)$ contains orthogonal functions that satisfy appropriate constraints similar to those in (5.65) and (5.66). Note that the dimensions of $\mathbf{P}$ in (5.34) and $\widetilde{\boldsymbol{\Omega}}$ in (5.36) are not related to the dimensions of $\acute{\boldsymbol{g}}(\tau) \in \mathbb{R}^{p_1}$.

## 5.5 Numerical Examples

We employ the proposed methods to various numerical examples to demonstrate their effectiveness. All examples were tested in MATLAB© in conjunction with the optimization parser Yalmip [603] and SDP numerical solver SDPT3 [512]. Note that we do not prescribe fixed positive eigenvalue margins for positive semi-definite (PSD) variables. Instead, the validity of a feasible result (this principle applies equally to minimization programs) is confirmed by verifying that all eigenvalues of the resulting PSD variables are strictly positive.

### 5.5.1 Stability analysis of an LDDS

Consider the following LDDS

$$\begin{aligned}\dot{x}(t) &= 0.33x(t) - 5\int_{-r}^{0}\sin(\cos(12\tau))x(t+\tau)\mathsf{d}\tau\\ &= 0.33x(t) - \begin{bmatrix}5 & \mathbf{0}^\top\end{bmatrix}\int_{-r}^{0}\begin{bmatrix}\varphi_1(\tau)\\ \boldsymbol{f}(\tau)\end{bmatrix}x(t+\tau)\mathsf{d}\tau, \quad t \geq t_0\end{aligned} \tag{5.75}$$

with any $t_0 \in \mathbb{R}$, where $\varphi_1(\tau) = \sin(\cos(12\tau))$. The corresponding state space matrices of (5.1) for (5.75) and (5.7) are $A_1 = 0.33$ and $A_3 = -\begin{bmatrix}5 & \mathbf{0}^\top\end{bmatrix}$ and the rest of the state space matrices in (5.1) is zero with $m = q = 0$.

For the example in (5.75), we consider two cases for $\boldsymbol{f}(\cdot)$. The first one is Legendre polynomials $\boldsymbol{f}(\tau) = \boldsymbol{\ell}_d(\tau) = \mathbf{Col}_{i=0}^{d}\,\ell_i(\tau)$ with

$$\ell_d(\tau) := \sum_{k=0}^{d}\binom{d}{k}\binom{d+k}{k}\left(\frac{\tau}{r}\right)^k \tag{5.76}$$

with $\mathsf{F}_1^{-1} = \int_{-r}^{0} \boldsymbol{\ell}_d(\tau)\boldsymbol{\ell}_d^{\top}(\tau)\mathrm{d}\tau = r^{-1}\bigoplus_{i=0}^{d} 2i+1$ and the corresponding $M_1$ in (5.3) can be easily constructed. The second one

$$\boldsymbol{f}(\tau) = \boldsymbol{h}_d(\tau) = \mathbf{Col}\left[1,\ \mathbf{Col}_{i=1}^{d/2} \sin 12i\tau,\ \mathbf{Col}_{i=1}^{d/2} \cos 12i\tau\right]$$

contains trigonometric functions paired with $M_1 = 0 \oplus \begin{bmatrix} \mathsf{O}_{d/2} & \bigoplus_{i=1}^{d/2} 12i \\ -\bigoplus_{i=1}^{d/2} 12i & \mathsf{O}_{d/2} \end{bmatrix}$ that satisfies the first relation in (5.3). Note that the $d > 0$ in $\boldsymbol{h}_d(\tau)$ must be even, and the functions in $\boldsymbol{h}_d(\tau)$ are not orthogonal over $[-r, 0]$. Therefore, the corresponding $\mathsf{F}$ for $\boldsymbol{h}_d(\tau)$ is not a diagonal matrix. Since $0.33 > 0$ in (5.75) that is non-Hurwitz, the method in [607] cannot be applied. Furthermore, since $\varphi_1(\tau) = \sin(\cos(12\tau))$ does not satisfy the differentiation closure property as in (5.3), the approach in [169] is not applicable to (5.75).

We utilized spectrum methods from [202] to (5.75) with discretization index $M = 200$. The resulting numerical results reveal that system (5.75) is stable within intervals: $[0.093, 0.169]$, $[0.617, 0.692]$, $[1.14, 1.216]$, $[1.664, 1.739]$, $[2.188, 2.263]$ and $[2.711, 2.787]$.

For this example, we employ a single delay version of Theorem 5.1 to (5.75), which is derived via the LKF

$$\mathsf{v}(\boldsymbol{x}(t), \boldsymbol{\psi}(\cdot)) = \boldsymbol{\eta}^{\top}(t)P\boldsymbol{\eta}(t) + \int_{-r}^{0} \boldsymbol{y}^{\top}(t+\tau)\left[Q + (\tau + r)R\right]\boldsymbol{y}(t+\tau)\mathrm{d}\tau \tag{5.77}$$

as a simplified version of (5.43), where $P \in \mathbb{S}^{n+(d+1)\nu}$, $Q; R \in \mathbb{S}^{\nu}$ and $\boldsymbol{\eta}(t) := \begin{bmatrix} \boldsymbol{x}(t) \\ \int_{-r}^{0} F_d(\tau)\boldsymbol{y}(t+\tau)\mathrm{d}\tau \end{bmatrix}$ with $F_d(\tau) = \boldsymbol{f}(\tau) \otimes I_\nu$. Furthermore, the corresponding $\boldsymbol{\vartheta}(t)$ in (5.18) and (5.53) is defined as $\boldsymbol{\vartheta}(t) := \mathbf{Col}\left[\boldsymbol{x}(t),\ \boldsymbol{y}(t-r),\ \int_{-r}^{0} F_d(\tau)\boldsymbol{y}(t+\tau)\mathrm{d}\tau\right]$. Now apply the corresponding stability condition derived by (5.77) with a one delay version congruent transformation (5.64) with $\eta_1 = 1$ to (5.75) with $\boldsymbol{f}(\tau) = \boldsymbol{\ell}_d(\tau)$ and $\boldsymbol{f}(\tau) = \boldsymbol{h}_d(\tau)$, respectively. The results concerning detectable delay margins are listed in Table 5.1 and Table 5.2. Note that the values of $N$ and $d$ in these tables are presented when the margins of the stable delay intervals can be determined by the numerical results produced by Theorem 5.1 or the method in [470]. Note that also the results in Table 5.1 and Table 5.2 associated with [202, 203] are calculated with $M = 200$. Finally, NoDVs in Table 5.1 and Table 5.2 refers to the number of decision variables.

| [202, 203] | $[0.093, 0.169]$ | $[0.617, 0.692]$ | $[1.14, 1.216]$ |
|---|---|---|---|
| Theorem 1 | $d = 3$ (NoDVs: 17) | $d = 6$ (NoDVs: 38) | $d = 10$ (NoDVs: 80) |
| $\boldsymbol{f}(\tau)$ | $\boldsymbol{\ell}_d(\tau)$ | $\boldsymbol{h}_d(\tau)$ | $\boldsymbol{h}_d(\tau)$ |
| [470] | $N = 3$ (NoDVs: 17) | $N = 11$ (NoDVs: 93) | $N = 23$ (NoDVs: 327) |

**Table 5.1:** Testing of stable delay margins

| [202] | $[1.664, 1.739]$ | $[2.188, 2.263]$ | $[2.711, 2.787]$ |
|---|---|---|---|
| Theorem 5.1 | $d = 10$ (NoDVs: 80) | $d = 10$ (NoDVs: 80) | $d = 10$ (NoDVs: 80) |
| $\boldsymbol{f}(\tau)$ | $\boldsymbol{h}_d(\tau)$ | $\boldsymbol{h}_d(\tau)$ | $\boldsymbol{h}_d(\tau)$ |
| [470] | – | – | – |

**Table 5.2:** Testing of stable delay margins

In Table 5.1 and Table 5.2, the results correspond to [470] were produced by Theorem 5.1 via (5.77) with $\boldsymbol{f}(\tau) = \boldsymbol{\ell}_d(\tau)$ and $d = N$, which is essentially equivalent to the method in [470]. At $N = 25$, the margins of the stable delay intervals $[1.664, 1.739]$, $[2.188, 2.263]$ and $[2.711, 2.787]$ are still undetectable using the polynomial approximation approach proposed in [470]. For $N > 25$, our experiments show that the computational time becomes too long for accurately computing the approximation coefficient and error term using the `vpaintegral` in MATLAB©. In addition, the function `integral` in MATLAB© does not offer a viable alternative for calculating the approximation coefficient and error term due to its limited numerical accuracy. The results in Table 5.1 and Table 5.2 can be explained by the fact that $\varphi_1(\tau) = \sin(\cos(12\tau))$, $\tau \in [a, b]$ is not easily approximated by polynomials when the length of $[a, b]$ becomes relatively large. As a consequence, we have shown the advantage of our method over the one in [470] when it comes to the stability analysis of (5.75).

### 5.5.2 Stability and dissipativity analysis with DDs

Consider a system in the form of (5.1) with $r_1 = 2$, $r_2 = 4.05$ and the state space parameters

$$\begin{gathered}
A_1 = \begin{bmatrix} 0.01 & 0 \\ 0 & -3 \end{bmatrix}, A_2 = \begin{bmatrix} 0 & 0.1 \\ 0.2 & 0 \end{bmatrix}, A_3 = \begin{bmatrix} -0.1 & 0 \\ 0 & -0.2 \end{bmatrix}, A_6 = I_2, A_7 = A_8 = \mathsf{O}_2,\ D_1 = \begin{bmatrix} 0.2 \\ 0.3 \end{bmatrix}, \\
A_4\left(\begin{bmatrix} \boldsymbol{\varphi}_1(\tau) \\ \acute{\boldsymbol{f}}(\tau) \end{bmatrix} \otimes I_2\right) = \begin{bmatrix} 3\sin(18\tau) & -0.3\mathsf{e}^{\cos(18\tau)} \\ 0 & 3\sin(18\tau) \end{bmatrix},\ A_5\left(\begin{bmatrix} \boldsymbol{\varphi}_2(\tau) \\ \grave{\boldsymbol{f}}(\tau) \end{bmatrix} \otimes I_2\right) = \begin{bmatrix} -10\cos(18\tau) & 0 \\ 0.5\mathsf{e}^{\sin(18\tau)} & -10\cos(18\tau) \end{bmatrix}, \\
C_1 = \begin{bmatrix} -0.1 & 0.2 \\ 0 & 0.1 \end{bmatrix},\ C_2 = \begin{bmatrix} -0.1 & 0 \\ 0 & 0.2 \end{bmatrix},\ C_3 = \begin{bmatrix} 0 & 0.1 \\ -0.1 & 0 \end{bmatrix},\ D_2 = \begin{bmatrix} 0.12 \\ 0.1 \end{bmatrix}, \\
C_4\left(\begin{bmatrix} \boldsymbol{\varphi}_1(\tau) \\ \acute{\boldsymbol{f}}(\tau) \end{bmatrix} \otimes I_\nu\right) = 0.1 \oplus 0,\ C_5\left(\begin{bmatrix} \boldsymbol{\varphi}_2(\tau) \\ \grave{\boldsymbol{f}}(\tau) \end{bmatrix} \otimes I_\nu\right) = 0.2 \oplus 0.1,\ C_6 = \begin{bmatrix} 0 & 0 \\ 0 & 0.1 \end{bmatrix}, C_7 = \begin{bmatrix} 0.2 & 0 \\ 0 & 0 \end{bmatrix},
\end{gathered} \tag{5.78}$$

with $\boldsymbol{\varphi}_1(\tau) = \boldsymbol{\varphi}_2(\tau) = \begin{bmatrix} \mathsf{e}^{\sin(18\tau)} & \mathsf{e}^{\cos(18\tau)} \end{bmatrix}^\top$ and $n = m = 2$, $q = 1$. We find out that the system with (5.78) is stable using the spectral method toolbox developed in [203]. Moreover, we aim to minimize $\mathcal{L}^2$ gain $\gamma$ as the performance criterion, which corresponds to

$$\gamma > 0, \quad J_1 = -\gamma I_2, \quad \widetilde{J} = I_2, \quad J_2 = \mathbf{0}_2, \quad J_3 = \gamma \tag{5.79}$$

for the supply rate function in (5.22).

Even if one assumes the method in [463] can be extended to handle systems with multiple delay channels, this method still remains inapplicable in this context, given that $A_1$ is not a Hurwitz matrix. In addition, since $\boldsymbol{\varphi}_1(\tau) = \boldsymbol{\varphi}_2(\tau)$ fails to satisfy the 'differentiation closure' property in (5.3), simply extending the corresponding conditions in [169] for a linear CDDS may not solve the problem of dissipativity and stability analysis, even if a multiple distinct delays version of the method in [169] is derivable.

Let

$$\acute{\boldsymbol{f}}(\tau) = \acute{\boldsymbol{\phi}}(\tau) = \begin{bmatrix} 1 \\ \mathsf{Col}_{i=1}^{d} \sin 18i\tau \\ \mathsf{Col}_{i=1}^{d} \cos 18i\tau \end{bmatrix}, \quad \grave{\boldsymbol{f}}(\tau) = \grave{\boldsymbol{\phi}}(\tau) = \begin{bmatrix} 1 \\ \mathsf{Col}_{i=1}^{\delta} \sin 18i\tau \\ \mathsf{Col}_{i=1}^{\delta} \cos 18i\tau \end{bmatrix} \tag{5.80}$$

in (5.78) and (5.3), which are associated with matrices

$$M_1 = M_3 = 0 \oplus \begin{bmatrix} \mathsf{O}_d & \bigoplus_{i=1}^{d} 18i \\ -\bigoplus_{i=1}^{d} 18i & \mathsf{O}_d \end{bmatrix}, \quad M_2 = M_4 = 0 \oplus \begin{bmatrix} \mathsf{O}_\delta & \bigoplus_{i=1}^{\delta} 18i \\ -\bigoplus_{i=1}^{\delta} 18i & \mathsf{O}_\delta \end{bmatrix} \tag{5.81}$$

for the differential equations in (5.3). By taking into account $\acute{\boldsymbol{f}}(\cdot)$, $\grave{\boldsymbol{f}}(\cdot)$ in (5.80) and $\varphi_1(\tau) = \varphi_2(\tau) = \begin{bmatrix} \mathsf{e}^{\sin(18\tau)} & \mathsf{e}^{\cos(18\tau)} \end{bmatrix}^\top$ with (2.1a) and (5.11), we obtain

$$
\begin{aligned}
A_4 &= \begin{bmatrix} \mathsf{O}_2 & \begin{matrix} 0 & -0.3 \\ 0 & 0 \end{matrix} & \mathsf{O}_2 & \begin{matrix} 3 & 0 \\ 0 & 3 \end{matrix} & \mathsf{O}_{2\times(4d-2)} \end{bmatrix} \\
A_5 &= \begin{bmatrix} \begin{matrix} 0 & 0 & 0 & 0 \\ 0.5 & 0 & 0 & 0 \end{matrix} & \mathsf{O}_{2,2\delta+2} & \begin{matrix} -10 & 0 \\ 0 & -10 \end{matrix} & \mathsf{O}_{2,2\delta-2} \end{bmatrix} \\
C_4 &= \begin{bmatrix} \mathsf{O}_{2,4} & 0.1 \oplus 0 & \mathsf{O}_{2,4d} \end{bmatrix}, \quad C_5 = \begin{bmatrix} \mathsf{O}_{2,4} & 0.2 \oplus 0.1 & \mathsf{O}_{2,4\delta} \end{bmatrix}
\end{aligned}
\tag{5.82}
$$

that can satisfy (5.7) for the DDs in (5.78).

Now apply Theorem 5.1 with (5.64)[4] with $\eta_1 = \eta_2 = 1$ to the system with the parameters in (5.78)–(5.82). In this case, matrices $\acute{\Gamma}_d$, $\grave{\Gamma}_\delta$, $\acute{\mathsf{E}}_d$, $\grave{\mathsf{E}}_\delta$ and $\acute{\mathsf{F}}_d$, $\grave{\mathsf{F}}_\delta$ are computed via the function `vpaintegral` in MATLAB© that can produce results with high-numerical precisions. With $d = \delta = 1$, a feasible result can be produced with $\min \gamma = 0.64655$, requiring 196 unknowns. With $d = \delta = 2$, we obtain feasible solutions with $\min \gamma = 0.32346$, requiring 376 variables. Finally, with $d = \delta = 10$ our method can produce feasible solutions with $\min \gamma = 0.31265$, requiring 4120 variables. It is worth remarking that even with $d = \delta = 10$, which is relatively large, the computational time of $\acute{\Gamma}_d$, $\grave{\Gamma}_\delta$, $\acute{\mathsf{E}}_d$, $\grave{\mathsf{E}}_\delta$ and $\acute{\mathsf{F}}_d$, $\grave{\mathsf{F}}_\delta$ by `vpaintegral` is still acceptable (about a minute).

In contrast, let $\acute{\boldsymbol{f}}(\tau) = \boldsymbol{j}_d^{0,0}(\tau)_{-r_1}^{0} = \mathbf{Col}_{i=0}^{d}\, j_k^{0,0}(\tau)_{-r_1}^{0}$ and $\grave{\boldsymbol{f}}(\tau) = \boldsymbol{j}_d^{0,0}(\tau)_{-r_2}^{-r_1} = \mathbf{Col}_{i=0}^{d}\, j_k^{0,0}(\tau)_{-r_2}^{-r_1}$ which are Legendre polynomials associated with $\acute{\mathsf{F}}_d = r_1^{-1}\mathsf{D}_d$ and $\grave{\mathsf{F}}_\delta = r_3^{-1}\mathsf{D}_\delta$. (see (4.6) for the definition of orthogonal polynomials) The characteristics of the functions in $\varphi_1(\tau) = \varphi_2(\tau)$ indicate that they may be difficult to approximate using polynomials. Indeed, with $d = \delta = 15$ and the corresponding $A_4$,$A_5$ and $C_4$,$C_5$, Theorem 5.1 yields no feasible solutions.

The preceding examples have demonstrated an important contribution of the proposed methodology in this chapter. Namely, the general structure of $\acute{\boldsymbol{f}}(\cdot)$ and $\grave{\boldsymbol{f}}(\cdot)$ in (5.1) and (5.43) can lead to the construction of less conservative dissipative conditions for Theorem 5.1, compared to existing results where only polynomial functions can be considered for the kernel functions.

[4]Note that here we do not apply (5.34) in Theorem 5.1, see Remark 5.12 for further details.

# Chapter 6

# Dissipative Delay Range Analysis of Coupled Differential-Difference Delay Systems with Distributed Delays

## 6.1 Introduction

In the previous chapters, we assumed that delay values are known. However, this hypothesis can be restrictive when it comes to analyzing practical systems, since the exact values of time delays can be unknown. In this chapter, our goal is to establish new methodologies for the delay range stability analysis of a linear differential-difference system subject to the dissipative constraint in (2.10).

For an LTDS, stability information can be obtained by analyzing its corresponding spectrum. Various frequency-domain approaches [178, 268] have been developed in the frequency-domain, which can provide almost complete stability characterization when the delay systems possess certain structures. For more complex delay structures, such as DDs with general kernels, the numerical schemes in [202, 203, 215] can produce reliable results verifying a systems stability with given delay values. These results suffer almost no conservatism if numerical complexities are ignored. Furthermore, the method in [327] allows one to calculate the value of $\mathcal{H}^{\infty}$ norm of a delay system with known point-wise delay values. However, to the best of our knowledge, none of the existing spectral-based approaches can handle the problem of delay range stability analysis subject to performance objectives [514] for LTDSs. IIn other words, it is not possible to test whether a delay system is stable and simultaneously dissipative with a supply rate function [452] for all $r \in [r_1, r_2]$, where the exact delay value $r$ is unknown but bounded by $r_1 \leq r \leq r_2$ with known constants $r_2 > r_1 > 0$.

As remarked in the previous chapters, constructing LKFs [178, 452] is a standard approach in the time-domain for analyzing the stability of delay systems. Numerous different functionals (see [178, 405, 452] and the references therein) have been proposed in the literature [169, 461] to address the problem of point-wise delay stability. Compared to their frequency-domain counterparts, time-domain approaches can be more adaptable for handling the problem of range stability analysis with performance objectives, although only sufficient conditions can be constructed. In [644, 645], methods for the range stability of linear pointwise-delay systems are presented based on the principle of quadratic separation. A solution to the same problem was proposed in [469, 621] based on constructing LKFs. However, no methodologies based on the LKF approach have been developed for the range stability analysis of LDDS.[1] Furthermore, almost all LKFs used in the literature are parameterized by matrices that are independent of delay values, which can be a conservative approach when applied to range stability analysis. This motivates the use of new functionals with matrix parameters that are dependent on delay values for constructing range stability conditions subject to performance objectives.

In this chapter, we propose methodologies which allow one to conduct delay range dissipativity and stability analysis for a linear CDDS [22]; [642, 646] where the delay value is unknown but bounded. The linear CDDS model considered in this chapter contains DDs with polynomial kernels, which can represent various models of TDSs such as the cell dynamics in [641]. We employ a novel LKF with delay-dependent matrix parameters in conjunction with a quadratic supply rate function to derive our dissipativity and stability conditions. The resulting sufficient conditions, which are expressed in terms of sum-of-squares constraints [647, 648], are obtained by equivalently converting robust LMIs into SoS conditions by using the relaxation technique described in [649]. The advantage of this transformation is that it does not introduce any potential conservatism theoretically. Furthermore, we extend our proposed method to address the problem of estimating the delay margins with dissipative constraints. Finally, we also demonstrate that a hierarchy of feasibility can be established for our proposed dissipativity and stability condition, analogous to that presented in [469].

The chapter is organized as follows. In Section 6.2 we formulate the linear CDDS model to be analyzed in this chapter. Subsequently, theoretical preliminaries are presented in Section 6.3, providing the necessary tools for deriving the main results in the following section. In section 4, the main results on range stability analysis under a dissipative constraint are presented, which include remarks and expositions. Finally, we tested several numerical examples in section 5 to demonstrate the advantage of our proposed framework.

## 6.2 Problem Formulation

In this chapter, the stability of the following linear coupled differential-difference system

$$
\begin{aligned}
\dot{\boldsymbol{x}}(t) &= A_1\boldsymbol{x}(t) + A_2\boldsymbol{y}(t-r) + \int_{-r}^{0} A_3(r)L_d(\tau)\boldsymbol{y}(t+\tau)\mathsf{d}\tau + D_1\boldsymbol{w}(t), \quad \forall t \geq t_0 \\
\boldsymbol{y}(t) &= A_4\boldsymbol{x}(t) + A_5\boldsymbol{y}(t-r) \\
\boldsymbol{z}(t) &= C_1\boldsymbol{x}(t) + C_2\boldsymbol{y}(t-r) + \int_{-r}^{0} C_3(r)L_d(\tau)\boldsymbol{y}(t+\tau)\mathsf{d}\tau + D_2\boldsymbol{w}(t) \\
\boldsymbol{x}(t_0) &= \boldsymbol{\xi}, \quad \forall \theta \in [-r, 0], \quad \boldsymbol{y}(t_0+\theta) = \boldsymbol{\phi}(\theta)
\end{aligned}
\tag{6.1}
$$

with DDs is considered, where $t_0 \in \mathbb{R}$ and $\boldsymbol{\xi} \in \mathbb{R}^n$ and $\boldsymbol{\phi}(\cdot) \in \widehat{\mathcal{C}}([-r, 0); \mathbb{R}^\nu)$. Moreover, $\boldsymbol{x}(t) \in \mathbb{R}^n$ and $\boldsymbol{y}(t) \in \mathbb{R}^\nu$ satisfy the equations in (6.1), $\boldsymbol{w}(\cdot) \in \mathcal{L}^2([t_0, \infty); \mathbb{R}^q)$ represents the disturbance, $\boldsymbol{z}(t) \in \mathbb{R}^m$ is the regulated output. Note that $\boldsymbol{\xi} \in \mathbb{R}^n$ and $\boldsymbol{\phi}(\cdot) \in \widehat{\mathcal{C}}([-r, 0); \mathbb{R}^n)$ are the initial conditions for the system at $t = t_0$. $\widehat{\mathcal{C}}(\mathcal{X}; \mathbb{R}^n)$ stands for the Banach space of bounded right piecewise continuous functions with a uniform norm $\|\boldsymbol{f}(\cdot)\|_\infty := \sup_{\tau\in\mathcal{X}} \|\boldsymbol{f}(\tau)\|_2$. The dimensions of the state space matrices in (6.1) are determined by the indexes $n; \nu \in \mathbb{N}$ and $m; q \in \mathbb{N}_0 := \mathbb{N} \cup \{0\}$. Moreover, $L_d(\tau) := \boldsymbol{\ell}_d(\tau) \otimes I_\nu$ with $\boldsymbol{\ell}_d(\tau) \in \mathbb{R}^{d+1}$ contains polynomials in each row up to degree $d \in \mathbb{N}_0$. $A_3(r) \in \mathbb{R}^{n\times\varrho}$ and $C_3(r) \in \mathbb{R}^{m\times\varrho}$ are functions of $r$ which satisfies that $rA_3(r) \in \mathbb{R}^{n\times\varrho}$ and $rC_3(r) \in \mathbb{R}^{m\times\varrho}$ are polynomial matrices of $r$ with $\varrho = (d+1)\nu$. $r$ is a constant but with unknown and bounded values as $r \in [r_1, r_2]$, where the values of $r_2 > r_1 > 0$ are known. Finally, it is assumed $\rho(A_5) < 1$ which ensures the input to state stability of $\boldsymbol{y}(t) = A_4\boldsymbol{x}(t) + A_5\boldsymbol{y}(t-r)$ [22] where $\rho(A_5)$ stands for the spectral radius of $A_5$. Since $\rho(A_5) < 1$ is independent of $r$, this condition ensures the input to state stability of $\boldsymbol{y}(t) = A_4\boldsymbol{x}(t) + A_5\boldsymbol{y}(t-r)$ for all $r > 0$.

[1] The methods proposed in [568] can handle DDs with polynomials kernels. Nevertheless, the approaches in [568] are derived not based on LKFs, but the principle of robust control (Quadratic Separation).

**Remark 6.1.** Many delay-related systems can be modeled using (6.1). See [22, 169] and the references therein. In comparison with the CDDS model in [22], the system in (6.1) takes disturbances into account and contains DD with polynomial kernels at both the state and output. In terms of real-time applications, the structures of $A_3(r)$ and $C_3(r)$ can be justified by the fact that the DD kernels of a system can be dependent on the numerical values of $r$.[2]

## 6.3 Preliminaries

We present in this section some mathematical preliminaries for the derivation of the stability condition in the later section. This includes an integral inequality and the foundation of matrix polynomials optimization [647–649]. Without loss of generality, we assume in this chapter that $\boldsymbol{\ell}_d(\tau) = \mathsf{Col}_{i=0}^{d}\, \ell_i(r,\tau)$ in (6.1) consists of Legendre polynomials [467, 469, 470, 572]

$$\ell_d(r,\tau) := \sum_{k=0}^{d} \binom{d}{k}\binom{d+k}{k}\left(\frac{\tau}{r}\right)^k = \sum_{k=0}^{d} \binom{d}{k}\binom{d+k}{k}\tau^k r^{-k},\ \forall d \in \mathbb{N}_0,\ \forall \tau \in [-r,0], \tag{6.2}$$

with $\int_{-r}^{0} \boldsymbol{\ell}_d(\tau)\boldsymbol{\ell}_d^\top(\tau)\mathsf{d}\tau = \bigoplus_{k=0}^{d} \frac{r}{(2k+1)}$. Note that the form of (6.2) is derived from the structure of Jacobi polynomials (4.6) with $\alpha = \beta = 0$ and $a = -r$, $b = 0$.

Some properties of Legendre polynomials are set out as follows.

**Property 6.1.** *Given $d \in \mathbb{N}_0$ and $\boldsymbol{m}_d(\tau) := \mathsf{Col}_{i=0}^{d}\, \tau^i$, then the following three properties hold for all $r > 0$.*

- $$\exists! \mathsf{L}_d(\cdot) \in \left(\mathbb{R}_{[d+1]}^{(d+1)\times(d+1)}\right)^{\mathbb{R}_+}, \exists! \Lambda_d \in \mathbb{R}_{[d+1]}^{(d+1)\times(d+1)} : \forall \tau \in \mathbb{R},\ \ \boldsymbol{\ell}_d(\tau) = \mathsf{L}_d(r)\boldsymbol{m}_d(\tau) = \Lambda_d\left(\bigoplus_{i=0}^{d} r^i\right)^{-1}\boldsymbol{m}_d(\tau) \tag{6.3}$$
- $$\mathsf{L}_d^{-1}(r) = \left(\bigoplus_{i=0}^{d} r^i\right)\Lambda_d^{-1} \tag{6.4}$$
- $$\exists! \acute{\mathsf{L}}_d \in \mathbb{R}^{(d+1)\times(d+1)},\ \forall \tau \in \mathbb{R},\ \frac{\mathsf{d}\boldsymbol{\ell}_d(\tau)}{\mathsf{d}\tau} = r^{-1}\acute{\mathsf{L}}_d\boldsymbol{\ell}_d(\tau) \tag{6.5}$$

*where $\mathbb{R}_{[n]}^{n\times n} := \{X \in \mathbb{R}^{n\times n} : \operatorname{rank}(X) = n\}$, and $\exists!$ stands for the symbol of unique existential quantification.*

*Proof.* Since $\boldsymbol{\ell}_d(\tau)$ contains polynomials with $\int_{-r}^{0} \boldsymbol{\ell}_d(\tau)\boldsymbol{\ell}_d^\top(\tau)\mathsf{d}\tau = \bigoplus_{k=0}^{d} \frac{r}{(2k+1)}$ which is of full rank, (6.3) can be easily derived based on the form of (6.2) together with the property of positive definite matrices. By (6.3) and $r > 0$, (6.4) can be obtained. Finally, by (6.3), we have $\frac{\mathsf{d}\boldsymbol{\ell}_d(\tau)}{\mathsf{d}\tau} = \Lambda_d\left(\bigoplus_{i=0}^{d} r^i\right)^{-1}\frac{\mathsf{d}\boldsymbol{m}_d(\tau)}{\mathsf{d}\tau}$. Now it is obvious that, $\frac{\mathsf{d}\boldsymbol{m}_d(\tau)}{\mathsf{d}\tau} = \begin{bmatrix} \mathbf{0}_d^\top & 0 \\ \bigoplus_{i=1}^{d} i & \mathbf{0}_d \end{bmatrix}\boldsymbol{m}_d(\tau)$ for all $d \in \mathbb{N}_0$ if we define $\mathbf{0}_0$ and $\bigoplus_{i=1}^{0} i$ as $0\times 1$ and $0\times 0$ empty matrices, respectively. Using this relation with (6.3) and (6.4) we can obtain that

$$\frac{\mathsf{d}\boldsymbol{\ell}_d(\tau)}{\mathsf{d}\tau} = \Lambda_d\left(\bigoplus_{i=0}^{d} r^i\right)^{-1}\frac{\mathsf{d}\boldsymbol{m}_d(\tau)}{\mathsf{d}\tau} = \Lambda_d\left(\bigoplus_{i=0}^{d} r^i\right)^{-1}\begin{bmatrix} \mathbf{0}_d^\top & 0 \\ \bigoplus_{i=1}^{d} i & \mathbf{0}_d \end{bmatrix}\boldsymbol{m}_d(\tau)$$

[2] See a representative example by Example 2 in [461].

$$= \mathsf{\Lambda}_d \left(\bigoplus_{i=0}^{d} r^i\right)^{-1} \begin{bmatrix} \mathbf{0}_d^\top & 0 \\ \bigoplus_{i=1}^{d} i & \mathbf{0}_d \end{bmatrix} \left(\bigoplus_{i=0}^{d} r^i\right) \mathsf{\Lambda}_d^{-1} \boldsymbol{\ell}_d(\tau) \tag{6.6}$$

Note that based on the final term in (6.6), (6.6) can be rewritten as

$$\begin{aligned} \frac{\mathsf{d}\boldsymbol{\ell}_d(\tau)}{\mathsf{d}\tau} &= \mathsf{\Lambda}_d \begin{bmatrix} 1 & \mathbf{0}_d^\top \\ \mathbf{0}_d & \bigoplus_{i=1}^{d} r^{-i} \end{bmatrix} \begin{bmatrix} \mathbf{0}_d^\top & 0 \\ \bigoplus_{i=1}^{d} i & \mathbf{0}_d \end{bmatrix} \begin{bmatrix} \bigoplus_{i=0}^{d-1} r^i & \mathbf{0}_d \\ \mathbf{0}_d^\top & r^d \end{bmatrix} \mathsf{\Lambda}_d^{-1} \boldsymbol{\ell}_d(\tau) \\ &= \mathsf{\Lambda}_d \begin{bmatrix} \mathbf{0}_d^\top & 0 \\ r^{-1} \left(\bigoplus_{i=1}^{d} i\right) & \mathbf{0}_d \end{bmatrix} \mathsf{\Lambda}_d^{-1} \boldsymbol{\ell}_d(\tau) \end{aligned} \tag{6.7}$$

which gives (6.5). ■

**Remark 6.2.** Consider DD kernels in standard polynomials such as $\mathbb{R}^{n\times\nu} \ni \widehat{A}(\tau) = AM_d(\tau) = A\left(\boldsymbol{m}_d(\tau) \otimes I_\nu\right)$ and $\mathbb{R}^{m\times\nu} \ni \widehat{C}(\tau) = CM_d(\tau) = C\left(\boldsymbol{m}_d(\tau) \otimes I_\nu\right)$, where $\boldsymbol{m}_d(\tau) := \mathbf{Col}_{i=0}^{d}\, \tau^i$ and the matrices $A \in \mathbb{R}^{n\times\varrho}$, $C \in \mathbb{R}^{m\times\varrho}$ can be easily determined by the structure of $M_d(\tau)$. By the definition of Legendre polynomials $\boldsymbol{\ell}_d(\tau)$ in (6.2) with (6.3) and (6.4), we have $AM_d(\tau) = A\left(\mathsf{L}_d^{-1}(r) \otimes I_\nu\right) L_d(\tau)$ and $CM_d(\tau) = C\left(\mathsf{L}_d^{-1}(r) \otimes I_\nu\right) L_d(\tau)$, where $A\left(\mathsf{L}_d^{-1}(r) \otimes I_\nu\right)$ and $C\left(\mathsf{L}_d^{-1}(r) \otimes I_\nu\right)$ are polynomial matrices with respect to $r$ corresponding to $A_3(r)$ and $C_3(r)$ in (6.1). This demonstrates that the choice of Legendre polynomials $\boldsymbol{\ell}_d(\tau)$ in (6.2) together with polynomial matrices $A_3(r)$ and $C_3(r)$ in (6.1) can handle DDs with polynomial kernels.

The following inequality was initially derived in [467, 469] with different notation.

**Lemma 6.1.** *Given $U(r) \in \mathbb{S}_{\succeq 0}^{n}$ for all $r > 0$, then the inequality*

$$\int_{-r}^{0} \boldsymbol{x}^\top(\tau) U \boldsymbol{x}(\tau) \mathsf{d}\tau \geq [*] \left(r^{-1}\mathsf{D}_d \otimes U(r)\right) \left[\int_{-r}^{0} \left(\boldsymbol{\ell}_d(\tau) \otimes I_n\right) \boldsymbol{x}(\tau) \mathsf{d}\tau\right] \tag{6.8}$$

*holds for all $\boldsymbol{x}(\cdot) \in \mathscr{L}^2([-r,0];\mathbb{R}^n)$ and for all $r > 0$, where $\boldsymbol{\ell}_d(\tau)$ is defined in (6.2) and $\mathsf{D}_d := \bigoplus_{i=0}^{d} 2i+1$.*

*Proof.* Given $U(r) \in \mathbb{S}_{\succeq 0}^{n}$ for all $r > 0$. Let $\mathcal{K} = [-r, 0]$ and $\boldsymbol{f}(\tau) = \boldsymbol{\ell}_d(\tau)$ in Equation (2.11), then it gives the form of the inequality (6.8) with a known $r$ since $\int_{-r}^{0} \boldsymbol{\ell}_d(\tau)\boldsymbol{\ell}_d^\top(\tau)\mathsf{d}\tau = r\mathsf{D}_d^{-1}$. Note that the result is naturally valid for all $r > 0$ which gives this lemma. ■

**Remark 6.3.** Note that since $\boldsymbol{x}(\cdot) \in \mathscr{L}^2([-r,0];\mathbb{R}^n)$ in (6.8) and all functions in $\boldsymbol{\ell}_d(\tau)$ are bounded, all the integrals in (6.8) are well defined.

In the following definition, we define the space of univariate polynomial matrices. For the expression of multivariate polynomial matrices, see [647].

**Definition 6.1.** The space containing polynomial matrices between $\mathbb{R}^n$ to $\mathbb{R}^{p\times q}$ is defined as

$$\mathbb{R}^{p\times q}\left[\mathbb{R}\right] := \left\{ F(\cdot) \in \left(\mathbb{R}^{p\times q}\right)^{\mathbb{R}} \;\middle|\; \begin{matrix} F(x) = \sum_{i=0}^{p} Q_i x^p \;\;\&\;\; p \in \mathbb{N}_0 \\ \&\;\; Q_i \in \mathbb{R}^{p\times q} \end{matrix} \right\}. \tag{6.9}$$

Furthermore, the degree of a polynomial matrix is defined as

$$\deg\left(\sum_{i=0}^{p} Q_i x^p\right) = \max_{i=0,\ldots,p} \left(\left[\mathbb{1}_{\mathbb{R}^{p\times q}\setminus\{\mathsf{O}_{p,q}\}}(Q_i)\right] i\right) \tag{6.10}$$

where $\mathbb{1}_{\mathcal{X}}(\cdot)$ is the standard indicator function. This also allows us to define

$$\mathbb{R}^{p\times q}\left[\mathbb{R}\right]_d := \left\{F(\cdot) \in \mathbb{R}^{p\times q}\left[\mathbb{R}\right] : \deg\left(F(\cdot)\right) = d\right\}, \text{ with } \; d \in \mathbb{N}_0 \tag{6.11}$$

which contains polynomials with degree $d$.

The following definition gives the space of univariate sum-of-squares polynomial matrices. For the definition of the structure of multivariate sum-of-squares polynomial matrices, see [649] for details.

**Definition 6.2.** A polynomial in $\mathbb{S}^m[\mathbb{R}]$ is classified as a sum-of-squares polynomial if and only if it belongs to the space

$$\boldsymbol{\Sigma}\left(\mathbb{R};\mathbb{S}^m_{\succeq 0}\right) := \left\{ F(\cdot)\in\mathbb{S}^m\left[\mathbb{R}\right] \,\middle|\, \begin{array}{c} F(x)=\Phi(x)^\top\Phi(x) \\ \exists\Phi(\cdot)\in\mathbb{R}^{p\times m}\left[\mathbb{R}\right] \quad \& \quad p\in\mathbb{N} \end{array} \right\}. \tag{6.12}$$

We also define $\boldsymbol{\Sigma}_d\left(\mathbb{R};\mathbb{S}^m_{\succeq 0}\right) := \left\{F(\cdot)\in\boldsymbol{\Sigma}\left(\mathbb{R};\mathbb{S}^m_{\succeq 0}\right) : \deg\left(F(\cdot)\right)=2d\right\}$ with $d\in\mathbb{N}_0$. Finally, it is obvious to see that $\boldsymbol{\Sigma}_0\left(\mathbb{R};\mathbb{S}^m_{\succeq 0}\right)=\mathbb{S}^m_{\succeq 0}$.

The following lemma allows one to solve SoS constraints numerically via SDP. Unlike the original [649, Lemma 1], we only need to consider the univariate case.

**Lemma 6.2.** *$P(\cdot)\in\boldsymbol{\Sigma}\left(\mathbb{R};\mathbb{S}^m_{\succeq 0}\right)$ if and only if there exists $Q\in\mathbb{S}^{(d+1)m}_{\succeq 0}$ such that*

$$\forall x\in\mathbb{R},\quad P(x)=\left(\boldsymbol{m}(x)\otimes I_m\right)^\top Q\left(\boldsymbol{m}(x)\otimes I_m\right), \tag{6.13}$$

*where $\boldsymbol{m}(x) := \mathbf{Col}_{i=0}^{d}\, x^i$ with $d\in\mathbb{N}_0$.*

*Proof.* Let $u(\cdot)=\mathbf{Col}_{i=0}^{d}\, x^i$ with $d\in\mathbb{N}_0$ in [649, Lemma 1], then Lemma 6.2 with $\boldsymbol{m}(x) := \mathbf{Col}_{i=0}^{d}\, x^i$ is obtained. ■

**Remark 6.4.** In numerical calculations, we can only enforce $Q\succ 0$ instead of $Q\succeq 0$. As a result, the membership certificate produced by numerical calculations in reality is $P(\cdot)\in\boldsymbol{\Sigma}\left(\mathbb{R};\mathbb{S}^m_{\succ 0}\right)\subset\boldsymbol{\Sigma}\left(\mathbb{R};\mathbb{S}^m_{\succeq 0}\right)$.

## 6.4 Main Results of Dissipativity and Stability Analysis

In this section, the main results on range dissipativity and stability analysis are detailed. The section is divided into five subsections, and it begins by presenting the criteria for determining range delay stability and dissipativity for (6.1).

### 6.4.1 Criteria for range delay stability and dissipativity

The following range stability criteria for (6.1) can be obtained by modifying [22, Theorem 3].

**Lemma 6.3.** *Given $r_2>r_1>0$, the origin of the system (6.1) with $\boldsymbol{w}(t)\equiv\mathbf{0}_q$ is globally uniformly asymptotically (exponentially) stable for all $r\in[r_1,r_2]$, if there exist $\epsilon_1;\epsilon_2;\epsilon_3>0$ and a differentiable functional $\mathsf{v}:\mathbb{R}_+\times\mathbb{R}^n\times\widehat{\mathcal{C}}([-r,0);\mathbb{R}^\nu)\to\mathbb{R}_{\geq 0}$ such that $\forall r\in[r_1,r_2]$, $\mathsf{v}(r,\mathbf{0}_n,\mathbf{0}_\nu)=0$ and*

$$\epsilon_1\left\|\boldsymbol{\xi}\right\|_2^2\leq\mathsf{v}(r,\boldsymbol{\xi},\boldsymbol{\phi}(\cdot))\leq\epsilon_2\left(\left\|\boldsymbol{\xi}\right\|_2\vee\left\|\boldsymbol{\phi}(\cdot)\right\|_\infty\right)^2 \tag{6.14}$$

$$\left.\frac{\mathsf{d}^+}{\mathsf{d}t}\mathsf{v}(r,\boldsymbol{x}(t),\mathbf{y}_t(\cdot))\right|_{t=t_0,\boldsymbol{x}(t_0)=\boldsymbol{\xi},\mathbf{y}_{t_0}(\cdot)=\boldsymbol{\phi}(\cdot)}\leq-\epsilon_3\left\|\boldsymbol{\xi}\right\|_2^2 \tag{6.15}$$

*hold for all $r\in[r_1,r_2]$ and for any $\boldsymbol{\xi}\in\mathbb{R}^n$ and $\boldsymbol{\phi}(\cdot)\in\widehat{\mathcal{C}}([-r,0);\mathbb{R}^\nu)$ in (6.1), where $t_0\in\mathbb{R}$ and $\frac{\mathsf{d}^+}{\mathsf{d}x}f(x)=\limsup_{\eta\downarrow 0}\frac{f(x+\eta)-f(x)}{\eta}$. Furthermore, $\mathbf{y}_t(\cdot)$ in (6.15) is defined by $\forall t\geq t_0$, $\forall\theta\in[-r,0)$, $\mathbf{y}_t(\theta)=\boldsymbol{y}(t+\theta)$ where $\boldsymbol{y}(t)$ here and $\boldsymbol{x}(t)$ in (6.15) satisfying (6.1) with $\boldsymbol{w}(t)\equiv\mathbf{0}_q$.*

*Proof.* Theorem 3 of [22] is proposed for the stability analysis of a system with a given $r > 0$. However, it can be easily extended point-wisely by treating $r$ in the system as an uncertain parameter belonging to an interval $[r_1, r_2]$ with $r_2 > r_1 > 0$. Moreover, the functions $\mathsf{v}(\cdot)$,$u(\cdot)$,$\mathsf{v}(\cdot)$ and $w(\cdot)$ in [22, Theorem 3] should be parameterized by $r$ in this case. Thus a corresponding range stability criteria can be obtained which can be applied to (6.1) with $\boldsymbol{w}(t) \equiv \mathbf{0}_q$. Following the aforementioned steps and letting the functions $u(r, \cdot), v(r, \cdot), w(r, \cdot)$ be the quadratic functions $\epsilon_i x^2$, $i = 1, 2, 3$, Lemma 6.3 can be obtained given the fact that (6.1) is a special case of the general system considered in [22, Theorem 3]. ■

**Definition 6.3** (Dissipativity)**.** Given $r_2 > r_1 > 0$, the system (6.1) with a supply rate function $\mathsf{s}(\boldsymbol{z}(t), \boldsymbol{w}(t))$ is dissipative for all $r \in [r_1, r_2]$, if there exists a differentiable functional $\mathsf{v} : \mathbb{R}_+ \times \mathbb{R}^n \times \widehat{\mathscr{C}}([-r, 0); \mathbb{R}^\nu) \to \mathbb{R}$ such that

$$\forall r \in [r_1, r_2],\ \forall t \geq t_0 :\ \dot{\mathsf{v}}(r, \boldsymbol{x}(t), \mathbf{y}_t(\cdot)) - \mathsf{s}(\boldsymbol{z}(t), \boldsymbol{w}(t)) \leq 0 \tag{6.16}$$

with $t_0 \in \mathbb{R}$, where $\mathbf{y}_t(\cdot)$ is defined by the equality $\forall t \geq t_0$, $\forall \theta \in [-r, 0)$, $\mathbf{y}_t(\theta) = \boldsymbol{y}(t + \theta)$, and $\boldsymbol{x}(t)$, $\boldsymbol{y}(t)$ and $\boldsymbol{z}(t)$ satisfy the equalities in (6.1) with $\boldsymbol{w}(\cdot) \in \mathscr{L}^2([t_0, \infty); \mathbb{R}^q)$.

To incorporate dissipitivity into the analysis of (6.1), we utilize the same quadratic supply rate function

$$\begin{gathered}
\mathsf{s}(\boldsymbol{z}(t), \boldsymbol{w}(t)) = \begin{bmatrix} \boldsymbol{z}(t) \\ \boldsymbol{w}(t) \end{bmatrix}^\top \mathbf{J} \begin{bmatrix} \boldsymbol{z}(t) \\ \boldsymbol{w}(t) \end{bmatrix} \quad \text{with} \quad \mathbf{J} = \begin{bmatrix} \widetilde{J}^\top J_1^{-1} \widetilde{J} & J_2 \\ * & J_3 \end{bmatrix} \in \mathbb{S}^{(m+q)} \\
\widetilde{J}^\top J_1^{-1} \widetilde{J} \preceq 0, \quad J_1^{-1} \prec 0
\end{gathered} \tag{6.17}$$

as in (2.10).

### 6.4.2 Conditions for range dissipativity and stability analysis

In this subsection, the main results on dissipativity and stability analysis are derived in the following theorem where the optimization constraints can be solved via the method of sum-of-squares programming.

**Theorem 6.1.** *Given $\lambda_1; \lambda_2; \lambda_3 \in \mathbb{N}_0$ and $\boldsymbol{\ell}_d(\tau)$ consisting of the Legendre polynomials in (6.2) with $d \in \mathbb{N}_0$, the system (6.1) with supply rate function (6.17) is dissipative for all $r \in [r_1, r_2]$, and the origin of (6.1) with $\boldsymbol{w}(t) \equiv \mathbf{0}_q$ is globally uniformly asymptotically stable for all $r \in [r_1, r_2]$, if there exist matrix polynomials*

$$P(\cdot) \in \mathbb{S}^{n+\varrho}[\mathbb{R}]_{\lambda_1}, \quad S(\cdot) \in \mathbb{S}^\nu[\mathbb{R}]_{\lambda_2} \quad U(\cdot) \in \mathbb{S}^\nu[\mathbb{R}]_{\lambda_3}, \quad \widehat{P}(\cdot) \in \mathbb{S}^{n+\varrho}[\mathbb{R}] \quad \widehat{S}(\cdot); \widehat{U}(\cdot) \in \mathbb{S}^\nu[\mathbb{R}]$$

*and $\delta_i \in \mathbb{N}_0$, $i = 1, \ldots, 8$ with $\delta_7 \neq 0$ such that*

$$P(\cdot) + \Big[\mathsf{O}_n \oplus (\mathsf{D}_d \otimes S(\cdot))\Big] + g(\cdot)\widehat{P}(\cdot) \in \boldsymbol{\Sigma}_{\delta_1}\Big(\mathbb{R}; \mathbb{S}^{n+\varrho}_{\succ 0}\Big) \ \ \widehat{P}(\cdot) \in \boldsymbol{\Sigma}_{\delta_2}\Big(\mathbb{R}; \mathbb{S}^{n+\varrho}_{\succeq 0}\Big) \tag{6.18}$$

$$S(\cdot) + g(\cdot)\widehat{S}(\cdot) \in \boldsymbol{\Sigma}_{\delta_3}\left(\mathbb{R}; \mathbb{S}^n_{\succeq 0}\right), \ \widehat{S}(\cdot) \in \boldsymbol{\Sigma}_{\delta_4}\left(\mathbb{R}; \mathbb{S}^n_{\succeq 0}\right) \tag{6.19}$$

$$U(\cdot) + g(\cdot)\widehat{U}(\cdot) \in \boldsymbol{\Sigma}_{\delta_5}\left(\mathbb{R}; \mathbb{S}^n_{\succeq 0}\right), \ \widehat{U}(\cdot) \in \boldsymbol{\Sigma}_{\delta_6}\left(\mathbb{R}; \mathbb{S}^n_{\succeq 0}\right) \tag{6.20}$$

$$-\begin{bmatrix} J_1 & \widetilde{J}\Sigma(\cdot) \\ * & \mathbf{\Phi}_d(\cdot) \end{bmatrix} + g(\cdot)\mathbf{\Psi}(\cdot) \in \boldsymbol{\Sigma}_{\delta_7}\Big(\mathbb{R}; \mathbb{S}^{m+q+2n+\varrho}_{\succ 0}\Big), \quad \mathbf{\Psi}(\cdot) \in \boldsymbol{\Sigma}_{\delta_8}\Big(\mathbb{R}; \mathbb{S}^{m+q+2n+\varrho}_{\succeq 0}\Big) \tag{6.21}$$

*where $\varrho = (d + 1)\nu$ and $g(r) = (r - r_1)(r - r_2)$ and*

$$
\boldsymbol{\Phi}_d(r) := \mathsf{Sy}\left(\begin{bmatrix} \mathsf{O}_{q,n} & \mathsf{O}_{q,\varrho} \\ I_n & \mathsf{O}_{n,\varrho} \\ \mathsf{O}_n & \mathsf{O}_{n,\varrho} \\ \mathsf{O}_{\varrho,n} & rI_\varrho \end{bmatrix} P(r) \begin{bmatrix} D_1 & A_1 & A_2 & rA_3(r) \\ \mathsf{O}_{\varrho,q} & L_d(0)A_4 & L_d(0)A_5 - L_d(-r) & -\widehat{\mathsf{L}}_d \end{bmatrix}\right)
$$

$$
+\Gamma^\top (rS(r) + rU(r))\,\Gamma - \Big[J_3 \oplus \mathsf{O}_n \oplus rS(r) \oplus (r\mathsf{D}_d \otimes U(r))\Big] - \mathsf{Sy}\left(\Big[\Sigma^\top J_2 \quad \mathsf{O}_{(n+\nu+\varrho+q),(n+\nu+\varrho)}\Big]\right), \tag{6.22}
$$

*with*

$$
\Gamma := \Big[\mathsf{O}_{\nu,q} \quad A_4 \quad A_5 \quad \mathsf{O}_{\nu,\varrho}\Big], \quad \Sigma(r) := \Big[D_2 \quad C_1 \quad C_2 \quad rC_3(r)\Big]. \tag{6.23}
$$

*and* $\widehat{\mathsf{L}}_d := \acute{\mathsf{L}}_d \otimes I_\nu$ *in which* $\acute{\mathsf{L}}_d$ *is given in* (6.5).

*Proof.* The proof of this theorem is based on the construction of the parameterized functional

$$
\begin{aligned}
\mathsf{v}(r, \boldsymbol{x}(t), \mathbf{y}_t(\cdot)) := \begin{bmatrix} \boldsymbol{x}(t) \\ \int_{-r}^{0} L_d(\tau)\boldsymbol{y}(t+\tau)\mathsf{d}\tau \end{bmatrix}^\top P(r) \begin{bmatrix} \boldsymbol{x}(t) \\ \int_{-r}^{0} L_d(\tau)\boldsymbol{y}(t+\tau)\mathsf{d}\tau \end{bmatrix} \\
+ \int_{-r}^{0} \boldsymbol{y}^\top(t+\tau)\Big[rS(r) + (\tau + r)U(r)\Big]\boldsymbol{y}(t+\tau)\mathsf{d}\tau,
\end{aligned} \tag{6.24}
$$

where $\mathbf{y}_t(\cdot)$ follows the definition in (6.16) and the functional satisfies $\mathsf{v}(r, \mathbf{0}_n, \mathbf{0}_\nu) = 0$ for all $r \in [r_1, r_2]$ with given $r_2 > r_1 > 0$. Furthermore, $L_d(\tau)$ in (6.24) is defined as $L_d(\tau) = \boldsymbol{\ell}_d(\tau) \otimes I_\nu$ with $\boldsymbol{\ell}_d(\tau)$ in (6.2), and the matrix parameters in (6.24) are $P(\cdot) \in \mathbb{S}^{n+\varrho}[\mathbb{R}]_{\lambda_1}$, $S(\cdot) \in \mathbb{S}^\nu[\mathbb{R}]_{\lambda_2}$ and $U(\cdot) \in \mathbb{S}^\nu[\mathbb{R}]_{\lambda_3}$ with degree indexes $\lambda_1; \lambda_2; \lambda_3 \in \mathbb{N}_0$.[3]

*First of all, we demonstrate that the feasible solutions to* (6.20)–(6.21) *guarantee the existence of* (6.24) *satisfying* (6.16) *and* (6.15). Differentiating $\mathsf{v}(r, \boldsymbol{x}(t), \mathbf{y}_t(\cdot))$ along the trajectory of (6.1) and considering the relation

$$
\begin{aligned}
\int_{-r}^{0} L_d(\tau)\dot{\boldsymbol{y}}(t+\tau)\mathsf{d}\tau &= L_d(0)\boldsymbol{y}(t) - L_d(-r)\boldsymbol{y}(t-r) - \widehat{\mathsf{L}}_d\frac{1}{r}\int_{-r}^{0} L_d(\tau)\boldsymbol{y}(t+\tau)\mathsf{d}\tau \\
&= L_d(0)A_4\boldsymbol{x}(t) + \left(L_d(0)A_5 - L_d(-r)\right)\boldsymbol{y}(t-r) - \widehat{\mathsf{L}}_d\frac{1}{r}\int_{-r}^{0} L_d(\tau)\boldsymbol{y}(t+\tau)\mathsf{d}\tau
\end{aligned} \tag{6.25}
$$

produces

$$
\begin{aligned}
&\forall t \geq t_0, \quad \dot{\mathsf{v}}(r, \boldsymbol{x}(t), \mathbf{y}_t(\cdot)) - \mathsf{s}(\boldsymbol{z}(t), \boldsymbol{w}(t)) \\
&= \boldsymbol{\chi}_d^\top(t)\,\mathsf{Sy}\left(\begin{bmatrix} \mathsf{O}_{q,n} & \mathsf{O}_{q,\varrho} \\ I_n & \mathsf{O}_{n,\varrho} \\ \mathsf{O}_{\nu,n} & \mathsf{O}_{\nu,\varrho} \\ \mathsf{O}_{\varrho,n} & rI_\varrho \end{bmatrix} P(r) \begin{bmatrix} D_1 & A_1 & A_2 & rA_3(r) \\ \mathsf{O}_{\varrho,q} & L_d(0)A_4 & L_d(0)A_5 - L_d(-r) & -\widehat{\mathsf{L}}_d \end{bmatrix}\right)\boldsymbol{\chi}_d(t) \\
&+ \boldsymbol{\chi}_d^\top(t)\Big[\Gamma^\top (rS(r) + rU(r))\,\Gamma - \Big(J_3 \oplus \mathsf{O}_n \oplus rS(r) \oplus \mathsf{O}_\varrho\Big)\Big]\boldsymbol{\chi}_d(t) \\
&- \boldsymbol{\chi}_d^\top(t)\Big(\Sigma^\top(r)\widetilde{J}^\top J_1^{-1}\widetilde{J}\Sigma(r) + \mathsf{Sy}\left(\Big[\Sigma^\top(r)J_2 \quad \mathsf{O}_{(n+\nu+\varrho+q),(n+\nu+\varrho)}\Big]\right)\Big]\boldsymbol{\chi}_d(t), \\
&- \int_{-r}^{0} \boldsymbol{y}^\top(t+\tau)U(r)\boldsymbol{y}(t+\tau)\mathsf{d}\tau,
\end{aligned} \tag{6.26}
$$

where

$$
\boldsymbol{\chi}_d(t) := \mathbf{Col}\left(\boldsymbol{w}(t),\ \boldsymbol{x}(t),\ \boldsymbol{y}(t-r),\ \frac{1}{r}\int_{-r}^{0} L_d(\tau)\boldsymbol{y}(t+\tau)\mathsf{d}\tau\right) \tag{6.27}
$$

[3] Note that $\mathbb{S}^\nu[\mathbb{R}]_0 = \mathbb{S}^\nu$

and $\Gamma$, $\Sigma(r)$ have been defined in (6.23), and $\widehat{\mathsf{L}}_d := \acute{\mathsf{L}}_d \otimes I_\nu$ in (6.25) can be obtained by (6.5) with (2.1a). Assume $U(r) \succeq 0$, $\forall r \in [r_1, r_2]$. Seeing that $\mathbf{y}_t(\cdot) \in \widehat{\mathcal{C}}([-r,0);\mathbb{R}^\nu) \subset \mathcal{L}^2([-r,0);\mathbb{R}^\nu)$, now apply (6.8) to the integral $\int_{-r}^{0} \boldsymbol{y}^\top(t+\tau)U(r)\boldsymbol{y}(t+\tau)\mathsf{d}\tau$ in (6.26). It produces

$$\forall r \in [r_1, r_2], \quad \int_{-r}^{0} \boldsymbol{y}^\top(t+\tau)U(r)\boldsymbol{y}(t+\tau)\mathsf{d}\tau \geq [*]\left(r\mathsf{D}_d \otimes U(r)\right)\left[\int_{-r}^{0} r^{-1}L_d(\tau)\boldsymbol{y}(t+\tau)\mathsf{d}\tau\right] \tag{6.28}$$

with $\mathsf{D}_d = \bigoplus_{i=0}^{d} 2i+1$. Moreover, applying (6.28) to (6.26) yields

$$\forall r \in [r_1, r_2],\ \forall t \geq t_0,\ \dot{\mathsf{v}}(r, \boldsymbol{x}(t), \mathbf{y}_t(\cdot)) - \mathsf{s}(\boldsymbol{z}(t), \boldsymbol{w}(t)) \leq \boldsymbol{\chi}_d^\top(t)\left(\boldsymbol{\Phi}_d(r) - \Sigma^\top(r)\widetilde{J}^\top J_1^{-1}\widetilde{J}\Sigma(r)\right)\boldsymbol{\chi}_d(t) \tag{6.29}$$

where $\boldsymbol{\Phi}_d(r)$ and $\boldsymbol{\chi}_d(t)$ have been defined in (6.22) and (6.27), respectively. Based on the structure of (6.29), it is evident that if

$$\forall r \in [r_1, r_2]: \quad \boldsymbol{\Phi}_d(r) - \Sigma^\top(r)\widetilde{J}^\top J_1^{-1}\widetilde{J}\Sigma(r) \prec 0, \quad U(r) \succeq 0 \tag{6.30}$$

is satisfied then the dissipative inequality in (6.16) : $\dot{\mathsf{v}}(r, \boldsymbol{x}(t), \mathbf{y}_t(\cdot)) - \mathsf{s}(\boldsymbol{z}(t), \boldsymbol{w}(t)) \leq 0$ holds $\forall r \in [r_1, r_2]$ and $\forall t \geq t_0$.

Furthermore, by considering $J_1^{-1} \prec 0$ and the structure of $\boldsymbol{\Phi}_d(r) - \Sigma^\top(r)\widetilde{J}^\top J_1^{-1}\widetilde{J}\Sigma(r) \prec 0$, $\forall r \in [r_1, r_2]$ with the properties of negative definite matrices, it is obvious that given (6.30) holds then there exists (6.24) and $\epsilon_3 > 0$ satisfying $\forall r \in [r_1, r_2]$, $\forall t \geq t_0$, $\dot{\mathsf{v}}(r, \boldsymbol{x}(t), \mathbf{y}_t(\cdot)) \leq -\epsilon_3 \left\|\boldsymbol{x}(t)\right\|_2^2$ along the trajectory of (6.1) with $\boldsymbol{w}(t) \equiv \mathbf{0}_q$. Now consider the case of $t = t_0$ for the previous inequality with the initial conditions in (6.1), it shows that the feasible solutions to (6.30) imply the existence of $\epsilon_3$ and (6.24) satisfying (6.15). On the other hand, given $J_1^{-1} \prec 0$, applying the Schur complement to (6.30) enables one to conclude that (6.30) holds if and only if

$$\forall r \in \mathcal{G}: \quad \boldsymbol{\Theta}_d(r) = \begin{bmatrix} J_1 & \widetilde{J}\Sigma(r) \\ * & \boldsymbol{\Phi}_d(r) \end{bmatrix} \prec 0, \quad U(r) \succeq 0 \tag{6.31}$$

with $\mathcal{G} := \left\{\rho \in \mathbb{R} : g(\rho) := (\rho - r_1)(\rho - r_2) \leq 0\right\} = [r_1, r_2]$. By applying the matrix sum-of-squares relaxation technique proposed in [649] to (6.31), given that $g(\cdot)$ naturally satisfies the qualification constraint in [649, Theorem 1]. Then we can conclude that (6.31) holds if and only if[4] (6.20) and (6.21) hold for some $\delta_i$, $i = 5, \ldots, 8$. This shows that the feasible solutions to (6.20)–(6.21) guarantee the existence of $\epsilon_3 > 0$ and (6.24) satisfying (6.16) and (6.15).

*Now we start to prove that* (6.18)–(6.20) *imply that* (6.24) *satisfies* (6.14) *with* $\epsilon_1 > 0$ *and* $\epsilon_2 > 0$. Given the structure of (6.24) and consider the situation of $t = t_0$ with the initial conditions in (6.1), it follows that there exists $\lambda > 0$ such that for all $r \in [r_1, r_2]$

$$\begin{aligned}
\mathsf{v}(r, \boldsymbol{\xi}, \boldsymbol{\phi}(\cdot)) &\leq \begin{bmatrix} \boldsymbol{\xi} \\ \int_{-r}^{0} L_d(\tau)\boldsymbol{\phi}(\tau)\mathsf{d}\tau \end{bmatrix}^\top \lambda \begin{bmatrix} \boldsymbol{\xi} \\ \int_{-r}^{0} L_d(\tau)\boldsymbol{\phi}(\tau)\mathsf{d}\tau \end{bmatrix} + \int_{-r}^{0} \boldsymbol{\phi}^\top(\tau)\lambda\boldsymbol{\phi}(\tau)\mathsf{d}\tau \\
&\leq \lambda \left\|\boldsymbol{\xi}\right\|_2^2 + \int_{-r}^{0} \boldsymbol{\phi}^\top(\tau)L_d^\top(\tau)\mathsf{d}\tau\lambda\int_{-r}^{0} L_d(\tau)\boldsymbol{\phi}(\cdot)\mathsf{d}\tau + \lambda r\left\|\boldsymbol{\phi}(\tau)\right\|_\infty^2 \\
&\leq \lambda \left\|\boldsymbol{\xi}\right\|_2^2 + \lambda r\left\|\boldsymbol{\phi}(\cdot)\right\|_\infty^2 + [*]\left(\lambda\mathsf{D}_d \otimes I_n\right)\int_{-r}^{0} L_d(\tau)\boldsymbol{\phi}(\tau)\mathsf{d}\tau \\
&\leq \lambda \left\|\boldsymbol{\xi}\right\|_2^2 + \lambda r\left\|\boldsymbol{\phi}(\cdot)\right\|_\infty^2 + r\int_{-r}^{0} \boldsymbol{\phi}^\top(\tau)\lambda\boldsymbol{\phi}(\tau)\mathsf{d}\tau \\
&\leq \lambda \left\|\boldsymbol{\xi}\right\|_2^2 + \left(\lambda r + \lambda r^2\right)\left\|\boldsymbol{\phi}(\cdot)\right\|_\infty^2
\end{aligned} \tag{6.32}$$

[4] See the results related to the equations (1) and (6) in [649]

$$\leq \left(\lambda + \lambda r^2\right) \|\boldsymbol{\xi}\|_2^2 + \left(\lambda r + \lambda r^2\right) \|\boldsymbol{\phi}(\cdot)\|_\infty^2$$

$$\leq 2\left(\lambda r + \lambda r^2\right) \left(\|\boldsymbol{\xi}\|_2 \vee \|\boldsymbol{\phi}(\cdot)\|_\infty\right)^2 \tag{6.33}$$

$$\leq 2\left(\lambda r_2 + \lambda r_2^2\right) \left(\|\boldsymbol{\xi}\|_2 \vee \|\boldsymbol{\phi}(\cdot)\|_\infty\right)^2 \tag{6.34}$$

holds for any $\boldsymbol{\xi} \in \mathbb{R}^n$ and $\boldsymbol{\phi}(\cdot) \in \widehat{\mathcal{C}}([-r,0);\mathbb{R}^\nu)$ in (6.1). As a consequence, it shows that there exist $\epsilon_2 > 0$ such that (6.24) satisfies the upper bound property in (6.14).

Now assume $S(r) \succeq 0, \forall r \in [r_1, r_2]$. Then applying (6.8) to the integral $\int_{-r}^{0} \boldsymbol{y}^\top(t+\tau)S(r)\boldsymbol{y}(t+\tau)\mathsf{d}\tau$ in (6.24) at $t = t_0$ yields

$$\forall r \in [r_1, r_2],\ \ r\int_{-r}^{0} \boldsymbol{\phi}^\top(\tau)S(r)\boldsymbol{\phi}(\tau)\mathsf{d}\tau \geq \int_{-r}^{0} \boldsymbol{\phi}^\top(\tau)L_d^\top(\tau)\mathsf{d}\tau \left(\mathsf{D}_d \otimes S(r)\right) \int_{-r}^{0} L_d(\tau)\boldsymbol{\phi}(\tau)\mathsf{d}\tau, \tag{6.35}$$

with $\mathsf{D}_d = \bigoplus_{i=0}^{d} 2i+1$. Applying (6.35) to (6.24) at $t = t_0$ produces that for all $r \in [r_1, r_2]$

$$\mathsf{v}(r,\boldsymbol{\xi},\boldsymbol{\phi}(\cdot)) \geq \begin{bmatrix} \boldsymbol{\xi} \\ \int_{-r}^{0} L_d(\tau)\boldsymbol{\phi}(\tau)\mathsf{d}\tau \end{bmatrix}^\top \left(P(r) + \left[\mathsf{O}_n \oplus \left(\mathsf{D}_d \otimes S(r)\right)\right]\right) \begin{bmatrix} \boldsymbol{\xi} \\ \int_{-r}^{0} L_d(\tau)\boldsymbol{\phi}(\tau)\mathsf{d}\tau \end{bmatrix} + \int_{-r}^{0} (\tau + r)\boldsymbol{\phi}^\top(\tau)U(r)\boldsymbol{\phi}(\tau)\mathsf{d}\tau. \tag{6.36}$$

holds for any $\boldsymbol{\xi} \in \mathbb{R}^n$ and $\boldsymbol{\phi}(\cdot) \in \widehat{\mathcal{C}}([-r,0);\mathbb{R}^\nu)$ in (6.1). By (6.36), we see that if

$$\forall r \in \mathcal{G}:\ \ \Pi_d(r) := P(r) + \left[\mathsf{O}_n \oplus \left(\mathsf{D}_d \otimes S(r)\right)\right] \succ 0,\quad S(r) \succeq 0,\ \ U(r) \succeq 0 \tag{6.37}$$

is satisfied, then there exists $\epsilon_1 > 0$ such that for all $r \in [r_1, r_2]$ we have $\epsilon_1 \|\boldsymbol{\xi}\|_2 \leq \mathsf{v}(r,\boldsymbol{\xi},\boldsymbol{\phi}(\cdot))$ for any $\boldsymbol{\xi} \in \mathbb{R}^n$ and $\boldsymbol{\phi}(\cdot) \in \widehat{\mathcal{C}}([-r,0);\mathbb{R}^\nu)$ in (6.1), where $\mathcal{G} := \left\{\rho \in \mathbb{R} : g(\rho) := (\rho - r_1)(\rho - r_2) \leq 0\right\} = [r_1, r_2]$. Now apply the matrix relaxation technique in [649] to the conditions in (6.37). Then one can conclude that (6.37) holds if and only if (6.18)–(6.20) hold for some $\delta_i$, $i = 1, \ldots, 6$. Considering the upper bound result derived in (6.34), one can see that the feasible solutions to (6.18)–(6.20) guarantee the existence of (6.24) and $\epsilon_1; \epsilon_2 > 0$ satisfying (6.14).

In conclusion, we have demonstrated that the feasible solutions to (6.18)–(6.21) imply the existence of $\epsilon_1; \epsilon_2; \epsilon_3 > 0$ and (6.24) satisfying(6.14)–(6.16). This finishes the proof. ∎

**Remark 6.5.** Note that the structure of (6.24) is inspired by the complete LKF proposed in [22]. Given that all matrix parameters in (6.24) are related to $r$ polynomially, it can be reasonably expected that less conservative results can be obtained for range delay stability analysis using (6.24) as compared to constructing an LKF with only constant matrix parameters.

**Remark 6.6.** All the SoS constraints in Theorem 6.1 can be solved numerically using the relation specified in (6.13). The dimension of the corresponding certificate variable $Q$ in (6.13) is specified by the values of $\delta_i$, $i = 1, \ldots, 8$ with $\lambda_1$, $\lambda_2$ and $\lambda_3$ in (6.24).

**Remark 6.7.** Various forms of $g(r)$ can be utilized to describe the set $\mathcal{G} = [r_1, r_2] = \{r \in \mathbb{R} : g(r) \leq 0\}$, provided that $g(r) \leq 0$ can equivalently characterize the interval $[r_1, r_2]$ and satisfy the qualification constraints in [649]. This implies that the feasibility of the corresponding sum-of-squares (SoS) constraints is not affected by the use of a valid $g(r)$ with different form, as they are ultimately equivalent to (6.31) and (6.37). Nevertheless, considering its low degree and the consequent reduction in computational burden to solve (6.31) and (6.37), the form $g(r) = (r - r_1)(r - r_2)$ the optimal choice for solving (6.18)–(6.21).

**Remark 6.8.** Pointwise delay stability analysis at $r = r_0 > 0$ can be evaluated by solving

$$S \succeq 0,\ U \succeq 0,\ \ \Pi \succ 0,\ \ \mathbf{\Theta}_d(r_0) \prec 0, \tag{6.38}$$

where the value of $r_0$ is given with $\lambda_1 = \lambda_2 = \lambda_3 = 0$ in (6.24). Since the value of $r_0$ is fixed, there is no need to consider non-constant polynomial matrix variables for (6.24). It should be emphasized that whenever (6.38) is referenced in this chapter, it is assumed that $\lambda_1 = \lambda_2 = \lambda_3 = 0$ in (6.24).

### 6.4.3 Reducing the computational burden of Theorem 6.1 for certain cases

The SoS constraints in (6.18)–(6.21) can be applied to handle any form of (6.1) with $\lambda_1, \lambda_2, \lambda_3$ of known values in (6.24), provided that appropriate selections of $\delta_i$, $i = 1, \ldots, 8$ are made. However, if any inequality in (6.31) or (6.37) is affine in $r$, it can be solved using the property of the convex hull, resulting in a significant reduction in numerical complexity compared to solving the equivalent SoS constraints in equations (6.18)–(6.21). Notwithstanding, this approach is only applicable to very specific cases, which are elaborated on in the following discussion.

**Case 1.** Let $\lambda_1 = \lambda_2 = \lambda_3 = 0$ in (6.24) and matrix-valued functions $rA_3(r), rC_3(r)$ that are affine in $r$, then $\mathbf{\Theta}_d(r) \prec 0$ in (6.39) becomes affine in $r$ and the matrix inequalities in (6.18)–(6.20) become standard LMIs

$$S \succeq 0,\ U \succeq 0,\ \ \Pi \succ 0, \tag{6.39}$$

with $\delta_i = 0$, $i = 1, \ldots, 6$ and $\widehat{P}(r) = \mathsf{O}_{\nu+\varrho}$, $\widehat{S}(r) = \widehat{U}(r) = \mathsf{O}_{\nu}$, where (6.39) can be obtained directly from (6.37) with $\lambda_1 = \lambda_2 = \lambda_3 = 0$ also. Since $\mathbf{\Theta}_d(r) \prec 0$ in (6.39) is affine in $r$, hence $\forall r \in [r_1, r_2]$, $\mathbf{\Theta}_d(r) \prec 0$ can be solved by convex hull properties instead of using (6.21). Meanwhile, $\forall r \in [r_1, r_2]$, $\mathbf{\Theta}_d(r) \prec 0$ can still be solved via the SoS condition (6.21), with degrees $\delta_7 = 1$ and $\delta_8 = 0$. However, since using SoS does not provide any additional feasibility, it is preferable in this case to solve $\forall r \in [r_1, r_2]$, $\mathbf{\Theta}_d(r) \prec 0$ using the property of the convex hull instead of SoS to reduce the number of decision variables.

**Case 2.** If any inequality in (6.37) is affine or convex[5] in $r$, then it can be directly solved using the property of the convex hull. Meanwhile, if unstructured matrix variables are considered in (6.24) without predefined sparsities, the only scenario where $\mathbf{\Theta}_d(r) \prec 0$ in (6.31) is an affine (convex) matrix inequality in $r$ is when $\lambda_1 = \lambda_2 = \lambda_3 = 0$ in (6.24) with $rA_3(r)$ being a constant[6] and $rC_3(r)$ being affine in $r$, based on the structure of $\mathbf{\Theta}_d(r) \prec 0$. Consequently, the convex hull's properties cannot be applied to solve $\forall r \in [r_1, r_2]$, $\mathbf{\Theta}_d(r) \prec 0$ if $rA_3(r)$ is non-constant and $rC_3(r)$ is not affine in $r$.

We have demonstrated that for certain situations, the parameter dependent LMIs in (6.37) and (6.31) can be solved using the convex hull's properties requiring fewer decision variables instead of solving the equivalent SoS constraints in (6.18)–(6.21). However, as illustrated in Case 1 and Case 2, the SoS constraint in (6.21) still necessitates the use of the proposed method if $rA_3(r)$ is non-constant and $rC_3(r)$ is not affine in $r$, which still holds true even if (6.24) is parameterized only by constant matrix parameters ($\lambda_1 = \lambda_2 = \lambda_3 = 0$).

[5] This may include the situation such as $S(r) = S_1 + r^4 S_2$. However, the handling of $\forall r \in [r_1, r_2]$, $S(r) \succ 0$ is the same as for an affine case. In addition, the variable structure such as $S(r) = r^3 S_1$ is not considered in this chapter, as it can always be equivalently transformed into a constant parameter.

[6] See also in Remark 5 of [621] for a relevant exposition of range stability without considering output and disturbance

### 6.4.4 Estimating delay margins subject to prescribed performance objectives

Given an initial $r_0$ along with a supply rate function (6.17) (assume no decision variables in (6.17)) which renders (6.38) to be feasible, we are interested in the following problem.

**Problem 1.** Finding the minimum $\grave{r}$ or maximum $\acute{r}$ which render (6.1) stable and dissipative with (6.17) over $[\grave{r}, r_0]$ or $[r_0, \acute{r}]$, where the matrices in (6.17) contain no decision variables.

The control interpretation of this problem is straightforward: Given a specific performance objective, compute the maximum stable delay interval of a delay system such that the system can consistently satisfy the given performance objective within this interval.

Problem 1 can be solved by the following optimization programs

$$\min \rho \quad \text{subject to} \quad (6.18)-(6.21) \quad \text{with} \quad g(r) = (r-\rho)(r-r_0) \tag{6.40}$$

or

$$\max \rho \quad \text{subject to} \quad (6.18)-(6.21) \quad \text{with} \quad g(r) = (r-r_0)(r-\rho), \tag{6.41}$$

with given $\lambda_1$, $\lambda_2$, $\lambda_3$ and $\delta_i$, $i = 1, \ldots, 8$. Specifically, the optimization problems in (6.40) and (6.41) can be efficiently handled using an iterative one-dimensional search scheme in conjunction with the bisection method [507]. Since both (6.40) and (6.41) represent delay range dissipativity and stability conditions, the use of the bisection method does not produce false feasible solutions even if (6.40) and (6.41) are not necessarily quasi-convex. Furthermore, as explicated in Section 6.4.3, if any inequality in (6.37) and (6.31) is affine (convex) in $r$, then it can be directly solved using the property of the convex hull to replace the SoS conditions in (6.40) and (6.41). Finally, it is vitally important to emphasize here that the result of dissipativity over $[\grave{r}, r_0]$ or $[r_0, \acute{r}]$, which is individually produced by (6.40) and (6.41), cannot be automatically merged due to the mathematical nature of range dissipativity analysis.

### 6.4.5 A hierarchy of the conditions in Theorem 6.1

Here we show that the feasibility of the range dissipativity and stability condition in Theorem 6.1 exhibits a hierarchy with respect to $d$.

**Theorem 6.2.** *Given Legendre polynomials $\boldsymbol{\ell}_d(\tau)$ in (6.2), we have*

$$\forall d \in \mathbb{N}_0, \quad \mathcal{F}_d \subseteq \mathcal{F}_{d+1} \tag{6.42}$$

*where*

$$\begin{aligned} \mathcal{F}_d &:= \left\{ (r_1, r_2) \,\middle|\, r_2 > r_1 > 0 \quad \& \quad (6.37) \text{ and } (6.31) \text{ hold} \right\} \\ &= \left\{ (r_1, r_2) \,\middle|\, r_2 > r_1 > 0 \quad \& \quad (6.18)\text{–}(6.21) \text{ hold} \quad \& \quad \delta_7 \in \mathbb{N} \quad \& \quad \delta_i; \delta_8 \in \mathbb{N}_0, i = 1, \ldots, 6 \right\}. \end{aligned} \tag{6.43}$$

*Proof.* Let $d \in \mathbb{N}_0$ and $(r_1, r_2) \in \mathcal{F}_d$ with $\mathcal{F}_d \neq \varnothing$ and suppose (6.37) is satisfied by $S(r), U(r)$ and $P_d(r)$ at $d$. By assuming $P_{d+1}(r) = P_d(r) \oplus \mathsf{O}_n$ for $P(r)$ in (6.24) at $d+1$ and considering the structure of (6.37), we have

$$\begin{gathered} \forall r \in [r_1, r_2], \quad \Pi_{d+1}(r) = P_{d+1}(r) + \left[ \mathsf{O}_n \oplus (\mathsf{D}_{d+1} \otimes S(r)) \right] = \Pi_d(r) \oplus [(2d+3) \otimes S(r)] \succ 0, \\ \forall r \in [r_1, r_2], \quad S(r) \succeq 0, \quad U(r) \succeq 0. \end{gathered} \tag{6.44}$$

Since $S(r) \succeq 0$ and $2d+3>0$, it can be concluded from (6.44) that the existence of feasible solutions to $\forall r \in [r_1, r_2]$, $\Pi_d(r) \succ 0$ implies the existence of feasible solutions to $\forall r \in [r_1, r_2]$, $\Pi_{d+1}(r) \succ 0$, given $\forall r \in [r_1, r_2]$, $S(r) \succeq 0$, $U(r) \succeq 0$.

Given $P_{d+1}(r) = P_d(r) \oplus \mathsf{O}_n$ and consider (6.31) with the structure of $\boldsymbol{\Theta}_d(r)$ and (6.27) and $\boldsymbol{\Phi}_d(r)$ in (6.22), it is obvious to see that

$$\forall r \in [r_1, r_2], \quad \boldsymbol{\Theta}_{d+1}(r) = \begin{bmatrix} J_1 & \begin{bmatrix} \Sigma(r) & \mathsf{O}_{m,n} \end{bmatrix} \\ * & \boldsymbol{\Phi}_{d+1}(r) \end{bmatrix} = \begin{bmatrix} J_1 & \begin{bmatrix} \Sigma(r) & \mathsf{O}_{m,n} \end{bmatrix} \\ * & \boldsymbol{\Phi}_d(r) \oplus [-r(2d+3)U(r)] \end{bmatrix} \prec 0. \tag{6.45}$$

Since $\forall r \in [r_1, r_2], U(r) \succeq 0$, we can conclude that the existence of the feasible solutions to $\forall r \in [r_1, r_2]$, $\boldsymbol{\Theta}_d(r) \prec 0$ implies the existence of a feasible solution to $\forall r \in [r_1, r_2]$, $\boldsymbol{\Theta}_{d+1}(r) \prec 0$. Consequently, we have shown that the existence of feasible solutions to (6.37) and (6.31) at $d$ implies the existence of a feasible solution to (6.37) and (6.31) at $d+1$. Finally, since (6.44) and (6.45) at $d+1$ are equivalently verified by the SoS conditions (6.18)–(6.21) with $d+1$ for some $\delta_7 \in \mathbb{N}$ and $\delta_i; \delta_8 \in \mathbb{N}_0$, $i = 1, \ldots, 6$. This finishes the proof. ∎

## 6.5 Numerical Examples

The effectiveness of the proposed methods is demonstrated through their application to several numerical examples. These examples were coded and computed using MATLAB© and Yalmip [603] as the optimization parser. Mosek [513] was used as the SDP numerical solver. Moreover, all SoS constraints are implemented using the `coefficient` function of Yalmip.

### 6.5.1 Delay range stability analysis

In this subsection, we aim to establish the delay range stability of

$$\begin{aligned} \dot{\boldsymbol{x}}(t) &= A_1 \boldsymbol{x}(t) + A_2 \boldsymbol{y}(t-r) + \int_{-r}^{0} A_3(\tau) \boldsymbol{y}(t+\tau) \mathsf{d}\tau \\ \boldsymbol{y}(t) &= A_4 \boldsymbol{x}(t) + A_5 \boldsymbol{y}(t-r) \end{aligned} \tag{6.46}$$

with different state space parameters presented in Table 6.1, where the delay margins $r_{\min}$ and $r_{\max}$ are calculated via the software package [203] designed for the spectral method in [202].

| Parameters | $A_1$ | $A_2$ | $A_3(r)L_d(\tau)$ | $A_4$ | $A_5$ | $r_{\min}$ | $r_{\max}$ |
|---|---|---|---|---|---|---|---|
| Example 1 | $\begin{bmatrix} 0 & 1 \\ -2 & 0.1 \end{bmatrix}$ | $\begin{bmatrix} 0 \\ 1 \end{bmatrix}$ | $\begin{bmatrix} 0 \\ 0 \end{bmatrix}$ | $\begin{bmatrix} 1 & 0 \end{bmatrix}$ | 0 | 0.10016827 | 1.71785 |
| Example 2 | $\begin{bmatrix} 0.2 & 0.01 \\ 0 & -2 \end{bmatrix}$ | $\mathsf{O}_2$ | $\begin{bmatrix} -1+0.3\tau & 0.1 \\ 0 & -0.1 \end{bmatrix}$ | $\begin{bmatrix} 1 & 0 \\ 0 & 1 \end{bmatrix}$ | $\mathsf{O}_2$ | 0.1944 | 1.7145 |

**Table 6.1:** Numerical Examples of (6.46)

The examples in Table 6.1 are taken from [645] and [568], respectively, which are denoted via equivalent CDDS representations. Note that Example 2 cannot be analyzed by the range stability approaches in [469, 621, 644].

Note that the matrix $A_3(r)$ of Example 1 and Example 2 in Table 6.1 are $A_3(r) = \mathsf{O}_{2,(d+1)}$ and

$$A_3(r) = \begin{bmatrix} -1 & 0.1 & 0.3 & 0 & \mathsf{O}_{2,(2d-2)} \\ 0 & -0.1 & 0 & 0 & \end{bmatrix} \left(\mathsf{L}_d^{-1}(r) \otimes I_2\right)$$
$$= \begin{bmatrix} -1-0.15r & 0.1 & 0.15r & 0 & \mathsf{O}_{2,(2d-2)} \\ 0 & -0.1 & 0 & 0 & \end{bmatrix}, \quad (6.47)$$

respectively. Table 6.2–6.3 presents the results of the detectable stable delay interval with the maximum length calculated using our methodology and compares them to the results reported in [645] and [568], respectively. The column labeled NoDVs represents the number of decision variables, and the values of $\delta_7$ and $\delta_8$ are the degrees of the SoS constraints in (6.21). As mentioned in Section 6.4.4 regarding the reduction of numerical complexity in Theorem 6.1, if any inequality in (6.37) is affine, then it is solved via the convex hull's properties. Lastly, a numerical solution to Example 1 with $r = 1$ produced by DDE23 [650] in MATLAB© is presented in Figure 6.1.[7]

| Solutions for Delay Range Stability | $[r_1, r_2]$ | NoDVs |
|---|---|---|
| [645] ($N = 5$) | $[0.10016829, 1.7178]$ | 294 |
| Theorem 6.1 ($\lambda_1 = 1, \lambda_2 = \lambda_3 = 0,\ d = 4,\ \delta_7 = 1,\ \delta_8 = 0$) | $[0.10016828, 1.71785]$ | 231 |
| Theorem 6.1 ($\lambda_1 = 1, \lambda_2 = \lambda_3 = 0,\ d = 5,\ \delta_7 = 1,\ \delta_8 = 0$) | $[0.10016827, 1.71785]$ | 291 |

**Table 6.2:** Detectable stable interval with the maximum length of Example 1 in Table 6.1.

| Solutions for Delay Range Stability | $[r_1, r_2]$ | NoDVs |
|---|---|---|
| [568] ($l = 1, r = 3$) | $[0.2, 1.29]$ | 5973 |
| [568] ($l = 2, r = 3$) | $[0.2, 1.3]$ | 14034 |
| Theorem 6.1 ($\lambda_1 = \lambda_2 = \lambda_3 = 0,\ d = 4,\ \delta_7 = 2,\ \delta_8 = 1$) | $[0.27, 1.629]$ | 1394 |
| Theorem 6.1 ($\lambda_1 = 1, \lambda_2 = \lambda_3 = 0,\ d = 4,\ \delta_7 = 2,\ \delta_8 = 1$) | $[0.1944, 1.7145]$ | 1472 |

**Table 6.3:** Largest detectable stable interval of Example 2 in Table 6.1

The results presented in Table 6.2–6.3 demonstrate the clear advantage of our proposed methodologies in detecting stable intervals for Example 1 and 2 [8] with fewer variables compared to the existing results reported in [568] and [645]. In addition, the results in Table 6.3 clearly show the benefits of constructing functional (6.24) with delay-dependent matrix parameters for addressing delay range stability problems. Finally, note that Theorem 6.1 does not necessitate the constraint that $A_1 + A_2$ being nonsingular, as is required by [645, Theorem 4].

**Remark 6.9.** Note that the number of decision variables of Theorem 6.1 in Table 6.2–6.3 may be further reduced by simplifying the SoS certificate variable in (6.13) for each instance where a SoS condition must be solved.

### 6.5.2 Range dissipativity and stability analysis

Consider the following neutral delay system

$$\begin{aligned} &\frac{\mathsf{d}}{\mathsf{d}t}\left(\boldsymbol{y}(t) - A_4\boldsymbol{y}(t-r)\right) = A_1\boldsymbol{y}(t) + A_2\boldsymbol{y}(t-r) + \int_{-r}^{0} A_3(r)L_d(\tau)\boldsymbol{y}(t+\tau)\mathsf{d}\tau + D_1\boldsymbol{w}(t) \\ &\boldsymbol{z}(t) = C_1\boldsymbol{y}(t) + C_2\boldsymbol{y}(t-r) + \int_{-r}^{0} C_3(r)L_d(\tau)\boldsymbol{y}(t+\tau) + D_2\boldsymbol{w}(t) \end{aligned} \quad (6.48)$$

[7] Figure 6.1 is generated via `matlab2tikz v1.1.0` by the original figure produced in MATLAB©

[8] Here the stable intervals refer to the ones whose delay margins are calculated by the method in [203] as listed in Table 6.1

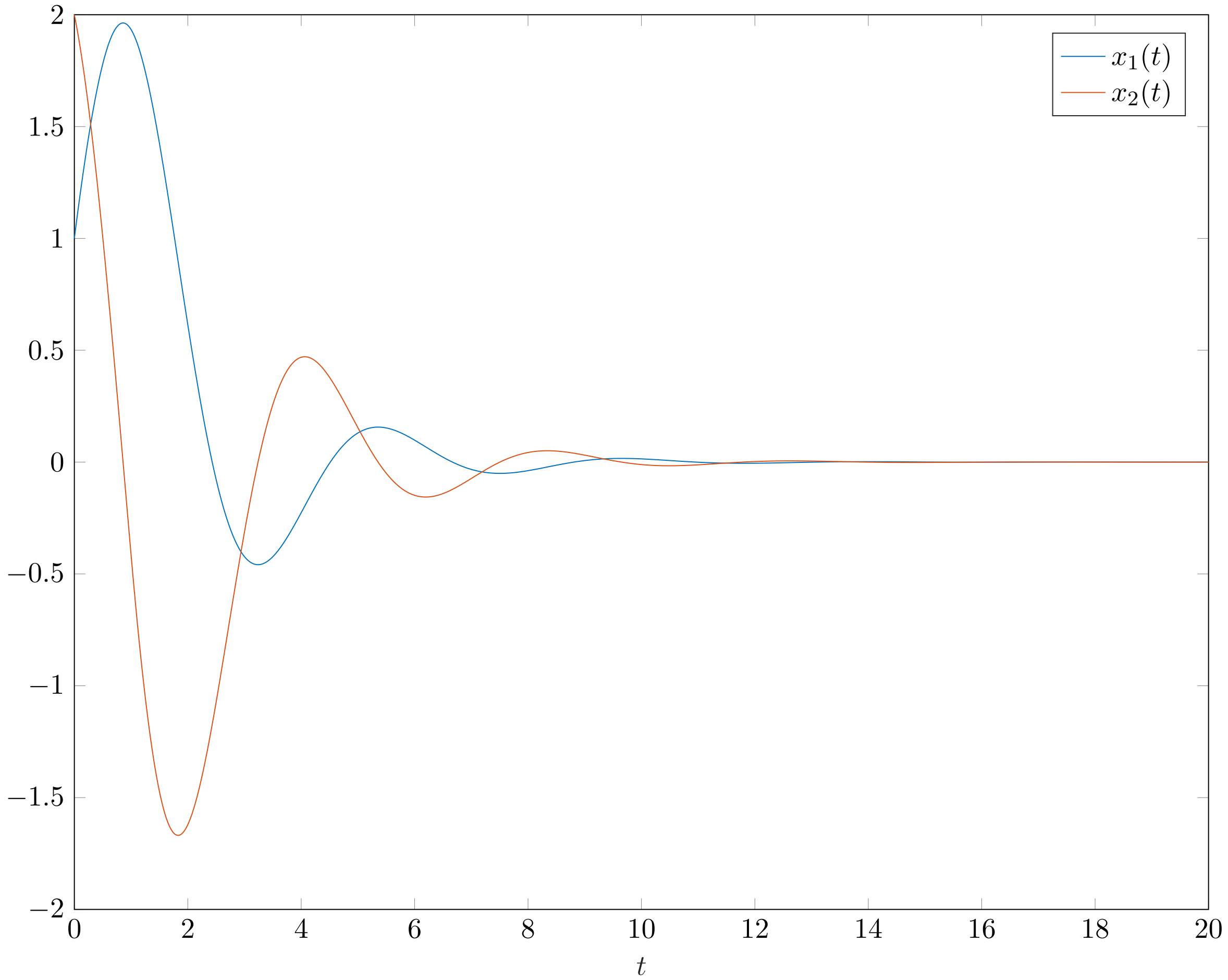

**Figure 6.1:** A numerical solution to Example 1 in Table 6.1

with DDs terms at the state and output. Let $\boldsymbol{x}(t) = \boldsymbol{y}(t) - A_4\boldsymbol{y}(t-r)$, then (6.48) can be rewritten as

$$
\begin{aligned}
\dot{\boldsymbol{x}}(t) &= A_1\boldsymbol{x}(t) + (A_1A_4 + A_2)\boldsymbol{y}(t-r) + \int_{-r}^{0} A_3(r)L_d(\tau)\boldsymbol{y}(t+\tau) + D_1\boldsymbol{w}(t), \\
\boldsymbol{y}(t) &= \boldsymbol{x}(t) + A_4\boldsymbol{y}(t-r), \\
\boldsymbol{z}(t) &= C_1\boldsymbol{x}(t) + (C_1A_4 + C_2)\boldsymbol{y}(t-r) + \int_{-r}^{0} C_3(r)L_d(\tau)\boldsymbol{y}(t+\tau) + D_2\boldsymbol{w}(t).
\end{aligned}
\tag{6.49}
$$

which is now in line with the CDDS form in (6.1). We first consider a linear neutral delay system (6.48) with the parameters

$$
\begin{aligned}
&A_1 = 100\begin{bmatrix} -2.103 & 1 & 2 \\ 3 & -9 & 0 \\ 1 & 2 & -6 \end{bmatrix}, \quad A_2 = 100\begin{bmatrix} 1 & 0 & -3 \\ -0.5 & -0.5 & -1 \\ -0.5 & -1.5 & 0 \end{bmatrix}, \quad A_4 = \frac{1}{72}\begin{bmatrix} -1 & 5 & 2 \\ 4 & 0 & 3 \\ -2 & 4 & 1 \end{bmatrix} \\
&D_1 = \begin{bmatrix} 0 \\ 0 \\ 0.1 \end{bmatrix}, \quad C_1 = \begin{bmatrix} -0.1 & 0.1 & 0.2 \\ 0.4 & 0.01 & 0 \\ 0.1 & 0.21 & 0.1 \end{bmatrix}, \quad C_2 = \begin{bmatrix} 0.1 & 0 & 0.2 \\ 0.4 & 0 & -0.1 \\ 0 & -0.5 & 0.3 \end{bmatrix}, \quad D_2 = \begin{bmatrix} 0 \\ 0.1 \\ 0 \end{bmatrix}
\end{aligned}
\tag{6.50}
$$

with $A_3 = C_3 = \mathsf{O}_{3,3(d+1)}$, which is modified based on the circuit model in [651]. Note that the system's model does not possess any DD. For the performance criteria, we aim to minimize the $\mathcal{L}^2$ gain $\gamma$ corresponding to

$$
J_1 = -\gamma I_m, \quad \widetilde{J} = I_m, \quad J_2 = \mathsf{O}_{m,q}, \quad J_3 = \gamma I_q \tag{6.51}
$$

in (6.17).

Assuming $\lambda_1 = 1$ and $\lambda_2 = \lambda_3 = 0$ in (6.24), with a delay range $[r_1, r_2] = [0.1, 0.5]$. Theorem 6.1 is applied to (6.51)–(6.49) in conjunction with (6.50) and the assignment of $A_3 = C_3 = \mathbf{O}_{3\times 3(d+1)}$. As all inequalities in (6.37) are either affine with respect to $r$ or simple LMIs, they are solved directly via the convex hull's properties rather than solving the SoS conditions in (6.18)–(6.20). The results of $\min \gamma$ over $r \in [0.1, 0.5]$ are enumerated in Table 6.4. Note that $\delta_7$ and $\delta_8$ are the degrees of the SoS constraints in (6.21). Finally, Figure 6.2 depicts a numerical solution to the system in this case at $r = 0.1$ produced by DDENSD [652] in MATLAB©.[9]

| Theorem 6.1 | $\delta_7 = 1,\ \delta_8 = 0$ | $\delta_7 = 2,\ \delta_8 = 1$ | $\delta_7 = 3,\ \delta_8 = 2$ |
|---|---|---|---|
| $d = 1$ | 0.441 | 0.441 | 0.441 |
| $d = 2$ | 0.364 | 0.364 | 0.364 |
| $d = 3$ | 0.361 | 0.361 | 0.361 |

**Table 6.4:** Values of $\min \gamma$ valid over $r \in [0.1, 0.5]$

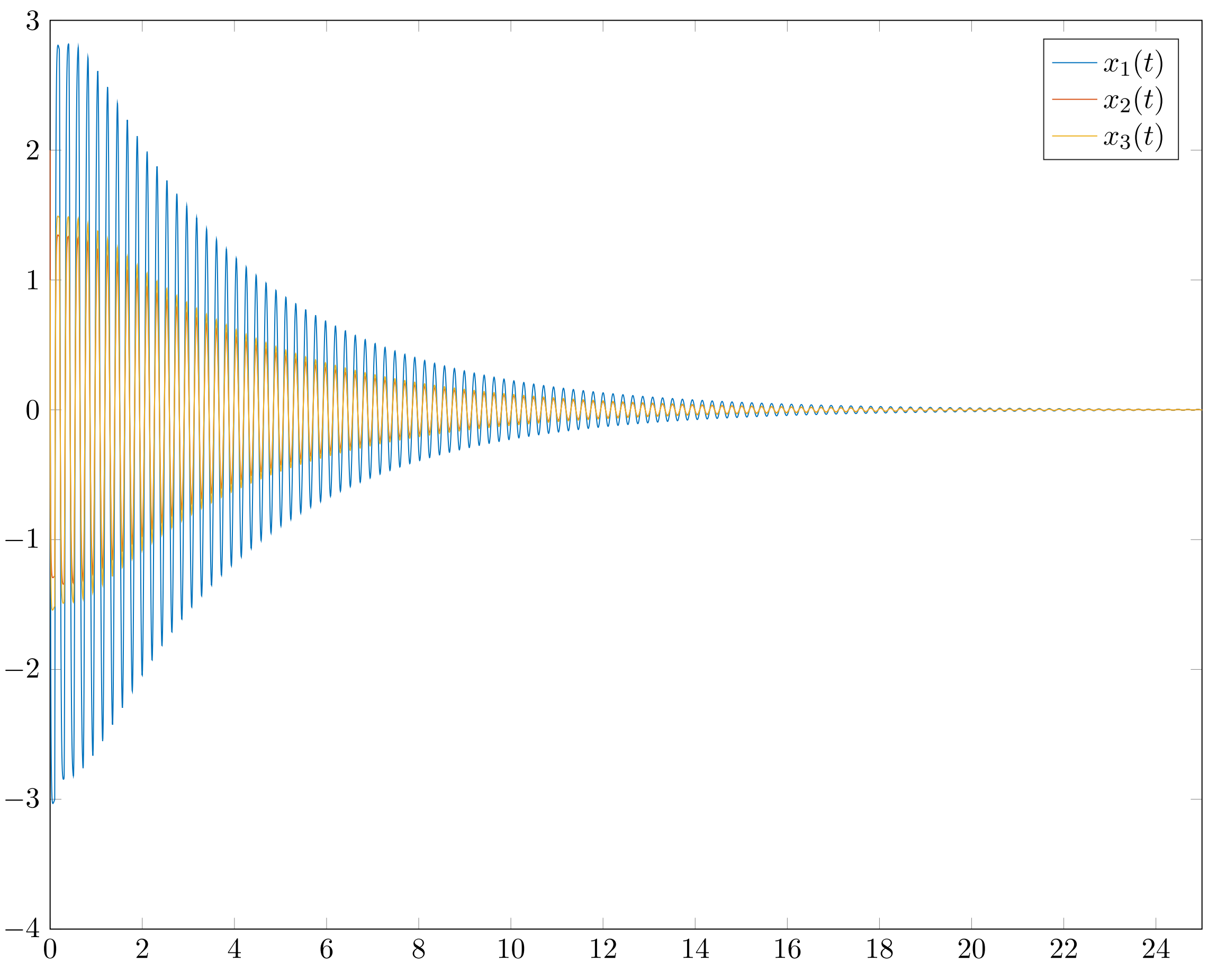

**Figure 6.2:** A numerical solution to (6.49) with (6.50) and $A_3 = C_3 = \mathbf{O}_{3,3(d+1)}$

By applying Theorem 6.1 with constant matrix parameters $\lambda_1 = \lambda_2 = \lambda_3 = 0$ in (6.24) to the same aforementioned model with $[r_1, r_2] = [0.1, 0.5]$, it is observed that the conditions in (6.37) are simple LMIs and $\mathbf{\Theta}_d(r) \prec 0$ in (6.31) can be solved by the convex hull's properties. However, even with $d = 10$, the corresponding range dissipativity and stability condition with $\lambda_1 = \lambda_2 = \lambda_3 = 0$

[9] Figure 6.2 is generated via `matlab2tikz v1.1.0` by the original figure produced in MATLAB©

fails to yield feasible solutions. This illustrates the advantage of applying an LKF with delay-dependent parameters when a functional with constant parameters is insufficient to derive conditions capable of detecting a stable delay interval.

To partially verify the results of $\min\gamma$ in Table 6.4, we apply [10] the `sigma` function in MATLAB© to compute the singular values ($\min\gamma$) of an LTDS system over a fixed frequency range. By extracting the peak value produced by `sigma`, it shows that the system (6.49) with (6.50) and $A_3 = C_3 = \mathsf{O}_{3,3(d+1)}$ guarantees $\min\gamma = 0.101074$ and $\min\gamma = 0.36064$ at $r = 0.1$ and $r = 0.5$, respectively. This demonstrates that the values of $\min\gamma$ in Table 6.4, which are valid over $r \in [0.1, 0.5]$, are consistent with the point-wise $\min\gamma$ values obtained via `sigma` function. Additionally, the minimum value $\min\gamma = 0.361$ in Table 6.4 is very close to the point-wise result $\min\gamma = 0.36064$ at $r = 0.5$.

Now consider new $A_3 L_d(\tau)$ and $C_3 L_d(\tau)$ with parameters

$$A_3(r)L_d(\tau) = \begin{bmatrix} 0.1\tau & 0.1 & 0.3 \\ 0.2 & 0.1 & 0.3-0.1\tau \\ -0.1 & -0.2+0.1\tau & 0.2 \end{bmatrix}, \quad C_3(r)L_d(\tau) = \begin{bmatrix} 0.1 & 0 & 0.2 \\ 0.4 & 0 & -1 \\ 0 & -0.5 & 0.3 \end{bmatrix} \tag{6.52}$$

which along with (6.50) constitute a linear neutral system with DDs. Note that we can easily obtain the corresponding matrix coefficients as

$$\begin{aligned} A_3(r) &= \begin{bmatrix} -0.05r & 0.1 & 0.3 & 0.05r & 0 & 0 & \\ 0.2 & 0.1 & 0.3+0.05r & 0 & 0 & -0.05r & \mathsf{O}_{3,(3d-3)} \\ -0.1 & -0.2-0.05r & 0.2 & 0 & 0.05r & 0 & \end{bmatrix}, \\ C_3(r) &= \begin{bmatrix} 0.1 & 0 & 0.2 & \\ 0.4 & 0 & -1 & \mathsf{O}_{3,3d} \\ 0 & -0.5 & 0.3 & \end{bmatrix}. \end{aligned} \tag{6.53}$$

To the best of our knowledge, there exist no results pertaining to delay range dissipative analysis concerning linear neutral systems with non-constant DD kernels. As such, it can be claimed that no existing schemes can address the problem at hand.

Assume $\lambda_1 = 1$ and $\lambda_2 = \lambda_3 = 0$ in (6.24) with delay range $[r_1, r_2] = [0.1, 0.5]$. Now apply Theorem 6.1 to the system (6.49) with (6.50) and (6.53). Once again, as all corresponding inequalities in (6.37) are affine in $r$, they are solved directly via the convex hull's properties rather than solving (6.18)–(6.20). The resulting values of $\min\gamma$ over $r \in [0.1, 0.5]$ are given in Table 6.5, where $\delta_7$ and $\delta_8$ are the degrees of the SoS constraint (6.21).

| Theorem 6.1 | $\delta_7 = 2$, $\delta_8 = 1$ | $\delta_7 = 3$, $\delta_8 = 2$ |
|---|---|---|
| $d = 1$ | 0.47 | 0.47 |
| $d = 2$ | 0.382 | 0.382 |
| $d = 3$ | 0.37822 | 0.37822 |

**Table 6.5:** values of $\min\gamma$ valid over $[0.1, 0.5]$

Unfortunately, even for the case of point-wise delays, the `sigma` function is currently unable to handle an LDDS.[11] Thus, we suggest using (6.38) to partially verify the results in Table 6.5. Specifically, apply (6.38) with $d = 3$ to the system (6.49) with (6.50) and (6.53) at $r_0 = 0.1$ and $r_0 = 0.5$, respectively. It yields $\min \gamma = 0.10101$ at $r_0 = 0.1$ and $\min \gamma = 0.37822$ at $r_0 = 0.5$, respectively. This verifies that the values of range $\min \gamma$ in Table 6.5 are consistent with the point-wise $\min \gamma$ values we presented.

[10] We first use the default range of `sigma` to determine which frequency range contains the peak singular value. Based on the previous information, we assigned `w = logspace(-1, 2, 2000000)` to `sigma` to ensure the accuracy of $\min\gamma$.

[11] We thank Dr. Suat Gumussoy for providing this information about the `sigma` function in MATLAB©.

Finally, the estimation problem described in Section 6.4.4 can be readily applied to system (6.49) with (6.50) and (6.53) assuming that the value of $\gamma$ in (6.17) is known. Specifically, for system (6.49) with (6.50) and (6.53), let $\gamma = 0.37822$ and $r_0 = 0.5$, with which feasible solutions can be produced by (6.38) with $d = 3$. One can then use (6.40)[12] with $\lambda_1 = 1$ and $\lambda_2 = \lambda_3 = 0$, $d = 3$ and $\delta_7 = 2$, $\delta_8 = 1$ to determine the minimum $r^*$, which renders the corresponding system stable and satisfies $\gamma = 0.37822$ over $[r^*, r_0]$. Given the results in Table 6.5, it is predictable that the optimal value of $r^*$ here is $r^* = 0.1$.

[12] Note that in this case instead of solving the SoS contraints in (6.18)–(6.20) in (6.40), the equivalent inequalities in (6.37) are directly solved for (6.40) via the convex hull's properties.

# Chapter 7

# Dissipative Stabilization of Linear Systems with Uncertain Bounded Time-Varying Distributed Delays

## 7.1 Introduction

In previous chapters, the delay parameters of linear systems were assumed as constants. However, in certain real-time applications such as those described in [60, 653–655], delays can be time-varying [656, 657]. A particular class of delays $r(\cdot) \in \mathcal{M}([r_1, r_2]; \mathbb{R})$ is of significant research interest as it can be applied to model sampled-data [545] or networked control systems (NCSs) [658], or even a delay that is bounded but non-deterministic [659]. This strongly motivates the investigation of solutions for stability analysis and synthesis for linear systems with time-varying delays in the form of $r(\cdot) \in \mathcal{M}([r_1, r_2]; \mathbb{R})$.

There exists a significant body of research pertaining to the stability analysis [539, 625, 660–672] and stabilization [673–677] of linear systems in the form of $\dot{\boldsymbol{x}}(t) = A_1\boldsymbol{x}(t) + A_2\boldsymbol{x}(t - r(t))$ with $r(\cdot) \in \mathcal{M}([r_1, r_2]; \mathbb{R})$, based on the construction of KFs. Furthermore, it has been shown in [541, 678] that the ideas of the LKF approach on $\dot{\boldsymbol{x}}(t) = A_1\boldsymbol{x}(t) + A_2\boldsymbol{x}(t - r(t))$ can be successfully applied to resolve the synthesis problems of NCSs. It is worth noting that, unlike the case of constant delays, frequency-domain-based approaches [8, 202, 203, 215, 224] may not be easily extended to handle systems with time-varying delays $r(\cdot) \in \mathcal{M}([r_1, r_2]; \mathbb{R})$ if the exact expression of $r(t)$ is unknown. This clearly illustrates the advantage and adaptability of the LKF approach when dealing with time-varying delays $r(\cdot) \in \mathcal{M}([r_1, r_2]; \mathbb{R})$ with unknown expressions.

It has been established in [159, 679] that the digital communication channel of Networked Control Systems (NCSs) with stochastic packet delays and loss can be modeled using DDs. However, to the best of our knowledge, no existing results address the problem of stability analysis and stabilization of systems with time-varying $r(\cdot) \in \mathcal{M}([r_1, r_2]; \mathbb{R})$ DDs at system state, input, and output, where the DD kernels can be non-constant. In the Theorem 2 of [680], a solution to stabilization is proposed for systems in the form of $\dot{\boldsymbol{x}}(t) = A\boldsymbol{x}(t) + \int_{-r(t)}^{0} B(\tau)\boldsymbol{u}(t+\tau)\mathrm{d}\tau$. However, all the poles of $A$ in [680] are assumed to be positioned on the imaginary axis, and the delay considered therein is $r(\cdot) \in (0, r_2]^{\mathbb{R}}$. Moreover, the stability of positive linear systems with distributed time-varying delays is investigated in [681, 682]. Although the method in [682] does include criteria that can determine the stability of non-positive linear systems, the delay structure $r(\cdot) \in [0, r_2]^{\mathbb{R}}$ considered in [682] is still rather restrictive. Additionally, the synthesis (stability analysis) results in [159, 169, 463, 470, 683], which are derived to handle LDDSs with constant delay values, may not be easily extended to cope with the case of systems with an unknown time-varying delay $r(\cdot) \in \mathcal{M}([r_1, r_2]; \mathbb{R})$. This is especially true for the approximation approaches in [159, 463, 470, 683] since the resulting approximation coefficients can become nonlinear with respect to $r(t)$ if DD kernels are approximated over the interval $[-r(t), 0]$. Because of this, it is crucial to develop solutions to the stabilization of linear systems with nontrivial DDs where delay function $r(\cdot) \in \mathcal{M}([r_1, r_2]; \mathbb{R})$ is unknown but bounded.

In this chapter, new approaches based on the LKF approach for the design of a state feedback controller for a linear system with DDs are developed where the time-varying delay $r(\cdot) \in \mathcal{M}([r_1, r_2]; \mathbb{R})$ is unknown but bounded by given values $0 \leq r_1 \leq r_2$. We do not consider a pointwise time-varying delay $\boldsymbol{x}(t - r(t))$ [684] at this stage as its presence, as its presence can significantly change the strategies of constructing LKFs. DDs can be found in the system's states, inputs and outputs, ensuring the generality of our model. Additionally, the DD kernels considered in this chapter follow the same class as in Chapter 2. To obtain numerically tractable synthesis (stability) conditions via the construction of an LKF, a novel integral inequality is derived which generalizes Lemma 2 of [685]. Consequently, sufficient conditions for the existence of a state feedback controller, which ensures that the system is stable and dissipative with a supply rate function, are constructed in terms of matrix inequalities set out in the first theorem of this chapter. For the conditions in our first theorem, there is a Bilinear Matrix Inequality (BMI) corresponding to the problem of dissipative synthesis, whereas that inequality becomes convex when non-stabilization problems are considered. To tackle the problem of non-convexity, we derive the second theorem of this chapter via the application of the Projection Lemma where the dissipative synthesis condition is denoted in terms of LMIs. Furthermore, an iterative algorithm is also derived to solve the bilinear condition in the first theorem and the algorithm can be initiated by feasible solutions to the second theorem. To the best of the author's knowledge, no existing methods can handle the synthesis problem considered in this chapter. Finally, two numerical examples are presented to demonstrate the effectiveness of our proposed methodologies.

The remainder of this chapter is organized as follows. The model of our CLS is first derived in Section 7.2. Secondly, some indispensable lemmas and definitions are presented in Section 7.3 which includes the derivation of a novel integral inequality in Lemma 7.3. Next, Section 7.4 sets forth the main results on dissipative synthesis for the existence of a state feedback controller, which are enunciated in Theorem 7.1,7.2 and Algorithm 4. The testing results of two numerical examples are given in Section 7.5 prior to the final conclusion.

## 7.2 Problem Formulation

Consider an LDDS

$$\begin{aligned}
\dot{\boldsymbol{x}}(t) &= A_1\boldsymbol{x}(t) + \int_{-r(t)}^{0} \widetilde{A}_2(\tau)\boldsymbol{x}(t+\tau)\mathrm{d}\tau + B_1\boldsymbol{u}(t) + \int_{-r(t)}^{0} \widetilde{B}_2(\tau)\boldsymbol{u}(t+\tau)\mathrm{d}\tau + D_1\boldsymbol{w}(t), \\
\boldsymbol{z}(t) &= C_1\boldsymbol{x}(t) + \int_{-r(t)}^{0} \widetilde{C}_2(\tau)\boldsymbol{x}(t+\tau)\mathrm{d}\tau + B_4\boldsymbol{u}(t) + \int_{-r(t)}^{0} \widetilde{B}_5(\tau)\boldsymbol{u}(t+\tau)\mathrm{d}\tau + D_2\boldsymbol{w}(t), \\
&\forall \theta \in [-r_2, 0], \quad \boldsymbol{x}(t_0 + \theta) = \boldsymbol{\phi}(\tau), \quad r(\cdot) \in \mathcal{M}(\mathbb{R}; [r_1, r_2])
\end{aligned} \tag{7.1}$$

to be stabilized, where $t_0 \in \mathbb{R}$ and $\boldsymbol{\phi}(\cdot) \in \mathcal{C}([-r_2, 0]; \mathbb{R}^n)$, and $r_2 > 0$, $r_2 \geq r_1 \geq 0$ are given constants. The differential equation in holds for almost all $\widetilde{\forall} t \geq t_0$ with respect to the Lebesgue measure. Furthermore, $\boldsymbol{x} : [t_0 - r_2, \infty) \to \mathbb{R}^n$ satisfies (7.1), $\boldsymbol{u}(t) \in \mathbb{R}^p$ denotes input signals, $\boldsymbol{w}(t) \in \mathcal{L}^2([t_0, +\infty); \mathbb{R}^q)$ represents the disturbance, and $\boldsymbol{z}(t) \in \mathbb{R}^m$ is the regulated output. The size of the given state space parameters in (7.1) is determined by the values of $n \in \mathbb{N}$ and $m; p; q \in \mathbb{N}_0 := \mathbb{N} \cup \{0\}$. Finally, the matrix-valued DDs in (7.1) satisfy

$$\begin{aligned} &\widetilde{A}_2(\cdot) \in \mathscr{L}^2([-r_2, 0]; \mathbb{R}^{n\times n}), \quad \widetilde{C}_2(\cdot) \in \mathscr{L}^2([-r_2, 0]; \mathbb{R}^{m\times n}), \\ &\widetilde{B}_2(\cdot) \in \mathscr{L}^2([-r_2, 0]; \mathbb{R}^{n\times p}), \quad \widetilde{B}_5(\cdot) \in \mathscr{L}^2([-r_2, 0]; \mathbb{R}^{m\times p}). \end{aligned} \tag{7.2}$$

**Remark 7.1.** Systems with DDs and a time-varying delay function can be found in the models of complex and neural networks [148, 149, 686–688]. Moreover, the results in [100, 114, 115] have demonstrated that event-triggered mechanisms can be described as a DD with constant delay values.

**Remark 7.2.** No point-wise time-varying delay $\boldsymbol{x}(t-r(t))$ is considered in (7.1) in this chapter, as its presence can cause significant ramifications to the derivations of synthesis conditions if the Krasovskiĭ functional approach is utilized. Moreover, numerous future works can be done for the stabilization of $\dot{\boldsymbol{x}}(t) = A_1\boldsymbol{x}(t) + A_2\boldsymbol{x}(t - r(t))$ given the limited existing results in the literature. Thus we leave the synthesis problem encompassing both $\boldsymbol{x}(t - r(t))$ and distributed time-varying delays to future research.

The DDs in (7.2) are infinite-dimensional. To construct synthesis constraints of finite dimensions for (7.1), we propose a decomposition scenario as follows.

**Proposition 1.** *The conditions in* (7.2) *hold if and only if there exist* $\boldsymbol{f}_1(\cdot) \in C^1([-r_2, 0]; \mathbb{R}^{d_1})$, $\boldsymbol{f}_2(\cdot) \in C^1([-r_2, 0]; \mathbb{R}^{d_2})$, $\boldsymbol{\varphi}_1(\cdot) \in \mathscr{L}^2([-r_2, 0]; \mathbb{R}^{\delta_1})$, $\boldsymbol{\varphi}_2(\cdot) \in \mathscr{L}^2([-r_2, 0]; \mathbb{R}^{\delta_2})$, $M_1 \in \mathbb{R}^{d_1\times\kappa_1}$, $M_2 \in \mathbb{R}^{d_2\times\kappa_2}$, $A_2 \in \mathbb{R}^{n\times\kappa_1 n}$, $A_3 \in \mathbb{R}^{n\times\kappa_2 n}$, $B_2 \in \mathbb{R}^{n\times\kappa_1 p}$, $B_3 \in \mathbb{R}^{n\times\kappa_2 p}$, $C_2 \in \mathbb{R}^{m\times\kappa_1 n}$, $C_3 \in \mathbb{R}^{m\times\kappa_2 n}$, $B_5 \in \mathbb{R}^{m\times\kappa_1 p}$ *and* $B_6 \in \mathbb{R}^{m\times\kappa_2 p}$ *such that*

$$\forall\tau \in [-r_1, 0], \quad \widetilde{A}_2(\tau) = A_2\left(\widehat{\boldsymbol{f}}_1(\tau) \otimes I_n\right), \quad \widetilde{B}_2(\tau) = B_2\left(\widehat{\boldsymbol{f}}_1(\tau) \otimes I_p\right), \tag{7.3}$$

$$\forall\tau \in [-r_2, -r_1], \quad \widetilde{A}_2(\tau) = A_3\left(\widehat{\boldsymbol{f}}_2(\tau) \otimes I_n\right), \quad \widetilde{B}_2(\tau) = B_3\left(\widehat{\boldsymbol{f}}_2(\tau) \otimes I_p\right), \tag{7.4}$$

$$\forall\tau \in [-r_1, 0], \quad \widetilde{C}_2(\tau) = C_2\left(\widehat{\boldsymbol{f}}_1(\tau) \otimes I_n\right), \quad \widetilde{B}_5(\tau) = B_5\left(\widehat{\boldsymbol{f}}_1(\tau) \otimes I_p\right), \tag{7.5}$$

$$\forall\tau \in [-r_2, -r_1], \quad \widetilde{C}_2(\tau) = C_3\left(\widehat{\boldsymbol{f}}_2(\tau) \otimes I_n\right), \quad \widetilde{B}_5(\tau) = B_6\left(\widehat{\boldsymbol{f}}_2(\tau) \otimes I_p\right), \tag{7.6}$$

$$\forall\tau \in [-r_2, 0], \quad \frac{\mathsf{d}\boldsymbol{f}_1(\tau)}{\mathsf{d}\tau} = M_1\widehat{\boldsymbol{f}}_1(\tau), \quad \frac{\mathsf{d}\boldsymbol{f}_2(\tau)}{\mathsf{d}\tau} = M_2\widehat{\boldsymbol{f}}_2(\tau) \tag{7.7}$$

$$\mathsf{G}_1 = [\,]_{0\times 0} \ \text{ or } \ \mathsf{G}_1 \succ 0, \quad \mathsf{G}_1 := \int_{-r_1}^{0} \widehat{\boldsymbol{f}}_1(\tau)\widehat{\boldsymbol{f}}_1^\top(\tau)\mathsf{d}\tau \tag{7.8}$$

$$\mathsf{G}_2 = [\,]_{0\times 0} \ \text{ or } \ \mathsf{G}_2 \succ 0, \quad \mathsf{G}_2 := \int_{-r_2}^{-r_1} \widehat{\boldsymbol{f}}_2(\tau)\widehat{\boldsymbol{f}}_2^\top(\tau)\mathsf{d}\tau \tag{7.9}$$

*where* $\kappa_1 = d_1 + \delta_1$, $\kappa_2 = d_2 + \delta_2$ *with* $d_1; d_2; \delta_1; \delta_2 \in \mathbb{N}_0$ *satisfying* $d_1 + d_2 > 0$, *and*

$$\widehat{\boldsymbol{f}}_1(\tau) = \begin{bmatrix} \boldsymbol{\varphi}_1(\tau) \\ \boldsymbol{f}_1(\tau) \end{bmatrix}, \quad \widehat{\boldsymbol{f}}_2(\tau) = \begin{bmatrix} \boldsymbol{\varphi}_2(\tau) \\ \boldsymbol{f}_2(\tau) \end{bmatrix}. \tag{7.10}$$

*Finally, the derivatives in* (7.7) *at* $\tau = 0$ *and* $\tau = -r_2$ *are one-sided derivatives. Note that if matrix multiplications in* (7.3)–(7.10) *involve any empty matrix, then it follows the definition and properties of empty matrices in MATLAB.*

*Proof.* First of all, it is straightforward to see that (7.2) is implied by (7.3)–(7.10) since $\boldsymbol{\varphi}_1(\cdot) \in \mathscr{L}^2([-r_2,0];\mathbb{R}^{\delta_1})$, $\boldsymbol{\varphi}_2(\cdot) \in \mathscr{L}^2([-r_2,0];\mathbb{R}^{\delta_2})$, $\boldsymbol{f}_1(\cdot) \in C^1([-r_2,0];\mathbb{R}^{d_1}) \subset \mathscr{L}^2([-r_2,0];\mathbb{R}^{d_1})$ and $\boldsymbol{f}_2(\cdot) \in C^1([-r_2,0];\mathbb{R}^{d_2}) \subset \mathscr{L}^2([-r_2,0];\mathbb{R}^{d_2})$.

We now proceed to prove that (7.2) implies the existence of the parameters in Proposition Proposition 1 satisfying (7.3)–(7.9). Given any $\boldsymbol{f}_1(\cdot) \in C^1([-r_2,0];\mathbb{R}^{d_1})$, $\boldsymbol{f}_2(\cdot) \in C^1([-r_2,0];\mathbb{R}^{d_2})$, one can always find appropriate $\boldsymbol{\varphi}_1(\cdot) \in \mathscr{L}^2([-r_2,0];\mathbb{R}^{\delta_1})$, $\boldsymbol{\varphi}_2(\cdot) \in \mathscr{L}^2([-r_2,0];\mathbb{R}^{\delta_2})$ with $M_1 \in \mathbb{R}^{d_1\times\kappa_1}$ and $M_2 \in \mathbb{R}^{d_2\times\kappa_2}$ such that (7.7)–(7.9) are satisfied with (7.10), where $\mathsf{G}_1 \succ 0$ and $\mathsf{G}_2 \succ 0$ in (7.8) implies that the functions in $\widehat{\boldsymbol{f}}_1(\cdot)$ and $\widehat{\boldsymbol{f}}_2(\cdot)$ in (7.10) are linearly independent[1] in a Lebesgue sense over $[-r_2,0]$ and $[-r_2,-r_1]$, respectively. This is true since $\frac{\mathrm{d}\boldsymbol{f}_1(\tau)}{\mathrm{d}\tau}(\cdot) \in \mathscr{L}^2([-r_1,0];\mathbb{R}^{d_1})$ and $\frac{\mathrm{d}\boldsymbol{f}_1(\tau)}{\mathrm{d}\tau}(\cdot) \in \mathscr{L}^2([-r_2,-r_2];\mathbb{R}^{d_2})$, and the dimensions of $\boldsymbol{\varphi}_1(\tau)$ and $\boldsymbol{\varphi}_2(\tau)$ can always be enlarged with more linearly independent functions. Note that if any vector-valued function $\boldsymbol{f}_1(\tau)$, $\boldsymbol{f}_2(\tau)$, $\boldsymbol{\varphi}_1(\tau)$ and $\boldsymbol{\varphi}_2(\tau)$ is $[\,]_{0\times 1}$, then it can be handled by using empty matrices, as reflected in (7.8) and (7.9).

Given any $\boldsymbol{f}_1(\cdot) \in C^1([-r_2,0];\mathbb{R}^{d_1})$, $\boldsymbol{f}_2(\cdot) \in C^1([-r_2,0];\mathbb{R}^{d_2})$, we have shown that one can always construct appropriate $\boldsymbol{\varphi}_1(\cdot) \in \mathscr{L}^2([-r_2,0];\mathbb{R}^{\delta_1})$, $\boldsymbol{\varphi}_2(\cdot) \in \mathscr{L}^2([-r_2,0];\mathbb{R}^{\delta_2})$ with $M_1$ and $M_2$ such that the conditions in (7.7)–(7.9) are satisfied with (7.10). For this reason, based on the definition of matrix-valued functions and the fact that the dimensions of $\boldsymbol{\varphi}_1(\tau)$ and $\boldsymbol{\varphi}_2(\tau)$ can be increased without limits, one can always construct appropriate constant matrices $A_{2,i}$, $A_{3,i}$, $C_{2,i}$, $C_{3,i}$, $B_{2,i}$, $B_{3,i}$, $B_{5,i}$, $B_{6,i}$ and $\boldsymbol{f}_1(\tau)$, $\boldsymbol{f}_2(\tau)$, $\boldsymbol{\varphi}_1(\tau)$ and $\boldsymbol{\varphi}_2(\tau)$ such that

$$\forall \tau \in [-r_1,0],\ \widetilde{A}_2(\tau) = \sum_{i=1}^{\kappa_1} A_{2,i}g_i(\tau),\ \ \widetilde{C}_2(\tau) = \sum_{i=1}^{\kappa_1} C_{2,i}g_i(\tau), \tag{7.11}$$

$$\forall \tau \in [-r_1,0],\ \widetilde{B}_2(\tau) = \sum_{i=1}^{\kappa_1} B_{2,i}g_i(\tau),\ \ \widetilde{B}_5(\tau) = \sum_{i=1}^{\kappa_1} B_{5,i}g_i(\tau) \tag{7.12}$$

$$\forall \tau \in [-r_2,-r_1],\ \widetilde{A}_2(\tau) = \sum_{i=1}^{\kappa_2} A_{3,i}h_i(\tau),\ \ \widetilde{C}_2(\tau) = \sum_{i=1}^{\kappa_2} C_{3,i}h_i(\tau), \tag{7.13}$$

$$\forall \tau \in [-r_2,-r_1],\ \widetilde{B}_2(\tau) = \sum_{i=1}^{\kappa_2} B_{3,i}h_i(\tau),\ \ \widetilde{B}_5(\tau) = \sum_{i=1}^{\kappa_2} B_{6,i}h_i(\tau) \tag{7.14}$$

$$\boldsymbol{g}^\top(\tau) = \widehat{\boldsymbol{f}}_1^\top(\tau) = \begin{bmatrix} \boldsymbol{\varphi}_1^\top(\tau) & \boldsymbol{f}_1^\top(\tau) \end{bmatrix}^\top,\quad \boldsymbol{h}(\tau) = \widehat{\boldsymbol{f}}_2^\top(\tau) = \begin{bmatrix} \boldsymbol{\varphi}_2^\top(\tau) & \boldsymbol{f}_2^\top(\tau) \end{bmatrix}^\top \tag{7.15}$$

with $\kappa_1;\kappa_2 \in \mathbb{N}_0$, where $\boldsymbol{f}_1(\tau)$, $\boldsymbol{f}_2(\tau)$, $\boldsymbol{\varphi}_1(\tau)$ and $\boldsymbol{\varphi}_2(\tau)$ satisfy (7.7)–(7.9) for some $M_1$ and $M_2$. Now (7.11)–(7.14) can be further rewritten as

$$\begin{aligned}
&\forall \tau \in [-r_1,0],\ \widetilde{A}_2(\tau) = \left[\operatorname*{\mathsf{Row}}_{i=1}^{\kappa_1} A_{2,i}\right]\left(\widehat{\boldsymbol{f}}_1(\tau)\otimes I_n\right), && \widetilde{C}_2(\tau) = \left[\operatorname*{\mathsf{Row}}_{i=1}^{\kappa_1} C_{2,i}\right]\left(\widehat{\boldsymbol{f}}_1(\tau)\otimes I_n\right)\\
&\forall \tau \in [-r_2,-r_1],\ \widetilde{A}_2(\tau) = \left[\operatorname*{\mathsf{Row}}_{i=1}^{\kappa_2} A_{3,i}\right]\left(\widehat{\boldsymbol{f}}_2(\tau)\otimes I_n\right), && \widetilde{C}_2(\tau) = \left[\operatorname*{\mathsf{Row}}_{i=1}^{\kappa_2} C_{3,i}\right]\left(\widehat{\boldsymbol{f}}_2(\tau)\otimes I_n\right)\\
&\forall \tau \in [-r_1,0],\ \widetilde{B}_2(\tau) = \left[\operatorname*{\mathsf{Row}}_{i=1}^{\kappa_1} B_{2,i}\right]\left(\widehat{\boldsymbol{f}}_1(\tau)\otimes I_p\right), && \widetilde{B}_5(\tau) = \left[\operatorname*{\mathsf{Row}}_{i=1}^{\kappa_1} B_{5,i}\right]\left(\widehat{\boldsymbol{f}}_1(\tau)\otimes I_p\right)\\
&\forall \tau \in [-r_2,-r_1],\ \widetilde{B}_2(\tau) = \left[\operatorname*{\mathsf{Row}}_{i=1}^{\kappa_2} B_{3,i}\right]\left(\widehat{\boldsymbol{f}}_2(\tau)\otimes I_p\right), && \widetilde{B}_5(\tau) = \left[\operatorname*{\mathsf{Row}}_{i=1}^{\kappa_2} B_{6,i}\right]\left(\widehat{\boldsymbol{f}}_2(\tau)\otimes I_p\right).
\end{aligned} \tag{7.16}$$

which are in line with the forms in (7.3)–(7.6). Given all the aforementioned statements we have presented, Proposition 1 is proved. ■

[1]See [566, Theorem 7.2.10] for more information

**Remark 7.3.** Proposition 1 can effectively handle the DDs in (7.1), as it uses a group of basis functions to decompose the DDs without appealing to the application of approximations. The potential choices of the functions in (7.3)–(7.6) are further discussed in the next section in light of the construction of an LKF related to $\boldsymbol{f}_1(\cdot)$ and $\boldsymbol{f}_2(\cdot)$.

### 7.2.1 The formulation of our CLS

In this chapter, our goal is to construct a state feedback controller $\boldsymbol{u}(t) = K\boldsymbol{x}(t)$ with $K \in \mathbb{R}^{p\times n}$ to stabilize open-loop system (7.1). Then the corresponding CLS is denoted as

$$
\begin{aligned}
&\dot{\boldsymbol{x}}(t) = A_1\boldsymbol{x}(t) + \int_{-r_1}^{0} A_2\left(\widehat{\boldsymbol{f}}_1(\tau)\otimes I_n\right)\boldsymbol{x}(t+\tau)\mathsf{d}\tau + \int_{-r(t)}^{-r_1} A_3\left(\widehat{\boldsymbol{f}}_2(\tau)\otimes I_n\right)\boldsymbol{x}(t+\tau)\mathsf{d}\tau + B_1K\boldsymbol{x}(t)\\
&+\int_{-r_1}^{0} B_2(I_{\kappa_1}\otimes K)\left(\widehat{\boldsymbol{f}}_1(\tau)\otimes I_n\right)\boldsymbol{x}(t+\tau)\mathsf{d}\tau + \int_{-r(t)}^{-r_1} B_3(I_{\kappa_2}\otimes K)\left(\widehat{\boldsymbol{f}}_2(\tau)\otimes I_n\right)\boldsymbol{x}(t+\tau)\mathsf{d}\tau + D_1\boldsymbol{w}(t)\\
&\boldsymbol{z}(t) = C_1\boldsymbol{x}(t) + \int_{-r_1}^{0} C_2\left(\widehat{\boldsymbol{f}}_1(\tau)\otimes I_n\right)\boldsymbol{x}(t+\tau)\mathsf{d}\tau + \int_{-r(t)}^{-r_1} C_3\left(\widehat{\boldsymbol{f}}_2(\tau)\otimes I_n\right)\boldsymbol{x}(t+\tau)\mathsf{d}\tau + B_4K\boldsymbol{x}(t)\\
&+\int_{-r_1}^{0} B_5(I_{\kappa_1}\otimes K)\left(\widehat{\boldsymbol{f}}_1(\tau)\otimes I_n\right)\boldsymbol{x}(t+\tau)\mathsf{d}\tau + \int_{-r(t)}^{-r_1} B_6(I_{\kappa_2}\otimes K)\left(\widehat{\boldsymbol{f}}_2(\tau)\otimes I_n\right)\boldsymbol{x}(t+\tau)\mathsf{d}\tau + D_2\boldsymbol{w}(t),\\
&\forall\theta\in[-r_2,0],\quad \boldsymbol{x}(t_0+\theta) = \boldsymbol{\phi}(\tau),\qquad r(\cdot)\in\mathcal{M}\left(\mathbb{R};[r_1,r_2]\right)
\end{aligned}
\tag{7.17}
$$

by Lemma 2.1 and Proposition 1, where the decompositions of the DDs are constructed via

$$
\begin{aligned}
\left(\widehat{\boldsymbol{f}}_i(\tau)\otimes I_p\right)K &= \left(\widehat{\boldsymbol{f}}_i(\tau)\otimes I_p\right)(1\otimes K)\\
&= \left(I_{\kappa_i}\widehat{\boldsymbol{f}}_i(\tau)\otimes KI_n\right) = (I_{\kappa_i}\otimes K)\left(\widehat{\boldsymbol{f}}_i(\tau)\otimes I_n\right),\quad i=1,2
\end{aligned}
\tag{7.18}
$$

by (2.1b). Note that (7.17) has different forms for the following three cases $r_2 > r_1 > 0$, and $r_1 = 0$; $r_2 > 0$, and $r_1 = r_2 > 0$.[2] This implies that each case may require a distinct formulation for the corresponding synthesis conditions for (7.17). To avoid deriving three separate synthesis conditions, we rewrite (7.17) as

$$
\begin{aligned}
&\dot{\boldsymbol{x}}(t) = \left(\mathbf{A} + \mathbf{B}_1\left[\left(I_{\widehat{3}+\kappa}\otimes K\right)\oplus \mathsf{O}_q\right]\right)\boldsymbol{\chi}(t),\quad \widetilde{\forall} t\geq t_0\\
&\boldsymbol{z}(t) = \left(\mathbf{C} + \mathbf{B}_2\left[\left(I_{\widehat{3}+\kappa}\otimes K\right)\oplus \mathsf{O}_q\right]\right)\boldsymbol{\chi}(t),\\
&\forall\theta\in[-r_2,0],\quad \boldsymbol{x}(t_0+\theta) = \boldsymbol{\phi}(\theta)
\end{aligned}
\tag{7.19}
$$

with $t_0$ and $\boldsymbol{\phi}(\cdot)$ in (7.1), where $\kappa = \kappa_1 + 2\kappa_2$ and

$$
\mathbf{A} = \begin{bmatrix}\widehat{\mathsf{O}}_{n,n} & A_1 & A_2\left(\sqrt{\mathsf{G}_1}\otimes I_n\right) & A_3\left(\sqrt{\mathsf{G}_2}\otimes I_n\right) & \mathsf{O}_{n,\kappa_2 n} & D_1\end{bmatrix}
\tag{7.20}
$$

$$
\mathbf{B}_1 = \begin{bmatrix}\widehat{\mathsf{O}}_{n,p} & B_1 & B_2\left(\sqrt{\mathsf{G}_1}\otimes I_p\right) & B_3\left(\sqrt{\mathsf{G}_2}\otimes I_p\right) & \mathsf{O}_{n,\kappa_2 p} & \mathsf{O}_{n,q}\end{bmatrix}
\tag{7.21}
$$

$$
\mathbf{C} = \begin{bmatrix}\widehat{\mathsf{O}}_{m,n} & C_1 & C_2\left(\sqrt{\mathsf{G}_1}\otimes I_n\right) & C_3\left(\sqrt{\mathsf{G}_2}\otimes I_n\right) & \mathsf{O}_{m,\kappa_2 n} & D_2\end{bmatrix}
\tag{7.22}
$$

$$
\mathbf{B}_2 = \begin{bmatrix}\widehat{\mathsf{O}}_{m,p} & B_4 & B_5\left(\sqrt{\mathsf{G}_1}\otimes I_p\right) & B_6\left(\sqrt{\mathsf{G}_2}\otimes I_p\right) & \mathsf{O}_{m,\kappa_2 p} & \mathsf{O}_{m,q}\end{bmatrix}
\tag{7.23}
$$

[2] Since (7.17) becomes a delay-free system with $r_1 = r_2 = 0$, hence such a case is not considered here.

$$
\boldsymbol{\chi}(t)=\begin{bmatrix}\widehat{\mathbb{1}}\boldsymbol{x}(t-r_1)\\ \mathbb{1}\boldsymbol{x}(t-r_2)\\ \boldsymbol{x}(t)\\ \int_{-r_1}^{0}\left(\sqrt{\mathsf{G}_1^{-1}}\widehat{\boldsymbol{f}}_1(\tau)\otimes I_n\right)\boldsymbol{x}(t+\tau)\mathsf{d}\tau\\ \int_{-r(t)}^{-r_1}\left(\sqrt{\mathsf{G}_2^{-1}}\widehat{\boldsymbol{f}}_2(\tau)\otimes I_n\right)\boldsymbol{x}(t+\tau)\mathsf{d}\tau\\ \int_{-r_2}^{-r(t)}\left(\sqrt{\mathsf{G}_2^{-1}}\widehat{\boldsymbol{f}}_2(\tau)\otimes I_n\right)\boldsymbol{x}(t+\tau)\mathsf{d}\tau\\ \boldsymbol{w}(t)\end{bmatrix},\quad \widehat{\mathsf{O}}_{n,p}=\begin{cases}\mathsf{O}_{n,2p} & \text{for } r_2>r_1>0\\ \mathsf{O}_{n,p} & \text{for } r_1=r_2>0\\ \mathsf{O}_{n,p} & \text{for } r_1=0;r_2>0\end{cases}\quad \widehat{3}=\begin{cases}3 & \text{for } r_2>r_1>0\\ 2 & \text{for } r_1=r_2>0\\ 2 & \text{for } r_1=0;r_2>0\end{cases} \tag{7.24}
$$

$$
\mathbb{1}=\begin{cases}I_n & \text{for } r_2>r_1\geq 0\\ [\,]_{0\times n} & \text{for } r_1=r_2>0,\end{cases}\qquad \widehat{\mathbb{1}}=\begin{cases}I_n & \text{for } r_2\geq r_1>0\\ [\,]_{0\times n} & \text{for } r_1=0;r_2>0.\end{cases} \tag{7.25}
$$

Note that $\sqrt{X}$ stands for the unique square root of $X\succ 0$ and the terms in (7.20)–(7.23) are obtained by the following relations for $i\in\{1,2\}$:

$$
\left(\widehat{\boldsymbol{f}}_i(\tau)\otimes I_n\right)=\sqrt{\mathsf{G}_i}\sqrt{\mathsf{G}_i^{-1}}\widehat{\boldsymbol{f}}_i(\tau)\otimes I_n=\left(\sqrt{\mathsf{G}_i}\otimes I_n\right)\left(\sqrt{\mathsf{G}_i^{-1}}\widehat{\boldsymbol{f}}_i(\tau)\otimes I_n\right), \tag{7.26}
$$

$$
(I_{\kappa_i}\otimes K)\left(\widehat{\boldsymbol{f}}_i(\tau)\otimes I_n\right)=\left(\sqrt{\mathsf{G}_i}\sqrt{\mathsf{G}_i^{-1}}\otimes K\right)\left(\widehat{\boldsymbol{f}}_i(\tau)\otimes I_n\right)=\left(\sqrt{\mathsf{G}_i}\otimes I_p\right)(I_{\kappa_i}\otimes K)\left(\sqrt{\mathsf{G}_i^{-1}}\widehat{\boldsymbol{f}}_i(\tau)\otimes I_n\right) \tag{7.27}
$$

which themselves can be obtained via (2.1a) with the fact that $\mathsf{G}_1$ and $\mathsf{G}_2$ in (7.8) are invertible[3]. Now the expressions of the CLS in (7.17) at $r_1=r_2>0$ and $r_1=0$; $r_2>0$ can be obtained by (7.19) with $r_1=r_2>0$, $d_2=\delta_2=0$, and $r_1=0$; $r_2>0$, $d_1=\delta_1=0$ in (7.20)–(7.24), respectively.

By introducing the terms $\widehat{\mathsf{O}}$, $\widehat{3}$, $\mathbb{1}$ and $\widehat{\mathbb{1}}$ in (7.24)–(7.25), we can equivalently denote CLS (7.17) in the form of (7.19). This form can characterize all cases of $r_2\geq r_1\geq 0$, $r_2>0$ without introducing redundant terms into the parameters (7.20)–(7.24). This is crucial for deriving well-posed synthesis conditions in this chapter.

[3]Note that $\sqrt{X^{-1}}=\left(\sqrt{X}\right)^{-1}$ for any $X\succ 0$, based on the application of the eigendecomposition of $X\succ 0$

**Remark 7.4.** The existence and uniqueness of the solution to the CLS in (7.19) are guaranteed by [9, Theorem 1.1, Chapter 6] which is developed for a general LTDS. Specifically, consider $\int_{-r(t)}^{0} G(\tau)\boldsymbol{\phi}(\tau)\mathrm{d}\tau$ with $r(\cdot) \in \mathcal{M}(\mathbb{R};[r_1,r_2])$, $r_2 > 0$, $r_2 \geq r_1 \geq 0$ and $G(\cdot) \in \mathcal{L}^2([-r_2,0];\mathbb{R}^{m\times n})$ and $\boldsymbol{\phi}(\cdot) \in \mathcal{C}([-r_2,0];\mathbb{R}^n)$. By using the Cauchy Schwartz inequality with the fact that $G(\cdot) \in \mathcal{L}^2([-r_2,0];\mathbb{R}^{m\times n})$ and $\boldsymbol{\phi}(\cdot) \in \mathcal{C}([-r_2,0];\mathbb{R}^n) \subset \mathcal{L}^2([-r_2,0];\mathbb{R}^n)$, we have

$$
\begin{aligned}
&\left\|\int_{-r(t)}^{0} G(\tau)\boldsymbol{\phi}(\tau)\mathrm{d}\tau\right\|_2 = \left\|\int_{-r(t)}^{0} \underset{i=1}{\overset{m}{\mathsf{Col}}}\, \boldsymbol{g}_i^\top(\tau)\boldsymbol{\phi}(\tau)\mathrm{d}\tau\right\|_2 = \sqrt{\sum_{i=1}^{m}\left(\int_{-r(t)}^{0} \boldsymbol{g}_i^\top(\tau)\boldsymbol{\phi}(\tau)\mathrm{d}\tau\right)^2} \\
&\leq \sqrt{\sum_{i=1}^{m}\left(\int_{-r(t)}^{0} \|\boldsymbol{g}_i(\tau)\|_2^2\,\mathrm{d}\tau \int_{-r(t)}^{0} \|\boldsymbol{\phi}(\tau)\|_2^2\,\mathrm{d}\tau\right)} \leq \sqrt{\sum_{i=1}^{m}\left(\int_{-r_2}^{0} \|\boldsymbol{g}_i(\tau)\|_2^2\,\mathrm{d}\tau \int_{-r_2}^{0} \|\boldsymbol{\phi}(\tau)\|_2^2\,\mathrm{d}\tau\right)} \\
&\leq \sqrt{\sum_{i=1}^{m}\left(\alpha\int_{-r_2}^{0} \|\boldsymbol{\phi}(\cdot)\|_\infty^2\,\mathrm{d}\tau\right)} = \sqrt{mr_2\alpha\,\|\boldsymbol{\phi}(\cdot)\|_\infty^2} = \sqrt{mr_2\alpha}\,\|\boldsymbol{\phi}(\cdot)\|_\infty
\end{aligned}
\tag{7.28}
$$

for some $\alpha > 0$, where $G(\tau) = \mathsf{Col}_{i=1}^{m}\, \boldsymbol{g}_i^\top(\tau)$. Now this shows that all the integrals in (7.19) satisfy the inequality below eq.(1.5) in Chapter 6 of [9]. This is because of identity $\int_{-r(t)}^{0} G(\tau)\boldsymbol{\phi}(\tau)\mathrm{d}\tau = \int_{-r_2}^{0} \mathit{1}(r(t)+\tau)G(\tau)\boldsymbol{\phi}(\tau)\mathrm{d}\tau$, where the function $\mathit{1}(r(t)+\tau)G(\tau)$ is integrable in $\tau$ for all $t \in \mathbb{R}$ and measurable in $t \in \mathbb{R}$ for all $\tau \in [-r_2,0]$.

## 7.3 Important Lemmas and Definition

In this section, we present several lemmas and definitions that are crucial for deriving the synthesis conditions in the following section. A novel integral inequality is also derived for the handling of time-varying delay in the context of constructing LKFs.

The properties of the commutation matrix [689] are critical for deriving the main results of this chapter.

**Lemma 7.1.**

$$
\begin{gathered}
\forall X \in \mathbb{R}^{d\times\delta},\ \ \forall Y \in \mathbb{R}^{n\times m}\ \ \mathsf{K}_{(n,d)}\,(X\otimes Y)\,\mathsf{K}_{(\delta,m)} = Y\otimes X \\
\forall m,n \in \mathbb{N},\quad \mathsf{K}_{(n,m)}^{-1} = \mathsf{K}_{(m,n)} = \mathsf{K}_{(n,m)}^\top
\end{gathered}
\tag{7.29}
$$

*where* $\mathsf{K}_{(n,d)}$ *is the commutation matrix defined by the identity*

$$
\forall A \in \mathbb{R}^{n\times d},\ \ \mathsf{K}_{(n,d)}\,\mathbf{vec}\,(A) = \mathbf{vec}\left(A^\top\right)
$$

*which follows the definition in [689], where* $\mathbf{vec}(\cdot)$ *stands for the vectorization of a matrix. See Section 4.2 of [690] for the definition and more details of* $\mathbf{vec}(\cdot)$.

**Remark 7.5.** Note that for $\mathsf{K}_{(n,d)}$, we have $\mathsf{K}_{(n,1)} = \mathsf{K}_{(1,n)} = I_n$, $\forall n \in \mathbb{N}$ which gives the identity

$$
\mathsf{K}_{(n,d)}\,(\boldsymbol{f}(\tau)\otimes I_n) = \mathsf{K}_{(n,d)}\,(\boldsymbol{f}(\tau)\otimes I_n)\,K_{(1,n)} = I_n\otimes\boldsymbol{f}(\tau)
\tag{7.30}
$$

with $\boldsymbol{f}(\tau) \in \mathbb{R}^d$. The commutation matrix $\mathsf{K}_{(n,d)}$ can be coded by $\mathsf{K}_{(n,d)} = \texttt{vecperm(d,n)}$ that is a function in The Matrix Computation Toolbox [691] for MATLAB©.

Two IIs are presented as follows. The first one can be derived from the proof of Theorem 4.1 in Chapter 4, and the second inequality is specifically derived for the handling of time-varying delays in the next section.

**Lemma 7.2.** *Let $\varpi(\cdot)$ in (4.1) be given with $d \in \mathbb{N}$ and let $U \in \mathbb{S}^n_{\succeq 0}$ with $n \in \mathbb{N}$ and $\boldsymbol{f}(\cdot) \in \mathscr{L}^2_{\varpi}(\mathcal{K};\mathbb{R}^d)$ satisfying*

$$\int_{\mathcal{K}} \varpi(\tau)\boldsymbol{f}(\tau)\boldsymbol{f}^{\top}(\tau)\mathrm{d}\tau \succ 0, \tag{7.31}$$

*then*

$$\begin{aligned}\forall \boldsymbol{x}(\cdot) \in \mathscr{L}^2_{\varpi}(\mathcal{K};\mathbb{R}^n),\ &\int_{\mathcal{K}} \varpi(\tau)\boldsymbol{x}^{\top}(\tau)U\boldsymbol{x}(\tau)\mathrm{d}\tau \\ &\geq \int_{\mathcal{K}} \varpi(\tau)\boldsymbol{x}^{\top}(\tau)F^{\top}(\tau)\mathrm{d}\tau \left(U \otimes \mathsf{F}^{-1}\right) \int_{\mathcal{K}} \varpi(\tau)F(\tau)\boldsymbol{x}(\tau)\mathrm{d}\tau\end{aligned} \tag{7.32}$$

*where $F(\tau) := I_n \otimes \boldsymbol{f}(\tau)$ and $\mathsf{F} = \int_{\mathcal{K}} \varpi(\tau)\boldsymbol{f}(\tau)\boldsymbol{f}^{\top}(\tau)\mathrm{d}\tau$.*

*Proof.* (7.32) can be easily proved by changing the order of the Kronecker product for relevant expressions in the proof of Theorem 4.1. Note that the definition of $\mathsf{F}$ in Lemma 7.2 is different from the definition of $\mathsf{F}$ in Theorem 4.1. Nevertheless, (7.32) is essentially equivalent to (4.3). ■

**Lemma 7.3.** *Let $\varpi(\cdot)$ in (4.1) with $d \in \mathbb{N}$ and $\mathcal{K} = [-r_2, -r_1]$ with $0 \leq r_1 < r_2$ be given. Assume $U \in \mathbb{S}^n_{\succeq 0}$ with $n \in \mathbb{N}$ and $\boldsymbol{f}(\tau) := \mathbf{Col}_{i=1}^{d} f_i(\tau) \in \mathscr{L}^2_{\varpi}([-r_2,-r_1];\mathbb{R}^d)$ satisfying*

$$\int_{-r_2}^{-r_1} \varpi(\tau)\boldsymbol{f}(\tau)\boldsymbol{f}^{\top}(\tau)\mathrm{d}\tau \succ 0, \tag{7.33}$$

*then we have*

$$\begin{aligned}\int_{-r_2}^{-r_1} \varpi(\tau)\boldsymbol{x}^{\top}(\tau)U\boldsymbol{x}(\tau)\mathrm{d}\tau &\geq [*]\left(\begin{bmatrix} U & Y \\ * & U\end{bmatrix} \otimes \mathsf{F}^{-1}\right)\begin{bmatrix}\int_{-\varrho}^{-r_1}\left(I_n \otimes \boldsymbol{f}(\tau)\right)\boldsymbol{x}(\tau)\varpi(\tau)\mathrm{d}\tau \\ \int_{-r_2}^{-\varrho}\left(I_n \otimes \boldsymbol{f}(\tau)\right)\boldsymbol{x}(\tau)\varpi(\tau)\mathrm{d}\tau\end{bmatrix} \\ &= [*]\left(\begin{bmatrix}\mathsf{K}_{(d,n)} & \mathsf{O}_{dn} \\ * & \mathsf{K}_{(d,n)}\end{bmatrix}\left(\begin{bmatrix} U & Y \\ * & U\end{bmatrix} \otimes \mathsf{F}^{-1}\right)\begin{bmatrix}\mathsf{K}_{(n,d)} & \mathsf{O}_{dn} \\ * & \mathsf{K}_{(n,d)}\end{bmatrix}\right)\begin{bmatrix}\int_{-\varrho}^{-r_1}\left(\boldsymbol{f}(\tau) \otimes I_n\right)\boldsymbol{x}(\tau)\varpi(\tau)\mathrm{d}\tau \\ \int_{-r_2}^{-\varrho}\left(\boldsymbol{f}(\tau) \otimes I_n\right)\boldsymbol{x}(\tau)\varpi(\tau)\mathrm{d}\tau\end{bmatrix}\end{aligned} \tag{7.34}$$

*for all $\boldsymbol{x}(\cdot) \in \mathscr{L}^2_{\varpi}(\mathcal{K};\mathbb{R}^n)$ and for any $Y \in \mathbb{R}^{n\times n}$ satisfying $\left[\begin{smallmatrix} U & Y \\ * & U\end{smallmatrix}\right] \succeq 0$, where $\varrho \in [r_1, r_2]$ and $\mathsf{F} := \int_{-r_2}^{-r_1} \varpi(\tau)\boldsymbol{f}(\tau)\boldsymbol{f}^{\top}(\tau)\mathrm{d}\tau$.*

*Proof.* The proof is based on the insights illustrated in Section 4.1 of [685]. Consider the equality

$$\begin{aligned}\int_a^b \varpi(\tau)\boldsymbol{x}^{\top}(\tau)U\boldsymbol{x}(\tau)\mathrm{d}\tau &= \int_{\varrho}^{b} \varpi(\tau)\begin{bmatrix}\boldsymbol{x}(\tau) \\ \mathbf{0}_n\end{bmatrix}^{\top}\begin{bmatrix} U & Y \\ * & U\end{bmatrix}\begin{bmatrix}\boldsymbol{x}(\tau) \\ \mathbf{0}_n\end{bmatrix}\mathrm{d}\tau \\ &+ \int_a^{\varrho} \varpi(\tau)\begin{bmatrix}\mathbf{0}_n \\ \boldsymbol{x}(\tau)\end{bmatrix}^{\top}\begin{bmatrix} U & Y \\ * & U\end{bmatrix}\begin{bmatrix}\mathbf{0}_n \\ \boldsymbol{x}(\tau)\end{bmatrix}\mathrm{d}\tau = \int_a^b \boldsymbol{y}^{\top}(\tau)\begin{bmatrix} U & Y \\ * & U\end{bmatrix}\boldsymbol{y}(\tau)\mathrm{d}\tau\end{aligned} \tag{7.35}$$

which holds for any $Y \in \mathbb{R}^{n\times n}$ with

$$\mathbb{R}^{2n} \ni \boldsymbol{y}(\tau) := \begin{cases}\begin{bmatrix}\boldsymbol{x}(\tau) \\ \mathbf{0}_n\end{bmatrix}, & \forall \tau \in [\varrho, b] \\ \begin{bmatrix}\mathbf{0}_n \\ \boldsymbol{x}(\tau)\end{bmatrix}, & \forall \tau \in [a, \varrho],\end{cases} \qquad \varrho \in [a,b]. \tag{7.36}$$

Let $Y \in \mathbb{R}^{n\times n}$ satisfying $\left[\begin{smallmatrix} U & Y \\ * & U \end{smallmatrix}\right] \succeq 0$, then one can apply (7.32) with (7.29)–(7.30) to the rightmost integral in (7.35) with $\mathcal{K} = [a,b]$ and $\mathsf{f}(\cdot) \in \mathscr{L}^2_{\varpi}(\mathcal{K};\mathbb{R}^d)$ satisfying (7.33). Then we have

$$\begin{aligned}
\int_a^b \varpi(\tau)\boldsymbol{x}^\top(\tau)U\boldsymbol{x}(\tau)\mathsf{d}\tau &= \int_a^b \varpi(\tau)\boldsymbol{y}^\top(\tau)\begin{bmatrix} U & Y \\ * & U \end{bmatrix}\boldsymbol{y}(\tau)\mathsf{d}\tau \\
&\geq [*]\left(\mathsf{F}^{-1}\otimes\begin{bmatrix} U & Y \\ * & U \end{bmatrix}\right)\left(\int_a^b \varpi(\tau)\left(\mathsf{f}(\tau)\otimes I_{2n}\right)\boldsymbol{y}(\tau)\mathsf{d}\tau\right) \\
&= [*]\left(\mathsf{F}^{-1}\otimes\begin{bmatrix} U & Y \\ * & U \end{bmatrix}\right)\left(\int_a^b \varpi(\tau)\mathsf{K}_{2n,d}\left(I_{2n}\otimes\mathsf{f}(\tau)\right)\boldsymbol{y}(\tau)\mathsf{d}\tau\right) \\
= \int_a^b \varpi(\tau)\boldsymbol{y}^\top(\tau)\left(I_{2n}\otimes\mathsf{f}^\top(\tau)\right)\mathsf{d}\tau &\left(\begin{bmatrix} U & Y \\ * & U \end{bmatrix}\otimes\mathsf{F}^{-1}\right)\int_a^b \varpi(\tau)\left(I_{2n}\otimes\mathsf{f}(\tau)\right)\boldsymbol{y}(\tau)\mathsf{d}\tau \qquad (7.37)
\end{aligned}$$

where $\mathsf{F} = \int_a^b \varpi(\tau)\mathsf{f}(\tau)\mathsf{f}^\top(\tau)\mathsf{d}\tau$. Furthermore, it follows that

$$\begin{aligned}
\int_a^b \left(I_{2n}\otimes\mathsf{f}(\tau)\right)\boldsymbol{y}(\tau)\varpi(\tau)\mathsf{d}\tau &= \int_\varrho^b \begin{bmatrix} I_n\otimes\mathsf{f}(\tau) & \mathsf{O}_{dn} \\ \mathsf{O}_{dn} & I_n\otimes\mathsf{f}(\tau) \end{bmatrix}\begin{bmatrix} \boldsymbol{x}(\tau) \\ \mathbf{0}_n \end{bmatrix}\varpi(\tau)\mathsf{d}\tau \\
+\int_a^\varrho \begin{bmatrix} I_n\otimes\mathsf{f}(\tau) & \mathsf{O}_{dn} \\ \mathsf{O}_{dn} & I_n\otimes\mathsf{f}(\tau) \end{bmatrix}&\begin{bmatrix} \mathbf{0}_n \\ \boldsymbol{x}(\tau) \end{bmatrix}\varpi(\tau)\mathsf{d}\tau = \begin{bmatrix} \int_\varrho^b \left[I_n\otimes\mathsf{f}(\tau)\right]\boldsymbol{x}(\tau)\varpi(\tau)\mathsf{d}\tau \\ \int_a^\varrho \left[I_n\otimes\mathsf{f}(\tau)\right]\boldsymbol{x}(\tau)\varpi(\tau)\mathsf{d}\tau \end{bmatrix} \qquad (7.38)
\end{aligned}$$

by the properties of the Kronecker product. Substituting the expressions in (7.38) for the last integral in (7.37) and using (7.30) yield (7.34). ■

**Remark 7.6.** Note that the value of $\mathsf{F}$ in (7.34) is related to the values of $a$ and $b$, but not to the value of $\varrho$. This means that $\varrho$ can be any function as long as its value is bounded by $[a,b]$. This property enables us to handle time-varying delays and derive tractable dissipative conditions in the following section.

Given that the differential equation in (7.17) is defined for almost all $t \geq t_0$ even in the case of $\boldsymbol{w}(t) \equiv \mathbf{0}_n$, the common Lyapunov-Krasovskiĭ stability theorem[4] cannot be applied to (7.17). A Lyapunov-Krasovskiĭ stability criterion is presented as follows which can analyze the stability of (7.17). See Theorem E.1 for the general Lyapunov-Krasovskiĭ stability criterion, which is derived for analyzing the stability of general FDEs subject to the Carathéodory conditions in [9, section 2.6].

**Corollary 7.1.** *Let $\boldsymbol{w}(t) \equiv \mathbf{0}_q$ in (7.19) and $r_2 \geq r_1 \geq 0$, $r_2 > 0$ be given, then the trivial solution $\boldsymbol{x}(t) \equiv \mathbf{0}_n$ of (7.19) is globally uniformly asymptotically stable if there exist $\epsilon_1;\epsilon_2;\epsilon_3 > 0$ and a differentiable functional $\mathsf{v} : \mathcal{C}([-r_2,0];\mathbb{R}^n) \to \mathbb{R}$ such that*

$$\forall\boldsymbol{\phi}(\cdot) \in \mathcal{C}([-r_2,0];\mathbb{R}^n),\ \epsilon_1\left\|\boldsymbol{\phi}(0)\right\|_2^2 \leq \mathsf{v}(\boldsymbol{\phi}(\cdot)) \leq \epsilon_2\left\|\boldsymbol{\phi}(\cdot)\right\|_\infty^2, \qquad (7.39)$$

$$\widetilde{\forall} t \geq t_0 \in \mathbb{R},\ \frac{\mathsf{d}}{\mathsf{d}t}\mathsf{v}(\mathbf{x}_t(\cdot)) \leq -\epsilon_3\left\|\boldsymbol{x}(t)\right\|_2^2 \qquad (7.40)$$

*where $\mathbf{x}_t(\cdot)$ in (7.40) is defined by the equality $\forall t \geq t_0$, $\forall\theta \in [-r_2,0]$, $\mathbf{x}_t(\theta) = \boldsymbol{x}(t+\theta)$ in which $\boldsymbol{x}(\cdot) \in \mathcal{C}\left(\mathbb{R}_{\geq t_0 - r_2};\mathbb{R}^n\right)$ satisfies $\dot{\boldsymbol{x}}(t) = \left(\mathbf{A}+\mathbf{B}_1\left[\left(I_{\widehat{3}+\kappa}\otimes K\right)\oplus\mathsf{O}_q\right]\right)\boldsymbol{\chi}(t)$ in (7.19) with $\boldsymbol{w}(t) \equiv \mathbf{0}_q$.*

*Proof.* Since (7.19) with $\boldsymbol{w}(t) \equiv \mathbf{0}_q$ is a linear system and $r(\cdot) \in \mathcal{M}\left(\mathbb{R};[r_1,r_2]\right)$, the system (7.19) with $\boldsymbol{w}(t) \equiv \mathbf{0}_q$ is a special case of the general time-varying system (E.1). Then (7.39) and (7.40) can be obtained by letting $\alpha_1(s) = \epsilon_1 s^2$, $\alpha_2(s) = \epsilon_2 s^2$, $\alpha_3(s) = \epsilon_3 s^2$ with $\epsilon_1;\epsilon_2;\epsilon_3 > 0$. ■

[4] See [9, Theorem 2, Section 5.1] and [178, Theorem 1.3]

**Definition 7.1.** Given $0 \neq r_2 \geq r_1 \geq 0$, the CLS in (7.19) with a supply rate function $\mathsf{s}(\boldsymbol{z}(t), \boldsymbol{w}(t))$ is dissipative if there exists a differentiable functional $\mathsf{v} : \mathcal{C}([-r_2, 0]; \mathbb{R}^n) \to \mathbb{R}$ such that

$$\widetilde{\forall} t \geq t_0, \quad \dot{\mathsf{v}}(\mathbf{x}_t(\cdot)) - \mathsf{s}(\boldsymbol{z}(t), \boldsymbol{w}(t)) \leq 0 \tag{7.41}$$

where $t_0$, $\boldsymbol{z}(t)$ and $\boldsymbol{w}(t)$ are given in (7.19) along with $\forall t \geq t_0$, $\forall \theta \in [-r_2, 0]$, $\mathbf{x}_t(\theta) = \boldsymbol{x}(t+\theta)$ where $\boldsymbol{x}(t)$ is the solution to the system (7.19).

Note that (7.41) implies the origin integral-based definition of dissipativity via the properties of Lebesgue integrals. To denote dissipativity, we apply the same quadratic supply function

$$\begin{gathered}
\mathsf{s}(\boldsymbol{z}(t), \boldsymbol{w}(t)) = \begin{bmatrix} \boldsymbol{z}(t) \\ \boldsymbol{w}(t) \end{bmatrix}^\top \begin{bmatrix} \widetilde{J}^\top J_1^{-1} \widetilde{J} & J_2 \\ * & J_3 \end{bmatrix} \begin{bmatrix} \boldsymbol{z}(t) \\ \boldsymbol{w}(t) \end{bmatrix}, \\
\mathbb{S}^m \ni \widetilde{J}^\top J_1^{-1} \widetilde{J} \preceq 0, \ \ \mathbb{S}^m \ni J_1^{-1} \prec 0, \ \ \widetilde{J} \in \mathbb{R}^{m \times m}
\end{gathered} \tag{7.42}$$

as in (2.10).

## 7.4 Dissipative Controller Synthesis

The results on dissipative controller synthesis are presented in this section enunciated in Theorem 7.1, 7.2 and Algorithm 4.

**Theorem 7.1.** *Let $r_2 > r_1 > 0$ and all the parameters in Proposition 1 be given, then the CLS in (7.19) with supply rate function (7.42) is dissipative and the trivial solution $\boldsymbol{x}(t) \equiv \mathbf{0}_n$ of (7.19) with $\boldsymbol{w}(t) \equiv \mathbf{0}_q$ is globally uniformly asymptotically stable if there exist $K \in \mathbb{R}^{p \times n}$ and $P_1 \in \mathbb{S}^n$, $P_2 \in \mathbb{R}^{n \times \varrho}$, $P_3 \in \mathbb{S}^\varrho$ with $\varrho = (d_1 + d_2)n$ and $Q_1; Q_2; R_1; R_2 \in \mathbb{S}^n$ and $Y \in \mathbb{R}^{n \times n}$ such that*

$$\begin{bmatrix} P_1 & P_2 \\ * & P_3 \end{bmatrix} + \left(\mathsf{O}_n \oplus [I_{d_1} \otimes Q_1] \oplus [I_{d_2} \otimes Q_2]\right) \succ 0, \tag{7.43}$$

$$Q_1 \succeq 0, \ \ Q_2 \succeq 0, \ \ R_1 \succeq 0, \ \begin{bmatrix} R_2 & Y \\ * & R_2 \end{bmatrix} \succeq 0, \tag{7.44}$$

$$\begin{bmatrix} \boldsymbol{\Psi} & \boldsymbol{\Sigma}^\top \widetilde{J}^\top \\ * & J_1 \end{bmatrix} = \mathsf{Sy}\left[\mathbf{P}^\top \boldsymbol{\Pi}\right] + \boldsymbol{\Phi} \prec 0 \tag{7.45}$$

*where $\boldsymbol{\Sigma} = \mathbf{C} + \mathbf{B}_2 \left[\left(I_{\widehat{3}+\kappa} \otimes K\right) \oplus \mathsf{O}_q\right]$ with $\mathbf{C}$ and $\mathbf{B}_2$ in (7.22) and (7.23), and*

$$\boldsymbol{\Psi} = \mathsf{Sy}\left( \begin{bmatrix} \widehat{\mathsf{O}}_{n,n}^\top & \widehat{\mathsf{O}}_{\varrho,n}^\top \\ I_n & \mathsf{O}_{n,\varrho} \\ \mathsf{O}_{\kappa n,n} & \widehat{I}^\top \\ \mathsf{O}_{q,n} & \mathsf{O}_{q,\varrho} \end{bmatrix} \begin{bmatrix} P_1 & P_2 \\ * & P_3 \end{bmatrix} \begin{bmatrix} \mathbf{A} + \mathbf{B}_1 \left[\left(I_{\widehat{3}+\kappa} \otimes K\right) \oplus \mathsf{O}_q\right] \\ \left[\widehat{\mathbf{F}} \otimes I_n \quad \mathsf{O}_{\varrho,q}\right] \end{bmatrix} - \begin{bmatrix} \widehat{\mathsf{O}}_{m,n}^\top \\ \mathsf{O}_{(n+\kappa n),m} \\ J_2^\top \end{bmatrix} \boldsymbol{\Sigma} \right) - \Xi \tag{7.46}$$

$$\widehat{I} = \left(\sqrt{\mathsf{F}_1^{-1}} \oplus \sqrt{\mathsf{F}_2^{-1}}\right) \begin{bmatrix} \mathsf{O}_{d_1,\delta_1} & I_{d_1} & \mathsf{O}_{d_1,\delta_2} & \mathsf{O}_{d_1,d_2} & \mathsf{O}_{d_1,\delta_2} & \mathsf{O}_{d_1,d_2} \\ \mathsf{O}_{d_2,\delta_1} & \mathsf{O}_{d_2,d_1} & \mathsf{O}_{d_2,\delta_2} & I_{d_2} & \mathsf{O}_{d_2,\delta_2} & I_{d_2} \end{bmatrix} \left(\sqrt{\mathsf{G}}_1 \oplus \sqrt{\mathsf{G}}_2 \oplus \sqrt{\mathsf{G}}_2\right) \otimes I_n \tag{7.47}$$

$$\begin{aligned}
\Xi = \Bigg[ & [Q_1 - Q_2 - r_3 R_2] \oplus [\mathbb{1} Q_2] \oplus \left[\widehat{\mathbb{1}}(-Q_1 - r_1 R_1)\right] \oplus (I_{\kappa_1} \otimes R_1) \\
& \oplus \left( \begin{bmatrix} \mathsf{K}_{(\kappa_2,n)} & \mathsf{O}_{\kappa_2 n} \\ * & \mathsf{K}_{(\kappa_2,n)} \end{bmatrix} \left( \begin{bmatrix} R_2 & Y \\ * & R_2 \end{bmatrix} \otimes I_{\kappa_2} \right) \begin{bmatrix} \mathsf{K}_{(n,\kappa_2)} & \mathsf{O}_{\kappa_2 n} \\ * & \mathsf{K}_{(n,\kappa_2)} \end{bmatrix} \right) \oplus J_3 \Bigg]
\end{aligned} \tag{7.48}$$

$$\widehat{\mathbf{F}} = \begin{bmatrix} -\sqrt{\mathsf{F}_1^{-1}} \boldsymbol{f}_1(-r_1) & \mathbf{0}_{d_1} & \sqrt{\mathsf{F}_1^{-1}} \boldsymbol{f}_1(0) & -\sqrt{\mathsf{F}_1^{-1}} M_1 \sqrt{\mathsf{G}_1} & \mathsf{O}_{d_1,\kappa_2} & \mathsf{O}_{d_1,\kappa_2} \\ \sqrt{\mathsf{F}_2^{-1}} \boldsymbol{f}_2(-r_1) & -\sqrt{\mathsf{F}_2^{-1}} \boldsymbol{f}_2(-r_2) & \mathbf{0}_{d_2} & \mathsf{O}_{d_2,\kappa_1} & -\sqrt{\mathsf{F}_2^{-1}} M_2 \sqrt{\mathsf{G}_2} & -\sqrt{\mathsf{F}_2^{-1}} M_2 \sqrt{\mathsf{G}_2} \end{bmatrix} \tag{7.49}$$

*with* $\mathbf{A}, \mathbf{B}_1$ *in* (7.20)–(7.21) *and* $\mathbb{1}, \widehat{\mathbb{1}}$ *in* (7.25) *and* $\mathsf{G}_1, \mathsf{G}_2$ *in* (7.8)–(7.9)*. Moreover,* $\mathsf{F}_1 = \int_{-r_1}^{0} \boldsymbol{f}_1(\tau)\boldsymbol{f}_1^\top(\tau)\mathsf{d}\tau$ *and* $\mathsf{F}_2 = \int_{-r_2}^{-r_1} \boldsymbol{f}_2(\tau)\boldsymbol{f}_2^\top(\tau)\mathsf{d}\tau$ *and the rest of the parameters in* (7.45) *is defined as*

$$\mathbf{P} = \begin{bmatrix} \widehat{\mathsf{O}}_{n,n} & P_1 & P_2\widehat{I} & \mathsf{O}_{n,q} & \mathsf{O}_{n,m} \end{bmatrix}, \quad \mathbf{\Pi} = \begin{bmatrix} \mathbf{A} + \mathbf{B}_1 \left[ \left( I_{\widehat{3}+\kappa} \otimes K \right) \oplus \mathsf{O}_q \right] & \mathsf{O}_{n,m} \end{bmatrix} \tag{7.50}$$

*and*

$$\mathbf{\Phi} = \mathsf{Sy}\left( \begin{bmatrix} \widehat{\mathsf{O}}_{\varrho,n}^\top \\ P_2 \\ \widehat{I}^\top P_3 \\ \mathsf{O}_{(q+m),\varrho} \end{bmatrix} \begin{bmatrix} \widehat{\mathbf{F}} \otimes I_n & \mathsf{O}_{\varrho,(q+m)} \end{bmatrix} + \begin{bmatrix} \widehat{\mathsf{O}}_{m,n}^\top \\ \mathsf{O}_{(n+\kappa n),m} \\ -J_2^\top \\ \widetilde{J} \end{bmatrix} \begin{bmatrix} \mathbf{\Sigma} & \mathsf{O}_m \end{bmatrix} \right) - \Xi \oplus (-J_1)\,. \tag{7.51}$$

*Furthermore, with* $r_1 = r_2$*,* $d_2 = \delta_2 = 0$ *and* $Q_2 = R_2 = Y = \mathsf{O}_n$*, then the inequalities in* (7.43)–(7.45) *are a dissipative synthesis condition for the CLS in* (7.19) *with* $r_1 = r_2 > 0$*. Finally, with* $r_2 > 0$*;* $r_1 = 0$*,* $d_1 = \delta_1 = 0$ *and* $Q_1 = R_1 = \mathsf{O}_n$*, then the inequalities in* (7.43)–(7.45) *are a dissipative synthesis condition for the CLS in* (7.19) *with* $r_2 > 0$*;* $r_1 = 0$*.*

*Proof.* The proof of Theorem 7.1 is via the construction of

$$\begin{aligned} \mathsf{v}(\mathbf{x}_t(\cdot)) &= \boldsymbol{\eta}^\top(t) \begin{bmatrix} P_1 & P_2 \\ * & P_3 \end{bmatrix} \boldsymbol{\eta}(t) + \int_{-r_1}^{0} \boldsymbol{x}^\top(t+\tau)\left[ Q_1 + (\tau + r_1)R_1 \right] \boldsymbol{x}(t+\tau)\mathsf{d}\tau \\ &+ \int_{-r_2}^{-r_1} \boldsymbol{x}^\top(t+\tau) \left[ Q_2 + (\tau + r_2)R_2 \right] \boldsymbol{x}(t+\tau)\mathsf{d}\tau, \end{aligned} \tag{7.52}$$

where $\mathbf{x}_t(\cdot)$ follows the same definition in (7.41), and $P_1 \in \mathbb{S}^n$, $P_2 \in \mathbb{R}^{n\times\varrho}$, $P_3 \in \mathbb{S}^\varrho$ with $\varrho = (d_1+d_2)n$, and $Q_1; Q_2; R_1; R_2 \in \mathbb{S}^n$ and

$$\boldsymbol{\eta}(t) := \mathbf{Col}\left[ \boldsymbol{x}(t), \int_{-r_1}^{0} \left( \sqrt{\mathsf{F}_1^{-1}}\boldsymbol{f}_1(\tau) \otimes I_n \right) \boldsymbol{x}(t+\tau)\mathsf{d}\tau, \int_{-r_2}^{-r_1} \left( \sqrt{\mathsf{F}_2^{-1}}\boldsymbol{f}_2(\tau) \otimes I_n \right) \boldsymbol{x}(t+\tau)\mathsf{d}\tau \right] \tag{7.53}$$

with $\mathsf{F}_1 = \int_{-r_1}^{0} \boldsymbol{f}_1(\tau)\boldsymbol{f}_1^\top(\tau)\mathsf{d}\tau$ and $\mathsf{F}_2 = \int_{-r_2}^{-r_1} \boldsymbol{f}_2(\tau)\boldsymbol{f}_2^\top(\tau)\mathsf{d}\tau$. Note that given the conditions in (7.8)–(7.9), both $\sqrt{\mathsf{F}_1^{-1}}$ and $\sqrt{\mathsf{F}_2^{-1}}$ are well defined.

We first prove this theorem for the case of $r_2 > r_1 > 0$. Then the synthesis conditions for the cases of $r_1 = r_2 > 0$ and $r_1 = 0$; $r_2 > 0$ can be easily obtained based on the synthesis condition for $r_2 > r_1 > 0$, respectively.

Now given $t_0 \in \mathbb{R}$ in (7.19) with $r_2 > r_1 > 0$, differentiating $\mathsf{v}(\mathbf{x}_t(\cdot))$ along the trajectory of (7.19) and consider (7.42) produces

$$\begin{aligned} &\widetilde{\forall} t \geq t_0, \quad \dot{\mathsf{v}}(\mathbf{x}_t(\cdot)) - \mathsf{s}(\boldsymbol{z}(t), \boldsymbol{w}(t)) \\ &= \boldsymbol{\chi}^\top(t)\,\mathsf{Sy}\left( \begin{bmatrix} \mathsf{O}_{2n,n} & \mathsf{O}_{2n,\varrho} \\ I_n & \mathsf{O}_{n,\varrho} \\ \mathsf{O}_{\kappa n,n} & \widehat{I}^\top \\ \mathsf{O}_{q,n} & \mathsf{O}_{q,\varrho} \end{bmatrix} \begin{bmatrix} P_1 & P_2 \\ * & P_3 \end{bmatrix} \begin{bmatrix} \mathbf{A} + \mathbf{B}_1 \left[ \left( I_{\widehat{3}+\kappa} \otimes K \right) \oplus \mathsf{O}_q \right] \\ \begin{bmatrix} \widehat{\mathbf{F}} \otimes I_n & \mathsf{O}_{\varrho,q} \end{bmatrix} \end{bmatrix} - \begin{bmatrix} \mathsf{O}_{(3n+\kappa n),m} \\ J_2^\top \end{bmatrix} \mathbf{\Sigma} \right) \boldsymbol{\chi}(t) \\ &+\boldsymbol{x}^\top(t)\left( Q_1 + r_1R_1 \right)\boldsymbol{x}(t) - \boldsymbol{x}^\top(t-r_2)Q_2\boldsymbol{x}(t-r_2) - \boldsymbol{x}^\top(t-r_1)\left( Q_1 - Q_2 - r_3R_2 \right)\boldsymbol{x}(t-r_1) \\ &- \boldsymbol{w}^\top(t)J_3\boldsymbol{w}(t) - \int_{-r_1}^{0} \boldsymbol{x}^\top(t+\tau)R_1\boldsymbol{x}(t+\tau)\mathsf{d}\tau - \int_{-r_2}^{-r_1} \boldsymbol{x}^\top(t+\tau)R_2\boldsymbol{x}(t+\tau)\mathsf{d}\tau \\ &- \boldsymbol{\chi}^\top(t)\mathbf{\Sigma}^\top\widetilde{J}^\top J_1^{-1}\widetilde{J}\mathbf{\Sigma}\boldsymbol{\chi}(t) \end{aligned} \tag{7.54}$$

where $\boldsymbol{\chi}(t)$ is given in (7.24) and $\mathbf{\Sigma}$, $\widehat{I}$ and $\widehat{\mathbf{F}}$ are defined in the statements of Theorem 7.1. Note that the expression of $\widehat{\mathbf{F}}$ in (7.49) is obtained by the relations

$$\int_{-r_1}^{0}\left(\sqrt{\mathsf{F}_1^{-1}}\boldsymbol{f}_1(\tau)\otimes I_n\right)\dot{\boldsymbol{x}}(t+\tau)\mathrm{d}\tau=\left(\sqrt{\mathsf{F}_1^{-1}}\boldsymbol{f}_1(0)\otimes I_n\right)\boldsymbol{x}(t)$$
$$-\left(\sqrt{\mathsf{F}_1^{-1}}\boldsymbol{f}_1(-r_1)\otimes I_n\right)\boldsymbol{x}(t-r_1)-\left(\sqrt{\mathsf{F}_1^{-1}}M_1\sqrt{\mathsf{G}_1}\otimes I_n\right)\int_{-r_1}^{0}\left(\sqrt{\mathsf{G}_1^{-1}}\widehat{\boldsymbol{f}}_1(\tau)\otimes I_n\right)\boldsymbol{x}(t+\tau)\mathrm{d}\tau, \quad (7.55)$$

$$\int_{-r_2}^{-r_1}\left(\sqrt{\mathsf{F}_2^{-1}}\boldsymbol{f}_2(\tau)\otimes I_n\right)\dot{\boldsymbol{x}}(t+\tau)\mathrm{d}\tau=\left(\sqrt{\mathsf{F}_2^{-1}}\boldsymbol{f}_2(-r_1)\otimes I_n\right)\boldsymbol{x}(t-r_1)-\left(\sqrt{\mathsf{F}_2^{-1}}\boldsymbol{f}_2(-r_2)\otimes I_n\right)\boldsymbol{x}(t-r_2),$$
$$-\left(\sqrt{\mathsf{F}_2^{-1}}M_2\sqrt{\mathsf{G}_2}\otimes I_n\right)\int_{-r(t)}^{-r_1}\left(\sqrt{\mathsf{G}_2^{-1}}\widehat{\boldsymbol{f}}_2(\tau)\otimes I_n\right)\boldsymbol{x}(t+\tau)\mathrm{d}\tau$$
$$-\left(\sqrt{\mathsf{F}_2^{-1}}M_2\sqrt{\mathsf{G}_2}\otimes I_n\right)\int_{-r_2}^{-r(t)}\left(\sqrt{\mathsf{G}_2^{-1}}\widehat{\boldsymbol{f}}_2(\tau)\otimes I_n\right)\boldsymbol{x}(t+\tau)\mathrm{d}\tau, \quad (7.56)$$

which are derived via (7.7)–(7.10) and (2.1a)–(2.1b). Moreover, the structure of $\widehat{I}$ in (7.54) and (7.47) is obtained based on the identities

$$\boldsymbol{f}_1(\tau)=\begin{bmatrix}\mathsf{O}_{d_1,\delta_1} & I_{d_1}\end{bmatrix}\widehat{\boldsymbol{f}}_1(\tau),\quad \boldsymbol{f}_2(\tau)=\begin{bmatrix}\mathsf{O}_{d_2,\delta_2} & I_{d_2}\end{bmatrix}\widehat{\boldsymbol{f}}_2(\tau), \quad (7.57)$$

$$\begin{bmatrix}\int_{-r_1}^{0}\left(\sqrt{\mathsf{F}_1^{-1}}\boldsymbol{f}_1(\tau)\otimes I_n\right)\boldsymbol{x}(t+\tau)\mathrm{d}\tau\\ \int_{-r_2}^{-r_1}\left(\sqrt{\mathsf{F}_2^{-1}}\boldsymbol{f}_2(\tau)\otimes I_n\right)\boldsymbol{x}(t+\tau)\mathrm{d}\tau\end{bmatrix}=\widehat{I}\begin{bmatrix}\int_{-r_1}^{0}\left(\sqrt{\mathsf{G}_1^{-1}}\widehat{\boldsymbol{f}}_1(\tau)\otimes I_n\right)\boldsymbol{x}(t+\tau)\mathrm{d}\tau\\ \int_{-r(t)}^{-r_1}\left(\sqrt{\mathsf{G}_2^{-1}}\widehat{\boldsymbol{f}}_2(\tau)\otimes I_n\right)\boldsymbol{x}(t+\tau)\mathrm{d}\tau\\ \int_{-r_2}^{-r(t)}\left(\sqrt{\mathsf{G}_2^{-1}}\widehat{\boldsymbol{f}}_2(\tau)\otimes I_n\right)\boldsymbol{x}(t+\tau)\mathrm{d}\tau\end{bmatrix} \quad (7.58)$$

in light of the structures of $\boldsymbol{\eta}(t)$ in (7.53) and $\boldsymbol{\chi}(t)$ in (7.24) and the property of the Kronecker product in (2.1b). Note that also the parameters $\mathbf{A}$, $\mathbf{B}_1$, $\mathbf{C}$ and $\mathbf{B}_2$ in (7.54) are given in (7.20)–(7.23).

Let $R_1\succeq 0$ and $\begin{bmatrix}R_2 & Y\\ * & R_2\end{bmatrix}\succeq 0$ with $Y\in\mathbb{R}^{n\times n}$. Now apply (7.32) and (7.34) with $\varpi(\tau)=1$ and $\mathsf{f}(\tau)=\sqrt{\mathsf{G}_1^{-1}}\widehat{\boldsymbol{f}}_1(\tau)$, $\mathsf{f}(\tau)=\sqrt{\mathsf{G}_2^{-1}}\widehat{\boldsymbol{f}}_2(\tau)$ to the integral terms $\int_{-r_1}^{0}\boldsymbol{x}^\top(t+\tau)R_1\boldsymbol{x}(t+\tau)\mathrm{d}\tau$ and $\int_{-r_2}^{-r_1}\boldsymbol{x}^\top(t+\tau)R_2\boldsymbol{x}(t+\tau)\mathrm{d}\tau$ in (7.54), respectively. Then we have

$$\int_{-r_1}^{0}\boldsymbol{x}^\top(t+\tau)R_1\boldsymbol{x}(t+\tau)\mathrm{d}\tau\geq[*]\left(I_{\kappa_1}\otimes R_1\right)\left[\int_{-r_1}^{0}\left(\sqrt{\mathsf{G}_1^{-1}}\widehat{\boldsymbol{f}}_1(\tau)\otimes I_n\right)\boldsymbol{x}(t+\tau)\mathrm{d}\tau\right] \quad (7.59)$$

$$\int_{-r_2}^{-r_1}\boldsymbol{x}^\top(t+\tau)R_2\boldsymbol{x}(t+\tau)\mathrm{d}\tau\geq[*]\left(\begin{bmatrix}R_2 & Y\\ * & R_2\end{bmatrix}\otimes I_{\kappa_2}\right)\begin{bmatrix}\int_{-r(t)}^{-r_1}\left(I_n\otimes\sqrt{\mathsf{G}_2^{-1}}\widehat{\boldsymbol{f}}_2(\tau)\right)\boldsymbol{x}(t+\tau)\mathrm{d}\tau\\ \int_{-r_2}^{-r(t)}\left(I_n\otimes\sqrt{\mathsf{G}_2^{-1}}\widehat{\boldsymbol{f}}_2(\tau)\right)\boldsymbol{x}(t+\tau)\mathrm{d}\tau\end{bmatrix}$$
$$=[*]\left(\begin{bmatrix}\mathsf{K}_{(\kappa_2,n)} & \mathsf{O}_{\kappa_2 n}\\ * & \mathsf{K}_{(\kappa_2,n)}\end{bmatrix}\left(\begin{bmatrix}R_2 & Y\\ * & R_2\end{bmatrix}\otimes I_{\kappa_2}\right)\begin{bmatrix}\mathsf{K}_{(n,\kappa_2)} & \mathsf{O}_{\kappa_2 n}\\ * & \mathsf{K}_{(n,\kappa_2)}\end{bmatrix}\right)\begin{bmatrix}\int_{-r(t)}^{-r_1}\left(\sqrt{\mathsf{G}_2^{-1}}\widehat{\boldsymbol{f}}_2(\tau)\otimes I_n\right)\boldsymbol{x}(t+\tau)\mathrm{d}\tau\\ \int_{-r_2}^{-r(t)}\left(\sqrt{\mathsf{G}_2^{-1}}\widehat{\boldsymbol{f}}_2(\tau)\otimes I_n\right)\boldsymbol{x}(t+\tau)\mathrm{d}\tau\end{bmatrix}. \quad (7.60)$$

Given the definition of $\mathbb{1}$ and $\widehat{\mathbb{1}}$ in (7.25) and $\widehat{\mathsf{O}}$ in (7.24) for the case of $r_2>r_1>0$, applying (7.59)–(7.60) to (7.54) with (7.44) produces

$$\widetilde{\forall} t\geq t_0,\quad \dot{\mathsf{v}}(\mathbf{x}_t(\cdot))-\mathsf{s}(\boldsymbol{z}(t),\boldsymbol{w}(t))\leq\boldsymbol{\chi}^\top(t)\left(\boldsymbol{\Psi}-\boldsymbol{\Sigma}^\top\widetilde{J}^\top J_1^{-1}\widetilde{J}\boldsymbol{\Sigma}\right)\boldsymbol{\chi}(t) \quad (7.61)$$

where $\boldsymbol{\Psi}$ is given in (7.46) and $\boldsymbol{\chi}(t)$ is given in (7.24). Now it is evident to conclude that if (7.44) and $\boldsymbol{\Psi}-\boldsymbol{\Sigma}^\top\widetilde{J}^\top J_1^{-1}\widetilde{J}\boldsymbol{\Sigma}\prec 0$ are true, then

$$\exists\epsilon_3>0:\ \widetilde{\forall} t\geq t_0,\quad \dot{\mathsf{v}}(\mathbf{x}_t(\cdot))-\mathsf{s}(\boldsymbol{z}(t),\boldsymbol{w}(t))\leq-\epsilon_3\left\|\boldsymbol{x}(t)\right\|_2. \quad (7.62)$$

Moreover, assuming $\boldsymbol{w}(t)\equiv\mathbf{0}_q$, one can also obtain

$$\exists\epsilon_3>0,\ \widetilde{\forall} t\geq t_0,\ \dot{\mathsf{v}}(\mathbf{x}_t(\cdot))\leq-\epsilon_3\left\|\boldsymbol{x}(t)\right\|_2 \quad (7.63)$$

by the structure of $\boldsymbol{\Psi}$ with the fact that $\boldsymbol{\Psi} \prec 0$ and the elements in $\boldsymbol{\chi}(t)$ along with the properties of quadratic forms.. Note that $\mathbf{x}_t(\cdot)$ in (7.63) is in line with the definition of $\mathbf{x}_t(\cdot)$ in (7.40). In consequence, there exists a functional in (7.52) satisfying (7.41) and (7.40) if (7.44) and $\boldsymbol{\Psi} - \boldsymbol{\Sigma}^\top \widetilde{J}^\top J_1^{-1} \widetilde{J} \boldsymbol{\Sigma} \prec 0$ are feasible for some matrices. Finally, applying the Schur complement to $\boldsymbol{\Psi} - \boldsymbol{\Sigma}^\top \widetilde{J}^\top J_1^{-1} \widetilde{J} \boldsymbol{\Sigma} \prec 0$ with (7.44) and $J_1^{-1} \prec 0$ gives the equivalent condition in (7.45). Therefore we have proved that the existence of the feasible solutions to (7.44) and (7.45) guarantees the existence of a functional (7.52) and $\epsilon_3 > 0$ satisfying (7.41) and (7.40).

We now demonstrate that if (7.43) and (7.44) are feasible for some matrices, then there exist $\epsilon_1 > 0$ and $\epsilon_2 > 0$ such that (7.52) satisfies (7.39). Let $\|\boldsymbol{\phi}(\cdot)\|_\infty^2 := \sup_{-r_2 \le \tau \le 0} \|\boldsymbol{\phi}(\tau)\|_2^2$ and consider the structure of (7.52) with $t = t_0$, it follows that there exists $\lambda > 0$ such that

$$
\begin{aligned}
&\mathsf{v}(\mathbf{x}_{t_0}(\cdot)) = \mathsf{v}(\boldsymbol{\phi}(\cdot)) \le \boldsymbol{\eta}^\top(t_0)\lambda\boldsymbol{\eta}(t_0) + \int_{-r_2}^{0} \boldsymbol{\phi}^\top(\tau)\lambda\boldsymbol{\phi}(\tau)\mathsf{d}\tau \le \lambda \|\boldsymbol{\phi}(0)\|_2^2 + \lambda r_2 \|\boldsymbol{\phi}(\cdot)\|_\infty^2 \\
&\quad + \int_{-r_1}^{0} \boldsymbol{\phi}^\top(\tau) \left(\sqrt{\mathsf{F}_1^{-1}}\boldsymbol{f}_1(\tau) \otimes I_n\right)^\top \mathsf{d}\tau \lambda \int_{-r_1}^{0} \left(\sqrt{\mathsf{F}_1^{-1}}\boldsymbol{f}_1(\tau) \otimes I_n\right) \boldsymbol{\phi}(\tau)\mathsf{d}\tau \\
&\quad + \int_{-r_2}^{-r_1} \boldsymbol{\phi}^\top(\tau) \left(\sqrt{\mathsf{F}_2^{-1}}\boldsymbol{f}_2(\tau) \otimes I_n\right)^\top \mathsf{d}\tau \lambda \int_{-r_2}^{-r_1} \left(\sqrt{\mathsf{F}_2^{-1}}\boldsymbol{f}_2(\tau) \otimes I_n\right) \boldsymbol{\phi}(\tau)\mathsf{d}\tau \\
&\le (\lambda + \lambda r_2) \|\boldsymbol{\phi}(\cdot)\|_\infty^2 + \lambda \int_{-r_2}^{0} \boldsymbol{\phi}^\top(\tau)\boldsymbol{\phi}(\tau)\mathsf{d}\tau \le (\lambda + 2\lambda r_2) \|\boldsymbol{\phi}(\cdot)\|_\infty^2
\end{aligned}
\tag{7.64}
$$

for any $\boldsymbol{\phi}(\cdot) \in \mathcal{C}([-r_2, 0]; \mathbb{R}^n)$ in (7.19), where (7.64) is derived via the property of quadratic forms: $\forall X \in \mathbb{S}^n, \exists \lambda > 0 : \forall \mathbf{x} \in \mathbb{R}^n \setminus \{\mathbf{0}\}, \mathbf{x}^\top (\lambda I_n - X)\mathbf{x} > 0$ together with the application of (7.32) with $\varpi(\tau) = 1$ and appropriate $\mathbf{f}(\tau)$. Consequently, the result in (7.64) implies that one can construct an upper bound for (7.52) that satisfies (7.39) with a $\epsilon_2 > 0$.

By applying (7.32) to (7.52) twice with $\varpi(\tau) = 1$ and $\mathbf{f}(\tau) = \sqrt{\mathsf{F}_1^{-1}}\boldsymbol{f}_1(\tau)$, $\mathbf{f}(\tau) = \sqrt{\mathsf{F}_2^{-1}}\boldsymbol{f}_2(\tau)$, we have

$$
\begin{aligned}
&\int_{-r_1}^{0} \boldsymbol{x}^\top(t+\tau)Q_1\boldsymbol{x}(t+\tau)\mathsf{d}\tau \ge [*]\left(I_{d_1} \otimes Q_1\right) \int_{-r_1}^{0} \left(\sqrt{\mathsf{F}_1^{-1}}\boldsymbol{f}_1(\tau) \otimes I_n\right) \boldsymbol{x}(t+\tau)\mathsf{d}\tau \\
&\int_{-r_2}^{-r_1} \boldsymbol{x}^\top(t+\tau)Q_2\boldsymbol{x}(t+\tau)\mathsf{d}\tau \ge [*]\left(I_{d_2} \otimes Q_2\right) \int_{-r_2}^{-r_1} \left(\sqrt{\mathsf{F}_2^{-1}}\boldsymbol{f}_2(\tau) \otimes I_n\right) \boldsymbol{x}(t+\tau)\mathsf{d}\tau
\end{aligned}
\tag{7.65}
$$

provided that (7.44) holds. Moreover, by applying (7.65) to (7.52) with (7.44) and (7.64), it is clear to see that the existence of the feasible solutions to (7.43) and (7.44) imply that (7.52) satisfies (7.39) with some $\epsilon_1; \epsilon_2 > 0$.

In conclusion, we have shown that if the conditions in (7.43)–(7.45) are feasible for some matrices, then there exists a functional (7.52) and $\epsilon_1; \epsilon_2 > 0$ satisfying dissipative condition (7.41) and the stability criteria in (7.39)–(7.40). This implies that the existence of feasible solutions to (7.43)–(7.45) means that the trivial solution to the CLS in (7.19) with $\boldsymbol{w}(t) \equiv \mathbf{0}_q$ is globally uniformly asymptotically stable, and the system (7.19) with (7.42) is dissipative.

Now let us consider the situation where $r_1 = r_2$ and the delay of the system (7.19) is of constant values. It is not difficult to show that the corresponding synthesis condition constructed via the functional in (7.52), following the procedures (7.52)–(7.65) with $r_1 = r_2$, can be obtained by choosing $d_2 = \delta_2 = 0$ in (7.43)–(7.45) with $Q_2 = R_2 = Y = \mathsf{O}_n$. Similarly, the corresponding synthesis condition for $r_1 = 0; r_2 > 0$ can be obtained by choosing $d_1 = \delta_1 = 0$ in (7.43)–(7.45) with $Q_1 = R_1 = \mathsf{O}_n$. Note that the use of $\mathbb{1}$, $\widehat{\mathbb{1}}$ in (7.25) and (7.48), and $\widehat{\mathsf{O}}$ in (7.24) allows (7.43)–(7.45) to cover the corresponding synthesis conditions for the cases of $r_1 = r_2$ and $r_1 = 0$; $r_2 > 0$, without introducing redundant matrices or matrices with ill-posed dimensions. ∎

**Remark 7.7.** Without using $\mathbb{1}$, $\widehat{\mathbb{1}}$ and $\widehat{\mathsf{O}}$, the synthesis condition derived for the case of $r_2 > r_1 > 0$ may not be be directly applicable to cases $r_1 = r_2$ or $r_1 = 0$; $r_2 > 0$. This is due to the changes of the mathematical structures of the CLS in (7.19) and the functional (7.52) corresponding to $r_1 = r_2$ or $r_1 = 0$; $r_2 > 0$.

**Remark 7.8.** Note that $\boldsymbol{f}_1(\cdot)$ and $\boldsymbol{f}_2(\cdot)$ in (7.53) can be any differentiable function since the decompositions in (7.3)–(7.6) are always constructible via proper choices of $\boldsymbol{\varphi}_1(\cdot)$ and $\boldsymbol{\varphi}_2(\cdot)$. This provides great flexibility in the structure of the Lyapunov-Krasovskiĭ functional in (7.52). Additionally, the functions inside of $\boldsymbol{f}_1(\cdot)$ and $\boldsymbol{f}_2(\cdot)$ can be chosen in view of the functions inside of the DDs in (7.1).

### 7.4.1 Some comments on (7.34)

The significance of the proposed inequality in (7.34) can be understood considering the procedures in the proof of Theorem 7.1. Indeed, assume that (7.32) is directly applied to the integrals $\int_{-r(t)}^{-r_1} \boldsymbol{x}^\top(t+\tau)Q_2\boldsymbol{x}(t+\tau)\mathrm{d}\tau$ and $\int_{-r_2}^{-r(t)} \boldsymbol{x}^\top(t+\tau)Q_2\boldsymbol{x}(t+\tau)\mathrm{d}\tau$ without using (7.34) at the step in (7.60), which gives the inequalities

$$
\begin{aligned}
\int_{-r(t)}^{-r_1} \boldsymbol{x}^\top(t+\tau)Q_2\boldsymbol{x}(t+\tau)\mathrm{d}\tau &\geq [*]\left(\widehat{\mathsf{F}}_1^{-1}(r(t))\otimes Q_2\right)\left[\int_{-r(t)}^{-r_1}\left(\widehat{\boldsymbol{f}}_2(\tau)\otimes I_n\right)\boldsymbol{x}(t+\tau)\mathrm{d}\tau\right]\\
\int_{-r_2}^{-r(t)} \boldsymbol{x}^\top(t+\tau)Q_2\boldsymbol{x}(t+\tau)\mathrm{d}\tau &\geq [*]\left(\widehat{\mathsf{F}}_2^{-1}(r(t))\otimes Q_2\right)\left[\int_{-r_2}^{-r(t)}\left(\widehat{\boldsymbol{f}}_2(\tau)\otimes I_n\right)\boldsymbol{x}(t+\tau)\mathrm{d}\tau\right]
\end{aligned}
\tag{7.66}
$$

where $\widehat{\mathsf{F}}_1(r(t)) = \int_{-r(t)}^{-r_1}\widehat{\boldsymbol{f}}_2(\tau)\widehat{\boldsymbol{f}}_2^\top(\tau)\mathrm{d}\tau$ and $\widehat{\mathsf{F}}_2(r(t)) = \int_{-r_2}^{-r(t)}\widehat{\boldsymbol{f}}_2(\tau)\widehat{\boldsymbol{f}}_2^\top(\tau)\mathrm{d}\tau$. Now combine the inequalities in (7.66), we have

$$
\begin{aligned}
\int_{-r_2}^{-r_1} \boldsymbol{x}^\top(t+\tau)Q_2\boldsymbol{x}(t+\tau)\mathrm{d}\tau \geq &\begin{bmatrix}\int_{-r(t)}^{-r_1}\left(\widehat{\boldsymbol{f}}_2(\tau)\otimes I_n\right)\boldsymbol{x}(t+\tau)\mathrm{d}\tau\\ \int_{-r_2}^{-r(t)}\left(\widehat{\boldsymbol{f}}_2(\tau)\otimes I_n\right)\boldsymbol{x}(t+\tau)\mathrm{d}\tau\end{bmatrix}^\top \times\\
&\begin{bmatrix}\widehat{\mathsf{F}}_1^{-1}(r(t))\otimes Q_2 & \mathsf{O}_{d_1n,d_2n}\\ \mathsf{O}_{d_2n,d_1n} & \widehat{\mathsf{F}}_2^{-1}(r(t))\otimes Q_2\end{bmatrix}\begin{bmatrix}\int_{-r(t)}^{-r_1}\left(\widehat{\boldsymbol{f}}_2(\tau)\otimes I_n\right)\boldsymbol{x}(t+\tau)\mathrm{d}\tau\\ \int_{-r_2}^{-r(t)}\left(\widehat{\boldsymbol{f}}_2(\tau)\otimes I_n\right)\boldsymbol{x}(t+\tau)\mathrm{d}\tau\end{bmatrix}
\end{aligned}
\tag{7.67}
$$

which also furnishes a LB for $\int_{-r_2}^{-r_1} \boldsymbol{x}^\top(t+\tau)Q_2\boldsymbol{x}(t+\tau)\mathrm{d}\tau$. Conventionally, the reciprocally convex combination lemma [661] or its derivatives [667, 692, 693] can be applied to a matrix in the form of $\begin{bmatrix}\frac{1}{1-\alpha}X & \mathsf{O}_n\\ \mathsf{O}_n & \frac{1}{\alpha}X\end{bmatrix}$ to construct a tractable LB with finite dimensions. However, the structure of $\begin{bmatrix}\frac{1}{1-\alpha}X & \mathsf{O}_n\\ \mathsf{O}_n & \frac{1}{\alpha}X\end{bmatrix}$ may not always be guaranteed by the matrix

$$
\begin{bmatrix}\mathsf{F}_1^{-1}(r(t))\otimes Q_2 & \mathsf{O}_{d_1n,d_2n}\\ \mathsf{O}_{d_2n,d_1n} & \mathsf{F}_2^{-1}(r(t))\otimes Q_2\end{bmatrix}
\tag{7.68}
$$

in (7.67), since $\mathsf{F}_1^{-1}(r(t))$ and $\mathsf{F}_2^{-1}(r(t))$ are nonlinear with respect to $r(t)$ in general.[5] Additionally, if (7.67) is applied directly to replace the step at (7.59) without the use of any kind of reciprocally convex combination lemmas, then the matrix in (7.68) does appear in the corresponding (7.45), where (7.45) becomes infinite-dimensional and also generally nonlinear with respect to $r(t)$. In contrast, the symmetric matrix in the LB in (7.60) is of finite-dimensional, which is constructed via the application of (7.34). This demonstrates the contribution of the integral inequality (7.34), by which a dissipative synthesis condition of finite dimensions can be derived via the Krasovskiĭ functional method.

[5] If $\widehat{\boldsymbol{f}}_2(\tau)$ only contains Legendre polynomials [468, 22.2.10] with appropriate structures, then the reciprocally convex combination lemma or its derivatives can be applied to (7.68). Notwithstanding, this is a very special case of $\widehat{\boldsymbol{f}}_2(\cdot) \in \mathscr{L}^2\left([-r_2, 0]; \mathbb{R}^{d_2+\delta_2}\right)$ considered in this chapter.

### 7.4.2 A convex dissipative synthesis condition

Similar to the obstacles we have encountered in Chapter 2 and Chapter 3, $\mathsf{Sy}\left[\mathbf{P}^\top\mathbf{\Pi}\right] + \mathbf{\Phi} \prec 0$ in (7.45) is bilinear with respect to the variables in $\mathbf{P}$ and $\mathbf{\Pi}$ if a synthesis problem is concerned, which cannot be solved directly via standard SDP solvers. To address this issue, we construct a convex dissipative synthesis condition in the following theorem by applying the Projection Lemma in Lemma 2.4 to (7.45).

**Theorem 7.2.** *Given $\{\alpha_i\}_{i=1}^{\widehat{3}+\kappa} \subset \mathbb{R}$ and $r_2 > r_1 > 0$ and the functions and parameters in Proposition 1, then the CLS in (7.19) with supply rate function (7.42) is dissipative and the trivial solution $\boldsymbol{x}(t) \equiv \mathbf{0}_n$ of (7.19) with $\boldsymbol{w}(t) \equiv \mathbf{0}_q$ is globally uniformly asymptotically stable if there exists $\acute{P}_1 \in \mathbb{S}^n$, $\acute{P}_2 \in \mathbb{R}^{n\times\varrho}$, $\acute{P}_3 \in \mathbb{S}^\varrho$ and $\acute{Q}_1; \acute{Q}_2; \acute{R}_1; \acute{R}_2; X \in \mathbb{S}^n$ and $\acute{Y} \in \mathbb{R}^{n\times n}$ and $V \in \mathbb{R}^{p\times n}$ such that*

$$\begin{bmatrix} \acute{P}_1 & \acute{P}_2 \\ * & \acute{P}_3 \end{bmatrix} + \left(\mathsf{O}_n \oplus \left[I_{d_1} \otimes \acute{Q}_1\right] \oplus \left[I_{d_2} \otimes \acute{Q}_2\right]\right) \succ 0, \tag{7.69}$$

$$\acute{Q}_1 \succeq 0,\ \ \acute{Q}_2 \succeq 0,\ \ \acute{R}_1 \succeq 0,\ \begin{bmatrix} \acute{R}_2 & \acute{Y} \\ * & \acute{R}_2 \end{bmatrix} \succeq 0, \tag{7.70}$$

$$\mathsf{Sy}\left(\begin{bmatrix} I_n \\ \mathbf{Col}_{i=1}^{\widehat{3}+\kappa}\, \alpha_i I_n \\ \mathsf{O}_{(q+m),n} \end{bmatrix}\begin{bmatrix} -X & \acute{\mathbf{\Pi}} \end{bmatrix}\right) + \begin{bmatrix} \mathsf{O}_n & \acute{\mathbf{P}} \\ * & \acute{\mathbf{\Phi}} \end{bmatrix} \prec 0 \tag{7.71}$$

*where* $\acute{\mathbf{\Pi}} = \left[\mathbf{A}\left[\left(I_{\widehat{3}+\kappa} \otimes X\right) \oplus I_q\right] + \mathbf{B}_1\left[\left(I_{\widehat{3}+\kappa} \otimes V\right) \oplus \mathsf{O}_q\right] \quad \mathsf{O}_{n,m}\right]$ *and* $\acute{\mathbf{P}} = \left[\widehat{\mathsf{O}}_{n,n} \quad \acute{P}_1 \quad \acute{P}_2\widehat{I} \quad \mathsf{O}_{n,(q+m)}\right]$ *with* $\widehat{I}$ *in* (7.47) *and*

$$\begin{aligned} \acute{\mathbf{\Phi}} = \mathsf{Sy}&\left(\begin{bmatrix} \widehat{\mathsf{O}}_{\varrho,n}^\top \\ \acute{P}_2 \\ \widehat{I}^\top\acute{P}_3 \\ \mathsf{O}_{(q+m),\varrho} \end{bmatrix}\begin{bmatrix} \widehat{\mathbf{F}} \otimes I_n & \mathsf{O}_{\varrho,(q+m)} \end{bmatrix} + \begin{bmatrix} \widehat{\mathsf{O}}_{m,n}^\top \\ \mathsf{O}_{(n+\kappa n),m} \\ -J_2^\top \\ \widetilde{J} \end{bmatrix}\begin{bmatrix} \acute{\mathbf{\Sigma}} & \mathsf{O}_m \end{bmatrix}\right) \\ &- \left(\left[\acute{Q}_1 - \acute{Q}_2 - r_3\acute{R}_2\right] \oplus \mathbb{1}\acute{Q}_2 \oplus \left[\widehat{\mathbb{1}}(-\acute{Q}_1 - r_1\acute{R}_1)\right] \oplus \left[I_{\kappa_1} \otimes \acute{R}_1\right]\right. \\ &\left.\oplus\left([*]\left(\begin{bmatrix} \acute{R}_2 & \acute{Y} \\ * & \acute{R}_2 \end{bmatrix} \otimes I_{\kappa_2}\right)\begin{bmatrix} \mathsf{K}_{(n,\kappa_2)} & \mathsf{O}_{\kappa_2 n} \\ * & \mathsf{K}_{(n,\kappa_2)} \end{bmatrix}\right) \oplus J_3 \oplus (-J_1)\right) \end{aligned} \tag{7.72}$$

*with* $\acute{\mathbf{\Sigma}} = \mathbf{C}\left[\left(I_{\widehat{3}+\kappa} \otimes X\right) \oplus I_q\right] + \mathbf{B}_2\left[\left(I_{\widehat{3}+\kappa} \otimes V\right) \oplus \mathsf{O}_q\right]$ *and* $\mathbf{A}$, $\mathbf{B}_1$, $\mathbf{B}_2$, $\mathbf{C}$ *are given in* (7.20)–(7.23). *The controller gain is computed via $K = VX^{-1}$. Furthermore, with $r_1 = r_2$, $d_2 = \delta_2 = 0$ and $\acute{Q}_2 = \acute{R}_2 = \acute{Y} = \mathsf{O}_n$, then the inequalities in* (7.69)–(7.71) *constitute a dissipative synthesis condition for the CLS with $r_1 = r_2 > 0$. Finally, with $r_2 > 0$; $r_1 = 0$, $d_1 = \delta_1 = 0$ and $\acute{Q}_1 = \acute{R}_1 = \mathsf{O}_n$, then the inequalities in* (7.69)–(7.71) *represent a dissipative synthesis condition for the CLS with $r_2 > 0$; $r_1 = 0$.*

*Proof.* Let us consider the case where $r_2 > r_1 > 0$. First of all, Firstly, note that $\mathsf{Sy}\left(\mathbf{P}^\top \mathbf{\Pi}\right) + \mathbf{\Phi} \prec 0$ in (7.45) can be reformulated as

$$\mathsf{Sy}\left(\mathbf{P}^\top \mathbf{\Pi}\right) + \mathbf{\Phi} = \begin{bmatrix} \mathbf{\Pi} \\ I_{3n+\kappa n+q+m} \end{bmatrix}^\top \begin{bmatrix} \mathsf{O}_n & \mathbf{P} \\ * & \mathbf{\Phi} \end{bmatrix} \begin{bmatrix} \mathbf{\Pi} \\ I_{3n+\kappa n+q+m} \end{bmatrix} \prec 0, \tag{7.73}$$

where the structure of (7.73) resembles one of the inequalities in (2.19b) as part of the statements of Lemma 2.4. Given that there are two matrix inequalities in (2.19b), a new matrix inequality must be constructed accordingly to utilize Lemma 2.4 in order to decouple the product between $\mathbf{P}$ and $\mathbf{\Pi}$ in (7.73). Now consider

$$\Upsilon^\top \begin{bmatrix} \mathsf{O}_n & \mathbf{P} \\ * & \mathbf{\Phi} \end{bmatrix} \Upsilon \prec 0 \tag{7.74}$$

with $\Upsilon^\top := \begin{bmatrix} \mathsf{O}_{(q+m),(4n+\kappa n)} & I_{q+m} \end{bmatrix}$. The inequality in (7.74) can be further simplified as

$$\Upsilon^\top \begin{bmatrix} \mathsf{O}_n & \mathbf{P} \\ * & \mathbf{\Phi} \end{bmatrix} \Upsilon = \begin{bmatrix} -J_3 - \mathsf{Sy}(D_2^\top J_2) & D_2^\top \widetilde{J} \\ * & J_1 \end{bmatrix} \prec 0, \tag{7.75}$$

where the left-hand side of the inequality in (7.75) is the $2 \times 2$ block matrix at the right bottom of $\mathsf{Sy}\left(\mathbf{P}^\top \mathbf{\Pi}\right) + \mathbf{\Phi}$ or $\mathbf{\Phi}$. As a result, we see that (7.75) is automatically satisfied if (7.73) or (7.45) holds. Hence, (7.75) and (7.45) hold if and only if (7.45) holds. In addition, the identities

$$\begin{aligned} &\begin{bmatrix} -I_n & \mathbf{\Pi} \end{bmatrix} \begin{bmatrix} \mathbf{\Pi} \\ I_{3n+\kappa n+q+m} \end{bmatrix} = \mathsf{O}_{n,(3n+\kappa n+q+m)}, \quad \begin{bmatrix} -I_n & \mathbf{\Pi} \end{bmatrix}_\perp = \begin{bmatrix} \mathbf{\Pi} \\ I_{3n+\kappa n+q+m} \end{bmatrix} \\ &\begin{bmatrix} I_{4n+\kappa n} & \mathsf{O}_{(4n+\kappa n),(q+m)} \end{bmatrix} \begin{bmatrix} \mathsf{O}_{(4n+\kappa n),(q+m)} \\ I_{q+m} \end{bmatrix} = \begin{bmatrix} I_{4n+\kappa n} & \mathsf{O}_{(4n+\kappa n),(q+m)} \end{bmatrix} \Upsilon = \mathsf{O}_{(4n+\kappa n),(q+m)} \\ &\begin{bmatrix} I_{4n+\kappa n} & \mathsf{O}_{(4n+\kappa n),(q+m)} \end{bmatrix}_\perp = \begin{bmatrix} \mathsf{O}_{(4n+\kappa n),(q+m)} \\ I_{q+m} \end{bmatrix} = \Upsilon, \end{aligned} \tag{7.76}$$

where $\operatorname{rank}\left(\begin{bmatrix} -I_n & \mathbf{\Pi} \end{bmatrix}\right) = n$ and $\operatorname{rank}\left(\begin{bmatrix} I_{4n+\kappa n} & \mathsf{O}_{(4n+\kappa n),(q+m)} \end{bmatrix}\right) = 4n + \kappa n$, which imply that Lemma 2.4 can be used with the terms in (7.76) given the rank nullity theorem.

By applying Lemma 2.4 to (7.73)–(7.75) with (7.76), it follows that (7.73) holds if and only if

$$\exists \mathbf{W} \in \mathbb{R}^{(4n+\kappa n+q+m)\times n} : \mathsf{Sy}\left(\begin{bmatrix} I_{4n+\kappa n} \\ \mathsf{O}_{(q+m),(4n+\kappa n)} \end{bmatrix} \mathbf{W} \begin{bmatrix} -I_n & \mathbf{\Pi} \end{bmatrix}\right) + \begin{bmatrix} \mathsf{O}_n & \mathbf{P} \\ * & \mathbf{\Phi} \end{bmatrix} \prec 0. \tag{7.77}$$

Now (7.77) is still bilinear due to the product between $\mathbf{W}$ and $\mathbf{\Pi}$. To convexify (7.77), we can select

$$\mathbf{W} := \mathbf{Col}\left[W,\ \mathbf{Col}_{i=1}^{\widehat{3}+\kappa} \alpha_i W\right] \tag{7.78}$$

with $W \in \mathbb{S}^n$ and $\{\alpha_i\}_{i=1}^{\widehat{3}+\kappa} \subset \mathbb{R}$. With (7.78) at hand, the inequality in (7.77) becomes

$$\mathbf{\Theta} = \mathsf{Sy}\left(\begin{bmatrix} W \\ \mathbf{Col}_{i=1}^{3+\kappa} \alpha_i W \\ \mathsf{O}_{(q+m),n} \end{bmatrix} \begin{bmatrix} -I_n & \mathbf{\Pi} \end{bmatrix}\right) + \begin{bmatrix} \mathsf{O}_n & \mathbf{P} \\ * & \mathbf{\Phi} \end{bmatrix} \prec 0 \tag{7.79}$$

which further implies (7.73). Note that using the structure in (7.78) suggests that (7.79) is no longer an equivalent but only a sufficient condition for (7.73) which is equivalent to (7.45). It is also important to stress that an invertible $W$ is automatically implied by (7.79) since the expression $-W - W^\top = -2W$ is the sole element in the first top left diagonal block of $\mathbf{\Theta}$.

Let $X^\top = W^{-1}$, and apply congruent transformations [578, page 12] to the matrix inequalities in (7.43),(7.44) and (7.79) taking into account that an invertible $W$ is implied by (7.79). Then, we can infer that

$$
\begin{gathered}
X^\top Q_1 X \succ 0, \quad X^\top Q_2 X \succ 0, \quad X^\top R_1 X \succ 0, \quad \begin{bmatrix} X^\top & \mathsf{O}_n \\ * & X^\top \end{bmatrix} \begin{bmatrix} R_2 & Y \\ * & R_2 \end{bmatrix} \begin{bmatrix} X & \mathsf{O}_n \\ * & X \end{bmatrix} \succ 0, \\
\left[ \left( I_{4+\kappa} \otimes X^\top \right) \oplus I_{q+m} \right] \mathbf{\Theta} \left[ (I_{4+\kappa} \otimes X) \oplus I_{q+m} \right] \prec 0, \\
[*] \left( \begin{bmatrix} P_1 & P_2 \\ * & P_3 \end{bmatrix} + \left( \mathsf{O}_n \oplus [I_{d_1} \otimes Q_1] \oplus [I_{d_2} \otimes Q_2] \right) \right) (I_{1+d_1+d_2} \otimes X) \succ 0
\end{gathered} \tag{7.80}
$$

hold if and only if (7.43),(7.44) and (7.79) hold. Moreover, considering (2.1a) and the definitions $\acute{Y} := X^\top Y X$ and

$$
\begin{bmatrix} \acute{P}_1 & \acute{P}_2 \\ * & \acute{P}_3 \end{bmatrix} := [*] \begin{bmatrix} P_1 & P_2 \\ * & P_3 \end{bmatrix} (I_{1+d_1+d_2} \otimes X), \quad \begin{bmatrix} \acute{Q}_1 & \acute{Q}_2 & \acute{R}_1 & \acute{R}_2 \end{bmatrix} := X^\top \begin{bmatrix} Q_1 X & Q_2 X & R_1 X & R_2 X \end{bmatrix}, \tag{7.81}
$$

then the inequalities in (7.80) can be rewritten as (7.69) and (7.70) and

$$
[*] \mathbf{\Theta} \left[ (I_{4+\kappa} \otimes X) \oplus I_{q+m} \right] = \acute{\mathbf{\Theta}} = \mathsf{Sy} \left( \begin{bmatrix} I_n \\ \mathbf{Col}_{i=1}^{3+\kappa} \alpha_i I_n \\ \mathsf{O}_{(q+m),n} \end{bmatrix} \begin{bmatrix} -X & \acute{\mathbf{\Pi}} \end{bmatrix} \right) + \begin{bmatrix} \mathsf{O}_n & \acute{\mathbf{P}} \\ * & \acute{\mathbf{\Phi}} \end{bmatrix} \prec 0, \tag{7.82}
$$

where $\acute{\mathbf{P}} = X\mathbf{P} \left[ (I_{3+\kappa} \otimes X) \oplus I_{q+m} \right] = \begin{bmatrix} \widehat{\mathsf{O}}_{n,n} & \acute{P}_1 & \acute{P}_2 \widehat{I} & \mathsf{O}_{n,q} & \mathsf{O}_{n,m} \end{bmatrix}$ and

$$
\begin{aligned}
\acute{\mathbf{\Pi}} = \mathbf{\Pi} \left[ (I_{3+\kappa} \otimes X) \oplus I_{q+m} \right] &= \begin{bmatrix} \mathbf{A} \left[ (I_{3+\kappa} \otimes X) \oplus I_q \right] + \mathbf{B}_1 \left[ (I_{3+\kappa} \otimes KX) \oplus \mathsf{O}_q \right] & \mathsf{O}_{n,m} \end{bmatrix} \\
&= \begin{bmatrix} \mathbf{A} \left[ (I_{3+\kappa} \otimes X) \oplus I_q \right] + \mathbf{B}_1 \left[ (I_{3+\kappa} \otimes V) \oplus \mathsf{O}_q \right] & \mathsf{O}_{n,m} \end{bmatrix}
\end{aligned} \tag{7.83}
$$

with $V = KX$ and $\acute{\mathbf{\Phi}}$ in (7.72). Note that (7.82)–(7.83) is equivalent to the statements in Theorem 7.2 given the definition of $\widehat{3}$ and $\widehat{\mathsf{O}}$ in (7.24). Note that also the form of $\acute{\mathbf{\Phi}}$ in (7.72) is derived via the relations $\widehat{I} (I_\kappa \otimes X) = (I_{d_1+d_2} \otimes X) \widehat{I}$ and

$$
\begin{aligned}
&\begin{bmatrix} \widehat{\mathbf{F}} \otimes I_n & \mathsf{O}_{\varrho,(q+m)} \end{bmatrix} \left[ (I_{3+\kappa} \otimes X) \oplus I_{q+m} \right] = \begin{bmatrix} I_{d_1+d_2} \widehat{\mathbf{F}} \otimes X I_n & \mathsf{O}_{\varrho,(q+m)} \end{bmatrix} \\
&\quad = \begin{bmatrix} (I_{d_1+d_2} \otimes X) \left( \widehat{\mathbf{F}} \otimes I_n \right) & \mathsf{O}_{\varrho,(q+m)} \end{bmatrix} = (I_{d_1+d_2} \otimes X) \begin{bmatrix} \widehat{\mathbf{F}} \otimes I_n & \mathsf{O}_{\varrho,(q+m)} \end{bmatrix},
\end{aligned} \tag{7.84}
$$

$$
\begin{aligned}
\begin{bmatrix} \mathsf{K}_{(n,\kappa_2)} & \mathsf{O}_{\kappa_2 n} \\ * & \mathsf{K}_{(n,\kappa_2)} \end{bmatrix} \begin{bmatrix} I_{\kappa_2} \otimes X & \mathsf{O}_{\kappa_2 n} \\ * & I_{\kappa_2} \otimes X \end{bmatrix} &= \begin{bmatrix} X \otimes I_{\kappa_2} & \mathsf{O}_{\kappa_2 n} \\ * & X \otimes I_{\kappa_2} \end{bmatrix} \begin{bmatrix} \mathsf{K}_{(n,\kappa_2)} & \mathsf{O}_{\kappa_2 n} \\ * & \mathsf{K}_{(n,\kappa_2)} \end{bmatrix} \\
&= \left( \begin{bmatrix} X & \mathsf{O}_n \\ * & X \end{bmatrix} \otimes I_{\kappa_2} \right) \begin{bmatrix} \mathsf{K}_{(n,\kappa_2)} & \mathsf{O}_{\kappa_2 n} \\ * & \mathsf{K}_{(n,\kappa_2)} \end{bmatrix}
\end{aligned} \tag{7.85}
$$

which are established from the properties of matrices with (2.1a),(2.1b) and (7.29). Furthermore, we see that $X$ is invertible if (7.71) holds since $-X - X^\top = -2X$ is the sole element in the first top left diagonal block of $\acute{\mathbf{\Theta}}$ in (7.71). This is consistent with the fact that an invertible $W$ is implied by the matrix inequality in (7.79).

Thus, we have shown the equivalence between (7.43)–(7.44) and (7.69)–(7.70) for the case of $r_2 > r_1 > 0$. Meanwhile, it has been shown that (7.71) is equivalent to (7.79) which implies (7.45). In consequence, (7.43)–(7.45) are satisfied if (7.69)–(7.71) hold with some $W \in \mathbb{S}^n$ and $\{\alpha_i\}_{i=1}^{\widehat{3}+\kappa} \subset \mathbb{R}$. Thus it demonstrates that the existence of the feasible solutions to (7.69)–(7.71) ensures that the trivial solution $\boldsymbol{x}(t) \equiv \mathbf{0}_n$ of the CLS in (7.19) with $\boldsymbol{w}(t) \equiv \mathbf{0}_q$ is globally uniformly asymptotically stable and (7.19) with (7.42) is dissipative.

Now for the case of $r_1 = r_2$, it is not difficult to show that a synthesis condition can be obtained by setting $d_2 = \delta_2 = 0$ in (7.69)–(7.71) with $\acute{Q}_2 = \acute{R}_2 = \acute{Y} = \mathsf{O}_n$ and $r_1 = r_2$, given the definition of $\widehat{3}$ and $\widehat{\mathsf{O}}$ in (7.24). The proof of such a synthesis condition for $r_1 = r_2$ follows the same procedures presented above, with substitutions $3 \leftarrow \widehat{3}$ and $4 \leftarrow \widehat{3}+1$ and $d_2 = \delta_2 = 0$ in (7.73)–(7.85). Similarly, a synthesis condition for the case of $r_1 = 0$; $r_2 > 0$ can be obtained by letting $d_1 = \delta_1 = 0$ in (7.69)–(7.71) with substitutions $3 \leftarrow \widehat{3}$ and $4 \leftarrow \widehat{3}+1$ and $\acute{Q}_1 = \acute{R}_1 = \mathsf{O}_n$ and $r_1 = 0$; $r_2 > 0$. ■

**Remark 7.9.** Please note that Theorem 7.2 has been specifically derived to address the synthesis problem of (7.19). In the case of an open-loop system where $B_1 = \widetilde{B}_2(\tau) = \mathsf{O}_{n,p}$ and $B_4 = B_5(\tau) = \mathsf{O}_{m,p}$, then Theorem 7.1 should be applied instead of Theorem 7.2. This is because Theorem 7.2 is more conservative compared to Theorem 7.1 for a specific problem of stability analysis.

**Remark 7.10.** For $\{\alpha_i\}_{i=1}^{\widehat{3}+\kappa} \subset \mathbb{R}$ in (7.71), some values of $\alpha_i$ may have a more significant impact on the feasibility of (7.71). For example, the value of $\alpha_{\widehat{3}}$ may have a significant impact on the feasibility of (7.71) since it may determine the feasibility of the very diagonal block related to $A_1$ in (7.71). A simple assignment for $\{\alpha_i\}_{i=1}^{\widehat{3}+\kappa} \subset \mathbb{R}$ can be $\alpha_i = 0$ for $i = 1, \ldots, \widehat{3}+\kappa$ with $i \neq \widehat{3}$ which allows one to only adjust the value of $\alpha_{\widehat{3}}$ to use Theorem 7.2.

### 7.4.3 An inner convex approximation solution to Theorem 7.1

For a dissipative synthesis problem, Theorem 7.2 offers a convex solution. Nevertheless, the simplification in (7.78) can render Theorem 7.2 more conservative than Theorem 7.1, while the BMI in Theorem 7.1 cannot be solved by standard SDP solvers. Following the strategies outlined in Chapter 2 and Chapter 3, an iterative algorithm is proposed in this subsection based on the method proposed in [565]. This algorithm provides an inner convex approximation solution to the BMI in (7.45), and can be initiated by a feasible solution to Theorem 7.2. Thus the advantage of both Theorem 7.1 and Theorem 7.2 are combined together in the proposed algorithm without the need to solve nonlinear optimization constraints.

First of all, it is important to note that (7.43) and (7.44) remain convex even when considering a synthesis problem. Now we can see that (7.45) can be rewritten as

$$\mathfrak{U}(\mathbf{H}, K) := \mathsf{Sy}\left[\mathbf{P}^\top \mathbf{\Pi}\right] + \mathbf{\Phi} = \mathsf{Sy}\left(\mathbf{P}^\top \mathbf{B}\left[\left(I_{\widehat{3}+\kappa} \otimes K\right) \oplus \mathsf{O}_{p+m}\right]\right) + \widehat{\mathbf{\Phi}} \prec 0 \tag{7.86}$$

with $\mathbf{B} := \begin{bmatrix} \mathbf{B}_1 & \mathsf{O}_{n,m} \end{bmatrix}$ and $\widehat{\mathbf{\Phi}} := \mathsf{Sy}\left(\mathbf{P}^\top \begin{bmatrix} \mathbf{A} & \mathsf{O}_{n,m} \end{bmatrix}\right) + \mathbf{\Phi}$, where $\mathbf{P}$ is given in (7.50), and $\mathbf{A}$ and $\mathbf{B}_1$ are given in (7.20)–(7.21), and $\mathbf{H} := \begin{bmatrix} P_1 & P_2 \end{bmatrix}$ with $P_1$ and $P_2$ in Theorem 7.1. It is worth mentioning $\widehat{\mathbf{\Phi}}$ is convex with respect to all the decision variables it contains. Given the conclusions of Example 3 in [565], one can conclude that the function $\Delta\left(\bullet, \widetilde{\mathbf{G}}, \bullet, \widetilde{\mathbf{\Gamma}}\right)$, which is defined as

$$\Delta\left(\mathbf{G}, \widetilde{\mathbf{G}}, \mathbf{\Gamma}, \widetilde{\mathbf{\Gamma}}\right) := \begin{bmatrix} \mathbf{G}^\top - \widetilde{\mathbf{G}}^\top & \mathbf{\Gamma}^\top - \widetilde{\mathbf{\Gamma}}^\top \end{bmatrix} \left[Z \oplus (I_n - Z)\right]^{-1} \begin{bmatrix} \mathbf{G} - \widetilde{\mathbf{G}} \\ \mathbf{\Gamma} - \widetilde{\mathbf{\Gamma}} \end{bmatrix}$$

$$+\,\mathsf{Sy}\left(\widetilde{\mathbf{G}}^{\top}\mathbf{\Gamma}+\mathbf{G}^{\top}\widetilde{\mathbf{\Gamma}}-\widetilde{\mathbf{G}}^{\top}\widetilde{\mathbf{\Gamma}}\right)+\mathbf{T} \tag{7.87}$$

with $Z \oplus (I_n - Z) \succ 0$ satisfying

$$\forall \mathbf{G};\widetilde{\mathbf{G}} \in \mathbb{R}^{n\times l},\ \forall \mathbf{\Gamma};\widetilde{\mathbf{\Gamma}} \in \mathbb{R}^{n\times l},\ \ \mathbf{T}+\mathsf{Sy}\left(\mathbf{G}^{\top}\mathbf{\Gamma}\right) \preceq \Delta\left(\mathbf{G},\widetilde{\mathbf{G}},\mathbf{\Gamma},\widetilde{\mathbf{\Gamma}}\right),\ \ \mathbf{T}+\mathsf{Sy}\left(\mathbf{G}^{\top}\mathbf{\Gamma}\right) = \Delta(\mathbf{G},\mathbf{G},\mathbf{\Gamma},\mathbf{\Gamma}), \tag{7.88}$$

is a PSD-convex overestimate [565] of $\acute{\Delta}(\mathbf{G},\mathbf{\Gamma}) = \mathbf{T}+\mathsf{Sy}\left[\mathbf{G}^{\top}\mathbf{\Gamma}\right]$ with respect to the parameterization

$$\begin{bmatrix}\mathbf{vec}(\widetilde{\mathbf{G}})\\ \mathbf{vec}(\widetilde{\mathbf{\Gamma}})\end{bmatrix} = \begin{bmatrix}\mathbf{vec}(\mathbf{G})\\ \mathbf{vec}(\mathbf{\Gamma})\end{bmatrix}. \tag{7.89}$$

Let

$$\begin{aligned}
&\mathbf{T}=\widehat{\mathbf{\Phi}},\ \ \mathbf{G}=\mathbf{P}=\begin{bmatrix}\widehat{\mathsf{O}}_{n,n} & P_1 & P_2\widehat{I} & \mathsf{O}_{n,q} & \mathsf{O}_{n,m}\end{bmatrix},\\
&\widetilde{\mathbf{G}}=\widetilde{\mathbf{P}}=\begin{bmatrix}\widehat{\mathsf{O}}_{n,n} & \widetilde{P}_1 & \widetilde{P}_2\widehat{I} & \mathsf{O}_{n,q} & \mathsf{O}_{n,m}\end{bmatrix},\\
&\mathbf{H}=\begin{bmatrix}P_1 & P_2\end{bmatrix},\ \widetilde{\mathbf{H}}:=\begin{bmatrix}\widetilde{P}_1 & \widetilde{P}_2\end{bmatrix},\ \widetilde{P}_1\in\mathbb{S}^n,\ \widetilde{P}_2\in\mathbb{R}^{n\times dn},\\
&\mathbf{\Gamma}=\mathbf{B}\mathbf{K},\ \ \mathbf{K}=\left[\left(I_{\widehat{3}+\kappa}\otimes K\right)\oplus \mathsf{O}_{p+m}\right],\ \ \ \widetilde{\mathbf{\Gamma}}=\mathbf{B}\widetilde{\mathbf{K}},\ \ \widetilde{\mathbf{K}}=\left[\left(I_{\widehat{3}+\kappa}\otimes \widetilde{K}\right)\oplus \mathsf{O}_{p+m}\right]
\end{aligned} \tag{7.90}$$

in (7.87) with $l=\widehat{3}n+\kappa n+q+m$ and $Z\oplus(I_n-Z)\succ 0$ and $\widehat{\mathbf{\Phi}}$, $\mathbf{H}$ and $K$ in line with the definition in (7.86), one can conclude that

$$\begin{aligned}
\mathfrak{U}(\mathbf{H},K) &= \widehat{\mathbf{\Phi}}+\mathsf{Sy}\left[\mathbf{P}^{\top}\mathbf{B}\left[\left(I_{\widehat{3}+\kappa}\otimes K\right)\oplus \mathsf{O}_{p+m}\right]\right] \preceq \mathfrak{S}\left(\mathbf{H},\widetilde{\mathbf{H}},K,\widetilde{K}\right)\\
&:= \widehat{\mathbf{\Phi}}+\mathsf{Sy}\left(\widetilde{\mathbf{P}}^{\top}\mathbf{B}\mathbf{K}+\mathbf{P}^{\top}\mathbf{B}\widetilde{\mathbf{K}}-\widetilde{\mathbf{P}}^{\top}\mathbf{B}\widetilde{\mathbf{K}}\right)+\begin{bmatrix}\mathbf{P}^{\top}-\widetilde{\mathbf{P}}^{\top} & \mathbf{K}^{\top}\mathbf{B}^{\top}-\widetilde{\mathbf{K}}^{\top}\mathbf{B}^{\top}\end{bmatrix}\left[Z\oplus(I_n-Z)\right]^{-1}[*]
\end{aligned} \tag{7.91}$$

by (7.88), where $\mathfrak{S}(\bullet,\widetilde{\mathbf{H}},\bullet,\widetilde{K})$ in (7.91) is a PSD-convex overestimate of $\mathfrak{U}(\mathbf{H},K)$ in (7.86) with respect to the parameterization

$$\begin{bmatrix}\mathbf{vec}(\widetilde{\mathbf{H}})\\ \mathbf{vec}(\widetilde{K})\end{bmatrix} = \begin{bmatrix}\mathbf{vec}(\mathbf{H})\\ \mathbf{vec}(K)\end{bmatrix}. \tag{7.92}$$

From (7.91), it is obvious that one can infer (7.86) from $\mathfrak{S}\left(\mathbf{H},\widetilde{\mathbf{H}},K,\widetilde{K}\right)\prec 0$. Moreover, it is also true that $\mathfrak{S}\left(\mathbf{H},\widetilde{\mathbf{H}},K,\widetilde{K}\right)\prec 0$ in (7.91) holds if and only if

$$\begin{bmatrix}
\widehat{\mathbf{\Phi}}+\mathsf{Sy}\left(\widetilde{\mathbf{P}}^{\top}\mathbf{B}\mathbf{K}+\mathbf{P}^{\top}\mathbf{B}\widetilde{\mathbf{K}}-\widetilde{\mathbf{P}}^{\top}\mathbf{B}\widetilde{\mathbf{K}}\right) & \mathbf{P}^{\top}-\widetilde{\mathbf{P}}^{\top} & \mathbf{K}^{\top}\mathbf{B}^{\top}-\widetilde{\mathbf{K}}^{\top}\mathbf{B}^{\top}\\
* & -Z & \mathsf{O}_n\\
* & * & Z-I_n
\end{bmatrix}\prec 0 \tag{7.93}$$

holds based on the application of the Schur complement given that $Z\oplus(I_n-Z)\succ 0$. Now we can infer (7.86) from (7.93) that can be solved by standard numerical solvers of SDP provided that the values of $\widetilde{\mathbf{H}}$ and $\widetilde{K}$ are known.

By compiling all the aforementioned procedures according to the expositions in [565], an iterative algorithm is constructed in Algorithm 4, where $\mathbf{x}$ comprises all the variables in $P_3$, $Q_1$, $Q_2$ $R_1$, $R_2$, $Y$ in Theorem 7.1 and $Z$ in (7.93). Furthermore, $\mathbf{H}$, $\widetilde{\mathbf{H}}$, $\mathbf{K}$ and $\widetilde{\mathbf{K}}$ in Algorithm 4 are defined in (7.90) and $\rho_1$, $\rho_2$ and $\varepsilon$ are given constants for regularizations and setting up error tolerance, respectively.

Based on the mathematic constructs in [565], certain initial data for $\widetilde{\mathbf{H}}$ and $\widetilde{K}$ must be obtained to initialize Algorithm 4. This data can be part of a feasible solution to (7.43)–(7.45) in Theorem 7.1. For this reason, $\widetilde{P}_1 \leftarrow P_1$, $\widetilde{P}_2 \leftarrow P_2$ and $\widetilde{K} \leftarrow K$ is used for the initial data of $\widetilde{\mathbf{H}}$ and $\widetilde{K}$ in Algorithm 4 if $P_1$, $P_2$ and $K$ are feasible solutions to (7.43)–(7.45). Generally speaking, acquiring a feasible solution to

**Algorithm 4:** An inner convex approximation solution to Theorem 7.1 with $r_2 > r_1 > 0$

**begin**

 **solve** Theorem 7.2 with given $\alpha_i$ to obtain a feasible $K$, and then **solve** Theorem 7.1 with the previous $K$ to obtain $\mathbf{H} = \begin{bmatrix} P_1 & P_2 \end{bmatrix}$.

 **update** $\widetilde{\mathbf{H}} \longleftarrow \mathbf{H}, \quad \widetilde{K} \longleftarrow K$,

 **solve** $\min_{\mathbf{x},\mathbf{H},K} \operatorname{tr}\left[\rho_1[*](\mathbf{H}-\widetilde{\mathbf{H}})\right] + \operatorname{tr}\left[\rho_2[*](K-\widetilde{K})\right]$ subject to (7.43)–(7.44) and (7.93) to obtain $\mathbf{H}$ and $K$

 **while** $\dfrac{\left\| \begin{bmatrix} \mathbf{vec}(\mathbf{H}) \\ \mathbf{vec}\,(K) \end{bmatrix} - \begin{bmatrix} \mathbf{vec}(\widetilde{\mathbf{H}}) \\ \mathbf{vec}(\widetilde{K}) \end{bmatrix} \right\|_\infty}{\left\| \begin{bmatrix} \mathbf{vec}(\widetilde{\mathbf{H}}) \\ \mathbf{vec}(\widetilde{K}) \end{bmatrix} \right\|_\infty + 1} \geq \varepsilon$ **do**

  **update** $\widetilde{\mathbf{H}} \longleftarrow \mathbf{H}, \quad \widetilde{K} \longleftarrow K$;

  **solve** $\min_{\mathbf{x},\mathbf{H},K} \operatorname{tr}\left[\rho_1[*](\mathbf{H}-\widetilde{\mathbf{H}})\right] + \operatorname{tr}\left[\rho_2[*](K-\widetilde{K})\right]$ subject to (7.43)–(7.44) and (7.93) to obtain $\mathbf{H}$ and $K$;

 **end**

**end**

Theorem 7.1 may not be a trivial task. Nevertheless, as proposed in Theorem 7.2, initial values of $\widetilde{P}_1$, $\widetilde{P}_2$ and $\widetilde{K}$ can be obtained by solving the constraints in (7.69)–(7.71) with given values[6] of $\{\alpha_i\}_{i=1}^{\widehat{3}+\kappa}$.

> **Remark 7.11.** If a convex objective function is considered in Theorem 7.1, for instance $\mathcal{L}^2$ gain $\gamma > 0$ minimization, a termination criterion [565] can be added to Algorithm 4 to track the progress of the objective function between each adjacent iteration.

> **Remark 7.12.** For the delay values $r_2 > 0$; $r_1 = 0$ or $r_2 = r_1 > 0$, Algorithm 4 can be utilized via the corresponding synthesis conditions with appropriate parameter assignments as stated in the statements of Theorem 7.1 and Theorem 7.2.

Since we have proposed plenty of technical results in this chapter, a summary concerning their relations is presented as follows:

- The first noteworthy result is the decomposition scenario in Proposition 1. This enables us to denote general DDs in terms of the products between constants and some appropriate functions.

- By using Proposition 1, one can derive the synthesis results in Theorem 7.1, where the condition is expressed as finite-dimensional optimization constraints through the application of the integral inequality proposed in (7.34).

- Theorem 7.2 has been proposed as a convexification of the BMI in Theorem 7.1 via the application of the Projection Lemma.

[6]Note that as demonstrated in Remark 7.10 that one may apply Theorem 7.2 with $\alpha_i = 0$ for $i = 1, \ldots, \widehat{3}+\kappa, i \neq \widehat{3}$ which allow users to only adjust the value of $\alpha_{\widehat{3}}$ to solve the conditions in Theorem 7.2

- Algorithm 4 has been further proposed to solve the BMI in Theorem 7.1 based on an inner convex approximation algorithm. The initial value of Algorithm 4 can be provided by solving the synthesis condition in Theorem 7.2.

## 7.5 Numerical Examples

We investigate two numerical examples for the demonstration of the effectiveness of our proposed methodologies. All numerical programs are coded and computed in MATLAB using Yalmip [603] as the optimization interface. Moreover, we use SDPT3 [512] as the numerical solver of SDP.

### 7.5.1 Stability and dissipative analysis of a linear system with a time-varying distributed delay

Consider a system in the form of (7.1) with any $r(\cdot) \in \mathcal{M}(\mathbb{R}; [r_1, r_2])$ and the state space matrices

$$
\begin{aligned}
&A_1 = \begin{bmatrix} 0.1 & 0 \\ 0 & -1 \end{bmatrix}, \ \widetilde{A}_2(\tau) = \begin{bmatrix} 0.3e^{\cos(5\tau)} - 0.1e^{\sin(5\tau)} - 0.4 & 0.01e^{\cos(5\tau)} - 0.1e^{\sin(5\tau)} + 1 \\ \ln(2-\tau) - 1 & 0.4 - 0.3e^{\cos(5\tau)} \end{bmatrix}, \\
&B_1 = \widetilde{B}_2(\tau) = B_4 = \widetilde{B}_5(\tau) = \begin{bmatrix} 0 \\ 0 \end{bmatrix}, D_1 = \begin{bmatrix} 0.1 \\ 0.2 \end{bmatrix}, C_1 = \begin{bmatrix} -0.1 & 0.2 \\ 0 & 0.1 \end{bmatrix}, \\
&\widetilde{C}_2(\tau) = \begin{bmatrix} 0.2e^{\sin(5\tau)} - 0.11 & 0.1 - 0.5\ln(2-\tau) \\ 0.1e^{\sin(5\tau)} & 0.14e^{\cos(5\tau)} - 0.2e^{\sin(5\tau)} \end{bmatrix}, \ D_2 = \begin{bmatrix} 0.12 \\ 0.1 \end{bmatrix}.
\end{aligned}
\tag{7.94}
$$

Moreover, let

$$
J_1 = -\gamma I_m, \ \ \widetilde{J} = I_m, \ \ J_2 = \mathsf{O}_{m,q}, \ \ J_3 = \gamma I_q \tag{7.95}
$$

for the supply rate function (7.42) where the objective is to calculate the minimum value of $\mathcal{L}^2$ gain $\gamma$. Note that all the controller gains in (7.94) are zeros, and the DDs in (7.94) contain different types of functions.

To the best of our knowledge, there are no existing approaches, whether based on time or frequency-domain methods, capable of analyzing the stability of (7.1) with the parameters in (7.94). Note that since $r(t)$ is time-varying and its expression is unknown, the DD kernels in (7.94) may not be approximated over $[-r(t), 0]$ via the approaches in [463, 470]. Similarly, the DDs may not be easily analyzed in the frequency domain using existing analytical methods such as those in [203, 215, 268]. Furthermore, no existing methods are capable of calculating the $\mathcal{L}^2$ gain of the system considered in this subsection.

By observing the functions inside of $\widetilde{A}_2(\cdot), \widetilde{C}_2(\cdot)$ in (7.94), we choose

$$
\begin{aligned}
&\boldsymbol{f}_1(\tau) = \boldsymbol{f}_2(\tau) = \begin{bmatrix} 1 \\ e^{\sin(5\tau)} \\ e^{\cos(5\tau)} \\ \ln(2-\tau) \end{bmatrix}, \ \ \boldsymbol{\varphi}_1(\tau) = \boldsymbol{\varphi}_2(\tau) = \begin{bmatrix} \cos(5\tau)e^{\sin(5\tau)} \\ \sin(5\tau)e^{\cos(5\tau)} \\ \dfrac{1}{\tau - 2} \end{bmatrix}, \\
&M_1 = M_2 = \begin{bmatrix} 0 & 0 & 0 & 0 & 0 & 0 & 0 \\ 5 & 0 & 0 & 0 & 0 & 0 & 0 \\ 0 & -5 & 0 & 0 & 0 & 0 & 0 \\ 0 & 0 & 1 & 0 & 0 & 0 & 0 \end{bmatrix}
\end{aligned}
\tag{7.96}
$$

for the functions $\boldsymbol{f}_1(\cdot), \boldsymbol{f}_2(\cdot)$ and $\boldsymbol{\varphi}_1(\cdot), \boldsymbol{\varphi}_2(\cdot)$ in Proposition 1, which corresponds to $d_1 = d_2 = 4$, $\delta_1 = \delta_2 = 3$, $n = m = 2$, $q = 1$, and

$$\begin{aligned} A_2 = A_3 &= \begin{bmatrix} 0 & 0 & 0 & 0 & 0 & 0 & -0.4 & 1 & -0.1 & -0.1 & 0.3 & 0.01 & 0 & 0 \\ 0 & 0 & 0 & 0 & 0 & 0 & -1 & 0.4 & 0 & 0 & 0 & -0.3 & 1 & 0 \end{bmatrix}, \quad B_2 = B_3 = \mathsf{O}_{2\times 7} \\ C_2 = C_3 &= \begin{bmatrix} 0 & 0 & 0 & 0 & 0 & 0 & -0.11 & 0.1 & 0.2 & 0 & 0 & 0 & 0 & -0.5 \\ 0 & 0 & 0 & 0 & 0 & 0 & 0 & 0 & 0.1 & -0.2 & 0 & 0.14 & 0 & 0 \end{bmatrix}, \quad B_5 = B_6 = \mathsf{O}_{2\times 7}. \end{aligned} \tag{7.97}$$

Now apply Theorem 7.1 to (7.19) with the parameters in (7.94)–(7.97), where the conditions in Theorem 7.1 are all convex in this case. This produces results in Table 7.1–7.2, where several detectable delay boundaries are presented with the corresponding $\min \gamma$.

| $[r_1, r_2]$ | $[0.98, 1.25]$ | $[1, 1.23]$ | $[1.02, 1.21]$ | $[1.04, 1.19]$ |
|---|---|---|---|---|
| $r_3 = r_2 - r_1$ | 0.27 | 0.23 | 0.19 | 0.15 |
| $\min \gamma$ | 0.5511 | 0.51356 | 0.48277 | 0.45692 |

**Table 7.1:** $\min \gamma$ produced with decreasing values of $r_3$

| $[r_1, r_2]$ | $[0.8, 1.07]$ | $[1, 1.27]$ | $[1.2, 1.47]$ | $[1.32, 1.59]$ |
|---|---|---|---|---|
| $r_3 = r_2 - r_1$ | 0.27 | 0.27 | 0.27 | 0.27 |
| $\min \gamma$ | 0.35556 | 0.59179 | 1.7935 | 25.9774 |

**Table 7.2:** $\min \gamma$ produced with a fixed value for $r_3$

The results of $\min \gamma$ in Table 7.1 indicate that smaller $r_3$ can lead to smaller $\min \gamma$ values. Indeed, it is more challenging to make the system dissipative for $r(\cdot) \in \mathcal{M}(\mathbb{R}; [r_1, r_2])$ with a large value of $r_3$ than for all $r(\cdot) \in \mathcal{M}(\mathbb{R}; [\acute{r}_1, \acute{r}_2])$ with a smaller value of $\acute{r}_3 = \acute{r}_2 - \acute{r}_1$ if $[\acute{r}_1, \acute{r}_2] \subset [r_1, r_2]$. In addition, the values of $\min \gamma$ in Table 7.2 show that the values of $r_1$ and $r_2$ can significantly affect the resulting $\min \gamma$ even with a fixed $r_3 = r_2 - r_1$.

In order to partially verify the results in Table 7.1 and Table 7.2, we apply the frequency domain method in [203] to (7.94) assuming that $r(\cdot) \in \mathcal{M}(\mathbb{R}; [r_1, r_2])$ is an unknown function with a constant value. (Note that an unknown $r(\cdot)$ with a constant value is an option for $r(\cdot) \in \mathcal{M}(\mathbb{R}; [r_1, r_2])$) The result shows that the system with a constant value of $r$ is stable over $[0.61, 1.64]$, which is consistent with the results in Table 7.1 and Table 7.2. This is because the results in Table 7.1 and Table 7.2 imply that the system with a constant delay value is stable over the intervals therein, which are all the subsets of $[0.61, 1.64]$.

**Remark 7.13.** It should be noted that the values of $\min \gamma$ in Table 7.1–7.2 are valid for any $r(\cdot) \in \mathcal{M}(\mathbb{R}; [r_1, r_2])$ with given $r_1$ and $r_2$ due to the nature of the proposed methods. This also holds true for other types of dissipative constraints.

### 7.5.2 Dissipative stabilization of a linear system with a time-varying distributed delay

Consider a system in the form of (7.1) with state space parameters

$$
\begin{aligned}
&A_1 = \begin{bmatrix} -1 & -1.9 \\ 0 & 0.1 \end{bmatrix}, \ \widetilde{A}_2(\tau) = \begin{bmatrix} 0.2\cos(e^\tau) + 0.1\sin(e^\tau) & 0.01\cos(e^\tau) - 0.1\sin(e^\tau) \\ 0 & -0.4\cos(e^\tau) \end{bmatrix}, \ \tau \in [-r_1, 0] \\
&\widetilde{A}_2(\tau) = \begin{bmatrix} 0.2\cos(e^\tau) + 0.1\sin(e^\tau) - 0.2 & 0.01\cos(e^\tau) - 0.1\sin(e^\tau) + 1 \\ \ln(2 - \cos(\tau)) - 1.2 & 1 - 0.4\cos(e^\tau) \end{bmatrix}, \quad \tau \in [-r(t), -r_1] \\
&B_1 = \begin{bmatrix} 0 \\ 1 \end{bmatrix}, \quad \widetilde{B}_2(\tau) = \begin{bmatrix} 0.1\sin(e^\tau) - 0.1 \\ 0.12\cos(e^\tau) + 0.1 \end{bmatrix}, D_1 = \begin{bmatrix} 0.01 \\ 0.02 \end{bmatrix}, C_1 = \begin{bmatrix} 0.1 & 0.15 \\ 0 & -0.2 \end{bmatrix}, \\
&\widetilde{C}_2(\tau) = \begin{bmatrix} 0.2\sin(e^\tau) + 0.1 & 0.1 \\ -0.2\sin(e^\tau) & 0.3\sin(e^\tau) - 0.1\cos(e^\tau) \end{bmatrix}, \quad B_4 = \begin{bmatrix} 0 \\ 0.1 \end{bmatrix} \\
&\widetilde{B}_5(\tau) = \begin{bmatrix} 0 \\ 0.1 - 0.1\sin(e^\tau) \end{bmatrix}, \quad D_2 = \begin{bmatrix} 0.1 \\ 0.2 \end{bmatrix}.
\end{aligned} \tag{7.98}
$$

and $r_1 = 0.5$, $r_2 = 1$. Moreover, let

$$
J_1 = -\gamma I_m, \quad \widetilde{J} = I_m, \quad J_2 = \mathsf{O}_{m,q}, \quad J_3 = \gamma I_q \tag{7.99}
$$

for the supply rate function (7.42) to calculate the minimum value of $\mathcal{L}^2$ gain $\gamma$.

To the best of our knowledge, no existing methods may handle a dissipative stabilization problem for open-loop system (7.1) with the parameters in (7.98).

By observing the functions inside of the DDs $\widetilde{A}_2(\cdot), \widetilde{B}_2(\cdot), \widetilde{C}_2(\cdot), \widetilde{B}_5(\cdot)$, we select

$$
\begin{aligned}
&\boldsymbol{f}_1(\tau) = \begin{bmatrix} 1 \\ \sin(e^\tau) \\ \cos(e^\tau) \end{bmatrix}, \quad \boldsymbol{f}_2(\tau) = \begin{bmatrix} 1 \\ \sin(e^\tau) \\ \cos(e^\tau) \\ \ln(2 - \cos\tau) \end{bmatrix}, \ \boldsymbol{\varphi}_1(\tau) = \begin{bmatrix} e^\tau \cos(e^\tau) \\ e^\tau \sin(e^\tau) \end{bmatrix}, \ \boldsymbol{\varphi}_2(\tau) = \begin{bmatrix} e^\tau \cos(e^\tau) \\ e^\tau \sin(e^\tau) \\ \dfrac{\sin\tau}{2 - \cos\tau} \end{bmatrix} \\
&M_1 = \begin{bmatrix} 0 & 0 & 0 & 0 & 0 \\ 1 & 0 & 0 & 0 & 0 \\ 0 & -1 & 0 & 0 & 0 \end{bmatrix}, \quad M_2 = \begin{bmatrix} 0 & 0 & 0 & 0 & 0 & 0 & 0 \\ 1 & 0 & 0 & 0 & 0 & 0 & 0 \\ 0 & -1 & 0 & 0 & 0 & 0 & 0 \\ 0 & 0 & 1 & 0 & 0 & 0 & 0 \end{bmatrix}
\end{aligned} \tag{7.100}
$$

and

$$
\begin{aligned}
&A_2 = \begin{bmatrix} 0 & 0 & 0 & 0 & 0 & 0 & 0.1 & -0.1 & 0.2 & 0.01 \\ 0 & 0 & 0 & 0 & 0 & 0 & 0 & 0 & 0 & -0.4 \end{bmatrix}, \ A_3 = \begin{bmatrix} 0 & 0 & 0 & 0 & 0 & 0 & -0.2 & 1 & 0.1 & -0.1 & 0.2 & 0.01 & 0 & 0 \\ 0 & 0 & 0 & 0 & 0 & 0 & -1.2 & 1 & 0 & 0 & 0 & -0.4 & 1 & 0 \end{bmatrix} \\
&B_2 = \begin{bmatrix} 0 & 0 & -0.1 & 0.1 & 0 \\ 0 & 0 & 0.1 & 0 & 0.12 \end{bmatrix}, \quad B_3 = \begin{bmatrix} 0 & 0 & 0 & -0.1 & 0.1 & 0 & 0 \\ 0 & 0 & 0 & 0.1 & 0 & 0.12 & 0 \end{bmatrix} \\
&C_2 = \begin{bmatrix} 0 & 0 & 0 & 0 & 0.1 & 0.1 & 0.2 & 0 & 0 & 0 \\ 0 & 0 & 0 & 0 & 0 & 0 & -0.2 & 0.3 & 0 & -0.1 \end{bmatrix}, \quad C_3 = \begin{bmatrix} 0 & 0 & 0 & 0 & 0 & 0 & 0.1 & 0.1 & 0.2 & 0 & 0 & 0 & 0 & 0 \\ 0 & 0 & 0 & 0 & 0 & 0 & 0 & 0 & -0.2 & 0.3 & 0 & -0.1 & 0 & 0 \end{bmatrix} \\
&B_5 = \begin{bmatrix} 0 & 0 & 0 & 0 & 0 \\ 0 & 0 & 0.1 & -0.1 & 0 \end{bmatrix}, \quad B_6 = \begin{bmatrix} 0 & 0 & 0 & 0 & 0 & 0 & 0 \\ 0 & 0 & 0 & 0.1 & -0.1 & 0 & 0 \end{bmatrix}.
\end{aligned} \tag{7.101}
$$

to utilize Proposition 1 with $d_1 = 3, d_2 = 4, \delta_1 = 2, \delta_3 = 3, n = m = 2, q = 1$.

Now apply Algorithm 4 to (7.19) with the parameters in (7.98)–(7.101) and $\alpha_3 = 0.5$, $\alpha_1 = \alpha_2 = \alpha_i = 0$, $i = 4, \ldots, 12$ for using Theorem 7.2. This produces controller gains that guarantee $\min \gamma$ as shown in Table 7.3, where Number of Iterations represents the number of iterations in the while loop within Algorithm 4.

| Controller gain $K$ | $\begin{bmatrix} 0.4182 \\ -2.7551 \end{bmatrix}^\top$ | $\begin{bmatrix} 0.5011 \\ -2.7108 \end{bmatrix}^\top$ | $\begin{bmatrix} 0.5787 \\ -2.6595 \end{bmatrix}^\top$ | $\begin{bmatrix} 0.6505 \\ -2.6021 \end{bmatrix}^\top$ |
|---|---|---|---|---|
| $\min \gamma$ | 0.36657 | 0.3607 | 0.3551 | 0.3498 |
| Number of Iterations | 10 | 20 | 30 | 40 |

**Table 7.3:** Controller gains with $\min \gamma$ produced with different iterations

Since the expression of $r(t)$ is unknown, existing approaches based on the frequency domain may not be directly applicable to analyzing the stability of the resulting CLSs obtained using our methods. To partially verify the results in Table 7.3, we assume $r(t)$ to be an unknown constant $\widehat{r} \in [r_1, r_2]$. This allows one to apply the spectral method in [203] to compute the spectral abscissa of the spectrum of the resulting CLSs with a constant delay. As our synthesis results indicate that any resulting CLS is stable for all $r(\cdot) \in \mathcal{M}\left(\mathbb{R}; [r_1, r_2]\right)$, the same CLSs with a constant delay $\widehat{r}$ are stable for $\widehat{r} \in [r_1, r_2]$ as the case of $r(t) = \widehat{r}$ is included by $\mathcal{M}\left(\mathbb{R}; [r_1, r_2]\right)$. The numerical results produced by [203] show that all the resulting CLSs are stable for $\widehat{r} \in [r_1, r_2]$ under the hypothesis that $r(t) = \widehat{r}$ is a constant delay.

We considered the CLSs stabilized $K = \begin{bmatrix} 0.6505 & -2.6021 \end{bmatrix}$ in Table 7.3 for numerical simulation. Specifically, we assumed $t_0 = 0$, $\boldsymbol{z}(t) = \mathbf{0}_2, t < 0$, and $\boldsymbol{\phi}(\tau) = \begin{bmatrix} 5 & 3 \end{bmatrix}^\top, \tau \in [-1, 0]$ as the initial condition, and $\boldsymbol{w}(t) = 50 \sin 20\pi t(\mathit{1}(t) - \mathit{1}(t-5))$ as the disturbance where $\mathit{1}(t)$ is the Heaviside step function. Furthermore, we considered a time-varying delay $r(t) = 0.75 + 0.25\cos(100t)$ which[7] exhibits strong oscillation. Numerical simulations were performed in Simulink with the aforementioned data via the ODE solver `ode8` with 0.0001 as the fundamental sampling time. Our simulation results are presented in Figure 7.1–7.3 concerning the trajectories of the states, outputs and the controller feedback. Note that the update method of the MATLAB function block in Simulink was set to discrete, and all DD integrals in the simulations were discretized using the trapezoidal rule

$$
\begin{aligned}
&\int_{-r_2}^{0} F(t,\tau)\boldsymbol{x}(t+\tau)\mathrm{d}\tau \approx \\
&\frac{r_2}{n}\left[\frac{F(t,-r_2)\boldsymbol{x}(t-r_2)}{2} + \sum_{k=1}^{n-1} F\left(t, r_2\frac{k}{n} - r_2\right)\boldsymbol{x}\left(t + r_2\frac{k}{n} - r_2\right) + \frac{F(t,0)\boldsymbol{x}(t)}{2}\right] \\
&= \frac{r_2}{n}\left[\frac{F(t,-r_2)\boldsymbol{x}(t-r_2)}{2} + \sum_{k=1}^{n-1} F\left(t, -r_2\frac{n-k}{n}\right)\boldsymbol{x}\left(t - r_2\frac{n-k}{n}\right) + \frac{F(t,0)\boldsymbol{x}(t)}{2}\right] \\
&= \frac{r_2}{n}\left[\frac{F(t,-r_2)\boldsymbol{x}(t-r_2)}{2} + \sum_{k=1}^{n-1} F\left(t, -r_2\frac{k}{n}\right)\boldsymbol{x}\left(t - r_2\frac{k}{n}\right) + \frac{F(t,0)\boldsymbol{x}(t)}{2}\right]
\end{aligned} \tag{7.102}
$$

given the fact that $\sum_{k=1}^{n-1} f(n-k) = \sum_{k=1}^{n-1} f(k)$, where

$$
F(t,\tau) := \begin{cases} \widetilde{F}(\tau) & \forall \tau \in [-r(t), 0] \\ 0 & \forall \tau \in [-r_2, r(t)), \end{cases} \tag{7.103}
$$

and $\widetilde{F}(\tau)$ is piecewise continuous on $[-r(t), 0]$.

**Remark 7.14.** Note that (7.103) enables one to discretize $\int_{-r(t)}^{0} \widetilde{F}(\tau)\boldsymbol{x}(t+\tau)\mathrm{d}\tau$ via (7.102), which avoids dealing with $\int_{-r(t)}^{0} \widetilde{F}(\tau)\boldsymbol{x}(t+\tau)\mathrm{d}\tau$ directly.

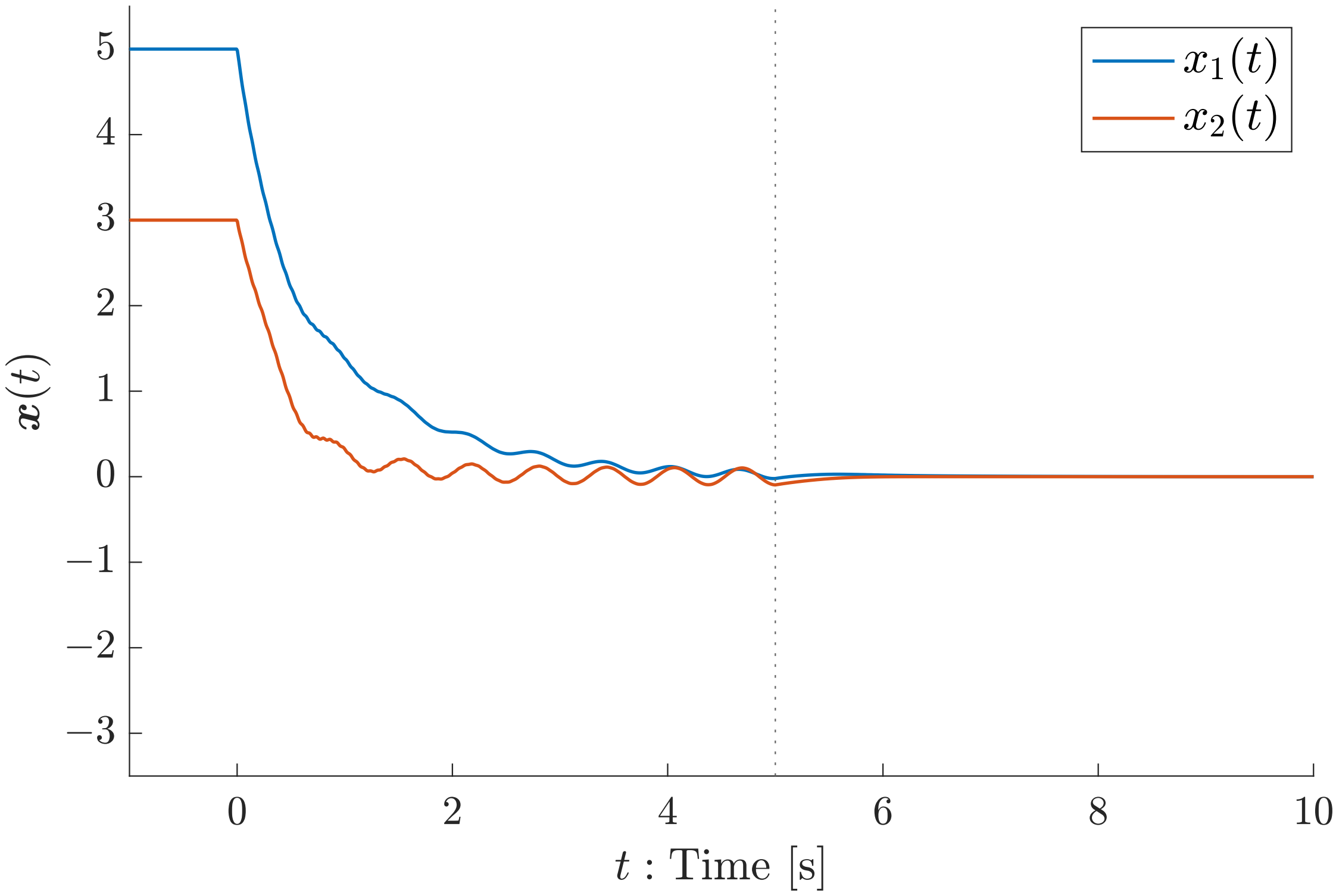

**Figure 7.1:** The close-loop system's trajectory $\boldsymbol{x}(t)$ with $K = \begin{bmatrix} 0.6505 & -2.6021 \end{bmatrix}$ in Table 7.3

**Remark 7.15.** Due to the lack of appropriate numerical solvers in Simulink for delay systems, we were limited to using an ODE solver (ode8) in Simulink to conduct our simulation. As we cannot predict the potential problems of utilizing an ODE solver on a delay system, the numerical results in Figure 7.1–7.3 only provide an estimation of the actual behavior of the system trajectories and output, and the numerical accuracy in this case may not be guaranteed.

**Remark 7.16.** The results in Figure 7.1–7.3 unequivocally demonstrate the effectiveness of the proposed stabilization method when considering a time-varying delay in the form of $r(t) = 0.75 + 0.25\cos(100t)$. It should be noted that the abrupt change around $t = 5$ in Figure 7.3 is due to the form of the disturbance signal $\boldsymbol{w}(t) = 50\sin 20\pi t(\mathit{1}(t) - \mathit{1}(t-5))$ which satisfies $\forall t > 5, \boldsymbol{w}(t) = 0$.

[7]Note that this function satisfying $\forall t \geq t_0, r_1 = 0.5 \leq r(t) \leq 1 = r_2$

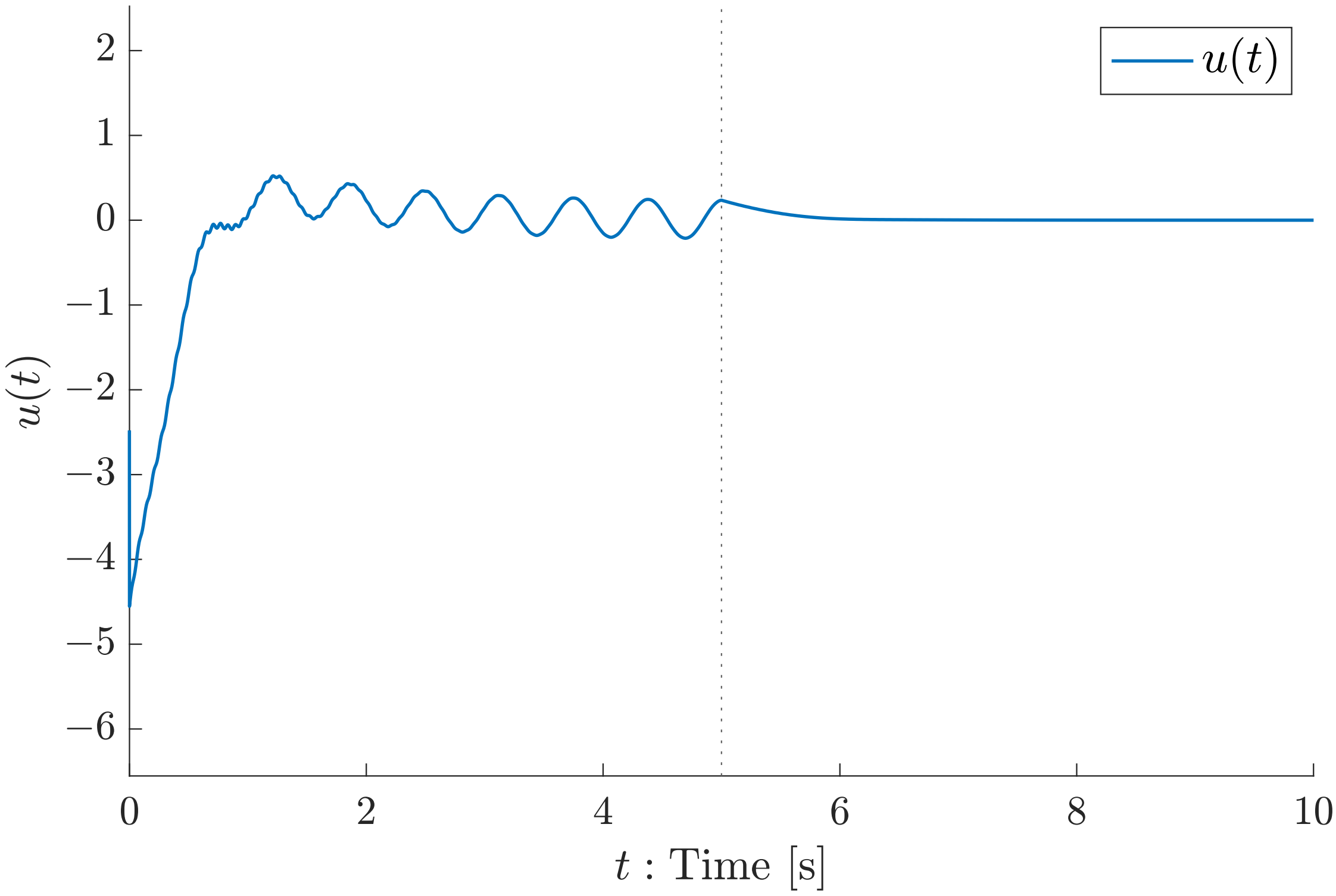

**Figure 7.2:** The trajectory of the controller effort $\boldsymbol{u}(t) = K\boldsymbol{x}(t)$ with $K = \begin{bmatrix} 0.6505 & -2.6021 \end{bmatrix}$ in Table 7.3

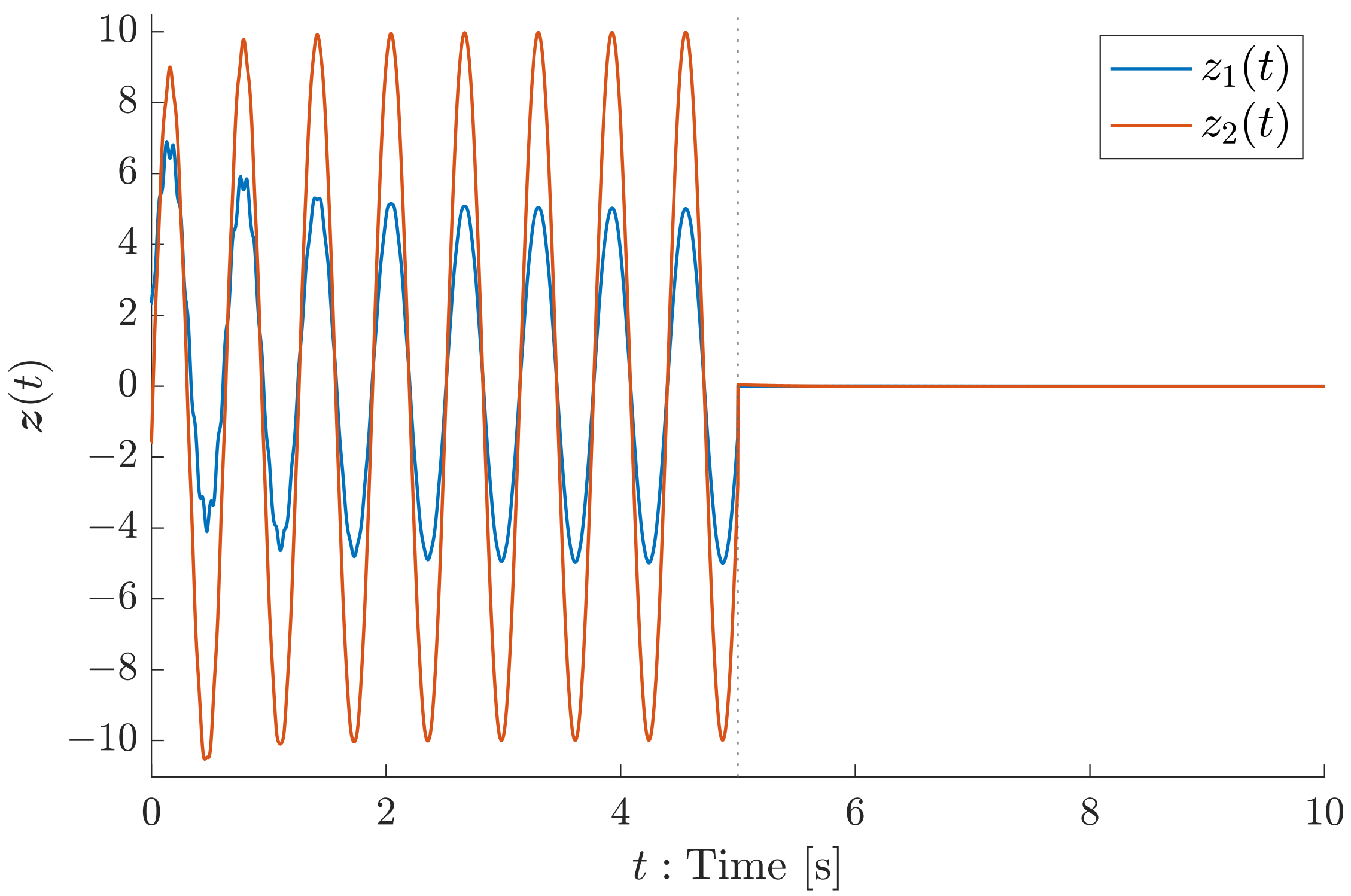

**Figure 7.3:** The output $\boldsymbol{z}(t)$ of the CLS with $K = \begin{bmatrix} 0.6505 & -2.6021 \end{bmatrix}$ in Table 7.3

# Chapter 8

# Summaries and Future Works

In this book, we have presented effective methods for the stability analysis and stabilization of linear systems with DDs possessing nontrivial kernels. Solutions to the stability (dissipativity) analysis and stabilization problems have been derived using the LKF approach via novel IIs. The conclusions of the results presented in this book are summarized in the following section.

## 8.1 Conclusions

- In Chapter 2, we investigated the problem of stabilizing an LDDS with DDs in the states, input, and output. The stabilization problem incorporated a dissipative constraint expressed by the quadratic supply rate function (2.10) to ensure the performance of the resulting controllers. By constructing the LKF in (2.20) with the newly proposed inequality (2.12), sufficient conditions for the existence of a dissipative stabilizing controller are derived in terms of matrix inequalities in Theorem 2.1. To solve the resulting BMI (2.21b) in Theorem 2.1, we developed convex synthesis conditions via the application of the Projection Lemma in Theorem 2.2, whose feasible solutions imply the existence of the feasible solutions for Theorem 2.1. To reduce the potential conservatism of Theorem 2.2 due to the simplification of $\mathbf{W}$ at (2.40), an iterative algorithm was proposed as Algorithm 1 based on the ideas developed in [565]. Given the generality of our LKF in (2.20) and the use of the novel integral inequality in (2.12), our proposed synthesis solutions can produce nonconservative results requiring fewer decision variables compared to existing methods in the literature. Two numerical examples we have investigated confirm this conclusion, which clearly demonstrate the advantages of our proposed methodology over existing approaches.

- In Chapter 3, we developed dissipative stabilization conditions for uncertain LDDSs with general uncertainties among the systems state space matrices. An instrumental mathematical device designed to tackle the general uncertainties in (3.1) was proposed as Lemma 3.1 using matrix inequalities. This allows us to derive dissipative synthesis conditions for uncertain CLS (3.5), as enunciated in Section 3.3. It has also been shown in Section 3.4 that the idea presented in Section 3.3 can be further modified to calculate the gains of a resilient dynamical state feedback controller for uncertain input delay system (3.36), where the controller itself is robust against uncertainties in general forms. It is crucial to point out that the design of a resilient controller is made possible by the mathematical structure of (3.41) with a dynamic state feedback controller. Namely, there are no matrix products among the uncertainties of the controller gains and input gain matrix, which allows us to handle the uncertainties in (3.41) in a manner similar to (3.5).

- In Chapter 4, we explored two distinct classes of IIs, delineated in Theorem 4.1, 4.2 and Theorem 4.3, where we also have addressed the relationships concerning inequality bound gaps. The inequalities (4.3),(4.17) and (4.23) serve to generalize the vast majority of existing IIs found in the literature, with many being essentially equivalent in terms of inequality bound gaps within the context of SDP. Moreover, we have demonstrated the utility of our proposed inequalities and their properties in deriving equivalent stability conditions for a linear CDDS with a DD. Our inequalities possess significant potential for the applications in broader contexts such as the stability analysis of PDE-related systems, sampled-data systems, or other types of IDSs whenever the contexts are suitable.

- In Chapter 5, a new method for the dissipativity and stability analysis of a linear CDDS with DDs in state and output equations has been proposed in Theorem 5.1 in terms of LMIs. Our approach is capable of handling DDs with $\mathcal{L}^2$ kernels and simultaneously includes approximation errors in the resulting conditions (5.34)–(5.36), owing to the innovative integral inequality proposed in (5.24). In contrast to existing approaches such as the one in [470] which relies on the application of Legendre polynomial approximations, our method allows for the application of a broader class of elementary functions $\acute{\boldsymbol{f}}(\cdot)$ and $\grave{\boldsymbol{f}}(\cdot)$ to approximate the DD kernels of (5.1). Since the generality of the LKF in (5.43) is also related to the structure of $\acute{\boldsymbol{f}}(\cdot)$ and $\grave{\boldsymbol{f}}(\cdot)$, our proposed methods derived from constructing (5.43) can produce less conservative results compared to a functional parameterized by Legendre polynomials such as the one considered in [470]. The testing results of numerical examples in Chapter 5 have clearly demonstrated the advantages of the proposed methodologies over existing approaches.A potential avenue for future exploration is to investigate whether the hierarchy conclusion can be derived without imposing an orthogonality constraint.

- In Chapter 6, we introduced a novel solution to the problem of delay-range dissipativity and stability analysis of CDDS (6.1) with polynomial-kernel-DD, as described in Theorem 6.1 using SoS constraints (6.18)–(6.21) based on the construction of an LKF (6.24). The superiority of our proposed methodologies lies in the form of functional (6.24) with delay-dependent matrix parameters, which mathematically results in less conservative conditions in terms of delay-parameter-dependent LMIs in (6.31) and (6.37). The challenge of numerically solving robust LMIs is circumvented by applying the matrix relaxation technique set forth in [649], resulting in SoS constraints (6.18)–(6.21) which are equivalent to (6.31) and (6.37). Meanwhile, it is shown in Section 6.4.3 that on certain occasions, it is unnecessary to solve SoS constraints (6.18)–(6.21). Rather, one may solve (6.31) and (6.37) directly utilizing the convex hull's properties. Moreover, a solution to the delay margins estimation problem with prescribed performance objectives is also presented in Section 6.4.4 which is followed by the feasibility hierarchy established in Theorem 6.2. Finally, the tests of numerical examples in Section 6.5 have confirmed the advantage of our proposed methods that less conservative results with reduced computational burden can be produced compared to existing approaches.

- In Chapter 7, solutions to the design of a dissipative state feedback controller of a linear system with DDs (7.1) have been proposed where the delay $r(\cdot) \in \mathcal{M}([r_1, r_2]; \mathbb{R})$ is time-varying, measurable and bounded. The key step in deriving the synthesis condition in Theorem 7.1 is the application of the novel inequality proposed in Lemma 7.3, which led to simple LMI terms as enunciated in Section 7.4.1. Although (7.45) in Theorem 7.1 is bilinear, it has been shown in Theorem 7.2 that convex conditions (7.69)–(7.71) are constructible by the application of the Projection Lemma. Moreover, Algorithm 4 is proposed to solve the BMI in Theorem 7.1 iteratively, which can be initiated through the feasible solutions to Theorem 7.2. Considering the generality of $r(\cdot) \in \mathcal{M}([r_1, r_2]; \mathbb{R})$, our proposed methodologies have the potential to handle a large class of delay scenarios and solve numerous delay-related problems.

## 8.2 Future works

The results presented in this book offer the potential for further extension to address new challenges concerning the control of TDSs or IDSs. Here, we propose several potential avenues for future exploration.

- Regarding the results in Chapter 2 and Chapter 3, we can broaden the scope of the kernel function $\boldsymbol{f}(\cdot)$ to include any differentiable functions. This would necessitate modifying Assumption 2.1 to cope with the new functions produced by $\frac{\mathrm{d}\boldsymbol{f}(\tau)}{\mathrm{d}\tau}$. Since the inequality (2.12) can address any $\boldsymbol{g}(\cdot) \in \mathcal{L}_2(\mathcal{K}; \mathbb{R}^d)$, then it can be utilized to handle $\boldsymbol{f}(\cdot) \in \mathcal{C}^1(\mathcal{K}; \mathbb{R}^d)$ in the construction of LKFs. Additionally, the challenges of designing filters and dynamical output feedback controllers may also be considered.

- The results in Chapter 3 demonstrate the feasibility of utilizing the LKF approach to construct a dynamical state feedback controller resembling the structure of a predictor controller for an uncertain linear system with input delays. It would be of significant research interest to extend this strategy to address the challenge of designing dissipative dynamical state feedback controllers for linear systems with pointwise-delays in both states and inputs. Existing results on predictor controllers in [495, 497, 498, 500–502] have indicated that the effects of a point-wise input delay in a linear system with a state delay can be totally eliminated by the use of a predictor controller possessing DDs. It would also be interesting to explore the challenge of designing a dissipative dynamical state feedback controller for a linear system with DDs at its inputs.

- It would be advantageous to extend the inequalities proposed in Chapter 4 to consider more general structures, such as inequalities in Hilbert space.

- In Chapter 5 the values of delays are assumed to be known. For real-time operations, it would be preferable to consider that delay values are uncertain, as considered in Chapter 6. If delay values are unknown, the approximation scheme and Lemma 5.2 in Chapter 5 would need to be modified to produce constant approximation coefficients that are anticipated to be independent of the uncertain delay $r$.

- It is natural to extend the results in Chapter 6 to address a dissipative range synthesis problem.

# Appendix A

# Proof of Lemma 3.1

*Proof.* The proof draws inspiration from the strategies illustrated in [452]. Let us first consider the scenario where $\mathcal{F} = \mathcal{D}$. Our initial task is to construct an equivalent condition that can guarantee the well-posedness of $(I_m - \Delta F)^{-1}$ for all $\Delta \in \mathcal{D}$ with $\mathcal{D}$ in (3.8). Assume first that $F \neq \mathsf{O}_{p,m}$, it is evident that $(I_m - \Delta F)^{-1}$ is well defined for all $\Delta \in \mathcal{D}$ if and only if $\forall \Delta \in \mathcal{D}$, $\operatorname{rank}(I_m - \Delta F) = m$. Moreover, we can deduce that $\forall \Delta \in \mathcal{D}$, $\operatorname{rank}(I_m - \Delta F) = m$ if and only if $\forall \Delta \in \mathcal{D}$, $(I_m - \Delta F)^\top (I_m - \Delta F) \succ 0$, based on the properties of ranks and Gramian matrices.

Let $\mathbb{R}^m \ni \boldsymbol{\mu} := \Delta F \boldsymbol{\theta}$ and $\mathbb{R}^p \ni \boldsymbol{\omega} := F\boldsymbol{\theta}$ with $\boldsymbol{\theta} \in \mathbb{R}^m$, we can deduce that $\forall \Delta \in \mathcal{D}$, $(I_m - \Delta F)^\top (I_m - \Delta F) \succ 0$ if and only if $\forall \boldsymbol{\theta} \in \mathbb{R}^m \backslash \{\mathbf{0}_m\}$, $\forall \Delta \in \mathcal{D}$, $\boldsymbol{\theta}^\top (I_m - \Delta F)^\top (I_m - \Delta F)\boldsymbol{\theta} > 0$ which is equivalent to

$$\begin{bmatrix} \boldsymbol{\theta} \\ \boldsymbol{\mu} \end{bmatrix}^\top \begin{bmatrix} I_m & -I_m \\ * & I_m \end{bmatrix} \begin{bmatrix} \boldsymbol{\theta} \\ \boldsymbol{\mu} \end{bmatrix} > 0, \ \ \forall \begin{bmatrix} \boldsymbol{\theta} \\ \boldsymbol{\mu} \end{bmatrix} \in \mathcal{M} \setminus \{\mathbf{0}_{2m}\}, \tag{A.1}$$

$$\text{with } \ \mathcal{M} \in \left\{ \begin{bmatrix} \acute{\boldsymbol{\theta}} \\ \acute{\boldsymbol{\mu}} \end{bmatrix} \middle| \, \acute{\boldsymbol{\mu}} = \acute{\Delta} F \acute{\boldsymbol{\theta}} \ \ \& \ \ \acute{\Delta} \in \mathcal{D} \ \ \& \ \ \acute{\boldsymbol{\theta}} \in \mathbb{R}^m \right\}. \tag{A.2}$$

Given the definition of $\mathcal{D}$ and the property of quadratic forms, it holds true that $\forall \boldsymbol{\theta} \in \mathbb{R}^m$ and $\forall \Delta \in \mathcal{D}$ we have

$$\boldsymbol{\theta}^\top F^\top \begin{bmatrix} I \\ \Delta \end{bmatrix}^\top \begin{bmatrix} \Theta_1^{-1} & \Theta_2 \\ * & \Theta_3 \end{bmatrix} \begin{bmatrix} I \\ \Delta \end{bmatrix} F\boldsymbol{\theta} = \begin{bmatrix} F\boldsymbol{\theta} \\ \Delta F\boldsymbol{\theta} \end{bmatrix}^\top \begin{bmatrix} \Theta_1^{-1} & \Theta_2 \\ * & \Theta_3 \end{bmatrix} \begin{bmatrix} F\boldsymbol{\theta} \\ \Delta F\boldsymbol{\theta} \end{bmatrix} \geq 0. \tag{A.3}$$

Now if $\begin{bmatrix} \boldsymbol{\theta}^\top & \boldsymbol{\mu}^\top \end{bmatrix}^\top \in \mathcal{M}$, then

$$\begin{bmatrix} F\boldsymbol{\theta} \\ \boldsymbol{\mu} \end{bmatrix}^\top \begin{bmatrix} \Theta_1^{-1} & \Theta_2 \\ * & \Theta_3 \end{bmatrix} \begin{bmatrix} F\boldsymbol{\theta} \\ \boldsymbol{\mu} \end{bmatrix} = \begin{bmatrix} F\boldsymbol{\theta} \\ \Delta F\boldsymbol{\theta} \end{bmatrix}^\top \begin{bmatrix} \Theta_1^{-1} & \Theta_2 \\ * & \Theta_3 \end{bmatrix} \begin{bmatrix} F\boldsymbol{\theta} \\ \Delta F\boldsymbol{\theta} \end{bmatrix} \geq 0 \tag{A.4}$$

for some $\Delta \in \mathcal{D}$ given the property in (A.3). Therefore,

$$\mathcal{M} \subseteq \mathcal{J} := \left\{ \begin{bmatrix} \acute{\boldsymbol{\theta}} \\ \acute{\boldsymbol{\mu}} \end{bmatrix} \in \mathbb{R}^{2m} \middle| \begin{bmatrix} \acute{\boldsymbol{\theta}} \\ \acute{\boldsymbol{\mu}} \end{bmatrix}^\top \begin{bmatrix} F & \mathsf{O}_{p,m} \\ \mathsf{O}_m & I_m \end{bmatrix}^\top \begin{bmatrix} \Theta_1^{-1} & \Theta_2 \\ * & \Theta_3 \end{bmatrix} \begin{bmatrix} F & \mathsf{O}_{p,m} \\ \mathsf{O}_m & I_m \end{bmatrix} \begin{bmatrix} \acute{\boldsymbol{\theta}} \\ \acute{\boldsymbol{\mu}} \end{bmatrix} \geq 0 \right\}. \tag{A.5}$$

By (A.5), it follows that

$$\begin{bmatrix} \boldsymbol{\theta} \\ \boldsymbol{\mu} \end{bmatrix}^\top \begin{bmatrix} I_m & -I_m \\ * & I_m \end{bmatrix} \begin{bmatrix} \boldsymbol{\theta} \\ \boldsymbol{\mu} \end{bmatrix} > 0, \ \ \forall \begin{bmatrix} \boldsymbol{\theta} \\ \boldsymbol{\mu} \end{bmatrix} \in \mathcal{J} \setminus \{\mathbf{0}_{2m}\} \tag{A.6}$$

holds then (A.1) holds since $\mathcal{M} \subseteq \mathcal{J}$. By applying S-procedure to (A.6), we can conclude that (A.6) is true if and only if [1]

$$\exists \alpha > 0 : \begin{bmatrix} I_m & -I_m - \alpha F^\top \Theta_2 \\ * & I_m - \alpha\Theta_3 \end{bmatrix} - \alpha \begin{bmatrix} F^\top \\ \mathsf{O}_{m,p} \end{bmatrix} \Theta_1^{-1} \begin{bmatrix} F & \mathsf{O}_{p,m} \end{bmatrix}$$

$$\begin{bmatrix} I_m & -I_m - \alpha F^\top \Theta_2 \\ * & I_m - \alpha\Theta_3 \end{bmatrix} - \begin{bmatrix} \alpha F^\top \\ \mathsf{O}_{m,p} \end{bmatrix} (\alpha\Theta_1)^{-1} \begin{bmatrix} \alpha F & \mathsf{O}_{p,m} \end{bmatrix} \succ 0. \quad \text{(A.7)}$$

Utilizing the Schur complement on (A.7) and taking into account that $\Theta_1^{-1} \succ 0$, we obtain (3.7). Given that $\forall \Delta \in \mathcal{D}$, $(I_m - \Delta F)^\top (I_m - \Delta F) \succ 0$ is equivalent to (A.1), which is implied by (3.7), we can conclude that $(I_m - \Delta F)^{-1}$, $F \neq \mathsf{O}_{p,m}$ is well defined for all $\Delta \in \mathcal{D}$ if (3.7) holds. Finally, it is evident that the well-posedness of $(I_m - \Delta F)^{-1}$ is automatically achieved if $F = \mathsf{O}_{p,m}$, where (3.7) need not be considered.

Now we shall proceed to prove the remaining results in Lemma 3.1. Let us assume that (3.7) holds true for a $\alpha > 0$ indicating that $(I_m - \Delta F)^{-1}$ is well defined. By the definition of negative definite matrices, we know that

$$\forall \Delta \in \mathcal{D}, \quad \Phi + \mathsf{Sy}\left[G(I_m - \Delta F)^{-1} \Delta H\right] \prec 0 \quad \text{(A.8)}$$

if and only if

$$\forall \Delta \in \mathcal{D}, \quad \forall \mathbf{x} \in \mathbb{R}^n \setminus \{\mathbf{0}_n\}, \quad \mathbf{x}^\top \left(\Phi + \mathsf{Sy}\left[G(I_m - \Delta F)^{-1} \Delta H\right]\right) \mathbf{x} < 0. \quad \text{(A.9)}$$

Now let $\mathbb{R}^m \ni \boldsymbol{\varrho} := (I_m - \Delta F)^{-1} \Delta H \mathbf{x}$ and consider the fact that $(I_m - \Delta F)^{-1}$ is well defined, we have $\boldsymbol{\varrho} = \Delta H \mathbf{x} + \Delta F \boldsymbol{\varrho}$ by which we can reformulate (A.9) as

$$\begin{bmatrix} \mathbf{x} \\ \boldsymbol{\varrho} \end{bmatrix}^\top \begin{bmatrix} \Phi & G \\ * & \mathsf{O}_m \end{bmatrix} \begin{bmatrix} \mathbf{x} \\ \boldsymbol{\varrho} \end{bmatrix} < 0, \quad \forall \begin{bmatrix} \mathbf{x} \\ \boldsymbol{\varrho} \end{bmatrix} \in \mathcal{X} \setminus \{\mathbf{0}_{n+m}\} \quad \text{(A.10)}$$

$$\text{with } \mathcal{X} \in \left\{ \begin{bmatrix} \acute{\mathbf{x}} \\ \acute{\boldsymbol{\varrho}} \end{bmatrix} \in \mathbb{R}^{n+m} \,\middle|\, \acute{\boldsymbol{\varrho}} := (I_m - \acute{\Delta} F)^{-1} \acute{\Delta} H \acute{\mathbf{x}} \;\;\&\;\; \acute{\Delta} \in \mathcal{D} \right\}.$$

Based on the definition of $\mathcal{D}$ and the property of quadratic forms, it holds true that $\forall \mathbf{x} \in \mathbb{R}^n$ and $\forall \Delta \in \mathcal{D}$ we have

$$(H\mathbf{x} + F\boldsymbol{\varrho})^\top \begin{bmatrix} I \\ \Delta \end{bmatrix}^\top \begin{bmatrix} \Theta_1^{-1} & \Theta_2 \\ * & \Theta_3 \end{bmatrix} \begin{bmatrix} I \\ \Delta \end{bmatrix} (H\mathbf{x} + F\boldsymbol{\varrho}) =$$

$$\begin{bmatrix} H\mathbf{x} + F\boldsymbol{\varrho} \\ \Delta\,(H\mathbf{x} + F\boldsymbol{\varrho}) \end{bmatrix}^\top \begin{bmatrix} \Theta_1^{-1} & \Theta_2 \\ * & \Theta_3 \end{bmatrix} \begin{bmatrix} H\mathbf{x} + F\boldsymbol{\varrho} \\ \Delta\,(H\mathbf{x} + F\boldsymbol{\varrho}) \end{bmatrix} = [*]^\top \begin{bmatrix} \Theta_1^{-1} & \Theta_2 \\ * & \Theta_3 \end{bmatrix} \begin{bmatrix} H\mathbf{x} + F\boldsymbol{\varrho} \\ \boldsymbol{\varrho} \end{bmatrix} \geq 0, \quad \text{(A.11)}$$

where $\boldsymbol{\varrho} = (I_m - \Delta F)^{-1} \Delta H \mathbf{x}$ that is equivalent to $\boldsymbol{\varrho} = \Delta H \mathbf{x} + \Delta F \boldsymbol{\varrho}$. It becomes evident that

$$\forall \begin{bmatrix} \mathbf{x} \\ \boldsymbol{\varrho} \end{bmatrix} \in \mathcal{X}, \quad \begin{bmatrix} \mathbf{x} \\ \boldsymbol{\varrho} \end{bmatrix}^\top \begin{bmatrix} H & F \\ \mathsf{O}_{m,p} & I_m \end{bmatrix}^\top \begin{bmatrix} \Theta_1^{-1} & \Theta_2 \\ * & \Theta_3 \end{bmatrix} \begin{bmatrix} H & F \\ \mathsf{O}_{m,p} & I_m \end{bmatrix} \begin{bmatrix} \mathbf{x} \\ \boldsymbol{\varrho} \end{bmatrix} \geq 0 \quad \text{(A.12)}$$

considering the property in (A.11) and the fact that

$$\forall \begin{bmatrix} \mathbf{x} \\ \boldsymbol{\varrho} \end{bmatrix} \in \mathcal{X}, \; \exists \Delta \in \mathcal{D}, \quad \boldsymbol{\varrho} = \Delta H \mathbf{x} + \Delta F \boldsymbol{\varrho}.$$

[1]Since (A.7) with $\alpha = 0$ cannot be feasible, we only need $\alpha > 0$ in (A.15) for the existence of $\alpha \geq 0$ introduced by the application of S-procedure.

Accordingly, one can derive the following relation

$$\mathcal{X} \subseteq \mathcal{Y} := \left\{ \begin{bmatrix} \acute{\mathbf{x}} \\ \acute{\varrho} \end{bmatrix} \;\middle|\; \begin{bmatrix} \acute{\mathbf{x}} \\ \acute{\varrho} \end{bmatrix}^\top \begin{bmatrix} H & F \\ \mathsf{O}_{m,p} & I_m \end{bmatrix}^\top \begin{bmatrix} \Theta_1^{-1} & \Theta_2 \\ * & \Theta_3 \end{bmatrix} \begin{bmatrix} H & F \\ \mathsf{O}_{m,p} & I_m \end{bmatrix} \begin{bmatrix} \acute{\mathbf{x}} \\ \acute{\varrho} \end{bmatrix} \geq 0 \right\}. \tag{A.13}$$

From (A.13), we see that (A.10) holds if

$$\forall \begin{bmatrix} \mathbf{x} \\ \varrho \end{bmatrix} \in \mathcal{Y} \setminus \{\mathbf{0}_{n+m}\}, \quad \begin{bmatrix} \mathbf{x} \\ \varrho \end{bmatrix}^\top \begin{bmatrix} \Phi & G \\ * & \mathsf{O}_m \end{bmatrix} \begin{bmatrix} \mathbf{x} \\ \varrho \end{bmatrix} < 0, \tag{A.14}$$

holds. By applying the S-procedure [167] to (A.14) with $\mathcal{Y}$, we see that (A.14) is true if and only if [2]

$$\begin{aligned} \exists \kappa > 0: \; & \begin{bmatrix} \Phi & G + \kappa H^\top \Theta_2 \\ * & \kappa F^\top \Theta_2 + \kappa \Theta_2^\top F + \kappa \Theta_3 \end{bmatrix} + \kappa \begin{bmatrix} H^\top \\ F^\top \end{bmatrix} \Theta_1^{-1} \begin{bmatrix} H & F \end{bmatrix} \\ & = \begin{bmatrix} \Phi & G + \kappa H^\top \Theta_2 \\ * & \kappa F^\top \Theta_2 + \kappa \Theta_2^\top F + \kappa \Theta_3 \end{bmatrix} + \begin{bmatrix} \kappa H^\top \\ \kappa F^\top \end{bmatrix} (\kappa \Theta_1)^{-1} \begin{bmatrix} \kappa H & \kappa F \end{bmatrix} \prec 0. \end{aligned} \tag{A.15}$$

It should be emphasized that for cases involving a single constraint such as (A.14), the S-procedure can yield an equivalent condition.

Utilizing the Schur complement on (A.15) with $\Theta_1^{-1} \succ 0$, it yields (3.9), which is equivalent to (A.14). Since (A.14) implies (A.10) which is equivalent to (3.8), we have demonstrated that (3.9) is a sufficient condition for (3.8) with $\mathcal{F} = \mathcal{D}$. Moreover, it is evident that (3.8) and the well-posedness of $(I_m - \Delta F)^{-1}$ with $\mathcal{F} \subseteq \mathcal{D}$ are ensured by (3.7) and (3.9), respectively. Additionally, if $\Theta_1^{-1} = \mathsf{O}_p$, then (A.15) and (A.7) become (3.11) and (3.10), respectively, where $\Theta_1$ does not need to be considered as the Schur complement need not be applied at the steps of (A.15) and (A.7). Finally, it is important to point out that (3.7) and (3.9) can accommodate an uncertainty term $G\Delta(I_p - F\Delta)^{-1}H$ via the Sylvester's determinant identity[3] $\det(I_p - F\Delta) = \det(I_m - \Delta F)$, and $\Delta(I_p - F\Delta)^{-1} = (I_m - \Delta F)^{-1}\Delta$ that is a special case of (A.3) in [452]. ∎

[2] Since (A.15) with $\kappa = 0$ cannot be feasible, $\kappa \geq 0$ is applied in (A.15) for the condition $\kappa \geq 0$ introduced by the application of S-procedure

[3] See B.1.16 in [694]

# Appendix B

# Proof of Theorem 4.1

*Proof.* The proof bears similarities to Lemma 2.3, with the exception that we have a weight function $\varpi(\cdot)$ here.

Let $\boldsymbol{\varepsilon}(\tau) := \boldsymbol{x}(\tau) - F^\top(\tau)(\mathsf{F} \otimes I_n) \int_{\mathcal{K}} \varpi(\theta) F(\theta) \boldsymbol{x}(\theta) \mathrm{d}\theta$, where $F(\tau) = \boldsymbol{f}(\tau) \otimes I_n$. By combining the expression $\boldsymbol{\varepsilon}(\cdot)$ with $\int_{\mathcal{K}} \varpi(\tau) \boldsymbol{\varepsilon}^\top(\tau) U \boldsymbol{\varepsilon}(\tau) \mathrm{d}\tau$, we have

$$
\begin{aligned}
\int_{\mathcal{K}} \varpi(\tau) \boldsymbol{\varepsilon}^\top(\tau) U \boldsymbol{\varepsilon}(\tau) \mathrm{d}\tau = \int_{\mathcal{K}} \varpi(\tau) \boldsymbol{x}^\top(\tau) U \boldsymbol{x}(\tau) \mathrm{d}\tau - 2 \int_{\mathcal{K}} \varpi(\tau) \boldsymbol{x}^\top(\tau) U F^\top(\tau) \mathrm{d}\tau (\mathsf{F} \otimes I_n) \boldsymbol{\vartheta} + \\
\boldsymbol{\vartheta}^\top \int_{\mathcal{K}} \varpi(\tau) (\mathsf{F} \otimes I_n)^\top F(\tau) U F^\top(\tau) (\mathsf{F} \otimes I_n) \mathrm{d}\tau \boldsymbol{\vartheta}, \quad \text{(B.1)}
\end{aligned}
$$

where $\boldsymbol{\vartheta} := \int_{\mathcal{K}} \varpi(\theta) F(\theta) \boldsymbol{x}(\theta) \mathrm{d}\theta$. Now apply (2.1a) to the term $UF^\top(\tau)$ with $F(\tau) = \boldsymbol{f}(\tau) \otimes I_n$ and $U = U^\top$, we have

$$
U(\boldsymbol{f}^\top(\tau) \otimes I_n) = \boldsymbol{f}^\top(\tau) \otimes U = \left(\boldsymbol{f}^\top(\tau) \otimes I_n\right) (I_d \otimes U) = F^\top(\tau)(I_d \otimes U). \quad \text{(B.2)}
$$

Apply (B.2) to some of the terms in (B.1). It follows that

$$
\begin{aligned}
\int_{\mathcal{K}} \varpi(\tau) \boldsymbol{x}^\top(\tau) U F^\top(\tau) \mathrm{d}\tau (\mathsf{F} \otimes I_n) \boldsymbol{\vartheta} = \int_{\mathcal{K}} \varpi(\tau) \boldsymbol{x}^\top(\tau) F^\top(\tau) \mathrm{d}\tau \, (I_d \otimes U) \, (\mathsf{F} \otimes I_n) \boldsymbol{\vartheta} \\
= \boldsymbol{\vartheta}^\top (I_d \otimes U)(\mathsf{F} \otimes I_n) \boldsymbol{\vartheta} = \boldsymbol{\vartheta}^\top (\mathsf{F} \otimes U) \boldsymbol{\vartheta}. \quad \text{(B.3)}
\end{aligned}
$$

By (B.2) and (2.1a) and the fact that $\mathsf{F} = \mathsf{F}^\top$, we have

$$
\begin{aligned}
\int_{\mathcal{K}} (\mathsf{F} \otimes I_n)^\top \varpi(\tau) F(\tau) U F^\top(\tau) (\mathsf{F} \otimes I_n) \mathrm{d}\tau = (\mathsf{F} \otimes I_n) \int_{\mathcal{K}} \varpi(\tau) F(\tau) F^\top(\tau) \mathrm{d}\tau (\mathsf{F} \otimes U) \\
(\mathsf{F} \otimes I_n) \left[ \int_{\mathcal{K}} \varpi(\tau) \boldsymbol{f}(\tau) \boldsymbol{f}^\top(\tau) \mathrm{d}\tau \otimes I_n \right] (\mathsf{F} \otimes U) = \mathsf{F} \otimes U \quad \text{(B.4)}
\end{aligned}
$$

where $F(\tau) = \boldsymbol{f}(\tau) \otimes I_n$ and $\mathsf{F}^{-1} = \int_{\mathcal{K}} \varpi(\tau) \boldsymbol{f}(\tau) \boldsymbol{f}^\top(\tau) \mathrm{d}\tau$. Substituting the result in (B.4) for the last integral in (B.1) and also considering the relation in (B.3) yields

$$
\begin{aligned}
\int_{\mathcal{K}} \varpi(\tau) \boldsymbol{\varepsilon}^\top(\tau) U \boldsymbol{\varepsilon}(\tau) \mathrm{d}\tau = \int_{\mathcal{K}} \varpi(\tau) \boldsymbol{x}^\top(\tau) U \boldsymbol{x}(\tau) \mathrm{d}\tau \\
- \int_{\mathcal{K}} \varpi(\tau) \boldsymbol{x}^\top(\tau) F^\top(\tau) \mathrm{d}\tau \, (\mathsf{F} \otimes U) \int_{\mathcal{K}} \varpi(\tau) F(\tau) \boldsymbol{x}(\tau) \mathrm{d}\tau. \quad \text{(B.5)}
\end{aligned}
$$

Given $U \succ 0$, (B.5) gives (4.3). This finishes the proof. ∎

# Appendix C

# Proof of Theorem 4.4

*Proof.* Let $U \succ 0$ and $\varpi(\cdot)$ and $\boldsymbol{f}(\cdot)$ satisfying (4.2) be given which gives $\mathsf{F}^{-1} = \int_{\mathcal{K}} \varpi(\tau)\boldsymbol{f}(\tau)\boldsymbol{f}^{\top}(\tau)\mathrm{d}\tau \in \mathbb{S}^{d}_{\succ 0}$. Using the Schur complement to (4.22) with $U \succ 0$ concludes that $Y \succeq X^{\top}U^{-1}X$ for any $Y$, $X = \mathbf{Row}_{i=1}^{d} X_i \in \mathbb{R}^{n\times \rho d n}$ satisfying (4.22). Now consider $\mathsf{W}$ in Theorem 4.3 with $Y \succeq X^{\top}U^{-1}X$ and (4.21), we have

$$
\begin{aligned}
\mathsf{W} \succeq \int_{\mathcal{K}} \varpi(\tau)(\boldsymbol{f}^{\top}(\tau)\otimes I_{nd})X^{\top}U^{-1}X(\boldsymbol{f}(\tau)\otimes I_{nd})\mathrm{d}\tau &= \widehat{X}^{\top}\int_{\mathcal{K}} \varpi(\tau)(\boldsymbol{f}(\tau)\otimes I_n)U^{-1}\left(\boldsymbol{f}^{\top}(\tau)\otimes I_n\right)\mathrm{d}\tau\widehat{X} \\
&= \widehat{X}^{\top}\left(\int_{\mathcal{K}} \varpi(\tau)\left(I_d\otimes U^{-1}\right)(\boldsymbol{f}(\tau)\boldsymbol{f}^{\top}(\tau)\otimes I_n)\mathrm{d}\tau\right)\widehat{X} = \widehat{X}^{\top}\left(\mathsf{F}^{-1}\otimes U^{-1}\right)\widehat{X} \quad \text{(C.1)}
\end{aligned}
$$

with $\widehat{X} = \mathbf{Col}_{i=1}^{d} X_i \in \mathbb{R}^{dn\times \rho n}$. By the structures in (C.1), one can also infer that $\mathsf{W}$ in Theorem 4.3 satisfies

$$
\mathsf{W} = \widehat{X}^{\top}\left(\mathsf{F}^{-1}\otimes U^{-1}\right)\widehat{X} \quad \text{(C.2)}
$$

with

$$
\widehat{X} = (\mathsf{F}\otimes U)\,\Upsilon, \quad Y = X^{\top}U^{-1}X \quad \text{(C.3)}
$$

for a given $U \succ 0$, where the values of $X$ for $Y$ in (C.3) can be determined by the structural relation $X = \mathbf{Row}_{i=1}^{d} X_i \in \mathbb{R}^{n\times \rho d n}$ with $\widehat{X} = (\mathsf{F}\otimes U)\,\Upsilon = \mathbf{Col}_{i=1}^{d} X_i \in \mathbb{R}^{dn\times \rho n}$.

By (C.1)–(C.3) with (4.14), we have

$$
\mathsf{Sy}\left(\Upsilon^{\top}\widehat{X}\right) - \mathsf{W} \preceq \mathsf{Sy}\left(\Upsilon^{\top}\widehat{X}\right) - \widehat{X}^{\top}\left(\mathsf{F}^{-1}\otimes U^{-1}\right)\widehat{X} \preceq \Upsilon^{\top}\left(\mathsf{F}\otimes U\right)\Upsilon \quad \text{(C.4)}
$$

holds for any $Y$ and $X = \mathbf{Row}_{i=1}^{d} X_i \in \mathbb{R}^{n\times \rho d n}$ satisfying (4.22), and

$$
\mathsf{Sy}\left(\Upsilon^{\top}\widehat{X}\right) - \mathsf{W} = \mathsf{Sy}\left(\Upsilon^{\top}\widehat{X}\right) - \widehat{X}^{\top}\left(\mathsf{F}^{-1}\otimes U^{-1}\right)\widehat{X} = \Upsilon^{\top}\left(\mathsf{F}\otimes U\right)\Upsilon \quad \text{(C.5)}
$$

if $Y$ and $X$ satisfy the equalities in (C.3). Given $\Upsilon\boldsymbol{z} = \int_{\mathcal{K}} \varpi(\tau)F(\tau)\boldsymbol{x}(\tau)\mathrm{d}\tau$ in Theorem 4.3 in light of the results in (C.4) and (C.5), it follows that

$$
\boldsymbol{z}^{\top}\left[\mathsf{Sy}\left(\Upsilon^{\top}\widehat{X}\right) - \mathsf{W}\right]\boldsymbol{z} \leq \boldsymbol{z}^{\top}\Upsilon^{\top}(\mathsf{F}\otimes U)\,\Upsilon\boldsymbol{z} = \int_{\mathcal{K}} \varpi(\tau)\boldsymbol{x}^{\top}(\tau)F^{\top}(\tau)\mathrm{d}\tau\,(\mathsf{F}\otimes U)\int_{\mathcal{K}} \varpi(\tau)F(\tau)\boldsymbol{x}(\tau)\mathrm{d}\tau \quad \text{(C.6)}
$$

holds for any $Y$ and $X$ satisfying (4.22) with $U \succ 0$, and

$$
\boldsymbol{z}^{\top}\left[\mathsf{Sy}\left(\Upsilon^{\top}\widehat{X}\right) - \mathsf{W}\right]\boldsymbol{z} = \boldsymbol{z}^{\top}\Upsilon^{\top}(\mathsf{F}\otimes U)\,\Upsilon\boldsymbol{z} = \int_{\mathcal{K}} \varpi(\tau)\boldsymbol{x}^{\top}(\tau)F^{\top}(\tau)\mathrm{d}\tau\,(\mathsf{F}\otimes U)\int_{\mathcal{K}} \varpi(\tau)F(\tau)\boldsymbol{x}(\tau)\mathrm{d}\tau \quad \text{(C.7)}
$$

holds with (C.3).

In summary, the aforementioned arguments demonstrate that, under the same $\varpi(\cdot)$, $U$ and $\boldsymbol{f}(\cdot)$, one can always find $X$ and $Y$ for (4.22) such that (4.23) becomes identical to (4.3), in which case it corresponds to the maximum LB value in (4.23). Since the maximum LB value in (4.17) is also identical to the one in (4.3), this completes the proof of this theorem. ∎

# Appendix D

# Proof of Lemma 5.2

*Proof.* The proof of Lemma 5.2 is inspired by the proofs of [470, Lemma 2] and [169, Lemma 5]. Firstly, it can be inferred that $\mathcal{E}_d$ in (5.25) is invertible for any $\mathsf{f}(\cdot) \in \mathcal{L}^2_{\varpi}(\mathcal{K};\mathbb{R}^d)$ , $\mathsf{g}(\cdot) \in \mathcal{L}^2_{\varpi}(\mathcal{K};\mathbb{R}^\delta)$ satisfying (5.23) since

$$\mathcal{E}_d = \int_{\mathcal{K}} \varpi(\tau)\mathbf{e}(\tau)\mathbf{e}^\top(\tau)\mathrm{d}\tau = \begin{bmatrix} I_\delta & -\mathsf{A}\end{bmatrix}\int_{\mathcal{K}}\varpi(\tau)\begin{bmatrix}\mathsf{g}(\tau)\\ \mathsf{f}(\tau)\end{bmatrix}\begin{bmatrix}\mathsf{g}^\top(\tau) & \mathsf{f}^\top(\tau)\end{bmatrix}\mathrm{d}\tau\begin{bmatrix} I_\delta & -\mathsf{A}\end{bmatrix}^\top \succ 0, \tag{D.1}$$

where the positive definite matrix inequality can be derived based on (5.23) and the properties of congruent transformations, taking into account that rank $\begin{bmatrix} I_\delta & -\mathsf{A}\end{bmatrix} = \delta$. For this reason, $\mathcal{E}_d^{-1}$ is always well defined.

Let $\boldsymbol{v}(\tau) := \boldsymbol{x}(\tau) - \mathsf{F}^\top(\tau)(\mathcal{F}_d \otimes I_n)\int_{\mathcal{K}}\varpi(\theta)\mathsf{F}(\theta)\boldsymbol{x}(\theta)\mathrm{d}\theta - \mathsf{E}^\top(\tau)\left(\mathcal{E}_d^{-1}\otimes I_n\right)\int_{\mathcal{K}}\varpi(\theta)\mathsf{E}(\theta)\boldsymbol{x}(\theta)\mathrm{d}\theta$, where $\mathsf{F}(\cdot)$, $\mathsf{E}(\cdot)$ are defined Lemma 5.2. By $\mathsf{A} = \int_{\mathcal{K}}\varpi(\tau)\mathsf{g}(\tau)\mathsf{f}^\top(\tau)\mathrm{d}\tau\mathcal{F}_d$ and $\mathbf{e}(\tau) = \mathsf{g}(\tau) - \mathsf{A}\mathsf{f}(\tau) \in \mathbb{R}^\delta$, we see that

$$\begin{aligned}\int_{\mathcal{K}}\varpi(\tau)\mathbf{e}(\tau)\mathsf{f}^\top(\tau)\mathrm{d}\tau &= \int_{\mathcal{K}}\varpi(\tau)\left[\mathsf{g}(\tau) - \mathsf{A}\mathsf{f}(\tau)\right]\mathsf{f}^\top(\tau)\mathrm{d}\tau = \int_{\mathcal{K}}\varpi(\tau)\mathsf{g}(\tau)\mathsf{f}^\top(\tau)\mathrm{d}\tau - \mathsf{A}\int_{\mathcal{K}}\varpi(\tau)\mathsf{f}(\tau)\mathsf{f}^\top(\tau)\mathrm{d}\tau\\ &= \int_{\mathcal{K}}\varpi(\tau)\mathsf{g}(\tau)\mathsf{f}^\top(\tau)\mathrm{d}\tau - \left(\int_{\mathcal{K}}\varpi(\tau)\mathsf{g}(\tau)\mathsf{f}^\top(\tau)\mathrm{d}\tau\right)\mathcal{F}_d\mathcal{F}_d^{-1} = \mathsf{O}_{\delta,d}.\end{aligned} \tag{D.2}$$

By combining the expression of $\boldsymbol{v}(\cdot)$ with $\int_{\mathcal{K}}\varpi(\tau)\boldsymbol{v}^\top(\tau)U\boldsymbol{v}(\tau)\mathrm{d}\tau$ and considering (D.2), it yields

$$\begin{aligned}\int_{\mathcal{K}}\varpi(\tau)\boldsymbol{v}^\top(\tau)U\boldsymbol{v}(\tau)\mathrm{d}\tau &= \int_{\mathcal{K}}\varpi(\tau)\boldsymbol{x}^\top(\tau)U\boldsymbol{x}(\tau)\mathrm{d}\tau - 2\int_{\mathcal{K}}\varpi(\tau)\boldsymbol{x}^\top(\tau)U\mathsf{F}^\top(\tau)\mathrm{d}\tau(\mathcal{F}_d\otimes I_n)\boldsymbol{\zeta}\\ &+ \boldsymbol{\zeta}^\top\int_{\mathcal{K}}\varpi(\tau)(\mathcal{F}_d\otimes I_n)^\top\mathsf{F}(\tau)U\mathsf{F}^\top(\tau)(\mathcal{F}_d\otimes I_n)\mathrm{d}\tau\boldsymbol{\zeta} - 2\int_{\mathcal{K}}\varpi(\tau)\boldsymbol{x}^\top(\tau)U\mathsf{E}^\top(\tau)\mathrm{d}\tau(\mathcal{E}_d^{-1}\otimes I_n)\boldsymbol{\omega}\\ &+ \boldsymbol{\omega}^\top\int_{\mathcal{K}}\varpi(\tau)(\mathcal{E}_d^{-1}\otimes I_n)^\top\mathsf{E}(\tau)U\mathsf{E}^\top(\tau)(\mathcal{E}_d^{-1}\otimes I_n)\mathrm{d}\tau\boldsymbol{\omega},\end{aligned} \tag{D.3}$$

where $\boldsymbol{\zeta} := \int_{\mathcal{K}}\varpi(\theta)\mathsf{F}(\theta)\boldsymbol{x}(\theta)\mathrm{d}\theta$ and $\boldsymbol{\omega} := \int_{\mathcal{K}}\varpi(\theta)\mathsf{E}(\theta)\boldsymbol{x}(\theta)\mathrm{d}\theta$. Apply (2.1a) to $U\mathsf{F}^\top(\tau)$ and $U\mathsf{E}^\top(\tau)$ and consider $\mathsf{F}(\tau) = \mathsf{f}(\tau)\otimes I_n$ and $\mathsf{E}(\tau) = \mathbf{e}(\tau)\otimes I_n$. Then we see that

$$U\mathsf{F}^\top(\tau) = \mathsf{F}^\top(\tau)(I_d\otimes U), \quad U\mathsf{E}^\top(\tau) = \mathsf{E}^\top(\tau)(I_\delta\otimes U) \tag{D.4}$$

given that $(X\otimes Y)^\top = X^\top\otimes Y^\top$. Additionally, it is true that

$$\begin{aligned}\int_{\mathcal{K}}\varpi(\tau)\mathsf{F}(\tau)\mathsf{F}^\top(\tau)\mathrm{d}\tau &= \left(\int_{\mathcal{K}}\varpi(\tau)\mathsf{f}(\tau)\mathsf{f}^\top(\tau)\mathrm{d}\tau\right)\otimes I_n = \mathcal{F}_d^{-1}\otimes I_n\\ \int_{\mathcal{K}}\varpi(\tau)\mathsf{E}(\tau)\mathsf{E}^\top(\tau)\mathrm{d}\tau &= \left(\int_{\mathcal{K}}\varpi(\tau)\mathbf{e}(\tau)\mathbf{e}^\top(\tau)\mathrm{d}\tau\right)\otimes I_n = \mathcal{E}_d\otimes I_n\end{aligned} \tag{D.5}$$

since $\mathsf{F}(\tau) = \mathbf{f}(\tau) \otimes I_n$ and $\mathsf{E}(\tau) = \boldsymbol{e}(\tau) \otimes I_n$. Utilizing (D.4) and (D.5) with (2.1a) on some of the terms in (D.3), it follows that

$$\begin{aligned}
&\int_{\mathcal{K}} \varpi(\tau)\boldsymbol{x}^\top(\tau)U\mathsf{F}^\top(\tau)\mathrm{d}\tau(\mathcal{F}_d \otimes I_n)\boldsymbol{\zeta} = \boldsymbol{\zeta}^\top(\mathcal{F}_d \otimes U)\boldsymbol{\zeta}, \\
&\int_{\mathcal{K}} \varpi(\tau)\boldsymbol{x}^\top(\tau)U\mathsf{E}^\top(\tau)\mathrm{d}\tau(\mathcal{E}_d^{-1} \otimes I_n)\boldsymbol{\omega} = \boldsymbol{\omega}^\top(\mathcal{E}_d^{-1} \otimes U)\boldsymbol{\omega},
\end{aligned} \tag{D.6}$$

and

$$\begin{aligned}
&\int_{\mathcal{K}} (\mathcal{F}_d \otimes I_n)^\top \varpi(\tau)\mathsf{F}(\tau)U\mathsf{F}^\top(\tau)(\mathcal{F}_d \otimes I_n)\mathrm{d}\tau = (\mathcal{F}_d \otimes I_n)\int_{\mathcal{K}} \varpi(\tau)\mathsf{F}(\tau)\mathsf{F}^\top(\tau)\mathrm{d}\tau(\mathcal{F}_d \otimes U) = \mathcal{F}_d \otimes U \\
&\int_{\mathcal{K}} (\mathcal{E}_d^{-1} \otimes I_n)^\top \varpi(\tau)\mathsf{E}(\tau)U\mathsf{E}^\top(\tau)(\mathcal{E}_d^{-1} \otimes I_n)\mathrm{d}\tau = (\mathcal{E}_d^{-1} \otimes I_n)\int_{\mathcal{K}} \varpi(\tau)\mathsf{E}(\tau)\mathsf{E}^\top(\tau)\mathrm{d}\tau\left(\mathcal{E}_d^{-1} \otimes U\right) \\
&\qquad = \mathcal{E}_d^{-1} \otimes U.
\end{aligned} \tag{D.7}$$

Substituting the result in (D.7) for the last integral in (D.3) and also considering (D.6), we obtain

$$\begin{aligned}
\int_{\mathcal{K}} \varpi(\tau)\boldsymbol{v}^\top(\tau)U\boldsymbol{v}(\tau)\mathrm{d}\tau = &\int_{\mathcal{K}} \varpi(\tau)\boldsymbol{x}^\top(\tau)U\boldsymbol{x}(\tau)\mathrm{d}\tau - \int_{\mathcal{K}} \varpi(\tau)\boldsymbol{x}^\top(\tau)\mathsf{F}^\top(\tau)\mathrm{d}\tau\left(\mathcal{F}_d \otimes U\right)\int_{\mathcal{K}} \varpi(\tau)\mathsf{F}(\tau)\boldsymbol{x}(\tau)\mathrm{d}\tau \\
&- \int_{\mathcal{K}} \varpi(\tau)\boldsymbol{x}^\top(\tau)\mathsf{E}^\top(\tau)\mathrm{d}\tau\left(\mathcal{E}_d^{-1} \otimes U\right)\int_{\mathcal{K}} \varpi(\tau)\mathsf{E}(\tau)\boldsymbol{x}(\tau)\mathrm{d}\tau.
\end{aligned} \tag{D.8}$$

Since $U \succeq 0$, (D.8) gives (5.24). This finishes the proof. ■

# Appendix E

# Lyapunov-Krasovskiĭ Stability Theorem for FDEs Subject to the Carathéodory Conditions

**Theorem E.1.** *Consider an FDE in the form of*

$$\begin{aligned} &\widetilde{\forall} t \in \mathcal{T},\ \ \dot{\boldsymbol{x}}(t) = \boldsymbol{f}(t, \mathbf{x}_t(\cdot)),\ \ \mathcal{T} = [t_0, +\infty) \cap \mathcal{U},\ \ t_0 \in \mathcal{U} \subseteq \mathbb{R} \\ &\forall \theta \in [-r, 0],\ \ \boldsymbol{x}(t_0+\theta) = \mathbf{x}_{t_0}(\theta) = \boldsymbol{\phi}(\theta),\ \ r>0, \\ &\forall t \in \mathbb{R},\ \mathbf{0}_n = \boldsymbol{f}(t, \mathbf{0}_n(\cdot)) \end{aligned} \tag{E.1}$$

*where* $\boldsymbol{f} : \mathbb{R} \times \mathcal{C}([-r,0];\mathbb{R}^n) \to \mathbb{R}^n$ *satisfies the Carathéodory conditions in [9, section 2.6] and*

$$\exists c(\cdot) \in \mathbb{R}_{>0}^{\mathbb{R}_{>0}},\ \forall\delta>0,\ \forall \boldsymbol{\phi}(\cdot) \in \mathcal{C}_\delta([-r,0];\mathbb{R}^n),\ \widetilde{\forall} t\in\mathbb{R},\ \ \left\|\boldsymbol{f}\left(t,\boldsymbol{\phi}(\cdot)\right)\right\|_1 < c(\delta). \tag{E.2}$$

*Then the trivial solution* $\boldsymbol{x}(t) \equiv \mathbf{0}_n$ *of* (E.1) *is globally uniformly asymptotically stable if there exist three functions* $\alpha_1(\cdot); \alpha_2(\cdot); \alpha_3(\cdot) \in \mathcal{K}_\infty$, *and a differentiable functional* $\mathsf{v} : \mathbb{R} \times \mathcal{C}([-r,0];\mathbb{R}^n) \to \mathbb{R}$ *such that*

$$\forall t \in \mathbb{R},\ \forall \boldsymbol{\phi}(\cdot) \in \mathcal{C}([-r,0];\mathbb{R}^n),\ \ \alpha_1\left(\|\boldsymbol{\phi}(0)\|_2\right) \le \mathsf{v}(t,\boldsymbol{\phi}(\cdot)) \le \alpha_2\left(\|\boldsymbol{\phi}(\cdot)\|_\infty\right), \tag{E.3}$$

$$\widetilde{\forall} t \ge t_0 \in \mathbb{R},\ \frac{\mathsf{d}}{\mathsf{d}t}\mathsf{v}(t,\mathbf{x}_t(\cdot)) \le -\alpha_3\left(\|\boldsymbol{x}(t)\|_2\right) \tag{E.4}$$

*where* $\|\boldsymbol{\phi}(\cdot)\|_\infty^2 := \max_{-r_2 \le \tau \le 0} \|\boldsymbol{\phi}(\tau)\|_2^2$, *and* $\mathbf{x}_t(\cdot)$, $\boldsymbol{x}(\cdot)$ *in* (E.4) *satisfy* $\dot{\boldsymbol{x}}(t) = \boldsymbol{f}(t,\mathbf{x}_t(\cdot))$ *in* (E.1) *for almost all* $t \ge t_0$. *Moreover,* $\mathcal{K}_\infty$ *follows the standard definition in [695]. Note that the notation* $\widetilde{\forall}$ *denotes for almost all with respect to the Lebesgue measure [571].*

*Proof.* The proof here is based on the procedure in [9, Theorem 2.1] and [178, Theorem 1.3] with certain variations. To establish the uniform stability of the trivial solution, let

$$\mathbb{R}_{\ge 0} \ni \epsilon \mapsto \delta(\epsilon) = 1/2 \min\left(\epsilon, \alpha_2^{-1}\left(\alpha_1(\epsilon)\right)\right) \tag{E.5}$$

where $\alpha_2^{-1}(\cdot)$ is well defined since $\alpha_2(\cdot) \in \mathcal{K}_\infty$. It is evident that $\delta(\cdot) \in \mathcal{K}_\infty$ and satisfies $\forall \epsilon > 0$, $0 < \delta(\epsilon) < \epsilon$ and $\delta(\epsilon) < \alpha_2^{-1}\left(\alpha_1(\epsilon)\right)$ that further implies

$$\forall \epsilon > 0,\ \ \alpha_2\left(\delta(\epsilon)\right) < \alpha_1(\epsilon) \tag{E.6}$$

since $\alpha_2(\cdot) \in \mathcal{K}_\infty$. By (E.4), it holds true that $\widetilde{\forall} t \ge t_0 \in \mathbb{R}$, $\dot{\mathsf{v}}(t,\mathbf{x}_t(\cdot)) \le 0$. Now, by applying the Fundamental Theorem of Calculus for Lebesgue integrals [571] to the previous proposition, we see that

$$\forall t_0 \in \mathbb{R}, \forall t \geq t_0, \forall \boldsymbol{\phi}(\cdot) \in \mathcal{C}([-r,0];\mathbb{R}^n), \int_{t_0}^{t} \dot{\mathsf{v}}(t, \mathbf{x}_t(\cdot)) \mathrm{d}\tau$$
$$= \mathsf{v}\left(t, \mathbf{x}_t(\cdot)\right) - \mathsf{v}\big(t_0, \mathbf{x}_{t_0}(\cdot)\big) = \mathsf{v}\left(t, \mathbf{x}_t(\cdot)\right) - \mathsf{v}\big(t_0, \boldsymbol{\phi}(\cdot)\big) \leq 0 \quad \text{(E.7)}$$

which further implies that $\forall t_0 \in \mathbb{R}$, $\forall t \geq t_0$, $\forall \epsilon > 0$, $\forall \boldsymbol{\phi}(\cdot) \in \mathcal{C}_{\delta(\epsilon)}([-r,0];\mathbb{R}^n)$:

$$\alpha_1(\|\boldsymbol{x}(t)\|_2) \leq \mathsf{v}(t, \mathbf{x}_t(\cdot)) \leq \mathsf{v}(t_0, \boldsymbol{\phi}(\cdot)) \leq \alpha_2\left(\|\boldsymbol{\phi}(\cdot)\|_\infty\right) < \alpha_2(\delta(\epsilon)) < \alpha_1(\epsilon) \quad \text{(E.8)}$$

by (E.3) and (E.6), where $\mathcal{C}_{\delta(\epsilon)}([-r,0];\mathbb{R}^n) := \left\{\boldsymbol{\phi}(\cdot) \in \mathcal{C}([-r,0];\mathbb{R}^n) : \|\boldsymbol{\phi}(\cdot)\|_\infty < \delta(\epsilon)\right\}$. Therefore,

$$\forall \epsilon > 0,\ \forall \boldsymbol{\phi}(\cdot) \in \mathcal{C}_{\delta(\epsilon)}([-r,0];\mathbb{R}^n),\ \forall t_0 \in \mathbb{R},\ \forall t \geq t_0,\ \|\boldsymbol{x}(t)\|_2 < \epsilon, \ \text{ since } \ \alpha_1(\cdot) \in \mathcal{K}_\infty \quad \text{(E.9)}$$

where $\delta(\epsilon) = 1/2 \min\left(\epsilon, \alpha_2^{-1}\left(\alpha_1(\epsilon)\right)\right)$ is independent of $t_0 \in \mathbb{R}$ and $\lim_{\epsilon \to +\infty} \delta(\epsilon) = +\infty$ since $\delta(\cdot) \in \mathcal{K}_\infty$. Now one can further infer

$$\forall \epsilon > 0,\ \exists \delta > 0,\ \forall \boldsymbol{\phi}(\cdot) \in \mathcal{C}_\delta([-r,0];\mathbb{R}^n),\ \forall t_0 \in \mathbb{R},\ \forall t \geq t_0,\ \|\mathbf{x}_t(\cdot)\|_\infty \leq \max_{\tau \geq t_0} \|\boldsymbol{x}(\tau)\| < \epsilon \quad \text{(E.10)}$$

form (E.9) which demonstrates uniform stability.

To prove global uniform asymptotic stability, we shall employ proof by contradiction. Note that the origin is globally uniform asymptotically stable if it is uniformly stable, as established earlier, and

$$\forall \eta > 0, \forall \delta > 0,\ \exists \beta \geq 0,\ \forall \boldsymbol{\phi}(\cdot) \in \mathcal{C}_\delta([-r,0];\mathbb{R}^n),\ \forall t_0 \in \mathbb{R},\ \forall t \geq t_0 + \beta,\ \|\mathbf{x}_t(\cdot)\|_\infty < \eta. \quad \text{(E.11)}$$

Now let us assume that

$$\exists \epsilon > 0,\ \exists \delta > 0,\ \exists \boldsymbol{\phi}(\cdot) \in \mathcal{C}_\delta([-r,0];\mathbb{R}^n),\ \exists t_0 \in \mathbb{R},\ \forall t \geq t_0,\ \|\mathbf{x}_t(\cdot)\|_\infty \geq \epsilon. \quad \text{(E.12)}$$

Given the definition $\|\mathbf{x}_t(\cdot)\|_\infty := \max_{\tau \in [-r,0]} \|\boldsymbol{x}(t+\tau)\|_2$ with (E.12), it implies

$$\exists \epsilon > 0,\ \exists \delta > 0,\ \exists \boldsymbol{\phi}(\cdot) \in \mathcal{C}_\delta([-r,0];\mathbb{R}^n),\ \exists t_0 \in \mathbb{R},\ \forall t \geq t_0,\ \exists \lambda \in [t-r, t],\ \|\boldsymbol{x}(\lambda)\|_2 \geq \epsilon. \quad \text{(E.13)}$$

Let $\epsilon > 0$, $\delta > 0$, $\boldsymbol{\phi}(\cdot) \in \mathcal{C}_\delta([-r,0];\mathbb{R}^n)$ and $t_0 \in \mathbb{R}$ in (E.13) be given, then we see that there exists a sequence $\mathbb{N} \ni k \to t_k \in \mathbb{R}_{\geq t_0}$ such that

$$\forall k \in \mathbb{N},\ (2k-1)r \leq t_k - t_0 \leq 2kr \ \ \& \ \ \|\boldsymbol{x}(t_k)\|_2 \geq \epsilon. \quad \text{(E.14)}$$

Additionally,

$$\begin{aligned}
\|\boldsymbol{x}(t)\|_2 = \left\|\boldsymbol{x}(t_k) + \int_{t_k}^{t} \dot{\boldsymbol{x}}(\tau)\mathrm{d}\tau\right\|_2 &\geq \|\boldsymbol{x}(t_k)\|_2 - \left\|\int_{t_k}^{t} \dot{\boldsymbol{x}}(\tau)\mathrm{d}\tau\right\|_2 = \|\boldsymbol{x}(t_k)\|_2 - \left\|\int_{t_k}^{t} \boldsymbol{f}(\tau, \mathbf{x}_\tau(\cdot))\mathrm{d}\tau\right\|_2 \\
&\geq \|\boldsymbol{x}(t_k)\|_2 - \left\|\int_{t_k}^{t} \boldsymbol{f}(\tau, \mathbf{x}_\tau(\cdot))\mathrm{d}\tau\right\|_1 = \|\boldsymbol{x}(t_k)\|_2 - \sum_{i=1}^{n} \left|\int_{t_k}^{t} f_i(\tau, \mathbf{x}_\tau(\cdot))\mathrm{d}\tau\right| \\
&\geq \|\boldsymbol{x}(t_k)\|_2 - \left|\sum_{i=1}^{n} \int_{t_k}^{t} |f_i(\tau, \mathbf{x}_\tau(\cdot))|\,\mathrm{d}\tau\right| = \|\boldsymbol{x}(t_k)\|_2 - \left|\int_{t_k}^{t} \sum_{i=1}^{n} |f_i(\tau, \mathbf{x}_\tau(\cdot))|\,\mathrm{d}\tau\right| \\
&\geq \|\boldsymbol{x}(t_k)\|_2 - \left|\int_{t_k}^{t} \|\boldsymbol{f}(\tau, \mathbf{x}_\tau(\cdot))\|_1\,\mathrm{d}\tau\right| \quad \text{(E.15)}
\end{aligned}$$

holds true for all $t \geq t_0$ and $k \in \mathbb{N}$ based on the properties of Lebesgue integrals and norms. Since $\forall t \geq t_0$, $\forall k \in \mathbb{N}$, $\left|\int_{t_k}^{t} \|\boldsymbol{f}(\tau, \mathbf{x}_\tau(\cdot))\|_1\,\mathrm{d}\tau\right| < \left|\int_{t_k}^{t} c(\delta)\ \mathrm{d}\tau\right| = c(\delta)|t - t_k|$ by (E.2) with a given $\delta > 0$ and $\boldsymbol{\phi}(\cdot) \in \mathcal{C}_\delta([-r,0];\mathbb{R}^n)$, it follows that

$$
\begin{aligned}
\forall k \in \mathbb{N},\ \ \forall t \in \mathcal{T}_k := \left[t_k - \frac{\epsilon}{2c(\delta)}, t_k + \frac{\epsilon}{2c(\delta)}\right],\ \ &\|\boldsymbol{x}(t)\|_2 \geq \|\boldsymbol{x}(t_k)\|_2 - \left|\int_{t_k}^{t} \|\boldsymbol{f}(\tau, \mathbf{x}_\tau(\cdot))\|_1 \,\mathsf{d}\tau\right| \\
> \|\boldsymbol{x}(t_k)\|_2 - \left|\int_{t_k}^{t} c(\delta)\mathsf{d}\tau\right| = \|\boldsymbol{x}(t_k\|_2 - c(\delta)|t - t_k| &\geq \epsilon - c(\delta)\frac{\epsilon}{2c(\delta)} = \frac{\epsilon}{2}. \qquad \text{(E.16)}
\end{aligned}
$$

Consequently,

$$
\widetilde{\forall} t \in \mathbb{R}_{\geq t_0} \cap \bigcup_{k\in\mathbb{N}} \mathcal{T}_k,\ \frac{\mathsf{d}}{\mathsf{d}t}\mathsf{v}(\mathbf{x}_t(\cdot)) \leq -\alpha_3\left(\epsilon/2\right). \quad \& \quad \widetilde{\forall} t \in \mathbb{R}_{\geq t_0},\ \frac{\mathsf{d}}{\mathsf{d}t}\mathsf{v}(\mathbf{x}_t(\cdot)) \leq 0. \tag{E.17}
$$

Since $c(\delta) > 0$ in $\mathcal{T}_k = [t_k - \epsilon/2c(\delta), t_k + \epsilon/2c(\delta)]$ can be made arbitrarily large for any $\delta > 0$, we can assume that $\bigcap_{k\in\mathbb{N}} \mathcal{T}_k = \varnothing$ and $t_1 - \epsilon/2c(\delta) \geq t_0$. As a result, we have

$$
\begin{aligned}
&\forall k \in \mathbb{N},\ \mathsf{v}(t_k, \mathbf{x}_{t_k}(\cdot)) - \mathsf{v}(t_0, \boldsymbol{\phi}(\cdot)) \\
&\quad= \int_{t_0}^{t_k} \frac{\mathsf{d}}{\mathsf{d}\tau}\mathsf{v}(\mathbf{x}_\tau(\cdot))\mathsf{d}\tau = \int_{\bigcup_{i=1}^{k-1}\mathcal{T}_i} \frac{\mathsf{d}}{\mathsf{d}\tau}\mathsf{v}(\mathbf{x}_\tau(\cdot))\mathsf{d}\tau + \underbrace{\int_{[t_k,t_0]\setminus\bigcup_{i=1}^{k-1}\mathcal{T}_i} \frac{\mathsf{d}}{\mathsf{d}\tau}\mathsf{v}(\mathbf{x}_\tau(\cdot))\mathsf{d}\tau}_{0} \\
&\qquad\leq -\int_{\bigcup_{i=1}^{k-1}\mathcal{T}_i} \alpha_3\left(\epsilon/2\right)\mathsf{d}\tau = -\sum_{i=1}^{k-1}\int_{\mathcal{T}_i} \alpha_3\left(\epsilon/2\right)\mathsf{d}\tau = -\alpha_3\left(\epsilon/2\right)\frac{\epsilon}{c(\delta)}(k-1) \qquad \text{(E.18)}
\end{aligned}
$$

by (E.17). This further indicates that

$$
\begin{aligned}
\forall k \in \mathbb{N},\ \mathsf{v}(t_k, \mathbf{x}_{t_k}(\cdot)) \leq \mathsf{v}(t_0, \boldsymbol{\phi}(\cdot)) - \alpha_3\left(\epsilon/2\right)\frac{\epsilon}{c(\delta)}(k-1) &\leq \alpha_2\left(\|\boldsymbol{\phi}(\cdot)\|_\infty\right) - \alpha_3\left(\epsilon/2\right)\frac{\epsilon}{c(\delta)}(k-1) \\
&< \alpha_2(\delta) - \alpha_3\left(\epsilon/2\right)\frac{\epsilon}{c(\delta)}(k-1) \qquad \text{(E.19)}
\end{aligned}
$$

by (E.3) and the fact that $\|\boldsymbol{\phi}(\cdot)\|_\infty < \delta$ and $\alpha_2(\cdot) \in \mathcal{K}_\infty$. Now it is easy to see that

$$
\alpha_2(\delta) - \alpha_3\left(\epsilon/2\right)\frac{\epsilon}{c(\delta)}(k-1) < 0 \iff \frac{\alpha_2(\delta)}{\alpha_3\left(\epsilon/2\right)}\frac{c(\delta)}{\epsilon} + 1 < k. \tag{E.20}
$$

Let $\kappa(\epsilon, \delta) = \left\lceil \frac{\alpha_2(\delta)}{\alpha_3\left(\epsilon/2\right)}\frac{c(\delta)}{\epsilon} \right\rceil + 1$, then we can obtain $\forall k > \kappa(\epsilon, \delta)$, $\mathsf{v}(t_k, \mathbf{x}_{t_k}(\cdot)) < 0$ by (E.19), which is a contradiction considering (7.39). In consequence, (E.12) cannot be true for $t_k$ with any $k > \kappa(\epsilon, \delta)$, which implies that $\exists k \leq \kappa(\epsilon, \delta)$, $\|\mathbf{x}_{t_k}(\cdot)\|_\infty < \epsilon$. This further indicates that

$$
\forall \epsilon > 0,\ \forall \delta > 0,\ \forall \boldsymbol{\phi}(\cdot) \in \mathcal{C}_\delta([-r, 0]; \mathbb{R}^n),\ \forall t_0 \in \mathbb{R},\ \exists \theta \in [t_0, t_0 + 2r\kappa(\epsilon, \delta)],\ \|\mathbf{x}_\theta(\cdot)\|_\infty < \epsilon \tag{E.21}
$$

considering (E.14).

Now let $\epsilon > 0$ in (E.21) be

$$
\epsilon(\eta) = 1/3 \min\left(\eta, \alpha_2^{-1}\left(\alpha_1(\eta)\right)\right) \tag{E.22}
$$

with a given $\eta > 0$, and assume $\boldsymbol{\phi}(\cdot)$, $t_0$, $\theta$ in (E.21) are also given. Note that the structure of $\epsilon(\cdot)$ in (E.22) guarantees $\epsilon(\cdot) \in \mathcal{K}_\infty$ and $\alpha_2(\epsilon(\eta)) < \alpha_1(\eta)$ for any $\eta > 0$ similar to the property in (E.6). From the computation in (??) and (E.4), we also know that for all $t \geq \theta$ and $\boldsymbol{\phi}(\cdot) \in \mathcal{C}([-r, 0]; \mathbb{R}^n)$

$$
\int_\theta^t \dot{\mathsf{v}}(t, \mathbf{x}_t(\cdot))\mathsf{d}\tau = \mathsf{v}\left(t, \mathbf{x}_t(\cdot)\right) - \mathsf{v}\left(\theta, \mathbf{x}_\theta(\cdot)\right) \leq 0. \tag{E.23}
$$

With (E.3) and (E.21)–(E.23) and $\alpha_2(\epsilon(\eta)) < \alpha_1(\eta)$, we find that

$$\forall \eta > 0, \forall \delta > 0, \quad \forall t_0 \in \mathbb{R}, \quad \forall \boldsymbol{\phi}(\cdot) \in \mathcal{C}_\delta([-r, 0]; \mathbb{R}^n), \quad \forall t \geq \theta \ \text{ such that}$$

$$\alpha_1(\|\boldsymbol{x}(t)\|_2) \leq \mathsf{v}(t, \mathbf{x}_t(\cdot)) \leq \mathsf{v}(\theta, \mathbf{x}_\theta(\cdot)) \leq \alpha_2(\|\mathbf{x}_\theta(\cdot)\|_\infty) < \alpha_2(\epsilon(\eta)) < \alpha_1(\eta). \quad \text{(E.24)}$$

Since $\theta \leq t_0 + 2\kappa(\epsilon, \delta)r$ in (E.21), relation (E.24) further implies

$$\forall t \geq t_0 + 2\kappa(\epsilon, \delta)r \geq \theta, \quad \alpha_1(\|\boldsymbol{x}(t)\|_2) < \alpha_1(\eta)$$

and then $\|\boldsymbol{x}(t)\|_2 < \eta$ since $\alpha_1(\cdot) \in \mathcal{K}_\infty$. Given that $2r\kappa(\epsilon(\eta), \delta)$ is independent of $t_0$, we have proved the global asymptotic stability in (E.11) with $\beta = 2\kappa(\epsilon, \delta)r$. This proves the global uniform asymptotic stability of $\boldsymbol{x}(t) \equiv \mathbf{0}_n$, as the uniform stability has been proved with $\delta(\cdot) \in \mathcal{K}_\infty$ in (E.5) satisfying $\lim_{\epsilon \to +\infty} \delta(\epsilon) = +\infty$. ■